\documentclass[a4paper,12pt,amsart,frenchb]{article}

\usepackage{amsmath,amsbsy,amsfonts,amssymb}
\usepackage[french]{babel}
\newcommand{\ass}[2]{\vskip0.3cm\noindent
{\bf {#1}}. { \sl {#2}}\vskip0.3cm\noindent
}
  
\begin{document}

 \title{   La formule des traces locale tordue}
\author{J.-L. Waldspurger}
\date{13 septembre 2012}
\maketitle

 {\bf Introduction}
 
 On se propose de g\'en\'eraliser au cas tordu les r\'esultats d'Arthur contenus dans les articles [A1] et [A7]. Soient $F$ un corps local, $G$ un groupe r\'eductif connexe d\'efini sur $F$ et $\tilde{G}$ un espace tordu sur $G$, au sens de Labesse (cf. 2.1). Nous imposons une condition \`a $\tilde{G}$ (2.1(2)) qui revient \`a dire qu'il existe un groupe alg\'ebrique non connexe $G^+$ d\'efini sur $F$, de composante neutre $G$, tel que $\tilde{G}$ soit une composante connexe de $G^+$. Mais la structure de groupe sur $G^+$ ne joue aucun r\^ole, seules importent les actions \`a droite et \`a gauche de $G$ sur $\tilde{G}$. Notons $Z_{G}$ le centre de $G$ et $Z_{G}(F)^{\theta}$ le sous-groupe des $z\in Z_{G}(F)$ tels que $z\gamma=\gamma z$ pour tout $\gamma\in \tilde{G}$. On fixe un caract\`ere unitaire $\omega$ de $G(F)$ dont la restriction \`a $Z_{G}(F)^{\theta}$ est triviale. On s'int\'eresse aux "distributions" $\omega$-\'equivariantes sur $\tilde{G}(F)$. Ce sont des formes lin\'eaires $l:C_{c}^{\infty}(\tilde{G}(F))\to {\mathbb C}$ telles que, pour tout $f\in C_{c}^{\infty}(\tilde{G}(F))$ et tout $g\in G(F)$, on ait l'\'egalit\'e $l(^gf)=\omega(g)^{-1}l(f)$, o\`u $^gf$ est la fonction $^gf(\gamma)=f(g^{-1}\gamma g)$. Il y a deux  types basiques de  telles distributions. D'abord les int\'egrales orbitales. On fixe $\gamma\in \tilde{G}(F)$, disons fortement r\'egulier. On note $Z_{G}(\gamma)$ son commutant dans $G$ et on munit le quotient $Z_{G}(\gamma,F)\backslash G(F)$ d'une mesure invariante \`a droite. Pour $f\in C_{c}^{\infty}(\tilde{G}(F))$,  l'int\'egrale orbitale de $f$ au point $\gamma$ est
 $$I_{\tilde{G}}(\gamma,\omega,f)=D^{\tilde{G}}(\gamma)^{1/2}\int_{Z_{G}(\gamma,F)\backslash G(F)}\omega(x)f(x^{-1} \gamma x)\,dx.$$
 La fonction $D^{\tilde{G}}(\gamma)$ est la variante tordue de la fonction habituelle.
 Il y a aussi les caract\`eres de repr\'esentations. Soit $\pi$ une repr\'esentation admissible et irr\'eductible de $G(F)$. Pour $\gamma\in \tilde{G}(F)$, notons $ad_{\gamma}$ l'automorphisme de $G$ tel que $\gamma g=ad_{\gamma}(g)\gamma$ pour tout $g\in G$. La classe d'\'equivalence de la repr\'esentation $\pi\circ ad_{\gamma}$ ne d\'epend pas de $\gamma$. Supposons que $\pi\circ ad_{\gamma}$ soit isomorphe \`a $\omega\otimes \pi$. On peut alors prolonger $\pi$ en une "$\omega$-repr\'esentation" $\tilde{\pi}$ de $\tilde{G}(F)$, c'est-\`a-dire une application $\tilde{\pi}$ de $\tilde{G}(F)$ dans le groupe des automorphismes de l'espace de $\pi$ qui v\'erifie la condition $\tilde{\pi}(g \gamma g')=\pi(g)\tilde{\pi}(\gamma)\pi(g')\omega(g')$ pour tous $\gamma\in \tilde{G}(F)$ et $g,g'\in G(F)$. Pour $f\in \tilde{G}(F)$, on d\'efinit l'op\'erateur $\tilde{\pi}(f)=\int_{\tilde{G}(F)}f(\gamma)\tilde{\pi}(\gamma)\,d\gamma$ (une mesure de Haar \'etant fix\'ee sur $G(F)$ et transport\'ee \`a $\tilde{G}(F)$), puis le caract\`ere
 $$I_{\tilde{G}}(\tilde{\pi},f)=trace(\tilde{\pi}(f)).$$
 
 La formule des traces locale tordue \'etablit une \'egalit\'e entre deux expressions, l'une contenant des int\'egrales orbitales, l'autre des caract\`eres de $\omega$-repr\'esentations temp\'er\'ees. Pr\'ecisons un tout petit peu. Soient $f_{1},f_{2}\in C_{c}^{\infty}(\tilde{G}(F))$. On d\'efinit une expression
 $$J^{\tilde{G}}_{g\acute{e}om}(\omega,f_{1},f_{2})=\sum_{\tilde{M}\in {\cal P}(\tilde{M}_{0})}\vert \tilde{W}^M\vert \vert \tilde{W}^G\vert ^{-1}J_{\tilde{M},g\acute{e}om}^{\tilde{G}}(\omega,f_{1},f_{2}).$$
 L'ensemble ${\cal P}(\tilde{M}_{0})$ est celui des "ensembles de Levi" $\tilde{M}$ de $\tilde{G}$ contenant un ensemble de Levi minimal fix\'e $\tilde{M}_{0}$. Si $\tilde{M}=\tilde{G}$,  $J_{\tilde{G},g\acute{e}om}^{\tilde{G}}(\omega,f_{1},f_{2})$ est une certaine int\'egrale d'expressions
 $$\overline{I_{\tilde{G}}(\gamma,\omega,f_{1})}I_{\tilde{G}}(\gamma,\omega,f_{2})$$
 en des points $\gamma\in \tilde{G}(F)$ qui sont fortement r\'eguliers et elliptiques. Si $\tilde{M}\not=\tilde{G}$, l'expression est plus compliqu\'ee: elle fait intervenir des int\'egrales orbitales pond\'er\'ees, qui g\'en\'eralisent les int\'egrales orbitales d\'efinies ci-dessus, mais ne sont plus $\omega$-\'equivariantes. On d\'efinit aussi une expression
  $$J^{\tilde{G}}_{spec}(\omega,f_{1},f_{2})=\sum_{\tilde{M}\in {\cal P}(\tilde{M}_{0})}\vert \tilde{W}^M\vert \vert \tilde{W}^G\vert ^{-1}J_{\tilde{M},spec}^{\tilde{G}}(\omega,f_{1},f_{2}).$$
  Dans le cas o\`u $\tilde{M}=\tilde{G}$ et o\`u  $Z_{G}(F)^{\theta}$ est compact, le terme $J_{\tilde{G},spec}^{\tilde{G}}(\omega,f_{1},f_{2})$ est une somme de produits
  $$\overline{I_{\tilde{G}}(\tilde{\pi},f_{1})}I_{\tilde{G}}(\tilde{\pi},f_{2}),$$
  o\`u $\tilde{\pi}$ d\'ecrit un certain ensemble de $\omega$-repr\'esentations de $\tilde{G}(F)$. M\^eme dans ce cas simple, la d\'efinition de ces caract\`eres doit \^etre un peu g\'en\'eralis\'ee, car les repr\'esentations sous-jacentes aux $\tilde{\pi}$ ne sont pas irr\'eductibles mais seulement de longueur finie. Dans le cas o\`u $Z_{G}(F)^{\theta}$ n'est plus compact, il faut int\'egrer les produits ci-dessus selon des param\`etres r\'eminiscents de l'existence de ce centre. Dans le cas o\`u $\tilde{M}\not=\tilde{G}$, l'expression fait intervenir des g\'en\'eralisations des caract\`eres, \`a savoir les caract\`eres pond\'er\'es qui, eux non plus, ne sont pas $\omega$-\'equivariants.
  
  La formule des traces locale tordue  (th\'eor\`eme 5.1) affirme l'\'egalit\'e
  $$J^{\tilde{G}}_{g\acute{e}om}(\omega,f_{1},f_{2})=J^{\tilde{G}}_{spec}(\omega,f_{1},f_{2}).$$
  
  La premi\`ere cons\'equence en est le "th\'eor\`eme 0" de Kazhdan: pour $f\in C_{c}^{\infty}(\tilde{G}(F))$, si $f$ v\'erifie $I_{\tilde{G}}(\tilde{\pi},f)=0$ pour toute $\omega$-repr\'esentation temp\'er\'ee $\tilde{\pi}$ de $\tilde{G}(F)$, alors $I_{\tilde{G}}(\gamma,\omega,f)=0$ pour tout $\gamma$ fortement r\'egulier (th\'eor\`eme 5.5).
  
  Pour aller plus loin, on doit transformer la formule pr\'ec\'edente en une formule invariante, comme dans [A7]. Pour cela, on doit utiliser le th\'eor\`eme de Paley-Wiener tordu. Celui-ci a \'et\'e d\'emontr\'e par Rogawski  ([R]) dans le cas o\`u $F$ est non-archim\'edien et $\omega=1$. Il est d\'emontr\'e toujours dans le cas non-archim\'edien, mais pour tout $\omega$, dans le travail en cours de Henniart et Lemaire ([HL]). Dans le cas o\`u $F$ est archim\'edien, il est d\'emontr\'e par Delorme et Mezo pour $\omega=1$ ([DM]).   Nous montrons en 6.3 que leur r\'esultat s'\'etend ais\'ement au cas $\omega$ quelconque.  A partir de ce point, les fonctions $f_{1}$ et $f_{2}$ sont suppos\'ees $K$-finies quand $F$ est archim\'edien. En utilisant le th\'eor\`eme de Paley-Wiener tordu, on transforme la formule en l'\'egalit\'e suivante (th\'eor\`eme 6.6):
  $$I^{\tilde{G}}_{g\acute{e}om}(\omega,f_{1},f_{2})=I^{\tilde{G}}_{disc}(\omega,f_{1},f_{2}).$$
  Le terme de gauche est de la forme
  $$\sum_{\tilde{M}\in {\cal P}(\tilde{M}_{0})}\vert \tilde{W}^M\vert \vert \tilde{W}^G\vert ^{-1}I_{\tilde{M},g\acute{e}om}^{\tilde{G}}(\omega,f_{1},f_{2}).$$
  Il ne contient que des distributions $\omega$-\'equivariantes (int\'egrales orbitales pond\'er\'ees invariantes). Mais le terme principal est le m\^eme que pr\'ec\'edemment:
  $$I_{\tilde{G},g\acute{e}om}^{\tilde{G}}(\omega,f_{1},f_{2})=J_{\tilde{G},g\acute{e}om}^{\tilde{G}}(\omega,f_{1},f_{2}).$$
  Du c\^ot\'e spectral, on a simplement
  $$I^{\tilde{G}}_{disc}(\omega,f_{1},f_{2})=J^{\tilde{G}}_{\tilde{G},spec}(\omega,f_{1},f_{2}).$$
  Il ne contient que d'honn\^etes caract\`eres de $\omega$-repr\'esentations temp\'er\'ees.
  
 Supposons pour simplifier que $Z_{G}(F)^{\theta}$ soit compact. L'ensemble  des $\omega$-repr\'esentations qui interviennent ici contient le sous-ensemble des repr\'esentations elliptiques au sens d'Arthur. Notons-le ici $E_{ell}(\tilde{G},\omega)$ (cette notation changera de sens dans le corps de l'article). On sait que le caract\`ere d'une $\omega$-repr\'esentation $\tilde{\pi}$ de longueur finie est localement int\'egrable, donc donn\'e par une fonction $\gamma\mapsto \Theta(\tilde{\pi},\gamma)$. D'apr\`es le th\'eor\`eme de Paley-Wiener, on peut associer \`a tout \'el\'ement $\tilde{\pi}\in E_{ell}(\tilde{G},\omega)$ un "pseudo-coefficient" $f_{\tilde{\pi}}$. La formule ci-dessus permet d'exprimer $\Theta(\tilde{\pi},\gamma)$  au moyen des int\'egrales orbitales pond\'er\'ees invariantes  de $f_{\tilde{\pi}}$ (th\'eor\`eme 7.2). Notons $\tilde{G}(F)_{ell}$ l'ensemble des \'el\'ements fortement r\'eguliers et elliptiques de $\tilde{G}(F)$. On peut munir l'espace des fonctions $\omega$-\'equivariantes sur $\tilde{G}(F)_{ell}$ d'un produit hermitien raisonnable. Notons $I_{cusp}(\tilde{G}(F),\omega)$ l'espace des fonctions sur $\tilde{G}(F)_{ell}$ de la forme $\gamma\mapsto I_{\tilde{G}}(\gamma,\omega,f)$, o\`u $f$ est une fonction cuspidale sur $\tilde{G}(F)$ (c'est-\`a-dire que $I_{\tilde{G}}(\gamma,\omega,f)=0$ si $\gamma$ est fortement r\'egulier et non elliptique).
 On montre que la famille des restrictions \`a $\tilde{G}(F)_{ell}$ des caract\`eres $\Theta(\tilde{\pi},\gamma)$ forme une base orthonormale  de $I_{cusp}(\tilde{G}(F),\omega)$ et on calcule la norme de chaque \'el\'ement de base (th\'eor\`eme 7.3).
 
 Une bonne partie de l'article n'est qu'un d\'ecalque des travaux d'Arthur. On s'est autoris\'e \`a passer rapidement sur les points dont les d\'emonstrations sont similaires \`a celles du cas non tordu. Le point qui m\'erite attention est l'\'etablissement de la partie spectrale de version non invariante de la formule des traces locale tordue (paragraphe 3).  Comme dans le cas de la formule des traces d'Arthur-Selberg, l'utilisation dans le cas tordu de la combinatoire du cas non tordu ne conduit \`a rien de fructueux. On doit profond\'ement modifier cette combinatoire. Heureusement pour nous, la bonne combinatoire a \'et\'e mise au point par les r\'edacteurs du Morning Seminar, principalement dans la lecture 15 de Langlands. Ces notes ont \'et\'e r\'edig\'ees r\'ecemment (cf. [LW]). On utilise ici exactement la m\^eme m\'ethode, dont la mise en oeuvre est \'evidemment beaucoup plus simple que dans le cas global.

\section{G\'en\'eralit\'es}

\bigskip

\subsection{Notations}

Soit $F$ un corps local  de caract\'eristique nulle. On note $\vert .\vert _{F}$ sa valeur absolue usuelle. Soit $G$ un groupe r\'eductif connexe d\'efini sur $F$. On note $Z_{G}$ le centre de $G$ et $A_{G}$ le plus grand tore  contenu dans $Z_{G}$ et d\'eploy\'e sur $F$. On pose ${\cal A}_{G}=X_{*}(A_{G})\otimes_{{\mathbb Z}}{\mathbb R}$, o\`u $X_{*}(A_{G})$ d\'esigne selon l'usage le groupe des sous-groupes \`a un param\`etre de $A_{G}$. On note $a_{G}$ la dimension de ${\cal A}_{G}$. Notons $X_{F}^*(G)$ le groupe des caract\`eres alg\'ebriques d\'efinis sur $F$ de $G$. Par restriction \`a $A_{G}$, il s'envoie dans le groupe $X^*(A_{G})$ des caract\`eres de $A_{G}$, donc dans le dual ${\cal A}_{G}^*$ de ${\cal A}_{G}$. De fa\c{c}on g\'en\'erale, on note $<.,.>$ l'accouplement naturel entre un espace vectoriel et son dual. On d\'efinit l'homomorphisme $H_{G}:G(F)\to {\cal A}_{G}$ par la condition: $e^{<\chi,H_{G}(g)>}=\vert \chi(g)\vert_{F}$ pour tout $g\in G(F)$ et tout $\chi\in X_{F}^*(G)$. On note ${\cal A}_{G,F}$ l'image de $G(F)$ par cet homomorphisme et ${\cal A}_{A_{G},F}$ l'image de $A_{G}(F)$. Si $F$ est archim\'edien, ${\cal A}_{G,F}={\cal A}_{A_{G},F}={\cal A}_{G}$. Si $F$ est non archim\'edien, ${\cal A}_{G,F}$ et ${\cal A}_{A_{G},F}$ sont des r\'eseaux dans ${\cal A}_{G}$. Plus pr\'ecis\'ement, ce sont des r\'eseaux dans $log(q)X_{*}(A_{G})\otimes_{{\mathbb Z}}{\mathbb Q}$, o\`u $q$ est le nombre d'\'el\'ements du corps r\'esiduel de $F$.  Pour tout sous-groupe ferm\'e ${\cal L}$ dans ${\cal A}_{G}$, on note ${\cal L}^{\vee}$ le groupe des  $\lambda\in {\cal A}_{G}^*$ tels que $<\lambda,H>\in 2\pi {\mathbb Z}$ pour tout $H\in{\cal L}$. Ainsi, via l'exponentielle, $i{\cal A}_{G}^*/i{\cal L}^{\vee}$ s'identifie au  groupe dual de ${\cal L}$. On pose ${\cal A}_{G,F}^*={\cal A}_{G}^*/{\cal A}_{G,F}^{\vee}$.   

Sauf mention expresse du contraire, tous les sous-groupes alg\'ebriques de $G$ que l'on consid\'erera seront suppos\'es d\'efinis sur $F$.   Un Levi de $G$ est une composante de Levi d'un sous-groupe parabolique de $G$.  Une paire parabolique est un couple $(P,M)$ o\`u $P$ est un sous-groupe parabolique et $M$ est une composante de Levi de $P$. L'expression "soit $P=MU_{P}$ un sous-groupe parabolique" signifie que $M$ est une composante de Levi de $P$ et que $U_{P}$ est le radical unipotent de $P$. Plus pr\'ecis\'ement, si une paire parabolique minimale  $(P_{0},M_{0})$ est fix\'ee, l'expression "soit $P=MU_{P}$ un sous-groupe parabolique semi-standard (ou standard)" signifie que $P$ contient $M_{0}$ (ou $P_{0}$), que $M$ est la composante de Levi de $P$ qui contient $M_{0}$ et que $U_{P}$ est le radical unipotent de $P$. On utilise les notations habituelles d'Arthur. Par exemple, pour un Levi $M$, on note ${\cal P}(M)$, resp. ${\cal F}(M)$, l'ensemble des sous-groupes paraboliques $P$ dont $M$ est une composante de Levi, resp. qui contiennent $M$. Pour un sous-groupe parabolique $P=MU_{P}$, on note $\delta_{P}$ le module usuel, qui est une fonction sur $P(F)$.

 Fixons une paire parabolique minimale $(P_{0},M_{0})$ et un sous-groupe compact maximal $K$ de $G(F)$. On suppose:
 
 - si $F$ est non-archim\'edien, $K$ est le fixateur d'un point sp\'ecial dans l'appartement attach\'e \`a $A_{M_{0}}$ de l'immeuble de $G$;
 
 - si $F$ est archim\'edien, les alg\`ebres de Lie de $K$ et de $A_{M_{0}}$ sont orthogonales pour la forme de Killing.
 
 On simplifie les notations en notant par un simple indice $0$ les objets relatifs \`a $M_{0}$ ou $P_{0}$. Par exemple $A_{0}=A_{M_{0}}$,  ${\cal A}_{0}={\cal A}_{M_{0}}$, $H_{0}=H_{M_{0}}$, $\delta_{0}=\delta_{P_{0}}$. On note $\Delta_{0}$ l'ensemble des racines simples de $A_{0}$ associ\'e \`a $P_{0}$. On note ${\cal A}_{0}^{\geq}$ l'ensemble des $H\in {\cal A}_{0}$ tels que $<\alpha,H>\geq 0$ pour tout $\alpha\in \Delta_{0}$. On note $M_{0}(F)^{\geq}$ l'ensemble des $m\in M_{0}(F)$ tels que $H_{0}(m)\in {\cal A}_{0}^{\geq}$. On a l'\'egalit\'e $G(F)=KM_{0}(F)^{\geq}K$. Plus pr\'ecis\'ement, pour $g\in G(F)$, si l'on \'ecrit $g=kmk'$, avec $k,k'\in K$ et $m\in M_{0}(F)^{\geq}$, l'\'el\'ement $H_{0}(m)$ est uniquement d\'etermin\'e par $g$, on le note $h_{0}(g)$. On note $W^G$ le groupe de Weyl de $G$ relatif au tore $A_{0}$, c'est-\`a-dire $W^G=Norm_{G(F)}(A_{0})/M_{0}(F)$ (de fa\c{c}on g\'en\'erale, si un groupe $H$ op\`ere sur un ensemble $X$ et si $Y\subset X$, on note $Norm_{H}(Y)$ le sous-groupe des $h\in H$ tels que $h(Y)=Y$). On fixe une forme quadratique  $(.,.)$ sur ${\cal A}_{0}$ d\'efinie positive et invariante par $W^G$. On d\'efinit la norme $\vert H\vert =(H,H)^{1/2}$ pour tout $H\in {\cal A}_{0}$. Par dualit\'e, ${\cal A}_{0}^*$ est aussi muni d'une norme.  
 
 Quand on remplace le groupe $G$ par un autre groupe r\'eductif, par exemple un Levi $M$, on ajoute des exposants $M$ dans les notations. Toutefois, si on remplace $G$ par une composante de Levi $M$ d'un sous-groupe parabolique $P$ fix\'e, il est parfois plus commode d'ajouter un exposant $P$ au lieu de $M$, ou de remplacer un indice $M$ par $P$. Par exemple, on pourra noter ${\cal A}_{P}$ au lieu de ${\cal A}_{M}$. Ou encore, si $P$ est standard, on notera $\Delta_{0}^M$ ou $\Delta_{0}^P$ l'ensemble des racines simples de $A_{0}$ dans $M$ associ\'e au parabolique minimal $P_{0}\cap M$. Soit $M$  un Levi de $G$. Choissons un \'el\'ement $x\in G(F)$ tel que $M'=x^{-1}Mx$ contienne $M_{0}$.  On munit le sous-espace ${\cal A}_{M'}$ de ${\cal A}_{0}$ de la restriction de la forme quadratique fix\'ee plus haut. De la conjugaison $ad_{x}$ se d\'eduit fonctoriellement un isomorphisme not\'e simplement $H\mapsto xH$ de ${\cal A}_{M'}$ sur ${\cal A}_{M}$, gr\^ace auquel on transporte \`a ${\cal A}_{M}$ la forme quadratique sur ${\cal A}_{M'}$. La forme obtenue ne d\'epend pas du choix de $x$.
 
  Si $M'\subset M$ sont deux Levi, l'espace ${\cal A}_{M}$ est un sous-espace de ${\cal A}_{M'}$ dont on note ${\cal A}^M_{M'}$ l'orthogonal. Pour $H\in {\cal A}_{M'}$, on note $H_{M}$ et $H^M$ ses projections sur chacun de ces sous-espaces. Remarquons que
  
    (1) l'application $H\mapsto H_{M}$ envoie surjectivement ${\cal A}_{M',F}$ sur ${\cal A}_{M,F}$.
  
  Preuve. En conjuguant $M'$, on peut supposer  $M_{0}\subset M'$. Si $H=H_{M'}(x)$, pour $x\in M'(F)$, on a $H_{M}=H_{M}(x)$. Inversement, si $H=H_{M}(y)$, pour $y\in M(F)$, on fixe $P'\in {\cal P}^M(M')$ et on  \'ecrit $y=xuk$, avec $x\in M'(F)$, $u\in U_{P'}(F)$ et $k\in K\cap M(F)$. Alors $H_{M}(y)=H_{M}(x)$ et $H_{M}=(H_{M'}(x))_{M}$. D'o\`u (1). $\square$
  
  Dualement \`a (1), on a
  
  (2) l'injection ${\cal A}_{M}^*\to {\cal A}_{M'}^*$ se quotiente en une injection ${\cal A}_{M,F}^*\to {\cal A}_{M',F}^*$. 
  
    Soit $P=MU_{P}$ est un sous-groupe parabolique semi-standard. En utilisant  l'\'egalit\'e $G(F)=P(F)K$, on prolonge la fonction $H_{M}$ en une fonction $H_{P}:G(F)\to {\cal A}_{M}$ par $H_{P}(muk)=H_{M}(m)$ pour tous $m\in M(F)$, $u\in U_{P}(F)$, $k\in K$. Si de plus $P$ est standard, on note $M_{0}(F)^{\geq,M}$ ou $M_{0}(F)^{\geq,P}$  l'ensemble des $m\in M_{0}(F)$ tels que $<\alpha,H_{0}(m)>\geq 0$ pour tout $\alpha\in \Delta_{0}^P$.    
   
 Si $X(x_{1},x_{2},...)$ et $Y(x_{1},x_{2},...)$ sont deux expressions r\'eelles positives ou nulles d\'ependant de variables $x_{1}$, $x_{2}$,..., on dit que $X(x_{1},x_{2},...)$ est essentiellement major\'ee par $Y(x_{1},x_{2},...)$ et on \'ecrit $X(x_{1},x_{2},...)<<Y(x_{1},x_{2},...)$ s'il existe $c>0$ tel que, pour tous $x_{1},x_{2},...$, on ait l'in\'egalit\'e $X(x_{1},x_{2},...)\leq cY(x_{1},x_{2},...)$. Cette terminologie est quelque peu impr\'ecise (la question \'etant  en pratique de savoir quelles donn\'ees sont "variables") mais \'evite d'introduire des kyrielles de constantes $c$ superflues.
 
 On note $C_{c}^{\infty}(G(F))$ l'espace des fonctions $f:G(F)\to {\mathbb C}$ qui sont lisses et \`a support compact. Lisse signifie localement constante si $F$ est non-archim\'edien, $C^{\infty}$ si $F$ est archim\'edien.

 \bigskip
 
 \subsection{Mesures}
  
 On munit $G(F)$ d'une mesure de Haar. Pour tout Levi $M$, on munit l'espace ${\cal A}_{M}$ d'une mesure de Haar. On suppose que, si $x\in G(F)$, l'isomorphisme $H\mapsto xH$ de ${\cal A}_{M}$ sur ${\cal A}_{xMx^{-1}}$ identifie les mesures sur ces espaces. L'espace $i{\cal A}_{M}^*$ s'identifie via l'exponentielle au groupe dual de ${\cal A}_{M}$ et on le munit de la mesure duale. Cela entra\^{\i}ne que, si ${\cal L}$ est un r\'eseau de ${\cal A}_{M}$, on a l'\'egalit\'e
 $$mes({\cal A}_{M}/{\cal L})mes(i{\cal A}_{M}^*/i{\cal L}^{\vee})=1,$$
 les mesures \'etant les quotients des mesures que l'on vient de fixer par les mesures de comptage sur les r\'eseaux. On pose $A_{M}(F)_{c}=Ker(H_{M})\cap A_{M}(F)$. Dans le cas o\`u $F$ est non-archim\'edien, on note $mes(i{\cal A}_{M,F}^*)$ la masse totale de $i{\cal A}_{M,F}^*$.
 Le groupe $ A_{M}(F)_{c}$ est un sous-groupe ouvert compact de $A_{M}(F)$. On munit $A_{M}(F)$ de la mesure de Haar pour laquelle
 $$mes( A_{M}(F)_{c})=mes({\cal A}_{M}/{\cal A}_{A_{M},F}).$$
 Si $F$ est archim\'edien,  on pose $mes(i{\cal A}_{M,F}^*)=1$. On a la suite exacte
 $$1\to  A_{M}(F)_{c}\to A_{M}(F)\to {\cal A}_{M}\to 0.$$
 On munit le groupe compact $ A_{M}(F)_{c}$ de la mesure de Haar de masse totale $1$ et le groupe $A_{M}(F)$ de la mesure de Haar compatible avec la suite ci-dessus.  En tout cas, soit $f$ une fonction de Schwartz sur $i{\cal A}_{M,F}^*$, d\'efinissons une fonction $\hat{f}$ sur ${\cal A}_{M,F}$ par 
  $$\hat{f}(H)=mes(i{\cal A}_{M,F}^*)^{-1}\int_{i{\cal A}_{M,F}^*}f(\Lambda)e^{-<\Lambda,H>}\,d\Lambda.$$
  On a la formule d'inversion
  $$f(\Lambda)=\int_{{\cal A}_{M,F}}\hat{f}(H)e^{<\Lambda,H>}\,dH.$$
  Remarquons qu'une int\'egrale sur ${\cal A}_{M,F}$ est en fait une s\'erie  si $F$ est non-archim\'edien.
  
  Pour deux Levi $M'\subset M$, on munit ${\cal A}_{M'}^M$ 
   de la mesure compatible aux mesures sur ${\cal A}_{M}$ et ${\cal A}_{M'}$ et \`a l'isomorphisme ${\cal A}_{M'}={\cal A}_{M}\oplus {\cal A}_{M'}^M$. On munit $i{\cal A}_{M'}^{M,*}$ de la mesure duale comme ci-dessus. Dans le cas non archim\'edien, la masse totale de $i{\cal A}_{M'}^{M,*}/\left((i{\cal A}_{M',F}^{\vee}+i{\cal A}_{M}^*)\cap i{\cal A}_{M'}^{M,*}\right)$ est $mes(i{\cal A}_{M',F}^*)mes(i{\cal A}_{M,F}^*)^{-1}$.
   
  Soit $\alpha$ une racine r\'eduite de $A_{0}$ dans $G$. Comme on le sait, il est attach\'e \`a $\alpha$ un Levi $M_{\alpha}$ qui est minimal parmi les \'el\'ements de ${\cal L}(M_{0})-\{M_{0}\}$ et un sous-groupe parabolique $P_{\alpha}=M_{0}U_{\alpha}$ de $M_{\alpha}$. On munit $U_{\alpha}(F)$ d'une mesure de Haar. On suppose que si $\alpha$ et $\beta$ sont deux telles racines r\'eduites et si $k\in K$ v\'erifie $kM_{0}k^{-1}=M_{0}$, $kP_{\alpha}k^{-1}=P_{\beta}$, alors la conjugaison par $k$ identifie les mesures sur $U_{\alpha}(F)$ et $U_{\beta}(F)$. Consid\'erons un sous-groupe unipotent $U$ de $G$ qui est produit de tels groupes $U_{\alpha_{i}}$, pour $i=1,...,n$. On munit $U(F)$ de la mesure pour laquelle le produit
  $$\prod_{i=1,...,n}U_{\alpha_{i}}(F)\to U(F)$$
  pr\'eserve les mesures. Cela ne d\'epend pas de l'ordre du produit. L'hypoth\`ese sur $U$ est v\'erifi\'ee si $U$ est le radical unipotent d'un sous-groupe parabolique semi-standard de $G$, ou si $U$ est le radical unipotent d'un sous-groupe parabolique semi-standard d'un Levi semi-standard de $G$, ou si $U$ est une intersection de tels groupes.
  
  Pour tout Levi $M\in {\cal L}(M_{0})$, on munit $M(F)$ d'une mesure. On impose les conditions suivantes, dont on v\'erifie ais\'ement qu'elles sont loisibles. D'abord, si $M$ et $M'$ sont deux \'el\'ements de ${\cal L}(M_{0})$ et si $x\in G(F)$ v\'erifie $xMx^{-1}=M'$, on suppose que la conjugaison par $x$ identifie les mesures sur $M(F)$ et $M'(F)$. Ensuite, si $M\subset L$ sont deux \'el\'ements de ${\cal L}(M_{0})$  et si $P=MU_{P}\in {\cal P}^L(M)$, on impose que l'on ait l'\'egalit\'e
 $$\int_{L(F)}f(l)\,dl\,=\int_{U_{\bar{P}}(F)}\int_{M(F)}\int_{U_{P}(F)}f(\bar{u}mu)\delta_{P}(m)\,du\,dm\,d\bar{u}$$
 pour toute $f\in C_{c}^{\infty}(L(F))$, o\`u $\bar{P}=MU_{\bar{P}}$ est le parabolique oppos\'e \`a $P$.   
  
  On munit le groupe $K$ de la mesure de Haar de masse totale $1$. Remarquons que, dans le cas o\`u $F$ est non-archim\'edien, on n'impose pas que la mesure sur $G(F)$ se restreigne en cette mesure sur $K$.   Plus g\'en\'eralement, pour tout Levi semi-standard $M$ de $G$, on pose $K^M=K\cap M(F)$ et on munit ce groupe de la mesure de Haar de masse totale $1$. On v\'erifie que, pour deux L\'evi semi-standard $M\subset L$, il existe un r\'eel $\gamma(L\vert M)>0$ tel que, pour tout $P=MU_{P}\in {\cal P}^L(M)$, on ait l'\'egalit\'e
  $$\int_{L(F)}f(l)\,dl\,=\gamma(L\vert M)\int_{M(F)}\int_{U_{P}(F)} \int_{K^L}f(muk)\,dk\,du\,dm$$
 pour toute $f\in C_{c}^{\infty}(L(F))$. On a l'\'egalit\'e $\gamma(G\vert M)=\gamma(G\vert L)\gamma(L\vert M)$.

 On introduit la fonction $D_{0}$ sur $M_{0}(F)^{\geq}$, \`a valeurs r\'eelles positives, telle que l'on ait l'\'egalit\'e
 $$\int_{G(F)}f(g)\,dg=\int_{K\times K}\int_{M_{0}(F)^{\geq}}D_{0}(m)f(kmk')\,dm\,dk\,dk'$$
 pour toute $f\in C_{c}^{\infty}(G(F))$.  Elle v\'erifie la majoration
 
 (1)  $D_{0}(m)<<\delta_{0}(m)^{1/2}$ pour tout $m\in M_{0}(F)^{\geq}$.
 
 Cf. [A1] corollaire 1.2. 
 
 {\bf Remarque.} Nos mesures ne sont pas normalis\'ees comme en [A1]: pour $M\in {\cal L}(M_{0})$, Arthur prend pour mesure sur $M(F)$ celle que l'on a d\'efinie multipli\'ee par $\gamma(G\vert M)^2$. Cela entra\^{\i}ne que notre fonction $D_{0}$ est \'egale \`a celle d'Arthur multipli\'ee par $\gamma(G\vert M_{0})^{2}$.

 \bigskip
 
 \subsection{D\'efinitions combinatoires}
 
 Les propri\'et\'es \'enonc\'ees dans ce paragraphe sont aujourd'hui bien connues. La plupart sont dues \`a Arthur, Langlands et Labesse. On donne pour r\'ef\'erence l'article r\'ecent [LW] qui les a rassembl\'ees.
 
 Soit $P=MU_{P}$ un sous-groupe parabolique semi-standard de $G$. On note $\Delta_{P}$ l'ensemble des racines simples de $A_{M}$ associ\'e \`a $P$. Il s'agit d'une base de ${\cal A}_{M}^{G,*}$. On note $\{\check{\varpi}_{\alpha}; \alpha\in \Delta_{P}\}$ sa base duale. A toute racine $\alpha\in \Delta_{P}$, on peut associer une coracine $\check{\alpha}\in {\cal A}_{M}^G$. Dans le cas o\`u $P=P_{0}$ et $M=M_{0}$, l'ensemble des racines de $A_{0}$ est un vrai syst\`eme de racines et on definit les coracines selon l'usage. Dans le cas g\'en\'eral, quitte \`a conjuguer $P$, on peut supposer $P$ standard. Alors $\Delta_{P}$ est l'ensemble des projections non nulles sur ${\cal A}_{M}$ d'\'el\'ements de $\Delta_{0}$ et on d\'efinit les coracines comme les projections non nulles des coracines associ\'ees aux \'el\'ements de $\Delta_{0}$. Cela ne d\'epend pas de la conjugaison.  
 On introduit la base $\{\varpi_{\alpha};\alpha \in \Delta_{P}\}$ de ${\cal A}_{M}^{G,*}$ duale de celle des coracines.   Une propri\'et\'e essentielle est que $(\alpha,\beta)\leq0$ et  $(\varpi_{\alpha},\varpi_{\beta})\geq0$ pour tous $\alpha\not=\beta \in \Delta_{P}$.
 
 Soient $P=MU_{P}\subset Q=LU_{Q}$ deux sous-groupes paraboliques semi-standard. On note $\Delta_{P}^Q\subset \Delta_{P}$ l'ensemble des racines simples de $A_{M}$ dans l'alg\`ebre de Lie de $L\cap U_{P}$. On d\'efinit les fonctions suivantes sur ${\cal A}_{M}$:
 
 $\tau_{P}^Q$ fonction caract\'eristique de l'ensemble des $H\in {\cal A}_{M}$ tels que $<\alpha,H>>0$ pour tout $\alpha\in \Delta_{P}^Q$;
 
 $\hat{\tau}_{P}^Q$ fonction caract\'eristique de l'ensemble des $H\in {\cal A}_{M}$ tels que $<\varpi_{\alpha}^L,H>>0$ pour tout $\alpha\in \Delta_{P}^Q$;
 
 $\phi_{P}^Q$ fonction caract\'eristique de l'ensemble des $H\in {\cal A}_{M}$ tels que $<\varpi_{\alpha}^L,H>\leq 0$ pour tout $\alpha\in \Delta_{P}^Q$;
 
 $\delta_{M}^Q$ fonction caract\'eristique du sous-ensemble ${\cal A}_{L}$ de ${\cal A}_{M}$.

 Si $P'=M'U_{P'}$ est un sous-groupe parabolique semi-standard contenu dans $P$, ces fonctions peuvent \^etre consid\'er\'ees comme des fonctions sur ${\cal A}_{M'}$, en les identifiant avec leurs compos\'ees avec la projection orthogonale de ${\cal A}_{M'}$ sur ${\cal A}_{M}$.  On supprimera souvent les indices $P$ quand $P=P_{0}$.
 
 On a
 $$(1) \qquad \sum_{R; M\subset R\subset Q}\delta_{M}^R(H)\tau_{R}^Q(H)=1.$$
 Cela traduit la d\'ecomposition de ${\cal A}_{M}$ en chambres positives relatives aux paraboliques $R$ indiqu\'es. On a l'\'egalit\'e (qui est un lemme de Langlands):
 
 $$(2) \qquad \sum_{R; P\subset R\subset Q}\phi_{P}^R(H)\tau_{R}^Q(H)=1.$$

 On d\'efinit la fonction $\Gamma_{P}^Q$ sur ${\cal A}_{M}\times {\cal A}_{M}$ par 
 $$\Gamma_{P}^Q(H,X)=\sum_{R; P\subset R\subset Q}(-1)^{a_{R}-a_{Q}}\tau_{P}^R(H)\hat{\tau}_{R}^Q(H-X).$$
 On rappelle que $a_{R}$, resp. $a_{Q}$, est la dimension de ${\cal A}_{R}$, resp. ${\cal A}_{Q}$.  Sur le support de cette fonction, on a une majoration $\vert H_{M}^L\vert <<\vert X_{M}^L\vert $.  On a
 $$(3) \qquad \sum_{R; P\subset R\subset Q}\Gamma_{R}^Q(H,X)\phi_{P}^R(H)=\phi_{P}^Q(H-X).$$
 
 Preuve. Il r\'esulte de la formule du bin\^ome que la d\'efinition de $\phi_{P}^Q(H)$ peut s'\'ecrire
 $$(4) \qquad \phi_{P}^Q(H)=\sum_{S; P\subset S\subset Q}(-1)^{a_{S}-a_{Q}}\hat{\tau}_{S}^Q(H).$$
 On applique cette d\'efinition en rempla\c{c}ant $Q$ par $R$ ainsi que la d\'efinition de $\Gamma_{R}^Q(H,X)$. Le membre de gauche de (3) devient
 $$\sum_{S,R,S'; P\subset S\subset R\subset S'\subset Q}(-1)^{a_{S}-a_{R}+a_{S'}-a_{Q}}\hat{\tau}_{S}^R(H)\tau_{R}^{S'}(H)\hat{\tau}_{S'}^Q(H-X).$$
 Pour $S,S'$ fix\'es, la somme
 $$\sum_{R; S\subset R\subset S'}(-1)^{a_{S}-a_{R}}\hat{\tau}_{S}^R(H)\tau_{R}^{S'}(H)$$
 vaut $1$ si $S=S'$ et $0$ sinon (cf. [LW] proposition 1.7.2). L'expression pr\'ec\'edente se simplifie en
 $$\sum_{S; P\subset S\subset Q}(-1)^{a_{S}-a_{Q}}\hat{\tau}_{S}^Q(H-X).$$
 En appliquant de nouveau (4), c'est $\phi_{P}^Q(H-X)$. Cela prouve (3). $\square$
 
 On a 
 $$(5)\qquad \sum_{R; P\subset R\subset Q}\Gamma_{P}^R(H,X)\tau_{R}^Q(H-X)=\tau_{P}^Q(H).$$
 La preuve est similaire \`a celle de (3).
 
 Soit $M$ un Levi semi-standard. Consid\'erons une famille  ${\cal Y}=(Y[P])_{P\in {\cal P}(M)}$ o\`u, pour tout $P$, $Y[P]$ est un \'el\'ement de ${\cal A}_{M}$.   On dit qu'elle est $(G,M)$-orthogonale si elle v\'erifie les conditions suivantes. Soient $P,P'\in {\cal P}(M)$ deux \'el\'ements adjacents. Soit $\alpha$ l'unique racine r\'eduite de $A_{M}$ qui est positive pour $P$ et n\'egative pour $P'$. Alors $Y[P]-Y[P']$ appartient \`a la droite port\'ee par $\check{\alpha}$.  On dit que ${\cal Y}$ est positive si $Y[P]-Y[P']$ appartient \`a la demi-droite port\'ee par $\check{\alpha}$. Si $F$ est non-archim\'edien, on dit que ${\cal Y}$ est $p$-adique si $Y[P]\in {\cal A}_{M,F}\otimes_{{\mathbb Z}}{\mathbb Q}$ pour tout $P\in {\cal P}(M)$. Pour une famille $(G,M)$-orthogonale ${\cal Y}$ et pour $Q=LU_{Q}\in {\cal F}(M)$, on pose $Y[Q]=Y[P]_{L}$, o\`u $P$ est un \'el\'ement de ${\cal P}(M)$ tel que $P\subset Q$. Cela ne d\'epend pas du choix de $P$. Si on fixe $L\in {\cal L}(M)$, la famille $(Y[Q])_{Q\in {\cal P}(L)}$ est encore $(G,L)$-orthogonale et est positive si la famille de d\'epart l'est. Dans le cas particulier o\`u $M=M_{0}$, on associe \`a tout \'el\'ement $T\in {\cal A}_{0}$ une famille $(G,M_{0})$-orthogonale ${\cal T}=(T[P])_{P\in {\cal P}(M_{0})}$ de la fa\c{c}on suivante. Pour $P\in {\cal P}(M_{0})$, soit $w\in W^G$ tel que $w(P_{0})=P$. On pose $T[P]=wT$. Cette famille est positive si et seulement si $T\in {\cal A}_{0}^{\geq}$.

Pour une famille $(G,M)$-orthogonale ${\cal Y}=(Y[P])_{P\in {\cal P}(M)}$ et pour $Q=LU_{Q}\in {\cal F}(M)$, on d\'efinit une fonction  $\Gamma_{M}^Q(.,{\cal Y})$ sur ${\cal A}_{M}$ par
 $$\Gamma_{M}^Q(H,{\cal Y})=\sum_{R\in {\cal F}(M); R\subset Q}\delta_{M}^R(H)\Gamma_{R}^Q(H,Y[R]).$$
 Sur le support de cette fonction, on a une majoration $\vert H^L\vert <<sup_{P\in {\cal P}(M); P\subset Q}\vert Y(P)^L\vert $. Dans le cas o\`u la famille est positive, c'est la fonction caract\'eristique de l'ensemble des $H\in {\cal A}_{M}$ tels que $H^L$ appartient \`a l'enveloppe convexe de la famille des points $Y[P]^L$ pour $P\in {\cal P}(M)$, $P\subset Q$. On a l'\'egalit\'e
 $$(6) \qquad \sum_{Q\in {\cal F}(M)}\Gamma_{M}^Q(H,{\cal Y})\tau_{Q}^G(H-Y[Q])=1.$$
 Cf. [LW] lemme 1.8.4(3).
 
 Fixons $P_{1}\in {\cal P}(M)$, ce qui d\'efinit un ordre sur l'ensemble des racines de $A_{M}$, not\'e $\alpha>_{P_{1}}0$. Pour $P\in {\cal P}(M)$, notons $\phi_{P,P_{1}}^G$ la fonction caract\'eristique de l'ensemble des $H\in {\cal A}_{M}$ tels que
 
 - pour $\alpha\in \Delta_{P}$ tel que $ \alpha>_{P_{1}}0$, $<\varpi_{\alpha},H>\leq 0$;
 
 - pour  $\alpha\in \Delta_{P}$ tel que $ \alpha<_{P_{1}}0$, $<\varpi_{\alpha},H>> 0$.
 
 Notons $s(P,P_{1})$ le nombre de $\alpha\in \Delta_{P}$ tel que $ \alpha<_{P_{1}}0$.   On a l'\'egalit\'e
 
 $$(7)\qquad \Gamma_{M}^G(H,{\cal Y})=\sum_{P\in {\cal P}(M)}(-1)^{s(P,P_{1})}\phi_{P,P_{1}}^G(H-Y[P]).$$
 Cf. [LW] proposition 1.8.7(2).

 Pour $H\in {\cal A}_{0}$, on note ${\cal C}^G(H)$ l'enveloppe convexe des $wH$ pour $w\in W^G$. Supposons  $H\in {\cal A}_{0}^{\geq}$. Alors, pour tout $H'\in {\cal A}_{0}$, $H'$ appartient \`a ${\cal C}^G(H)$ si et seulement si on a $\phi^G_{w(P_{0})}(H'-wH)=1$ pour tout $w\in W^G$. Pour $H'\in {\cal A}_{0}^{\geq}$, $H'$ appartient \`a ${\cal C}^G(H)$ si et seulement si $\phi_{P_{0}}^G(H'-H)=1$. 
 
 \ass{Lemme}{Pour tout $m\in M_{0}(F)$, tout $k\in K$ et tout $P\in {\cal P}(M_{0})$, $H_{P}(km)$ appartient \`a ${\cal C}^G(H_{0}(m))$.}

Preuve. On utilise les deux ingr\'edients suivants.

(8) Pour tous $g\in G(F)$ et tous $P,P'\in {\cal P}(M_{0})$, on a $\phi_{P}^G(H_{P}(g)-H_{P'}(g))=1$.

 Autrement dit, la famille $(-H_{P}(g))_{P\in {\cal P}(M_{0})}$ est $(G,M_{0})$-orthogonale et positive, cf. [A2] p. 40.

(9) Supposons $m\in M_{0}(F)^{\geq}$; alors, pour tout $k\in K$, $\phi_{P_{0}}^G(H_{\bar{P}_{0}}(km)-H_{0}(m))=1$.

Ceci est un ingr\'edient de la th\'eorie de la transform\'ee de Satake, cf. [HC] lemme 35 et [BT] corollaire 4.3.17.

Venons-en \`a la preuve du lemme. Conjuguer $m$ par un \'el\'ement de $Norm_{K}(M_{0})$ ne change pas le probl\`eme. On peut donc supposer $m\in M_{0}(F)^{\geq}$. Soit $w\in W^G$. On doit prouver que $\phi^G_{w(P_{0})}(H_{P}(km)-wH_{0}(m))=1$. Ou encore que $\phi^G_{P_{0}}(w^{-1}H_{P}(km)-H_{0}(m))=1$. Or $w^{-1}H_{P}(km)=H_{w^{-1}(P)}(\dot{w}^{-1}km)$, o\`u $\dot{w}$ est un repr\'esentant de $w$ dans $K$. Quitte \`a remplacer $k$ par $\dot{w}k$ et $P$ par $w(P)$, on est ramen\'e \`a prouver $\phi_{P_{0}}^G(H_{P}(km)-H_{0}(m))=1$ pour tous $k$, $P$.  D'apr\`es (8), on a $\phi_{\bar{P}_{0}}^G(H_{\bar{P}_{0}}(km)-H_{P}(km))=1$. C'est \'equivalent \`a $\phi_{P_{0}}^G(H_{P}(km)-H_{\bar{P}_{0}}(km))=1$. En utilisant (9) et le fait que $\phi_{P_{0}}^G$ est la fonction caract\'eristique d'un c\^one (qui est stable par addition), on en d\'eduit la relation $\phi_{P_{0}}^G(H_{P}(km)-H_{0}(m))=1$  cherch\'ee. $\square$
 
 \bigskip
 
 \subsection{$(G,M)$-familles}
 Soit $M$ un Levi de $G$. D'apr\`es Arthur, une $(G,M)$-famille est une famille $(\varphi(\Lambda,P))_{P\in {\cal P}(M)}$ v\'erifiant les conditions suivantes:
 
 - pour tout $P$, $\Lambda\mapsto \varphi(\Lambda,P)$ est une fonction $C^{\infty}$ sur $i{\cal A}_{M}^*$;
 
 - soient $P,P'\in {\cal P}(M)$ deux paraboliques adjacents, soit ${\cal H}\subset {\cal A}_{M}^*$ le mur s\'eparant les chambres positives associ\'ees \`a $P$ et $P'$; alors $\varphi(\Lambda,P)=\varphi(\Lambda,P')$ pour $\Lambda\in i{\cal H}$.

 Pour $P\in {\cal P}(M)$, le r\'eseau engendr\'e par les coracines associ\'ees aux \'el\'ements de $\Delta_{P}$ ne d\'epend pas de $P$. On le note ${\mathbb Z}(\check{\Delta}_{M})$. On d\'efinit la fonction  m\'eromorphe $\epsilon_{P}^G$ sur ${\cal A}_{M,{\mathbb C}}^*={\cal A}_{M}^*\otimes_{{\mathbb R}}{\mathbb C}$ par 
 $$\epsilon_{P}^G(\Lambda)=mes({\cal A}_{M}^G/{\mathbb Z}[\check{\Delta}_{M}])\prod_{\alpha\in \Delta_{P}}<\Lambda,\check{\alpha}>^{-1}.$$

 {\bf Remarque.} Arthur note cette fonction $\theta_{P}^G(\Lambda)^{-1}$. On r\'eserve la lettre $\theta$ pour un autre usage.
 
  Pour une $(G,M)$-famille comme ci-dessus, on d\'efinit la fonction $\varphi_{M}^G$ sur $i{\cal A}_{M}^*$ par
 $$\varphi_{M}^G(\Lambda)=\sum_{P\in {\cal P}(M)}\varphi(\Lambda,P)\epsilon_{P}^G(\Lambda).$$
 Le point essentiel est qu'il s'agit d'une fonction $C^{\infty}$ (les singularit\'es des fonctions $\epsilon_{P}^G$ disparaissent). Dans le cas o\`u $F$ est archim\'edien, on a aussi:
 
 \ass{Lemme}{Supposons $F$ archim\'edien. Si les fonctions $\varphi(\Lambda,P)$ sont de Schwartz, $\varphi_{M}^G$ l'est aussi. Si les fonctions $\varphi(\Lambda,P)$ ainsi que leurs d\'eriv\'ees sont \`a croissance mod\'er\'ee, $\varphi_{M}^G$ et ses d\'eriv\'ees le sont aussi.}
 
 Cela r\'esulte par exemple du lemme 13.2.2 de [LW].
 
 \bigskip
 
 \subsection{Variantes des fonctions $\epsilon_{P}^G$ }

  Soient $M$ un Levi semi-standard et  $Y\in {\cal A}_{M}$. Si $F$ est non-archim\'edien, on suppose plus pr\'ecis\'ement que $Y\in {\cal A}_{M,F}\otimes_{{\mathbb Z}}{\mathbb Q}$. Soit $X\in {\cal A}_{G,F}$. On pose
  $${\cal A}_{M,F}^G(X)=\{H\in {\cal A}_{M,F}; H_{G}=X\}.$$
  Si $F$ est non-archim\'edien, on munit cet ensemble de la mesure de comptage. Si $F$ est archim\'edien, cet ensemble est \'egal \`a $X+{\cal A}_{M}^G$ et on le munit de la mesure d\'eduite de celle sur ${\cal A}_{M}^G$.
  
  Pour $P\in {\cal P}(M)$ et $\Lambda\in {\cal A}_{M,{\mathbb C}}^*$, consid\'erons l'int\'egrale
 $$\epsilon_{P}^{G,Y}(X;\Lambda)=\int_{{\cal A}_{M,F}^G(X)}\phi_{P}^G(H-Y)e^{<\Lambda,H>}\,dH.$$
 Elle est absolument convergente si $<\Lambda,\check{\alpha}>>0$ pour tout $\alpha\in \Delta_{P}$.  Dans ce domaine, elle ne d\'epend que de la projection de $\Lambda$ dans ${\cal A}_{M,{\mathbb C}}^{*}/i{\cal A}_{M,F}^{\vee}$. Si $F$ est archim\'edien, on calcule
 $$\epsilon_{P}^{G,Y}(X;\Lambda)=e^{<\Lambda,X+Y^G>}\epsilon_{P}^G(\Lambda).$$
 Supposons $F$ non archim\'edien. On peut fixer $X'\in {\cal A}_{M,F}$ tel que $X'_{G}=X$ (cf. 1.1(1)). Notons ${\cal A}_{M,F}^G={\cal A}_{M,F}^G(0)={\cal A}_{M,F}\cap {\cal A}_{M}^G$. Alors ${\cal A}_{M,F}^G(X)=X'+{\cal A}_{M,F}^G$. Par le changement de variables $H\mapsto H+X'$,
 $$\epsilon_{P}^{G,Y}(X;\Lambda)=e^{<\Lambda,X'>}\epsilon_{P}^{G,Y-X'}(0;\Lambda).$$
  Fixons un entier $k\geq1$, posons ${\cal L}_{k}=\frac{1}{k}log(q){\mathbb Z}[\check{\Delta}_{M}]$ o\`u $q$ est le nombre d'\'el\'ements du corps r\'esiduel de $F$.   Si $k$ est assez grand, ce r\'eseau contient ${\cal A}_{M,F}^G$ et $Y^G-(X')^G$. On fixe $k$ de sorte qu'il en soit  ainsi. Par inversion de Fourier, on a
 $$\epsilon_{P}^{G,Y}(X;\Lambda)= [{\cal L}_{k}:{\cal A}_{M,F}^G]^{-1}e^{<\Lambda,X'>}\sum_{\nu\in i{\cal A}_{M,F}^{G,\vee}/i{\cal L}_{k}^{\vee}}\sum_{H\in {\cal L}_{k}}\phi_{P}^G(H+X'-Y)e^{<\Lambda+\nu,H>}$$
  $$=[{\cal L}_{k}:{\cal A}_{M,F}^G]^{-1}e^{<\Lambda,X'>}\sum_{\nu\in i{\cal A}_{M,F}^{G,\vee}/i{\cal L}_{k}^{\vee}}\sum_{H\in {\cal L}_{k}}\phi_{P}^G(H)e^{<\Lambda+\nu,H+Y^G-(X')^G>}.$$
Les ensembles ${\cal A}_{M,F}^{G,\vee}$ et ${\cal L}_{k}^{\vee}$ sont des r\'eseaux dans  ${\cal A}_{M}^{G,*}$. L'ensemble des $H\in {\cal L}_{k}$ tels que $\phi_{P}^G(H)=1$ est celui des $\sum_{\alpha\in \Delta_{P}}\kappa n_{\alpha}\check{\alpha}$ pour des entiers $n_{\alpha}\leq0$, o\`u $\kappa=\frac{1}{k}log(q)$. La s\'erie
$$\sum_{H\in {\cal L}_{k}}\phi_{P}^G(H)e^{<\Lambda,H>}$$
 se calcule. Elle vaut $\epsilon_{P,k}^G(\Lambda)$, o\`u
 $$\epsilon^G_{P,k}(\Lambda)=\prod_{\alpha\in \Delta_{P}}(1-e^{-<\Lambda,\kappa\check{\alpha}>})^{-1}.$$
 On obtient
 $$(1) \qquad  \epsilon_{P}^{G,Y}(X;\Lambda)=[{\cal L}_{k}:{\cal A}_{M,F}^G]^{-1}e^{<\Lambda,X'>}\sum_{\nu\in i{\cal A}_{M,F}^{G,\vee}/i{\cal L}_{k}^{\vee}}e^{<\Lambda+\nu,Y^G-(X')^G>} \epsilon^G_{P,k}(\Lambda+\nu).$$
 Cette expression se prolonge m\'eromorphiquement \`a tout $\Lambda\in {\cal A}_{M,{\mathbb C}}^*$. 
 
 Fixons $P_{1}\in {\cal P}(M)$. Pour $P\in {\cal P}(M)$ et $\Lambda\in {\cal A}_{M,{\mathbb C}}^*$, posons
 $$\epsilon_{P,P_{1}}^{G,Y}(X;\Lambda)=\int_{{\cal A}_{M,F}^G(X)}\phi_{P,P_{1}}^G(H-Y)e^{<\Lambda,H>}\,dH.$$
 Elle est absolument convergente si $<\Lambda,\check{\alpha}>>0$ pour tout $\alpha\in \Delta_{P_{1}}$. Montrons que
 
 (2) cette fonction se prolonge m\'eromorphiquement \`a tout $\Lambda\in {\cal A}_{M,{\mathbb C}}^*$; on a l'\'egalit\'e $\epsilon_{P,P_{1}}^{G,Y}(X;\Lambda)=(-1)^{s(P,P_{1})}\epsilon_{P}^{G,Y}(X;\Lambda)$. 
 
 On traite le cas o\`u $F$ est non-archim\'edien, le cas archim\'edien \'etant plus facile. Le m\^eme calcul qui a conduit \`a l'\'egalit\'e (1) conduit \`a une \'egalit\'e similaire exprimant $\epsilon_{P,P_{1}}^{G,Y}(X;\Lambda)$. La fonction $\epsilon_{P,k}^G(\Lambda)$ y est remplac\'ee par une fonction $\epsilon_{P,P_{1},k}^G(\Lambda)$. Supposons $<\Lambda,\check{\alpha}>>0$ pour tout $\alpha\in \Delta_{P_{1}}$. Alors cette fonction est d\'efinie par 
 $$\epsilon_{P,P_{1},k}^G(\Lambda)=\sum_{H\in {\cal L}_{k}}\phi_{P,P_{1}}^G(H)e^{<\Lambda,H>}.$$
 L'ensemble des $H\in {\cal L}_{k}$ tels que $\phi_{P,P_{1}}^G(H)=1$ est celui des $\sum_{\alpha\in \Delta_{P}}\kappa n_{\alpha}\check{\alpha}$, o\`u $\kappa$ est comme ci-dessus et les entiers $n_{\alpha}$ v\'erifient
 
 - $n_{\alpha}\leq0$ si $<\Lambda,\check{\alpha}>>0$;
 
 - $n_{\alpha}>0$ si $<\Lambda,\check{\alpha}><0$.
 
 Pour un r\'eel $t\not=0$, on a les \'egalit\'es \'el\'ementaires 
 
  $\sum_{n>0}e^{nt}=-(1-e^{-t})^{-1}$, si $t<0$;
  
  $\sum_{n\leq0}e^{nt}=(1-e^{-t})^{-1}$, si $t>0$.
  
  On calcule alors $\epsilon_{P,P_{1},k}^G(\Lambda)=(-1)^{s(P,P_{1})}\epsilon_{P,k}^G(\Lambda)$. Cette \'egalit\'e se prolonge \`a tout $\Lambda$. Cela entra\^{\i}ne (2).

 \bigskip
 
 \subsection{Variantes des fonctions $\varphi_{M}^G$}
 Soient $M$ un Levi semi-standard, $X\in {\cal A}_{G,F}$,  ${\cal Y}=(Y[P])_{P\in {\cal P}(M)}$ une famille $(G,M)$-orthogonale et $(\varphi(\Lambda,P))_{P\in {\cal P}(M)}$ une $(G,M)$-famille.
 
 Supposons d'abord $F$ archim\'edien.   D\'efinissons la fonction
 $$(1) \qquad \varphi_{M}^{G,{\cal Y}}(X;\Lambda)=\sum_{P\in {\cal P}(M)}\varphi(\Lambda,P)\epsilon_{P}^{G,Y[P]}(X;\Lambda).$$
 Si l'on pose $\varphi({\cal Y};\Lambda,P)=\varphi(\Lambda,P)e^{<\Lambda,Y[P]^G>}$, la famille $(\varphi({\cal Y};\Lambda,P))_{P\in {\cal P}(M)}$ est une $(G,M)$-famille et on a l'\'egalit\'e $\varphi_{M}^{G,{\cal Y}}(X;\Lambda)=e^{<\Lambda,X>}\varphi_{M}^G({\cal Y};\Lambda)$. Cette fonction est donc $C^{\infty}$. 
 
 Supposons maintenant $F$ non-archim\'edien. On suppose que la famille ${\cal Y}$ est $p$-adique, cf. 1.3.  On suppose aussi que la famille $(\varphi(\Lambda,P))_{P\in {\cal P}(M)}$ est $p$-adique, notion que l'on d\'efinit de la fa\c{c}on suivante: on suppose que

 - pour tout $P$, la fonction $\Lambda\mapsto \varphi(\Lambda,P)$ est invariante par $ i{\cal A}_{M,F}^{\vee}$, autrement dit se descend en une fonction sur $i{\cal A}_{M,F}^*$. 
 
 On d\'efinit la fonction $\varphi_{M}^{G,{\cal Y}}(X;\Lambda)$ par l'\'egalit\'e (1).

 \ass{Lemme}{On suppose que $F$ est non-archim\'edien, que ${\cal Y}$ est une famille $(G,M)$-orthogonale $p$-adique et que $(\varphi(\Lambda,P))_{P\in {\cal P}(M)}$ est une $(G,M)$-famille $p$-adique. Alors la fonction $\varphi_{M}^{G,{\cal Y}}(X;\Lambda)$ est $C^{\infty}$ sur $i{\cal A}_{M}^*$ et invariante par $i{\cal A}_{M,F}^{\vee}$.}
 
 Preuve. L'invariance est \'evidente puisque toutes nos fonctions le sont. Reprenons le r\'eseau ${\cal L}_{k}$ de 1.5. On utilise le lemme 10.2 de [A1]. Il affirme l'existence d'une $(G,M)$-famille  $(u_{k}(\Lambda,P))_{P\in {\cal P}(M)}$, o\`u $\Lambda\mapsto u_{k}(\Lambda,P)$ appartient \`a $ C_{c}^{\infty}(i{\cal A}_{M}^*/i{\cal A}_{G}^*)$, v\'erifiant les propri\'et\'es suivantes:
 
 - pour tout $P$ et tout $\Lambda$, $\sum_{\mu\in i{\cal L}_{k}^{\vee}}u_{k}(\Lambda+\mu,P)=1$;
 
 - pour tout $P$, la fonction $v_{k}(\Lambda,P)=u_{k}(\Lambda,P)\epsilon_{P,k}^G(\Lambda)\epsilon_{P}^G(\Lambda)^{-1}$ est lisse sur $i{\cal A}_{M}^*/i{\cal A}_{G}^*$.
 
 On peut glisser la somme
$$ \sum_{\mu\in i{\cal L}_{k}^{\vee}}u_{k}(\Lambda+\nu+\mu,P)$$
dans la formule 1.5(1). En regroupant les deux sommes de la formule obtenue, on obtient
$$\epsilon_{P}^{G,Y}(X;\Lambda)= [{\cal L}_{k}:{\cal A}_{M,F}^G]^{-1}e^{<\Lambda,X'>}\sum_{\nu\in i{\cal A}_{M,F}^{G,\vee}}e^{<\Lambda+\nu,Y^G-(X')^G>}\epsilon_{P,k}^G(\Lambda+\nu)u_{k}(\Lambda+\nu,P)$$
$$= [{\cal L}_{k}:{\cal A}_{M,F}^G]^{-1}e^{<\Lambda,X'>}\sum_{\nu\in i{\cal A}_{M,F}^{G,\vee}}e^{<\Lambda+\nu,Y^G-(X')^G>}v_{k}(\Lambda+\nu,P)\epsilon_{P}^G(\Lambda+\nu).$$
Notons que la somme est localement finie d'apr\`es la compacit\'e du support de $u_{k}(\Lambda,P)$.
On en d\'eduit
$$(2) \qquad \varphi_{M}^{G,{\cal Y}}(X;\Lambda)=[{\cal L}_{k}:{\cal A}_{M,F}^G]^{-1}e^{<\Lambda,X>}\sum_{\nu\in i{\cal A}_{M,F}^{G,\vee}}e^{-<\nu,(X')^G>}$$
$$\sum_{P\in {\cal P}(M)}\varphi(\Lambda,P)e^{<\Lambda+\nu,Y[P]^G>}v_{k}(\Lambda+\nu,P)\epsilon_{P}^G(\Lambda+\nu).$$
Posons
$$\varphi(\nu,{\cal Y};\Lambda,P)=\varphi(\Lambda-\nu,P)e^{<\Lambda,Y[P]^G>}v_{k}(\Lambda,P).$$
Montrons que 

(3) $(\varphi(\nu,{\cal Y};\Lambda,P))_{P\in {\cal P}(M)}$ est une $(G,M)$-famille. 

La famille $(e^{<\Lambda,Y[P]^G>})_{P\in {\cal P}(M)}$ est une $(G,M)$-famillle. La famille $(v_{k}(\Lambda,P))_{P\in {\cal P}(M)}$ aussi: cela r\'esulte ais\'ement du fait que $(u_{k}(\Lambda,P))_{P\in {\cal P}(M)}$ en est une. Il reste \`a montrer que $(\varphi(\Lambda-\nu,P))_{P\in {\cal P}(M)}$ en est une. Or $i{\cal A}_{M,F}^{G,\vee}=(i{\cal A}_{M,F}^{\vee}+i{\cal A}_{G}^*)\cap i{\cal A}_{M}^{G,*}$. Il suffit de montrer que $(\varphi(\Lambda-\nu,P))_{P\in {\cal P}(M)}$ est une $(G,M)$-famille pour $\nu\in i{\cal A}_{M,F}^{\vee}+i{\cal A}_{G}^*$. Par hypoth\`ese, nos fonctions sont invariantes par $i{\cal A}_{M,F}^{\vee}$. On peut donc se limiter \`a $\nu\in i{\cal A}_{G}^*$. Mais il est clair que les conditions d\'efinissant une $(G,M)$-famille sont invariantes par translation par $i{\cal A}_{G}^*$. D'o\`u (3). 

Alors (2) se r\'ecrit
$$\varphi_{M}^{G,{\cal Y}}(X;\Lambda)=[{\cal L}_{k}:{\cal A}_{M,F}^G]^{-1}e^{<\Lambda,X>}\sum_{\nu\in i{\cal A}_{M,F}^{G,\vee}}e^{-<\nu,(X')^G>}\varphi_{M}^G(\nu,{\cal Y};\Lambda+\nu)$$
et ces fonctions sont $C^{\infty}$ d'apr\`es les propri\'et\'es des $(G,M)$-familles habituelles. $\square$
 
 {\bf Remarque} Pour $Z\in {\cal A}_{A_{G},F}$, on a l'\'egalit\'e
 $$\varphi_{M}^{G,{\cal Y}}(X+Z;\Lambda)=e^{<\Lambda,Z>}\varphi_{M}^{G,{\cal Y}}(X;\Lambda).$$
 En particulier, si l'on restreint la variable $\Lambda$ \`a l'ensemble $(i{\cal A}_{M}^{G,*}+i{\cal A}_{M,F}^{\vee})/i{\cal A}_{M,F}^{\vee}$, la fonction $X\mapsto \varphi_{M}^{G,{\cal Y}}(X;\Lambda).$ devient invariante par ${\cal A}_{A_{G},F}$.
 
 \bigskip
 
 \subsection{D\'eveloppement en fonction d'un param\`etre $T$}
 Si $F$ est archim\'edien, notons $PolExp$ l'ensemble des fonctions $f:{\cal A}_{0}\to {\mathbb C}$ pour lesquelles il existe une famille $(p_{\Lambda})_{\Lambda\in i{\cal A}_{0}^*}$ de polyn\^omes sur ${\cal A}_{0}$ de sorte que
 
 - l'ensemble  des $\Lambda$ tels que $p_{\Lambda}\not=0$ est fini;
 
 - pour tout $T$, on a l'\'egalit\'e
 $$f(T)=\sum_{\Lambda\in i{\cal A}_{0}^*}e^{<\Lambda,T>}p_{\Lambda}(T).$$
 
 La famille $(p_{\Lambda})_{\Lambda\in i{\cal A}_{0}^*}$ est uniquement d\'etermin\'ee. Plus pr\'ecis\'ement, la connaissance de $f$ dans un ouvert non vide de ${\cal A}_{0}$ suffit \`a d\'eterminer cette famille. En particulier, le polyn\^ome $p_{0}$ est bien d\'etermin\'e. On pose $c_{0}(f)=p_{0}(0)$.
 
 Si $\Xi\subset i{\cal A}_{0}^*$ est un ensemble fini et $N$ est un entier naturel, on note plus pr\'ecis\'ement $PolExp_{\Xi,N}$ l'ensemble des $f\in PolExp$ tels que, avec les notations ci-dessus, le degr\'e des $p_{\Lambda}$ soit inf\'erieur ou \'egal \`a $N$ et  $p_{\Lambda}$ soit nul  si $\Lambda\not\in \Xi$.

 Si $F$ est non-archim\'edien, notons $PolExp$ l'ensemble des fonctions $f:{\cal A}_{M_{0},F}\otimes_{{\mathbb Z}}{\mathbb Q}\to {\mathbb C}$ qui v\'erifient la condition suivante. Pour tout r\'eseau  ${\cal R}\subset {\cal A}_{M_{0},F}\otimes_{{\mathbb Z}}{\mathbb Q}$, il existe une famille $(p_{{\cal R},\Lambda})_{\Lambda\in i{\cal A}_{0}^*/i{\cal R}^{\vee}}$ telle que
 
 - l'ensemble  des $\Lambda$ tels que $p_{{\cal R},\Lambda}\not=0$ est fini;
 
  - pour tout $T\in {\cal R}$, on a l'\'egalit\'e
 $$f(T)=\sum_{\Lambda\in i{\cal A}_{0}^*/i{\cal R}^{\vee}}e^{<\Lambda,T>}p_{{\cal R},\Lambda}(T).$$
 
 De nouveau, la famille $(p_{{\cal R},\Lambda})_{\Lambda\in i{\cal A}_{0}^*/i{\cal R}^{\vee}}$ est uniquement d\'etermin\'ee. Plus pr\'ecis\'ement, la connaissance de $f$ dans l'intersection de ${\cal A}_{M_{0},F}\otimes_{{\mathbb Z}}{\mathbb Q} $ et d'un c\^one  ouvert non vide de ${\cal A}_{0}$ suffit \`a d\'eterminer cette famille. On pose $c_{{\cal R},0}(f)=p_{{\cal R},0}(0)$.
 
 Soit $\boldsymbol{\Xi}=(\Xi_{{\cal R}})_{{\cal R}}$ une famille index\'ee par les r\'eseaux dans ${\cal A}_{M_{0},F}\otimes_{{\mathbb Z}}{\mathbb Q}$, o\`u $\Xi_{{\cal R}}$ est un sous-ensemble fini de  $ i{\cal A}_{0}^*/i{\cal R}^{\vee}$. Soit $N$ un entier naturel. On note plus pr\'ecis\'ement $PolExp_{\boldsymbol{\Xi},N}$ l'ensemble des $f\in PolExp$ tels que, avec les notations ci-dessus et pour tout ${\cal R}$, le degr\'e des $p_{{\cal R},\Lambda}$ soit inf\'erieur ou \'egal \`a $N$ et que l'on ait $p_{{\cal R},\Lambda}=0$ si $\Lambda\not\in \Xi_{{\cal R}}$.
 
 Soient $M$ un Levi semi-standard, ${\cal Y}=(Y[P])_{P\in {\cal P}(M)}$ une famille $(G,M)$-orthogonale, $(\varphi(\Lambda,P))_{P\in {\cal P}(M)}$ une $(G,M)$-famille, $T\in {\cal A}_{0}$ et $X\in {\cal A}_{G,F}$. Dans le cas o\`u $F$ est non-archim\'edien, on suppose que les deux familles sont $p$-adiques et que $T\in {\cal A}_{M_{0},F}\otimes_{{\mathbb Z}}{\mathbb Q}$. On d\'eduit de $T$ une famille $(G,M)$-orthogonale $(T[P])_{P\in {\cal P}(M)}$, cf. 1.3. On d\'efinit la famille ${\cal Y}(T)=(Y[P]+T[P])_{P\in {\cal P}(M)}$. Elle est encore $(G,M)$-orthogonale, $p$-adique si $F$ est non-archim\'edien. On a d\'efini dans le paragraphe pr\'ec\'edent un terme $\varphi_{M}^{G,{\cal Y}(T)}(X;\Lambda)$. D'autre part, pour $P\in {\cal P}(M)$, posons $\varphi({\cal Y};\Lambda,P)=\varphi(\Lambda,P)e^{<\Lambda,Y[P]^G>}$. La famille $(\varphi({\cal Y};\Lambda,P))_{P\in {\cal P}(M)}$ est une $(G,M)$-famille. On a d\'efini en 1.4 la fonction
 $\varphi_{M}^G({\cal Y};\Lambda)$.
 
 {\bf Remarque.} On appliquera souvent ces constructions \`a la famille ${\cal Y}$ nulle, c'est-\`a-dire $Y[P]=0$ pour tout $P$. Dans ce cas on a $\varphi_{M}^G({\cal Y};\Lambda)= \varphi_{M}^G(\Lambda)$ et on note simplement $\varphi_{M}^{G,T}(X;\Lambda)=\varphi_{M}^{G,{\cal Y}(T)}(X;\Lambda)$.

 Supposons $F$ non-archim\'edien. Fixons une base de l'espace des  op\'erateurs diff\'erentiels \`a coefficients constants sur $i{\cal A}_{M}^*$, de degr\'e born\'e par $a_{M}-a_{G}$.
 Appelons norme de la $(G,M)$-famille $(\varphi(\Lambda,P))_{P\in {\cal P}(M)}$ le sup des $\vert D\varphi(\Lambda,P)\vert $, quand $\Lambda$ parcourt $i{\cal A}_{M,F}^*$, $P$ parcourt ${\cal P}(M)$ et $D$ parcourt la base fix\'ee. 
  
 \ass{Lemme}{(i) Supposons $F$ archim\'edien. Pour tout $\Lambda_{0}\in i{\cal A}_{M}^*$, la fonction
 $f:T\mapsto \varphi^{G,{\cal Y}(T)}_{M}(X;\Lambda_{0})$ appartient \`a $PolExp$. Plus pr\'ecis\'ement, il existe un entier $N$ et un sous-ensemble fini $\Xi\subset i{\cal A}_{0}^*$ ne d\'ependant que de $\Lambda_{0}$ tels que $f\in PolExp_{\Xi,N}$. On a
 $$c_{0}(f)=\left\lbrace\begin{array}{cc}e^{<\Lambda_{0},X>}\varphi_{M}^G({\cal Y};\Lambda_{0}),&\text{ si }\Lambda_{0}\in i{\cal A}_{G}^*,\\ 0,&\text{ sinon.}\\ \end{array}\right.$$
 
 (ii) Supposons $F$ non-archim\'edien. Pour tout $\Lambda_{0}\in i{\cal A}_{M}^*$, la fonction
 $f:T\mapsto \varphi^{G,{\cal Y}(T)}_{M}(X;\Lambda_{0})$ appartient \`a $PolExp$.  Plus pr\'ecis\'ement, il existe un entier $N$ et  une famille d'ensembles finis $\boldsymbol{\Xi}$ ne d\'ependant que de $\Lambda_{0}$ tels que $f\in PolExp_{\boldsymbol{\Xi},N}$. Soit ${\cal R}\subset {\cal A}_{M_0,F}\otimes_{{\mathbb Z}}{\mathbb Q}$ un r\'eseau.  Si $\Lambda_{0}\not\in i{\cal A}_{M,F}^{\vee}+i{\cal A}_{G}^*$, il existe un entier $k_{1}$ ne d\'ependant que de ${\cal R}$ tel que $c_{\frac{1}{k}{\cal R},0}(f)=0$ pour tout entier $k\geq k_{1}$. Supposons maintenant $\Lambda_{0}\in \Lambda_{1}+i{\cal A}_{M,F}^{\vee}$, o\`u $\Lambda_{1}\in i{\cal A}_{G}^*$.    Alors, il existe un r\'eel $c>0$  ne d\'ependant  que de ${\cal R}$ tel que, pour tout entier $k\geq1$, on ait la majoration
$$\vert c_{\frac{1}{k}{\cal R},0}(f)- mes(i{\cal A}_{M,F}^*)mes(i{\cal A}_{G,F}^*)^{-1}e^{<\Lambda_{1},X>}\varphi^G_{M}({\cal Y};\Lambda_{1})\vert \leq cNk^{-1},$$
o\`u $N$ est la norme de la $(G,M)$-famille $(\varphi({\cal Y};\Lambda,P))_{P\in {\cal P}(M)}$. }
  
  Preuve.  Rappelons comment on calcule un terme tel que $\varphi_{M}^G(\Lambda_{0})$. Pour tout $P\in {\cal P}(M)$, notons $\Delta_{P}(\Lambda_{0})$ l'ensemble des $\alpha\in \Delta_{P}$ tels que $<\Lambda_{0},\check{\alpha}>=0$. Notons $n_{P}(\Lambda_{0})$ le nombre d'\'el\'ements de cet ensemble. On fixe un \'el\'ement $\xi\in i{\cal A}_{M}^{G,*}$ en position g\'en\'erale. Pour $t\in {\mathbb R}$, posons 
  $$\varphi(t,\Lambda_{0},P)=\varphi(\Lambda_{0}+t\xi,P)\left(\prod_{\alpha\in \Delta_{P}-\Delta_{P}(\Lambda_{0})}<\Lambda_{0}+t\xi,\check{\alpha}>^{-1}\right)\left(\prod_{\alpha\in \Delta_{P}(\Lambda_{0})}<\xi,\check{\alpha}>^{-1}\right).$$
  Cette fonction est $C^{\infty}$ en $t=0$.   Alors
$$(1)\qquad \varphi_{M}^G(\Lambda_{0})=mes({\cal A}_{M}^G/{\mathbb Z}[\check{\Delta}_{M}])\sum_{P\in {\cal P}(M)}(n_{P}(\Lambda_{0})!)^{-1}\left(\frac{d^{n_{P}(\Lambda_{0})}}{dt^{n_{P}(\Lambda_{0})}}\varphi(t,\Lambda_{0},P)\right) _{t=0}.$$ 
Cela r\'esulte simplement de l'\'egalit\'e
$$\varphi_{M}^G(\Lambda_{0})=lim_{t\to 0}\sum_{P\in {\cal P}(M)}\varphi(\Lambda_{0}+t\xi,P)\epsilon_{P}^G(\Lambda_{0}+t\xi).$$

  Preuve de (i). On a  dit  en 1.6 que l'on avait l'\'egalit\'e $ f(T)=e^{<\Lambda,X>}\varphi_{M}^G({\cal Y}(T);\Lambda_{0})$.   La formule (1) montre que $\varphi_{M}^G({\cal Y}(T);\Lambda_{0})$  est combinaison lin\'eaire de termes $D\varphi({\cal Y}(T);\Lambda_{0},P)$ o\`u $P$ parcourt ${\cal P}(M)$ et $D$ parcourt les op\'erateurs diff\'erentiels sur $i{\cal A}_{M}^*$, \`a coefficients constants et de degr\'e born\'e par $a_{M}-a_{G}$. Comme fonction de $T$,  un tel terme est de la forme $e^{<\Lambda_{0},T[P]^G>}p(T)$, o\`u $p$ est un polyn\^ome de degr\'e inf\'erieur ou \'egal \`a $a_{M}-a_{G}$. Cette fonction appartient \`a $PolExp$, donc $f$ aussi. Plus pr\'ecis\'ement, $f\in PolExp_{\Xi,a_{M}-a_{G}}$ o\`u $\Xi=\{w\Lambda_{0}^G; w\in W^G\}$.  Le coefficient $c_{0}(f)$ ne voit que les termes dont la partie  exponentielle est triviale.   Il n'y en a que si que si $\Lambda_{0}^G=0$. Supposons cette condition v\'erifi\'ee. Alors toutes les exponentielles sont triviales et $c_{0}(f)$ est simplement $f(0)$. Mais  on a $f(0)=e^{<\Lambda_{0},X>}\varphi_{M}^G({\cal Y};\Lambda_{0})$. D'o\`u (i).
  
  Preuve de (ii). Soit  ${\cal R}$ un r\'eseau dans ${\cal A}_{M_0,F}\otimes_{{\mathbb Z}}{\mathbb Q}$. On a d\'efini en 1.5 des r\'eseaux ${\cal L}_{k}$ pour tout entier $k\geq1$. On peut fixer $k_{0}$ de sorte que $T[P]^G\in {\cal L}_{k_{0}}$ pour tout $T\in {\cal R}$. On a alors l'\'egalit\'e 1.6(2) que l'on r\'ecrit
  $$f(T) =[{\cal L}_{k_{0}}:{\cal A}_{M,F}^G]^{-1}e^{<\Lambda_{0},X>}\sum_{\nu\in i{\cal A}_{M,F}^{G,\vee}}e^{-<\nu,(X')^G>}\varphi_{M}^G(\nu,{\cal Y}(T);\Lambda_{0}+\nu),$$
  o\`u, pour tout $P\in {\cal P}(M)$, 
  $$\varphi(\nu,{\cal Y}(T);\Lambda,P)=\varphi(\Lambda-\nu,P)e^{<\Lambda,Y[P]^G+T[P]^G>}v_{k_{0}}(\Lambda,P).$$
  La somme en $\nu$ est finie. En appliquant (1), on obtient
  $$(2) \qquad f(T)=mes({\cal A}_{M}^G/{\mathbb Z}[\check{\Delta}_{M}])[{\cal L}_{k_{0}}:{\cal A}_{M,F}^G]^{-1}e^{<\Lambda_{0},X>}\sum_{\nu\in i{\cal A}_{M,F}^{G,\vee}}e^{-<\nu,(X')^G>}$$
  $$\sum_{P\in {\cal P}(M)}(n_{P}(\Lambda_{0}+\nu)!)^{-1}\left(\frac{d^{n_{P}(\Lambda_{0}+\nu)}}{dt^{n_{P}(\Lambda_{0}+\nu)}}\varphi(\nu,{\cal Y}(T);t,\Lambda_{0}+\nu,P)\right) _{t=0}.$$ 
  Comme dans le cas archim\'edien une fonction 
  $$T\mapsto \left(\frac{d^{n_{P}(\Lambda_{0}+\nu)}}{dt^{n_{P}(\Lambda_{0}+\nu)}}\varphi(\nu,{\cal Y}(T);t,\Lambda_{0}+\nu,P)\right) _{t=0}$$
  est produit de  $e^{<\Lambda_{0}+\nu,T[P]^G>}$ et d'un polyn\^ome en $T$ de degr\'e inf\'erieur ou \'egal \`a $a_{M}-a_{G}$. Une telle fonction appartient \`a $PolExp$, donc la fonction $f$  appartient aussi \`a cet espace.  Plus pr\'ecis\'ement, $f\in PolExp_{\boldsymbol{\Xi},a_{M}-a_{G}}$, o\`u $\boldsymbol{\Xi}=(\Xi_{{\cal R}})_{{\cal R}}$ est la famille telle que $\Xi_{{\cal R}}$ soit l'ensemble des projections dans $i{\cal A}_{0}^*/i{\cal R}^{\vee}$ des $w(\Lambda_{0}^G+\nu)$ pour $w\in W^G$ et $\nu\in i{\cal A}_{M,F}^{G,\vee}$. Pour tout $P\in {\cal P}(M)$, notons ${\cal S}(P)$ l'image de ${\cal R}$ par l'application $T\mapsto T[P]^G$. C'est un r\'eseau de ${\cal A}_{M}^G$.  Le coefficient $c_{{\cal R},0}(f)$ s\'electionne les termes de la formule (2) pour lesquels les exponentielles sont triviales pour tout $T\in {\cal R}$, c'est-\`a-dire les couples $(\nu,P)$ tels que $\Lambda_{0}^G+\nu\in i{\cal S}(P)^{\vee}$.  Supposons $\Lambda_{0}\not\in i{\cal A}_{M,F}^{\vee}+i{\cal A}_{G}^*= i{\cal A}_{M,F}^{G,\vee}\oplus i{\cal A}_{G}^*$.  Si  ${\cal R}$ est assez grand, on a ${\cal S}(P)^{\vee}\subset {\cal A}_{M,F}^{G,\vee}$  pour tout $P$ et aucun $\nu\in i{\cal A}_{M,F}^{G,\vee}$ ne contribue. Donc $c_{{\cal R},0}(f)=0$ ce qui prouve la premi\`ere assertion du (ii). Supposons maintenant $\Lambda_{0}\in\Lambda_{1}+i{\cal A}_{M,F}^{\vee}$, o\`u $\Lambda_{1}\in i{\cal A}_{G}^*$. On a les \'egalit\'es 
   $$\varphi_{M}^{G,{\cal Y}(T)}(X;\Lambda_{0})=\varphi_{M}^{G,{\cal Y}(T)}(X;\Lambda_{1})=e^{<\Lambda_{1},X>}(\varphi')_{M}^{G,{\cal Y}(T)}(X;0),$$
   o\`u $(\varphi'(\Lambda,P))_{P\in {\cal P}(M)}$ est la $(G,M)$-famille d\'eduite de la famille initiale par translation par $\Lambda_{1}$. Cela nous ram\`ene au cas o\`u $\Lambda_{0}=0$, ce que l'on suppose d\'esormais. Le calcul pr\'ec\'edent conduit \`a l'\'egalit\'e
  $$c_{{\cal R},0}(f)=mes({\cal A}_{M}^G/{\mathbb Z}[\check{\Delta}_{M}])[{\cal L}_{k_{0}}:{\cal A}_{M,F}^G]^{-1}\sum_{P\in {\cal P}(M)}$$
  $$\sum_{\nu\in i{\cal S}(P)^{\vee}}e^{-<\nu,(X')^G>}(n_{P}(\nu)!)^{-1}\left(\frac{d^{n_{P}(\nu)}}{dt^{n_{P}(\nu)}}\varphi(\nu,{\cal Y};t,\nu,P)\right) _{t=0}.$$ 
 Rempla\c{c}ons ${\cal R}$ par $\frac{1}{k}{\cal R}$. Cela  remplace  ${\cal S}(P)^{\vee}$ par $k{\cal S}(P)^{\vee}$.  On remarque que $\Delta_{P}(k\nu)=\Delta_{P}(\nu)$ et  $n_{P}(k\nu)=n_{P}(\nu)$. On peut remplacer l'entier $k_{0}$ par $kk_{0}$. On peut choisir $u_{kk_{0}}(\Lambda,P)=u_{k_{0}}(\frac{\Lambda}{k},P)$. Alors $v_{kk_{0}}(\Lambda,P)=k^{a_{M}-a_{G}}v_{k_{0}}(\frac{\Lambda}{k},P)$. Puisqu'on a aussi
  $$[{\cal L}_{kk_{0}}:{\cal A}_{M,F}^G]=k^{a_{M}-a_{G}}[{\cal L}_{k_{0}}:{\cal A}_{M,F}^G],$$
  l'\'egalit\'e ci-dessus devient
   $$(3) \qquad c_{\frac{1}{k}{\cal R},0}(f)=mes({\cal A}_{M}^G/{\mathbb Z}[\check{\Delta}_{M}])[{\cal L}_{k_{0}}:{\cal A}_{M,F}^G]^{-1}\sum_{P\in {\cal P}(M)}$$
   $$\sum_{\nu\in i{\cal S}(P)^{\vee}}e^{-<k\nu,(X')^G>}(n_{P}(\nu)!)^{-1}\left(\frac{d^{n_{P}(\nu)}}{dt^{n_{P}(\nu)}} f_{k}(t,\nu,P)\right) _{t=0},$$
   o\`u 
   $$f_{k}(t,\nu,P)=\varphi(t\xi,P)e^{<k\nu+t\xi,Y[P]^G>}v_{k_{0}}(\nu+\frac{t}{k}\xi,P)\left(\prod_{\alpha\in \Delta_{P}-\Delta_{P}(\nu)}<k\nu+t\xi,\check{\alpha}>^{-1}\right)$$
   $$\left(\prod_{\alpha\in \Delta_{P}(\nu)}<\xi,\check{\alpha}>^{-1}\right).$$
   Notons que puisque $v_{k_{0}}(\Lambda,P)$ est \`a support compact, la somme est limit\'ee \`a un ensemble fini ind\'ependant de $k$. Fixons $\nu$ et $P$. Si $k$ est assez grand, on a $e^{-<k\nu,(X')^G>}=e^{<k\nu,Y[P]^G>}=1$ et ces termes disparaissent. Dans la d\'efinition de $f_{k}(t,\nu,P)$, le produit $\varphi(t\xi,P)e^{<k\nu+t\xi,Y[P]^G>}$ devient $\varphi({\cal Y};t\xi,P)$. Le terme $\left(\frac{d^{n_{P}(\nu)}}{dt^{n_{P}(\nu)}} f_{k}(t,\nu,P)\right) _{t=0}$ est combinaison lin\'eaire de produits
   $$\left(\frac{d^{a}}{dt^{a}}\varphi({\cal Y};t\xi,P)\right)_{t=0}\left(\frac{d^b}{dt^b}v_{k_{0}}(\nu+\frac{t}{k}\xi,P)\right)_{t=0}\left(\frac{d^{c}}{dt^c}\prod_{\alpha\in \Delta_{P}-\Delta_{P}(\nu)}<k\nu+t\xi,\check{\alpha}>^{-1}\right)_{t=0}$$
   pour des entiers naturels $a,b,c$ tels que $a+b+c=n_{P}(\nu)$. On voit qu'un tel produit est born\'e par $CNk^{-m}$, o\`u $N$ est la norme de la $(G,M)$-famille $(\varphi({\cal Y};\Lambda,P))_{P\in {\cal P}(M)}$, $C$ ne d\'epend  que de $m$ et $v_{k_{0}}$, et $m=b+c+\vert \Delta_{P}-\Delta_{P}(\nu)\vert $. Les termes pour lesquels $m$ est strictement positif sont n\'egligeables pour la conclusion du lemme. Si $\nu\not=0$, on a $\Delta_{P}\not=\Delta_{P}(\nu)$ et tous les termes sont donc n\'egligeables. Pour $\nu=0$, on ne conserve que les termes o\`u les d\'erivations ne s'appliquent qu'\`a la fonction $\varphi({\cal Y};t\xi,P)$.  Si on \'elimine les termes n\'egligeables, le membre de droite de (3) devient    
   $$mes({\cal A}_{M}^G/{\mathbb Z}[\check{\Delta}_{M}])[{\cal L}_{k_{0}}:{\cal A}_{M,F}^G]^{-1}\sum_{P\in {\cal P}(M)}v_{k_{0}}(0,P)n!^{-1}\left(\frac{d^n}{dt^n}\varphi({\cal Y};t\xi,P)\right)_{t=0}\prod_{\alpha\in \Delta_{P}}<\xi,\check{\alpha}>^{-1},$$
   o\`u $n=a_{M}-a_{G}$.  On va montrer que
  
  (4) on a l'\'egalit\'e $ [{\cal L}_{k_{0}}:{\cal A}_{M,F}^G]^{-1}v_{k_{0}}(0,P)= mes(i{\cal A}_{M,F}^*)mes(i{\cal A}_{G,F}^*)^{-1}$. 
  
  En admettant cette relation et en utilisant (1), l'expression ci-dessus n'est autre que $ \varphi_{M}^G({\cal Y};0)$ et on obtient la conclusion de (ii).
  
  Il reste \`a prouver (4). Pour $\nu\in i{\cal L}_{k_{0}}^{\vee}$, $\nu\not=0$, la fonction $t\mapsto \epsilon_{P,{\cal L}_{k_{0}}}^G(\nu+t\xi)$ a un p\^ole d'ordre $n=a_{M}-a_{G}$ en $t=0$. La fonction $t\mapsto \epsilon_{P}^G(\nu+t\xi)^{-1}$ a un z\'ero d'ordre strictement inf\'erieur. Puisque $v_{k_{0}}(\nu+t\xi,P)=u_{k_{0}}(\nu+t\xi,P)\epsilon_{P,{\cal L}_{k_{0}}}^G(\nu+t\xi)\epsilon_{P}^G(\nu+t\xi)^{-1}$ est r\'eguli\`ere, on a n\'ecessairement $u_{k_{0}}(\nu,P)=0$. Puisque $\sum_{\nu\in i{\cal L}_{k_{0}}^{\vee}}u_{k_{0}}(\nu,P)=1$, on en d\'eduit $u_{k_{0}}(0,P)=1$. On calcule alors
  $$v_{k_{0}}(0,P)=mes({\cal A}_{M}^G/{\mathbb Z}[\check{\Delta}_{M}])^{-1}\lim_{\Lambda\to 0}\prod_{\alpha\in \Delta_{P}}\frac{<\Lambda,\check{\alpha}>}{1-e^{-<\Lambda,\kappa\check{\alpha}>}}=mes({\cal A}_{M}^G/{\mathbb Z}[\check{\Delta}_{M}])^{-1}\kappa^{-n},$$
  o\`u $\kappa=\frac{1}{k_{0}}log(q)$. Or ${\cal L}_{k_{0}}=\kappa {\mathbb Z}[\check{\Delta}_{M}]$, d'o\`u
  $$\kappa^n mes({\cal A}_{M}^G/{\mathbb Z}[\check{\Delta}_{M}])=mes({\cal A}_{M}^G/\kappa{\mathbb Z}[\check{\Delta}_{M}])=mes({\cal A}_{M}^G/{\cal L}_{k_{0}}).$$
 D'apr\`es les relations de compatibilit\'e  de nos mesures, on a 
 $$mes({\cal A}_{M}^G/{\cal A}_{M,F}^G)=mes(i{\cal A}_{M,F}^*)^{-1}mes(i{\cal A}_{G,F}^*),$$
  d'o\`u $mes({\cal A}_{M}^G/{\cal L}_{k_{0}})=[{\cal L}_{k_{0}}:{\cal A}_{M,F}^G]^{-1}mes(i{\cal A}_{M,F}^*)^{-1}mes(i{\cal A}_{G,F}^*)$. En assemblant ces calculs, on obtient (4). $\square$

 \bigskip
 
 \subsection{$(G,M)$-familles et inversion de Fourier}
 Soient $M$ un Levi semi-standard et $X$ un \'el\'ement de ${\cal A}_{G,F}$. Appliquons les constructions du paragraphe 1.6 \`a la famile ${\cal Y}=(Y[P])_{P\in {\cal P}(M)}$ nulle, c'est-\`a-dire $Y[P]=0$ pour tout $P$, et \`a la famille $(\varphi(\Lambda,P))_{P\in {\cal P}(M)}$ form\'ee des fonctions constantes \'egales \`a $1$. On en d\'eduit une fonction que l'on note
 $$\epsilon_{M}^{G,0}(X;\Lambda)=\sum_{P\in {\cal P}(M)}\epsilon_{P}^{G,0}(X;\Lambda).$$ 
 On a l'\'egalit\'e
 $$(1)\qquad \epsilon_{M}^{G,0}(X;\Lambda)=\left\lbrace\begin{array}{cc}e^{<\Lambda,X>},& \text{ si }F \text{ est non-archim\'edien et }X\in {\cal A}_{M,F}\cap {\cal A}_{G},\\&\text{ ou si }F\text{ est archim\'edien et }M=G,\\ 0,&\text{ sinon.}\end{array}\right.$$.
 
 Preuve. On peut fixer un \'el\'ement $P_{1}\in {\cal P}(M)$ et supposer $<\Lambda,\check{\alpha}>>0$ pour tout $\alpha\in \Delta_{P_{1}}$. En utilisant 1.5(2), on a
 $$\epsilon_{M}^{G,0}(X;\Lambda)=\sum_{P\in {\cal P}(M)}(-1)^{s(P,P_{1})}\epsilon_{P,P_{1}}^{G,0}(X;\Lambda)$$
 $$=\sum_{P\in {\cal P}(M)}(-1)^{s(P,P_{1})}\int_{{\cal A}_{M,F}^G(X)}\phi_{P,P_{1}}^G(H)e^{<\Lambda,H>}\,dH$$
 $$=\int_{{\cal A}_{M,F}^G(X)}e^{<\Lambda,H>}\sum_{P\in {\cal P}(M)}(-1)^{s(P,P_{1})}\phi_{P,P_{1}}^G(H) dH.$$
 D'apr\`es 1.3(7), la somme int\'erieure est \'egale \`a $\Gamma_{M}^G(H,{\cal Y})$, o\`u  ${\cal Y}$ est la famille nulle. Ce n'est autre que $\delta_{M}^G(H)$.  Si $F$ est archim\'edien et $M\not=G$, l'int\'egrale est nulle. Si $F$ est non archim\'edien, elle vaut $e^{<\Lambda,X>}$ s'il existe $H\in {\cal A}_{M,F}^G(X)$ avec $H^G=0$ et elle vaut $0$ sinon. Cette condition \'equivaut \`a $X\in {\cal A}_{M,F}\cap {\cal A}_{G}$. $\square$

 Soit $T$ un \'el\'ement de ${\cal A}_{0}$ et soit  $(\varphi(\Lambda,P))_{P\in {\cal P}(M)}$ une $(G,M)$-famille. Dans le cas o\`u $F$ est non-archim\'edien, on suppose que $T\in {\cal A}_{M_{0},F}\otimes_{{\mathbb Z}}{\mathbb Q}$ et que la $(G,M)$-famille est $p$-adique. Dans le cas archim\'edien, on suppose que les fonctions $\Lambda\mapsto \varphi(\Lambda,P)$ sont de Schwartz. Pour un sous-groupe parabolique $Q=LU_{Q}\in {\cal F}(M)$, on d\'efinit $\varphi(\Lambda,Q)$ pour $\Lambda\in i{\cal A}_{L,F}^*$ comme la restriction \`a $ i{\cal A}_{L,F}^*$ de $\varphi(\Lambda,P)$ pour n'importe quel $P\in {\cal P}(M)$ tel que $P\subset Q$.   On d\'efinit sa transform\'ee de Fourier $\hat{\varphi}(H,Q)$, qui est une fonction de $H\in {\cal A}_{L,F}$, par
 $$\hat{\varphi}(H,Q)=mes(i{\cal A}_{L,F}^*)^{-1}\int_{i{\cal A}_{L,F}^*}\varphi(\Lambda,Q)e^{-<\Lambda,H>}\,d\Lambda.$$
 C'est une fonction de Schwartz sur ${\cal A}_{L,F}$.

 \ass{Lemme}{L'expression 
 $$\sum_{Q=LU_{Q}\in {\cal F}(M)}\int_{ {\cal A}_{L,F}}\int_{{\cal A}_{M,F}^G(X+H_{G})}\delta_{M}^Q(H')\Gamma_{Q}^G(H',H+T[Q])\hat{\varphi}(H,Q) e^{<\Lambda,H'>}\,dH'\,dH$$
 est convergente et \'egale \`a $\varphi_{M}^{G,T}(X;\Lambda)$.}

 Preuve.     La condition $\Gamma_{Q}^G(H',H+T[Q])\not=0$ implique une majoration $\vert (H')^G_{L}\vert <<1+\vert H^G\vert $ ($T$ est ici consid\'er\'e comme une constante). Jointe aux conditions $\delta_{M}^Q(H')=1$ et  $H_{G}'-H_{G}=X$, cela implique $\vert H'\vert <<1+\vert H\vert $. L'int\'egrale  en $H'$ est donc essentiellement born\'ee par $(1+\vert H\vert )^D$ pour un entier $D$ convenable. Puisque $H\mapsto \hat{\varphi}(H,Q)$ est de Schwartz, l'int\'egrale en $H$ est convergente. Cela prouve la convergence \'enonc\'ee.
 
 Rappelons que, pour $\Lambda$ en position g\'en\'erale,
 $$\varphi_{M}^{G,T}(X;\Lambda)=\sum_{P\in {\cal P}(M)}\varphi(\Lambda,P)\epsilon^{G,T[P]}_{P}(X;\Lambda).$$
 Fixons $P$. Par inversion de Fourier
 $$\varphi(\Lambda,P)=\int_{{\cal A}_{M,F}}\hat{\varphi}(H,P)e^{<\Lambda,H>}\,dH.$$
 On va prouver l'\'egalit\'e
 $$(2) \qquad \epsilon^{G,T[P]}_{P}(X;\Lambda)e^{<\Lambda,H>}=\sum_{Q=LU_{Q}; P\subset Q}\int_{{\cal A}_{L,F}^G(X+H_{G})}\Gamma_{Q}^G(H',H+T[Q])\epsilon_{P}^{Q,0}(H';\Lambda)\,dH'.$$
 En admettant cela, on obtient
 $$\varphi_{M}^{G,T}(X;\Lambda)=\sum_{P\in {\cal P}(M)} \int_{{\cal A}_{M,F}}\hat{\varphi}(H,P)\epsilon^{G,T[P]}_{P}(X;\Lambda)e^{<\Lambda,H>}\,dH$$
 $$=\sum_{P\in {\cal P}(M)} \int_{{\cal A}_{M,F}}\hat{\varphi}(H,P)\sum_{Q=LU_{Q}; P\subset Q}\int_{{\cal A}_{L,F}^G(X+H_{G})}\Gamma_{Q}^G(H',H+T[Q])\epsilon_{P}^{Q,0}(H';\Lambda)\,dH'\,dH.$$
 Cette expression est absolument convergente pour $\Lambda$ en position g\'en\'erale pour les m\^emes raisons que ci-dessus, la fonction $H'\mapsto \epsilon_{P}^{Q,0}(H';\Lambda)$ \'etant born\'ee. On sort la somme en $Q$ des int\'egrales et  on la permute avec celle en $P$.   On obtient
$$\varphi_{M}^{G,T}(X;\Lambda)=\sum_{Q=LU_{Q}\in {\cal F}(M)}\sum_{P\in {\cal P}(M); P\subset Q} \int_{{\cal A}_{M,F}}\hat{\varphi}(H,P)$$
$$\int_{{\cal A}_{L,F}^G(X+H_{G})}\Gamma_{Q}^G(H',H+T[Q])\epsilon_{P}^{Q,0}(H';\Lambda)\,dH'\,dH.$$
On remarque que $H$ n'intervient que par sa projection $H_{L}$ dans l'int\'egrale int\'erieure. D\'ecomposons la premi\`ere int\'egrale en une int\'egrale sur $H\in {\cal A}_{L,F}$ et une int\'egrale sur ${\cal A}_{M,F}^L(H)$. On obtient
$$(3) \qquad \varphi_{M}^{G,T}(X;\Lambda)=\sum_{Q=LU_{Q}\in {\cal F}(M)}\sum_{P\in {\cal P}(M); P\subset Q} \int_{{\cal A}_{L,F}}\int_{{\cal A}_{M,F}^L(H)}\hat{\varphi}(H'',P)$$
$$\int_{{\cal A}_{L,F}^G(X+H_{G})}\Gamma_{Q}^G(H',H+T[Q])\epsilon_{P}^{Q,0}(H';\Lambda)\,dH'\,dH''\,dH.$$
On voit appara\^{\i}tre l'int\'egrale
$$\int_{{\cal A}_{M,F}^L(H)}\hat{\varphi}(H'',P)\,dH''.$$
Par inversion de Fourier, ceci n'est autre que
$$mes(i{\cal A}_{L,F}^*)^{-1}\int_{i{\cal A}_{L,F}^*}\varphi(\Lambda,P)e^{-<\Lambda,H>}\,d\Lambda.$$
Puisque le domaine d'int\'egration est restreint \`a $i{\cal A}_{L,F}^*$, on peut aussi bien remplacer $\varphi(\Lambda,P)$ par $\varphi(\Lambda,Q)$. L'int\'egrale ci-dessus devient $\hat{\varphi}(H,Q)$, en particulier est ind\'ependante de $P$. Dans la formule (3), la somme en $P$ ne porte plus que sur les termes $\epsilon_{P}^{Q,0}(H';\Lambda)$. Leur somme vaut $\epsilon_{M}^{Q,0}(H';\Lambda)$. Supposons $F$ non archim\'edien. L'\'egalit\'e (1) dit que $\epsilon_{M}^{Q,0}(H';\Lambda)$ vaut $e^{<\Lambda,X>}$ si $H'\in {\cal A}_{M,F}$, $0$ sinon. La formule (3)  devient
$$\varphi_{M}^{G,T}(X;\Lambda)=\sum_{Q=LU_{Q}\in {\cal F}(M)}\int_{ {\cal A}_{L,F}}\int_{{\cal A}_{L,F}^G(X+H_{G})\cap {\cal A}_{M,F}}\Gamma_{Q}^G(H',H+T[Q])\hat{\varphi}(H,Q) e^{<\Lambda,H'>}\,dH'\,dH.$$
Mais ${\cal A}_{L,F}^G(X+H_{G})\cap {\cal A}_{M,F}$ est l'ensemble des $H'\in {\cal A}_{M,F}^G(X+H_{G})$ tels que $\delta_{M}^Q(H')=1$. La formule ci-dessus est donc \'equivalente \`a celle de l'\'enonc\'e. Supposons maintenant $F$ archim\'edien. L'\'egalit\'e (1) transforme directement la formule (3) en celle de l'\'enonc\'e, \`a ceci pr\`es que l'on ne somme que sur les $Q\in {\cal P}(M)$. Mais si $Q\not\in {\cal P}(M)$, le support de la fonction $\delta_{M}^Q$ est de mesure nulle et la contribution de $Q$ \`a la formule de l'\'enonc\'e est nulle.

  Prouvons (2). Puisqu'il s'agit de fonctions m\'eromorphes en $\Lambda$, on peut supposer $Re<\Lambda,\check{\alpha}>>0$ pour tout $\alpha\in \Delta_{P}$. Alors
 $$\epsilon^{G,T}_{P}(X;\Lambda)e^{<\Lambda,H>}=\int_{ {\cal A}_{M,F}^G(X)}\phi_{P}^G(H'-T)e^{<\Lambda,H'+H>}\,dH'$$
 $$=\int_{ {\cal A}_{M,F}^G(X+H_{G})}\phi_{P}^G(H'-H-T)e^{<\Lambda,H'>}\,dH'.$$
 On utilise 1.3(3), d'o\`u
 $$\epsilon^{G,T}_{P}(X;\Lambda)=\int_{ {\cal A}_{M,F}^G(X+H_{G})}\sum_{Q=LU_{Q}; P\subset Q}\Gamma_{Q}^G(H',H+T)\phi_{P}^Q(H')e^{<\Lambda,H'>}\,dH'.$$
 On sort la somme en $Q$ de l'int\'egrale et on d\'ecompose celle-ci en une int\'egrale sur $H'\in {\cal A}_{L,F}^{G}(X+H_{G})$  et une int\'egrale sur $H''\in {\cal A}_{M,F}^{L}(H')$. On obtient
 $$\epsilon^{G,T}_{P}(X;\Lambda)=\sum_{Q=LU_{Q}; P\subset Q}\int_{{\cal A}_{L,F}^{G}(X+H_{G})}\Gamma_{Q}^G(H',H+T)\int_{{\cal A}_{M,F}^L(H')}\phi_{P}^Q(H'')e^{<\Lambda,H''>}\,dH''\,dH'.$$
 La derni\`ere int\'egrale n'est autre que $\epsilon_{P}^{Q,0}(H';\Lambda)$ et on obtient la formule (2). $\square$

 \bigskip
 
 \subsection{Repr\'esentations}
 Si $\pi$ est une repr\'esentation admissible de $G(F)$ et si $F$ est non-archim\'edien, on note $V_{\pi}$ un espace dans lequel elle se r\'ealise. Dans le cas o\`u $F$ est archim\'edien, pour \^etre correct, il faudrait introduire deux espaces: un espace dans lequel $\pi$ se r\'ealise et le sous-espace des vecteurs $K$-finis, qui est celui qui nous int\'eresse mais qui n'est pas invariant par l'action de $G(F)$. Pour simplifier, on n'introduit que ce dernier, que l'on note $V_{\pi}$, et on convient qu'il est plong\'e dans un espace implicite plus gros o\`u $\pi$ se r\'ealise.   Si $\pi$ est irr\'eductible et unitaire, on fixe un produit hermitien $(.,.)$ sur $V_{\pi}$ invariant par $\pi(g)$ pour tout $g\in G(F)$ (avec la convention ci-dessus dans le cas archim\'edien). Pour nous, un tel produit est lin\'eaire en la seconde variable et antilin\'eaire en la premi\`ere. On note $\Pi_{disc}(G(F))$ l'ensemble des classes de repr\'esentations irr\'eductibles de la s\'erie discr\`ete de $G(F)$. Pour $\sigma\in \Pi_{disc}(G(F))$, son degr\'e formel $d(\sigma)$ est d\'efini par l'\'egalit\'e
 $$ \int_{A_{G}(F)\backslash G(F)}(e_{1},\sigma(g)e_{2})(\sigma(g)e'_{1},e'_{2})\,dg\,=d(\sigma)^{-1}(e_{1},e'_{2})(e'_{1},e_{2})$$
 pour tous $e_{1},e_{2},e'_{1},e'_{2}\in V_{\sigma}$.
 
 Le groupe  ${\cal A}_{G,{\mathbb C}}^*$   agit sur l'ensemble des classes d'isomorphisme de repr\'esentations admissibles de $G(F)$: \`a une repr\'esentation $\pi$ et \`a $\lambda\in {\cal A}_{G,{\mathbb C}}^* $, on associe la repr\'esentation $\pi_{\lambda}$ d\'efinie par $\pi_{\lambda}(g)=e^{<\lambda,H_{G}(g)>}\pi(g)$. Cette action se quotiente en une action de 
 ${\cal A}_{G,{\mathbb C}}^*/i{\cal A}_{G,F}^{\vee}$. Pour deux repr\'esentations $\pi$, $\pi'$, consid\'erons  l'ensemble des $\lambda\in i{\cal A}_{G}^*$ tels que $\pi\simeq \pi'_{-\lambda}$. Il est stable par $i{\cal A}_{G,F}^{\vee}$ agissant par addition et on note $[\pi,\pi']$ son quotient par cette action. Si $\pi=\pi'$, c'est un groupe et  il est plus suggestif de le noter $Stab(i{\cal A}_{G,F}^*,\pi)$. Ce groupe ne d\'epend que de l'orbite de $\pi$ pour l'action de $i{\cal A}_{G,F}^*$. L'action du groupe $i{\cal A}_{G,F}^*$ conserve  l'ensemble $\Pi_{disc}(G(F))$. On note $\Pi_{disc}(G(F))/i{\cal A}_{G,F}^*$ l'ensemble des orbites.  
  
 Soit $P=MU_{P}$ un sous-groupe parabolique semi-standard de $G$ et soit $\sigma$ une repr\'esentation admissible de $M(F)$. Pour $\lambda\in {\cal A}_{M,{\mathbb C}}^*/i{\cal A}_{M,F}^{\vee}$, introduisons la repr\'esentation induite $\pi_{\lambda}=Ind_{P}^G(\sigma_{\lambda})$ de $G(F)$.   Elle se r\'ealise naturellement dans un espace $V_{\pi_{\lambda}}$ de fonctions $e:G(F)\to V_{\sigma}$ se transformant \`a gauche par $P(F)$ selon la formule usuelle.  On d\'efinit l'espace $V_{\sigma,P}$ des fonctions $e:K\to V_{\sigma}$ telles que
 
 - $e(muk)=\sigma(m)e(k)$ pour tous $m\in M(F)\cap K$, $u\in U(F)\cap K$, $k\in K$;
 
 - si $F$ est non archim\'edien, $e$ est localement constante;
 
 - si $F$ est archim\'edien, $e$ est lisse et $K$-finie.
 
 Par l'homomorphisme  $V_{\pi_{\lambda}}\to V_{\sigma,P}$ de restriction \`a $K$, $\pi_{\lambda}$ se r\'ealise  dans $V_{\sigma,P}$ pour tout $\lambda$. Supposons $\sigma$ irr\'eductible et unitaire.  Si $\lambda\in i{\cal A}_{M,F}^*$,  $\pi_{\lambda}$ est unitaire. Plus pr\'ecis\'ement, $\pi_{\lambda}$ conserve le produit hermitien sur $V_{\sigma,P}$  d\'efini par 
 $$(e_{1},e_{2})=\int_{K}(e_{1}(k),e_{2}(k))\,dk.$$
 
  Revenons \`a $\sigma$ admissible et $\lambda\in {\cal A}_{M,{\mathbb C}}^*/i{\cal A}_{M,F}^{\vee}$. Soit $P'=MU_{P'}\in {\cal P}(M)$. Posons $\pi'_{\lambda}=Ind_{P'}^G(\sigma_{\lambda})$. On d\'efinit   l'op\'erateur d'entrelacement
  $$J_{P'\vert P}(\sigma_{\lambda}):V_{\sigma,P}\to V_{\sigma,P'}.$$
  Modulo les isomorphismes 
  $$\begin{array}{ccccccc}V_{\sigma,P}&\to&V_{\pi_{\lambda}}&&V_{\sigma,P'}&\to&V_{\pi'_{\lambda}}\\e&\mapsto &e_{\lambda}&&e'&\mapsto&e'_{\lambda},\\ \end{array}$$
   il existe un r\'eel $c$ tel que, pour $Re(<\lambda,\check{\alpha}>)>c$ pour tout $\alpha\in \Delta_{P}$,  cet op\'erateur soit  donn\'e par la formule
  $$<\epsilon,(J_{P'\vert P}(\sigma_{\lambda})(e))_{\lambda}(g)>=\int_{(U_{P}(F)\cap U_{P'}(F))\backslash U_{P'}(F)}<\epsilon,e_{\lambda}(ug)>\,du$$
  pour tout $e\in V_{\sigma,P}$ et tout $\epsilon\in V_{\check{\sigma}}$, o\`u $\check{\sigma}$ est la contragr\'ediente de $\sigma$. Il se prolonge m\'eromorphiquement \`a tout ${\cal A}_{M,{\mathbb C}}^*/i{\cal A}_{M,F}^{\vee}$. Les op\'erateurs d'entrelacement v\'erifient de multiples propri\'et\'es, cf. par exemple [A3] paragraphe 1.

  Soit $x\in G(F)$.  Posons $M'=xMx^{-1}$, $P'=xPx^{-1}$ et supposons $M'$ semi-standard. On note $x\sigma$ la repr\'esentation $m'\mapsto \sigma(x^{-1}m'x)$ de $M'(F)$ dans $V_{\sigma}$. On note $\lambda\mapsto x\lambda$ l'application de ${\cal A}_{M,{\mathbb C}}^*$ dans ${\cal A}_{M',{\mathbb C}}^*$ d\'eduite par fonctorialit\'e de la conjugaison par $x$.  Posons $\pi'_{x\lambda}=Ind_{P'}^G((x\sigma)_{x\lambda})$. On r\'ealise $\pi_{\lambda}$ et $\pi'_{x\lambda}$ dans leurs mod\`eles $V_{\pi_{\lambda}}$ et $V_{\pi'_{x\lambda}}$.  Notons $\partial_{P}(x)$ le jacobien de l'application $ad_{x}:U_{P}\to U_{xPx^{-1}}$ (en fait, on peut \'ecrire $x=m'k'=km$, avec $m'\in M'(F)$, $m\in M(F)$, $k,k'\in K$; on a $\partial_{P}(x)=\delta_{P}(m)=\delta_{P'}(m')$). D\'efinissons l'op\'erateur
  $$\gamma(x):V_{\pi_{\lambda}}\to V_{\pi'_{x\lambda}}$$
  par $(\gamma(x)e)(g)=\partial_{P}(x)^{1/2}e(x^{-1}g)$. Il v\'erifie la relation $\gamma(x)\circ\pi_{\lambda} (g)=\pi'_{x\lambda}(g)\circ\gamma(x)$.  
    
   Supposons $\sigma$ irr\'eductible et unitaire et limitons-nous \`a $\lambda\in i{\cal A}_{M,F}^*$. En identifiant $V_{\pi_{\lambda}}$ \`a $V_{\sigma,P}$ et $V_{\pi'_{x\lambda}}$ \`a $V_{x\sigma,P'}$, ces espaces sont munis de produits hermitiens d\'efinis positifs pour lesquels les repr\'esentations $\pi_{\lambda}$ et $\pi'_{x\lambda}$ sont unitaires. Alors l'application $\gamma(x)$ est  unitaire.
  
   On appliquera ces constructions \`a des \'el\'ements du groupe de Weyl $W^G$. Soit $w\in W^G$. On choisit un repr\'esentant $\dot{w}$ de $w$ dans $K$. Le Levi $M'=w(M)$ est semi-standard.   Par abus de notation, on pose simplement $w\sigma=\dot{w}\sigma$ et $\gamma(w)=\gamma(\dot{w})$.  Notons que  $\gamma(w)$, vu comme  un op\'erateur de $V_{\sigma,P}$ dans $V_{w\sigma,w(P)}$, ne d\'epend pas de $\lambda$. On pourra aussi appliquer cette construction \`a des quotients tels que $W^G/W^M$. On pourra v\'erifier que nos formules ne d\'ependent pas des choix des rel\`evements $\dot{w}$.
    
  \bigskip
  
  \subsection{Op\'erateurs normalis\'es}
  Soient $M$ un Levi semi-standard et $\sigma\in \Pi_{disc}(M(F))$.   Pour $P\in {\cal P}(M)$,  l'op\'erateur $J_{P\vert \bar{P}}(\sigma_{\lambda})\circ J_{\bar{P}\vert P}(\sigma_{\lambda})$ est une homoth\'etie qui ne d\'epend pas de $P$. On d\'efinit $m^G(\sigma_{\lambda})$  comme le produit du degr\'e formel  $d(\sigma)$ et de l'inverse du rapport de cette homoth\'etie.   Remarquons qu'en vertu des relations de compatibilit\'e de nos mesures, $m^G(\sigma_{\lambda})$ ne d\'epend que des  mesures sur $G(F)$ et $A_{M}(F)$.
    
  On introduit des op\'erateurs normalis\'es $R_{P'\vert P}(\sigma_{\lambda})$ et des fonctions de normalisation $r_{P'\vert P}(\sigma_{\lambda})$, pour $P,P'\in {\cal P}(M)$, comme en [A3] paragraphe 2. On a l'\'egalit\'e $J_{P'\vert P}(\sigma_{\lambda})=r_{P'\vert P}(\sigma_{\lambda})R_{P'\vert P}(\sigma_{\lambda})$. Pour $\lambda\in i{\cal A}_{M,F}^*$, $R_{P'\vert P}(\sigma_{\lambda})$ est unitaire. On sait que, pour toute racine r\'eduite $\alpha$ de $A_{M}$, on peut d\'efinir une fonction $r_{\alpha}(\sigma_{\lambda})$ de sorte que
  $$(1) \qquad r_{P'\vert P}(\sigma_{\lambda})=\prod_{\alpha>_{P}0,\alpha<_{P'}0}r_{\alpha}(\sigma_{\lambda}),$$
  la notation signifiant que $\alpha$ parcourt les racines r\'eduites qui sont positives pour $P$ et n\'egatives pour $P'$.
  
   Il est utile de rappeler quelle est la forme des fonctions $r_{\alpha}(\sigma_{\lambda})$ dans le cas o\`u $F={\mathbb R}$. Fixons un sous-tore fondamental $S$ de $M$. Son alg\`ebre de Lie $\mathfrak{s}$  se d\'ecompose de fa\c{c}on unique en une somme directe $\mathfrak{a}_{M}\oplus \mathfrak{s}^{M}$ stable par l'action galoisienne. On pose $\mathfrak{h}_{M}=\mathfrak{a}_{M}({\mathbb R})$, $\mathfrak{h}^M=i\mathfrak{s}^M({\mathbb R})$, on d\'efinit l'alg\`ebre de Lie r\'eelle $\mathfrak{h}= \mathfrak{h}_{M}\oplus \mathfrak{h}^M$ et sa complexifi\'ee $\mathfrak{h}_{{\mathbb C}}$.   Remarquons que ${\cal A}_{M}$ s'identifie \`a $\mathfrak{h}_{M}$.  A la s\'erie discr\`ete  $\sigma$ est associ\'e un caract\`ere infinit\'esimal qui est une orbite dans $i\mathfrak{h}^*$ pour l'action du groupe de Weyl absolu de $S$. Choisissons $\mu_{\sigma}$ dans cette orbite. On sait que la projection $\mu_{\sigma}^M$ de $\mu_{\sigma}$ sur $i\mathfrak{h}^{M,*}$ appartient \`a un r\'eseau fixe de cet espace, plus pr\'ecis\'ement \`a un r\'eseau de $(X^*(S)\otimes_{{\mathbb Z}}{\mathbb Q})\cap i\mathfrak{h}^{M,*}$. Soit $\alpha$ une racine r\'eduite de $A_{M}$. D'apr\`es [A3] paragraphe 3, pour un choix convenable du point $\mu_{\sigma}$, la fonction $r_{\alpha}(\sigma_{\lambda})$ est produit d'une constante uniform\'ement born\'ee (c'est-\`a-dire born\'ee ind\'ependamment de $\sigma$ et $\lambda$), d'un nombre fini de fonctions
 $$\frac{1}{<\mu_{\sigma}+\lambda,\check{\beta}>},$$
 o\`u $\check{\beta}$ est une coracine de $S$ se restreignant \`a $\mathfrak{h}_{M}$ en un multiple positif et non nul de $\check{\alpha}$, et, \'eventuellement, d'une fonction
 $$\Gamma(\frac{<\mu_{\sigma}+\lambda,\check{\beta}>+N}{2})\Gamma(\frac{<\mu_{\sigma}+\lambda,\check{\beta}>+N+1}{2})^{-1}$$
 o\`u $\check{\beta}$ est une coracine de $S$ \'egale \`a un multiple positif et non nul de $\check{\alpha}$ et $N\in \{0,1\}$.

  On utilisera les propri\'et\'es suivantes, pour $\lambda\in i{\cal A}_{M,F}^*$.
  
  (2)  La fonction $\lambda\mapsto m^G(\lambda)$ est r\'eguli\`ere, \`a croissance mod\'er\'ee et \`a valeurs positives ou nulles.
  
 Les premi\`ere et  troisi\`eme propri\'et\'es sont dues \`a Harish-Chandra. La deuxi\`eme (qui ne concerne que le cas o\`u $F$ est archim\'edien) r\'esulte de la formule explicite [A3] proposition 3.1.
  
  (3) On a l'\'egalit\'e $r_{P\vert P'}(\sigma_{\lambda})=\overline{r_{P'\vert P}(\sigma_{\lambda})}$.
  
  Cela r\'esulte des relations d'adjonction v\'erifi\'ees tant par les op\'erateurs d'entrelacement que par les op\'erateurs normalis\'es.
  
  (4) On a l'\'egalit\'e $\vert r_{\bar{P}\vert P}(\sigma_{\lambda})\vert ^2=d(\sigma)m^G(\sigma_{\lambda})^{-1}$.
  
  Cela r\'esulte de (3) et de la d\'efinition de $m^G(\sigma_{\lambda})$.
  
  (5) Pour $P,P'\in {\cal P}(M)$, la fonction
  $$\lambda\mapsto r_{P\vert \bar{P}}(\sigma_{\lambda})r_{P'\vert \bar{P}'}(\sigma_{\lambda})^{-1}$$
  est r\'eguli\`ere de valeur absolue $1$. Si $F$ est archim\'edien, ses d\'eriv\'ees sont \`a croissance mod\'er\'ee.
  
  La premi\`ere assertion  r\'esulte de (4). Pour la seconde, on peut supposer $F={\mathbb R}$ car dans le cas o\`u $F={\mathbb C}$, les fonctions de normalisation sont les m\^emes que pour le groupe sur ${\mathbb R}$ d\'eduit de $G$ par restriction des scalaires.  La d\'ecomposition (1) nous ram\`ene \`a consid\'erer une fonction $r_{\alpha}(\sigma_{\lambda})r_{-\alpha}(\sigma_{\lambda})^{-1}$, ou encore, d'apr\`es (3), $r_{\alpha}(\sigma_{\lambda})\overline{r_{\alpha}(\sigma_{\lambda})}^{-1}$. D'apr\`es ce que l'on a dit ci-dessus, une telle fonction est produit fini de termes
  $$\frac{<\mu^M,\check{\beta}>-<\mu_{M}+\lambda,\check{\beta}>} {<\mu^M,\check{\beta}>+<\mu_{M}+\lambda,\check{\beta}>} $$
  et \'eventuellement d'un terme
  $$\frac{\Gamma(\frac{<\mu_{M}+\lambda,\check{\beta}>}{2})\Gamma(\frac{-<\mu_{M}+\lambda,\check{\beta}>+1}{2})}{\Gamma(\frac{-<\mu_{M}+\lambda,\check{\beta}>}{2})\Gamma(\frac{<\mu_{M}+\lambda,\check{\beta}>+1}{2})}.$$
  C'est un simple exercice de montrer que les d\'eriv\'ees de telles fonctions sont \`a croissance mod\'er\'ee en $\lambda$.

  (6) La fonction $\lambda\mapsto r_{P'\vert P}(\sigma_{\lambda})^{-1}$ est r\'eguli\`ere et \`a croissance mod\'er\'ee.
  
  Par la d\'ecomposition (1), on peut supposer $P'$ et $P$ adjacents. Il y a alors un parabolique $Q=LU_{Q}\in {\cal F}(M)$, qui est minimal parmi les paraboliques contenant strictement $M$, de sorte que $P,P'\subset Q$ et $P'\cap L=\overline{P\cap L}$. Par les propri\'et\'es habituelles, $r_{P'\vert P}(\sigma_{\lambda})=r^L_{\overline{P\cap L}\vert P\cap L}(\sigma_{\lambda})$. Alors la propri\'et\'e cherch\'ee r\'esulte de (2) et (4) appliqu\'ees au groupe $L$.
  
  (7) Pour quatre \'el\'ements $P_{1},P_{2},P_{3},P_{4}\in {\cal P}(M)$, la fonction
  $$\lambda\mapsto m^G(\sigma_{\lambda})r_{P_{1}\vert P_{2}}(\sigma_{\lambda})r_{P_{3}\vert P_{4}}(\sigma_{\lambda})$$
  est r\'eguli\`ere et \`a croissance mod\'er\'ee.
  
  On peut l'\'ecrire comme produit de
  $$r_{\bar{P}_{2}\vert P_{1}}(\sigma_{\lambda})^{-1}r_{\bar{P}_{4}\vert P_{3}}(\sigma_{\lambda})^{-1}$$
  et de
  $$m^G(\sigma_{\lambda})r_{\bar{P}_{2}\vert P_{2}}(\sigma_{\lambda})r_{\bar{P}_{4}\vert P_{4}}(\sigma_{\lambda}).$$
  La premi\`ere fonction v\'erifie les propri\'et\'es requises d'apr\`es (6). La seconde est de valeur absolue constante d'apr\`es (4).

  \bigskip
  
  \subsection{$R$-groupes}
  Soient $M$ un Levi semi-standard et $\sigma$ une repr\'esentation irr\'eductible de $M(F)$ de la s\'erie discr\`ete. On note ${\cal N}^G(\sigma)$ l'ensemble des couples $(A,n)$ o\`u $A$ est un automorphisme unitaire de $V_{\sigma}$ et $n\in Norm_{G(F)}(M)$ (le normalisateur de $M$ dans $G(F)$) qui v\'erifient la condition
  $$A\circ n\sigma(x)=\sigma(x)\circ A$$
  pour tout $x\in M(F)$. C'est un groupe pour le produit $(A,n)(A',n')=(AA',nn')$. Il y a un homomorphisme naturel $M(F)\to {\cal N}^G(\sigma)$ qui, \`a $x\in M(F)$, associe $(\sigma(m),m)$. L'image est un sous-groupe distingu\'e de ${\cal N}^G(\sigma)$. On note ${\cal W}^G(\sigma)$ le quotient. On a une suite exacte
  $$(1) \qquad 1\to {\mathbb U}\to {\cal W}^G(\sigma)\to W^G(\sigma)\to 1,$$
  o\`u ${\mathbb U}$ est le groupe des nombres complexes de module $1$ et
  $$W^G(\sigma)=\{w\in W^G; w(M)=M, w\sigma\simeq \sigma\}/W^M.$$

  Soit $P\in {\cal P}(M)$, posons $\pi=Ind_{P}^G(\sigma)$ que l'on r\'ealise dans l'espace naturel $V_{\pi}$. Pour $(A,n)\in {\cal N}^G(\sigma)$, on d\'efinit l'automorphisme $r_{P}(A,n)$ de $V_{\pi}$ par $r_{P}(A,n)=R_{P\vert nPn^{-1}}(\sigma)\circ\gamma(n)\circ A$. Expliquons que $A$ d\'esigne ici l'op\'erateur d\'eduit de $A$ par fonctorialit\'e, c'est-\`a-dire celui qui, \`a une fonction $e:G(F)\to V_{\sigma}$ associe la fonction $g\mapsto Ae(g)$. On v\'erifie que $r_{P}(A,n)$ est un entrelacement, c'est-\`a-dire que
  $$r_{P}(A,n)\circ \pi(g)=\pi(g)\circ r_{P}(A,n)$$
  pour tout $g\in G(F)$. On v\'erifie que l'application $r_{P}$ ainsi d\'efinie est une repr\'esentation unitaire de ${\cal N}^G(\sigma)$. Si on remplace $P$ par un autre $P'\in {\cal P}(M)$, les repr\'esentations $r_{P}$ et $r_{P'}$ sont \'equivalentes: on a $r_{P}(A,n)\circ R_{P\vert P'}(\sigma)=R_{P\vert P'}(\sigma)\circ r_{P'}(A,n)$. Le sous-groupe des $(A,n)\in {\cal N}^G(\sigma)$ tels que $r_{P}(A,n)$ soit l'identit\'e de $V_{\pi}$ est donc ind\'ependant de $P$. Il contient $M(F)$ (identifi\'e comme ci-dessus \`a un sous-groupe de ${\cal N}^G(\sigma)$). On note $W^G_{0}(\sigma)$ le quotient de ce sous-groupe par $M(F)$. On pose
  $${\cal R}^G(\sigma)={\cal W}^G(\sigma)/W^G_{0}(\sigma).$$
  
  Par la suite (1), $W^G_{0}(\sigma)$ s'identifie \`a un sous-groupe distingu\'e de $W^G(\sigma)$. On pose $R^G(\sigma)=W^G(\sigma)/W^G_{0}(\sigma)$. C'est le $R$-groupe de $\sigma$. On a une suite exacte
  $$1\to {\mathbb U}\to {\cal R}^G(\sigma)\to R^G(\sigma)\to 1.$$
  La repr\'esentation $r_{P}$ se quotiente en une repr\'esentation de ${\cal R}^G(\sigma)$ telle que $z\in {\mathbb U}$ agisse par l'homoth\'etie de rapport $z$.
  
  {\bf Dans la suite, toutes les repr\'esentations de ${\cal R}^G(\sigma)$ seront suppos\'ees v\'erifier cette condition.}
  
  On note $Irr({\cal R}^G(\sigma))$ l'ensemble des classes d'isomorphisme de repr\'esentations irr\'eductibles de ${\cal R}^G(\sigma)$. Le r\'esultat principal de la th\'eorie du $R$-groupe est qu'il existe une bijection $\rho\mapsto \pi_{\rho}$ de $Irr({\cal R}^G(\sigma))$ sur l'ensemble des classes de repr\'esentations irr\'eductibles de $G(F)$ qui sont des sous-repr\'esentations de $\pi$, de sorte que la repr\'esentation $r_{P}\otimes \pi$ de ${\cal R}^G(\sigma)\times G(F)$ dans $V_{\pi}$ se d\'ecompose en
  $$(2) \qquad \oplus_{\rho\in Irr({\cal R}^G(\sigma))}\rho\otimes \pi_{\rho}.$$
  Notons que la correspondance $\rho\mapsto \pi_{\rho}$ ne d\'epend pas de $P$: si on remplace $P$ par $P'$, l'op\'erateur $R_{P\vert P'}(\sigma)$ entrelace les repr\'esentations de ${\cal R}^G(\sigma)\times G(F)$ dans $V_{\pi}$ et dans $V_{\pi'}$, o\`u $\pi'=Ind_{P'}^G(\sigma)$.

   Harish-Chandra a d\'ecrit les groupes $W^G_{0}(\sigma)$ et $R^G(\sigma)$. Soit $\alpha$ une racine r\'eduite de $A_{M}$ dans $G$. Il lui est associ\'e un groupe de Levi $M_{\alpha}\in {\cal L}(M)$ qui est minimal dans ${\cal L}(M)-\{M\}$. Suposons que la mesure de Plancherel $m^{M_{\alpha}}(\sigma)$ soit nulle. On montre qu'alors $W^{M_{\alpha}}(M)$ a deux \'el\'ements: l'identit\'e et une sym\'etrie $s_{\alpha}$. L'ensemble des $\alpha$  v\'erifiant les conditions ci-dessus forme un syst\`eme de racines et $W^G_{0}(\sigma)$ est le groupe de Weyl de ce syst\`eme. En particulier,  $W^G_{0}(\sigma)$ est engendr\'e par ces sym\'etries $s_{\alpha}$. Fixons un ensemble de racines positives dans ce syst\`eme de racines.   On peut identifier $R^G(\sigma)$ au sous-ensemble des $w\in W^G(\sigma)$ qui conservent l'ensemble fix\'e de racines positives. On a alors la d\'ecomposition
   $$W^G(\sigma)=W^{G}_{0}(\sigma)\rtimes R^G(\sigma).$$
   
   {\bf Remarque.} Les d\'efinitions ci-dessus d\'ependent de choix d'op\'erateurs d'entrelacement normalis\'es. Si on change de choix, on obtient une nouvelle repr\'esentation, notons-la $\underline{r}_{P}$. On v\'erifie qu'il existe un caract\`ere unitaire $\chi$ de $W^G(\sigma)$, que l'on remonte  en un caract\`ere de ${\cal N}^G(\sigma)$, de sorte que $\underline{r}_{P}(A,n)=\chi(A,n)r_{P}(A,n)$. On peut introduire l'automorphisme $\boldsymbol{\chi}$ de ${\cal N}^G(\sigma)$ qui \`a $(A,n)$ associe $(\chi(A,n)A,n)$. On a alors $\underline{r}_{P}=r_{P}\circ\boldsymbol{\chi}$. A l'automorphisme $\boldsymbol{\chi}$ pr\`es, nos d\'efinitions sont donc ind\'ependantes des choix d'op\'erateurs normalis\'es. En tout cas, dans la suite, ces choix sont fix\'es.
   
   Soit $\lambda\in i{\cal A}_{G,F}^*$. L'application
   $$\begin{array}{ccc}{\cal N}^G(\sigma)&\to &{\cal N}^G(\sigma_{\lambda})\\ (A,n)&\mapsto& (Ae^{<\lambda,H_{G}(n)>},n)\\ \end{array}$$
   est un isomorphisme. A cause de la torsion par $e^{<\lambda,H_{G}(n)>}$, cet isomorphisme est  compatible avec les plongements de $M(F)$ dans les deux groupes ${\cal N}^G(\sigma)$ et ${\cal N}^G(\sigma_{\lambda})$. Il se quotiente en un isomorphisme ${\cal R}(\sigma)\simeq {\cal R}(\sigma_{\lambda})$ compatible avec l'identit\'e $R(\sigma)= R(\sigma_{\lambda})$.      
  
  \bigskip
  
  \subsection{Coefficients}
  Rappelons qu'Harish-Chandra a d\'efini une fonction $\Xi$ sur $G(F)$. Elle v\'erifie les propri\'et\'es:
 
 (1) la fonction $\Xi$ est biinvariante par $K$.
 
 (2) il existe $d\in {\mathbb N}$ tel que l'on ait les majorations
 $$\delta_{0}(m)^{-1/2}<<\Xi(m)<<\delta_{0}(m)^{-1/2}(1+\vert H_{0}(m)\vert )^d$$
 pour tout  $m\in M_{0}(F)^{\geq}$.
 
  Soient $P=MU_{P}$ un sous-groupe parabolique semi-standard et $\sigma\in \Pi_{disc}(M(F))$. Pour $\lambda\in i{\cal A}_{M,F}^*$, posons $\pi_{\lambda}=Ind_{P}^G(\sigma_{\lambda})$, que l'on r\'ealise dans l'espace $V_{\sigma,P}$. Pour $u,v\in V_{\sigma,P}$, on d\'efinit une fonction coefficient sur $G(F)$ qui, \`a $g\in G(F)$, associe le produit scalaire  $(u,\pi_{\lambda}(g)v)$. On a
  
  (3) $\vert (u,\pi_{\lambda}(g)v)\vert <<\Xi(g)$ pour tout $g\in G(F)$ et tout $\lambda\in i{\cal A}_{M,F}^*$;
  
  plus pr\'ecis\'ement

(4) $\int_{K}\vert (u(k),(\pi_{\lambda}(g)v)(k))\vert \,dk\,<<\Xi(g)$ pour tout $g\in G(F)$ et tout $\lambda\in i{\cal A}_{M,F}^*$.

Cela r\'esulte de la preuve du lemme VI.2.2 de [W].  
  
  Soit $Q=LU_{Q}$ un sous-groupe parabolique standard. Notons $W^G(L\vert P)=\{s\in W^G; s(M)\subset L, P_{0}\cap L\subset s(P)\cap L\}/W^M$. Rappelons que pour tout $s$ dans cet ensemble, on fixe un rel\`evement de $s$ dans $K$, que l'on note encore $s$. On pose $Q_{s}=(s(P)\cap L)U_{Q}$ et $\underline{Q}_{s}=(s(P)\cap L)U_{\bar{Q}}$. On introduit la repr\'esentation $\pi^L_{s\lambda}=Ind_{s(P)\cap L}^L((s\sigma)_{s\lambda})$, que l'on r\'ealise dans l'espace $V^L_{s\sigma,s(P)\cap L}$. Le terme
  $$J_{Q_{s}\vert s(P)}((s\sigma)_{s\lambda})\circ\gamma(s)u$$
 est une fonction sur $K$ appartenant \`a l'espace de la repr\'esentation $Ind_{Q_{s}}^G((s\sigma)_{s\lambda})$. Sa restriction \`a $K\cap L(F)$ appartient \`a $V^L_{s\sigma,s(P)\cap L}$. On note $u_{s}(\lambda)$ cette restriction. De m\^eme, on note $v_{s}(\lambda)$ la restriction de
$$J_{\underline{Q}_{s}\vert s(P)}((s\sigma)_{s\lambda})\circ\gamma(s)v$$ 
  \`a $K\cap L(F)$. Pour $l\in L(F)$, posons
  $$E_{\flat}(l,\lambda)=\gamma(G\vert L)^{-1}\sum_{s\in W^G(L\vert P)}(u_{s}(\lambda),\pi^L_{s\lambda}(l)v_{s}(\lambda)).$$
  On appelle cette fonction le terme constant faible du coefficient $(u,\pi_{\lambda}(g)v)$.
  
  \ass{Proposition}{Soit $\nu>0$ un r\'eel. Il existe $c>0$ et une fonction $C$ sur $i{\cal A}_{M,F}^*$, lisse, \`a croissance mod\'er\'ee et \`a valeurs positives de sorte que l'on ait la majoration
  $$\vert (u,\pi_{\lambda}(m)v)-\delta_{Q}(m)^{-1/2}E_{\flat}(m,\lambda)\vert \leq C(\lambda)\delta_{Q}(m)^{-1/2}\Xi^L(m)e^{-c\vert H_{0}(m)\vert }$$
  pour tout $\lambda\in i{\cal A}_{M,F}^*$ et tout $m\in M_{0}(F)^{\geq}$ v\'erifiant la condition
  $$<\alpha,H_{0}(m)>\geq \nu\vert H_{0}(m)\vert $$
  pour tout $\alpha\in \Delta_{0}-\Delta_{0}^Q$.}
  
  Cf. [A1] lemme 7.1 et [W] lemme VI.2.3 dans le cas non-archim\'edien.

  \subsection{La formule de Plancherel}
   
  Cette formule exprime toute fonction  $K$-finie $f\in C_{c}^{\infty}(G(F))$ \`a l'aide de ses actions dans les repr\'esentations induites unitaires de repr\'esentations de la s\'erie discr\`ete des Levi de $G$ (la condition que $f$ est $K$-finie n'est une restriction que si $F$ est archim\'edien).  Dans la formule ci-dessous, pour tout $M_{disc}\in {\cal L}(M_{0})$, on fixe un parabolique $S\in {\cal P}(M_{disc})$ et, pour tout \'el\'ement de   $\Pi_{disc}(M_{disc}(F))/i{\cal A}_{M_{disc},F}^*$, on fixe un point-base $\sigma$ dans cette orbite.  Alors, pour toute fonction $K$-finie $f\in C_{c}^{\infty}(G(F))$ et tout $g\in G(F)$, on a l'\'egalit\'e
  $$f(g)=\sum_{M_{disc}\in {\cal L}(M_{0})}\vert W^{M_{disc}}\vert \vert W^G\vert ^{-1}\sum_{\sigma\in \Pi_{disc}(M_{disc}(F))/i{\cal A}_{M_{disc},F}^*} \vert Stab(i{\cal A}_{M_{disc},F}^*,\sigma)\vert ^{-1}$$
  $$\int_{i{\cal A}_{M_{disc},F}^*}m^G(\sigma_{\lambda})trace(Ind_{S}^G(\sigma_{\lambda},g^{-1})Ind_{S}^G(\sigma_{\lambda},f))\,d\lambda.$$

\bigskip

\subsection{La formule du produit scalaire d'Arthur}
  Pour exprimer cette formule  sous une forme suffisamment g\'en\'erale, on fixe un sous-tore $D\subset A_{G}$ et un caract\`ere unitaire $\omega$ de $G(F)$ dont la restriction \`a $D(F)$ est triviale.  On adapte nos notations en notant par exemple ${\cal A}_{0}^D$ l'orthogonal du sous-espace ${\cal A}_{D}\subset {\cal A}_{0}$ et $H_{D}$, $H^D$ les projections d'un \'el\'ement $H\in {\cal A}_{0}$ sur ${\cal A}_{D}$ et ${\cal A}^D$.  On consid\`ere un \'el\'ement $T\in {\cal A}_{0}$ que l'on traite comme une variable. Mais on suppose qu'il reste dans un sous-ensemble fixe de ${\cal A}_{0}$ d\'efini par des relations

$<\alpha,T>>0$ pour tout $\alpha\in \Delta_{0}$;

$<\alpha,T>> c_{\star}\vert T\vert $ pour tout $\alpha\in \Delta_{0}$, o\`u $c_{\star}>0$ est un r\'eel fix\'e;

$T\in {\cal A}_{M_{0},F}\otimes_{{\mathbb Z}}{\mathbb Q}$ si $F$ est non-archim\'edien.

Remarquons que, dans un tel domaine, les fonctions $\vert T\vert $ et $<\alpha,T>$, pour $\alpha\in \Delta_{0}$, sont \'equivalentes. On note $\kappa^T$ la fonction caract\'eristique de l'ensemble des $g\in G(F)$ tels que $\phi^G(h_{0}(g)-T)=1$, cf. 1.1 pour la d\'efinition de $h_{0}(g)$.

Soient $P=MU_{P}$ et $P'=M'U_{P'}$ deux sous-groupes paraboliques semi-standard et soient $\sigma\in \Pi_{disc}(M(F))$, $\sigma'\in \Pi_{disc}(M'(F))$. On suppose que:

- les restrictions \`a $D(F)$ des caract\`eres centraux de $\sigma$ et $\sigma'$ co\"{\i}ncident. 

Fixons $u,v\in V_{\sigma,P}$ et $u',v'\in V_{\sigma',P'}$. Pour $\lambda\in i{\cal A}_{M,F}^*$ et $\lambda'\in i{\cal A}_{M',F}^*$ tels que $\lambda_{D}-\lambda'_{D}\in i{\cal A}_{D,F}^{\vee}$, posons $\pi_{\lambda}=Ind_{P}^G(\sigma_{\lambda})$ et $\pi'_{\lambda'}=Ind_{P'}^G(\sigma'_{\lambda'})$.  Pour $X\in {\cal A}_{G,F}$, notons $G(F;X)$ l'ensemble des $g\in G(F)$ tels que $H_{G}(g)= X$.   On munit cet ensemble d'une mesure de sorte que l'on ait l'\'egalit\'e
$$\int_{ G(F)}f(g)\,dg\,=\int_{ {\cal A}_{G,F}}\int_{G(F;X)}f(g)\,dg\,dX$$
pour toute $f\in C_{c}^{\infty}( G(F))$. On note  $D(F)_{c}=Ker(H_{G})\cap D(F)$. On munit ce groupe d'une mesure de Haar, et l'ensemble $D(F)_{c}\backslash G(F;X)$ de la mesure quotient. Remarquons que $D(F)_{c}\backslash G(F;X)=D(F)\backslash D(F)G(F;X)$, ce dernier quotient est donc muni d'une mesure.

Posons
$$\omega^T(X;\lambda,\lambda')=\int_{D(F)\backslash D(F)G(F;X)}(v,\pi_{\lambda}(g)u)(\pi'_{\lambda'}(g)v',u')\omega(g)\kappa^T(g)\,dg.$$
Cette expression ne d\'epend que de la classe $X+{\cal A}_{D,F}$.

  La formule d'Arthur calcule une valeur approch\'ee de $\omega^T(X;\lambda,\lambda')$. Rappelons ce r\'esultat.
 
 Soit $S=M_{S}U_{S}$ un sous-groupe parabolique standard. Notons $W^G(M_{S}\vert M)=\{s\in W^G; s(M)=M_{S}\}/W^{M}$. On d\'efinit de m\^eme $W^G(M_{S}\vert M')$. Soient $s\in W^G(M_{S}\vert M)$, $s'\in W^G(M_{S}\vert M')$, $k,k'\in K$.  D\'efinissons les \'el\'ements $u_{s}(k',\lambda), v_{s}(k,\lambda)\in V_{\sigma}$ et $u'_{s'}(k,\lambda'),v'_{s'}(k',\lambda')\in V_{\sigma'}$ par les formules suivantes
 $$u_{s}(k',\lambda)=(J_{\bar{S}\vert s(P)}((s\sigma)_{s\lambda})\circ\gamma(s)u)(k'),$$
 $$v_{s}(k,\lambda)=(J_{S\vert s(P)}((s\sigma)_{s\lambda})\circ\gamma(s)v)(k),$$
 $$u'_{s'}(k,\lambda')=(J_{S\vert s'(P')}((s'\sigma')_{s'\lambda'})\circ\gamma(s')u')(k),$$
 $$v'_{s'}(k',\lambda')=(J_{\bar{S}\vert s'(P')}((s'\sigma')_{s'\lambda'})\circ\gamma(s')v')(k').$$
 On rappelle que les termes entre parenth\`eses sont des fonctions sur $K$; cela a donc un sens de les \'evaluer en $k$ ou $k'$. Pour des \'el\'ements $\epsilon,\eta\in V_{\sigma}$ et $\epsilon',\eta'\in V_{\sigma'}$, et pour $\Lambda\in {\cal A}_{M_{S},{\mathbb C}}^*$ tel que $\Lambda_{D}\in i{\cal A}_{D,F}^{\vee}$, posons
 $$r^T_{S,s,s'}(X;\epsilon,\eta,\epsilon',\eta';\Lambda)= \int_{D(F)\backslash D(F) M_{S}^G(F;X)}(\eta,s\sigma(x)\epsilon)(s'\sigma'(x)\eta',\epsilon')$$
 $$\phi_{S}^G(H_{0}(x)-T)e^{<\Lambda,H_{0}(x)>}\omega(x)\,dx,$$
 o\`u $M_{S}^G(F;X)=M_{S}(F)\cap G(F;X)$ et la mesure est d\'efinie de fa\c{c}on analogue \`a celle sur $D(F)\backslash D(F)G(F;X)$. Cette int\'egrale est absolument convergente si    $Re<\Lambda,\check{\alpha}>>0$ pour tout $\alpha\in \Delta_{S}$ et  se prolonge en une fonction m\'eromorphe d\'efinie pour tout $\Lambda\in{\cal A}_{M_{S},{\mathbb C}}^*$ tel que $\Lambda_{D}\in i{\cal A}_{D,F}^{\vee}$.  Posons
 $$r^T_{S,s,s'}(X;\lambda,\lambda';\Lambda)=\int_{K\times K}\omega(k'k^{-1})$$
 $$r^T_{S,s,s'}(X,u_{s}(k',\lambda),v_{s}(k,\lambda),u'_{s'}(k,\lambda'),v'_{s'}(k',\lambda');\Lambda)\,dk\,dk'.$$
Pour des points $\lambda$ et $\lambda'$  o\`u les op\'erateurs d'entrelacement intervenant ci-dessus n'ont pas de p\^oles, cette expression est holomorphe au point $\Lambda=s\lambda-s'\lambda'$.  On pose 
 
 $$r^T(X;\lambda,\lambda')=\sum_{S; P_{0}\subset S}\sum_{s\in W^G(M_{S}\vert M),s'\in W^G(M_{S}\vert M')}r^T_{S,s,s'}(X;\lambda,\lambda'; s\lambda-s'\lambda').$$

On obtient ainsi une fonction m\'eromorphe en $\lambda$ et $\lambda'$. 

\ass{Th\'eor\`eme (Arthur)}{Il existe $c>0$ et une fonction lisse et \`a croissance mod\'er\'ee $C$ sur $i{\cal A}_{M,F}^*\times i{\cal A}_{M',F}^*$, \`a valeurs positives, de sorte que, pour tous $\lambda\in i{\cal A}_{M,F}^*$, $\lambda'\in i{\cal A}_{M',F}^*$ pour tous $X$, $T$, $u,v,u',v'$, on ait la majoration
$$\vert \omega^T(X;\lambda,\lambda')-r^T(X;\lambda,\lambda')\vert \leq m^G(\sigma_{\lambda})^{-1/2}m^G(\sigma_{\lambda'})^{-1/2}C(\lambda,\lambda')e^{-c\vert T\vert }\vert u\vert \vert v\vert \vert u'\vert \vert v'\vert .$$}

On peut exprimer   plus explicitement le terme $r^T_{S,s,t}(X;\lambda,\lambda';\Lambda)$.  Commen\c{c}ons par calculer $r^T_{S,s,s'}(X;\epsilon,\eta,\epsilon',\eta';\Lambda)$.  Pour simplifier les notations, on fixe $\epsilon,\eta,\epsilon',\eta'$. Pour $Y\in {\cal A}_{M_{S},F}$, posons
$$r_{S,s,s'}(Y )=\int_{D(F)\backslash D(F)M_{S}(F;Y)}(\eta,s\sigma(x)\epsilon)(s'\sigma'(x)\eta',\epsilon')\omega(x)\,dx,$$ 
o\`u $M_{S}(F;Y)=\{x\in M_{S}(F); H_{M_{S}}(x)= Y\}$. Posons comme en 1.6
$${\cal A}_{M_{S},F}^{G}(X)=\{Y\in {\cal A}_{M_{S},F}; Y_{G}= X\} .$$ 
  On a les \'egalit\'es
 $$r^T_{S,s,s'}(X;\epsilon,\eta,\epsilon',\eta';\Lambda)=\int_{({\cal A}_{M_{S},F}^{G}(X)+{\cal A}_{D,F})/{\cal A}_{D,F}}\phi_{S}^G(Y-T)e^{<\Lambda,Y>}r_{S,s,s'}(Y)\,dY$$
 $$=\int_{{\cal A}_{M_{S},F}^{G}(X)}\phi_{S}^G(Y-T)e^{<\Lambda,Y>}r_{S,s,s'}(Y)\,dY.$$
  Notons $\omega_{\sigma}$ et $\omega_{\sigma'}$ les caract\`eres centraux de $\sigma$ et $\sigma'$.  Il est clair que, pour $Y\in {\cal A}_{M_{S},F}$, $r_{S,s,s'}(Y)=0$ si la restriction de $s(\omega_{\sigma})s'(\omega_{\sigma'})^{-1}\omega$ \`a $A_{M_{S}}(F)_{c}$ est non triviale. Supposons cette condition v\'erifi\'ee.  Il existe alors un \'el\'ement $\Lambda_{s,s'}\in i{\cal A}_{M_{S}}^*$ tel que
 $$(s(\omega_{\sigma})s'(\omega_{\sigma'})^{-1}\omega)(a)=e^{<\Lambda_{s,s'},H_{M_{S}}(a)>}
$$
pour tout $a\in A_{M_{S}}(F)$. Cet \'el\'ement est uniquement d\'etermin\'e modulo $i{\cal A}_{A_{M_{S}},F}^{\vee}$ et v\'erifie $(\Lambda_{s,s'})_{D}\in i{\cal A}_{D,F}^{\vee}$. La fonction 
$$Y\mapsto e^{-<\Lambda_{s,s'},Y>}r_{S,s,s'}( Y)$$
ne d\'epend que de la classe $Y+{\cal A}_{A_{M_{S}},F}$.  L'ensemble $A_{M_{S}}(F)M_{S}(F;Y)$ est ouvert dans $M_{S}(F)$ et donc muni d'une mesure. On v\'erifie que, pour une fonction $f$ sur le quotient $A_{M_{S}}(F)\backslash A_{M_{S}}(F)M_{S}(F;Y)$, on a la formule d'int\'egration
$$\int_{D(F)\backslash D(F)M_{S}(F;Y)}f(x)\,dx\,=C\int_{A_{M_{S}}(F)\backslash A_{M_{S}}(F)M_{S}(F;Y)}f(x)\,dx,$$
o\`u
$$C=mes(D(F)_{c})^{-1}mes(i{\cal A}_{M_{S},F}^*)^{-1}[{\cal A}_{M_{S},F}:{\cal A}_{A_{M_{S}},F}].$$
De ces consid\'erations r\'esulte  l'\'egalit\'e

$e^{-<\Lambda_{s,s'},Y>}r_{S,s,s'}( Y)=$
$$C\int_{A_{M_{S}}(F)\backslash A_{M_{S}}(F)M_{S}(F;Y)}(\eta,s\sigma(x)\epsilon)(s'\sigma'(x)\eta',\epsilon')\omega(x)e^{-<\Lambda_{s,s'},H_{M_{S}}(x)>}\,dx.$$

Pour $\mu\in i{\cal A}_{A_{M_{S}},F}^{\vee}/i{\cal A}_{M_{S},F}^{\vee}$, posons
$$ r^T_{S,s,s'}(\mu)=\sum_{Y\in{\cal A}_{M_{S},F}/{\cal A}_{A_{M_{S}},F}}e^{-<\Lambda_{s,s'}+\mu,Y>}r^T_{S,s,s'}(Y).$$
Par inversion de Fourier, on a l'\'egalit\'e
$$(1) \qquad  r_{S,s,s'}( Y)= [{\cal A}_{M_{S},F}:{\cal A}_{A_{M_{S}},F}]^{-1}\sum_{\mu\in  i{\cal A}_{A_{M_{S}},F}^{\vee}/i{\cal A}_{M_{S},F}^{\vee}}e^{<\Lambda_{s,s'}+\mu,Y>}r^T_{S,s,s'}(\mu ).$$
L'\'egalit\'e $s(M)=M_{S}$ entra\^{\i}ne $mes(i{\cal A}_{M_{S},F}^*)=mes(i{\cal A}_{M,F}^*)$. D'autre part, on a
$$r^T_{S,s,s'}(\mu )=C \int_{A_{M_{S}}(F)\backslash M_{S}(F)}(\eta,s\sigma(x)\epsilon)(s'\sigma'(x)\eta',\epsilon')e^{-<\Lambda_{s,s'}+\mu,H_{M_{S}}(x)>}\omega(x)\,dx.$$
Cette derni\`ere expression est le produit scalaire de deux coefficients de repr\'esentations de la s\'erie discr\`ete. Il est nul si la repr\'esentation $s'\sigma'$ n'est pas \'equivalente \`a la repr\'esentation $(\omega s\sigma)_{-\Lambda_{s,s'}-\mu}$, autrement dit si $\Lambda_{s,s'}+\mu\not\in [s'\sigma',\omega s\sigma]$. Soit $\nu\in [s'\sigma',\omega s\sigma]$.  Introduisons un isomorphisme unitaire $A_{\nu}:V_{\sigma}\to V_{\sigma'}$ tel que $ (s'\sigma')(x)\circ A_{\nu}=\omega(x)A_{\nu}\circ(s\sigma)_{-\nu}(x) $ pour tout $x\in M'(F)$. Alors
$$r^T_{S,s,s'}(\nu-\Lambda_{s,s'} )= Cd(\sigma)^{-1}(A_{\nu}\eta,\epsilon')(\eta',A_{\nu}\epsilon).$$
 La formule (1) se r\'ecrit
 $$r_{S,s,s'}(Y)= mes(D(F)_{c})^{-1}mes(i{\cal A}_{M,F}^*)^{-1}d(\sigma)^{-1}\sum_{\nu\in [s'\sigma',\omega s\sigma]}e^{<\nu,Y>}(A_{\nu}\eta,\epsilon')(\eta',A_{\nu}\epsilon).$$
 On a d\'efini $\epsilon_{S}^{G,T}(X;\Lambda)$ en 1.6. On obtient finalement
 $$ r^T_{S,s,s'}(X;\epsilon,\eta,\epsilon',\eta';\Lambda)=mes(D(F)_{c})^{-1}mes(i{\cal A}_{M,F}^*)^{-1} d(\sigma)^{-1}$$
 $$\sum_{\nu\in [s'\sigma',\omega s\sigma]}(A_{\nu}\eta,\epsilon')(\eta',A_{\nu}\epsilon)\epsilon_{S}^{G,T}(X;\Lambda+\nu).$$
 
 Pour calculer  $r^T_{S,s,s'}(X;\lambda,\lambda';\Lambda)$, on doit remplacer $\epsilon,\eta,\epsilon',\eta'$ par $u_{s}(k',\lambda)$, $v_{s}(k,\lambda)$, $u'_{s'}(k,\lambda')$, $v'_{s'}(k',\lambda')$, multiplier par $\omega(k'k^{-1})$ puis int\'egrer en $k,k'\in K$. Il appara\^{\i}t des int\'egrales
 $$\int_{K}(v'_{s'}(k',\lambda'),A_{\nu}u_{s}(k',\lambda) )\omega(k')\,dk',$$
 $$\int_{K}(A_{\nu}v_{s}(k,\lambda),u'_{s'}(k,\lambda'))\omega(k)^{-1}\,dk.$$
 Consid\'erons la premi\`ere, la seconde \'etant analogue. L'op\'erateur $A_{\nu}$ d\'efinit par fonctorialit\'e un op\'erateur de $V_{ s\sigma,\bar{S}}$ dans $V_{\omega^{-1} s'\sigma',\bar{S}}$, que nous notons encore $A_{\nu}$. Notons $\underline{\omega}$ l'op\'erateur qui, \`a une fonction $f$ sur $K$, associe la fonction $k\mapsto \omega(k)f(k)$. Alors $\underline{\omega}\circ A_{\nu}$ envoie $V_{s\sigma,\bar{S}}$ dans $V_{s'\sigma',\bar{S}}$. L'int\'egrale ci-dessus se r\'ecrit
 $$\int_{K}((J_{\bar{S}\vert s'(P')}((s'\sigma')_{s'\lambda'})\circ\gamma(s')v')(k'), (\underline{\omega}\circ A_{\nu}\circ J_{\bar{S}\vert s(P)}((s\sigma)_{s\lambda})\circ \gamma(s)u)(k'))\,dk'.$$
 Par d\'efinition du produit scalaire dans $V_{s'\sigma',\bar{S}}$, ce n'est autre que
 $$(J_{\bar{S}\vert s'(P')}((s'\sigma')_{s'\lambda'})\circ\gamma(s')v', \underline{\omega}\circ A_{\nu}\circ J_{\bar{S}\vert s(P)}((s\sigma)_{s\lambda})\circ \gamma(s)u).$$
 On obtient alors
 $$(2) \qquad r^T_{S,s,s'}(X;\lambda,\lambda';\Lambda)=mes(D(F)_{c})^{-1}mes(i{\cal A}_{M,F}^*)^{-1}d(\sigma)^{-1}$$
 $$\sum_{\nu\in [s'\sigma',\omega s\sigma]}(J_{\bar{S}\vert s'(P')}((s'\sigma')_{s'\lambda'})\circ\gamma(s')v', \underline{\omega}\circ A_{\nu}\circ J_{\bar{S}\vert s(P)}((s\sigma)_{s\lambda})\circ \gamma(s)u)$$
 $$(\underline{\omega}\circ A_{\nu}\circ J_{S\vert s(P)}((s\sigma)_{s\lambda})\circ\gamma(s)v,J_{S\vert s'(P')}((s'\sigma')_{s'\lambda'})\circ\gamma(s')u')\epsilon_{S}^{G,T}(X;\Lambda+\nu).$$

 \bigskip
 
 \section{Espaces tordus}
 
 \subsection{Notations}
 
 Soit $\tilde{G}$ un espace tordu sous $G$. Rappelons qu'un tel espace est une vari\'et\'e alg\'ebrique sur $F$ munie d'actions alg\'ebriques de $G$ \`a droite et \`a gauche
 $$\begin{array}{ccc}G\times \tilde{G}\times G&\to&\tilde{G}\\ (g,\gamma,g')&\mapsto &g\gamma g'\\ \end{array}$$
 pour chacune desquelles $\tilde{G}$ est un espace principal homog\`ene. On impose la condition 
 
(1)  $\tilde{G}(F)\not=\emptyset$. 
 
 Pour $\gamma\in \tilde{G}$, on note $ad_{\gamma}$ l'automorphisme de $G$ tel que $\gamma g=ad_{\gamma}(g)\gamma$ pour tout $g\in G$. Il arrive que $ad_{\gamma}$ induise sur certains objets des automorphismes ind\'ependants de $\gamma$. On note alors $\theta$ ces automorphismes. Par exemple, $\tilde{G}$ d\'etermine un automorphisme $\theta$ de $Z_{G}$ ou de $A_{G}$.  On impose la condition
 
(2)  l'automorphisme $\theta$ de $Z_{G}$ est d'ordre fini.

 {\bf Remarque.} Les hypoth\`eses (1) et  (2) impliquent qu'il existe un groupe alg\'ebrique non connexe $G^+$ d\'efini sur $F$, de composante neutre $G$, tel que $\tilde{G}$ s'identifie \`a une composante connexe de $G^+$, cf. [L] paragraphe 3.4.
 Cela nous permet d'utiliser pour notre espace $\tilde{G}$ les r\'esultats d\'emontr\'es dans la litt\'erature pour les groupes non connexes. Il est n\'eanmoins tr\`es probable que ces r\'esultats restent vrais sous une hypoth\`ese plus faible que (2).

 On note $A_{\tilde{G}}$ le sous-tore de $A_{G}$ tel que $X_{*}(A_{\tilde{G}})=X_{*}(A_{G})^{\theta}$, l'exposant signifiant selon l'usage l'ensemble des invariants par $\theta$. On pose ${\cal A}_{\tilde{G}}=X_{*}(A_{\tilde{G}})\otimes_{{\mathbb Z}}{\mathbb R}={\cal A}_{G}^{\theta}$. On note $a_{\tilde{G}}$ la dimension de ${\cal A}_{\tilde{G}}$. On note $H_{\tilde{G}}:G(F)\to {\cal A}_{\tilde{G}}$ la compos\'ee de $H_{G}$ et de la projection orthogonale de ${\cal A}_{G}$ sur ${\cal A}_{\tilde{G}}$. On note ${\cal A}_{\tilde{G},F}$ l'image de $G(F)$ par l'application $H_{\tilde{G}}$.
 
 Soit $(P,M)$ une paire parabolique de $G$. Notons $\tilde{P}$ le normalisateur de $P$ dans $\tilde{G}$, c'est-\`a-dire l'ensemble des $\gamma\in \tilde{G}$ tels que $ad_{\gamma}$ conserve $P$. Notons $\tilde{M}$ le normalisateur commun de $P$ et $M$ dans $\tilde{G}$. Si $\tilde{P}$ n'est pas vide, les ensembles $\tilde{M}$, $\tilde{P}(F)$, $\tilde{M}(F)$ ne le sont pas non plus. On dit alors que $\tilde{P}$ est un espace parabolique et $\tilde{M}$ est un espace de Levi. Dans le cas de la paire minimale $(P_{0},M_{0})$, il est toujours vrai que $\tilde{P}_{0}$ n'est pas vide. Soit $\gamma_{0}\in \tilde{M}_{0}(F)$. Un sous-groupe parabolique standard $P$ donne naissance \`a un espace parabolique $\tilde{P}$ si et seulement si $ad_{\gamma_{0}}(P)=P$. On peut supposer et on suppose que la forme quadratique fix\'ee sur ${\cal A}_{0}$ est invariante par l'automorphisme  d\'eduit de $ad_{\gamma_{0}}$ par fonctorialit\'e.

 On adopte pour les espaces de Levi et les espaces paraboliques des notations similaires \`a celles introduites pour les groupes. Par exemple, on note ${\cal P}(\tilde{M})$ l'ensemble des espaces paraboliques de composante de Levi $\tilde{M}$. Ou bien, pour $H\in {\cal A}_{M}$, on note $H_{\tilde{M}}$ et $H^{\tilde{M}}$ ses projections sur l'espace ${\cal A}_{\tilde{M}}$, resp. sur son orthogonal ${\cal A}_{M}^{\tilde{M}}$ dans ${\cal A}_{M}$.  On munit les espaces ${\cal A}_{\tilde{M}}$, $i{\cal A}_{\tilde{M}}^*$ et le groupe $A_{\tilde{M}}(F)$ de mesures de Haar v\'erifiant des conditions similaires \`a celles de 1.2. 
 
 Soient $L$ un Levi de $G$ et $\gamma\in \tilde{G}(F)$ tel que $ad_{\gamma}(L)=L$. Notons $\theta$ l'automorphisme de ${\cal A}_{L}$ d\'eduit fonctoriellement de $ad_{\gamma}$. Posons ${\cal A}_{L}^{\theta}=\{H\in {\cal A}_{L}; \theta H=H\}$. Notons $M$ le plus grand Levi contenant $L$ tel que ${\cal A}_{L}^{\theta}\subset {\cal A}_{M}$. Posons $\tilde{M}=M\gamma$. On a
 
 (3) l'ensemble $\tilde{M}$ est un ensemble de Levi;
 
 (4) on a l'\'egalit\'e ${\cal A}_{L}^{\theta}={\cal A}_{\tilde{M}}$;
 
 (5) l'ensemble des espaces paraboliques $\tilde{P}$ tels que $L\subset P$ et $ad_{\gamma}(P)=P$ est \'egal \`a ${\cal F}(\tilde{M})$.
 
 Preuve. Notons $T$ le sous-tore de $A_{L}$ tel que $X_{*}(T)=X_{*}(A_{L})\cap {\cal A}_{L}^{\theta}$. Le Levi $M$ est le commutant de ce tore dans $G$. Par d\'efinition de $T$ et $\theta$, $ad_{\gamma}$ fixe tout point de $T$. Donc $\tilde{M}$ est le commutant de $T$ dans $\tilde{G}$ (c'est-\`a-dire l'ensemble des $\gamma'\in \tilde{G}$ tels que $ad_{\gamma'}$ fixe tout point de $T$). On sait qu'un tel commutant est un espace de Levi, pourvu qu'il ne soit pas vide. Cette derni\`ere condition est v\'erifi\'ee. D'o\`u (1). Evidemment $T$ est inclus dans $A_{\tilde{M}}$, donc ${\cal A}_{L}^{\theta}\subset {\cal A}_{\tilde{M}}$. Inversement, un \'el\'ement $H\in {\cal A}_{\tilde{M}}$ appartient \`a ${\cal A}_{M}$ donc aussi \`a ${\cal A}_{L}$. Puisque $\gamma\in \tilde{M}$, on a aussi $\theta H=H$, donc $H\in {\cal A}_{L}^{\theta}$. Cela prouve (2). Soit $\tilde{P}=\tilde{M}_{P}U_{P}$ un espace parabolique tel que $L\subset M_{P}$ et $ad_{\gamma}(P)=P$. Puisque $ad_{\gamma}$ conserve $L$ et que $M_{P}$ est l'unique composante de Levi de $P$ contenant $L$, on a aussi $ad_{\gamma}(M_{P})=M_{P}$. Donc $\tilde{M}_{P}=M_{P}\gamma$. Alors ${\cal A}_{\tilde{M}_{P}}$ est un sous-espace de ${\cal A}_{L}$ sur lequel $\theta$ agit trivialement. D'o\`u ${\cal A}_{\tilde{M}_{P}}\subset {\cal A}_{L}^{\theta}={\cal A}_{\tilde{M}}$. D'o\`u $\tilde{M}\subset \tilde{M}_{P}$, ce qui \'equivaut \`a $\tilde{P}\in {\cal F}(\tilde{M})$. La r\'eciproque est claire. $\square$
 
  On a fix\'e une paire parabolique minimale $(P_{0},M_{0})$ et un groupe compact maximal $K$. On a
 
 (6) si $F$ est archim\'edien, il existe $\gamma_{0}\in \tilde{M}_{0}(F)$ tel que $ad_{\gamma_{0}}$ conserve la composante neutre de $K$.
 
 Cf. [DM] lemme 1.

 \bigskip
 
 \subsection{D\'efinitions combinatoires}
 Soit $\tilde{P}=\tilde{M}U$ un espace parabolique tel que $M_{0}\subset M$. Pour $\alpha\in \Delta_{P}$, on note $\tilde{\alpha}$ la restriction de $\alpha$ \`a $A_{\tilde{M}}$.   On note $\Delta_{\tilde{P}}$ l'ensemble de ces $\tilde{\alpha}$ pour $\alpha\in \Delta_{P}$. Si on fixe $\gamma\in \tilde{M}$, l'automorphisme $ad_{\gamma}$ d\'efinit par fonctorialit\'e une permutation de $\Delta_{P}$ ind\'ependante de $\gamma$, que l'on peut noter $\theta$. Alors $\Delta_{\tilde{P}}$ s'identifie \`a l'ensemble des orbites dans $\Delta_{P}$ pour l'action du groupe de permutations engendr\'e par $\theta$. Pour $\alpha\in \Delta_{P}$, on note $\varpi_{\tilde{\alpha}}$ la restriction de $\varpi_{\alpha}$ \`a $A_{\tilde{M}}$. Cette restriction ne d\'epend en effet que de la restriction $\tilde{\alpha}$ de $\alpha$. Les ensembles $\Delta_{\tilde{P}}$ et $\{\varpi_{\tilde{\alpha}};\tilde{\alpha}\in \Delta_{\tilde{P}}\}$ sont des bases de ${\cal A}_{\tilde{M}}^{\tilde{G},*}$.
 
 Soient $\tilde{P}=\tilde{M}U_{P}\subset \tilde{Q}=\tilde{L}U_{Q}$ deux espaces paraboliques tels que $M_{0}\subset M$. On d\'efinit naturellement le sous-ensemble $\Delta_{\tilde{P}}^{\tilde{Q}}\subset \Delta_{\tilde{P}}$. On d\'efinit les fonctions $\tau_{\tilde{P}}^{\tilde{Q}}$, $\hat{\tau}_{\tilde{P}}^{\tilde{Q}}$, $\phi_{\tilde{P}}^{\tilde{Q}}$, $\delta_{\tilde{P}}^{\tilde{Q}}$ sur ${\cal A}_{\tilde{M}}$ et $\Gamma_{\tilde{P}}^{\tilde{Q}}$ sur ${\cal A}_{\tilde{M}}\times {\cal A}_{\tilde{M}}$ en rempla\c{c}ant $\alpha$ par $\tilde{\alpha}$ et $\varpi_{\alpha}$ par $\varpi_{\tilde{\alpha}}$ dans les d\'efinitions de 1.3. On d\'efinit aussi la notion de famille de points $(\tilde{G},\tilde{M})$-orthogonale. Pour une telle famille ${\cal Y}$, on d\'efinit la fonction $\Gamma_{\tilde{M}}^{\tilde{Q}}(.,{\cal Y})$. Toutes les relations \'enonc\'ees dans le paragraphe 1.3 restent valables pour ces nouvelles fonctions. Evidemment, toutes ces fonctions peuvent \^etre consid\'er\'ees comme des fonctions sur ${\cal A}_{M}$ ou m\^eme sur ${\cal A}_{0}$, par composition avec la projection orthogonale sur ${\cal A}_{\tilde{M}}$.
 
 Soient maintenant $Q=LU_{Q}\subset R$ deux sous-groupes paraboliques tels que $M_{0}\subset L$ et soit  $\tilde{P}=\tilde{M}U_{P}$ un espace parabolique tel que $Q\subset P\subset R$. Notons $\tilde{\sigma}_{Q}^R$ la fonction caract\'eristique du sous-ensemble des $H\in {\cal A}_{L}$ qui v\'erifient les conditions suivantes:
 
 $ <\alpha,H>>0$ pour tout $\alpha\in \Delta_{Q}^R$;
 
 $<\alpha,H>\leq0$ pour tout $\alpha\in \Delta_{Q}-\Delta_{Q}^R$;
 
 $<\varpi_{\tilde{\alpha}},H_{\tilde{M}}>>0$ pour tout $\tilde{\alpha}\in \Delta_{\tilde{P}}$.
 
 On d\'emontre que cette fonction ne d\'epend pas de $\tilde{P}$ v\'erifiant les conditions ci-dessus, cf. [LW] lemme 2.10.3.  Soulignons que cette fonction n'est d\'efinie que s'il existe au moins un $\tilde{P}$ v\'erifiant ces conditions. Pour un sous-groupe parabolique $Q=LU_{Q}$ et un espace parabolique $\tilde{P}$ tel que $Q\subset P$, on a l'\'egalit\'e
 
 (1) $\tau_{Q}^P(H)\hat{\tau}_{\tilde{P}}(H)=\sum_{R; P\subset R}\tilde{\sigma}_{Q}^R(H)$ pour tout $H\in {\cal A}_{L}$,
 
 cf. [LW] lemme 2.10.5.

 \bigskip
 
 \subsection{$(\tilde{G},\tilde{M})$-familles}
 Les d\'efinitions et r\'esultats des paragraphes 1.4 \`a 1.8 s'adaptent aux espaces tordus. Pour un espace de Levi $\tilde{M}$, on d\'efinit une notion de $(\tilde{G},\tilde{M})$-famille: c'est une famille $(\varphi(\Lambda,\tilde{P}))_{\tilde{P}\in {\cal P}(\tilde{M})}$ o\`u $\Lambda\mapsto\varphi(\Lambda,\tilde{P})$ est une fonction $C^{\infty}$ sur $i{\cal A}_{\tilde{M}}^*$; pour deux espaces paraboliques adjacents $\tilde{P}$ et $\tilde{P}'$, les fonctions $\varphi(\Lambda,\tilde{P})$ et $\varphi(\Lambda,\tilde{P}')$ se recollent sur le mur s\'eparant les chambres positives associ\'ees aux deux espaces paraboliques. Pour une telle famille, on d\'efinit la fonction $\varphi_{\tilde{M}}^{\tilde{G}}(\Lambda)$.  Pour $X\in {\cal A}_{\tilde{G},F}$ et pour une famille $(\tilde{G},\tilde{M})$-orthogonale  ${\cal Y}$, on d\'efinit la fonction $\varphi_{\tilde{M}}^{\tilde{G},{\cal Y}}(X;\Lambda)$. Ces fonctions v\'erifient les m\^emes propri\'et\'es que dans le cas non tordu.
 
 Consid\'erons un Levi $L$ contenu dans $\tilde{M}$. Soit $(\varphi(\Lambda,Q))_{Q\in {\cal P}(L)}$ une $(G,L)$-famille.  Pour $\tilde{P}\in {\cal P}(\tilde{M})$, choisissons $Q\in {\cal P}(L)$ tel que $Q\subset P$. La restriction de $\Lambda\mapsto \varphi(\Lambda,Q)$ \`a $i{\cal A}_{\tilde{M}}^*$ ne d\'epend pas du choix de $Q$. On la note $\varphi(\Lambda,\tilde{P})$. Alors la famille $(\varphi(\Lambda,\tilde{P}))_{\tilde{P}\in {\cal P}(\tilde{M})}$ est une $(\tilde{G},\tilde{M})$-famille. 
 
 \bigskip

\subsection{Une construction auxiliaire}
Pour un caract\`ere $\mu$ de $G(F)$ et pour $\gamma\in \tilde{G}(F)$, le caract\`ere $\mu\circ ad_{\gamma}$ ne d\'epend pas de $\gamma$, on le note $\mu\circ \theta$. On note $\mu\circ(1-\theta)$ le caract\`ere $\mu(\mu\circ\theta)^{-1}$. Dans la suite, on consid\'erera un caract\`ere unitaire $\omega$ de $G(F)$ dont la restriction \`a $Z_{G}(F)^{\theta}$ est triviale. Certains aspects de la th\'eorie se simplifient quand ce caract\`ere   est de la forme $\mu\circ(1-\theta)$. Ce n'est pas toujours le cas mais nous allons voir que l'on peut se ramener \`a ce cas en effectuant des extensions.

Soient $(G',\tilde{G}')$ un couple v\'erifiant les m\^emes hypoth\`eses que notre couple $(G,\tilde{G})$. Soit $p:G'\to G$ un homomorphisme de groupes alg\'ebriques et $\tilde{p}:\tilde{G}'\to \tilde{G}$ un homomorphisme de vari\'et\'es alg\'ebriques (tous deux d\'efinis sur $F$). On dit que $(p,\tilde{p})$ est un homomorphisme d'espaces tordus si et seulement si on a l'\'egalit\'e
$$\tilde{p}(x'\gamma'y')=p(x')\tilde{p}(\gamma')p(y')$$
pour tous $\gamma'\in \tilde{G}'$ et $x',y'\in G'$.

\ass{Proposition}{Soit $\omega$ un caract\`ere unitaire de $G(F)$ dont la restriction \`a $Z_{G}(F)^{\theta}$ est triviale. Il existe un espace tordu $(G',\tilde{G}')$ v\'erifiant les hypoth\`eses (1) et (2) de 2.1, et un homomorphisme d'espaces  tordus $(p,\tilde{p}):(G',\tilde{G}')\to (G,\tilde{G})$ v\'erifiant les conditions suivantes:

(i) l'homomorphisme $p$ s'ins\`ere dans une suite exacte
$$1\to C\to G'\stackrel{p}{\to}G\to 1$$
o\`u $C$ est un sous-tore central de $G'$;

(ii) pour toute extension $F'$ de $F$, l'homomorphisme $p:G'(F')\to G(F')$ est surjectif;

(iii)  il existe un caract\`ere unitaire $\mu'$ de $G'(F)$ tel que $\omega\circ p=\mu'\circ(1-\theta')$ (en notant $\theta'$ l'analogue de $\theta$ pour $(G',\tilde{G}')$).}

Preuve. On commence par montrer 

(1) il existe $(G',\tilde{G}')$ v\'erifiant (i) et (ii) et tel que le groupe d\'eriv\'e de $G'$ soit simplement connexe.

C'est une simple adaptation de la th\'eorie des $z$-extensions. Fixons $\gamma\in \tilde{G}(F)$, posons $\theta=ad_{\gamma}$. Notons $G_{SC}$ le rev\^etement simplement connexe du groupe d\'eriv\'e de $G$. On sait que l'on a un isomorphisme
$$(G_{SC}\times Z_{G}^0)/\Xi\simeq G,$$
o\`u $\Xi$ est un sous-goupe ab\'elien fini de $G_{SC}\times Z_{G}^0$. L'automorphisme $\theta$ se rel\`eve en un automorphisme de $G_{SC}$, donc aussi  en un automorphisme de $G_{SC}\times Z_{G}^0$
qui pr\'eserve $\Xi$. On note tous ces automorphismes $\theta$. Fixons une cl\^oture alg\'ebrique $\bar{F}$ de $F$, posons $\Gamma_{F}=Gal(\bar{F}/F)$. Fixons une extension galoisienne finie $F'$ de $F$ telle que $\Gamma_{F'}$ agisse trivialement sur $\Xi$. Notons $\hat{\Xi}$ le groupe $Hom(\Xi,\bar{F}^{\times})$ des caract\`eres alg\'ebriques  de $\Xi$ et  ${\mathbb Z}[\hat{\Xi}]$ le ${\mathbb Z}$-module libre de base $\hat{\Xi}$. Pour plus de clart\'e, si $\mu\in \hat{\Xi}$, on note $\underline{\mu}$ l'\'el\'ement de base de ${\mathbb Z}[\hat{\Xi}]$. Le module  ${\mathbb Z}[\hat{\Xi}]$ est muni d'une action de $\Gamma_{F}$ triviale sur $\Gamma_{F'}$. Notons $\Lambda$ l'ensemble des fonctions $f:\Gamma_{F'}\backslash \Gamma_{F}\to{\mathbb Z}[\hat{\Xi}]$. On le munit de l'action galoisienne  par translations \`a droite. On note $C$ le tore sur $F$ tel que $X_{*}(C)=\Lambda$, muni de cette action galoisienne.
Par ailleurs $\theta$ se transporte en une permutation de $\hat{\Xi}$, puis en un automorphisme de  ${\mathbb Z}[\hat{\Xi}]$, puis en un automorphisme de $\Lambda$ (par action sur l'espace d'arriv\'ee), enfin en un automorphisme $\theta_{C}$ de $C$. Il est clair que $\theta_{C}$ est d\'efini sur $F$. On d\'efinit deux applications
$$\Xi\to {\mathbb Z}[\hat{\Xi}]\otimes_{{\mathbb Z}}\bar{F}^{\times}$$
$${\mathbb Z}[\hat{\Xi}]\otimes_{{\mathbb Z}}\bar{F}^{\times}\to C.$$
La premi\`ere envoie $\xi\in \Xi$ sur $\sum_{\mu\in \hat{\Xi}}\underline{\mu}\otimes \mu(\xi)$. La seconde est d\'eduite de l'application ${\mathbb Z}[\hat{\Xi}]\to \Lambda$ qui \`a $x\in {\mathbb Z}[\hat{\Xi}]$ associe la fonction $f$ d\'efinie par $f(\sigma)=\sigma x$. Notons $\iota:\Xi\to C$ la compos\'ee des deux applications. Elle est injective, \'equivariante pour les actions galoisiennes et v\'erifie $\iota\circ \theta=\theta_{C}\circ\iota$. Posons
$$G'=(G_{SC}\times Z_{G}^0\times C)/\{(\xi,\iota(\xi)); \xi\in \Xi\}.$$
On a bien la suite exacte
$$1\to C\to G'\stackrel{p}{\to} G\to 1.$$
La condition (ii) est v\'erifi\'ee parce que $C$ est un tore induit, donc $H^1(\Gamma_{F'},C)=0$ pour toute extension $F'$ de $F$. Enfin, le groupe d\'eriv\'e de $G'$ est simplement connexe parce que l'homomorphisme naturel $G_{SC}\to G'$ est injectif. Le produit des  automorphismes $\theta$ de $G_{SC}\times Z_{G}^0$ et $\theta_{C}$ de $C$ se descend en un automorphisme $\theta'$ de $G'$. Sa restriction \`a $Z_{G'}$ est d'ordre fini. On  introduit une vari\'et\'e $\tilde{G}'$ muni d'un isomorphisme d\'efini sur $F$ de $G'$ sur $\tilde{G}'$. On note $g\mapsto g\boldsymbol{\theta}'$ cet isomorphisme. On d\'efinit les actions \`a gauche et \`a droite de $G'$ sur $\tilde{G}'$ par $(x,g\boldsymbol{\theta},y)\mapsto xg\theta'(y)\boldsymbol{\theta}'$. On d\'efinit $\tilde{p}:\tilde{G}'\to \tilde{G}$ par $\tilde{p}(g\boldsymbol{\theta}')=g\gamma$. Le couple $(p,\tilde{p})$ est un homomorphisme d'espaces tordus. Cela prouve (1).

On montre ensuite

(2) les assertions de la proposition sont v\'erifi\'ees si le groupe d\'eriv\'e de $G$ est simplement connexe.

On fixe de nouveau $\gamma\in \tilde{G}(F)$ et on en d\'eduit un automorphisme $\theta$ de $G$. On introduit le tore $D=Z_{G}^0/(Z_{G_{SC}}\cap Z_{G}^0)$. On a la suite exacte
$$1\to G_{SC} \to G\stackrel{\underline{d}}{\to}D\to 1,$$
d'o\`u
$$1\to G_{SC}(F)\to G(F)\stackrel{\underline{d}}{\to}  D(F).$$
De $\theta$ se d\'eduit un automorphisme encore not\'e $\theta$ de $D$. Tout caract\`ere de $G_{SC}(F)$ \'etant trivial,  $\omega$ se factorise en un caract\`ere du quotient $G_{SC}(F)\backslash G(F)$, que l'on peut prolonger en un caract\`ere de $D(F)$. On note $\omega_{D}$ ce prolongement. Notons $Z_{G}^{\theta,0}$ et $D^{\theta,0}$ les composantes neutres des sous-groupes de points fixes par $\theta$ dans $Z_{G}$ et $D$. L'homomorphisme $Z_{G}^{\theta,0}\to D^{\theta,0}$ est surjectif, donc $\underline{d}(Z_{G}^{\theta,0}(F))$ est d'indice fini dans $D^{\theta,0}(F)$. Par hypoth\`ese, $\omega$ est trivial sur $Z_{G}(F)^{\theta}$, a fortiori sur $Z_{G}^{\theta,0}(F)$. Donc la restriction de $\omega_{D}$ \`a $D^{\theta,0}(F)$ est d'ordre fini. Notons $n_{1}$ l'ordre de $\omega_{D}$, $n_{2}$ celui de l'automorphisme $\theta$ de $D$ et posons $N=n_{1}n_{2}$. Si $N=1$, on a $D=D^{\theta,0}$ puis $\omega=1$. On prend $G'=G$, $\tilde{G}'=\tilde{G}$ et c'est termin\'e. Supposons $N>1$. Posons
$$C=\{(x_{1},...,x_{N})\in D^N; \prod_{i=1,...,N}\theta_{D}^{i}x_{i}=1\}.$$
C'est un tore isomorphe \`a $D^{N-1}$. Posons $G'=G\times C$. L'extension
$$1\to C\to G'\stackrel{p}{\to }G\to 1$$
v\'erifie les conditions (i) et (ii). D\'efinissons un automorphisme $\theta'$ de $G'$ par
$$\theta'(g,x_{1},...,x_{N})=(\theta(g),\underline{d}(g)x_{N},x_{1},...,x_{N-2},\underline{d}\circ\theta(g^{-1})x_{N-1}).$$
On v\'erifie que
$$(\theta')^N(g,x_{1},...,x_{N})=(\theta^N(g),x_{1},...,x_{N}).$$
Il en r\'esulte que la restriction de $\theta'$ \`a $Z_{G'}$ est d'ordre fini. On d\'efinit l'espace tordu $\tilde{G}'$ et l'application $\tilde{p}$ comme dans la preuve de (1). D\'efinissons un caract\`ere $\mu'$ de $G'(F)$ par
$$\mu'(g,x_{1},...,x_{N})=\prod_{i=1,...,N-1}\prod_{j=0,...,i-1}\omega_{D}^{-1}\circ\theta^j(x_{i}).$$
On calcule
$$\mu'\circ(1-\theta')(g,x_{1},...,x_{N})=\omega_{D}^{-1}(x_{1}x_{N}^{-1}\underline{d}(g)^{-1})\prod_{i=2,...,N-1}\prod_{j=0,...,i-1}\omega_{D}^{-1}\circ\theta^j(x_{i}x_{i-1}^{-1})$$
$$=\omega_{D}\circ\underline{d}(g)\omega_{D}(x_{1}^{-1}x_{N})\left(\prod_{i=2,...,N-1}\prod_{j=0,...,i-1}\omega_{D}^{-1}\circ\theta^j(x_{i})\right)\left(\prod_{i=1,...,N-2}\prod_{j=0,...,i}\omega_{D}^{-1}\circ\theta^j(x_{i}^{-1})\right)$$
$$=\omega(g)\omega_{D}(x_{N})\left(\prod_{i=1,...,N-2}\omega_{D}\circ\theta^{i}(x_{i})\right)\prod_{j=0,...,N-2}\omega_{D}^{-1}\circ\theta^j(x_{N-1}).$$
On utilise la relation
$$x_{N}=\prod_{i=1,...,N-1}\theta^{i}(x_{i}^{-1})$$
et on obtient
$$\mu'\circ(1-\theta')(g,x_{1},...,x_{N})=\omega(g)\prod_{j=0,...,N-1}\omega_{D}^{-1}\circ\theta^{j}(x_{N-1}).$$
L'application $x\mapsto \prod_{j=0,...,n_{2}-1}\theta^j(x)$ envoie $D$ dans $D^{\theta,0}$, donc $D(F)$ dans $D^{\theta,0}(F)$. Donc, pour tout $x\in D(F)$, $\prod_{j=0,...,N-1}\theta^j(x)$ est la puissance $n_{1}$-i\`eme d'un \'el\'ement de $D^{\theta,0}(F)$. Il en r\'esulte que
$$\prod_{j=0,...,N-1}\omega_{D}\circ\theta^j(x)=1$$
pour tout $x\in D(F)$. L'\'egalit\'e ci-dessus devient
$$\mu'\circ(1-\theta')(g,x_{1},...,x_{N})=\omega(g),$$
comme on le voulait. Cela prouve (2).

Dans le cas g\'en\'eral, on commence par construire  un couple v\'erifiant (1). On affecte les objets relatifs \`a ce couple d'un indice $1$: $G'_{1}$, $\tilde{G}'_{1}$, $p_{1}$, $\tilde{p}_{1}$. On applique (2) \`a ce couple et au caract\`ere $\omega\circ p_{1}$. On obtient des objets $G'$, $\tilde{G}'$, $\mu'$ et des applications $p_{2}:G'\to G'_{1}$, $\tilde{p}_{2}:\tilde{G}'\to \tilde{G}'_{1}$. On pose $p=p_{1}\circ p_{2}$, $\tilde{p}=\tilde{p}_{1}\circ\tilde{p}_{2}$ et on note $C$ l'image r\'eciproque de $C_{1}$ dans $G'$. Il est clair que les objets $G'$, $\tilde{G}'$, $C$, $p$ et $\tilde{p}$ v\'erifient les conditions de l'\'enonc\'e. $\square$

 \bigskip
 
 \subsection{Repr\'esentations}
Soit $\omega$ un caract\`ere unitaire de $G(F)$. On suppose que sa restriction \`a $Z_{G}(F)^{\theta}$ est triviale. On appelle $\omega$-repr\'esentation lisse de $\tilde{G}(F)$ un couple $(\pi,\tilde{\pi})$, o\`u $\pi$ est une repr\'esentation lisse de $G(F)$ dans un espace $V_{\pi}$ et $\tilde{\pi} $ est une application de $\tilde{G}(F)$ dans le groupe des automorphismes de $V_{\pi}$ telle que $\tilde{\pi}(g\gamma g')=\pi(g)\tilde{\pi}(\gamma)\pi(g')\omega(g')$ pour tous $g,g'\in G(F)$ et $\gamma\in \tilde{G}(F)$. 

{\bf Remarque.} Dans le cas o\`u $F$ est archim\'edien, il faudrait distinguer, comme on l'a dit en 1.9, l'espace de la repr\'esentation et son sous-espace $V_{\pi}$ des vecteurs $K$-finis.
\bigskip

En pratique, on parlera plut\^ot de la $\omega$-repr\'esentation $\tilde{\pi}$, en occultant la repr\'esentation $\pi$ sous-jacente. Pour une $\omega$-repr\'esentation lisse $\tilde{\pi}$ et pour $z\in {\mathbb C}^{\times}$, on d\'efinit $z\tilde{\pi}$ par  $(z\tilde{\pi})(\gamma)=z\tilde{\pi}(\gamma)$.
 
 {\bf Variante.} On dira que $\tilde{\pi}$ est unitaire s'il existe un produit hermitien d\'efini positif sur $V_{\pi}$ tel que $\tilde{\pi}(\gamma)$ soit un op\'erateur unitaire pour tout $\gamma\in \tilde{G}(F)$. Notons ${\mathbb U}$ le groupe des nombres complexes de module $1$. Pour une $\omega$-repr\'esentation lisse unitaire $\tilde{\pi}$ et pour $z\in {\mathbb U}$, $z\tilde{\pi}$ est encore unitaire.
 \bigskip
 
  On dit que $\tilde{\pi}$ est admissible si et seulement si $\pi$ l'est. On dit que $\tilde{\pi}$ est temp\'er\'ee si et seulement si $\tilde{\pi}$ est unitaire et $\pi$ est temp\'er\'ee. Nos hypoth\`eses que  $\omega$ unitaire et que l'automorphisme $\theta$ de $A_{G}$ est d'ordre fini impliquent que $\tilde{\pi}$ est de longueur finie si et seulement si $\pi$ l'est. Il y a une notion naturelle d'irr\'eductibilit\'e et une notion plus fine de $G$-irr\'eductibilit\'e: $\tilde{\pi}$ est $G$-irr\'eductible si et seulement si $\pi$ est irr\'eductible. Notons $Irr(G(F);\theta,\omega)$ l'ensemble des classes de repr\'esentations lisses irr\'eductibles  $\pi$ de $G(F)$ telles que $\omega\otimes \pi \simeq \pi\circ \theta$, o\`u $\pi\circ\theta$ est la classe de $\pi\circ ad_{\gamma}$ pour un quelconque $\gamma\in \tilde{G}(F)$. L'application $\tilde{\pi}\mapsto \pi$ est une surjection de l'ensemble des $\omega$-repr\'esentations lisses $G$-irr\'eductibles de $\tilde{G}(F)$ sur $Irr(G(F);\theta,\omega)$. Les fibres sont isomorphes \`a ${\mathbb C}^{\times}$.
  
  La mesure de Haar fix\'ee sur $G(F)$ d\'etermine une mesure sur $\tilde{G}(F)$.   Pour $f\in C_{c}^{\infty}(\tilde{G}(F))$, on d\'efinit l'op\'erateur
$$\tilde{\pi}(f)=\int_{\tilde{G}(F)}f(\gamma)\tilde{\pi}(\gamma)d\gamma.$$
Supposons $\tilde{\pi}$ admissible. Alors la trace de cet op\'erateur est bien d\'efinie. A $\tilde{\pi}$ est associ\'e son caract\`ere, qui est  la distribution sur $C_{c}^{\infty}(\tilde{G}(F))$ d\'efinie par
$$I_{\tilde{G}}(\tilde{\pi},f)=trace(\tilde{\pi}(f)).$$
Si $g\in G(F)$ et si l'on note $^gf$ la fonction d\'efinie par $^gf(\gamma)=f(g^{-1}\gamma g)$, on a
$$I_{\tilde{G}}(\tilde{\pi},{^gf})=\omega(g)^{-1}I_{\tilde{G}}(\tilde{\pi},f).$$ 
On note $D_{spec}(\tilde{G}(F);\omega)$ l'espace de distributions engendr\'e par les caract\`eres de $\omega$-repr\'esentations admissibles de longueur finie. Il est \'evidemment engendr\'e par les caract\`eres de $\omega$-repr\'esentations irr\'eductibles. Une $\omega$-repr\'esentation irr\'eductible  a un caract\`ere non  nul si et seulement si elle est $G$-irr\'eductible, cf. [L] proposition A.4.1. Tout \'el\'ement $\pi\in Irr(G(F);\theta,\omega)$ d\'etermine une droite $D_{\pi}$ dans $D_{spec}(\tilde{G}(F);\omega)$, \`a savoir la droite port\'ee par le caract\`ere de $\tilde{\pi}$, o\`u $\tilde{\pi}$ est un prolongement quelconque de $\pi$ en une repr\'esentation de $\tilde{G}(F)$. On a 
$$(1)\qquad D_{spec}(\tilde{G}(F);\omega)=\oplus_{\pi\in Irr(G(F);\theta,\omega)}D_{\pi};$$
cf. [L] proposition A.4.1 dans le cas non archim\'edien; la preuve est similaire dans le cas archim\'edien;
 
(2) l'espace  $D_{spec}(\tilde{G}(F);\omega)$ est form\'e de distributions localement int\'egrables sur $\tilde{G}(F)$ et lisses sur les \'el\'ements fortement r\'eguliers. 

Preuve. Notons ${\bf 1}$ le caract\`ere trivial de $G(F)$. Si $\omega={\bf 1}$, l'assertion est due \`a Bouaziz dans le cas archim\'edien ([B] th\'eor\`eme 2.1.1) et \`a Clozel dans le cas non-archim\'edien ([C]). Supposons qu'il existe un caract\`ere unitaire $\mu$ de $G(F)$ tel que $\omega=\mu\circ(1-\theta)$ (cf. 2.4). Soit $\tilde{\pi}$ une $\omega$-repr\'esentation de $\tilde{G}(F)$ de repr\'esentation sous-jacente $\pi$. Posons $\pi_{1}=\pi\otimes \mu$. Fixons $\gamma_{0}\in \tilde{G}(F)$ et d\'efinissons $\tilde{\pi}_{1}$ par $\tilde{\pi}_{1}(g\gamma_{0})=\mu(g)\tilde{\pi}(g\gamma_{0})$. On v\'erifie que $\tilde{\pi}_{1}$ est une ${\bf 1}$-repr\'esentation de $\tilde{G}(F)$, de repr\'esentation sous-jacente $\pi_{1}$. Son caract\`ere est localement int\'egrable, donc associ\'e \`a  une fonction localement int\'egrable $\gamma\mapsto\Theta(\tilde{\pi}_{1},\gamma)$ sur $\tilde{G}(F)$. Le caract\`ere de $\tilde{\pi}$ est alors associ\'e \`a la fonction $g\gamma_{0}\mapsto \mu(g)^{-1}\Theta(\tilde{\pi}_{1},g\gamma_{0})$, qui est elle-aussi localement int\'egrable. Dans le cas g\'en\'eral, on introduit des objets $G'$, $\tilde{G}'$, $C$, $p$, $\tilde{p}$ v\'erifiant la proposition 2.4. On pose $\pi'=\pi\circ p$, $\tilde{\pi}'=\tilde{\pi}\circ\tilde{p}$, $\omega'=\omega\circ p$. Le terme $\tilde{\pi}'$ est une $\omega'$-repr\'esentation de $\tilde{G}'(F)$. L'hypoth\`ese du cas pr\'ec\'edent est v\'erifi\'ee: $\omega'$ est de la forme $\mu'\circ(1-\theta')$. Donc le caract\`ere de $\tilde{\pi}'$ est localement int\'egrable, associ\'e \`a une fonction localement int\'egrable $\gamma'\mapsto \Theta(\tilde{\pi}',\gamma')$ sur $\tilde{G}'(F)$. Le caract\`ere central de la repr\'esentation $\pi'$ est trivial sur $C(F)$ par construction. Il en r\'esulte que la fonction ci-dessus est invariante par translations (\`a droite ou \`a gauche) par $C(F)$. Il existe donc une fonction $\gamma\mapsto \Theta(\tilde{\pi},\gamma)$ sur $\tilde{G}(F)$ telle que l'on ait l'\'egalit\'e $\Theta(\tilde{\pi}',\gamma')=\Theta(\tilde{\pi},\tilde{p}(\gamma'))$ pour tout $\gamma'\in \tilde{G}'(F)$. On v\'erifie ais\'ement que cette fonction est elle-aussi localement int\'egrable et que sa distribution associ\'ee est le caract\`ere de $\tilde{\pi}$. $\square$

Soit $\tilde{M}$ un espace de Levi de $\tilde{G}$ tel que $M_{0}\subset M$ et soit $\tilde{\sigma}$ une $\omega$-repr\'esentation lisse de $\tilde{M}(F)$. Pour $\tilde{P}\in {\cal P}(\tilde{M})$, on d\'efinit l'induite habituelle $Ind_{P}^G(\sigma)$. On d\'efinit une repr\'esentation $\tilde{\pi}=Ind_{\tilde{P}}^{\tilde{G}}(\tilde{\sigma})$ de $\tilde{G}(F)$ dans l'espace de cette induite par la formule
$$(3) \qquad (\tilde{\pi}(\gamma)f)(g)=\omega(g')\delta_{P}(\gamma')^{1/2}\tilde{\sigma}(\gamma')f(g'),$$
o\`u:

- on a \'ecrit $g\gamma =\gamma'g'$, avec $\gamma'\in \tilde{M}(F)$;

- $\delta_{P}(\gamma')$ est la valeur absolue du d\'eterminant de l'action de $ad_{\gamma'}$ dans l'alg\`ebre de Lie du radical unipotent de $P$.

{\bf Remarque.} Il peut arriver que $\tilde{\sigma}$ soit $M$-irr\'eductible mais que $Ind_{\tilde{P}}^{\tilde{G}}(\tilde{\sigma})$ ait un caract\`ere nul. Cela se produit quand $Ind_{P}^G(\sigma)$ n'est pas irr\'eductible et que l'action de $\tilde{G}(F)$ permute sans point fixe ses composantes irr\'eductibles. Donnons un exemple. Supposons $F$ non-archim\'edien et $G=SL_{2}$,   prenons pour $M$ le tore diagonal et pour $P$ le Borel sup\'erieur. Soient $\chi$ un caract\`ere quadratique non trivial de $F^{\times}$ et $D\in F^{\times}$ tel que $\chi(D)=-1$. Supposons $\tilde{G}=\{x\in GL_{2}; det(x)=D\}$. Prenons pour $\sigma$ le caract\`ere 
$$\left(\begin{array}{cc}a&0\\0&a^{-1}\\ \end{array}\right)\mapsto \chi(a)$$
de $M(F)$ et pour $\tilde{\sigma}$ la repr\'esentation
$$\left(\begin{array}{cc}a&0\\0&d\\ \end{array}\right)\mapsto \chi(a)$$
de $\tilde{M}(F)$.
  On v\'erifie que le caract\`ere de $Ind_{\tilde{P}}^{\tilde{G}}(\tilde{\sigma})$ co\"{\i}ncide  avec 
  la restriction \`a $\tilde{G}(F)$ du caract\`ere de la repr\'esentation $\pi$ de $GL_{2}(F)$ induite du caract\`ere $\chi\times 1$ du tore diagonal. Or $\pi\simeq (\chi\circ det)\pi$, donc le caract\`ere de $\pi$ est nul sur les \'el\'ements de $GL_{2}(F)$ de d\'eterminant $D$.

\bigskip

 \subsection{Torsion par un caract\`ere}
 Notons $G(F)^1$ le noyau de l'homomorphisme $H_{\tilde{G}}$. Il est invariant par $ad_{\gamma}$ pour tout $\gamma\in \tilde{G}(F)$. Posons $\tilde{{\cal A}}_{\tilde{G},F}=G(F)^1\backslash \tilde{G}(F)=\tilde{G}(F)/G(F)^1$ et notons $\tilde{H}_{\tilde{G}}:\tilde{G}(F)\to \tilde{{\cal A}}_{\tilde{G},F}$ la projection naturelle. L'ensemble $\tilde{{\cal A}}_{\tilde{G},F}$ est un espace principal homog\`ene sous l'action de ${\cal A}_{\tilde{G},F}$ (les actions \`a gauche et \`a droite co\"{\i}ncident). On pose $\tilde{{\cal A}}_{\tilde{G}}={\cal A}_{\tilde{G}}\otimes_{{\cal A}_{\tilde{G},F}}\tilde{{\cal A}}_{\tilde{G},F}$. C'est un espace affine sous ${\cal A}_{\tilde{G}}$. Par abus de notation, on note $\tilde{{\cal A}}_{\tilde{G}}^*$ l'espace des fonctions affines sur $\tilde{{\cal A}}_{\tilde{G}}$, c'est-\`a-dire les fonctions 
 $$\begin{array}{cccc}\tilde{\lambda}:&\tilde{{\cal A}}_{\tilde{G}}&\to& {\mathbb R}\\ &\tilde{H}&\mapsto &<\tilde{\lambda},\tilde{H}>\\ \end{array}$$
  telles qu'il existe $\lambda\in {\cal A}_{\tilde{G}}^*$ de sorte que $<\tilde{\lambda},H+\tilde{H}>=<\lambda,H>+<\tilde{\lambda},\tilde{H}>$ pour tout $H\in {\cal A}_{\tilde{G}}$ et $\tilde{H}\in \tilde{{\cal A}}_{\tilde{G}}$. On  d\'efinit de fa\c{c}on \'evidente le complexifi\'e $\tilde{{\cal A}}_{\tilde{G},{\mathbb C}}^*$ et son sous-espace imaginaire $i\tilde{{\cal A}}_{\tilde{G}}^*$. On note $\tilde{{\cal A}}_{\tilde{G},F}^{\vee}$ le sous-groupe des $\tilde{\lambda}\in \tilde{{\cal A}}_{\tilde{G}}^*$ tels que $\tilde{\lambda}(\tilde{{\cal A}}_{\tilde{G},F})\subset 2\pi{\mathbb Z}$. On pose $i\tilde{{\cal A}}_{\tilde{G},F}^*=(i\tilde{{\cal A}}_{\tilde{G}}^*)/(i\tilde{{\cal A}}_{\tilde{G},F}^{\vee})$. On a des suites exactes
 $$0\to {\mathbb R}\to \tilde{{\cal A}}_{\tilde{G}}^*\to {\cal A}_{\tilde{G}}^*\to 0,$$
 $$0\to {\mathbb C}\to \tilde{{\cal A}}_{\tilde{G},{\mathbb C}}^*\to {\cal A}_{\tilde{G},{\mathbb C}}^*\to 0,$$
 $$0\to i{\mathbb R}/2\pi i{\mathbb Z}\to i\tilde{{\cal A}}_{\tilde{G},F}^*\to i{\cal A}_{\tilde{G},F}^*\to 0.$$
 Le choix d'un point-base $\gamma\in \tilde{G}(F)$ scinde ses suites: on identifie par exemple $\tilde{\cal A}_{\tilde{G}}^*$ au sous-espace des fonctions affines qui valent $0$ au point $\tilde{H}_{\tilde{G}}(\gamma)$.
 Pour $\tilde{\lambda}\in \tilde{{\cal A}}_{\tilde{G}}^*$ (ou $\tilde{\lambda}\in  i\tilde{{\cal A}}_{\tilde{G},F}^*$), on note sans plus de commentaire $\lambda$ son image dans ${\cal A}_{\tilde{G}}^*$ (ou $i{\cal A}_{\tilde{G},F}^*$).

Soit $\tilde{\pi}$ une $\omega$-repr\'esentation lisse de $\tilde{G}(F)$. Pour $\tilde{\lambda}\in \tilde{{\cal A}}_{\tilde{G},{\mathbb C}}^*$, on d\'efinit la repr\'esentation  $\tilde{\pi}_{\tilde{\lambda}}$ par $\tilde{\pi}_{\tilde{\lambda}}(\gamma)=e^{<\tilde{\lambda},\tilde{H}_{\tilde{G}}(\gamma)>}\tilde{\pi}(\gamma)$. Sa repr\'esentation sous-jacente de $G(F)$ est $\pi_{\lambda}$. Evidemment, $\tilde{\pi}_{\tilde{\lambda}}$ ne d\'epend que de l'image de $\tilde{\lambda}$ dans $\tilde{{\cal A}}_{\tilde{G},{\mathbb C}}^*/i\tilde{{\cal A}}_{\tilde{G},F}^{\vee}$. Remarquons que:

(1)  si $\pi$ est irr\'eductible, le stabilisateur de $\tilde{\pi}$ dans $i\tilde{{\cal A}}_{\tilde{G},F}^*$ s'identifie au stabilisateur $Stab(i{\cal A}_{\tilde{G},F}^*,\pi)$  de $\pi$ dans $i{\cal A}_{\tilde{G},F}^*$. 

Preuve. Il est clair que le premier groupe se projette injectivement dans le second. Inversement, soit $\lambda\in Stab(i{\cal A}_{\tilde{G},F}^*,\pi)$ et fixons arbitrairement $\tilde{\lambda}\in i\tilde{{\cal A}}_{\tilde{G},F}^*$ au-dessus de $\lambda$. Alors $\tilde{\pi}_{\tilde{\lambda}}$ est isomorphe \`a une repr\'esentation prolongeant $\pi$, donc \`a $z\tilde{\pi}$ pour un certain $z\in {\mathbb C}^{\times}$. On en d\'eduit $\tilde{\pi}_{n\tilde{\lambda}}\simeq z^n\tilde{\pi}$ pour tout $n\in {\mathbb N}$. On choisit $n$ tel que $n\lambda=0$ dans $i{\cal A}_{\tilde{G},F}^*$. Alors $n\tilde{\lambda}\in i{\mathbb R}/2\pi i{\mathbb Z}$ et $\tilde{\pi}_{n\tilde{\lambda}}=e^{n\tilde{\lambda}}\tilde{\pi}$, donc $z^n=e^{n\tilde{\lambda}}$ et $z$ est de module $1$. En \'ecrivant $z=e^{x}$, avec $x\in i{\mathbb R}/2\pi i{\mathbb Z}$, l'\'el\'ement $\tilde{\lambda}-x$ appartient au stabilisateur de $\tilde{\pi}$. $\square$

 On devra consid\'erer des fonctions \`a valeurs complexes $\varphi$ d\'efinies sur $\tilde{{\cal A}}_{\tilde{G},{\mathbb C}}$, resp. $i\tilde{{\cal A}}_{\tilde{G},F}^*$. La plupart v\'erifieront la condition
 
 (2) pour tout $z\in {\mathbb C}\subset \tilde{{\cal A}}_{\tilde{G},{\mathbb C}}^*$ et tout $\tilde{\lambda}\in \tilde{{\cal A}}_{\tilde{G},{\mathbb C}}^*$, $\varphi(z+\tilde{\lambda})=e^{z}\varphi(\tilde{\lambda})$
 
\noindent  resp.  une condition similaire. Modulo le choix d'un point-base permettant de scinder les suites exactes ci-dessus, une telle fonction s'identifie \`a une fonction sur ${\cal A}_{\tilde{G},{\mathbb C}}^*$, resp. $i{\cal A}_{\tilde{G},F}^*$. 

Une fonction $\varphi:i\tilde{{\cal A}}_{\tilde{G},F}^*\to {\mathbb C}$ est dite de Paley-Wiener si et seulement s'il existe une fonction $b:\tilde{{\cal A}}_{\tilde{G},F}\to {\mathbb C}$, lisse et \`a support compact, de sorte que
$$\varphi(\tilde{\lambda})=\int_{\tilde{{\cal A}}_{\tilde{G},F}}b(X)e^{<\tilde{\lambda},X>}\,dX.$$
Il revient au m\^eme de demander que $\varphi$ v\'erifie (2) et que la fonction sur $i{\cal A}_{\tilde{G},F}^*$ d\'eduite de $\varphi$ comme ci-dessus soit de Paley-Wiener.

\bigskip

\subsection{Caract\`eres pond\'er\'es}
Soit $\tilde{M}$ un espace de Levi de $\tilde{G}$ tel que $M_{0}\subset M$. Soit $\tilde{\pi}$ une $\omega$-repr\'esentation admissible et de longueur finie de $\tilde{M}(F)$. Fixons $\tilde{P}\in {\cal P}(\tilde{M})$, introduisons la repr\'esentation $Ind_{\tilde{P}}^{\tilde{G}}(\tilde{\pi})$, que l'on r\'ealise dans l'espace $V_{\pi,P}$. Supposons dans un premier temps que $\tilde{\pi}$ est en position g\'en\'erale, en ce sens que les op\'erateurs d'entrelacement intervenant ci-dessous sont bien d\'efinis. Pour $\tilde{Q}\in {\cal P}(\tilde{M})$, l'op\'erateur $J_{P\vert Q}(\pi)J_{Q\vert P}(\pi)$ est un automorphisme de $V_{\pi,P}$. Notons $\mu_{Q\vert P}(\pi)$ son inverse. Pour $\Lambda\in i{\cal A}_{\tilde{M}}^*$, posons
$${\cal M}(\pi;\Lambda,\tilde{Q})=\mu_{Q\vert P}(\pi)^{-1}\mu_{Q\vert P}(\pi_{\Lambda/2})J_{Q\vert P}(\pi)^{-1}J_{Q\vert P}(\pi_{\Lambda}).$$
On souligne la pr\'esence d'un indice $\Lambda/2$ dans cette formule.  La famille $({\cal M}(\pi;\Lambda,\tilde{Q}))_{\tilde{Q}\in {\cal P}(\tilde{M})}$ est une $(\tilde{G},\tilde{M})$-famille \`a valeurs op\'erateurs. On en d\'eduit un op\'erateur ${\cal M}_{\tilde{M}}^{\tilde{G}}(\pi;\Lambda)$. On pose ${\cal M}_{\tilde{M}}^{\tilde{G}}(\pi)={\cal M}_{\tilde{M}}^{\tilde{G}}(\pi;0)$. Si $F$ est archim\'edien, notons $C_{c}^{\infty}(\tilde{G}(F),K)$ le sous-espace des \'el\'ements de $C_{c}^{\infty}(\tilde{G}(F))$ qui sont $K$-finis \`a droite et \`a gauche. Pour unifier les notations, posons $C_{c}^{\infty}(\tilde{G}(F),K)=C_{c}^{\infty}(\tilde{G}(F))$ dans le cas o\`u $F$ est non-archim\'edien.  Le caract\`ere pond\'er\'e de $\tilde{\pi}$ est la distribution sur $C_{c}^{\infty}(\tilde{G}(F),K)$ d\'efinie par 
$$J_{\tilde{M}}^{\tilde{G}}(\tilde{\pi},f)=trace({\cal M}_{\tilde{M}}^{\tilde{G}}(\pi)Ind_{\tilde{P}}^{\tilde{G}}(\tilde{\pi},f)).$$
On a choisi un parabolique $\tilde{P}$, mais on v\'erifie que cette trace ne d\'epend pas de ce choix.

{\bf Remarque.} On verra en 5.4 que cette d\'efinition s'\'etend \`a tout $C_{c}^{\infty}(\tilde{G}(F))$ au moins si $\tilde{\pi}$ est temp\'er\'ee.
\bigskip 

Levons l'hypoth\`ese que $\tilde{\pi}$ est en position g\'en\'erale. Pour $\tilde{\pi}$ quelconque et pour $\tilde{\lambda}\in \tilde{{\cal A}}_{\tilde{M},{\mathbb C}}^*$, $\tilde{\pi}_{\tilde{\lambda}}$ est en position g\'en\'erale pour $\tilde{\lambda}$ hors d'un sous-ensemble ferm\'e de mesure nulle. On peut donc d\'efinir l'op\'erateur ${\cal M}_{\tilde{M}}^{\tilde{G}}(\pi_{\lambda})$ et la distribution $f\mapsto J_{\tilde{M}}^{\tilde{G}}(\tilde{\pi}_{\tilde{\lambda}},f)$ pour $\tilde{\lambda}$ dans un ouvert dense. Il est clair que ces termes v\'erifient la condition (2) de 2.6 et s'\'etendent en des fonctions m\'eromorphes de $\tilde{\lambda}$ d\'efinies pour tout $\tilde{\lambda}$. Si ces fonctions sont r\'eguli\`eres en $\tilde{\lambda}=0$, on d\'efinit ${\cal M}_{\tilde{M}}^{\tilde{G}}(\pi)$ et $J_{\tilde{M}}^{\tilde{G}}(\tilde{\pi},f)$ comme leurs valeurs en ce point $\tilde{\lambda}=0$.

\ass{Proposition}{Si $\tilde{\pi}$ est unitaire, les fonctions ci-dessus sont r\'eguli\`eres en $\tilde{\lambda}=0$.}

Cf. [A4] proposition 2.3. Arthur traite le cas d'une repr\'esentation irr\'eductible d'un groupe connexe mais sa preuve s'\'etend \`a notre situation. 

Remarquons que, dans notre construction, on n'a pas impos\'e \`a $\tilde{\pi}$ d'\^etre irr\'eductible. On v\'erifie que l'op\'erateur ${\cal M}_{\tilde{M}}^{\tilde{G}}(\pi)$ d\'epend fonctoriellement de $\tilde{\pi}$, pourvu qu'il soit d\'efini, et que, si l'on a une suite exacte
$$1\to \tilde{\pi}_{1}\to \tilde{\pi}_{2}\to \tilde{\pi}_{3}\to 1,$$
on a l'\'egalit\'e
$$(1) \qquad J_{\tilde{M}}^{\tilde{G}}(\tilde{\pi}_{2},f)=J_{\tilde{M}}^{\tilde{G}}(\tilde{\pi}_{1},f)+J_{\tilde{M}}^{\tilde{G}}(\tilde{\pi}_{3},f),$$
pourvu que tous les termes soient d\'efinis.

On a:

(2) supposons que $\tilde{\pi}$ soit  irr\'eductible mais pas $M$-irr\'eductible; alors $J_{\tilde{M}}^{\tilde{G}}(\tilde{\pi},f)=0$ pourvu que ce terme soit d\'efini.

Preuve. Puisque $\tilde{\pi}_{\tilde{\lambda}}$ v\'erifie les m\^emes hypoth\`eses que $\tilde{\pi}$, il suffit de consid\'erer le cas o\`u $\tilde{\pi}$ est en position g\'en\'erale. La d\'ecomposition de $\pi$ en composantes irr\'eductibles pour $M(F)$ induit une d\'ecomposition de $Ind_{P}^G(\pi)$. L'op\'erateur ${\cal M}_{\tilde{M}}^{\tilde{G}}(\pi)$ pr\'eserve chaque composante tandis  que $Ind_{\tilde{P}}^{\tilde{G}}(\tilde{\pi},f)$ les permute sans point fixe. $\square$

Plus g\'en\'eralement,

(3) supposons que le caract\`ere de $\tilde{\pi}$ soit nul;  alors $J_{\tilde{M}}^{\tilde{G}}(\tilde{\pi},f)=0$ pourvu que ce terme soit d\'efini.

Preuve. On peut encore supposer $\tilde{\pi}$ en position g\'en\'erale. Gr\^ace \`a (1), on peut supposer que $\tilde{\pi}$ est somme d'irr\'eductibles. On peut supprimer les irr\'eductibles non $M$-irr\'eductibles: cela ne modifie pas l'hypoth\`ese puisque leur caract\`ere est nul, ni la conclusion d'apr\`es (2).  Notons $\rho_{1},...,\rho_{k}$ les diff\'erentes repr\'esentations irr\'eductibles de $M(F)$ intervenant dans $\pi$. Pour chaque $\rho_{i}$, fixons un prolongement $\tilde{\rho}_{i}$ de $\rho_{i}$ \`a $\tilde{M}(F)$. Alors on a une \'egalit\'e
$$\tilde{\pi}=\oplus_{i=1,...,k}\oplus_{j=1,...,l_{i}}z_{i,j}\tilde{\rho}_{i},$$
pour des familles $(z_{i,j})_{j=1,...,l_{i}}$ de nombres complexes. On a alors
$$J_{\tilde{M}}^{\tilde{G}}(\tilde{\pi},f)=\sum_{i=1,...,k}(\sum_{j=1,...,l_{j}}z_{i,j})J_{\tilde{M}}^{\tilde{G}}(\tilde{\rho}_{i},f).$$
D'apr\`es 2.5(1), l'hypoth\`ese implique que $\sum_{j=1,...,l_{j}}z_{i,j}=0$ pour tout $i$. D'o\`u la conclusion. $\square$

{\bf Remarque.} Le terme $J_{\tilde{M}}^{\tilde{G}}(\tilde{\sigma},f)$ d\'epend de la mesure sur $G(F)$ (n\'ecessaire pour d\'efinir $Ind_{\tilde{P}}^{\tilde{G}}(\tilde{\sigma},f)$) et de la mesure sur ${\cal A}_{\tilde{M}}^{\tilde{G}}$ (n\'ecessaire pour d\'efinir le terme ${\cal M}_{\tilde{M}}^{\tilde{G}}(\sigma)$). Il ne d\'epend d'aucune autre mesure.

\bigskip

\subsection{$R$-groupes tordus}

Soient $M$ un Levi semi-standard de $G$ et $\sigma$ une repr\'esentation irr\'eductible et de la s\'erie discr\`ete de $M(F)$. 
 On note ${\cal N}^{\tilde{G}}(\sigma)$ l'ensemble des couples $(A,\gamma)$ o\`u $A$ est un automorphisme unitaire de $V_{\sigma}$ et $\gamma\in Norm_{\tilde{G}(F)}(M)$ (c'est-\`a-dire $\gamma\in \tilde{G}(F)$ et $ad_{\gamma}(M)=M$) qui v\'erifient la condition
 $$\sigma( ad_{\gamma}(x))\circ A=\omega(x)A\circ \sigma(x)$$
 pour tout $x\in M(F)$. Le groupe ${\cal N}^G(\sigma)$ agit \`a gauche et \`a droite sur ${\cal N}^{\tilde{G}}(\sigma)$ par 
 $$\begin{array}{ccc}{\cal N}^G(\sigma)\times {\cal N}^{\tilde{G}}(\sigma)\times {\cal N}^G(\sigma)&\to&{\cal N}^{\tilde{G}}(\sigma)\\ ((A',n'),(A,\gamma),(A'',n''))&\mapsto &(A'AA''\omega(n''),n'\gamma n'').\\ \end{array}$$
 L'ensemble ${\cal N}^{\tilde{G}}(\sigma)$ peut \^etre vide. Supposons qu'il ne l'est pas. C'est alors un espace principal homog\`ene sous ${\cal N}^G(\sigma)$, pour l'une ou l'autre des deux actions.
 
 Soit $P\in {\cal P}(M)$, posons $\pi=Ind_{P}^G(\sigma)$. On d\'efinit une application $\tilde{\nabla}_{P}$ de ${\cal N}^{\tilde{G}}(\sigma)$ dans le groupe des op\'erateurs unitaires de $V_{\pi}$ de la fa\c{c}on suivante. Soit $(A,\gamma)\in {\cal N}^{\tilde{G}}(\sigma)$. Alors $\tilde{\nabla}_{P}(A,\gamma)$ est la compos\'ee des op\'erateurs suivants:
 
 - l'op\'erateur $e\mapsto A\circ e$ de $V_{\pi}$ dans $V_{\pi_{1}}$, o\`u $\pi_{1}=Ind_{P}^G(\omega^{-1}(\sigma\circ ad_{\gamma}))$;
 
 - l'op\'erateur $\boldsymbol{\omega}$ de $V_{\pi_{1}}$ dans $V_{\pi_{2}}$, o\`u $\pi_{2}=Ind_{P}^G(\sigma\circ ad_{\gamma})$, qui \`a $e\in V_{\pi_{1}}$, associe la fonction $g\mapsto \omega(g)e(g)$;
 
 - l'op\'erateur $e\mapsto \partial_{P}(\gamma)^{1/2}e\circ ad_{\gamma}^{-1}$ de $V_{\pi_{2}}$ dans $V_{\pi'}$, o\`u $\pi'=Ind_{ad_{\gamma}(P)}^G(\sigma)$ et $\partial_{P}(\gamma)$ est le jacobien de l'application $ad_{\gamma}:U_{P}(F)\to U_{ad_{\gamma}(P)}(F)$;
 
 - l'op\'erateur $R_{P\vert ad{\gamma}(P)}(\sigma):V_{\pi'}\to V_{\pi}$.

 On v\'erifie les relations 
 
 $(1) \qquad \pi(ad_{\gamma}(g))\circ \tilde{\nabla}_{P}(A,\gamma)=\omega(g)\tilde{\nabla}_{P}(A,\gamma)\circ \pi(g)$ pour tout $g\in G(F)$;
 
 $(2) \qquad \pi(n')r_{P}(A',n')\tilde{\nabla}_{P}(A,\gamma)r_{P}(A'',n'')\pi(n'')\omega(n'')=\tilde{\nabla}_{P}(A'AA''\omega(n''),n'\gamma n'')$ pour tous $(A',n'), (A'',n'')\in {\cal N}^G(\sigma)$ et $(A,\gamma)\in {\cal N}^{\tilde{G}}(\sigma)$.
 
 Il s'en d\'eduit que si $(A',n')(A,\gamma)=(A,\gamma)(A'',n'')$, on a l'\'egalit\'e
 $$(3) \qquad r_{P}(A',n')\tilde{\nabla}_{P}(A,\gamma)=\tilde{\nabla}_{P}(A,\gamma)r_{P}(A'',n'').$$
 
 On v\'erifie directement que les orbites dans ${\cal N}^{\tilde{G}}(\sigma)$ pour l'action du sous-groupe $M(F)\subset {\cal N}^G(\sigma)$ sont les m\^emes, que l'on consid\`ere l'action \`a gauche ou l'action \`a droite. Notons  ${\cal W}^{\tilde{G}}(\sigma)$ l'ensemble de ces orbites. Gr\^ace \`a (3), les orbites dans ${\cal W}^{\tilde{G}}(\sigma)$ pour l'action de $W^G_{0}(\sigma)$ sont les m\^emes, que l'on consid\`ere l'action \`a gauche ou l'action \`a droite. Notons ${\cal R}^{\tilde{G}}(\sigma)$ l'ensemble de ces orbites. L'ensemble ${\cal R}^{\tilde{G}}(\sigma)$ est un espace principal homog\`ene sous l'action \`a droite ou \`a gauche de ${\cal R}^G(\sigma)$. Posons
 $$W^{\tilde{G}}(\sigma)=\{\gamma\in \tilde{G}(F); ad_{\gamma}(M)=M, \sigma\circ ad_{\gamma}\simeq \sigma\otimes \omega\}/M(F),$$
 et $R^{\tilde{G}}(\sigma)=W^{\tilde{G}}(\sigma)/W_{0}^G(\sigma)$. On a une application surjective
 $${\cal R}^{\tilde{G}}(\sigma)\to R^{\tilde{G}}(\sigma)$$
 dont les fibres sont isomorphes \`a ${\mathbb U}$.
 
 Soit $\tilde{\lambda}\in i\tilde{{\cal A}}_{\tilde{G},F}^*$. L'application
 $$\begin{array}{ccc}{\cal N}^{\tilde{G}}(\sigma)&\to &{\cal N}^{\tilde{G}}(\sigma_{\lambda})\\ (A,\gamma)&\mapsto&(Ae^{<\tilde{\lambda},\tilde{H}_{\tilde{G}}(\gamma)>},\gamma)\\ \end{array}$$
 est bijective. Il s'en d\'eduit une bijection ${\cal R}^{\tilde{G}}(\sigma)\simeq {\cal R}^{\tilde{G}}(\sigma_{\lambda})$ compatible avec l'identit\'e $R^{\tilde{G}}(\sigma)= R^{\tilde{G}}(\sigma_{\lambda})$. Ces bijections  sont aussi  compatibles avec les isomorphismes ${\cal N}^G(\sigma)\simeq {\cal N}^G(\sigma_{\lambda})$ et ${\cal R}^G(\sigma)\simeq {\cal R}^G(\sigma_{\lambda})$. 
 
 Soit $g\in G(F)$ tel que $M'=gMg^{-1}$ soit semi-standard. Posons $\sigma'=g\sigma=\sigma\circ ad_{g}^{-1}$. On d\'efinit une application
 $$\begin{array}{ccc}{\cal N}^{\tilde{G}}(\sigma)&\to&{\cal N}^{\tilde{G}}(\sigma')\\ (A,\gamma)&\mapsto& (A\omega(g),g\gamma g^{-1}).\\ \end{array}$$
 Elle est bijective et se quotiente en des bijections ${\cal R}^{\tilde{G}}(\sigma)\simeq {\cal R}^{\tilde{G}}(\sigma')$, $R^{\tilde{G}}(\sigma)\simeq R^{\tilde{G}}(\sigma')$. On note toutes ces applications $ad_{g}$. Ce sont ces applications qui servent \`a d\'efinir les notions de conjugaison par $G(F)$ utilis\'ees dans la suite. Par exemple, pour $\boldsymbol{\tilde{r}}\in {\cal R}^{\tilde{G}}(\sigma)$, on dit que les triplets $(M,\sigma,\boldsymbol{\tilde{r}})$ et $(M',\sigma',ad_{g}(\boldsymbol{\tilde{r}}))$ sont conjugu\'es.

 Appelons repr\'esentation de ${\cal R}^{\tilde{G}}(\sigma)$ un couple $(\rho,\tilde{\rho})$, o\`u $\rho$ est une repr\'esentation (disons unitaire et de dimension finie) de ${\cal R}^G(\sigma)$ et $\tilde{\rho}$ est un homomorphisme de ${\cal R}^{\tilde{G}}(\sigma)$ dans le groupe des automorphismes de $V_{\rho}$ v\'erifiant la condition
 $$\rho({\bf r}')\tilde{\rho}(\boldsymbol{\tilde{r}})\rho({\bf r}'')=\tilde{\rho}({\bf r}'\boldsymbol{\tilde{r}}{\bf r}'')$$
 pour tous ${\bf r}',{\bf r}''\in {\cal R}^G(\sigma)$, $\boldsymbol{\tilde{r}}\in {\cal R}^{\tilde{G}}(\sigma)$. On dit que $\tilde{\rho}$ est ${\cal R}^G(\sigma)$-irr\'eductible si $\rho$ est irr\'eductible. Le groupe ${\mathbb U}$ agit par multiplication sur ces repr\'esentations. Par ailleurs, tout \'el\'ement $\boldsymbol{\tilde{r}}\in {\cal R}^{\tilde{G}}(\sigma)$ d\'etermine un automorphisme $\theta_{\boldsymbol{\tilde{r}}}$ de ${\cal R}^G(\sigma)$ par l'\'egalit\'e $\boldsymbol{\tilde{r}}{\bf r}=\theta_{\boldsymbol{\tilde{r}}}({\bf r})\boldsymbol{\tilde{r}}$. Notons $\theta_{{\cal R}}$ sa classe modulo automorphismes int\'erieurs, qui ne d\'epend pas du choix de $\boldsymbol{\tilde{r}}$. Alors les orbites de ${\mathbb U}$ dans l'ensemble des classes de  repr\'esentations ${\cal R}^G(\sigma)$-irr\'eductibles de ${\cal R}^{\tilde{G}}(\sigma)$ sont en bijection avec l'ensemble $Irr({\cal R}^G(\sigma);\theta_{{\cal R}})$ des repr\'esentations irr\'eductibles $\rho$ de ${\cal R}^G(\sigma)$ telles que $\rho\circ\theta_{{\cal R}}\simeq \rho$.
 
 Consid\'erons la d\'ecomposition 1.11(2) de l'espace $V_{\pi}$. Soit $\rho\in Irr({\cal R}^G(\sigma))$. La relation (3) implique que, pour $(A,\gamma)\in {\cal N}^{\tilde{G}}(\sigma)$, $\tilde{\nabla}_{P}(A,\gamma)$ envoie la composante $\rho$-isotypique $\rho\otimes \pi_{\rho}$ sur la composante isotypique $(\rho\circ\theta_{{\cal R}})\otimes \pi_{\rho\circ \theta_{{\cal R}}}$. La relation (1) entra\^{\i}ne alors que $\pi_{\rho\circ\theta_{{\cal R}}}\circ\theta\simeq \omega\pi_{\rho}$. En particulier, $\pi_{\rho}\circ\theta\simeq \omega\pi_{\rho}$ si et seulement si $\rho\in Irr({\cal R}^G(\sigma),\theta_{{\cal R}})$. Supposons cette relation v\'erifi\'ee et choisissons une repr\'esentation $\tilde{\rho}$ de ${\cal R}^{\tilde{G}}(\sigma)$ prolongeant $\rho$. Alors il existe une unique repr\'esentation $\tilde{\pi}_{\tilde{\rho}}$ de $\tilde{G}(F)$  prolongeant $\pi_{\rho}$ de sorte que, pour tout $(A,\gamma)\in {\cal N}^{\tilde{G}}(\sigma)$, la restriction de $\tilde{\nabla}_{P}(A,\gamma)$ \`a $\rho\otimes \pi_{\rho}$ soit \'egale \`a $\tilde{\rho}(\boldsymbol{\tilde{r}})\otimes \tilde{\pi}_{\tilde{\rho}}(\gamma)$, o\`u $\boldsymbol{\tilde{r}}$ est l'image de $(A,\gamma)$ dans ${\cal R}^{\tilde{G}}(\sigma)$. 
 
 Ainsi, \`a un triplet $(M,\sigma,\rho)$, o\`u $M$ est un Levi semi-standard, $\sigma$ est une repr\'esentation irr\'eductible de la s\'erie discr\`ete de $M(F)$ telle que ${\cal N}^{\tilde{G}}(\sigma)\not=\emptyset$ et $\rho\in Irr({\cal R}^G(\sigma);\theta_{{\cal R}})$, on a associ\'e une repr\'esentation $G$-irr\'eductible et temp\'er\'ee $\tilde{\pi}_{\tilde{\rho}}$ de $\tilde{G}(F)$, uniquement d\'efinie \`a multiplication pr\`es par ${\mathbb U}$.  Comme en 1.11, on voit que cette correspondance ne d\'epend pas du sous-groupe parabolique $P$ choisi. Sur l'ensemble des triplets $(M,\sigma,\rho)$, on d\'efinit de fa\c{c}on \'evidente la relation de conjugaison par $G(F)$. La construction ci-dessus se quotiente en une bijection entre l'ensemble des classes de conjugaison par $G(F)$ de triplets $(M,\sigma,\rho)$ et l'ensemble des orbites de ${\mathbb U}$ dans l'ensemble des classes d'isomorphisme de repr\'esentations $G$-irr\'eductibles et temp\'er\'ees de $\tilde{G}(F)$.

 \bigskip
 
 \subsection{L'ensemble $E(\tilde{G},\omega)$}
 Soient $M$ et $\sigma$ comme dans le paragraphe pr\'ec\'edent. On suppose ${\cal N}^{\tilde{G}}(\sigma)\not=\emptyset$.  On fixe $P\in {\cal P}(M)$ et on note $\pi=Ind_{P}^G(\sigma)$. Soit $(A,\gamma)\in {\cal N}^{\tilde{G}}(\sigma)$. Gr\^ace \`a 2.7(1), on peut d\'efinir une $\omega$-repr\'esentation $\tilde{\pi}$ de $\tilde{G}(F)$ par la formule
 $$\tilde{\pi}(g\gamma)=\pi(g)\tilde{\nabla}_{P}(A,\gamma)$$
 pour tout $g\in G(F)$. Notons que cette repr\'esentation n'est pas irr\'eductible en g\'en\'eral.  La relation 2.7(2) montre que $\tilde{\pi}$ ne d\'epend que de l'image $\boldsymbol{\tilde{r}}$ de $(A,\gamma)$ dans ${\cal R}^{\tilde{G}}(\sigma)$. Elle ne d\'epend donc que du triplet $\boldsymbol{\tau}=(M,\sigma,\boldsymbol{\tilde{r}})$ et on peut la noter $\tilde{\pi}_{\boldsymbol{\tau}}$. On voit aussi que sa classe ne d\'epend pas de $P$.    On v\'erifie que la classe de $\tilde{\pi}_{\boldsymbol{\tau}}$ ne d\'epend que de la classe de conjugaison par $G(F)$ du triplet $ \boldsymbol{\tau}$ (cf. 2.8 pour cette notion de conjugaison). Il est tout aussi clair que la correspondance est \'equivariante pour les actions de ${\mathbb U}$: pour $z\in {\mathbb U}$, la repr\'esentation correspondant \`a $(M,\sigma,z\boldsymbol{\tilde{r}})$ est $z\tilde{\pi}_{\boldsymbol{\tau}}$. Le groupe ${\cal R}^G(\sigma)$ agissant \`a gauche et \`a droite sur ${\cal R}^{\tilde{G}}(\sigma)$, il agit aussi par conjugaison. Pour $\boldsymbol{\tilde{r}}\in {\cal R}^{\tilde{G}}(\sigma)$, il peut exister $z\in {\mathbb U}$, $z\not=1$, tel que $z\boldsymbol{\tilde{r}}$ soit conjugu\'e \`a $\boldsymbol{\tilde{r}}$. Dans ce cas,  $z\tilde{\pi}\simeq \tilde{\pi}$ donc le caract\`ere de $\tilde{\pi}$ est nul.
 On dit que le triplet $(M,\sigma,\boldsymbol{\tilde{r}})$ est essentiel si la classe de conjugaison de $\boldsymbol{\tilde{r}}$ par ${\cal R}^G(\sigma)$ ne coupe ${\mathbb U}\boldsymbol{\tilde{r}}$ qu'en le point $\boldsymbol{\tilde{r}}$. On note ${\cal E}(\tilde{G},\omega)$ l'ensemble des triplets $(M,\sigma,\boldsymbol{\tilde{r}})$ qui sont essentiels.
 
 Soit $\boldsymbol{\tau}=(M,\sigma,\boldsymbol{\tilde{r}})$ un triplet comme ci-dessus et soit $\tilde{\lambda}\in i\tilde{{\cal A}}_{\tilde{G},F}^*$. A l'aide de $\tilde{\lambda}$, on a d\'efini en 2.8 un isomorphisme ${\cal R}^{\tilde{G}}(\sigma)\simeq {\cal R}^{\tilde{G}}(\sigma_{\lambda})$. En identifiant ces deux ensembles, le triplet $\boldsymbol{\tau}_{\tilde{\lambda}}=(M,\sigma_{\lambda},\boldsymbol{\tilde{r}})$ v\'erifie encore les conditions requises. On v\'erifie que $\boldsymbol{\tau}$ est essentiel si et seulement si $\boldsymbol{\tau}_{\tilde{\lambda}}$ l'est et qu'on a l'\'egalit\'e $\tilde{\pi}_{\boldsymbol{\tau}_{\tilde{\lambda}}}=(\tilde{\pi}_{\boldsymbol{\tau}})_{\tilde{\lambda}}$.
 
 On note $E(\tilde{G},\omega)$ l'ensemble des triplets $(M,\sigma,\tilde{r})$, o\`u $(M,\sigma)$ est comme ci-dessus, $\tilde{r}\in R^{\tilde{G}}(\sigma)$ et il existe $\boldsymbol{\tilde{r}}\in {\cal R}^{\tilde{G}}(\sigma)$ relevant $\tilde{r}$ de sorte que le triplet $(M,\sigma,\boldsymbol{\tilde{r}})$ soit essentiel. Autrement dit, $E(\tilde{G},\omega)$ est le quotient de ${\cal E}(\tilde{G},\omega)$ par l'action naturelle de ${\mathbb U}$. On note $E(\tilde{G},\omega)/conj$ l'ensemble des classes de conjugaison par $G(F)$ dans l'ensemble $E(\tilde{G},\omega)$. Pour $\boldsymbol{\tau}\in {\cal E}(\tilde{G},\omega)$, l'espace port\'e par le caract\`ere de $\tilde{\pi}_{\boldsymbol{\tau}}$, c'est-\`a-dire $\{0\}$ si ce caract\`ere est nul et une droite sinon, ne d\'epend que de l'image $\tau$ de $\boldsymbol{\tau}$ dans $E(\tilde{G},\omega)$. On le note $D_{\tau}$. Il ne d\'epend aussi que de la classe de conjugaison de $\tau$.  On a d\'efini l'espace $D_{spec}(\tilde{G}(F);\omega)$ en 2.5. On note $D_{temp}(\tilde{G}(F);\omega)$ le sous-espace engendr\'e par les caract\`eres de $\omega$-repr\'esentations temp\'er\'ees.

   \ass{Proposition}{(i) Pour tout $\boldsymbol{\tau}\in {\cal E}(\tilde{G},\omega)$, le caract\`ere de $\tilde{\pi}_{\boldsymbol{\tau}}$ est non nul. 
  
  (ii) On a l'\'egalit\'e
  $$D_{temp}(\tilde{G}(F);\omega)=\oplus_{\tau\in E(\tilde{G},\omega)/conj}D_{\tau}.$$}
  
  Preuve. Pour tout couple $(M,\sigma)$ comme en 2.8 et pour tout $\rho\in Irr({\cal R}^G(\sigma),\theta_{{\cal R}})$, on fixe une extension $\tilde{\rho}$ de $\rho$ \`a ${\cal R}^{\tilde{G}}(\sigma)$. Pour tout $\tilde{r}\in R^{\tilde{G}}(\sigma)$, on fixe aussi un rel\`evement $\boldsymbol{\tilde{r}}\in {\cal R}^{\tilde{G}}(\sigma)$. Pour tout triplet $\beta=(M,\sigma,\rho)$ comme en 2.8, on pose $\tilde{\pi}_{\beta}=\tilde{\pi}_{\tilde{\rho}}$ et, pour tout $\tau=(M,\sigma,\tilde{r})\in E(\tilde{G},\omega)$, on pose $\tilde{\pi}_{\tau}=\tilde{\pi}_{\boldsymbol{\tau}}$, o\`u $\boldsymbol{\tau}=(M,\sigma,\boldsymbol{\tilde{r}})$. On a d\'ecrit les $\omega$-repr\'esentations temp\'er\'ees et $G$-irr\'eductibles en 2.8, \`a l'aide de triplets $\beta=(M,\sigma,\rho)$.  L'espace $D_{temp}(\tilde{G}(F);\omega)$ est somme directe des $D_{\tilde{\pi}_{\beta}}$ quand $\beta$ d\'ecrit les classes de conjugaison de triplets $\beta$. Il nous suffit de prouver l'assertion suivante. Fixons un couple $(M,\sigma)$. Notons $(\tau_{i})_{i=1,...,n}$ les diff\'erents \'el\'ements de $E(\tilde{G},\omega)$ de la forme $(M,\sigma,\tilde{r}_{i})$ et notons $(\beta_{j})_{j=1,...,m}$ les differents triplets comme ci-dessus de la forme $(M,\sigma,\rho_{j})$. Notons $X$ la matrice colonne \`a $n$ lignes dont les coefficients sont les caract\`eres $trace(\tilde{\pi}_{\tau_{i}})$ des $\tilde{\pi}_{\tau_{i}}$ et notons $Y$ la matrice colonne \`a $m$ lignes dont les coefficients sont les caract\`eres $trace(\tilde{\pi}_{\beta_{j}})$ des $\tilde{\pi}_{\beta_{j}}$. Alors
  
  (1)  il existe une matrice \`a coefficients complexes $M$ telle que $X=MY$ et $M$ est inversible.
  
    Soit $(A_{i},\gamma_{i})\in {\cal N}^{\tilde{G}}(\sigma)$ se projetant sur $\boldsymbol{\tilde{r}}_{i}$. On a vu en 2.8 que l'op\'erateur $\tilde{\nabla}_{P}(A_{i},\gamma_{i})$ se restreint \`a une composante $\rho_{j}\otimes \pi_{\rho_{j}}$ en l'op\'erateur $\tilde{\rho}_{j}(\boldsymbol{\tilde{r}}_{i})\otimes \tilde{\pi}_{\beta_{j}}(\gamma)$ et qu'il permute sans point fixe les composantes $\rho\otimes \pi_{\rho}$ pour $\rho\not\in Irr({\cal R}^G(\sigma),\theta_{{\cal R}})$. Il r\'esulte alors de la d\'efinition de $\tilde{\pi}_{\tau_{i}}$ que l'on a l'\'egalit\'e
  $$trace(\tilde{\pi}_{\tau_{i}})=\sum_{j=1,...,m}trace(\tilde{\rho}_{j}(\boldsymbol{\tilde{r}}_{i}))trace(\tilde{\pi}_{\beta_{j}}).$$
  Autrement dit, la matrice $M=(trace(\tilde{\rho}_{j}(\boldsymbol{\tilde{r}}_{i})))_{i=1,...,n;j=1,...,m}$ v\'erifie $X=MY$.  On va montrer que $M$ est inversible. Le groupe ${\cal R}^G(\sigma)$ est une extension de $R^G(\sigma)$ par   ${\mathbb U}$. Le groupe $R^G(\sigma)$ \'etant fini, l'extension provient d'une extension par un sous-groupe fini $Z\subset {\mathbb U}$. Substituons ce groupe $Z$ \`a ${\mathbb U}$ dans toutes les constructions. Fixons $\boldsymbol{\tilde{r}}\in {\cal R}^{\tilde{G}}(\sigma)$. On en d\'eduit un automorphisme $\theta_{\boldsymbol{\tilde{r}}}$ de ${\cal R}^G(\sigma)$. Fixons un entier $N\geq2$ tel que l'ordre de $\theta_{\boldsymbol{\tilde{r}}}$ divise $N$. Introduisons le produit semi-direct $H={\cal R}^G(\sigma)\rtimes({\mathbb Z}/N{\mathbb Z})$ o\`u $k\in {\mathbb Z}/N{\mathbb Z}$ agit par $(\theta_{\boldsymbol{\tilde{r}}})^k$. On identifie ${\cal R}^{\tilde{G}}(\sigma)$ \`a la composante ${\cal R}^G(\sigma)\times\{1\}$ en envoyant ${\bf r}\boldsymbol{\tilde{r}}$ sur $({\bf r},1)$, pour tout ${\bf r}\in {\cal R}^G(\sigma)$.  On est ramen\'e \`a un probl\`eme concernant les repr\'esentations irr\'eductibles du groupe fini $H$. Sa solution est bien connue, indiquons simplement le r\'esultat. Pour tout $\tilde{r}\in R^{\tilde{G}}(\sigma)$, notons $Stab(R^G(\sigma),\tilde{r})$ le stabilisateur de $\tilde{r}$ dans $R^G(\sigma)$ (agissant par conjugaison dans $R^{\tilde{G}}(\sigma)$). Introduisons la matrice 
  $$M'=(\vert Stab(R^G(\sigma),\tilde{r}_{i})\vert ^{-1}\overline{trace(\tilde{\rho}_{j}(\boldsymbol{\tilde{r}}_{i}))})_{i=1,...,n; j=1,...,m}.$$
   Alors la transpos\'ee de $M'$ est l'inverse de $M$. $\square$
    
  Le groupe $i{\cal A}_{\tilde{G},F}^*$ agit naturellement sur $E(\tilde{G},\omega)$: pour $\lambda\in i{\cal A}_{\tilde{G},F}^*$ et pour $\tau=(M,\sigma,\tilde{ r})$, on pose $\tau_{\lambda}=(M,\sigma_{\lambda},\tilde{r})$.   Ainsi, on dispose d'une action de $W^G\times i{\cal A}_{\tilde{G},F}^*$ sur $E(\tilde{G},\omega)$. Pour $\tau\in E(\tilde{G},\omega)$, on note $Stab(W^G\times i{\cal A}_{\tilde{G},F}^*,\tau)$ son stabilisateur dans $W^G\times i{\cal A}_{\tilde{G},F}^*$. C'est un groupe fini qui  contient le stabilisateur $ Stab(i{\cal A}_{\tilde{G},F}^*,\sigma)$ de $\sigma$ dans $i{\cal A}_{\tilde{G},F}^*$. Il contient aussi le groupe $W^M$ comme sous-groupe distingu\'e et son quotient $Stab(W^G\times i{\cal A}_{\tilde{G},F}^*,\tau)/W^M$ contient $W_{0}^G(\sigma)$ comme sous-groupe distingu\'e. On pose
  $${\bf Stab}(W^G\times i{\cal A}_{\tilde{G},F}^*,\tau)=(Stab(W^G\times i{\cal A}_{\tilde{G},F}^*,\tau)/W^M)/W_{0}^G(\sigma).$$
  Remarquons que ${\bf Stab}(W^G\times i{\cal A}_{\tilde{G},F}^*,\tau_{\lambda})={\bf Stab}(W^G\times i{\cal A}_{\tilde{G},F}^*,\tau)$ pour  $\lambda\in i{\cal A}_{\tilde{G},F}^*$ et que $\vert {\bf Stab}(W^G\times i{\cal A}_{\tilde{G},F}^*,\tau)\vert $ ne d\'epend que de l'image de $\tau$ dans $E(\tilde{G},\omega)/conj$.
  
  {\bf Remarque.} Le groupe $Stab(W^G\times i{\cal A}_{\tilde{G},F}^*,\tau)$ peut \^etre plus gros que le produit des stabilisateurs de $\tau$ dans chacun des groupes $W^G$ et $i{\cal A}_{\tilde{G},F}^*$: on peut avoir une relation $w\sigma\simeq \sigma_{\lambda}\not\simeq \sigma$.
  
  \bigskip
  
  \subsection{Triplets et induction}

  Soit $\tilde{L}$ un espace de Levi de $\tilde{G}$ contenant $\tilde{M}_{0}$, soit $M$ un Levi semi-standard de $L$ et soit $\sigma$ une repr\'esentation irr\'eductible et de la s\'erie discr\`ete de $M(F)$. Le groupe ${\cal N}^L(\sigma)$ est contenu dans ${\cal N}^G(\sigma)$ et l'ensemble ${\cal N}^{\tilde{L}}(\sigma)$ est contenu dans l'ensemble ${\cal N}^{\tilde{G}}(\sigma)$. Fixons $\tilde{Q}\in {\cal P}(\tilde{L})$ et $P\in {\cal P}(M)$ tel que $P\subset Q$. Posons $\Pi=Ind_{P}^G(\sigma)$, $\pi=Ind_{P\cap L}^L(\sigma)$.  On sait que $Ind_{Q}^G(\pi)\simeq \Pi$. Dans les mod\`eles naturels, l'isomorphisme envoie une fonction $e$ du premier espace sur la fonction $g\mapsto (e(g))(1)$. Soit $(A,n)\in {\cal N}^L(\sigma)$. On dispose de l'entrelacement $r_{P\cap L}^L(A,n)$ de $\pi$ et de l'entrelacement $r_{P}(A,n)$ de $\Pi$. On a
  
  (1) $r_{P}(A,n)$ s'identifie \`a l'op\'erateur d\'eduit par fonctorialit\'e de $r^L_{P\cap L}(A,n)$.
  
  C'est un simple calcul utilisant le fait que  l'op\'erateur $R_{P\vert  ad_{n}(P)}(\sigma)$ se d\'eduit par fonctorialit\'e de $R^L_{P\vert  ad_{n}(P)}(\sigma)$.
  
  Puisque $W^G_{0}(\sigma)$ est par d\'efinition    le noyau de $(A,n)\mapsto r_{P}(A,n)$  dans ${\cal N}^G(\sigma)/M(F)$ et de m\^eme pour $W^L_{0}(\sigma)$ , on d\'eduit de (1) que $W^L_{0}(\sigma)=W^G_{0}(\sigma)\cap ({\cal N}^L(\sigma)/M(F))$. Alors, d'apr\`es les d\'efinitions,  les plongements ${\cal N}^L(\sigma)\to {\cal N}^G(\sigma)$ et ${\cal N}^{\tilde{L}}(\sigma)\to {\cal N}^{\tilde{G}}(\sigma)$ se quotientent en des plongements ${\cal R}^L(\sigma)\to {\cal R}^G(\sigma)$ et ${\cal R}^{\tilde{L}}(\sigma)\to {\cal R}^{\tilde{G}}(\sigma)$. Evidemment, on a aussi des plongements $R^L(\sigma)\to R^G(\sigma)$ et $R^{\tilde{L}}(\sigma)\to R^{\tilde{G}}(\sigma)$. Supposons ${\cal R}^{\tilde{L}}(\sigma)\not=\emptyset$.
  
  {\bf Remarque.} Cette hypoth\`ese entra\^{\i}ne que la restriction de $\omega$ \`a $Z_{L}(F)^{\theta}$ est triviale.
  
  Soit $\boldsymbol{\tilde{r}}\in {\cal R}^{\tilde{L}}(\sigma)$, posons $\boldsymbol{\tau}=(M,\sigma,\boldsymbol{\tilde{r}})$. On peut consid\'erer ce triplet relativement \`a chacun des espaces ambiants $\tilde{L}$ ou $\tilde{G}$. Remarquons que la notion de triplet essentiel d\'epend de l'espace ambiant. Si le triplet est essentiel relativement \`a $\tilde{G}$, il l'est relativement \`a $\tilde{L}$, mais la r\'eciproque semble fausse en g\'en\'eral (je dis semble car je n'ai pas d'exemple).  Notons qu'associer une $\omega$-repr\'esentation \`a un triplet ne n\'ecessite pas que le triplet soit essentiel (s'il ne l'est pas, le caract\`ere de cette repr\'esentation est nul). Ainsi, on associe \`a $\boldsymbol{\tau}$ une repr\'esentation $\tilde{\pi}_{\boldsymbol{\tau}}$ de $\tilde{L}(F)$, resp. $\tilde{\Pi}_{\boldsymbol{\tau}}$ de $\tilde{G}(F)$.    
  \ass{Lemme}{La repr\'esentation $\tilde{\Pi}_{\boldsymbol{\tau}}$ est isomorphe \`a l'induite $Ind_{\tilde{Q}}^{\tilde{G}}(\tilde{\pi}_{\boldsymbol{\tau}})$.}
  
  Preuve.   Fixons $(A,\gamma)\in {\cal N}^{\tilde{G}}(\sigma)$ se projetant sur $\boldsymbol{\tilde{r}}$. On construit les repr\'esentations $\tilde{\Pi}_{\boldsymbol{\tau}}$ et $\tilde{\pi}_{\boldsymbol{\tau}}$ \`a l'aide de ces \'el\'ements comme en 2.8. Leurs repr\'esentations sous-jacentes de $G(F)$ sont respectivement  $ \Pi$ et $\pi$.  Comme on l'a dit ci-dessus, $\Pi\simeq Ind_{Q}^G(\pi)$.  L'action de $\gamma$ sur $Ind_{\tilde{Q}}^{\tilde{G}}(\tilde{\pi}_{\boldsymbol{\tau}})$ se d\'eduit de $\tilde{\nabla}_{P}^{\tilde{L}}(A,\gamma)$  par la formule 2.5(3). Il suffit de v\'erifier que, par l'isomorphisme pr\'ec\'edent, cette action co\"{\i}ncide avec $\tilde{\Pi}_{\boldsymbol{\tau}}(\gamma)=\tilde{\nabla}_{P}(A,\gamma)$.  Comme pour (1), c'est un simple calcul reposant sur le fait que l'op\'erateur $R_{P\vert  ad_{\gamma}(P)}(\sigma)$ se d\'eduit par fonctorialit\'e de $R^L_{P\vert  ad_{\gamma}(P)}(\sigma)$. $\square$
  
    \bigskip
  
  \subsection{Les ensembles $E_{disc}(\tilde{G},\omega)$ et $E_{ell}(\tilde{G},\omega)$}
  Soient $M$ un Levi semi-standard de $G$ et  $\sigma$ une repr\'esentation irr\'eductible et de la s\'erie discr\`ete de $M(F)$. Un \'el\'ement $\tilde{w}\in W^{\tilde{G}}(\sigma)$ agit naturellement sur ${\cal A}_{M}$. Notons ${\cal A}_{M}^{\tilde{w}}$ le sous-espace des points fixes par cette action. Il contient  ${\cal A}_{\tilde{G}}$. On note $W^{\tilde{G}}_{reg}(\sigma)$ l'ensemble des  $\tilde{w}\in W^{\tilde{G}}(\sigma)$ tels que ${\cal A}_{M}^{\tilde{w}}={\cal A}_{\tilde{G}}$. Rappelons que $R^{\tilde{G}}(\sigma)=W^G_{0}(\sigma)\backslash W^{\tilde{G}}(\sigma)$, en particulier $R^{\tilde{G}}(\sigma)=W^{\tilde{G}}(\sigma)$ si $W^G_{0}(\sigma)=\{1\}$.
  
  \ass{Lemme}{Soit $\tilde{r}\in R^{\tilde{G}}(\sigma)$. Les conditions suivantes sont \'equivalentes:
  
  (a) il n'existe pas d'espace de Levi  $\tilde{L}$ tel que $M\subset L\subsetneq G$ et $\tilde{r}\in R^{\tilde{L}}(\sigma)$;
  
  (b) $W^G_{0}(\sigma)=\{1\}$ et $\tilde{r}\in W^{\tilde{G}}_{reg}(\sigma)$.}
  
  Preuve. Fixons $\tilde{w}\in W^{\tilde{G}}(\sigma)$ d'image $\tilde{r}$. Pour un espace de Levi $\tilde{L}$ tel que $M\subset L\subsetneq G$, l'image de $R^{\tilde{L}}(\sigma)$ dans $R^{\tilde{G}}(\sigma)$ est $W^G_{0}(\sigma)\backslash(W^G_{0}(\sigma)W^{\tilde{L}}(\sigma))$. La condition (a) \'equivaut donc \`a
  
  (1) pour tout $w'\in W^{G}_{0}(\sigma)$, il n'existe pas de $\tilde{L}$ comme ci-dessus tel que $w'\tilde{w}\in W^{\tilde{L}}(\sigma)$.
  
  D'apr\`es 2.1(4), la condition $w'\tilde{w}\in W^{\tilde{L}}(\sigma)$ \'equivaut \`a ${\cal A}_{\tilde{L}}\subset{\cal A}_{M}^{w'\tilde{w}}$. Donc (1) est \'equivalent \`a
  
  (2) pour tout $w'\in W^{G}_{0}(\sigma)$, $w'\tilde{w}\in W^{\tilde{G}}_{reg}(\sigma)$.
  
  Si $W^G_{0}(\sigma)=\{1\}$, (2) est \'equivalente \`a (b). Si $W^G_{0}(\sigma)\not=\{1\}$, (b) n'est pas v\'erifi\'ee et on doit prouver que (2) ne l'est pas non plus. Comme on l'a dit en 1.11, $W^G_{0}(\sigma)$ est le groupe de Weyl d'un syst\`eme de racines inclus dans l'ensemble des racines de $A_{M}$ dans $G$. Il r\'esulte des d\'efinitions que l'action de $\tilde{w}$ conserve ce syst\`eme de racines. Fixons une base de ce syst\`eme. On peut alors trouver $w'\in W^G_{0}(\sigma)$ tel que $w'\tilde{w}$ conserve cette base (en permutant ses \'el\'ements). La somme des coracines associ\'ees aux \'el\'ements de cette base est alors un \'el\'ement de ${\cal A}_{M}^{w'\tilde{w}}$ qui n'appartient pas \`a ${\cal A}_{\tilde{G}}$. Donc $w'\tilde{w}\not\in W^{\tilde{G}}_{reg}(\sigma)$ et (2) n'est pas v\'erifi\'ee. $\square$
  
Soit   $\tilde{w}\in W^{\tilde{G}}(\sigma)$. Comme on vient de le dire, l'action de $\tilde{w}$ conserve l'ensemble de racines associ\'e \`a $W_{0}^G(\sigma)$, cf.1.11. En fixant un sous-ensemble positif, on pose $\epsilon_{\sigma}(\tilde{w})=(-1)^{n(\tilde{w})}$, o\`u $n(\tilde{w})$ est le nombre de racines positives $\alpha$ telles que $\tilde{w}(\alpha)$ soit n\'egative. Ce signe ne d\'epend pas de l'ensemble positif choisi.    Soit $\tau=(M,\sigma,\tilde{r})\in E(\tilde{G},\omega)$. On choisit $\tilde{w}\in {\cal W}^{\tilde{G}}(\sigma)$ se projetant sur $\tilde{r}$ et on pose
 $$\iota(\tau)=\vert W^G_{0}(\sigma)\vert ^{-1}\sum_{\tilde{w}'\in W^G_{0}(\sigma)\tilde{w}\cap W^{\tilde{G}}_{reg}(\sigma)}\epsilon_{\sigma}(\tilde{w}')\vert det((1-\tilde{w}')_{\vert {\cal A}_{M}^{\tilde{G}}})\vert ^{-1}.$$

 On dit que le triplet $\tau\in E(\tilde{G},\omega) $ est discret si $W^G_{0}(\sigma)\tilde{w}\cap W^{\tilde{G}}_{reg}(\sigma)\not=\emptyset$.   On dit que $\tau$ est elliptique s'il est discret et de plus $W^G_{0}(\sigma)=\{1\}$.  On note   $E_{disc}(\tilde{G},\omega)$, resp. $E_{ell}(\tilde{G},\omega)$, l'ensemble des triplets $(M,\sigma,r)$ qui sont  discrets, resp. elliptiques. On note   $E_{disc}(\tilde{G},\omega)/conj$, resp. $E_{ell}(\tilde{G},\omega)/conj$, l'ensemble des classes de conjugaison par $G(F)$ dans   $E_{disc}(\tilde{G},\omega)$, resp. $E_{ell}(\tilde{G},\omega)$.  
 
 \bigskip
 
 \subsection{Repr\'esentations elliptiques}
 Posons, avec les notations de 2.9,
 $$D_{ell}(\tilde{G}(F);\omega)=\oplus_{\tau\in E_{ell}(\tilde{G},\omega)/conj}D_{\tau}.$$
 Pour tout ensemble de Levi $\tilde{L}\in {\cal L}(\tilde{M}_{0})$, posons $W^G(\tilde{L})=Norm_{G(F)}(\tilde{L})/L(F)$.   Supposons que $\omega$ soit trivial sur $Z_{L}(F)^{\theta}$.  Le groupe $Norm_{G(F)}(\tilde{L})$ agit sur $D_{temp}(\tilde{L}(F);\omega)$: \`a $n\in Norm_{G(F)}(\tilde{L})$ et $d\in D_{temp}(\tilde{L}(F);\omega)$, on associe la distribution $f\mapsto \omega(n)^{-1}d(f^n)$, o\`u $f^n(\gamma)=f(n\gamma n^{-1})$.
 Cette action se descend en une action de $W^G(\tilde{L})$  sur $D_{temp}(\tilde{L}(F);\omega)$. Cette action conserve le sous-espace $D_{ell}(\tilde{L}(F);\omega)$. Selon l'usage, on note $D_{ell}(\tilde{L}(F);\omega)^{W^G(\tilde{L})}$ le sous-espace des invariants. 

D\'efinissons une action du groupe $Norm_{G(F)}(\tilde{L})$   sur l'ensemble ${\cal E}(\tilde{L})/conj$ des classes de conjugaison par $L(F)$ dans ${\cal E}(\tilde{L})$. Soient $n\in Norm_{G(F)}(\tilde{L})$ et $\boldsymbol{\tau}=(M,\sigma,\boldsymbol{\tilde{r}})\in {\cal E}(\tilde{L})$. Posons $M'=nMn^{-1}$, $\sigma'=n\sigma$. Quitte \`a multiplier $n$ \`a gauche par un \'el\'ement de $L(F)$, on peut supposer $M'$ semi-standard. On d\'efinit une bijection $ad_{n}:{\cal N}^{\tilde{L}}(\sigma)\simeq {\cal N}^{\tilde{L}}(\sigma')$ comme en 2.8: \`a $(A,\gamma)$, on associe $(A\omega(n),n\gamma n^{-1})$. Cette bijection se descend en une bijection $ad_{n}:{\cal R}^{\tilde{G}}(\sigma)\to {\cal R}^{\tilde{G}}(\sigma')$. La classe de conjugaison par $L(F)$ du triplet $\boldsymbol{\tau}'=(M',\sigma',ad_{n}(\boldsymbol{\tilde{r}}))$ ne d\'epend pas de la modification  de $n$ faite ci-dessus et ne d\'epend que de la classe de conjugaison par $L(F)$ de $\boldsymbol{\tau}$. L'action cherch\'ee est celle qui associe \`a la classe de $\boldsymbol{\tau}$ celle de $\boldsymbol{\tau}'$. Cette action se quotiente en une action de $W^G(\tilde{L})$.

Notons ${\cal L}(\tilde{M}_{0};\omega)$ l'ensemble des $\tilde{L}\in {\cal L}(\tilde{M}_{0})$ tels que 
$\omega$ soit trivial sur $Z_{L}(F)^{\theta}$. Soit $\tilde{L}\in {\cal L}(\tilde{M}_{0};\omega)$. Pour $\tilde{Q}\in {\cal P}(\tilde{M})$, l'op\'eration d'induction $Ind_{\tilde{Q}}^{\tilde{G}}$ d\'efinit une application de $D_{temp}(\tilde{L}(F);\omega)$ dans $ D_{temp}(\tilde{G}(F);\omega)$ qui ne d\'epend pas du  choix de $\tilde{Q}$. On note
$$ Ind_{\tilde{L}}^{\tilde{G}}:D_{temp}(\tilde{L}(F);\omega)\to D_{temp}(\tilde{G}(F);\omega)$$
cette application. 
Posons simplement $\tilde{W}^G=W^G(\tilde{M}_{0})$. Ce groupe agit sur ${\cal L}(\tilde{M}_{0})$ et conserve ${\cal L}(\tilde{M}_{0},\omega)$.

\ass{Proposition}{ (i) Pour tout $\tilde{L}\in  {\cal L}(\tilde{M}_{0},\omega)$, l'application $Ind_{\tilde{L}}^{\tilde{G}}$ est injective sur $D_{ell}(\tilde{L}(F);\omega)^{W^G(\tilde{L})}$.

(ii)On a l'\'egalit\'e
$$D_{temp}(\tilde{G}(F),\omega)=\oplus_{\tilde{L}\in  {\cal L}(\tilde{M}_{0},\omega)/\tilde{W}^G}Ind_{\tilde{L}}^{\tilde{G}}(D_{ell}(\tilde{L}(F);\omega)^{W^G(\tilde{L})}).$$}

Preuve. D'apr\`es la proposition 2.9, l'espace $D_{temp}(\tilde{G}(F);\omega)$ a une base param\'etr\'ee par $E(\tilde{G},\omega)/conj$. On note ici $(e^{\tilde{G}}(\tau))_{\tau\in E(\tilde{G},\omega)/conj}$ une telle base. Alors $(e^{\tilde{G}}(\tau))_{\tau\in E_{ell}(\tilde{G},\omega)/conj}$ est une base de $D_{ell}(\tilde{G}(F);\omega)$. Soit $\tilde{L}\in  {\cal L}(\tilde{M}_{0},\omega)$. Le groupe $W^G(\tilde{L})$ agit tant sur $D_{ell}(\tilde{L}(F);\omega)$ que sur $E_{ell}(\tilde{L})/conj$. Pour $w\in W^G(\tilde{L})$ et $\tau\in E_{ell}(\tilde{L})/conj$, on a une \'egalit\'e $w(e^{\tilde{L}}(\tau))=z(w,\tau)e^{\tilde{L}}(w\tau)$, avec $z(w,\tau)\in {\mathbb C}^{\times}$.  Posons $\underline{e}^{\tilde{L}}(\tau)=\sum_{w\in W^G(\tilde{L})}w(e^{\tilde{L}}(\tau))$. Cet \'el\'ement peut \^etre nul: la restriction de l'application $w\mapsto z(w,\tau)$ au stabilisateur de $\tau$ est un caract\`ere de ce groupe; $\underline{e}^{\tilde{L}}(\tau)$ est nul si et seulement si ce caract\`ere est non trivial. Notons $\underline{E}_{ell}(\tilde{L})$ l'ensemble des $\tau\in E_{ell}(\tilde{L})$ dont la classe v\'erifie $\underline{e}^{\tilde{L}}(\tau)\not=0$. L'action de $W^G(\tilde{L})$ respecte cet ensemble. Identifions l'ensemble quotient $(\underline{E}_{ell}(\tilde{L})/conj)/W^G(\tilde{L})$ \`a un ensemble de repr\'esentants (on fera d'autres identifications similaires dans la suite).  Alors $(\underline{e}^{\tilde{L}}(\tau))_{\tau\in (\underline{E}_{ell}(\tilde{L})/conj)/W^G(\tilde{L})}$ est une base de $D_{ell}(\tilde{L}(F);\omega)^{W^G(\tilde{L})}$.

 Pour ce qui est des triplets $(M,\sigma,\tilde{r})$, il y a une correspondance naturelle entre triplets pour $\tilde{L}$ et triplets pour $\tilde{G}$: c'est l'identit\'e, modulo l'identification de   ${\cal R}^{\tilde{L}}(\sigma)$ \`a un sous-ensemble de ${\cal R}^{\tilde{G}}(\sigma)$. On a d\'ej\`a dit qu'un triplet pouvait \^etre essentiel pour $\tilde{L}$ mais pas pour $\tilde{G}$. Notons $E^{\tilde{G}}(\tilde{L})$ l'ensemble  des triplets pour $\tilde{L}$ qui sont essentiels dans $\tilde{G}$. De la correspondance pr\'ec\'edente se d\'eduit une application $\iota_{\tilde{L}}^{\tilde{G}}:E^{\tilde{G}}(\tilde{L})/conj\to E(\tilde{G},\omega)/conj$. Celle-ci est invariante par l'action de $W^G(\tilde{L})$ sur l'espace de d\'epart.  D'apr\`es le lemme 2.10, l'application    $Ind_{\tilde{L}}^{\tilde{G}}$ se d\'ecrit ainsi: pour $\tau\in E(\tilde{L},\omega)/conj$, elle envoie $e^{\tilde{L}}(\tau)$ sur $0$ si $\tau\not\in E^{\tilde{G}}(\tilde{L})/conj$ et sur $z_{\tilde{L}}^{\tilde{G}}(\tau)e^{\tilde{G}}(\iota_{\tilde{L}}^{\tilde{G}}(\tau))$ si $\tau\in E^{\tilde{G}}(\tilde{L})/conj$, o\`u $z_{\tilde{L}}^{\tilde{G}}(\tau)\in {\mathbb C}^{\times}$. Soit $\tau\in E^{\tilde{G}}_{ell}(\tilde{L})/conj$ ($=(E^{\tilde{G}}(\tilde{L})/conj)\cap (E_{ell}(\tilde{L})/conj)$). Puisque l'induction est insensible \`a l'action de $W^G(\tilde{L})$, $Ind_{\tilde{L}}^{\tilde{G}}$ envoie $\underline{e}^{\tilde{L}}(\tau)$ sur $\vert W^G(\tilde{L})\vert z_{\tilde{L}}^{\tilde{G}}(\tau) e^{\tilde{G}}(\iota_{\tilde{L}}^{\tilde{G}}(\tau))$. A fortiori $\underline{e}^{\tilde{L}}(\tau)\not=0$, ce qui d\'emontre l'inclusion $E^{\tilde{G}}_{ell}(\tilde{L})\subset \underline{E}_{ell}(\tilde{L})$.
 
 On a d\'ecrit des bases de chacun de nos espaces et les matrices des applications d'induction dans ces bases. Le lemme r\'esulte alors des deux assertions suivantes:
 
 (1) pour tout $\tilde{L}\in  {\cal L}(\tilde{M}_{0},\omega)$, on a l'\'egalit\'e $E^{\tilde{G}}_{ell}(\tilde{L})= \underline{E}_{ell}(\tilde{L})$;
 
 (2) notons 
 $$\iota:\bigsqcup_{\tilde{L}\in  {\cal L}(\tilde{M}_{0},\omega)/\tilde{W}^G}(E^{\tilde{G}}_{ell}(\tilde{L})/conj)/W^G(\tilde{L})\to E(\tilde{G},\omega)/conj$$
  l'application qui co\"{\i}ncide avec  $\iota_{\tilde{L}}^{\tilde{G}}$ sur le sous-ensemble index\'e par $\tilde{L}$ de l'ensemble de d\'epart; alors $\iota$ est bijective.

 Prouvons (2). Soit $\tau=(M,\sigma,\tilde{r})\in E(\tilde{G},\omega)$.  Parmi les espaces de Levi $\tilde{L}$ tels que $M\subset L$ et $\tilde{r}\in {\cal R}^{\tilde{L}}(\sigma)$, consid\'erons un \'el\'ement minimal $\tilde{L}$. Quitte \`a conjuguer $\tau$ par un \'el\'ement de $G(F)$, on peut supposer $\tilde{L}\in  {\cal L}(\tilde{M}_{0},\omega)$. Notons $\tau^{\tilde{L}}$ le m\^eme triplet vu comme un \'el\'ement de $E(\tilde{L},\omega)$. On a $\tau=\iota_{\tilde{L}}^{\tilde{G}}(\tau^{\tilde{L}})$. Le lemme 2.11 et la minimalit\'e de $\tilde{L}$ assurent que $\tau^{\tilde{L}}\in E_{ell}(\tilde{L})$ (et forc\'ement $\tau^{\tilde{L}}\in E^{\tilde{G}}_{ell}(\tilde{L})$).  Cela prouve la surjectivit\'e de $\iota$. D\'emontrons l'injectivit\'e. On identifie $ {\cal L}(\tilde{M}_{0},\omega)/\tilde{W}^G$ \`a un ensemble de repr\'esentants dans $ {\cal L}(\tilde{M}_{0},\omega)$.  Soient $\tilde{L},\tilde{L}' $ dans cet ensemble, $\tau=(M,\sigma,\tilde{r})\in E^{\tilde{G}}_{ell}(\tilde{L})$, $\tau'=(M',\sigma',\tilde{r}')\in E^{\tilde{G}}_{ell}(\tilde{L}')$, supposons les triplets $(M,\sigma,\tilde{r})$ et $(M',\sigma',\tilde{r}')$ conjugu\'es par un \'el\'ement de $G(F)$. On doit prouver qu'alors $\tilde{L}=\tilde{L}'$ et que $\tau$ et $\tau'$ sont conjugu\'es par l'action du normalisateur $Norm_{G(F)}(\tilde{L})$. Soit $g\in G(F)$ qui conjugue le premier triplet en le second. On a $gMg^{-1}=M'$, $g\sigma=\sigma'$. Notons que, puisque $\tau$ est elliptique, $\tilde{r}$ est au d\'epart un \'el\'ement de $W^{\tilde{L}}_{reg}(\sigma)$. Mais, apr\`es passage \`a $\tilde{G}$,  le triplet n'est plus elliptique et $\tilde{r}$ devient une classe modulo $W^{G}_{0}(\sigma)$. De m\^eme pour $\tilde{r}'$. En continuant \`a consid\'erer $\tilde{r}$ et $\tilde{r}'$ comme des \'el\'ements de $W^{\tilde{L}}_{reg}(\sigma)$ et $W^{\tilde{L}'}_{reg}(\sigma')$,  la condition de conjugaison est $ad_{g}(W^G_{0}(\sigma)\tilde{r})=W_{0}^G(\sigma')\tilde{r}'$. On a dit que $W^G_{0}(\sigma)$ est le groupe de Weyl d'un sous-syst\`eme de racines $\Sigma$ de l'ensemble des racines de $A_{M}$ dans $G$. Puisque $W^L_{0}(\sigma)=\{1\}$ (condition d'ellipticit\'e), aucune de ces racines n'intervient dans $L$.  On peut donc trouver un \'el\'ement $H\in {\cal A}_{\tilde{L}}$ qui n'annule aucune de ces racines. Fixons un tel \'el\'ement. Il d\'etermine un sous-ensemble   positif $\Sigma_{+}\subset\Sigma$: l'ensemble des $\alpha\in \Sigma$ tels que $<\alpha,H>>0$. Puisque $\tilde{r}\in W^{\tilde{L}}$ fixe $H$, l'action de $\tilde{r}$ conserve $\Sigma_{+}$. On a  des objets analogues pour $\tau'$, que l'on affecte d'un $'$. La conjugaison par $g$ envoie $\Sigma$ sur $\Sigma'$. Quitte \`a multiplier $g$ \`a gauche par un \'el\'ement de $W^G_{0}(\sigma')$, on peut supposer qu'elle envoie $\Sigma_{+} $ sur $\Sigma'_{+}$. Ecrivons $ad_{g}(\tilde{r})=w\tilde{r}'$, avec $w\in W_{0}^G(\sigma)$. Puisque $\tilde{r}$ conserve $\Sigma_{+}$, $ad_{g}(\tilde{r})$ conserve $\Sigma'_{+}$. Il en est de m\^eme de $\tilde{r}'$. Donc aussi de $w$. Un \'el\'ement du groupe de Weyl qui conserve un ensemble positif est l'identit\'e. D'o\`u $w=1$. Puisque  $\tilde{r}\in W^{\tilde{L}}(\sigma)$, $\tilde{L}$ est le plus petit espace de Levi contenant \`a la fois $M$ et $\tilde{r}$. De m\^eme pour $\tilde{L}'$. Alors $g$ conjugue $\tilde{L}$ en $\tilde{L}'$. Deux \'el\'ements de ${\cal L}(\tilde{M}_{0})$ qui sont conjugu\'es par un \'el\'ement de $G(F)$ le sont par un \'el\'ement de $\tilde{W}^G$. Puisque les deux espaces de Levi  $\tilde{L}$ et $\tilde{L}'$ appartiennent \`a notre ensemble de repr\'esentants, on a $\tilde{L}=\tilde{L}'$ et $g$ appartient \`a $Norm_{G(F)}(\tilde{L})$.  C'est ce que l'on voulait prouver.
 
 Prouvons (1). Soient $\tilde{L}\in  {\cal L}(\tilde{M}_{0},\omega)$ et $\tau\in E_{ell}(\tilde{L})-E^{\tilde{G}}_{ell}(\tilde{L})$. On veut prouver que $\tau\not\in \underline{E}_{ell}(\tilde{L})$, autrement dit que la fonction $w\mapsto z(w,\tau)$ n'est pas constante sur le stabilisateur $Stab(W^G(\tilde{L}),\tau)$ de $\tau$ dans $W^G(\tilde{L})$ (en notant encore $\tau$ la classe de conjugaison de $\tau$ par $L(F)$). Relevons $\tau$ en un \'el\'ement $\boldsymbol{\tau}=(M,\sigma,\boldsymbol{\tilde{r}})\in {\cal E}(\tilde{L})$. Il r\'esulte des d\'efinitions que, pour $w\in Stab(W^G(\tilde{L}),\tau)$, $w(\boldsymbol{\tau})$ est \'egal   \`a  $(M,\sigma,z(w,\tau)\boldsymbol{\tilde{r}})$, \`a conjugaison pr\`es par un \'el\'ement de $L(F)$. Il s'agit donc de trouver un \'el\'ement $g\in Norm_{G(F)}(\tilde{L})$ tel que $ad_{g}(\boldsymbol{\tau})=(M,\sigma,z\boldsymbol{\tilde{r}})$, avec $z\not=1$. Puisque $\tau\not\in E^{\tilde{G}}_{ell}(\tilde{L})$, il existe en tout cas un \'el\'ement $g\in G(F)$ qui v\'erifie cette derni\`ere relation. La preuve de l'injectivit\'e de $\iota$ montre que, quitte \`a modifier $g$ par un \'el\'ement de $W^G_{0}(\sigma)$ (ce qui ne change par $ad_{g}(\boldsymbol{\tau})$), on peut supposer $g\in Norm_{G(F)}(\tilde{L})$. C'est ce qu'on voulait. $\square$
 
   \bigskip
 
 \section{Le calcul spectral}
 
 \bigskip
 
 \subsection{Position du probl\`eme}
   On consid\`ere un \'el\'ement $T\in {\cal A}_{0}$, qui intervient dans ce qui suit comme un param\`etre. On lui impose de v\'erifier les propri\'et\'es suivantes:
 
 $\theta(T)=T$;
 
 $<\alpha,T>>0$ pour tout $\alpha\in \Delta_{0}$;
 
 $<\alpha,T>\geq c_{\star}\vert T\vert $, o\`u $c_{\star}>0$ est un r\'eel fix\'e;
 
  si $F$ est non-archim\'edien, $T\in  {\cal A}_{M_0,F}\otimes_{\mathbb Z}{\mathbb Q}$.

 On note ${\tilde{\kappa}}^T$ la fonction caract\'eristique du sous-ensemble des $g\in G(F)$ tels que $\phi^{\tilde{G}}(h_{0}(g)-T)=1$. Ce sous-ensemble est invariant par $A_{G}(F)$ et sa projection dans $A_{G}(F)\backslash G(F)$ est compact.
 
 On d\'efinit le sous-espace $C_{c}^{\infty}(\tilde{G}(F);K)$ de $C_{c}^{\infty}(\tilde{G}(F))$ comme en 2.7.  
Soient $f_{1},f_{2}\in C_{c}^{\infty}(\tilde{G}(F),K)$.  On pose
$$J^T(\omega,f_{1},f_{2})=\int_{A_{\tilde{G}}(F)\backslash G(F)}\int_{\tilde{G}(F)} \bar{f}_{1}(\gamma)f_{2}(g^{-1}\gamma g) \omega(g){\tilde{\kappa}}^T(g)\,d\gamma\,dg.$$
C'est une int\'egrale \`a support compact. En effet, l'int\'egration en $\gamma$ est \`a support compact puisque $f_{1}$ l'est.  A cause de la fonction ${\tilde{\kappa}}^T$, on peut \'ecrire $g=ay$, o\`u $y$ reste dans un compact et $a\in A_{G}(F)$. La condition $f_{2}(g^{-1}\gamma g)\not=0$ impose alors que $a^{-1}\theta(a)$ reste dans un compact ce qui entra\^{\i}ne que l'image de $a$ dans $A_{\tilde{G}}(F)\backslash A_{G}(F)$ reste dans un compact.

On se propose dans cette section de montrer que la fonction $T\mapsto J^T(\omega,f_{1},f_{2})$ est asymptote \`a un \'el\'ement de $PolExp$ et de calculer une expression "spectrale" du terme constant  de cet \'el\'ement.

\bigskip

\subsection{Utilisation de la formule de Plancherel}
On fixe un \'el\'ement $\gamma_{0}\in \tilde{M}_{0}(F)$, v\'erifiant la condition 2.1(6) quand $F$ est archim\'edien. On note simplement $\theta=ad_{\gamma_{0}}$.

On d\'efinit les fonctions $\varphi_{1}$ et $\varphi_{2}$ sur $G(F)$ par $\varphi_{i}(x)=f_{i}(x\gamma_{0})$ pour $i=1,2$. Alors
$$J^T(\omega,f_{1},f_{2})=\int_{A_{\tilde{G}}(F)\backslash G(F)}\int_{G(F)} \bar{\varphi}_{1}(x)\varphi_{2}(g^{-1}x\theta(g)) \omega(g)\tilde{\kappa}^T(g)\,dx\,dg.$$
On utilise la formule de Plancherel-Harish-Chandra pour exprimer $\varphi_{2}$. C'est-\`a-dire
$$\varphi_{2}(g)=\sum_{M_{disc}\in {\cal L}(M_{0})}\vert W^{M_{disc}}\vert \vert W^G\vert ^{-1}\sum_{\sigma\in  \Pi_{disc}(M_{disc}(F))/i{\cal A}_{M_{disc},F}^*}  \vert Stab(i{\cal A}_{M_{disc},F}^*,\sigma)\vert ^{-1}$$
$$\int_{i{\cal A}_{M_{disc},F}^*}m^G(\sigma_{\lambda})trace(Ind_{S}^G(\sigma_{\lambda},g^{-1})Ind_{S}^G(\sigma_{\lambda},\varphi_{2}))\,d\lambda.$$

Pour tous $M$, $\sigma$, l'expression
$$\int_{A_{\tilde{G}}(F)\backslash G(F)}\int_{G(F)} \vert \varphi_{1}(x)\vert $$
$$\vert \int_{i{\cal A}_{M_{disc},F}^*}m^G(\sigma_{\lambda})trace(Ind_{S}^G(\sigma_{\lambda},\theta(^{-1}x^{-1}g)Ind_{S}^G(\sigma_{\lambda},\varphi_{2}))\,d\lambda\vert {\tilde{\kappa}}^T(g)\,dx\,dg$$
est convergente. En effet, l'int\'egrale en $x$ est \`a support compact puisque $\varphi_{1}$ l'est. L'int\'egrale int\'erieure d\'efinit une fonction de Schwartz-Harish-Chandra en $\theta(g)^{-1}x^{-1}g$.   Comme dans le paragraphe pr\'ec\'edent, on peut \'ecrire $g=ay$ o\`u $y$ reste dans un compact. L'int\'egrale restante en $a$ est celle d'une fonction essentiellement born\'ee par $(1+\vert \theta(H_0(a))-H_0(a)\vert )^{-r}$ pour tout r\'eel $r$. Une telle int\'egrale est convergente. On pose
$$J^T_{M_{disc},\sigma}(\omega,f_{1},f_{2})=  \int_{A_{\tilde{G}}(F)\backslash G(F)}\int_{G(F)}  \bar{\varphi}_{1}(x) $$
$$\int_{i{\cal A}_{M_{disc},F}^*}m^G(\sigma_{\lambda})trace(Ind_{S}^G(\sigma_{\lambda},\theta(g)^{-1}x^{-1}g)Ind_{S}^G(\sigma_{\lambda},\varphi_{2}))\,d\lambda \omega(g){\tilde{\kappa}}^T(g)\,dx\,dg$$
et on a
$$(1) \qquad J^T(\omega,f_{1},f_{2})=\sum_{M_{disc}\in {\cal L}(M_{0})}\vert W^{M_{disc}}\vert \vert W^G\vert ^{-1}$$
$$\sum_{\sigma\in  \Pi_{disc}(M_{disc}(F))/i{\cal A}_{M_{disc},F}^*}\vert Stab(i{\cal A}_{M_{disc},F}^*,\sigma)\vert ^{-1}J^T_{M_{disc},\sigma}(\omega,f_{1},f_{2}).$$

\bigskip

\subsection{Apparition d'int\'egrales de coefficients}
  Nous fixons maintenant un Levi semi-standard $M_{disc}$ et $\sigma\in \Pi_{disc}(M_{disc}(F))$. On fixe aussi $S\in {\cal P}(M_{disc})$. On pose  $\pi_{\lambda}=Ind_{S}^G(\sigma_{\lambda})$ que l'on r\'ealise dans l'espace  $V_{\sigma,S}$. On pose $S'=\theta^{-1}(S)$, $\sigma'=\sigma\circ \theta$, $\lambda'=\theta^{-1}\lambda$, $\pi'_{\lambda'}=Ind_{S'}^G(\sigma'_{\lambda'})$, que l'on r\'ealise dans l'espace $V_{\sigma',S'}$. 
  On fixe des bases orthonorm\'ees ${\cal B}$ de $V_{\sigma,S}$ et ${\cal B}'$ de $V_{\sigma',S'}$,  r\'eunions de bases des diff\'erents $K$-types intervenant. On introduit l'op\'erateur $V_{\pi'_{\lambda'}}\to V_{\pi_{\lambda}}$ qui \`a $e\in V_{\pi'_{\lambda'}}$ associe la fonction $g\mapsto e(\theta^{-1}(g))$.  Par restriction \`a $K$, on obtient un op\'erateur unitaire $U_{\theta,\sigma_{\lambda}}:V_{\sigma',S'}\to V_{\sigma,S}$, qui v\'erifie $U_{\theta,\sigma_{\lambda}}\pi'_{\lambda'}(g)=\pi_{\lambda}(\theta(g))U_{\theta,\sigma_{\lambda}}$. Si $K$ est stable par $\theta$, il est ind\'ependant de $\lambda$. Mais on n'a pas pos\'e cette hypoth\`ese sur $K$ et l'op\'erateur peut d\'ependre de $\lambda$.

   Pour $x,g\in G(F)$, on a
$$trace(Ind_{S}^G(\sigma_{\lambda},\theta(g)^{-1}x^{-1}g)Ind_{S}^G(\sigma_{\lambda},\varphi_{2}))=trace(U_{\theta,\sigma_{\lambda}}^{-1}Ind_{S}^G(\sigma_{\lambda},\theta(g)^{-1}x^{-1}g)Ind_{S}^G(\sigma_{\lambda},\varphi_{2})U_{\theta,\sigma_{\lambda}})$$
$$=\sum_{v'\in{\cal B}'}(v',U_{\theta,\sigma_{\lambda}}^{-1}\pi_{\lambda}(\theta(g)^{-1}x^{-1}g)\pi_{\lambda}(\varphi_{2})U_{\theta,\sigma_{\lambda}}v')$$
$$=\sum_{v'\in {\cal B}'}(\pi_{\lambda}(x\theta(g))U_{\theta,\sigma_{\lambda}}v',\pi_{\lambda}(g)\pi_{\lambda}(\varphi_{2})U_{\theta,\sigma_{\lambda}}v').$$
Cette somme est finie: $\varphi_{2}$ est $K$-finie \`a droite donc $\pi_{\lambda}(\varphi_{2})U_{\theta,\sigma_{\lambda}}$ annule presque tout $v'\in {\cal B}'$. 
L'expression
$$\int_{G(F)}  \bar{\varphi}_{1}(x) \int_{i{\cal A}_{M_{disc},F}^*}m^G(\sigma_{\lambda})\sum_{v'\in {\cal B}'}(\pi_{\lambda}(x\theta(g))U_{\theta,\sigma_{\lambda}}v',\pi_{\lambda}(g)\pi_{\lambda}(\varphi_{2})U_{\theta,\sigma_{\lambda}}v')\,d\lambda\,dx$$
est absolument convergente. En permutant les int\'egrales, on voit qu'elle vaut
$$\int_{i{\cal A}_{M_{disc},F}^*}m^G(\sigma_{\lambda})\sum_{v'\in {\cal B}'}(\pi_{\lambda}(\varphi_{1})\pi_{\lambda}(\theta(g))U_{\theta,\sigma_{\lambda}}v',\pi_{\lambda}(g)\pi_{\lambda}(\varphi_{2})U_{\theta,\sigma_{\lambda}}v')\,d\lambda.$$
 On a l'\'egalit\'e
$$(\pi_{\lambda}(\varphi_{1})\pi_{\lambda}(\theta(g))U_{\theta,\sigma_{\lambda}}v',\pi_{\lambda}(g)\pi_{\lambda}(\varphi_{2})U_{\theta,\sigma_{\lambda}}v')=( \pi_{\lambda}(\varphi_{1})U_{\theta,\sigma_{\lambda}}\pi'_{\lambda'}(g)v',\pi_{\lambda}(g) \pi_{\lambda}(\varphi_{2})U_{\theta,\sigma_{\lambda}}v').$$
On exprime matriciellement tous les op\'erateurs  intervenant. L'expression ci-dessus devient
 $$\sum_{u,v\in {\cal B},u'\in {\cal B}'}( \pi_{\lambda}(\varphi_{1})U_{\theta,\sigma_{\lambda}}u',v)(u,\pi_{\lambda}(\varphi_{2})U_{\theta,\sigma_{\lambda}}v')(v,\pi_{\lambda}(g)u)(\pi'_{\lambda'}(g)v',u').$$
 Cette somme est en fait finie d'apr\`es les propri\'et\'es de $K$-finitude des fonctions $\varphi_{1}$ et $\varphi_{2}$ et de l'op\'erateur $U_{\theta,\sigma_{\lambda}}$. Posons
 $$B_{u,v,u',v'}(\lambda)= ( \pi_{\lambda}(\varphi_{1})U_{\theta,\sigma_{\lambda}}u',v)(u,\pi_{\lambda}(\varphi_{2})U_{\theta,\sigma_{\lambda}}v').$$
 C'est une fonction de Schwartz sur $i{\cal A}_{M_{disc},F}^*$. On obtient l'\'egalit\'e
 $$(1) \qquad J^T_{M_{disc},\sigma}(\omega,f_{1},f_{2})=\sum_{u,v\in {\cal B},u',v'\in {\cal B}'}  \int_{A_{\tilde{G}}(F)\backslash G(F)}$$
 $$\int_{i{\cal A}_{M_{disc},F}^*}m^G(\sigma_{\lambda})B_{u,v,u',v'}(\lambda)(v,\pi_{\lambda}(g)u)(\pi'_{\lambda'}(g)v',u')\,d\lambda\,dx\,\omega(g){\tilde{\kappa}}^T(g)\,dg.$$
 
 \bigskip
 
 \subsection{Une premi\`ere approximation d'une int\'egrale de coefficients}
 Fixons $u,v\in {\cal B}$, $u',v'\in {\cal B}'$ et une fonction de Schwartz $B$ sur $i{\cal A}_{M_{disc},F}^*$. On pose
 $$j^T=
\int_{A_{\tilde{G}}(F)\backslash G(F)}\int_{i{\cal A}_{M_{disc},F}^*}m^G(\sigma_{\lambda})B(\lambda) (v,\pi_{\lambda}(g)u)(\pi'_{\lambda'}(g)v',u')\,d\lambda\,\omega(g)\tilde{\kappa}^T(g)\,dg,$$
o\`u on rappelle que $\lambda'=\theta^{-1}\lambda$.
Cette expression est convergente dans l'ordre indiqu\'e pour les m\^emes raisons que pr\'ec\'edemment.  La formule (1) du paragraphe pr\'ec\'edent exprime $J^T_{M_{disc},\sigma}(\omega,f_{1},f_{2})$ comme combinaison lin\'eaire de telles expressions $j^T$. 

On d\'efinit une fonction $\Omega_{G}$ sur $M_{0}(F)$ par
$$\Omega_{G}(m)=\int_{K\times K}\int_{i{\cal A}_{M_{disc},F}^*}m^G(\sigma_{\lambda})B(\lambda) (v,\pi_{\lambda}(kmk')u)(\pi'_{\lambda'}(kmk')v',u')\,d\lambda\,\omega(kmk')\,dk\,dk'.$$
Alors
$$j^T=\int_{A_{\tilde{G}}(F)\backslash M_{0}(F)^{\geq}}\Omega_{G}(m)D_{0}(m)\tilde{\kappa}^T(m)\,dm.$$

 Soit $Q=LU_{Q}$ un sous-groupe parablique standard de $G$.  Comme en 1.12, notons $W^G(L\vert  S)$ l'ensemble des $w\in W^G/W^{M_{disc}}$ tels que $w(M_{disc})\subset L$, $w(S)\cap L\supset P_{0}\cap L$. On identifiera souvent cet ensemble \`a un sous-ensemble de $W^G$ form\'e d'\'el\'ements $w$ de longueur minimale dans leur classe $W^Lw$. Pour un \'el\'ement $w$ de cet ensemble, on note $Q_{w}=(w(S)\cap L)U_{Q}$, $\underline{Q}_{w}=(w(S)\cap L)U_{\bar{Q}}$.  Notons $\pi_{w}=Ind_{Q_{w}}^G(w\sigma)$ et $\underline{\pi}_{w}=Ind_{\underline{Q}_{w}}^G(w\sigma)$, r\'ealis\'ees dans leurs espaces habituels $V_{w\sigma,Q_{w}}$ et $V_{w\sigma,\underline{Q}_{w}}$. Pour $v\in V_{w\sigma,Q_{w}}$ et $\underline{v}\in V_{w\sigma,\underline{Q}_{w}}$, on pose
$$(v,\underline{v})^L=\int_{K\cap L(F)}(v(k),\underline{v}(k))\,dk.$$
Le produit int\'erieur est le produit hermitien sur $V_{w\sigma}$. On pose des d\'efinitions analogues en rempla\c{c}ant $S$ par $S'=\theta^{-1}(S)$.  On prendra garde que les d\'efinitions de $Q_{w}$ et $\underline{Q}_{w}$ changent: pour $w\in W^G(L\vert  S')$, on a $Q_{w}=(w(S')\cap L)U_{Q}$.

Soient $w\in W^G(L\vert  S)$ et $w'\in W^G(L\vert  S')$.    On d\'efinit une fonction $\omega_{Q,w,w'}$ sur $G(F)\times K\times K\times {\cal A}_{M_{disc},{\mathbb C}}^*$ par
$$\omega_{Q,w,w'}(g,k,k',\lambda)=(J_{Q_{w}\vert w(S)}((w\sigma)_{w\lambda})\circ \gamma(w)\circ \pi_{\lambda}(k)v,J_{\underline{Q}_{w}\vert w(S)}((w\sigma)_{w\lambda})\circ\gamma(w)\circ\pi_{\lambda}(gk')u)^L$$
$$( J_{\underline{Q}_{w'}\vert w'(S')}((w'\sigma')_{w'\lambda'})\circ \gamma(w')\circ \pi'_{\lambda'}(gk')v',J_{Q_{w'}\vert w'(S')}((w'\sigma')_{w'\lambda'})\circ\gamma(w')\circ\pi'_{\lambda'}(k)u')^L.$$
  Cette fonction est m\'eromorphe en $\lambda$. Pour comprendre les questions de r\'egularit\'e et de croissance en $\lambda$, il est  commode de r\'ecrire cette d\'efinition en utilisant les op\'erateurs normalis\'es et les facteurs de normalisation. Pour $\lambda\in i{\cal A}_{M_{disc},F}^*$, on a l'\'egalit\'e
$$(1) \qquad \omega_{Q,w,w'}(g,k,k',\lambda)={\bf r}_{w,w'}(\sigma_{\lambda})$$
$$( R_{Q_{w}\vert w(S)}((w\sigma)_{w\lambda})\circ \gamma(w)\circ \pi_{\lambda}(k)v,R_{\underline{Q}_{w}\vert w(S)}((w\sigma)_{w\lambda})\circ\gamma(w)\circ\pi_{\lambda}(gk')u,)^L$$
$$( R_{\underline{Q}_{w'}\vert w'(S')}((w'\sigma')_{w'\lambda'})\circ \gamma(w')\circ \pi'_{\lambda'}(gk')v',R_{Q_{w'}\vert w'(S')}((w'\sigma')_{w'\lambda'})\circ\gamma(w')\circ\pi'_{\lambda'}(k)u')^L,$$
o\`u
$${\bf r}_{w,w'}(\sigma_{\lambda})=r_{ \underline{Q}_{w}\vert Q_{w}}((w\sigma)_{w\lambda})r_{Q_{w'}\vert \underline{Q}_{w'} }((w'\sigma')_{w'\lambda'}).$$
Les op\'erateurs normalis\'es sont holomorphes et unitaires sur $i{\cal A}_{M_{disc},F}^*$. Les produits scalaires de l'expression ci-dessus  sont donc  holomorphes et  born\'es en $\lambda$. On a

(2) le produit $m^G(\sigma_{\lambda}){\bf r}_{w,w'}(\sigma_{\lambda})$ est holomorphe et \`a croissance mod\'er\'ee sur $i{\cal A}_{M_{disc},F}^*$.

En effet, par transport de structure, on a aussi
$${\bf r}_{w,w'}(\sigma_{\lambda})=r_{S_{1}\vert S_{2}}(\sigma_{\lambda})r_{S_{3}\vert S_{4}}(\sigma_{\lambda}),$$
o\`u $S_{1}=w^{-1}(\underline{Q}_{w})$, $S_{2}=w^{-1}(Q_{w})$, $S_{3}=\theta((w')^{-1}(Q_{w'}))$, $S_{4}=\theta((w')^{-1}(\underline{Q}_{w'}))$. Ces quatre paraboliques appartiennent \`a ${\cal P}(M_{disc})$ et il reste \`a appliquer 1.10(7).

 On pose
$$\omega_{Q,w,w'}(g)=\int_{K\times K}\int_{i{\cal A}_{M_{disc},F}^*}\omega_{Q,w,w'}(g,k,k',\lambda)m^G(\sigma_{\lambda})B(\lambda)\omega(k^{-1}gk')\,d\lambda\,dk\,dk'.$$

Pour tout sous-groupe parabolique standard $Q'=L'U_{Q'}\supset Q$, on note ${\cal W}_{Q}^{Q'}$ l'ensemble des $(w,w')\in W^G(L\vert  S)\times W^G(L\vert  S')$ tels que $\theta(W^{L'}w')\cap W^{L'}w\not=\emptyset$. On pose simplement ${\cal W}_{Q}={\cal W}_{Q}^Q$.

Soit $R$ un parabolique standard contenant $Q$.  Pour $w\in W^G(L\vert  S)$ et $w'\in W^G(L\vert  S')$, notons $s_{Q}^R(w,w')$ la somme des $(-1)^{a_{\tilde{P}}-a_{\tilde{G}}}$ sur les ensembles paraboliques $\tilde{P}$ tels que $Q\subset P\subset R$ et  $(w,w')\in {\cal W}_{Q}^P$. Remarquons que, si cet ensemble d'ensembles paraboliques n'est pas vide, il existe $\tilde{P}_{-}\subset \tilde{P}_{+}$ de sorte que cet ensemble soit simplement celui des $\tilde{P}$ tels que $\tilde{P}_{-}\subset \tilde{P}\subset \tilde{P}_{+}$. Donc $s_{Q}^R(w,w')$ est nul si l'ensemble est vide ou si $\tilde{P}_{-}\not=\tilde{P}_{+}$ et est \'egal \`a $(-1)^{a_{\tilde{P}_{+}}-a_{\tilde{G}}}$ si l'ensemble est non vide et $\tilde{P}_{-}=\tilde{P}_{+}$. Posons
$$j^T_{\star}=\sum_{Q=LU_{Q},R; P_{0}\subset Q\subset R}\sum_{w\in W^G(L\vert  S),w'\in W^G(L\vert  S')}s_{Q}^R(w,w')$$
$$\int_{A_{\tilde{G}}(F)\backslash M_{0}(F)^{\geq,Q}}\delta_{Q}(m)D_{0}^L(m)\tilde{\sigma}_{Q}^R(H_0(m)-T)\phi^Q(H_0(m)-T)\omega_{Q,w,w'}(m)\,dm.$$

\ass{Proposition}{L'expression ci-dessus est absolument convergente. Pour tout r\'eel $r$, on a la majoration
$$\vert j^T-j^T_{\star}\vert <<\vert T\vert ^{-r}$$
pour tout $T$.}

La preuve de cette proposition sera donn\'ee en 3.15. 

\bigskip

\subsection{Un lemme de majoration}
Soit $Q=LU_{Q}$ un parabolique standard. Pour $H\in {\cal A}_{0}$, on note ${\cal C}^L(H)$ l'enveloppe convexe des $sH$ pour $s\in W^L$. On pose
$$N^L(H)=1+inf\{\vert\theta  H'-H''\vert ; H',H''\in {\cal C}^L(H)\}.$$
Soient $w\in W^G(L\vert  S)$ et $w'\in W^G(L\vert  S')$, que l'on rel\`eve  en des \'el\'ements de $W^G$.  On pose
$$N^L_{w,w'}(H)=1+inf\{\vert w\circ\theta\circ (w')^{-1} H'-H''\vert ; H',H''\in {\cal C}^L(H)\}.$$
Ce terme ne d\'epend pas des rel\`evements choisis: changer de rel\`evement multiplie $w$ et $w'$ \`a gauche par des \'el\'ements de $W^L$ dont l'action conserve ${\cal C}^L(H)$. Remarquons que, dans le cas o\`u $W^Lw\cap \theta(W^Lw')\not=\emptyset$, on peut pour la m\^eme raison remplacer $w'$ et $w$ par des \'el\'ements tels que $w=\theta(w')$ et on obtient $N^L_{w,w'}(H)=N^L(H)$.

\ass{Lemme}{Soit $Q=LU_{Q}$ un sous-groupe parabolique standard.

(i) Soient $w\in W^G(L\vert  S)$ et $w'\in W^G(L\vert  S')$ et soit $\tilde{P}=\tilde{M}U_{P}$ le plus petit espace parabolique standard tel que $Q\subset P$ et $(w,w')\in {\cal W}_{Q}^P$. Alors on a une majoration
$$\vert (H-T)_{L}^M\vert <<N_{w,w'}^L(H)$$
pour tout $T$ et tout $H\in {\cal A}_{0}$ tel que $\phi^Q(H-T)\tau_{Q}^P(H-T)=1$ et $<\alpha,H>\geq0$ pour tout $\alpha\in \Delta_{0}^Q$.

(ii) Soit $\tilde{P}=\tilde{M}U_{P}$ le plus petit espace parabolique standard tel que $Q\subset P$. Alors on a une majoration
$$\vert (H-T)_{L}^M\vert <<N^L(H)$$
pour tout $T$ et tout $H\in {\cal A}_{0}$ tel que $\phi^Q(H-T)\tau_{Q}^P(H-T)=1$ et $<\alpha,H>\geq0$ pour tout $\alpha\in \Delta_{0}^Q$.

(iii) Soient $w\in W^G(L\vert  S)$ et $w'\in W^G(L\vert  S')$ et soit $\tilde{P}=\tilde{M}U_{P}$ le plus petit espace parabolique standard tel que $Q\subset P$ et $(w,w')\in {\cal W}_{Q}^P$. Supposons $(w,w')\not\in {\cal W}_{Q}$. Alors on a une majoration
$$\vert T\vert +\vert (H-T)^M_{L}\vert <<N^L_{w,w'}(H)$$
pour tout $T$ et tout $H\in {\cal A}_{0}$ tel que $\phi^Q(H-T)\tau_{Q}^P(H-T)=1$ et $<\alpha,H>\geq0$ pour tout $\alpha\in \Delta_{0}^Q$.

(iv) Soient $w\in W^G(L\vert  S)$ et $w'\in W^G(L\vert  S')$ et soit $R$ un sous-groupe parabolique contenant $Q$. Supposons $s_{Q}^R(w,w')\not=0$. Alors   on a une majoration
$$\vert (H-T)_{L}^{\tilde{G}}\vert <<N^L_{w,w'}(H)$$
pour tout $T$ et tout $H\in {\cal A}_{0}$ tel que $\phi^Q(H-T)\tilde{\sigma}_{Q}^R(H-T)=1$ et $<\alpha,H>\geq0$ pour tout $\alpha\in \Delta_{0}^Q$.}

Preuve de (i). Les hypoth\`eses $\phi^Q(H-T)=1$ et  $<\alpha,H>\geq0$ pour tout $\alpha\in \Delta_{0}^Q$ entra\^{\i}nent que $H^L$ appartient \`a ${\cal C}^L(T^L)$. Donc ${\cal C}^L(H)\subset {\cal C}^L(H_{L}+T^L)$. Posons $H_{\star}=H_{L}-T_{L}$. On obtient ${\cal C}^L(H)\subset {\cal C}^L(T+H_{\star})$. Soient $H',H''\in {\cal C}^L(T+H_{\star})$. Posons $s=w\theta(w')^{-1}$.  On veut minorer $\vert H'-s\theta H''\vert $. On a 
$$(1) \qquad \vert H'-s\theta H''\vert \geq \vert X\vert ,$$
 o\`u $X=(H'-s\theta H'')_{L}^M$. On a $H_{L}^{_{'}M}=T_{L}^M+H_{\star}^M$ parce que $H'\in {\cal C}^L(T+H_{\star})$, d'o\`u $X= T_{L}^M+H_{\star}^M-(s\theta H'')_{L}^M $. L'hypoth\`ese que $H''\in {\cal C}^L(T+H_{\star})$ signifie que l'on peut \'ecrire $H''=H_{\star}+\sum_{u\in W^L}y_{u}uT$, avec des $y_{u}\geq0$ tels que $\sum_{u\in W^L}y_{u}=1$. Donc
 $$X=T_{L}^M+H^M_{\star}-(s\theta H_{\star})_{L}^M-\sum_{u\in W^L}y_{u}(s\theta uT)_{L}^M.$$
  Posons $Y_{u}=T_{L}^M-(s\theta uT)_{L}^M$ et $Y=\sum_{u\in W^L}y_{u}Y_{u}$. Alors $ X=H_{\star}^M-(s\theta H_{\star})_{L}^M+Y$. Pour tout $u\in W^L$, on a $Y_{u}=(T-s\theta uT)_{L}^M=(T-u' T)_{L}^M$, o\`u $u'=s\theta(u)$, puisque $\theta T=T$. L'hypoth\`ese $(w,w')\in {\cal W}_{Q}^P$ entra\^{\i}ne que $s\in W^M$ donc aussi $u'\in W^M$. Puisque $T$ est dominant, $T-u'T$ est combinaison lin\'eaire \`a coefficients positifs ou nuls de $\check{\alpha}$ pour $\alpha\in \Delta_{0}^P$. Donc $Y_{u}$ est combinaison lin\'eaire \`a coefficients positifs ou nuls de $\check{\alpha}_{L}$ pour $\alpha\in \Delta_{0}^P-\Delta_{0}^Q$.  L'\'el\'ement $Y$ v\'erifie la m\^eme propri\'et\'e. Par ailleurs, l'hypoth\`ese $\tau_{Q}^P(H-T)=1$ signifie que $H_{\star}^M$ est combinaison lin\'eaire \`a coefficients strictement positifs de $\check{\varpi}_{\alpha}^M$ pour $\alpha\in \Delta_{0}^P-\Delta_{0}^Q$. Soit $U=\sum_{\alpha\in \Delta_{0}^P-\Delta_{0}^Q}u_{\alpha}\check{\varpi}_{\alpha}^M$ et $V=\sum_{\alpha\in \Delta_{0}^P-\Delta_{0}^Q}v_{\alpha}\check{\alpha}_{L}$ avec des coefficients $u_{\alpha}$ et $v_{\alpha}$ positifs ou nuls. Montrons que l'on a une majoration
$$(2)\qquad \vert U\vert +\vert V\vert << \vert U-(s\theta U)_{L}+V\vert .$$
Il suffit de prouver que le c\^one engendr\'e par les $\check{\varpi}_{\alpha}^M-(s\theta\check{\varpi}_{\alpha}^M)_{L}$ et les $\check{\alpha}_{L}$ est un "vrai" c\^one, c'est-\`a-dire ne contient pas d'espace vectoriel non nul. Il revient au m\^eme de prouver que, pour $U$ et $V$ comme ci-dessus, l'\'egalit\'e $U-(s\theta U)_{L}+V=0$ entra\^{\i}ne $U=0$, $V=0$. Or le produit scalaire $(U,V)$ est positif ou nul. Par produit scalaire avec $U$, l'\'egalit\'e $U-(s\theta U)_{L}+V=0$ entra\^{\i}ne donc $(U,U-(s\theta U)_{L})\leq0$, d'o\`u $(U,U)\leq (U,(s\theta U)_{L})$. Par Cauchy-Schwartz, cela implique $U=(s\theta U)_{L}
$ puis $U=s\theta U$. Si $U\not=0$, notons $\Delta_{0}(U)$ l'ensemble des $\alpha\in \Delta_{0}^P$ tels que $<\alpha,U>>0$. Cet ensemble  est non vide et inclus dans $\Delta_{0}^P-\Delta_{0}^Q$. Introduisons le sous-groupe parabolique standard $P'=M'U_{P'}$ tel que $\Delta_{0}^{M'}=\Delta_{0}^P-\Delta_{0}(U)$.On a $Q\subset P'\subsetneq P$. L'\'el\'ement  $U$ appartient \`a la chambre positive relative au sous-groupe parabolique $M\cap P'$ de $M$. L'\'egalit\'e $U=s\theta U$ entra\^{\i}ne $P'=s\theta(P')$. Puisque $P'$ et $s\theta(P')$ sont standard, cela implique $P'=\theta(P') $ et $s\in W^{M'}$. Mais alors $(w,w')\in {\cal W}_{Q}^{P'}$, ce qui contredit l'hypoth\`ese de minimalit\'e de $\tilde{P}$. Cette contradiction prouve que $U=0$. L'\'egalit\'e $U-(s\theta U)_{L}+V=0$ entra\^{\i}ne alors que $V=0$, ce qui prouve (2).

On applique (2) \`a $U=H_{\star}^M$ et $V=Y$. En abandonnant le terme $\vert Y\vert $, on obtient la majoration $\vert H_{\star}^M\vert <<\vert X\vert $. Gr\^ace \`a (1), cela entra\^{\i}ne la majoration du (i) de l'\'enonc\'e. 

La preuve du (ii) est similaire, il suffit de supprimer le terme $s$ des calculs.

Preuve de (iii).   On reprend la preuve de (i). A la fin de cette preuve, on avait abandonn\'e le terme $\vert Y\vert $. R\'etablissons-le. Pour obtenir le (iii) de l'\'enonc\'e, il suffit de prouver que, pour tout $u\in W^L$, on a une majoration 
$$(3)\qquad \vert T\vert <<\vert Y_{u}\vert .$$
Notons $\Sigma^M$ l'ensemble des racines de $A_{0}$ dans l'alg\`ebre de Lie de $M$, muni de la positivit\'e d\'efinie par $P_{0}\cap M$. On sait plus pr\'ecis\'ement que $T-u'T$ est combinaison lin\'eaire des $\check{\alpha}$ pour toutes les racines $\alpha\in \Sigma^M$ telles que $\alpha>0$ et $(u')^{-1}\alpha<0$. Les coefficients sont de la forme $<\beta,T>$ pour des $\beta\in \Sigma^M$, $\beta>0$, donc sont essentiellement minor\'es par $\vert T\vert $.  L'assertion (3) s'en d\'eduit pourvu qu'il y ait au moins une racine $\alpha$ telle que $\alpha>0$, $(u')^{-1}\alpha<0$ et $\check{\alpha}_{L}\not=0$. S'il n'en est pas ainsi,  on a $u'\in W^L$. L'\'egalit\'e $u'=w\theta(w')^{-1}\theta(u)$, jointe au fait que $u\in W^L$, entra\^{\i}ne alors que $W^Lw\cap \theta(W^Lw')\not=\emptyset$, contrairement \`a l'hypoth\`ese. Cela prouve (3) et ach\`eve la preuve du (iii) de l'\'enonc\'e.

Preuve de (iv). L'hypoth\`ese $s_{Q}^R(w,w')\not=0$ signifie qu'il existe un unique espace parabolique $\tilde{P}=\tilde{M}U_{P}$ tel que $Q\subset P\subset R$ et $(w,w')\in {\cal W}_{Q}^P$.
On reprend la preuve du (i). On a encore $s\in W^M$. Pour $H',H''\in {\cal C}^L(T+H_{\star})$, on a
$$\vert H'-s\theta H''\vert >>\vert (H'-s\theta H'')_{L}^M\vert +\vert (H'-s\theta H'')_{M}\vert .$$
La preuve de (i) s'applique au premier terme: on a une majoration
$$(4) \qquad \vert (H'-s\theta H'')_{L}^M\vert >> \vert H_{\star}^M\vert ,$$
o\`u $H_{\star}=(H-T)_{L}$. On a $(H'-s\theta H'')_{M}=H_{\star,M}-\theta(H_{\star,M})$ puisque $s\in W^M$.  Ecrivons $H_{\star, M}=H_{1}+H_{2}$ o\`u $H_{1}\in {\cal A}_{\tilde{M}}$ et $H_{2}$ appartient \`a l'orthogonal ${\cal A}_{M}^{\tilde{M}}$ de ce sous-espace dans ${\cal A}_{M}$. On a $H_{\star,M}-\theta(H_{\star,M})=H_{2}-\theta H_{2}$ et $1-\theta$ est injective sur ${\cal A}_{M}^{\tilde{M}}$, d'o\`u une majoration
$$(5)\qquad \vert H_{2}\vert <<\vert H_{\star,M}-\theta( H_{\star,M})\vert.$$
La condition $\tilde{\sigma}_{Q}^R(H_{\star})=1$ entra\^{\i}ne que

(6) $<\varpi_{\alpha},H_{1}>>0$ pour tout $\alpha\in \Delta_{0}-\Delta_{0}^P$.

Elle entra\^{\i}ne aussi que $ <\alpha,H_{\star}>\leq 0$ pour tout $\alpha\in \Delta_{0}-\Delta_{0}^R$. On remarque que, pour tout $\alpha\in \Delta_{0}-\Delta_{0}^P$, il existe $i\geq1$ tel que $\theta^{i}\alpha\not\in \Delta_{0}^R$: sinon l'espace  parabolique standard $\tilde{P}'$ associ\'e \`a la r\'eunion de $\tilde{\Delta}_{0}^P$ et de $\{\tilde{\alpha}\}$ v\'erifierait $P\subsetneq P'\subset R$, ce qui contredirait l'hypoth\`ese d'unicit\'e de $\tilde{P}$. Pour $\alpha\in \Delta_{0}-\Delta_{0}^P$ et  $i\geq1$ comme ci-dessus, on a 

$$(7) \qquad<\alpha,H_{1}>=<\theta^{i}\alpha,H_{1}>=<\theta^{i}\alpha,H_{\star}-H_{2}-H_{\star}^M>\leq -<\theta^{i}\alpha,H_{2}+H_{\star}^M>.$$
  Les conditions (6) et (7) bornent $\vert H^{\tilde{G}}_{1}\vert $: on a une majoration
$$\vert H_{1}^{\tilde{G}}\vert <<\vert H_{2}\vert +\vert H_{\star}^M\vert .$$
Gr\^ace \`a (5), on en d\'eduit
$$\vert H_{\star,M}^{\tilde{G}}\vert <<\vert H_{\star,M}-\theta (H_{\star,M})\vert+\vert H_{\star}^M\vert .$$
 D'o\`u, gr\^ace \`a (4),
$$\vert H_{\star}^{\tilde{G}}\vert << \vert H_{\star,M}-\theta (H_{\star,M})\vert+ \vert (H'-s\theta H'')_{L}^M\vert <<\vert H'-s\theta H''\vert .$$
Cela prouve l'assertion (iv) de l'\'enonc\'e. $\square$

\bigskip

\subsection{Majoration de coefficients}
 \ass{Lemme}{Soient $Q=LU_{Q}$ un parabolique standard, $w\in W^G(L\vert  S)$ et $w'\in W^G(L\vert S')$. Quel que soit le r\'eel $r$, il existe $c>0$ tel que l'on ait la majoration
$$\vert \omega_{Q,w,w'}(m)\vert\leq c\delta_{Q}(m)^{-1}\Xi^L(m)^2N^L_{w,w'}(H_0(m))^{-r}$$
pour tout $m\in M_{0}(F)$. }

Preuve. Notons $\rho_{w\lambda}^G=Ind_{\underline{Q}_{w}}^G((w\sigma)_{w\lambda})$, $\rho^{_{'}G}_{w'\lambda'}=Ind_{\underline{Q}_{w'}}((w'\sigma')_{w'\lambda'})$, $\rho_{w\lambda}=Ind_{w(S)\cap L}^L((w\sigma)_{w\lambda})$, $\rho'_{w'\lambda'}=Ind_{w'(S')\cap L}^L((w'\sigma')_{w'\lambda'})$ que l'on r\'ealise dans leurs espaces habituels  $V_{w\sigma,\underline{Q}_{w}}$, $V_{w'\sigma',\underline{Q}_{w'}}$, $V^L_{w\sigma,w(S)\cap L}$, $V^L_{w'\sigma',w'(S')\cap L}$. D'apr\`es 3.4(1) et (2), $\omega_{Q,w,w'}(m,k,k',\lambda)m^G(\sigma_{\lambda})$ est combinaison lin\'eaire de termes
$$(v_{0},\rho^G_{w\lambda}(m)u_{0})^L(\rho^{_{'}G}_{w'\lambda'}(m)v'_{0},u'_{0})^L$$
pour des \'el\'ements $u_{0},v_{0}\in V_{w\sigma,\underline{Q}_{w}}$ et $u'_{0},v'_{0}\in V_{w'\sigma',\underline{Q}_{w'}}$. Les coefficients sont des fonctions $C^{\infty}$ de $\lambda$, $k$ et $k'$ et sont \`a croissance mod\'er\'ee en $\lambda$. On dispose d'applications
$$V_{w\sigma,\underline{Q}_{w}}\to V^L_{w\sigma,w(S)\cap L},\,\,V_{w'\sigma',\underline{Q}_{w'}}\to V^L_{w'\sigma',w'(S')\cap L}$$
qui sont les restrictions \`a $K\cap L(F)$. En notant $e,f,e',f'$ les images de $u_{0},v_{0},u'_{0}, v'_{0}$ par ces applications, on a les \'egalit\'es
$$(v_{0,}\rho^G_{w\lambda}(m)u_{0})^L=\delta_{Q}(m)^{-1/2}(f,\rho_{w\lambda}(m)e),$$
$$(\rho^{_{'}G}_{w'\lambda'}(m)v'_{0},u'_{0})^L=\delta_{Q}(m)^{-1/2}(\rho'_{w'\lambda'}(m)f',e').$$
Quitte \`a remplacer la fonction  de Schwartz $B$ par son produit avec une  fonction \`a croissance mod\'er\'ee, on est ramen\'e \`a \'evaluer
$$(1) \qquad \delta_{Q}(m)^{-1}\int_{i{\cal A}_{M_{disc},F}^*} (f,\rho_{w\lambda}(m)e)(\rho'_{w'\lambda'}(m)f',e')B(\lambda)\,d\lambda.$$
On a
$$(f,\rho_{w\lambda}(m)e)=\int_{K\cap L(F)}(f(k),(\rho_{w\lambda}(m)e)(k))\,dk,$$
o\`u le produit int\'erieur est le produit hermitien sur $V_{w\sigma}$. On a
$$(\rho_{w\lambda}(m)e)(k)=e^{<w\lambda,H_{\underline{Q}_{w}}(km)>}(\rho_{w}(m)e)(k),$$
o\`u on a not\'e $\rho_{w}$ la repr\'esentation $\rho_{w\lambda}$ pour $\lambda=0$.  En utilisant des formules similaires pour le terme $(\rho'_{w'\lambda'}(m)f',e')$, on obtient que (1) est \'egal \`a
$$\delta_{Q}(m)^{-1}\int_{(K\cap L(F))\times (K\cap L(F))}(f(k),(\rho_{w}(m)e)(k))((\rho'_{w'}(m)f')(k'),e'(k'))\beta(m,k,k') \,dk\,dk',$$
o\`u
$$\beta(m,k,k')=\int_{i{\cal A}_{M_{disc},F}^*}B(\lambda)e^{<w\lambda,H_{\underline{Q}_{w}}(km)>-<w'\lambda',H_{\underline{Q}_{w'}}(k'm)>}\,d\lambda.$$
On va prouver que, pour tout r\'eel $r$, il existe $c>0$ tel que
$$(2)\qquad \vert \beta(m,k,k')\vert \leq cN^L_{w,w'}(H_0(m))^{-r}.$$
Admettons cela et finissons la d\'emonstration. De (2) r\'esulte que (1) est essentiellement major\'e par
$$\delta_{Q}(m)^{-1}N^L_{w,w'}(H_0(m))^{-r}$$
$$\int_{(K\cap L(F))\times (K\cap L(F))}\vert (f(k),(\rho_{w}(m)e)(k))((\rho'_{w'}(m)f')(k'),e'(k'))\vert  \,dk\,dk'.$$
En utilisant 1.12(4), la double int\'egrale est essentiellement major\'ee par $\Xi^L(m)^2$ et on obtient la majoration de l'\'enonc\'e.

D\'emontrons (2). On commence par prouver

(3) pour tout parabolique standard $Q'=L'U_{Q'}\subset Q$ et tout $k\in K\cap L(F)$, $H_{Q'}(km)$ appartient \`a ${\cal C}^L(H_0(m))$.

Quitte \`a conjuguer $m$ par un \'el\'ement de $K\cap Norm_{L(F)}(M_{0})$, ce qui ne change pas notre probl\`eme, on peut supposer $<\alpha,H_0(m)>\geq0$ pour tout $\alpha\in \Delta_{0}^Q$. D'apr\`es le lemme 1.3,  $H_{0}(km)$ appartient \`a ${\cal C}^L(H_0(m))$.   L'\'el\'ement $H_{Q'}(km)$ est la projection orthogonale  de $H_{0}(km)$ sur ${\cal A}_{L'}$. On est ramen\'e \`a montrer que, pour $H\in {\cal C}^L(H_0(m))$, alors la projection $H_{L'}$ de $H$ sur ${\cal A}_{L'}$ appartient aussi \`a ${\cal C}^L(H_0(m))$. L'espace ${\cal A}_{L'}^L$ est r\'eunion des chambres positives  associ\'ees aux paraboliques $Q''\in {\cal P}^L(L')$. On fixe un tel $Q''$ tel que $H_{L'}^L$ appartienne \`a la chambre associ\'ee \`a $Q''$ et on fixe $s\in W^L$ tel que $s(Q'')$ soit standard. Alors $H_{L'}=s^{-1}((sH)_{s(L')})$. Puisque ${\cal C}^L(H_{0}(m))$ est invariant par l'action de $s$, on peut aussi bien remplacer $H$ par $sH$ et $L'$ par $s(L')$. En oubliant cette construction, on est ramen\'e au cas o\`u  $<\alpha,H_{L'}>\geq 0$ pour tout $\alpha\in \Delta_{0}^Q$. Dans ce cas, l'appartenance \`a ${\cal C}^L(H_0(m))$ \'equivaut \`a   ce que $ H_{0}(m)-H_{L'}$ soit une combinaison lin\'eaire \`a coefficients positifs ou nuls  de $\check{\alpha}$ pour $\alpha\in\Delta_{0}^Q$. On a
$$H_{0}(m)-H_{L'}=H_{0}(m)-H_{0}(m)_{L'}+H_{0}(m)_{L'}-H_{L'}.$$
La premi\`ere diff\'erence est \'egale \`a $H_{0}(m)^{L'}$ qui appartient \`a la chambre positive ferm\'ee de ${\cal A}_{0}^{L'}$, a fortiori est combinaison lin\'eaire \`a coefficients positifs ou nuls de $\check{\alpha}$ pour $\alpha\in \Delta_{0}^{Q'}$. La deuxi\`eme diff\'erence est la projection de $H_{0}(m)-H$. Puisque $H\in {\cal C}^L(H_{0}(m))$, $H_{0}(m)-H$ est une  combinaison lin\'eaire \`a coefficients positifs ou nuls  de $\check{\alpha}$ pour $\alpha\in\Delta_{0}^Q$. Donc $H_{0}(m)_{L'}-H_{L'}$ est une telle combinaison lin\'eaire de $\check{\alpha}_{L'}$, pour $\alpha\in \Delta_{0}^Q$. Il suffit de voir qu'une telle projection est de la forme voulue. C'est clair si $\alpha\in \Delta_{0}^{Q'}$ puisqu'alors $\check{\alpha}_{L'}=0$. Soit $\alpha\in \Delta_{0}^Q-\Delta_{0}^{Q'}$. Alors $\check{\alpha}_{L'}=\check{\alpha}-\check{\alpha}^{L'}$. Parce que $\{\check{\alpha}; \alpha\in\Delta_{0}^Q\}$ est une base obtuse de ${\cal A}_{0}^L$, $-\check{\alpha}^{L'}$ appartient \`a la chambre positive ferm\'ee de ${\cal A}_{0}^{L'}$. A fortiori, $-\check{\alpha}^{L'}$ est combinaison lin\'eaire \`a coefficients positifs ou nuls de $\check{\beta}$ pour $\beta\in\Delta_{0}^{Q'}$. Cela prouve (3).

Le terme  $\beta(m,k,k')$ est la transform\'ee de Fourier de $B$ \'evalu\'ee en $\theta((w')^{-1}H') -w^{-1} H$, o\`u $H=H_{\underline{Q}_{w}}(km)$ et $H'=H_{\underline{Q}_{w'}}(k'm)$. Il est donc essentiellement major\'e par $(1+\vert \theta((w')^{-1}H') -w^{-1} H\vert )^{-r}$ pour tout r\'eel $r$, ou encore par $ (1+\vert (w\circ\theta\circ(w')^{-1})H'- H\vert )^{-r}$. Puisque $k'm$ et $km$ appartiennent \`a $L(F)$, on a les \'egalit\'es $H=H_{Q_{w}}(km)$ et $H'=H_{Q_{w'}}(k'm)$. Les paraboliques $Q_{w}$ et $Q_{w'}$  sont standard. Gr\^ace \`a (3), $H$ et $H'$ appartiennent \`a ${\cal C}^L(H_{0}(m))$. Donc
$1+\vert \theta((w')^{-1}H') -w^{-1} H\vert \geq N^L_{w,w'}(H_{0}(m))$ et la majoration (2) s'ensuit. $\square$

Pour tout parabolique standard $Q=LU_{Q}$ et pour $m\in M_{0}(F)$, posons
$$\Omega_{Q}(m)=\sum_{(w,w')\in {\cal W}_{Q}}\omega_{Q,w,w'}(m).$$
  Remarquons que, dans le cas o\`u $Q=G$, on retrouve la fonction $\Omega_{G}$ d\'ej\`a d\'efinie.

\ass{Corollaire}{Pour tout parabolique standard $Q=LU_{Q}$ et pour tout r\'eel $r$, il existe $c>0$ tel que l'on ait la majoration
$$\vert \Omega_{Q}(m)\vert \leq c\delta_{Q}(m)^{-1}\Xi^L(m)^2N^L(H_{0}(m))^{-r}$$
pour tout $m\in M_{0}(F)$.}

\bigskip

\subsection{Un  lemme d'\'equivalence}
 Soient $x^T(m)$ et $y^T(m)$ deux fonctions des variables $T$ et $m\in M_{0}(F)^{\geq}$. On dit que ces fonctions sont \'equivalentes si et seulement si elles v\'erifient la condition suivante: pour tout r\'eel $\nu>0$ et pour tout r\'eel $r$, il existe $c>0$ de sorte que l'on ait l'in\'egalit\'e
$$\vert x^T(m)-y^T(m)\vert \leq c \vert T\vert ^{-r}$$
pour tout $T$ et tout $m\in M_{0}(F)^{\geq}$ tel que $\vert H_0(m)\vert \leq \nu\vert T\vert $. 

\ass{Lemme}{Soit $\tilde{P}=\tilde{M}U_{P}$ un espace parabolique contenant $\tilde{P}_{0}$, soit $Q=LU_{Q}$ un sous-groupe parabolique tel que $P_{0}\subset Q\subset P$ et soit $\epsilon>0$. 

(i) Les fonctions $\phi^Q(H_{0}(m)-\epsilon T)\tau_{Q}^P(H_{0}(m)-\epsilon T)\Omega_{P}(m)\delta_{P_{0}}(m)$ et $\gamma(M\vert L)^{2}\phi^Q(H_{0}(m)-\epsilon T)\tau_{Q}^P(H_{0}(m)-\epsilon T)\Omega_{Q}(m)\delta_{P_{0}}(m)$ sont \'equivalentes.

(ii) Soient $w\in W^G(L\vert  S)$ et $w'\in W^G(L\vert  S')$. Supposons $(w,w')\in {\cal W}_{Q}^P-{\cal W}_{Q}$. Alors la fonction  $\phi^Q(H_{0}(m)-\epsilon T)\tau_{Q}^P(H_{0}(m)-\epsilon T)\omega_{Q,w,w'}(m)\delta_{P_{0}}(m)$ est \'equivalente \`a $0$.}

Preuve. On fixe $\nu>0$. Les $m$ \`a consid\'erer v\'erifient

- $m\in M_{0}(F)^{\geq}$;

-  $\phi^Q(H_{0}(m)-\epsilon T)\tau_{Q}^P(H_{0}(m)-\epsilon T)=1$, ce qui entra\^{\i}ne $<\alpha,H_{0}(m)>><\alpha,\epsilon T>$ pour tout $\alpha\in \Delta_{0}^P-\Delta_{0}^Q$;

- $\vert H_{0}(m)\vert\leq \nu\vert T\vert $. 

En cons\'equence

(1) il existe  un r\'eel $\nu'$ ne d\'ependant que de $\nu$ tel que l'on ait $<\alpha,H_{0}(m)>>\nu'\vert H_{0}(m)\vert $ pour $\alpha\in \Delta_{0}^P-\Delta_{0}^Q$. 

On peut donc approximer les produits scalaires intervenant dans $\Omega_{P}(m)$ par leurs termes constants faibles. Pr\'ecis\'ement, 
soient  $(w,w')\in {\cal W}_{P}$.   Consid\'erons le terme
$$X_{w}(m,k,k',\lambda)=(R_{P_{w}\vert w(S)}((w\sigma)_{w\lambda})\circ \gamma(w)\circ \sigma_{\lambda}(k)v,R_{\underline{P}_{w}\vert w(S)}((w\sigma)_{w\lambda})\circ\gamma(w)\circ\sigma_{\lambda}(mk')u)^M$$
qui intervient dans $\omega_{P,w,w'}(m,k,k',\lambda)$. On a vu dans la preuve pr\'ec\'edente que l'on pouvait  aussi l'\'ecrire
$$\delta_{P}(m)^{-1/2}(v_{w}(k,\lambda),\rho_{P,w\lambda}(m)u_{w}(k',\lambda)),$$
 o\`u $\rho_{P,w\lambda}=Ind_{w(S)\cap M}^M((w\sigma)_{w\lambda})$, $u_{w}(k',\lambda)$ est la restriction \`a $K\cap M(F)$ de 
$$R_{\underline{P}_{w}\vert w(S)}((w\sigma)_{w\lambda})\circ\gamma(w)\circ\sigma_{\lambda}(k')u,$$
et $v_{w}(k,\lambda)$ est la restriction \`a $K\cap M(F)$ de
$$R_{P_{w}\vert w(S)}((w\sigma)_{w\lambda})\circ \gamma(w)\circ \sigma_{\lambda}(k)v.$$
Notons $Y_{Q,w}(m,k,k',\lambda)$ le  terme constant faible relatif \`a $Q$ de
$$(v_{w}(k,\lambda),\rho_{P,w\lambda}(m)u_{w}(k',\lambda)),$$
 multipli\'e par $\delta_{Q}(m)^{-1/2}$. Remarquons que les \'el\'ements $u_{w}(k',\lambda)$ et $v_{w}(k,\lambda)$ restent dans des espaces de dimension finie et, quand on les \'ecrit dans une base de ces espaces, leurs coefficients sont des fonctions born\'ees de $k,k'$ et de $\lambda$. D'apr\`es (1), on peut appliquer la proposition 1.12 o\`u l'on remplace $G$ par $M$. Il existe donc $c>0$ et une fonction $C_{1}$ sur $i{\cal A}_{M_{disc},F}^*$, lisse, \`a croissance mod\'er\'ee et \`a valeurs positives, de sorte que l'on ait la majoration
 $$\vert X_{w}(m,k,k',\lambda)-Y_{Q,w}(m,k,k',\lambda)\vert \leq C_{1}(\lambda)\delta_{Q}(m)^{-1/2}\Xi^{L}(m)e^{-c\vert H_{0}(m)\vert }$$
pour tous $k$, $k'$, $\lambda$ et tout $m$ dans le domaine d\'ecrit ci-dessus. Notons que la condition $\tau_{Q}^P(H_{0}(m)-\epsilon T)=1$ entra\^{\i}ne une minoration $\vert H_{0}(m)\vert >> \vert T\vert $ (sauf dans le cas o\`u $P=Q$, mais alors le lemme est tautologique). Pour $m\in M_{0}(F)^{\geq}$, on a aussi 
$$\Xi^L(m)<<(1+\vert H_{0}(m)\vert )^{D_{1}}\delta_{P_{0}}^Q(m)^{-1/2}<<\vert T\vert ^{\nu D_{1}}\delta_{P_{0}}^Q(m)^{-1/2}$$
 pour un entier $D_{1}$ convenable. Quitte \`a r\'eduire $c$, on a donc
$$(2)\qquad \vert X_{w}(m,k,k',\lambda)-Y_{Q,w}(m,k,k',\lambda)\vert \leq C_{1}(\lambda) \delta_{P_{0}}(m)^{-1/2}e^{-c\vert T\vert }.$$
D'apr\`es 1.2(2) et 1.12(3), on a aussi
$$(3)\qquad \vert X_{w}(m,k,k',\lambda)\vert<<(1+\vert H_{0}(m)\vert )^{D_{2}}\delta_{P_{0}}(m)^{-1/2}<<\vert T\vert ^{\nu D_{2}}\delta_{P_{0}}(m)^{-1/2} ,$$
pour un autre entier $D_{2}$. Par diff\'erence, on en d\'eduit
$$(4) \qquad \vert Y_{Q,w}(m,k,k',\lambda)\vert <<C_{2}(\lambda)\vert T \vert ^{\nu D_{2}}\delta_{P_{0}}(m)^{-1/2},$$
pour une autre fonction $C_{2}$. 

On traite de la m\^eme fa\c{c}on le terme
$$X_{w'}(m,k,k',\lambda)=$$
$$(R_{\underline{P}_{w'}\vert w'(S')}((w'\sigma')_{w'\lambda'})\circ\gamma(w')\circ \pi'_{\lambda'}(mk')v',R_{P_{w'}\vert w'(S')}((w'\sigma')_{w'\lambda'})\circ\gamma(w')\circ \pi'_{\lambda'}(k)u')^M$$
qui intervient dans
$\omega_{P,w,w'}(m,k,k',\lambda)$. On d\'efinit $Y_{Q,w'}(m,k,k',\lambda)$, qui v\'erifie des majorations analogues \`a celles ci-dessus.

Rappelons l'\'egalit\'e
$$\omega_{P,w,w'}(m)=\int_{K\times K}\int_{i{\cal A}_{M_{disc},F}^*}m^G(\sigma_{\lambda})B(\lambda){\bf r}_{w,w'}(\sigma_{\lambda})X_{w}(m,k,k',\lambda)$$
$$X_{w'}(m,k,k',\lambda)\omega(k^{-1}mk')\,d\lambda\,dk\,dk'.$$
 Posons
$$\omega_{Q,w,w'}(m)=\gamma(M\vert L)^2\int_{K\times K}\int_{i{\cal A}_{M_{disc},F}^*}m^G(\sigma_{\lambda})B(\lambda){\bf r}_{w,w'}(\sigma_{\lambda})Y_{Q,w}(m,k,k',\lambda)$$
$$Y_{Q,w'}(m,k,k',\lambda)\omega(k^{-1}mk')\,d\lambda\,dk\,dk'.$$
En utilisant les majorations pr\'ec\'edentes, on obtient
$$\vert \omega_{P,w,w'}(m)-\gamma(M\vert L)^{-2}\omega_{Q,w,w'}(m)\vert \delta_{P_{0}}(m)\leq  e^{-c\vert T\vert }\int_{i{\cal A}_{M_{disc},F}^*}\vert C(\lambda)B(\lambda)m^G(\sigma_{\lambda}){\bf r}_{w,w'}(\sigma_{\lambda})\vert \,d\lambda$$
pour une certaine fonction $C$ \`a croissance mod\'er\'ee. L'int\'egrale ci-dessus est convergente, donc les fonctions $\phi^Q(H_{0}(m)-\epsilon T)\tau_{Q}^P(H_{0}(m)-\epsilon T)\omega_{P,w,w'}(m)\delta_{P_{0}}(m)$ et $\gamma(M\vert L)^{-2}\phi^Q(H_{0}(m)-\epsilon T)\tau_{Q}^P(H_{0}(m)-\epsilon T)\omega_{Q,w,w'}(m)\delta_{P_{0}}(m)$ sont \'equivalentes.

On calcule $Y_{Q,w}(m,k,k',\lambda)$ en utilisant les d\'efinitions de 1.12. On doit y remplacer $G$ par $M$, $P$ par $w(S)\cap M$, $Q$ par $Q\cap M$ et $\sigma_{\lambda}$ par $(w\sigma)_{w\lambda}$. Notons que, pour $s\in W^M(L\vert w(S)\cap M)$, les op\'erateurs $J^M_{(Q\cap M)_{s}\vert sw(S)\cap M}((sw\sigma)_{sw\lambda})$ et $J^M_{\underline{(Q\cap M)}_{s}\vert sw(S)\cap M}((sw\sigma)_{sw\lambda})$ qui interviennent  en 1.12 d\'efinissent par induction les op\'erateurs $J_{Q_{sw}\vert s(P_{w})}((sw\sigma)_{sw\lambda})$ et $J_{\underline{Q}_{sw}\vert s(\underline{P}_{w})}((sw\sigma)_{sw\lambda})$, avec les m\^emes d\'efinitions qu'en 3.4. En r\'etablissant les op\'erateurs d'entrelacement non normalis\'es dans la d\'efinition de $X_{w}(m,k,k',\lambda)$ et en utilisant les propri\'et\'es habituelles de ces op\'erateurs, on obtient l'\'egalit\'e
$$Y_{Q,w}(m,k,k',\lambda)=\gamma(M\vert L)^{-1}r_{ \underline{P}_{w}\vert w(S)}((w\sigma)_{w\lambda})^{-1}r_{w(S)\vert P_{w}}((w\sigma)_{w\lambda})^{-1}$$
$$\sum_{s\in W^M(L\vert w(S)\cap M)}( J_{Q_{sw}\vert sw(S)}((sw\sigma)_{sw\lambda})\circ\gamma(sw)\circ\pi_{\lambda}(k)v,J_{\underline{Q}_{sw}\vert sw(S)}((sw\sigma)_{sw\lambda})\circ\gamma(sw)\circ\pi_{\lambda}(mk')u)^L.$$
Remarquons que l'appication $s\mapsto sw$ est une bijection de $W^M(L\vert w(S)\cap M)$ sur l'ensemble des $s\in W^G(L\vert S)$ tels que $s\in W^Mw$. On calcule de m\^eme  $Y_{Q,w'}(m,k,k',\lambda)$.

Notons ${\cal W}_{Q,w,w'}^P$ l'ensemble des $(s,s')\in W^G(L\vert S)\times W^G(L\vert  S')$ tels que $s\in W^Mw$ et $s'\in W^Mw'$. Les calculs pr\'ec\'edents conduisent \`a l'\'egalit\'e
$$\omega_{Q,w,w'}(m)=\sum_{(s,s')\in {\cal W}_{Q,w,w'}^P}\omega_{Q,s,s'}(m).$$
Le terme $\Omega_{P}(m)$ est la somme des $\omega_{P,w,w'}(m)$ sur les $(w,w')\in {\cal W}_{P}$. On a $\cup_{(w,w')\in {\cal W}_{P}}{\cal W}_{Q,w,w'}={\cal W}_{Q}^P$.  Posons 
$$\Omega_{Q}^P(m)=\sum_{(s,s')\in {\cal W}_{Q}^P}\omega_{Q,s,s'}(m).$$
On obtient que les fonctions $\phi^Q(H_{0}(m)-\epsilon T)\tau_{Q}^P(H_{0}(m)-\epsilon T)\Omega_{P}(m)\delta_{P_{0}}(m)$ et $\gamma(M\vert L)^{-2}\phi^Q(H_{0}(m)-\epsilon T)\tau_{Q}^P(H_{0}(m)-\epsilon T)\Omega_{Q}^P(m)\delta_{P_{0}}(m)$ sont \'equivalentes. Pour obtenir le (i) de l'\'enonc\'e, il reste \`a prouver que la derni\`ere fonction est \'equivalente \`a $\gamma(M\vert L)^{-2}\phi^Q(H_{0}(m)-\epsilon T)\tau_{Q}^P(H_{0}(m)-\epsilon T)\Omega_{Q}(m)\delta_{P_{0}}(m)$. Or il  suffit pour cela de prouver le (ii) de l'\'enonc\'e.

 Soient maintenant $w$, $w'$ comme en (ii). En utilisant le lemme 3.6 et la majoration $\Xi^L(m)^2<<(1+\vert H_{0}(m)\vert )^D\delta_{P_{0}}^Q(m)^{-1}$ pour $m\in M_{0}(F)^{\geq,Q}$, cf. 1.2(2), il suffit pour prouver (ii) que, pour $m$ dans le domaine qui nous int\'eresse, on ait une minoration
$$N_{w,w'}^L(H_{0}(m))>>\vert T\vert .$$
Mais soit $\tilde{P}'=\tilde{M}'U_{P'}$ le plus petit espace parabolique standard tel que $Q\subset P'$ et $(w,w')\in {\cal W}_{Q}^{P'}$. On a $P'\subset P$. La condition $\tau_{Q}^P(H_{0}(m)-\epsilon T)=1$ entra\^{\i}ne $\tau_{Q}^{P'}(H_{0}(m)-\epsilon T)=1$. On peut appliquer le lemme 3.5(iii) qui nous fournit la minoration
$$N_{w,w'}^L(H_{0}(m))>>\epsilon\vert T\vert+\vert (H_{0}(m)-T)_{L}^{M'}\vert $$
plus forte que celle dont on a besoin. $\square$

\bigskip

\subsection{Un deuxi\`eme lemme d'\'equivalence}

 \ass{Lemme}{Soit $\tilde{P}=\tilde{M}U_{P}$ un espace parabolique contenant $\tilde{P}_{0}$, soit $Q=LU_{Q}$ un sous-groupe parabolique tel que $P_{0}\subset Q\subset P$ et soit $\epsilon>0$. 

(i) Les fonctions $\phi^Q(H_{0}(m)-\epsilon T)\tau_{Q}^P(H_{0}(m)-\epsilon T)\Omega_{P}(m)\delta_{P}(m)D_{0}^M(m)$ et $\phi^Q(H_{0}(m)-\epsilon T)\tau_{Q}^P(H_{0}(m)-\epsilon T)\Omega_{Q}(m)\delta_{Q}(m)D_{0}^L(m)$ sont \'equivalentes. 

(ii) Soient $w\in W^G(L\vert  S)$ et $w'\in W^G(L\vert  S')$. Supposons  $(w,w')\in {\cal W}_{Q}^P-{\cal W}_{Q}$.  Alors la fonction  $\phi^Q(H_{0}(m)-\epsilon T)\tau_{Q}^P(H_{0}(m)-\epsilon T)\omega_{Q,w,w'}(m)\delta_{Q}(m)D_{0}^L(m)$ est \'equivalente \`a $0$.}

Preuve. On a une majoration $\delta_{P}(m)D_{0}^M(m)<<\delta_{P_{0}}(m)$ pour tout $m\in M_{0}(F)^{\geq}$, cf. 1.2(1). On peut donc multiplier les deux fonctions du (i) l'\'enonc\'e pr\'ec\'edent par la fonction born\'ee $\delta_{P_{0}}(m)^{-1}\delta_{P}(m)D_{0}^M(m)$, elles restent \'equivalentes. On obtient que la premi\`ere fonction du (i) du pr\'esent \'enonc\'e est \'equivalente \`a la fonction 
$$(1) \qquad \gamma(M\vert L)^{-2} \phi^Q(H_{0}(m)-\epsilon T)\tau_{Q}^P(H_{0}(m)-\epsilon T)\Omega_{Q}(m)\delta_{P}(m)D_{0}^M(m).$$
Comme on l'a vu dans la preuve pr\'ec\'edente, une fois fix\'e un r\'eel $\nu$, on a une minoration $<\alpha,H_{0}(m)>>>\nu'\vert H_{0}(m)\vert $ pour les $m$ intervenant. D'apr\`es [A1] lemme 1.1 et en tenant compte de la remarque 1.2, il existe $c>0$ tel que l'on ait une majoration
$$\vert \gamma(M\vert L)^{-2}\delta_{P}(m)D_{0}^M(m)-\delta_{Q}(m)D_{0}^L(m)\vert <<\delta_{P_{0}}(m)e^{-c\vert H_{0}(m)\vert }.$$
De nouveau, on a vu dans la preuve pr\'ec\'edente que l'on avait une minoration $\vert H_{0}(m)\vert >>\vert T\vert$, sauf si $P=Q$ auquel cas le lemme est tautologique.  On peut remplacer la majoration ci-dessus par
$$\vert\gamma(M\vert L)^{-2} \delta_{P}(m)D_{0}^M(m)-\delta_{Q}(m)D_{0}^L(m)\vert <<\delta_{P_{0}}(m)e^{-c\vert T\vert }$$
pour un autre $c>0$. Il r\'esulte du corollaire 3.6 et de 1.2(2) que l'on a une majoration
$$\vert \Omega_{Q}(m)\vert<< (1+\vert H_{0}(m)\vert )^D\delta_{P_{0}}(m)^{-1}$$
pour un entier $D$ convenable. De ces deux derni\`eres majorations r\'esulte que la fonction (1) est \'equivalente \`a la deuxi\`eme fonction du (i) de l'\'enonc\'e. Cela d\'emontre ce (i). Le (ii) r\'esulte du (ii) du lemme pr\'ec\'edent par multiplication par la fonction born\'ee $\delta_{P_{0}}(m)^{-1}\delta_{P}(m)D_{0}^M(m)$.  $\square$ 

\bigskip

\subsection{Un troisi\`eme lemme d'\'equivalence}

\ass{Lemme}{Soient $Q=LU_{Q}$ un sous-groupe parabolique standard, $\tilde{P}=\tilde{M}U_{P}$ et $\tilde{P}'=\tilde{M}'U_{P'}$ deux espaces paraboliques standard et soit $\epsilon>0$. Notons $\tilde{P}_{-}=\tilde{M}_{-}U_{P_{-}}$ le plus grand espace parabolique standard tel que $P_{-}\subset Q$.   On suppose $Q\subset P$ et $P_{-}\subset P'\subset P$. Alors les fonctions 
$$\phi^Q(H_{0}(m-\epsilon T)\tau_{Q}^P(H_{0}(m)-\epsilon T)\Omega_{P'}(m)\delta_{P'}(m)D_{0}^{M'}(m)$$
 et  
 $$\phi^Q(H_{0}(m-\epsilon T)\tau_{Q}^P(H_{0}(m)-\epsilon T)\Omega_{P_{-}}(m)\delta_{P_{-}}(m)D_{0}^{M_{-}}(m)$$
  sont \'equivalentes.}
  
  Preuve. Fixons $\epsilon'>0$ que l'on pr\'ecisera plus tard. On utilise l'\'egalit\'e
  $$\sum_{Q'; P_{0}\subset Q'\subset P'}\phi^{Q'}(H_{0}(m)-\epsilon' T)\tau_{Q'}^{P'}(H_{0}(m)-\epsilon'T)=1$$
  cf. 1.3(2).  On peut fixer $Q'=L'U_{Q'}$ et prouver que les produits de nos fonctions avec $\phi^{Q'}(H_{0}(m)-\epsilon' T)\tau_{Q'}^{P'}(H_{0}(m)-\epsilon'T)$ sont  \'equivalentes. D'apr\`es le lemme pr\'ec\'edent, les fonctions
  $$\phi^{Q'}(H_{0}(m)-\epsilon' T)\tau_{Q'}^{P'}(H_{0}(m)-\epsilon'T)\Omega_{P'}(m)\delta_{P'}(m)D_{0}^{M'}(m)$$
  et
  $$\phi^{Q'}(H_{0}(m)-\epsilon' T)\tau_{Q'}^{P'}(H_{0}(m)-\epsilon'T)\Omega_{Q'}(m)\delta_{Q'}(m)D_{0}^{L'}(m)$$
  sont \'equivalentes.  La premi\`ere fonction de l'\'enonc\'e, multipli\'ee par   $\phi^{Q'}(H_{0}(m)-\epsilon' T)\tau_{Q'}^{P'}(H_{0}(m)-\epsilon'T)$, est alors \'equivalente \`a
  $$(1) \qquad \phi^Q(H_{0}(m-\epsilon T)\tau_{Q}^P(H_{0}(m)-\epsilon T) \phi^{Q'}(H_{0}(m)-\epsilon' T)$$
  $$\tau_{Q'}^{P'}(H_{0}(m)-\epsilon'T)\Omega_{Q'}(m)\delta_{Q'}(m)D_{0}^{L'}(m).$$ 
  Supposons que $Q'\subset P_{-}$. Alors l'\'egalit\'e $\phi^{Q'}(H_{0}(m)-\epsilon' T)\tau_{Q'}^{P'}(H_{0}(m)-\epsilon'T)=1$ entra\^{\i}ne $\phi^{Q'}(H_{0}(m)-\epsilon' T)\tau_{Q'}^{P_{-}}(H_{0}(m)-\epsilon'T)=1$. On peut appliquer le raisonnement ci-dessus en rempla\c{c}ant $P'$ par $P_{-}$. On obtient que la deuxi\`eme fonction de l'\'enonc\'e, multipli\'ee par $\phi^{Q'}(H_{0}(m)-\epsilon' T)\tau_{Q'}^{P'}(H_{0}(m)-\epsilon'T)$, est \'equivalente \`a (1), d'o\`u la conclusion cherch\'ee.
  
   Supposons maintenant que $Q'\not\subset P_{-}$. Il nous suffit de prouver que la fonction (1) est \'equivalente \`a $0$ et que  la deuxi\`eme fonction de l'\'enonc\'e, multipli\'ee par $\phi^{Q'}(H_{0}(m)-\epsilon' T)\tau_{Q'}^{P'}(H_{0}(m)-\epsilon'T)$, est elle-aussi \'equivalente \`a $0$. En utilisant le corollaire 3.6 et les majorations maintenant famili\`eres portant sur les fonctions $D_{0}$ et $\Xi$, on voit qu'il suffit de prouver que, pour $m\in M_{0}(F)^{\geq}$ tel que 
  $$\phi^Q(H_{0}(m-\epsilon T)\tau_{Q}^P(H_{0}(m)-\epsilon T) \phi^{Q'}(H_{0}(m)-\epsilon' T)\tau_{Q'}^{P'}(H_{0}(m)-\epsilon'T)=1,$$
  on a des minorations
  $$(2)\qquad N^{L'}(H_{0}(m))>>\vert T\vert ,\,\,N^{M_{-}}(H_{0}(m))>>\vert T\vert .$$
 
 Remarquons que la condition $\phi^{Q'}(H_{0}(m)-\epsilon' T)=1$ entra\^{\i}ne $<\alpha,H_{0}(m)>\leq c_{1}\epsilon'\vert T\vert $ pour un certain $c_{1}>0$ et pour tout $\alpha\in \Delta_{0}^{Q'}$, tandis que la condition $\phi^Q(H_{0}(m)-\epsilon T)\tau_{Q}^P(H_{0}(m)-\epsilon T)=1$ entra\^{\i}ne $<\alpha,H_{0}(m)>> c_{2}\vert T\vert $ pour un certain $c_{2}>0$ et pour tout $\alpha\in \Delta_{0}^P-\Delta_{0}^{Q}$. En supposant $c_{1}\epsilon'<c_{2}$, nos fonctions sont nulles sauf si $Q'\subset Q$. On suppose donc $Q'\subset Q$. L'hypoth\`ese $Q'\not\subset P_{-}$ signifie qu'il existe $\alpha\in \Delta_{0}^{Q'}$ tel que $\alpha\not\in \Delta_{0}^{P_{-}}$. Fixons un tel $\alpha$, qui appartient \`a $\Delta_{0}^Q$ par l'hypoth\`ese $Q'\subset Q$ que l'on vient de poser. Par d\'efinition de $P_{-}$, $\Delta_{0}^{P_{-}}$ est l'ensemble des \'el\'ements de $\Delta_{0}$ dont les images par les puissances de $\theta$ appartiennent toutes \`a $\Delta_{0}^Q$. Puisque $\alpha\not\in \Delta_{0}^{P_{-}}$, on peut fixer un entier $k>0$ tel que $\theta^{-k}\alpha\not\in \Delta_{0}^Q$.
 
 Notons $\tilde{P}''=\tilde{M}''U_{P''}$ le plus petit espace parabolique tel que $Q'\subset P''$. Il est inclus dans $\tilde{P}'$ et l'hypoth\`ese $\phi^{Q'}(H_{0}(m)-\epsilon' T)\tau_{Q'}^{P'}(H_{0}(m)-\epsilon'T)=1$ entra\^{\i}ne $\phi^{Q'}(H_{0}(m)-\epsilon' T)\tau_{Q'}^{P''}(H_{0}(m)-\epsilon'T)=1$. Appliquons le lemme 3.5(i). On obtient
 $$(3)\qquad N^{L'}(H_{0}(m))>>\vert (H_{0}(m)-\epsilon'T)_{L'}^{M''}\vert .$$
 Des hypoth\`eses $m\in M_{0}(F)^{\geq}$ et $ \phi^{Q'}(H_{0}(m)-\epsilon' T)=1$ r\'esulte une majoration
 $$\vert (H_{0}(m)-\epsilon'T)_{L'}^{M''}\vert\geq \vert H_{0}(m)^{M''}\vert -c_{3}\epsilon'\vert T\vert $$
 pour une constante $c_{3}>0$ convenable. Puisque $\alpha\in \Delta_{0}^{Q'}$, on a $\theta^{-k}\alpha\in \Delta_{0}^{P''}$. Donc
 $$<\theta^{-k}\alpha,H_{0}(m)^{M''}>=<\theta^{-k}\alpha,H_{0}(m)>.$$
  Puisque $\theta^{-k}\alpha\in \Delta_{0}^P-\Delta_{0}^Q$, l'hypoth\`ese $\tau_{Q}^P(H_{0}(m)-\epsilon T)=1$ entra\^{\i}ne 
  $$<\theta^{-k}\alpha,H_{0}(m)>> <\theta^{-k}\alpha, \epsilon T>.$$ On en  d\'eduit une majoration $\vert H_{0}(m)^{M''}\vert \geq c_{4}\epsilon \vert T\vert $ pour une constante $c_{4}>0$ convenable. On obtient
$$ \vert (H_{0}(m)-\epsilon'T)_{L'}^{M''}\vert\geq (c_{4}\epsilon-c_{3}\epsilon')\vert T\vert .$$
En supposant $c_{3}\epsilon'<c_{4}\epsilon$, on en d\'eduit 
$$\vert (H_{0}(m)-\epsilon'T)_{L'}^{M''}\vert>> \vert T\vert .$$
Jointe \`a (3), cette minoration entra\^{\i}ne la premi\`ere relation de (2).

  Soient maintenant $H',H''\in {\cal C}^{M_{-}}(H_{0}(m))$. Parce que $P_{-}$ est stable par $\theta$, on a $(\theta H')_{M_{-}}=\theta(H'_{M_{-}})$. On a aussi $H'_{M_{-}}=H''_{M_{-}}=H_{M_{-}}(m)$. D'o\`u
  $$\vert \theta H'-H''\vert \geq\vert (\theta H')_{M_{-}}-H''_{M_{-}}\vert =\vert \theta H_{M_{-}}(m)-H_{M_{-}}(m)\vert\geq\vert \theta(H_{M_{-}}(m)^M)-H_{M_{-}}(m)^M\vert  .$$
  Parce que $\phi^Q(H_{0}(m-\epsilon T)\tau_{Q}^P(H_{0}(m)-\epsilon T)=1$, on peut \'ecrire $H_{0}(m)^M=\epsilon T^M+X-Y$, avec  $X=\sum_{\beta\in \Delta_{0}^P-\Delta_{0}^{Q}}x_{\beta}\check{\varpi}_{\beta}^M$, $Y=\sum_{\beta\in \Delta_{0}^Q}y_{\beta}\check{\beta}$, les coefficients $x_{\beta}$ et $y_{\beta}$ \'etant positifs ou nuls. Donc $H_{M_{-}}(m)^M=\epsilon T_{M_{-}}^M+X-Y_{M_{-}}$, avec $Y_{M_{-}}=\sum_{\beta\in \Delta_{0}^L-\Delta_{0}^M}y_{\beta}\check{\beta}_{M}$. Remarquons que $T_{M_{-}}^M$ est fixe par $\theta$ puisque $T$ l'est. D'o\`u
  $$\vert \theta H_{M_{-}}(m)^M-H_{M_{-}}(m)^M\vert^2=\vert \theta X-X+Y_{M_{-}}-\theta(Y_{M_{_{}}})\vert ^2$$
  $$=\vert \theta X-X\vert ^2+\vert\theta(Y_{M_{-}})-Y_{M_{-}}\vert ^2+2(\theta X-X,Y_{M_{-}}-\theta(Y_{M_{-}})).$$
  Les \'el\'ements $X$ et $Y_{M_{-}}$ sont orthogonaux, et les \'el\'ements $\theta X$ et $\theta(Y_{M_{-}})$ aussi. L'\'el\'ement $\theta X$ est combinaison lin\'eaire \`a coefficients positifs ou nuls de $\check{\varpi}_{\beta}$ pour $\beta\in \Delta_{0}^P-\Delta_{0}^{P_{-}}$. Donc $(\theta X,Y_{M_{-}})\geq0$. Pour une raison analogue, $(X,\theta(Y_{M_{-}}))\geq0$. On obtient alors
  $$\vert \theta H_{M_{-}}(m)-H_{M_{-}}(m)\vert \geq\vert \theta(Y_{M_{-}})-Y_{M_{-}}\vert .$$
  Puisque $\alpha\in \Delta_{0}^Q$, on a 
  $$<\alpha,H_{0}(m)>=\epsilon<\alpha,T>-2y_{\alpha}-\sum_{\beta\in \Delta_{0}^Q, \beta\not=\alpha}y_{\beta}<\alpha,\check{\beta}>\geq c_{5}\epsilon\vert T\vert -2y_{\alpha}$$
  pour un certain $c_{5}>0$, puisque les $<\alpha,\check{\beta}>$ sont n\'egatifs ou nuls. Puisque $\alpha\in \Delta_{0}^{Q'}$, on a aussi $<\alpha,H_{0}(m)>\leq c_{1}\epsilon'\vert T\vert $. D'o\`u
  $$y_{\alpha}\geq \frac{1}{2}(c_{5}\epsilon-c_{1}\epsilon')\vert T\vert.$$
  En supposant $c_{1}\epsilon' <c_{5}\epsilon$, on obtient une minoration $y_{\alpha}>>\vert T\vert $. Puisque $\theta^{-k}\alpha\not\in \Delta_{0}^Q$, on a $<\theta^{-k}\varpi_{\alpha},Y_{M_{-}}>=0$, d'o\`u $<\varpi_{\alpha},\theta^k(Y_{M_{-}})>=0$. Alors $<\varpi_{\alpha},Y_{M_{-}}-\theta^k(Y_{M_{-}})>=y_{\alpha}>>\vert T\vert $, d'o\`u $\vert Y_{M_{-}}-\theta^k(Y_{M_{-}})\vert >>\vert T\vert $. L'\'egalit\'e $1-\theta^k=(1-\theta)(1+\theta+...+\theta^{k-1})$ et le fait que $\theta$ soit une isom\'etrie entra\^{\i}ne $\vert Y_{M_{-}}-\theta^k(Y_{M_{-}})\vert\leq k\vert Y_{M_{-}}-\theta(Y_{M_{-}})\vert $. D'o\`u  $\vert Y_{M_{-}}-\theta(Y_{M_{-}})\vert >>\vert T\vert $. En rassemblant nos calculs, on obtient la minoration
  $$\vert \theta H'-H''\vert >>\vert T\vert ,$$
  qui d\'emontre la seconde assertion de (2) et le lemme. $\square$

\bigskip

\subsection{Une fonction ${\cal E}^T$}

Pour $m\in M_{0}(F)$, posons
$${\cal E}^T(m)=\sum_{\tilde{P}=\tilde{M}U_{P}; \tilde{P}_{0}\subset \tilde{P}}(-1)^{dim(\mathfrak{a}_{\tilde{M}})-dim(\mathfrak{a}_{\tilde{G}})}\Omega_{P}(m)\hat{\tau}_{\tilde{P}}(H_{0}(m)-T)\delta_{P}(m)D_{0}^M(m).$$

\ass{Proposition}{Les fonctions ${\cal E}^T(m)$ et $\Omega_{G}(m)\phi^{\tilde{G}}(H_{0}(m)-T)D_{0}^G(m)$ sont \'equivalentes.}

Preuve. Fixons $\epsilon>0$ que l'on pr\'ecisera plus tard. Pour $m\in M_{0}(F)^{\geq}$, on a l'\'egalit\'e
$$\sum_{Q; P_{0}\subset Q}\phi^Q(H_{0}(m)-\epsilon T)\tau_{Q}(H_{0}(m)-\epsilon T)=1,$$
cf. 1.3(2). On peut donc fixer $Q$ et montrer que les fonctions $\phi^Q(H_{0}(m)-\epsilon T)\tau_{Q}(H_{0}(m)-\epsilon T){\cal E}^T(m)$ et $\phi^Q(H_{0}(m)-\epsilon T)\tau_{Q}(H_{0}(m)-\epsilon T)\phi^{\tilde{ G}}(H_{0}(m)-T)\Omega_{G}(m)D_{0}^G(m)$ sont \'equivalentes. Soit $\tilde{P}_{-}=\tilde{M}_{-}U_{P_{-}}$ le plus grand espace parabolique tel que $P_{0}\subset P \subset Q$. On introduit la fonction $\phi_{\tilde{P}_{-}}^{\tilde{G}}$, cf. 2.2. On d\'ecompose encore le probl\`eme en deux.   On va prouver:

(1) les fonctions $\phi_{\tilde{P}_{-}}^{\tilde{G}}(H_{0}(m)-T)\phi^Q(H_{0}(m)-\epsilon T)\tau_{Q}(H_{0}(m)-\epsilon T){\cal E}^T(m)$ et $\phi_{\tilde{P}_{-}}^{\tilde{G}}(H_{0}(m)-T)\phi^Q(H_{0}(m)-\epsilon T)\tau_{Q}(H_{0}(m)-\epsilon T)\phi^{\tilde{ G}}(H_{0}(m)-T)\Omega_{G}(m)D_{0}^G(m)$ sont \'egales sur $M_{0}(F)^{\geq}$;

(2) les fonctions $(1-\phi_{\tilde{P}_{-}}^{\tilde{G}}(H_{0}(m)-T))\phi^Q(H_{0}(m)-\epsilon T)\tau_{Q}(H_{0}(m)-\epsilon T){\cal E}^T(m)$ et $(1-\phi_{\tilde{P}_{-}}^{\tilde{G}}(H_{0}(m)-T))\phi^Q(H_{0}(m)-\epsilon T)\tau_{Q}(H_{0}(m)-\epsilon T)\phi^{\tilde{ G}}(H_{0}(m)-T)\Omega_{G}(m)D_{0}^G(m)$ sont toutes deux \'equivalentes \`a $0$.

Traitons (1). On se limite \`a consid\'erer des $m\in M_{0}(F)^{\geq}$ tels que $\phi_{\tilde{P}_{-}}^{\tilde{G}}(H_{0}(m)-T)\phi^Q(H_{0}(m)-\epsilon T)\tau_{Q}(H_{0}(m)-\epsilon T)=1$. Pour $\alpha\in \Delta^{P_{-}}_{0}$, on a une \'egalit\'e
$$\varpi_{\alpha}=\varpi_{\alpha}^{M_{-}}+\sum_{\beta\in \Delta_{0}^G-\Delta_{0}^{P_{-}}}x_{\beta}\varpi_{\beta}.$$
Parce que les poids fondamentaux forment une base aigu\" e de ${\cal A}_{0}^*$, les coefficients $x_{\beta}$ sont positifs ou nuls. On en d\'eduit une formule analogue
$$\varpi_{\tilde{\alpha}}=\varpi_{\tilde{\alpha}}^{M_{-}}+\sum_{\tilde{\beta}\in \Delta_{0}^{\tilde{G}}-\Delta_{0}^{\tilde{P}_{-}}}x_{\tilde{\beta}}\varpi_{\tilde{\beta}}.$$
Puisque $\phi_{\tilde{P}_{-}}^{\tilde{G}}(H_{0}(m)-T)=1$, on a $<\varpi_{\tilde{\beta}},H_{0}(m)-T>\leq 0$ pour tous les $\tilde{\beta}$ intervenant, donc 
$$(3) \qquad <\varpi_{\tilde{\alpha}},H_{0}(m)-T>\leq <\varpi_{\tilde{\alpha}}^{M_{-}},H_{0}(m)-T>.$$ 
Les hypoth\`eses $m\in M_{0}(F)^{\geq}$ et $\phi^Q(H_{0}(m)-\epsilon T)=1$ entra\^{\i}nent une majoration $\vert H_{0}(m)^L\vert\leq c\epsilon\vert T\vert $ pour un certain $c>0$, a fortiori
$\vert H_{0}(m)^{M_{-}}\vert \leq c\epsilon\vert T\vert $. D'o\`u $\vert <\varpi_{\tilde{\alpha}}^{M_{-}},H_{0}(m)>\vert \leq c_{\tilde{\alpha}} \epsilon \vert T\vert$ pour un certain $c_{\tilde{\alpha}}>0$. On pr\'ecise $\epsilon$ en  imposant que pour tout $T$, on ait l'in\'egalit\'e $c_{\tilde{\alpha}}\epsilon\vert T\vert \leq <\varpi_{\tilde{\alpha}}^{M_{-}},T>$ ( pour tout $\tilde{\alpha}$). Gr\^ace \`a (3), cela entra\^{\i}ne
$$<\varpi_{\tilde{\alpha}},H_{0}(m)-T>\leq 0.$$
On avait suppos\'e $\tilde{\alpha}\in \Delta_{0}^{\tilde{P}_{-}}$, mais cette relation est aussi vraie pour $\tilde{\alpha}\in \Delta_{0}^{\tilde{G}}-\Delta_{0}^{\tilde{P}_{-}}$ gr\^ace \`a l'hypoth\`ese $\phi_{\tilde{P}_{-}}^{\tilde{G}}(H_{0}(m)-T)=1$. Mais alors on a $\hat{\tau}_{\tilde{P}}(H_{0}(m)-T)=0$ pour tout $\tilde{P}$ propre. Dans la somme d\'efinissant ${\cal E}^T(m)$, il ne reste que la contribution de $\tilde{P}=\tilde{G}$, qui est simplement $\Omega_{G}(m)D_{0}^G(m)$. Comme on vient de le voir, on a $\phi^{\tilde{ G}}(H_{0}(m)-T)=1$ donc $\Omega_{G}(m)=\phi^{\tilde{ G}}(H_{0}(m)-T)\Omega_{G}(m)$ et la conclusion de (1) s'ensuit.

Traitons (2). On peut supposer $P_{-}\not=G$, sinon $1-\phi_{\tilde{P}_{-}}^{\tilde{G}}(H_{0}(m)-T)=0$. Pour la deuxi\`eme fonction, l'assertion est \'evidente car $(1-\phi_{\tilde{P}_{-}}^{\tilde{G}}(H_{0}(m)-T))\phi^{\tilde{ G}}(H_{0}(m)-T)=0$. Pour tout $\tilde{\alpha}\in \Delta_{0}^{\tilde{G}}-\Delta_{0}^{\tilde{P}_{-}}$, notons simplement $\hat{\tau}_{\tilde{\alpha}}$ la fonction caract\'eristique de l'ensemble des $H\in {\cal A}_{0}$ tels que $\varpi_{\tilde{\alpha}}(H)>0$. Le support de $1-\phi_{\tilde{P}_{-}}^{\tilde{G}}(H_{0}(m)-T)$ est contenu dans la r\'eunion des supports des $\hat{\tau}_{\tilde{\alpha}}(H_{0}(m)-T)$ quand $\tilde{\alpha}$ d\'ecrit $ \Delta_{0}^{\tilde{G}}-\Delta_{0}^{\tilde{P}_{-}}$. On peut fixer $\tilde{\alpha}$ dans cet ensemble et se contenter de prouver que la fonction
$$\hat{\tau}_{\tilde{\alpha}}(H_{0}(m)-T)\phi^Q(H_{0}(m)-\epsilon T)\tau_{Q}(H_{0}(m)-\epsilon T){\cal E}^T(m)$$
est \'equivalente \`a $0$. Dans la somme d\'efinissant ${\cal E}^T(m)$, on regroupe les $\tilde{P}$ en paires $(\tilde{P}',\tilde{P})$ de sorte que $\tilde{\Delta}_{0}^{P}=\tilde{\Delta}_{0}^{P'}\sqcup \{\tilde{\alpha}\}$. Il nous suffit de prouver que, pour une telle paire, 

(4) la fonction 
$$\hat{\tau}_{\tilde{\alpha}}(H_{0}(m)-T)\phi^Q(H_{0}(m)-\epsilon T)\tau_{Q}(H_{0}(m)-\epsilon T)$$
$$\left(\Omega_{P}(m)\hat{\tau}_{\tilde{P}}(H_{0}(m)-T)\delta_{P}(m)D_{0}^{M}(m)-\Omega_{P'}(m)\hat{\tau}_{\tilde{P}'}(H_{0}(m)-T)\delta_{P'}(m)D_{0}^{M'}(m)\right)$$
est \'equivalente \`a $0$. 

On fixe $\epsilon'>0$ que l'on pr\'ecisera par la suite. En utilisant le d\'ecoupage d\'ej\`a maintes fois utilis\'e, on peut fixer un sous-groupe parabolique standard $Q'\subset P$ et montrer que  la fonction (4) multipli\'ee par $\phi^{Q'}(H_{0}(m)-\epsilon' T)\tau_{Q'}^P(H_{0}(m)-\epsilon' T)$ est \'equivalente \`a $0$. Par hypoth\`ese, $\tilde{\alpha}\not\in \Delta_{0}^{\tilde{P}_{-}}$. Par d\'efinition de $\tilde{P}_{-}$, il  existe donc une racine $\alpha\in \Delta_{0}-\Delta_{0}^Q$ telle que sa restriction \`a ${\cal A}_{\tilde{M}_{0}}$ soit $\tilde{\alpha}$. Fixons une telle racine. L'hypoth\`ese $\phi^Q(H_{0}(m)-\epsilon T)\tau_{Q}(H_{0}(m)-\epsilon T)=1$ entra\^{\i}ne une minoration $<\alpha,H_{0}(m)>>> \epsilon\vert T\vert $. Si $\alpha\in \Delta_{0}^{Q'}$, l'hypoth\`ese $\phi^{Q'}(H_{0}(m)-\epsilon' T)=1$ entra\^{\i}ne une majoration $<\alpha,H_{0}(m)><< \epsilon'\vert T\vert $. On impose que $\epsilon'$ soit assez petit pour que ces deux in\'egalit\'es soient contradictoires. Alors $\alpha\not\in \Delta_{0}^{Q'}$. Notons $\tilde{P}'_{-}=\tilde{M}'_{-}U_{P'_{-}}$ le plus grand espace parabolique standard tel que $P'_{-}\subset Q'$. Alors $\tilde{\alpha}\not\in \Delta^{\tilde{P}'_{-}}$. Cette relation, jointe aux inclusions $P'_{-}\subset Q'\subset P$ et \`a l'\'egalit\'e $\Delta^{\tilde{P}}=\Delta^{\tilde{P}'}\sqcup \{\tilde{\alpha}\}$, entra\^{\i}ne que $P'_{-}\subset P'$. On peut alors appliquer le lemme 3.9:  les fonctions
$$\phi^{Q'}(H_{0}(m)-\epsilon' T)\tau_{Q'}^P(H_{0}(m)-\epsilon' T)\Omega_{P}(m)\delta_{P}(m)D_{0}^M(m)$$
et
$$\phi^{Q'}(H_{0}(m)-\epsilon' T)\tau_{Q'}^P(H_{0}(m)-\epsilon' T)\Omega_{P'}(m)\delta_{P'}(m)D_{0}^{M'}(m)$$
sont toutes deux \'equivalentes \`a
$$\phi^{Q'}(H_{0}(m)-\epsilon' T)\tau_{Q'}^P(H_{0}(m)-\epsilon' T)\Omega_{P'_{-}}(m)\delta_{P'_{-}}(m)D_{0}^{M'_{-}}(m).$$
Leur diff\'erence est donc \'equivalente \`a $0$. Par ailleurs, on a l'\'egalit\'e
$$\hat{\tau}_{\tilde{\alpha}}(H_{0}(m)-T)\hat{\tau}_{\tilde{P}}(H_{0}(m)-T)=\hat{\tau}_{\tilde{\alpha}}(H_{0}(m)-T)\hat{\tau}_{\tilde{P}'}(H_{0}(m)-T)$$
d'apr\`es l'hypoth\`ese sur  $P$ et $P'$. Donc la fonction (4), multipli\'ee par $\phi^{Q'}(H_{0}(m)-\epsilon' T)\tau_{Q'}^P(H_{0}(m)-\epsilon' T)$, est \'egale \`a la diff\'erence des deux fonctions ci-dessus, multipi\'ee par une fonction born\'ee. Cela prouve (4) et la proposition. $\square$

\bigskip

\subsection{Contr\^ole de la partie centrale}
  On sait que ${\cal A}_{A_{G},F}$ est un sous-groupe d'indice fini dans ${\cal A}_{G,F}$ (ces groupes sont tous deux \'egaux \`a ${\cal A}_{G}$ si $F$ est archim\'edien). Fixons un ensemble de repr\'esentants ${\bf b}$ du groupe quotient. Pour tout sous-ensemble ${\cal X}\subset G(F)$, notons ${\cal X}^{{\bf b}}$ l'ensemble des $g\in X$ tels que $H_{G}(g)\in {\bf b}$.
 Notons $H_{G}^{\tilde{G}}$ la compos\'ee de $H_{G}$ et de la projection sur l'orthogonal de ${\cal A}_{\tilde{G}}$.

\ass{Lemme}{Soient $Q$ un sous-groupe parabolique standard, $w\in W^G(L\vert  S)$ et $w'\in W^G(L\vert  S')$. Il existe un entier $D$ et, quel que soit le r\'eel $r$, il existe $c>0$ tel que l'on ait la majoration
$$\vert \omega_{Q,w,w'}(am)\vert \delta_{Q}(am)D^L(am)\leq c(1+\vert H_{0}(m)\vert)^D(1+\vert H_{G}^{\tilde{G}} (a)\vert )^{-r}$$
pour tout $m\in M_{0}(F)^{\geq,{\bf b}}$ et tout $a\in A_{\tilde{G}}(F)\backslash A_{G}(F)$.}

Preuve.    Il r\'esulte des d\'efinitions que, pour tout $\lambda\in i{\cal A}_{M_{disc},F}^*$ et tous $k,k'\in K$, on a l'\'egalit\'e
$$\omega_{Q,w,w'}(am,k,k',\lambda)=e^{-<\lambda',H_{G}(a)>+<\lambda,H_{G}(a)>}\omega_{Q,w,w'}(m,k,k',\lambda)$$
$$=e^{<\lambda,H_{G}(a)-\theta H_{G}(a)>}\omega_{Q,w,w'}(m,k,k',\lambda).$$
Alors 
$$\omega_{Q,w,w'}(am)=\int_{K\times K}\int_{i{\cal A}_{M_{disc}}^*/i({\cal A}_{G}^*+{\cal A}_{M_{disc},F}^{\vee})}\omega_{Q,w,w'}(m,k,k',\lambda)m^G(\sigma_{\lambda})$$
$$B'(\lambda,a)\omega(k^{-1}amk')\,d\lambda\,dk\,dk',$$
o\`u 
$$B'(\lambda,a)=\int_{i{\cal A}_{G,F}^*}B(\lambda+\mu)e^{<\lambda+\mu,H_{G}(a)-\theta H_{G}(a)>}\,d\mu.$$
On peut remplacer l'ensemble d'int\'egration en $\lambda$ par un domaine fondamental ${\cal X}$ dans $i{\cal A}_{M_{disc}}^*$. Alors $B'(\lambda,a)$ est de Schwartz en les deux variables. On remarque que les fonctions $\vert H_{G}(a)-\theta H_{G}(a)\vert $ et $\vert H_{G}^{\tilde{G}}(a)\vert $ sont \'equivalentes. On peut alors  trouver une fonction $B''$ de Schwartz sur ${\cal X}$, \`a valeurs positives, de sorte que l'on ait l'in\'egalit\'e $\vert B'(\lambda,a)\vert \leq B''(\lambda)(1+\vert H_{G}^{\tilde{G}} (a)\vert )^{-r}$. Alors $\omega_{Q,w,w'}(am)$ est major\'ee par le produit de $(1+\vert H_{G}^{\tilde{G}} (a)\vert )^{-r}$ et d'une fonction qui a essentiellement la m\^eme forme que $\omega_{Q,w,w'}(m)$, sauf que l'on a supprim\'e l'int\'egration centrale. Cette derni\`ere fonction est essentiellement major\'ee par $\delta_{Q}(m)^{-1}\Xi^L(m)^2$ par le m\^eme argument qu'au lemme 3.6. On en d\'eduit l'\'enonc\'e. $\square$

\bigskip

\subsection{Un lemme de convergence}
Soient $Q=LU_{Q}\subset R$ deux sous-groupes paraboliques standard et soient $w\in W^G(L\vert  S)$ et $w'\in W^G(L\vert  S')$. On suppose $s_{Q}^R(w,w')\not=0$. Il existe donc un unique espace parabolique $\tilde{P}=\tilde{M}U_{P}$ tel que $Q\subset P\subset R$ et $ (w,w')\in {\cal W}_{Q}^P$. Pour $m\in M_{0}(F)$, on pose
$$E^T_{Q,R,w,w'}(m)=\phi^Q(H_{0}(m)-T)\tilde{\sigma}_{Q}^R(H_{0}(m)-T)\delta_{Q}(m)D_{0}^L(m)\omega_{Q,w,w'}(m).$$

\ass{Lemme}{Sous ces hypoth\`eses, l'int\'egrale
$$\int_{A_{\tilde{G}}(F)\backslash M_{0}(F)^{\geq,Q}}\vert E^T_{Q,R,w,w'}(m)\vert \,dm$$
est convergente.}

Preuve. Le lemme 3.5(iv) nous fournit  une majoration
$$\vert (H_{0}(m)-T)^{\tilde{G}}_{L}\vert << N^L_{w,w'}(H_{0}(m))$$
pour tout  $T$ et tout $m\in M_{0}(F)^{\geq,Q}$ tel que $\phi^Q(H_{0}(m)-T)\tilde{\sigma}_{Q}^R(H_{0}(m)-T)=1$. On a $(H_{0}(m)-T)^{\tilde{G}}_{L}=H_{0}^{\tilde{G}}(m)-H_{0}^L(m)-T_{L}^{\tilde{G}}$. Les hypoth\`eses $\phi^Q(H_{0}(m)-T)=1$ et $m\in M_{0}(F)^{\geq,Q}$ entra\^{\i}nent une majoration $\vert H_{0}^L(m)\vert <<\vert T\vert $. D'o\`u 
$$ \vert H_{0}^{\tilde{G}}(m)\vert <<\vert T\vert+\vert (H_{0}(m)-T)^{\tilde{G}}_{L}\vert,$$
puis
$$(1) \qquad \vert H_{0}^{\tilde{G}}(m)\vert <<\vert T\vert +N^L_{w,w'}(H_{0}(m)).$$
  Remarquons que, si l'on consid\`ere $T$ comme fix\'e (comme on le peut ici), la majoration (1) entra\^{\i}ne plus simplement
$$\vert H_{0}^{\tilde{G}}(m)\vert <<N^L_{w,w'}(H_{0}(m))$$
dans le domaine consid\'er\'e. En utilisant cette majoration, le lemme 3.6 et les majorations famili\`eres concernant les fonctions $D_{0}^L$ et $\Xi^L$, on voit que l'int\'egrale de l'\'enonc\'e est essentiellement major\'ee par
$$\int_{A_{\tilde{G}}F)\backslash M_{0}(F)^{\geq,Q}}(1+\vert H_{0}^{\tilde{G}}(m)\vert )^{-r}\,dm$$
pour tout r\'eel $r$. Cette int\'egrale est convergente, d'o\`u le lemme. $\square$

\bigskip

\subsection{Comparaison de deux int\'egrales}
On conserve les hypoth\`eses du paragraphe pr\'ec\'edent. Soit $\nu$ un r\'eel strictement positif. Notons ${\bf 1}_{\nu T}$ la fonction caract\'eristique de l'ensemble des $m\in M_{0}(F)$ tels que $\vert H_{0}^{\tilde{G}}(m)\vert\leq \nu\vert T\vert $.

\ass{Lemme}{Sous ces hypoth\`eses, si $\nu$ est assez grand, la diff\'erence entre les deux int\'egrales
$$\int_{A_{\tilde{G}}(F)\backslash M_{0}(F)^{\geq,Q}} E^T_{Q,R,w,w'}(m)\,dm$$
et
$$\int_{A_{\tilde{G}}(F)\backslash M_{0}(F)^{\geq}}E^T_{Q,R,w,w'}(m){\bf 1}_{\nu T}(m) \,dm$$
est essentiellement major\'ee par $\vert T\vert ^{-r}$ pour tout r\'eel $r$.}

Preuve. On peut consid\'erer que le domaine d'int\'egration de la seconde int\'egrale  est l'ensemble des $m\in A_{\tilde{G}}(F)\backslash M_{0}(F)^{\geq}$ tels que ${\bf 1}_{\nu T}(m)=1$. Ce domaine est contenu dans le domaine d'int\'egration  $A_{\tilde{G}}(F)\backslash M_{0}(F)^{\geq,Q}$ de la premi\`ere int\'egrale. Notons ${\cal D}^T$ le compl\'ementaire du premier domaine dans le second. On va prouver plus pr\'ecis\'ement:

(1) pour $\nu$ assez grand, on a pour tout r\'eel $r$ une majoration
$$\int_{{\cal D}^T}\vert E^T_{Q,R,w,w'}(m)\vert \,dm\,<<\vert T\vert ^{-r}.$$

L'ensemble ${\cal D}^T$ est r\'eunion disjointe des deux sous-ensembles

${\cal D}^T_{1}=A_{\tilde{G}}(F)\backslash \{m\in M_{0}(F)^{\geq,Q}; {\bf 1}_{\nu T}(m)=0\}$;

${\cal D}^T_{2}=A_{\tilde{G}}(F)\backslash \{m\in M_{0}(F)^{\geq,Q}-M_{0}(F)^{\geq}; {\bf 1}_{\nu T}(m)=1\}$.

Compte tenu de la d\'efinition de $E^T_{Q,R,w,w'}(m)$, on peut ajouter la condition $\phi^Q(H_{0}(m)-T)\tilde{\sigma}_{Q}^R(H_{0}(m)-T)=1$. Sur ${\cal D}^T_{1}$, on utilise la relation (1) du paragraphe pr\'ec\'edent, que l'on \'ecrit plus pr\'ecis\'ement
$$\vert H_{0}^{\tilde{G}}(m)\vert \leq c_{1}(\vert T\vert +N^L_{w,w'}(H_{0}(m))) .$$
Jointe \`a la condition ${\bf 1}_{\nu T}(m)=0$, elle nous dit que
$$(\nu-c_{1})\vert T\vert \leq c_{1}N^L_{w,w'}(H_{0}(m))$$
et
$$(1-\frac{c_{1}}{\nu})\vert H_{0}^{\tilde{G}}(m)\vert \leq c_{1}N^L_{w,w'}(H_{0}(m)).$$
On suppose $\nu>c_{1}$. On en d\'eduit une majoration
$$\vert H_{0}^{\tilde{G}}(m)\vert +\vert T\vert <<N^L_{w,w'}(H_{0}(m)).$$
En utilisant le lemme 3.6 et  les majorations habituelles des fonctions $D_{0}^L$ et $\Xi^L$, on voit que  l'int\'egrale
$$ \int_{{\cal D}_{1}^T}\vert E^T_{Q,R,w,w'}(m)\vert \,dm$$
est essentiellement major\'ee par
$$\vert T\vert ^{-r}\int_{A_{\tilde{G}}(F)\backslash M_{0}(F)}(1+\vert H_{0}^{\tilde{G}}(m)\vert )^{-r}\,dm$$
pour tout r\'eel $r$. Cette derni\`ere  int\'egrale est convergente, d'o\`u la majoration 
$$ \int_{{\cal D}_{1}^T}\vert E^T_{Q,R,w,w'}(m)\vert \,dm<<\vert T\vert ^{-r}.$$

Traitons maintenant l'int\'egrale sur ${\cal D}^T_{2}$. Le lemme 3.5(iv) nous fournit la majoration
$$\vert (H_{0}(m)-T)_{L}^{\tilde{G}}\vert <<N^L_{w,w'}(H_{0}(m))$$
pour tout $T$ et tout $m\in M_{0}(F)^{\geq,Q}$ tel que $\phi^Q(H_{0}(m)-T)\tilde{\sigma}_{Q}^R(H_{0}(m)-T)=1$. Supposons de plus $m\not\in M_{0}(F)^{\geq}$. Il existe donc $\alpha\in \Delta_{0}$ tel que $<\alpha,H_{0}(m)>\leq0$. N\'ecessairement, $\alpha\not\in \Delta_{0}^Q$ (puisque $m\in M_{0}(F)^{\geq,Q}$). On a
$$<\alpha,(H_{0}(m)-T)_{L}^{\tilde{G}}>=<\alpha,H_{0}(m)-T>+<\alpha,(T-H_{0}(m))^L> .$$
Le premier terme est inf\'erieur ou \'egal \`a $-<\alpha,T>$. Le second est n\'egatif ou nul: l'hypoth\`ese $\phi^Q(H_{0}(m)-T)=1$ entra\^{\i}ne que $(T-H_{0}(m))^L$ est combinaison lin\'eaire \`a coefficients positifs ou nuls de $\check{\beta}$ pour $\beta\in \Delta_{0}^Q$ et on a $<\alpha,\check{\beta}>\leq0$ pour tout tel $\beta$. Donc
$$<\alpha,(H_{0}(m)-T)_{L}^{\tilde{G}}>\leq -<\alpha,T>.$$
A fortiori $\vert (H_{0}(m)-T)_{L}^{\tilde{G}}\vert >>\vert T\vert$, ce qui prouve la majoration
$$ \vert T\vert <<N^L_{w,w'}(H_{0}(m))$$
pour tout $m\in M_{0}(F)^{\geq,Q}-M_{0}(F)^{\geq}$ tel que $\phi^Q(H_{0}(m)-T)\tilde{\sigma}_{Q}^R(H_{0}(m)-T)=1$.   En utilisant  encore le lemme 3.6 et  les majorations habituelles des fonctions $D_{0}^L$ et $\Xi^L$, on voit que  l'int\'egrale
$$ \int_{{\cal D}_{2}^T}\vert E^T_{Q,R,w,w'}(m)\vert \,dm$$
est essentiellement major\'ee par
$$\vert T\vert ^{-r}\int_{A_{\tilde{G}}(F)\backslash M_{0}(F)}(1+\vert H_{0}^{\tilde{G}}(m)\vert )^D{\bf 1}_{\nu T}(m)\,dm$$
pour un certain entier $D$ et pour tout r\'eel $r$. La derni\`ere int\'egrale est convergente et essentiellement born\'ee par $\vert T\vert ^{D'}$ pour un certain entier $D'$.  On en d\'eduit
$$ \int_{{\cal D}_{2}^T}\vert E^T_{Q,R,w,w'}(m)\vert \,dm<<\vert T\vert ^{-r}$$
pour tout r\'eel $r$. Cela prouve (1) et le lemme. $\square$

\bigskip

\subsection{Un lemme d'\'equivalence}
 Posons
 $$E^T(m)=\sum_{Q,R; P_{0}\subset Q\subset R}\sum_{w\in W^G(L\vert  S),w'\in W^G(L\vert  S')}s_{Q}^R(w,w')E^T_{Q,R,w,w'}(m)$$
$$=\sum_{Q,R; P_{0}\subset Q\subset R}\sum_{w\in W^G(L\vert  S),w'\in W^G(L\vert  S')}s_{Q}^R(w,w')\phi^Q(H_{0}(m)-T)\tilde{\sigma}_{Q}^R(H_{0}(m)-T)$$
$$\delta_{Q}(m)D_{0}^L(m)\omega_{Q,w,w'}(m).$$

\ass{Lemme}{Les fonctions $E^T(m)$ et $\Omega_{G}(m)\phi^{\tilde{G}}(H_{0}(m)-T)D_{0}^G(m)$ sont \'equivalentes.}

Preuve. Dans la somme d\'efinissant ${\cal E}^T(m)$, on peut glisser une sous-somme
$$\sum_{Q; P_{0}\subset Q\subset P}\phi^Q(H_{0}(m)-T)\tau_{Q}^P(H_{0}(m)-T)$$
puisque celle-ci vaut $1$.  En utilisant le lemme 3.7 (i) et (ii), on peut alors, \`a \'equivalence pr\`es, remplacer la fonction $\Omega_{P}(m)\delta_{P}(m)D_{0}^M(m)$ par 
$$\sum_{(w,w')\in {\cal W}_{Q}^P}\omega_{Q,w,w'}(m)\delta_{Q}(m)D_{0}^L(m).$$
  On a aussi l'\'egalit\'e
$$\tau_{Q}^P(H_{0}(m)-T)\hat{\tau}_{\tilde{P}}(H_{0}(m)-T)=\sum_{R; P\subset R}\tilde{\sigma}_{Q}^R(H_{0}(m)-T),$$
cf. 2.2(1). A ce point, on a montr\'e que ${\cal E}^T(m)$ est \'equivalente \`a
$$\sum_{Q,R; P_{0}\subset Q\subset R}\sum_{w\in W^G(L\vert  S),w'\in W^G(L\vert  S')}\phi^Q(H_{0}(m)-T)\tilde{\sigma}_{Q}^R(H_{0}(m)-T)$$
$$\delta_{Q}(m)D_{0}^L(m)\omega_{Q,w,w'}(m)\sum_{\tilde{P}; Q\subset P\subset R, (w,w')\in {\cal W}_{Q}^P}(-1)^{dim(\mathfrak{a}_{\tilde{P}})-dim(\mathfrak{a}_{\tilde{G}})}.$$
La derni\`ere somme est \'egale par d\'efinition \`a $s_{Q}^R(w,w')$, donc l'expression ci-dessus est \'egale \`a $E^T(m)$. Donc $E^T(m)$ est \'equivalente \`a $ {\cal E}^T(m)$. La proposition 3.10 nous dit que ${\cal E}^T(m)$ est \'equivalente \`a la deuxi\`eme fonction de l'\'enonc\'e. $\square$

\bigskip

\subsection{Preuve de la proposition 3.4}
 Que l'expression d\'efinissant $j_{\star}^T$ soit absolument convergente r\'esulte du lemme 3.12. Fixons un r\'eel $\nu>0$ que l'on pr\'ecisera plus tard. Introduisons l'expression
 $$j_{1}^T= \int_{A_{\tilde{G}}(F)\backslash M_{0}(F)^{\geq}}E^T(m){\bf 1}_{\nu T}(m)\,dm.$$
  Elle aussi est absolument convergente. Le lemme 3.13 nous dit que
$$lim_{T\to \infty}j_{\star}^T-j_{1}^T=0.$$
Pour $m\in M_{0}(F)^{\geq}$, posons
$$\psi(m)=\Omega_{G}(m)D_{0}^G(m)\tilde{\kappa}^T(m)-E^T(m){\bf 1}_{\nu T}(m).$$
Remarquons que toute int\'egrale sur $A_{\tilde{G}}(F)\backslash M_{0}(F)^{\geq}$ se d\'ecompose en produit d'une int\'egrale sur $M_{0}(F)^{\geq,{\bf b}}$ et d'une int\'egrale sur $A_{\tilde{G}}(F)\backslash A_{G}(F)$, cf. 3.11 pour la d\'efinition de $M_{0}(F)^{\geq,{\bf b}}$. On d\'ecompose $j^T-j_{1}^T$ en $j_{2}^T+j_{3}^T$, o\`u
$$j_{2}^T= \int_{M_{0}(F)^{\geq,{\bf b}}}\int_{a\in A_{\tilde{G}}(F)\backslash A_{G}(F); \vert H_{G}^{\tilde{G}}(a)\vert >\vert T\vert }  \psi(am)\,da\,dm,$$
$$j_{3}^T= \int_{M_{0}(F)^{\geq,{\bf b}}}\int_{a\in A_{\tilde{G}}(F)\backslash A_{G}(F); \vert H_{G}^{\tilde{G}}(a)\vert \leq\vert T\vert }  \psi(am)\,da\,dm.$$
Dans la premi\`ere, on utilise le lemme 3.11. On peut l'appliquer \`a chaque fonction intervenant dans la d\'efinition de $E^T(m)$, ainsi qu'\`a la fonction $\Omega_{G}(m)D_{0}^G(m)$ qui en est un cas particulier. Ce lemme nous fournit une majoration
$$\vert \psi(am)\vert <<(1+\vert H_{0}(m)\vert)^D(1+\vert H_{G}^{\tilde{G}}(a)\vert )^{-r}$$
pour un certain entier $D$ et pour tout r\'eel $r$. L'int\'egrale en $a$ est essentiellement major\'ee par $\vert T\vert ^{-r}$ pour tout r\'eel $r$. L'int\'egrale en $m$ porte sur un domaine o\`u $\vert H_{0}(m)\vert <<\vert T\vert $ car une telle in\'egalit\'e est v\'erifi\'ee sur l'intersection de $M_{0}(F)^{\geq,{\bf b}}$ et des supports des fonctions $\tilde{\kappa}^T(am)$ ou ${\bf 1}_{\nu T}(am)$.  L'int\'egrale en $m$ est donc essentiellement born\'ee par $\vert T\vert ^{D'}$ pour un certain entier $D'$. Il en r\'esulte que $lim_{T\to \infty}j_{2}^T=0$. Consid\'erons $j_{3}^T$. Si $\nu$ est assez grand, les conditions $m\in M_{0}(F)^{\geq,{\bf b}}$, $\vert H_{G}^{\tilde{G}}(a)\vert \leq \vert T\vert $ et $\tilde{\kappa}^T(am)=1$ impliquent $\vert H_{0}(am)^{\tilde{G}}\vert \leq \nu\vert T\vert $. On ne change alors rien en multipliant la fonction $\Omega_{G}(am)D_{0}^G(am)\tilde{\kappa}^T(am)$ par ${\bf 1}_{\nu T}(am)$. Pour $m\in M_{0}(F)^{\geq}$, on a $\tilde{\kappa}^T(m)=\phi^{\tilde{ G}}(H_{0}(m)-T)$. Le lemme 3.14 nous fournit pour tout $r$ une majoration
$$\vert \psi(am)\vert<<  \vert T\vert ^{-r}$$
sur le domaine d'int\'egration de $j_{3}^T$, pour un certain entier $D$ et pour tout $r$. De nouveau, l'int\'egrale de cette fonction sur le domaine ${\bf 1}_{\nu T}(am)=1$ est major\'ee par $\vert T\vert ^{-r}$ pour tout $r$. Donc $lim_{T\to \infty}j_{3}^T=0$. Cela ach\`eve la preuve. $\square$
 
\bigskip

\subsection{D\'efinition de $(G,M)$-familles}
Posons $M'_{disc}=\theta^{-1}(M_{disc})$. Soit $t\in W^G(M_{disc}\vert M'_{disc})$, c'est-\`a-dire que $t\in W^G/W^{M'_{disc}}$ et $t(M'_{disc})=M_{disc}$. Soit $\nu\in [t\sigma',\omega \sigma]$. On fixe un automorphisme unitaire $A_{\nu}$ de $V_{\sigma}$ tel que $(t\sigma')(x)\circ A_{\nu}=\omega(x)A_{\nu}\circ \sigma_{-\nu}(x)$ pour tout $x\in M_{disc}(F)$. Par fonctorialit\'e, il d\'efinit des homomorphismes entre diff\'erentes repr\'esentations induites. Notons $\underline{\omega}$ l'op\'erateur qui, \`a une fonction $\varphi$ sur $G(F)$, associe la fonction $g\mapsto \omega(g)\varphi(g)$. Lui-aussi d\'efinit des homomorphismes entre certaines repr\'esentations induites. Pour $\Lambda\in i{\cal A}_{M_{disc},F}^*$, introduisons l'op\'erateur
$$A(t,\nu;\Lambda):\pi'_{t^{-1}(\Lambda+\nu)}=Ind_{S'}^G(\sigma'_{t^{-1}(\lambda+\nu)})\to \pi_{\Lambda}=Ind_{S}^G(\sigma_{\Lambda})$$
d\'efini par
$$A(t,\nu;\Lambda)=R_{S\vert t(S')}(\sigma_{\Lambda})\circ\gamma(t)\circ A_{\nu}^{-1}\underline{\omega}^{-1}.$$
Il  v\'erifie la relation d'entrelacement  $\omega(g)\pi_{\Lambda}(g)\circ A(t,\nu;\Lambda)=A(t,\nu;\Lambda)\circ \pi'_{t^{-1}(\Lambda+\nu)}(g)$.

Pour $S''\in {\cal P}(M_{disc})$, on d\'efinit la fonction $(\lambda,\Lambda)\mapsto \phi(t,\nu;\lambda,\Lambda,S'')$ des deux variables $\lambda,\Lambda\in i{\cal A}_{M_{disc},F}^*$ par
$$\phi(t,\nu;\lambda,\Lambda,S'')=(A(t,\nu;t\lambda'-\nu)v',J_{\bar{S}''\vert S}(\sigma_{t\lambda'-\nu})^{-1}\circ J_{\bar{S}''\vert S}(\sigma_{t\lambda'-\nu+\Lambda})u )$$
$$( J_{S''\vert S}(\sigma_{t\lambda'-\nu})^{-1}\circ J_{S''\vert S}(\sigma_{t\lambda'-\nu+\Lambda})v,A(t,\nu;t\lambda'-\nu)u').$$
 Pour $\lambda$ fix\'e, la famille $(\phi(t,\nu;\lambda,\Lambda,S''))_{S''\in {\cal P}(M_{disc})}$ de fonctions en $\Lambda$ est presque une $(G,M_{disc})$-famille et m\^eme presque une $(G,M_{disc})$-famille $p$-adique dans le cas o\`u $F$ est non-archim\'edien. "Presque" parce que les fonctions ne sont pas forc\'ement $C^{\infty}$: il peut y avoir des singularit\'es. Etudions cette question de r\'egularit\'e. En convertissant les op\'erateurs d'entrelacement en op\'erateurs normalis\'es, on peut r\'ecrire la d\'efinition de $\phi(t,\nu;\lambda,\Lambda,S'')$ sous la forme
 $$\phi(t,\nu;\lambda,\Lambda,S'')=r_{\bar{S}''\vert S''}(\sigma_{t\lambda'-\nu})^{-1}r_{ \bar{S}''\vert S''}(\sigma_{t\lambda'-\nu+\Lambda})$$
 $$(A(t,\nu;t\lambda'-\nu)v',R_{ \bar{S}''\vert S}(\sigma_{t\lambda'-\nu})^{-1}\circ R_{\bar{S}''\vert S}(\sigma_{t\lambda'-\nu+\Lambda})u)$$
 $$( R_{ S''\vert S}(\sigma_{t\lambda'-\nu})^{-1}\circ R_{S''\vert S}(\sigma_{t\lambda'-\nu+\Lambda})v,A(t,\nu;t\lambda'-\nu)u').$$
   Posons
 $$r_{S'',reg}(\sigma_{\lambda})=r_{ \bar{S}''\vert S''}(\sigma_{\lambda})r_{ \bar{S}\vert S}(\sigma_{\lambda})^{-1}.$$
 On peut \'ecrire
 $$\phi(t,\nu;\lambda,\Lambda,S'')=r_{ \bar{S}\vert S}(\sigma_{t\lambda'-\nu})^{-1}r_{ \bar{S}\vert S}(\sigma_{t\lambda'-\nu+\Lambda})\phi_{reg}(t,\nu;\lambda,\Lambda,S''),$$
 o\`u
 $$\phi_{reg}(t,\nu;\lambda,\Lambda,S'')=r_{S'',reg}(\sigma_{t\lambda'-\nu})^{-1}r_{S'',reg}(\sigma_{t\lambda'-\nu+\Lambda})$$
 $$(A(t,\nu;t\lambda'-\nu)v',R_{ \bar{S}''\vert S}(\sigma_{t\lambda'-\nu})^{-1}\circ R_{\bar{S}''\vert S}(\sigma_{t\lambda'-\nu+\Lambda})u)$$
 $$( R_{ S''\vert S}(\sigma_{t\lambda'-\nu})^{-1}\circ R_{S''\vert S}(\sigma_{t\lambda'-\nu+\Lambda})v,A(t,\nu;t\lambda'-\nu)u').$$

 Gr\^ace \`a 1.10(5), la fonction $\phi_{reg}(t,\nu;\lambda,\Lambda,S'')$ est r\'eguli\`ere en $\lambda$ et $\Lambda$. La famille $(\phi_{reg}(t,\nu;\lambda,\Lambda,S''))_{S''\in {\cal P}(M_{disc})}$ est vraiment une $(G,M)$-famille, $p$-adique dans le cas o\`u $F$ est non-archim\'edien. Toutes les singularit\'es de la famille de d\'epart se concentrent dans la fonction $r_{ \bar{S}\vert S}(\sigma_{t\lambda'-\nu})^{-1}r_{ \bar{S}\vert S}(\sigma_{t\lambda'-\nu+\Lambda})$  en facteur.
 
 Soit $Q=LU_{Q}\in {\cal F}(M_{disc})$. On d\'eduit de la famille $(\phi(t,\nu;\lambda,\Lambda,S''))_{S''\in {\cal P}(M_{disc})}$ une famille $(\phi(t,\nu;\lambda,\Lambda,S^{_{''}L}))_{S^{_{''}L}\in {\cal P}^L(M_{disc})}$, qui n'est autre que $(\phi(t,\nu;\lambda,\Lambda,S''))_{S''\in {\cal P}(M_{disc}),S''\subset Q}$.
   
 Pour $X\in {\cal A}_{L,F}$ et $S''\in {\cal P}(M_{disc})$ tel que $S''\subset Q$, on a d\'efini en 1.6 la fonction
 $$\epsilon_{S''}^{Q,T[S'']}(X;\Lambda)=\int_{{\cal A}_{M_{disc},F}^{L}(X)}\phi_{S''}^L(Y-T[S''])e^{<\Lambda,Y>}\,dY$$
 ($\phi_{S''}^L$ est la fonction combinatoire de 1.3, qui n'a rien \`a voir avec les fonctions de la $(G,M_{disc})$-famille).
 On pose
 $$\phi_{M_{disc}}^{Q,T}(t,\nu,X;\lambda,\Lambda)=\sum_{S''\in {\cal P}(M_{disc}),S''\subset Q}\phi(t,\nu;\lambda,\Lambda,S'')\epsilon_{S''}^{Q,T[S'']}(X;\Lambda).$$
 On d\'efinit de m\^eme $\phi_{reg,M_{disc}}^{Q,T}(t,\nu,X;\lambda,\Lambda)$.
 Rappelons que l'on a pos\'e $\lambda'=\theta^{-1}\lambda$. On pose
 $$\epsilon(t,\nu;\lambda)=r_{\bar{S}\vert S}(\sigma_{\lambda})r_{\bar{S}\vert S}(\sigma_{t\lambda'-\nu})^{-1} .$$
 
 \ass{Lemme}{(i) Pour $\lambda$ en position g\'en\'erale, la fonction $\Lambda\mapsto \phi_{M_{disc}}^{Q,T}(t,\nu,X;\lambda,\Lambda)$ est r\'eguli\`ere en $\Lambda=\lambda-t\lambda'+\nu$.
 
 (ii) Les fonctions
  $$\lambda\mapsto \phi_{M_{disc}}^{Q,T}(t,\nu,X;\lambda, \lambda-t\lambda'+\nu),$$
  $$\lambda\mapsto \phi_{reg,M_{disc}}^{Q,T}(t,\nu,X;\lambda, \lambda-t\lambda'+\nu),$$
  $$\lambda\mapsto  \epsilon(t,\nu;\lambda)$$
   sont $C^{\infty}$ sur $i{\cal A}_{M,F}^*$.   Si $F$ est archim\'edien, toutes leurs d\'eriv\'ees sont \`a croissance lente.
   
   (iii) On a l'\'egalit\'e
 $$\phi_{M_{disc}}^{Q,T}(t,\nu,X;\lambda, \lambda-t\lambda'+\nu)= \epsilon(t,\nu;\lambda) \phi_{reg,M_{disc}}^{Q,T}(t,\nu,X;\lambda, \lambda-t\lambda'+\nu) .$$ 
 
 (iv) Comme fonction de $X$, $ \phi_{M_{disc}}^{Q,T}(t,\nu,X;\lambda,\lambda-t\lambda'+\nu)$ ne d\'epend que de la classe $X+{\cal A}_{A_{\tilde{G}},F}$.}
 
 Preuve. Supposons d'abord $F$ non-archim\'edien. La discussion ci-dessus montre que l'on a l'\'egalit\'e 
 $$\phi_{M_{disc}}^{Q,T}(t,\nu,X;\lambda,\Lambda)=r_{ \bar{S}\vert S}(\sigma_{t\lambda'-\nu+\Lambda})r_{ \bar{S}\vert S}(\sigma_{t\lambda'-\nu})^{-1}\phi_{reg,M_{disc}}^{Q,T}(t,\nu,X;\lambda,\Lambda).$$
  La  derni\`ere fonction  est $C^{\infty}$ en $\lambda$ et $\Lambda$ d'apr\`es le lemme 1.6. Le facteur $r_{ \bar{S}\vert S}(\sigma_{t\lambda'-\nu+\Lambda})r_{ \bar{S}\vert S}(\sigma_{t\lambda'-\nu})^{-1}$ est \'evidemment r\'egulier en $\Lambda=\lambda-t\lambda'+\nu$ pour $\lambda$ en position g\'en\'erale. Il vaut alors $r_{ \bar{S}\vert S}(\sigma_{ \lambda})r_{ \bar{S}\vert S}(\sigma_{t\lambda'-\nu})^{-1}$. Pour prouver (ii) et (iii), il reste \`a montrer que  ce dernier terme est r\'egulier en $\lambda$. Par l'isomorphisme $\omega\sigma_{-\nu}\simeq  t\sigma'$, on a $\sigma_{t\lambda'-\nu}\simeq \omega^{-1}t(\sigma_{\lambda}\circ\theta) $. Par transport de structure et parce que les facteurs de normalisation sont insensibles \`a la torsion par le caract\`ere $\omega$, on a $r_{ \bar{S}\vert S}(\sigma_{t\lambda'-\nu})=r_{\theta t^{-1}(\bar{S})\vert \theta t^{-1}(S)}(\sigma_{\lambda})$. Mais alors, le quotient
 $$r_{\bar{S}\vert S}(\sigma_{\lambda})r_{ \bar{S}\vert S}(\sigma_{t\lambda'-\nu})^{-1}=r_{\bar{S}\vert S}(\sigma_{\lambda})r_{\theta t^{-1}(\bar{S})\vert \theta t^{-1}(S)}(\sigma_{\lambda})^{-1}$$
 est r\'egulier d'apr\`es 1.10(5).
 
 Supposons maintenant $F$ archim\'edien. Le m\^eme raisonnement prouve les assertions de r\'egularit\'e. Pour montrer l'assertion de croissance lente, le lemme 1.4 nous ram\`ene \`a prouver que  que
 
 - les d\'eriv\'ees des fonctions intervenant dans la d\'efinitions des termes $\phi_{reg}(t,\nu;\lambda,\Lambda,S'')$  sont \`a croissance lente;
 
 - les d\'eriv\'ees de la fonction $r_{ \bar{S}\vert S}(\sigma_{\lambda})r_{ \bar{S}\vert S}(\sigma_{t\lambda'-\nu})^{-1}$ sont \`a croissance lente.
 
 Les op\'erateurs d'entrelacement normalis\'es ont des d\'eriv\'ees \`a croissance lente. Les autres fonctions intervenant sont toutes de la forme $r_{ \bar{S}_{1}\vert S_{1}}(\sigma_{\mu})r_{\bar{S}_{2}\vert S_{2}}(\sigma_{\mu})^{-1}$ o\`u $S_{1},S_{2}\in {\cal P}(M_{disc})$  et $\mu$ d\'epend lin\'eairement de $\lambda$. L'assertion r\'esulte de 1.10(5). 
 
 Preuve de (iv). Pour $Y\in {\cal A}_{A_{\tilde{G}},F}$, l'\'egalit\'e suivante r\'esulte des d\'efinitions
 $$ \phi_{M_{disc}}^{Q,T}(t,\nu,X+Y;\lambda,\lambda-t\lambda'+\nu)=e^{<\lambda-t\lambda'+\nu,Y>} \phi_{M_{disc}}^{Q,T}(t,\nu,X;\lambda,\lambda-t\lambda'+\nu).$$
 Evidemment $e^{<\lambda-t\lambda',Y>}=1$. En comparant les restrictions \`a $A_{\tilde{G}}(F)$ des caract\`eres centraux des repr\'esentations $\omega\sigma$ et $ t\sigma'$, l'isomorphisme $\omega\sigma_{-\nu}\simeq  t\sigma'$ et le fait que $\omega$ soit trivial sur $A_{\tilde{G}}(F)$ entra\^{\i}ne que $\nu_{\tilde{G}}\in i{\cal A}_{A_{\tilde{G}},F}^{\vee}$. Donc $e^{<\nu,Y>}=1$. Cela prouve (iii). $\square$

\bigskip 
\subsection{D\'efinition d'une nouvelle int\'egrale} 
Pour deux sous-groupes paraboliques $Q=LU_{Q},R\in {\cal F}(M_{disc})$ tels que $Q\subset R$ et pour $t\in W^G(M_{disc}\vert M'_{disc})$, notons $s_{Q}^{R}(t)$ la somme des
$$(-1)^{dim(\mathfrak{a}_{\tilde{M}})-dim(\mathfrak{a}_{\tilde{G}})}$$
sur les espaces paraboliques $\tilde{P}=\tilde{M}U_{P}$ tels que $Q\subset P\subset R$ et $t\theta^{-1}(P)=P$.   Cette derni\`ere condition \'equivaut \`a $\gamma_{0}t^{-1}\in \tilde{P}(F)$ (o\`u on identifie $t$ \`a un rel\`evement dans $K$). La condition $s_{Q}^{R}(t)\not=0$ \'equivaut \`a ce qu'il existe un et un seul $\tilde{P}$ v\'erifiant ces conditions. Soient $t\in W^G(M_{disc}\vert M'_{disc})$ et $\nu\in [t\sigma',\omega\sigma]$. Posons
$${\bf E}^T_{Q,R,t,\nu}=mes(i{\cal A}_{M_{disc},F}^*)^{-1}\int_{{\cal A}_{L,F}/{\cal A}_{A_{\tilde{G}},F}}\tilde{\sigma}_{Q}^{R}(X-T[Q])$$
$$\int_{i{\cal A}_{M_{disc},F}^*}B(\lambda) \phi_{M_{disc}}^{Q,T}(t,\nu,X;\lambda,\lambda-t\lambda'+\nu)\,d\lambda\,dX.$$

\ass{Lemme}{Si $s_{Q}^R(t)\not=0$, l'expression ${\bf E}^T_{Q,R,t,\nu}$ est convergente dans l'ordre indiqu\'e.}

Preuve. La convergence absolue de l'int\'egrale int\'erieure r\'esulte du (ii) du lemme pr\'ec\'edent et du fait que $B$ est de Schwartz. Puisque ${\cal A}_{A_{L},F}$ est d'indice fini dans ${\cal A}_{L,F}$, on peut fixer $X_{0}\in {\cal A}_{L,F}$ et prouver que l'expression
$$\int_{{\cal A}_{A_{L},F}/{\cal A}_{A_{\tilde{G}},F}}\tilde{\sigma}_{Q}^{R}(X+X_{0}-T[Q])\vert \int_{i{\cal A}_{M_{disc},F}^*}B(\lambda) \phi_{M_{disc}}^{Q,T}(t,\nu,X+X_{0};\lambda,\lambda-t\lambda'+\nu)\,d\lambda\vert \,dX$$
est convergente. Il r\'esulte des d\'efinitions que
$$\phi_{M_{disc}}^{Q,T}(t,\nu,X+X_{0};\lambda,\lambda-t\lambda'+\nu)=e^{<\lambda-t\lambda'+\nu,X>} \phi_{M_{disc}}^{Q,T}(t,\nu,X_{0};\lambda,\lambda-t\lambda'+\nu).$$
Ainsi l'expression ci-dessus est de la forme
$$\int_{{\cal A}_{A_{L},F}/{\cal A}_{A_{\tilde{G}},F}}\tilde{\sigma}_{Q}^{R}(X+X_{0}-T[Q])\vert \int_{i{\cal A}_{M_{disc},F}^*}B'(\lambda)e^{<\lambda-t\lambda',X>}\,d\lambda\vert \,dX,$$
o\`u $B'$ est une fonction de Schwartz. En utilisant l'\'egalit\'e
$$<\lambda-t\lambda',X>=<\lambda,(1-\theta t^{-1})X>,$$
on voit que l'int\'egrale int\'erieure est essentiellement major\'ee par $(1+\vert ( 1-\theta t^{-1})X\vert )^{-r}$ pour tout r\'eel $r$. Pour d\'emontrer la convergence cherch\'ee, et puisque $T$ peut \^etre ici consid\'er\'e comme une constante, il suffit de prouver que l'on a une majoration
$$(1)\qquad \vert X^{\tilde{G}}\vert <<1+\vert T\vert +\vert ( 1-\theta t^{-1})X\vert $$
pour tout $T$ et tout $X\in {\cal A}_{A_{L},F}$ tel que $\tilde{\sigma}_{Q}^R(X+X_{0}-T[Q])=1$.  Puisque $X_{0}$ est fix\'e, cela r\'esulte par translation par $-X_{0}$ d'une majoration
$$ \vert (X-T[Q])^{\tilde{G}}_{L}\vert <<1+\vert ( 1-\theta t^{-1})X\vert $$
pour tout $X\in {\cal A}_{L}$ tel que $\tilde{\sigma}_{Q}^R(X-T[Q])=1$. Pour mieux comprendre la situation, introduisons un \'el\'ement $s\in W^G$ tel que $s(Q)$ soit standard. Posons $Q'=s(Q)=L'U_{Q'}$, 
$R'=s(R)$. La majoration pr\'ec\'edente  r\'esulte de la majoration

(2) $ \vert (X-T)^{\tilde{G}}_{L'}\vert <<1+\vert (s\theta t^{-1}s^{-1}-1)X\vert $ pour tout $X\in {\cal A}_{L'}$ tel que $\tilde{\sigma}_{Q'}^{R'}(X)=1$.

La condition $s_{Q}^R(t)\not=0$ se traduit ainsi: il existe un unique espace parabolique $\tilde{P}=\tilde{M}U_{P}$ tel que $Q'\subset P\subset R'$ et $ \theta(st)s^{-1}\in W^M$, autrement dit $s_{Q'}^{R'}(s,st)\not=0$. On applique le lemme 3.5(iv) en y rempla\c{c}ant les termes $Q,R,S,w,w',H$ par $Q',R',P_{0},s,st,X$ (ces donn\'ees v\'erifient les conditions de ce lemme). Ce lemme implique
$$\vert (X-T)^{\tilde{G}}_{L'}\vert <<N_{s,st}^{L'}(X).$$
Mais, par d\'efinition de ce dernier terme, on a 
$$N_{s,st}^{L'}(X)<<1+\vert (s\theta t^{-1}s^{-1}-1)X\vert .$$
Cela prouve (2) et le lemme. $\square$

\bigskip

\subsection{Apparition des $(G,M)$-familles}

Soient $Q=LU_{Q}\subset R$ deux sous-groupes paraboliques standard et soient $w\in W^G(L\vert S)$, $w'\in W^G(L\vert S')$ deux \'el\'ements tels que $s_{Q}^R(w,w')\not=0$. On note $\tilde{P}=\tilde{M}U_{P}$ l'unique espace parabolique tel que $Q\subset P\subset R$ et $(w,w')\in {\cal W}_{Q}^P$. On pose
$$E^T_{Q,R,w,w'}=\int_{A_{\tilde{G}}(F)\backslash M_{0}(F)^{\geq,Q}}E^T_{Q,R,w,w'}(m)\,dm$$
$$=\int_{A_{\tilde{G}}(F)\backslash M_{0}(F)^{\geq,Q}}\tilde{\sigma}_{Q}^R(H_{0}(m)-T)\phi^Q(H_{0}(m)-T)\delta_{Q}(m)D_{0}^L(m)\omega_{Q,w,w'}(m)\,dm.$$
 Le but du paragraphe est de d\'efinir  une approximation plus explicite de $E^T_{Q,R,w,w'}$. On pose
 $$E^T_{\star;Q,R,w,w'}= \sum_{t\in W^G(M_{disc}\vert M'_{disc})\cap w^{-1}W^Lw'}\sum_{\nu\in[t\sigma',\omega \sigma]}{\bf E}^T_{w^{-1}(Q),w^{-1}(R),t,\nu}.$$
 On identifie comme toujours  un \'el\'ement $t\in W^G(M_{disc}\vert M'_{disc})$ \`a un rel\`evement dans $W^G$. La condition que $t$ appartienne \`a $w^{-1}W^Lw'$ ne d\'epend pas du rel\`evement choisi.  D'autre part, on voit que, pour $t$ dans l'ensemble de sommation, on a $s_{w^{-1}(Q)}^{w^{-1}(R)}(t)=s_{Q}^R(w,w')$. Notre hypoth\`ese est que ce nombre est non nul. Les termes ${\bf E}^T_{w^{-1}(Q),w^{-1}(R),t,\nu}$ sont donc bien d\'efinis d'apr\`es le lemme pr\'ec\'edent.
 
 \ass{Lemme}{On a une majoration
 $$\vert E^T_{Q,R,w,w'}-mes(A_{\tilde{G}}(F)_{c})^{-1}E^T_{\star;Q,R,w,w'}\vert <<\vert T\vert ^{-r}$$
 pour tout r\'eel $r$.}
 
 Preuve. Fixons un r\'eel $\zeta>0$. Notons ici ${\bf 1}_{\zeta T}$ la fonction caract\'eristique de l'ensemble des $H\in {\cal A}_{0}$ tels que $\vert H\vert \leq \zeta\vert T\vert $ (ce n'est pas la m\^eme fonction qu'en 3.13). Notons $E^{T,\zeta}_{Q,R,w,w'}$ la variante de l'expression $E^T_{Q,R,w,w'}$ o\`u on glisse la fonction ${\bf 1}_{\zeta T}(H_{L}^{\tilde{G}}(m))$ dans l'int\'egrale. Notons $E^{T,\zeta}_{\star;Q,R,w,w'}$ la variante de $E^T_{\star;Q,R,w,w'}$ o\`u on glisse la fonction ${\bf 1}_{\zeta T}(X^{\tilde{G}})$ dans les int\'egrales d\'efinissant chaque ${\bf E}^T_{w^{-1}(Q),w^{-1}(R),t,\nu}$. Il r\'esulte de 3.13(1) que, si $\zeta$ est assez grand, on a la majoration
 $$\vert E^T_{Q,R,w,w'}-E^{T,\zeta}_{Q,R,w,w'}\vert<<\vert T\vert ^{-r}$$
 pour tout r\'eel $r$. De m\^eme, en reprenant la preuve du lemme 3.17, il r\'esulte de 3.17(1) que, si $\zeta$ est assez grand, on a la majoration
 $$\vert E^T_{\star;Q,R,w,w'}-E^{T,\zeta}_{\star;Q,R,w,w'}\vert<<\vert T\vert ^{-r}$$
 pour tout r\'eel $r$. On fixe $\zeta$ tel qu'il en soit ainsi. Il nous suffit alors de majorer
 $$\vert E^{T,\zeta}_{Q,R,w,w'}-mes(A_{\tilde{G}}(F)_{c})^{-1}E^{T,\zeta}_{\star;Q,R,w,w'}\vert .$$

 Pour $m\in M_{0}(F)^{\geq,Q}$, on a l'\'egalit\'e $\phi^Q(H_{0}(m)-T)=\kappa^{L,T}(m)$, ce dernier terme \'etant l'analogue de $\kappa^T(m)$ quand on remplace $G$ par $L$, cf. 1.14. On a aussi trivialement  $\tilde{\sigma}_{Q}^R(H_{0}(m)-T)=\tilde{\sigma}_{Q}^R(H_{L}(m)-T)$. Alors $E^{T,\zeta}_{Q,R,w,w'}$ est l'int\'egrale sur $A_{\tilde{G}}(F)\backslash M_{0}(F)^{\geq,Q}$ du produit de $D_{0}^L(m)$ et de
  $${\bf 1}_{\zeta T}(H_{L}^{\tilde{G}}(m))\tilde{\sigma}_{Q}^R(H_{L}(m)-T)\kappa^{L,T}(m)\delta_{Q}(m)\omega_{Q,w,w'}(m)  .$$
  Cette derni\`ere fonction  est d\'efinie non seulement pour $m\in M_{0}(F)$, mais en tout point de $L(F)$. Elle est biinvariante par $K\cap L(F)$ (cela r\'esulte de la pr\'esence d'une int\'egrale sur $K\times K$ dans la d\'efinition de $\omega_{Q,w,w'}(g)$, cf. 3.4). En se rappelant la d\'efinition de la fonction $D_{0}^L(m)$, cf. 1.2, on voit que
$$E^{T,\zeta}_{Q,R,w,w'}=\int_{A_{\tilde{G}}(F)\backslash L(F)}{\bf 1}_{\zeta T}(H_{L}^{\tilde{G}}(l))\tilde{\sigma}_{Q}^R(H_{L}(l)-T)\kappa^{L,T}(l)\delta_{Q}(l)\omega_{Q,w,w'}(l)   \,dl.$$
 L'int\'egrale sur $A_{\tilde{G}}(F)\backslash L(F)$ se d\'ecompose en une int\'egrale sur $ X\in{\cal A}_{L,F}/{\cal A}_{A_{\tilde{G}},F}$ d'int\'egrales portant sur $A_{\tilde{G}}(F)\backslash A_{\tilde{G}}(F)L(F;X)$ (on rappelle que  $L(F;X)$ est l'ensemble des $l\in L(F)$ tels que $H_{L}(l)= X $). Pr\'ecis\'ement
$$E^{T,\zeta}_{Q,R,w,w'}=\int_{{\cal A}_{L,F}/{\cal A}_{A_{\tilde{G}},F}}{\bf 1}_{\zeta T}(X^{\tilde{G}})\tilde{\sigma}_{Q}^R(X-T) e^T_{Q,w,w'}(X)\,dX,$$
o\`u  
$$e^T_{Q,w,w'}(X)=\int_{ A_{\tilde{G}}(F)\backslash A_{\tilde{G}}(F)L(F;X)}\delta_{Q}(l)\kappa^{L,T}(l)\omega_{Q,w,w'}(l)\,dl.$$

 Fixons $X$. Revenons \`a la d\'efinition de $\omega_{Q,w,w'}(l)$, cf. 3.4. C'est une int\'egrale portant sur $K\times K\times i{\cal A}_{M_{disc},F}^*$. Ces int\'egrales commutent entre elles et commutent \`a l'int\'egrale  ci-dessus sur $A_{\tilde{G}}(F)\backslash A_{\tilde{G}}(F) L(F;X)$ car celle-ci est \`a support compact \`a cause de la fonction $\kappa^{L,T}(l)$. On obtient
 $$e^T_{Q,w,w'}(X)=\int_{i{\cal A}_{M_{disc},F}^*}B(\lambda)m^G(\sigma_{\lambda})\omega^T_{Q,w,w'}(X,\lambda)\,d\lambda,$$
 o\`u
 $$\omega^T_{Q,w,w'}(X,\lambda)=\int_{K\times K}\int_{ A_{\tilde{G}}(F)\backslash A_{\tilde{G}}(F) L(F;X)}\omega_{Q,w,w'}(l,k,k',\lambda)\omega(l)\delta_{Q}(l)\kappa^{L,T}(l) \,dl\,\omega(k^{-1}k')\,dk\,dk'.$$
 
 Fixons  $\lambda$ que l'on suppose provisoirement  en position g\'en\'erale. Introduisons les repr\'esentations $\rho_{w\lambda}=Ind_{w(S)\cap L}^L((w\sigma)_{w\lambda})$ et $\rho'_{w' \lambda'}=Ind_{w'(S')\cap L}^L((w'\sigma')_{w'\lambda'})$ de $L(F)$, que l'on r\'ealise dans les espaces $V^L_{w\sigma,w(S)\cap L}$ et $V^L_{w'\sigma',w'(S')\cap L}$. Pour $k,k'\in K$, introduisons des \'el\'ements $u_{w}(k',\lambda),v_{w}(k,\lambda)\in V^L_{w\sigma,w(S)\cap L}$ et $u'_{w'}(k,\lambda'),v'_{w'}(k',\lambda')\in V^L_{w'\sigma',w'(S')\cap L}$ d\'efinis par 
 $$u_{w}(k',\lambda)=(R_{\underline{Q}_{w}\vert w(S)}((w\sigma)_{w\lambda})\circ\gamma(w)\circ\pi_{\lambda}(k')u)_{K\cap L(F)},$$
 $$v_{w}(k,\lambda)=(R_{Q_{w}\vert w(S)}((w\sigma)_{w\lambda})\circ\gamma(w)\circ\pi_{\lambda}(k)v)_{K\cap L(F)},$$
 $$u'_{w'}(k,\lambda')=(R_{Q_{w'}\vert w'(S')}((w'\sigma')_{w'\lambda'})\circ\gamma(w')\circ\pi'_{\lambda'}(k)u')_{K\cap L(F)},$$
 $$v'_{w'}(k',\lambda')=(R_{\underline{Q}_{w'}\vert w'(S')}((w'\sigma')_{w'\lambda'})\circ\gamma(w')\circ\pi'_{\lambda'}(k')v')_{K\cap L(F)}.$$
 Les termes entre parenth\`eses  de ces expressions sont des \'el\'ements de $V_{w\sigma,\underline{Q}_{w}}$, resp. $V_{w\sigma,Q_{w}}$, $V_{w'\sigma',Q_{w'}}$, $V_{w'\sigma',\underline{Q}_{w'}}$. Par d\'efinition de ces espaces, ce sont des fonctions sur $K$. L'indice final $K\cap L(F)$ signifie que l'on prend leurs restrictions \`a $K\cap L(F)$. On obtient des \'el\'ements des espaces indiqu\'es. Avec ces notations, la d\'efinition de $\omega^T_{Q,w,w'}(X,\lambda)$ se r\'ecrit
 $$\omega^T_{Q,w,w'}(X,\lambda)={\bf r}_{w,w'}(\sigma_{\lambda})\int_{K\times K}\int_{A_{\tilde{G}}(F)\backslash A_{\tilde{G}}(F)L(F;X)}(v_{w}(k,\lambda),\rho_{w\lambda}(l)u_{w}(k',\lambda))$$
 $$( \rho'_{w'\lambda'}(l)v'_{w'}(k',\lambda'),u'_{w'}(k,\lambda'))\omega(l)\kappa^{L,T}(l)\,dl\,\omega(k'k^{-1})\,dk\,dk'.$$
 Notons que le $\delta_{Q}(l)$ dispara\^{\i}t dans la transition entre induites pour $G(F)$ et induites pour $L(F)$. L'int\'egrale int\'erieure est de la forme de celles consid\'er\'ees en 1.14, le groupe ambiant $G$ de ce paragraphe \'etant remplac\'e par $L$ et le tore $D$ \'etant $A_{\tilde{G}}$.   On d\'eduit des constructions de  1.14  une valeur approch\'ee de $\omega^T_{Q,w,w'}(X,\lambda)$, notons-la $r^T_{Q,w,w'}(X,\lambda)$. On l'\'etudiera plus loin. Le th\'eor\`eme 1.14 entra\^{\i}ne l'existence  d'un r\'eel $c>0$ et d'une fonction lisse et \`a croissance mod\'er\'ee $C$ sur $i{\cal A}_{M_{disc},F}^*$ de sorte que l'on ait la majoration
 $$\vert \omega^T_{Q,w,w'}(X,\lambda)-r^T_{Q,w,w'}(X,\lambda)\vert \leq {\bf r}_{w,w'}(\sigma_{\lambda})m^L((w\sigma)_{w\lambda})^{-1/2}m^L((w'\sigma')_{w'\lambda'})^{-1/2}C(\lambda)e^{-c\vert T\vert }.$$
 D\'efinissons au moins formellement
 $$r^T_{Q,w,w'}(X)=\int_{i{\cal A}_{M_{disc},F}^*}B(\lambda)m^G(\sigma_{\lambda})r^T_{Q,w,w'}(X,\lambda)\,d\lambda,$$
 $$R^{T,\zeta}_{Q,R,w,w'}=\int_{{\cal A}_{L,F}/{\cal A}_{A_{\tilde{G}},F}}{\bf 1}_{\zeta T}(X^{\tilde{G}})\tilde{\sigma}_{Q}^R(X-T) r^T_{Q,w,w'}(X)\,dX.$$
 On obtient
 $$\vert E^{T,\zeta}_{Q,R,w,w'}-R^{T,\zeta}_{Q,R,w,w'}\vert \leq e^{-c\vert T\vert }\int_{{\cal A}_{L,F}/{\cal A}_{A_{\tilde{G}},F}}{\bf 1}_{\zeta T}(X^{\tilde{G}})\tilde{\sigma}_{Q}^R(X-T)\,dX$$
 $$\int_{i{\cal A}_{M_{disc},F}^*}\vert B(\lambda)m^G(\sigma_{\lambda}){\bf r}_{w,w'}(\sigma_{\lambda})m^L((w\sigma)_{w\lambda})^{-1/2}m^L((w'\sigma')_{w'\lambda'})^{-1/2}C(\lambda)\vert \,d\lambda.$$
 Ces calculs sont justifi\'es par le r\'esultat suivant:
 
 (1) le membre de droite de l'expression ci-dessus est convergent; il est essentiellement major\'e par $\vert T\vert ^{-r}$ pour tout r\'eel $r$.
 
 A cause de la fonction ${\bf 1}_{\zeta T}(X^{\tilde{G}})$, l'int\'egrale en $X$ est convergente et essentiellement major\'ee par $\vert T\vert ^D$ pour un certain entier $D$. Il suffit donc de prouver la convergence de l'int\'egrale en $\lambda$. Puisque $B$ est de Schwartz, celle-ci r\'esulte de:
 
 (2) $\vert m^G(\sigma_{\lambda}){\bf r}_{w,w'}(\sigma_{\lambda})m^L((w\sigma)_{w\lambda})^{-1/2}m^L((w'\sigma')_{w'\lambda'})^{-1/2}\vert =1$.
 
 On rappelle que
 $${\bf r}_{w,w'}(\sigma_{\lambda})=r_{ \underline{Q}_{w}\vert Q_{w}}((w\sigma)_{w\lambda})r_{Q_{w'}\vert \underline{Q}_{w'}}((w'\sigma')_{w'\lambda'}).$$
  On a l'\'egalit\'e
 $$r_{ \overline{Q_{w}}\vert Q_{w}}((w\sigma)_{w\lambda})=r_{\overline{Q_{w}}\vert \underline{Q}_{w}}((w\sigma)_{w\lambda})r_{\underline{Q}_{w}\vert Q_{w}}((w\sigma)_{w\lambda}).$$
 D'apr\`es la d\'efinition de $\underline{Q}_{w}$ et les propri\'et\'es usuelles des facteurs de normalisation, on a
 $$r_{\overline{Q_{w}}\vert \underline{Q}_{w}}((w\sigma)_{w\lambda})=r^L_{ (\overline{w(S)\cap L})\vert (w(S)\cap L)}((w\sigma)_{w\lambda}).$$
 On  utilise 1.10(4) dans $G$ et dans $L$. Parce que $d(\sigma)=d(w\sigma)$ et $m^G(\sigma_{\lambda})=m^G((w\sigma)_{w\lambda})$, on en d\'eduit
 $$\vert r_{ \underline{Q}_{w}\vert Q_{w}}((w\sigma)_{w\lambda})\vert =m^G(\sigma_{\lambda})^{-1/2}m^L((w\sigma)_{w\lambda})^{1/2}.$$
 De m\^eme, on a
 $$\vert r_{Q_{w'}\vert \underline{Q}_{w'}}((w'\sigma')_{w'\lambda'})\vert =m^G(\sigma_{\lambda})^{-1/2}m^L((w'\sigma')_{w'\lambda'})^{1/2}.$$
 L'assertion (2) en r\'esulte, puis (1).
 
 Pour d\'emontrer le lemme, il nous reste \`a prouver l'\'egalit\'e
 $$(3) \qquad R^{T,\zeta}_{Q,R,w,w'}=mes(A_{\tilde{G}}(F)_{c})^{-1}{\bf E}^{T,\zeta}_{Q,R,w,w'}.$$
 Pour cela, il faut revenir \`a la d\'efinition du membre de gauche, et d'abord du terme $r^T_{Q,w,w'}(X,\lambda)$. On suppose $\lambda$ en position g\'en\'erale et on d\'edouble la variable $\lambda$. Plus pr\'ecis\'ement, on conserve inchang\'es les termes o\`u intervient la variable $\lambda'=\theta^{-1}\lambda$ et, dans les autres termes, on remplace $\lambda$ par $\mu\in i{\cal A}_{M_{disc},F}^*$.  En particulier, on effectue cette substitution dans la d\'efinition de ${\bf r}_{w,w'}(\sigma_{\lambda})$, qui devient un terme d\'ependant de $\lambda'$ et $\mu$, notons-le ${\bf r}_{w,w'}(\lambda',\mu)$. D'apr\`es les d\'efinitions ci-dessus, d'apr\`es celles de 1.14 et la formule (2) de ce paragraphe, $r^T_{Q,w,w'}(X,\lambda)$ est la valeur en $\mu=\lambda$ d'une expression compos\'ee de:
 
 - une int\'egrale sur $(k,k')\in K\times K$;
 
 - une somme sur les sous-groupes paraboliques $S_{1}=M_{1}U_{1}$ tels que $P_{0}\subset S_{1}\subset Q$;
 
 - une somme sur $s\in W^L(M_{1}\vert w(M_{disc}))$, $s'\in W^L(M_{1}\vert w'(M'_{disc}))$;
  
 - une somme sur $\nu\in [s'w'\sigma',\omega sw\sigma]$;
 
 le terme que l'on somme est
 
 $$(4)\qquad  C {\bf r}_{w,w'}(\lambda',\mu)\omega(k'k^{-1}) d(\sigma)^{-1}\epsilon_{S_{1}}^{Q,T}(X;sw\mu-s'w'\lambda'+\nu)$$
 $$(J^L_{(\bar{S}_{1}\cap L)\vert (s'w'(S')\cap L)}((s'w'\sigma')_{s'w'\lambda'})\circ \gamma(s')v'_{w'}(k',\lambda'),\underline{\omega}\circ A_{\nu}\circ J^L_{(\bar{S}_{1}\cap L)\vert (sw(S)\cap L)}((sw\sigma)_{sw\mu})\circ\gamma(s)u_{w}(k',\mu) )$$
 $$( \underline{\omega}\circ A_{\nu}\circ J^L_{(S_{1}\cap L)\vert  (sw(S)\cap L)}((sw\sigma)_{sw\mu})\circ\gamma(s)v_{w}(k,\mu), J^L_{(S_{1}\cap L)\vert (s'w'(S')\cap L)}((s'w'\sigma')_{s'w'\lambda'})\circ\gamma(s')u'_{w'}(k,\lambda')) ,$$
 o\`u
 $$C=mes(i{\cal A}_{M_{disc},F}^*)^{-1}mes(A_{\tilde{G}}(F)_{c})^{-1}.$$
 
  On fait dispara\^{\i}tre le terme ${\bf r}_{w,w'}( \lambda',\mu)$ en r\'etablissant les op\'erateurs d'entrelacement non normalis\'es dans les d\'efinitions des \'el\'ements $u_{w}(k',\mu)$, etc... Notons ${\bf u}_{w}(k',\mu)$ etc... les \'el\'ements d\'efinis \`a l'aide des op\'erateurs non normalis\'es. La propri\'et\'e usuelle de compatibilit\'e des op\'erateurs d'entrelacement \`a l'induction conduit \`a l'\'egalit\'e
$$J^L_{(\bar{S}_{1}\cap L)\vert (sw(S)\cap L)}((sw\sigma)_{sw\mu})\circ\gamma(s){\bf u}_{w}(k',\mu)=$$
$$(J_{\bar{S}_{1}\vert s(\underline{Q}_{w})}((sw\sigma)_{sw\mu})\circ\gamma(s) \circ J_{\underline{Q}_{w}\vert w(S)}((w\sigma)_{w\mu})\circ\gamma(w)\circ\pi_{\mu}(k')u)_{K\cap L(F)}.$$
D'apr\`es les d\'efinitions, la distance entre les deux sous-groupes paraboliques $\bar{S}_{1}$ et $sw(S)$ est la somme des distances entre $\bar{S}_{1}$ et $s(\underline{Q}_{w})$ et entre $s(\underline{Q}_{w})$ et $sw(S)$. Par composition des op\'erateurs d'entrelacement, on obtient 
$$J^L_{(\bar{S}_{1}\cap L)\vert (sw(S)\cap L)}((sw\sigma)_{sw\mu})\circ\gamma(s){\bf u}_{w}(k',\mu)=(J_{\bar{S}_{1}\vert sw(S)}((sw\sigma)_{sw\mu})\circ\gamma(sw)\circ\pi_{\mu}(k')u)_{K\cap L(F)},$$
puis

$\underline{\omega}\circ A_{\nu}\circ J^L_{(\bar{S}_{1}\cap L)\vert (sw(S)\cap L)}((sw\sigma)_{sw\mu})\circ\gamma(s){\bf u}_{w}(k',\mu)=$
$$(\underline{\omega}\circ A_{\nu}\circ J_{\bar{S}_{1}\vert sw(S)}((sw\sigma)_{sw \mu})\circ\gamma(sw)\circ\pi_{\mu}(k')u)_{K\cap L(F)}.$$
Il y a l\`a un abus d'\'ecriture: le premier $\underline{\omega}$ porte sur des fonctions sur $K\cap L(F)$, le second sur des fonctions sur $K$.
De m\^eme
$$  J^L_{(\bar{S}_{1}\cap L)\vert (s'w'(S')\cap L)}((s'w'\sigma')_{s'w'\lambda'})\circ \gamma(s'){\bf v}'_{w'}(k',\lambda')=$$
$$( J_{\bar{S}_{1}\vert s'w'(S')}((s'w'\sigma')_{s'w'\lambda'})\circ\gamma(s'w')\circ\pi'_{\lambda'}(k')v')_{K\cap L(F)}$$
  Le premier produit scalaire de l'expression (4) (transform\'e comme on l'a dit ci-dessus) s'\'ecrit donc
$$\int_{K\cap L(F)}(( J_{\bar{S}_{1}\vert s'w'(S')}((s'w'\sigma')_{s'w'\lambda'})\circ\gamma(s'w')\circ\pi'_{\lambda'}(k')v')(h),$$
$$(\underline{\omega}\circ A_{\nu}\circ J_{\bar{S}_{1}\vert sw(S)}((sw\sigma)_{sw\mu})\circ\gamma(sw)\circ\pi_{\mu}(k')u)(h))\,dh.$$
On peut commuter l'op\'erateur $\pi'_{\lambda'}(k')$ aux op\'erateurs qui le pr\'ec\`edent et remplacer le premier terme par 
$$( J_{\bar{S}_{1}\vert s'w'(S')}((s'w'\sigma')_{s'w'\lambda'})\circ\gamma(s'w')v')(hk').$$
 On peut commuter l'op\'erateur $\pi_{\mu}(k')$ aux op\'erateurs qui le pr\'ec\`edent, mais, d'apr\`es la d\'efinition de $\underline{\omega}$,  la commutation \`a cet op\'erateur fait sortir un terme $\omega(k')^{-1}$. On peut donc remplacer le deuxi\`eme terme par
 $$\omega(k')^{-1}(\underline{\omega}\circ A_{\nu}\circ J_{\bar{S}_{1}\vert sw(S)}((sw\sigma)_{sw\mu})\circ\gamma(sw)u)(hk').$$
  On se rappelle que dans (4) figure une multiplication par $\omega(k')$ et que l'on doit int\'egrer tout cela en $k'\in K$. Apr\`es ces op\'erations, le premier produit scalaire de (4) devient plus simplement le produit scalaire
 $$(J_{\bar{S}_{1}\vert s'w'(S')}((s'w'\sigma')_{s'w'\lambda'})\circ\gamma(s'w')v' ,\underline{\omega}\circ A_{\nu}\circ J_{\bar{S}'\vert sw(S)}((sw\sigma)_{sw\mu})\circ\gamma(sw)u) ,$$
 ou encore
 $$(A_{\nu}^{-1}\circ \underline{\omega}^{-1}\circ J_{\bar{S}_{1}\vert s'w'(S')}((s'w'\sigma')_{s'w'\lambda'})\circ\gamma(s'w')v' ,  J_{\bar{S}'\vert sw(S)}((sw\sigma)_{sw\mu})\circ\gamma(sw)u) .$$
 
On traite de m\^eme le second produit scalaire. On obtient que  $r^T_{Q,w,w'}(X,\lambda)$ est la valeur en $\mu=\lambda$ d'une expression compos\'ee de sommes sur les $S_{1}$, $s$, $s'$, $\nu$ de termes
$$(5)\qquad  Cd(\sigma)^{-1}\epsilon_{S_{1}}^{Q,T}(X;sw\mu-s'w'\lambda'+\nu)$$
$$(A_{\nu}^{-1}\circ\underline{\omega}^{-1}\circ J_{\bar{S}_{1}\vert s'w'(S')}((s'w'\sigma')_{s'w'\lambda'})\circ\gamma(s'w')v' , J_{\bar{S}'\vert sw(S)}((sw\sigma)_{sw\mu})\circ\gamma(sw)u) $$
 
$$( J_{S_{1}\vert sw(S)}((sw\sigma)_{sw\mu})\circ\gamma(sw)v,A_{\nu}^{-1}\circ\underline{\omega}^{-1}\circ J_{S_{1}\vert s'w'(S')}((s'w'\sigma')_{s'w'\lambda'})\circ\gamma(s'w')u').$$
 L'ensemble de sommation $[s'w'\sigma';\omega sw\sigma]$  est celui des $\nu\in i{\cal A}_{M_{1},F}^{*}$ tels que $s'w'\sigma'\simeq \omega(sw\sigma)_{-\nu}$. Posons $t=w^{-1}s^{-1}s'w'$. C'est un \'el\'ement de $W^G(M_{disc}\vert M'_{disc})$. L'application
$$\begin{array}{ccc}[t\sigma';\omega \sigma]&\to&[s'w'\sigma';\omega sw\sigma]\\ \nu&\mapsto &sw\nu\\ \end{array}$$
est bijective. Plus pr\'ecis\'ement, remarquons que toutes les repr\'esentations $\sigma$, $t\sigma'$, $sw\sigma$ et $s'w'\sigma'$ se r\'ealisent naturellement dans  le m\^eme espace $V_{\sigma}$.  Soit $\nu\in [t \sigma';\omega \sigma]$, introduisons un automorphisme unitaire $A_{\nu}$ de $V_{\sigma}$ tel que $(t \sigma')(x)\circ A_{\nu}=\omega(x)A_{\nu}\circ \sigma_{-\nu}(x)$ pour tout $x\in M_{disc}(F)$. Alors $A_{\nu}$ v\'erifie aussi $  s'w'\sigma'(x)\circ A_{\nu}=\omega(x) A_{\nu}\circ(sw\sigma)_{-sw\nu}(x)$ pour tout $x\in M_{1}(F)$. On peut donc remplacer  l'ensemble de sommation $[s'w'\sigma';\omega sw\sigma]$ par $[t \sigma';\omega \sigma]$, en rempla\c{c}ant   $\nu$ par $sw\nu$, tout en prenant $A_{sw\nu}=A_{\nu}$.  Rappelons que les op\'erateurs intervenant dans (5) sont en fait d\'eduits par fonctorialit\'e de ces op\'erateurs $A_{\nu}$. On v\'erifie l'\'egalit\'e
$$A_{\nu}^{-1}\circ\underline{\omega}^{-1}\circ J_{\bar{S}_{1}\vert s'w'(S')}((s'w'\sigma')_{s'w'\lambda'})\circ\gamma(s'w')=J_{\bar{S}_{1}\vert s'w'(S')}((sw\sigma)_{s'w'\lambda'-sw\nu})\circ\gamma(s'w')\circ A_{\nu}^{-1}\circ\underline{\omega}^{-1}.$$
 
Comme en 3.16, on introduit l'op\'erateur
$$A(t,\nu;\Lambda)=R_{S\vert t(S')}(\sigma_{\Lambda})\circ\gamma(t)\circ A_{\nu}^{-1}\circ\underline{\omega}^{-1}:Ind_{S'}^G(\sigma'_{t^{-1}(\Lambda+\nu)})\to Ind_{S}^G(\sigma_{\Lambda}).$$
On calcule
$$J_{\bar{S}_{1}\vert s'w'(S')}((sw\sigma)_{s'w'\lambda'-sw\nu})\circ\gamma(s'w')\circ A_{\nu}^{-1}\circ\underline{\omega}^{-1}=r_{\bar{S}_{1}\vert s'w'(S')}((sw\sigma)_{s'w'\lambda'-sw\nu})$$
$$R_{\bar{S}_{1}\vert s'w'(S')}((sw\sigma)_{s'w'\lambda'-sw\nu})\circ\gamma(s'w')\circ A_{\nu}^{-1}\circ\underline{\omega}^{-1}$$
$$=r_{\bar{S}_{1}\vert s'w'(S')}((sw\sigma)_{s'w'\lambda'-sw\nu})R_{\bar{S}_{1}\vert sw(S)}((sw\sigma)_{s'w'\lambda'-sw\nu})\circ \gamma(sw)\circ A(t,\nu;t\lambda'-\nu)$$
$$=r_{\bar{S}_{1}\vert s'w'(S')}((sw\sigma)_{s'w'\lambda'-sw\nu})R_{ sw(S)\vert \bar{S}_{1}}((sw\sigma)_{s'w'\lambda'-sw\nu})^{-1}\circ \gamma(sw)\circ A(t,\nu;t\lambda'-\nu)$$
$$=r_{\bar{S}_{1}\vert s'w'(S')}((sw\sigma)_{s'w'\lambda''-sw\nu})r_{ sw(S)\vert \bar{S}_{1}}((sw\sigma)_{s'w'\lambda'-sw\nu})$$
$$J_{ sw(S)\vert \bar{S}_{1}}((sw\sigma)_{s'w'\lambda'-sw\nu})^{-1}\circ \gamma(sw)\circ A(t,\nu;t\lambda'-\nu).$$
La deuxi\`eme \'egalit\'e n\'ecessite que nos rel\`evements v\'erifient $t=w^{-1}s^{-1}s'w$, ce que l'on peut supposer.  Un calcul analogue vaut si l'on remplace $\bar{S}'$ par $S'$. En utilisant ensuite les propri\'et\'es d'adjonction des op\'erateurs d'entrelacement, on voit que $r^T_{Q,w,w'}(X,\lambda)$ est la valeur en $\mu=\lambda$ d'une expression compos\'ee de sommes sur les $S_{1}$, $s$, $s'$ comme pr\'ec\'edemment et sur $\nu\in [t\sigma',\omega \sigma]$ de termes
$$(6)\qquad  C\epsilon_{S_{1}}^{Q,T}(X;sw\mu-s'w'\lambda'+sw\nu)r(S_{1},s,s';s'w'\lambda'-sw\nu)$$
$$(A(t,\nu;t\lambda'-\nu)v',\gamma(w^{-1}s^{-1})\circ J_{ \bar{S}_{1}\vert sw(S)}((sw\sigma)_{s'w'\lambda'-sw\nu})^{-1}\circ J_{\bar{S}_{1}\vert sw(S)}((sw\sigma)_{sw\mu})\circ\gamma(sw)u )$$
$$(   \gamma(w^{-1}s^{-1})\circ J_{S_{1}\vert sw(S)}((sw\sigma)_{s'w'\lambda'-sw\nu})^{-1}\circ J_{S'\vert sw(S)}((sw\sigma)_{sw\mu})\circ\gamma(sw)v,A(t,\nu,t\lambda'-\nu)u'),$$
o\`u, pour $\xi\in i{\cal A}_{M_{S'},F}^*$, on a pos\'e
$$r(S_{1},s,s';\xi)=d(\sigma)^{-1}r_{S_{1}\vert  s'w'S'}((sw\sigma)_{\xi})r_{ sw(S)\vert S_{1}}((sw\sigma)_{\xi})\overline{r_{\bar{S}_{1}\vert s'w'S'}((sw\sigma)_{\xi})}\overline{r_{sw(S)\vert \bar{S}_{1}}((sw\sigma)_{\xi})}.$$
En utilisant 1.10(3), on obtient
$$r(S_{1},s,s';\xi)=d(\sigma)^{-1}r_{S_{1}\vert \bar{S}_{1}}((sw\sigma)_{\xi})r_{\bar{S}_{1}\vert S_{1}}((sw\sigma)_{\xi}).$$
Puisque $d(\sigma)=d(sw\sigma)$, ceci n'est autre que $m^G((sw\sigma)_{\xi})^{-1}$, lui-m\^eme \'egal \`a $m^G(\sigma_{w^{-1}s^{-1}\xi})$. Pour $\xi=s'w'\lambda'-sw\nu$, on a $\sigma_{w^{-1}s^{-1}\mu}=\sigma_{t\lambda'-\nu}\simeq \omega^{-1}(t\sigma')_{t\lambda'}=\omega^{-1} t((\sigma_{\lambda})\circ\theta)$. La mesure de Plancherel est invariante par automorphisme et on v\'erifie facilement qu'elle est aussi invariante par torsion par un caract\`ere unitaire. On obtient $r(S_{1},s,s';s'w'\lambda'-sw\nu)=m^G(\sigma_{\lambda})^{-1}$.
 
On peut  remplacer la sommation sur $s'$ par une sommation sur $t=w^{-1}s^{-1}s'w'$. Cet \'el\'ement d\'ecrit l'ensemble $W^G(M_{disc}\vert M'_{disc})\cap w^{-1}W^Lw'$, en identifiant pour simplifier le premier ensemble \`a un ensemble de repr\'esentants dans $W^G$.  On peut remplacer la somme en $S_{1}$ et $s$ par une somme sur $S''=w^{-1}s^{-1}(S_{1})$. Ce parabolique d\'ecrit l'ensemble ${\cal P}^{w^{-1}(Q)}(M_{disc})$ des $S''\in {\cal P}(M_{disc})$ tels que $S''\subset w^{-1}(Q)$.  On a les \'egalit\'es
$$\gamma(w^{-1}s^{-1})\circ J_{ \bar{S}_{1}\vert sw(S)}((sw\sigma)_{s'w'\lambda'-sw\nu})^{-1}\circ J_{\bar{S}_{1}\vert sw(S)}((sw\sigma)_{sw\mu})\circ\gamma(sw)=$$
$$J_{ \bar{S}''\vert S}(\sigma_{t\lambda'-\nu})^{-1}\circ J_{\bar{S}''\vert S}(\sigma_{\mu}),$$
$$\gamma(w^{-1}s^{-1})\circ J_{ S_{1}\vert sw(S)}((sw\sigma)_{s'w'\lambda'-sw\nu})^{-1}\circ J_{S_{1}\vert sw(S)}((sw\sigma)_{sw\mu})\circ\gamma(sw)=$$
$$J_{ S''\vert S}(\sigma_{t\lambda'-\nu})^{-1}\circ J_{S''\vert S}(\sigma_{\mu}),$$
$$\epsilon_{S_{1}}^{Q,T}(X,sw\mu-s'w'\lambda'+sw\nu)=\epsilon_{S''}^{w^{-1}(Q),T[S'']}(w^{-1}X,\mu-t\lambda'+\nu).$$
En utilisant les d\'efinitions de 3.16, l'expression (6) se r\'ecrit
$$ C m^G(\sigma_{\lambda})^{-1}\epsilon_{S''}^{w^{-1}(Q),T[S'']}(w^{-1}X,\mu-t\lambda'+\nu)\phi(t,\nu;\lambda,\mu-t\lambda'+\nu,S'').$$
En sommant sur $S''$, on obtient
$$ C m^G(\sigma_{\lambda})^{-1}\phi_{M_{disc}}^{w^{-1}(Q),T}(t,\nu,w^{-1}(X);\lambda,\mu-t\lambda'+\nu).$$
A ce point, le lemme 3.16 nous autorise \`a \'egaler $\mu$ \`a $\lambda$. On obtient finalement
$$r^T_{Q,w,w'}(X,\lambda)=C m^G(\sigma_{\lambda})^{-1}\sum_{t\in W^G(M_{disc}\vert M'_{disc})\cap w^{-1}W^Lw'}$$
$$\sum_{\nu\in [t\sigma',\omega \sigma]}\phi_{M_{disc}}^{w^{-1}(Q),T}(t,\nu,w^{-1}(X);\lambda,\lambda-t\lambda'+\nu).$$
En appliquant les d\'efinitions, on obtient ensuite
$$R^{T,\zeta}_{Q,R,w,w'}=C \sum_{t\in W^G(M_{disc}\vert M'_{disc})\cap w^{-1}W^Lw'}\sum_{\nu\in [t\sigma',\omega \sigma]}
\int_{{\cal A}_{L,F}/{\cal A}_{A_{\tilde{G}},F}}{\bf 1}_{\zeta T}(X^{\tilde{G}})\tilde{\sigma}_{Q}^R(X-T)$$
$$\int_{i{\cal A}_{M_{disc},F}^*}B(\lambda)\phi_{M_{disc}}^{w^{-1}(Q),T}(t,\nu,w^{-1}(X);\lambda,\lambda-t\lambda'+\nu)\,d\lambda\,dX.$$
On a $\tilde{\sigma}_{Q}^R(X-T)=\tilde{\sigma}_{w^{-1}(Q)}^{w^{-1}(R)}(w^{-1}X-T[w^{-1}(Q)])$. Par le changement de variables $X\mapsto wX$, la double int\'egrale multipli\'ee par $C$ devient $mes(A_{\tilde{G}}(F)_{c})^{-1}{\bf E}^T_{w^{-1}(Q),w^{-1}(R),t,\nu}$ et on obtient l'\'egalit\'e cherch\'ee
$$R^{T,\zeta}_{Q,R,w,w'}= mes(A_{\tilde{G}}(F)_{c})^{-1}E^T_{\star;Q,R,w,w'}.$$
Cela ach\`eve la preuve. $\square$
 
 \bigskip
 
 \subsection{Une nouvelle expression approchant $j^T$}
 Pour $Q=LU_{Q}\in {\cal F}(M_{disc})$ et $t\in W^G(M_{disc}\vert M'_{disc})$, d\'efinissons une fonction $H\mapsto S_{Q}(t;H)$ sur ${\cal A}_{0}$ par
 $$S_{Q}(t;H)=\sum_{R; Q\subset R}s_{Q}^R(t)\tilde{\sigma}_{Q}^R(H).$$
 En d\'evissant la d\'efinition de $s_{Q}^R(t)$, on a aussi
 $$S_{Q}(t;H)=\sum_{\tilde{P}=\tilde{M}U_{P}; Q\subset P, t\theta^{-1}(P)=P}(-1)^{dim(\mathfrak{a}_{\tilde{M}})-dim(\mathfrak{a}_{\tilde{G}})}\sum_{R; P\subset R}\tilde{\sigma}_{Q}^R(H),$$
 ou encore, en utilisant 2.2(1),
 $$S_{Q}(t;H)=\sum_{\tilde{P}=\tilde{M}U_{P}; Q\subset P,t\theta^{-1}(P)=P}(-1)^{dim(\mathfrak{a}_{\tilde{M}})-dim(\mathfrak{a}_{\tilde{G}})}\tau_{Q}^P(H)\hat{\tau}_{\tilde{P}} (H).$$
 
 Pour $\nu\in [t\sigma',\omega \sigma]$, posons 
 $${\bf E}_{t,\nu}^T=mes(i{\cal A}_{M_{disc},F}^*)^{-1}\sum_{Q=LU_{Q}\in {\cal F}(M_{disc})}\int_{{\cal A}_{L,F}/{\cal A}_{A_{\tilde{G}},F}}S_Q(t;X-T[Q])$$
 $$\int_{i{\cal A}_{M_{disc},F}^*}B(\lambda)\phi_{M_{disc}}^{Q,T}(t,\nu,X;\lambda;\lambda-t\lambda'+\nu)\,d\lambda\,dX.$$
 Ceci n'est autre que
 $$\sum_{Q,R\in {\cal F}(M_{disc}); Q\subset R}s_{Q}^R(t){\bf E}_{Q,R,t,\nu}^T$$
 qui est une expression convergente d'apr\`es le lemme  3.18. Posons
 $${\bf E}^T= \sum_{t\in W(M_{disc}\vert M'_{disc})}\sum_{\nu\in [t\sigma',\omega \sigma]}{\bf E}_{t,\nu}^T.$$
 
  \ass{Lemme}{On a la majoration
 $$\vert j^T-mes(A_{\tilde{G}}(F)_{c})^{-1}{\bf E}^T\vert <<\vert T\vert ^{-r}$$
 pour tout r\'eel $r$.}
 
 Preuve. La proposition 3.4 nous dit que l'on peut aussi bien remplacer $j^T$ par $j^T_{\star}$, c'est-\`a-dire par
 $$ \sum_{Q=LU_{Q},R; P_{0}\subset Q\subset R}\sum_{w\in W^G(L\vert S),w'\in W^G(L\vert S')}s_{Q}^R(w,w')E^T_{Q,R,w,w'}.$$
 Le lemme 3.18 nous dit que l'on peut aussi bien remplacer $E^T_{Q,R,w,w'}$ par $mes(A_{\tilde{G}}(F)_{c})^{-1}E^T_{\star;Q,R,w,w'}$. L'expression ci-dessus est alors remplac\'ee par
 $$mes(A_{\tilde{G}}(F)_{c})^{-1}\sum_{Q=LU_{Q},R; P_{0}\subset Q\subset R}\sum_{w\in W^G(L\vert S),w'\in W^G(L\vert S')}s_{Q}^R(w,w')$$
 $$\sum_{t\in W^G(M_{disc}\vert M'_{disc})\cap w^{-1}W^Lw'}\sum_{\nu\in [t\sigma',\omega \sigma]}{\bf E}^T_{w^{-1}(Q),w^{-1}(R),t,\nu}.$$
 Comme on l'a remarqu\'e en 3.18, pour $Q,...,t$ intervenant ci-dessus, on a l'\'egalit\'e $s_{Q}^R(w,w')=s_{w^{-1}(Q)}^{w^{-1}(R)}(t)$. On peut donc r\'ecrire l'expression pr\'ec\'edente sous la forme
 $$mes(A_{\tilde{G}}(F)_{c})^{-1}\sum_{t\in W(M_{disc}\vert M'_{disc})}\sum_{\nu\in [t\sigma',\omega \sigma]}E^T_{\star;t,\nu},$$
 o\`u
 $$E^T_{\star;t,\nu}=\sum_{Q=LU_{Q},R; P_{0}\subset Q\subset R}\sum_{w\in W^G(L\vert S),w'\in W^G(L\vert S'); t\in w^{-1}W^Lw'}s_{w^{-1}(Q)}^{w^{-1}(R)}(t){\bf E}^T_{w^{-1}(Q),w^{-1}(R),t,\nu}.$$
 Fixons $t$ et $\nu$. Pour tous $Q,R$ et tout $w\in W^G(L\vert S)$ il y a un et un seul $w'\in W^G(L\vert S')$ tel que $t\in w^{-1}W^Lw'$: en identifiant tous ces \'el\'ements \`a des rel\`evements dans $W^G$, $w'$ est l'\'el\'ement de longueur minimale dans la classe $W^Lwt$. On peut donc supprimer la somme en $w'$ et la condition $t\in w^{-1}W^Lw'$. L'application
 $$\begin{array}{ccc}\{Q=LU_{Q},R,w; P_{0}\subset Q\subset R, w\in W^G(L\vert S)\}&\to& \{Q',R'; M_{disc}\subset Q'\subset R'\}\\ (Q,R,w)&\mapsto& (w^{-1}(Q),w^{-1}(R))\\ \end{array}$$
 est bijective. On voit alors que $E^T_{\star;t,\nu}={\bf E}^T_{t,\nu}$. Cela d\'emontre le lemme. $\square$
 
 Dans toutes nos expressions intervient la fonction $B$ fix\'ee en 3.4. Il peut \^etre utile de pr\'eciser la d\'ependance en $B$ de la majoration du lemme. Pour cela, notons  pour quelques instants $j^T(B)$ et ${\bf E}^T(B)$ les termes not\'es pr\'ec\'edemment $j^T$ et ${\bf E}^T$. Introduisons un ensemble ${\cal N}$ de semi-normes sur l'espace des fonctions de Schwartz sur $i{\cal A}_{M_{disc},F}^*$. Dans le cas o\`u $F$ est non-archim\'edien, ces semi-normes sont
 $$B\mapsto sup\{ \vert (XB)(\lambda)\vert ;\lambda\in i{\cal A}_{M_{disc},F}^*\},$$
 o\`u $X$ parcourt les op\'erateurs diff\'erentiels \`a coefficients constants sur $i{\cal A}_{M_{disc}}^*$. Dans le cas o\`u $F$ est archim\'edien, ce sont les
 $$B\mapsto sup\{\vert (XB)(\lambda)\vert (1+\vert \lambda\vert )^N; \lambda\in i{\cal A}_{M_{disc}}^*\},$$
 o\`u $X$ parcourt les m\^emes op\'erateurs que ci-dessus et $N$ parcourt les entiers naturels. Le lemme se pr\'ecise en
 
 (1) pour tout r\'eel $r\geq1$, il existe $c_{r}>0$ et un sous-ensemble fini ${\cal N}_{r}\subset {\cal N}$ de sorte que l'on ait la majoration
 $$\vert j^T(B)-mes(A_{\tilde{G}}(F)_{c})^{-1}{\bf E}^T(B)\vert \leq c_{r}sup\{\underline{n}(B); \underline{n}\in {\cal N}_{r}\}\vert T\vert ^{-r}$$
 pour tout $T$ et toute fonction de Schwartz $B$ sur $i{\cal A}_{M_{disc},F}^*$.
 
 Il suffit de reprendre patiemment toutes nos majorations. A chaque fois que $B$ y intervient, c'est par l'interm\'ediaire d'une semi-norme appartenant \`a l'ensemble ${\cal N}$ et il est clair que, pour $r$ fix\'e, il n'y a qu'un nombre fini de telles semi-normes qui interviennent.
 \bigskip
 
 \subsection{ Un lemme d'inversion de Fourier}
  Dans ce paragraphe et jusqu'en 3.23, on fixe un \'element $t\in W^G(M_{disc}\vert M'_{disc})$ et un \'el\'ement $\nu\in [t\sigma',\omega \sigma]$. On d\'efinit un Levi $M_{t}$: c'est le plus grand Levi contenant $M_{disc}$ tel que ${\cal A}_{M_{t}}$ contienne l'ensemble
  $$\{H\in {\cal A}_{M_{disc}}; t\theta^{-1}H=H\}.$$
  On pose $\tilde{M}_{t}=M_{t}\gamma_{0}t^{-1}$, o\`u on rel\`eve $t$ en un \'el\'ement de $G$. D'apr\`es 2.1(3) et (5), $\tilde{M}_{t}$ est un ensemble de Levi et l'ensemble des espaces paraboliques $\tilde{P}$ tels que $M_{disc}\subset P$ et $t\theta^{-1}(P)=P$ n'est autre que ${\cal F}(\tilde{M}_{t})$.

Consid\'erons l'ensemble des couples $(\lambda,\Lambda)\in i{\cal A}_{M_{disc}}^*\times i{\cal A}_{\tilde{M}_{t}}^*$ tels que $\Lambda-\lambda+t\lambda'-\nu\in i{\cal A}_{M_{disc},F}^{\vee}$. Il est invariant par translation par $(i{\cal A}_{M_{disc},F}^{\vee}+i{\cal A}_{\tilde{M}_{t}}^*)\times i{\cal A}_{\tilde{M}_{t},F}^{\vee}$.  Consid\'erons  son quotient par l'action de ce groupe. Parce que $1-t\theta^{-1}$ se restreint en un automorphisme de $i{\cal A}_{M_{disc}}^{\tilde{M}_{t},*}$, on voit que ce quotient est fini, r\'eduit \`a un \'el\'ement si $F$ est archim\'edien. On voit aussi  que la projection $(\lambda,\Lambda)\to \lambda$ est injective. Autrement dit, si on note $\{\nu\}_{t}$ l'image de cette projection, il existe une application $\lambda\mapsto \Lambda(\lambda)$ de $\{\nu\}_{t}$ dans $i{\cal A}_{\tilde{M}_{t},F}^*$ de sorte que le quotient ci-dessus soit \'egal \`a $\{(\lambda,\Lambda(\lambda)); \lambda\in \{\nu\}_{t}\}$.

\ass{Lemme}{Soit $f$ une fonction de Schwartz sur $i{\cal A}_{M_{disc},F}^*$ et soit $X\in {\cal A}_{\tilde{M}_{t},F}$. L'int\'egrale
$$mes(i{\cal A}_{M_{disc},F}^*)^{-1}\int_{{\cal A}_{M_{disc},F}^{\tilde{M}_{t}}(X)}\int_{i{\cal A}_{M_{disc},F}^*}e^{<\lambda-t\lambda'+\nu,Y>}f(\lambda)\,d\lambda\,dY$$
est convergente dans cet ordre. Elle est \'egale \`a
$$mes(i{\cal A}_{\tilde{M}_{t},F}^*)^{-1}\vert det((1-t\theta^{-1})_{\vert {\cal A}_{M_{disc}}^{\tilde{M}_{t},*}})\vert ^{-1}\sum_{\lambda\in \{\nu\}_{t}}e^{<\Lambda(\lambda),X>}\int_{i{\cal A}_{\tilde{M}_{t},F}^*}f(\lambda+\xi)\,d\xi.$$}

Preuve. On suppose $F$ non-archim\'edien, le cas archim\'edien \'etant plus simple. On fixe un rel\`evement de $\nu$ dans $i{\cal A}_{M_{disc}}^*$, que l'on note encore $\nu$. On fixe $\mu\in i{\cal A}_{M_{disc}}^*$ tel que $\nu^{\tilde{M}_{t}}=t\theta^{-1}\mu-\mu$. On fixe $X'\in {\cal A}_{M_{disc},F}$ tel que $(X')^{\tilde{M}_{t}}=X$. On a ${\cal A}_{M_{disc},F}^{\tilde{M}_{t}}(X)=X'+{\cal A}_{M_{disc},F}^{\tilde{M}_{t}}$ (o\`u ${\cal A}_{M_{disc},F}^{\tilde{M}_{t}}={\cal A}_{M_{disc},F}\cap {\cal A}_{M_{disc}}^{\tilde{M}_{t}}$) et l'int\'egrale de l'\'enonc\'e se r\'ecrit
$$mes(i{\cal A}_{M_{disc},F}^*)^{-1}\sum_{Y\in{\cal A}_{M_{disc},F}^{\tilde{M}_{t}}}e^{<\nu_{\tilde{M}_{t}},X'>}\int_{i{\cal A}_{M_{disc},F}^*}e^{<(\lambda-\mu)-t\theta^{-1}(\lambda-\mu),Y+X'>}f(\lambda)\,d\lambda.$$
 Toutes les int\'egrations se font sur des groupes compacts. Une int\'egrale sur $i{\cal A}_{M_{disc},F}^*$ se d\'ecompose en le produit d'une int\'egrale sur $
 i{\cal A}_{M_{disc}}^*/(i{\cal A}_{M_{disc},F}^{\vee}+i{\cal A}_{\tilde{M}_{t}}^*)$ et d'une int\'egrale sur $i{\cal A}_{\tilde{M}_{t},F}^*$ (avec des pond\'erations provenant de nos choix de mesures).  L'expression ci-dessus se transforme en
$$(1) \qquad mes(i{\cal A}_{M_{disc},F}^*)^{-1}mes(i{\cal A}_{\tilde{M}_{t},F}^*)\sum_{Y\in {\cal A}_{M_{disc},F}^{\tilde{M}_{t}}}e^{<\nu_{\tilde{M}_{t}},X'>}$$
$$\int_{i{\cal A}_{M_{disc}}^*/(i{\cal A}_{M_{disc},F}^{\vee}+i{\cal A}_{\tilde{M}_{t}}^*)}e^{<(\lambda-\mu)-t\theta^{-1}(\lambda-\mu),Y+X'>} f'(\lambda)\,d\lambda,$$
o\`u
$$f'(\lambda)=mes(i{\cal A}_{\tilde{M}_{t},F}^*)^{-1}\int_{i{\cal A}_{\tilde{M}_{t},F}^*}f(\lambda+\xi)\,d\xi.$$
La fonction $f'$ est lisse sur $ i{\cal A}_{M_{disc}}^*/(i{\cal A}_{M_{disc},F}^{\vee}+i{\cal A}_{\tilde{M}_{t}}^*)$. L'accouplement
$$\begin{array}{ccc}( i{\cal A}_{M_{disc}}^*/(i{\cal A}_{M_{disc},F}^{\vee}+i{\cal A}_{\tilde{M}_{t}}^*))\times {\cal A}_{M_{disc},F}^{\tilde{M}_{t}}&\to &{\mathbb C}^{\times}\\ (\lambda,Y)&\mapsto&e^{<\lambda-t\theta^{-1}\lambda,Y>}\\ \end{array}$$
n'est pas parfait. Son noyau ${\cal K}$ est l'ensemble des $\lambda\in i{\cal A}_{M_{disc}}^*/(i{\cal A}_{M_{disc},F}^{\vee}+i{\cal A}_{\tilde{M}_{t}}^*)$  tels que $\lambda-t\theta^{-1}\lambda\in i{\cal A}_{M_{disc},F}^{\vee}+i{\cal A}_{\tilde{M}_{t}}^*$.  A ce d\'efaut pr\`es, la formule (1) est une inversion de Fourier. Elle est donc convergente dans l'ordre indiqu\'e. On la calcule selon la formule usuelle. Elle vaut
$$\vert {\cal K}\vert ^{-1}e^{<\nu_{\tilde{M}_{t}},X'>}\sum_{\lambda\in \mu+{\cal K}}e^{<(\lambda-\mu)-t\theta^{-1}(\lambda-\mu),X'>}f'(\lambda).$$
Notons $p:i{\cal A}_{M_{disc}}^*\to i{\cal A}_{M_{disc}}^{\tilde{M}_{t},*}$ la projection orthogonale. Via cette projection, on v\'erifie que ${\cal K}$ s'identifie au noyau de $1-t\theta^{-1}$ agissant dans $i{\cal A}_{M_{disc}}^{\tilde{M}_{t},*}/p(i{\cal A}_{M_{disc},F}^{\vee})$. Donc $\vert {\cal K}\vert =\vert det((1-t\theta^{-1})_{\vert {\cal A}_{M_{disc}}^{\tilde{M}_{t},*}})\vert $. D'autre part, en comparant les d\'efinitions, on voit que $\mu+{\cal K}=\{\nu\}_{t}$. Notre expression vaut donc
$$(2) \qquad \vert det((1-t\theta^{-1})_{\vert {\cal A}_{M_{disc}}^{\tilde{M}_{t},*}})\vert ^{-1}e^{<\nu_{\tilde{M}_{t}},X'>}\sum_{\lambda\in \{\nu\}_{t}}e^{<(\lambda-\mu)-t\theta^{-1}(\lambda-\mu),X'>}f'(\lambda).$$
Pour $\lambda\in \{\nu\}_{t}$, on introduit comme ci-dessus $\Lambda(\lambda)\in {\cal A}_{\tilde{M}_{t},F}^*$ tel que $\Lambda(\lambda)-\lambda+t\lambda'-\nu\in i{\cal A}_{M_{disc},F}^{\vee}$, ou encore $\Lambda(\lambda)-\lambda+\mu+t\theta^{-1}(\lambda-\mu)-\nu_{\tilde{M}_{t}}\in i{\cal A}_{M_{disc},F}^{\vee}$. Puisque $X'\in {\cal A}_{M_{disc},F}$, on a 
 $$e^{<(\lambda-\mu)-t\theta^{-1}(\lambda-\mu),X'>}=e^{<\Lambda(\lambda)-\nu_{\tilde{M}_{t}},X'>},$$
 d'o\`u
 $$e^{<\nu_{\tilde{M}_{t}},X'>}e^{<(\lambda-\mu)-t\theta^{-1}(\lambda-\mu),X'>}=e^{<\Lambda(\lambda),X'>}=e^{<\Lambda(\lambda),X>}.$$
 L'expression (2) se transforme en celle de l'\'enonc\'e. $\square$

 \bigskip

 \subsection{Transformation de ${\bf E}^T_{t,\nu}$}
  Dans la d\'efinition de ${\bf E}^T_{t,\nu}$ intervient la $(G,M_{disc})$-famille $(\phi(t,\nu;\lambda,\Lambda,S''))_{S''\in {\cal P}(M_{disc})}$ de 3.16. Les fonctions qui composent celle-ci ont le d\'efaut d'avoir des singularit\'es et, quand $F$ est archim\'edien, de ne pas \^etre de Schwartz en $\Lambda$: elles sont seulement \`a croissance mod\'er\'ee. Modifions cette famille. On fixe une fonction $C$ sur $i{\cal A}_{M_{disc},F}^*$, qui est $C^{\infty}$ et \`a support compact et telle que $C(0)=1$. On pose
 $$\phi_{\star}(t,\nu;\lambda,\Lambda,S'')=r_{\bar{S}\vert S}(\sigma_{t\lambda'-\nu+\Lambda})^{-1}r_{\bar{S}\vert S}(\sigma_{\lambda})C(\Lambda-\lambda'+t\lambda'-\nu)B(\lambda)\phi(t,\nu;\lambda,\Lambda,S'').$$
 On a aussi
 $$\phi_{\star}(t,\nu;\lambda,\Lambda,S'')= \epsilon(t,\nu;\lambda)C(\Lambda-\lambda'+t\lambda'-\nu)B(\lambda)\phi_{reg}(t,\nu;\lambda,\Lambda,S'').$$
  Les lemmes 1.6 et  3.16(ii)  montrent que $\phi_{\star}(t,\nu;\lambda,\Lambda,S'')$ est $C^{\infty}$ en les deux variables $\lambda,\Lambda\in i{\cal A}_{M_{disc},F}^*$. Les propri\'et\'es de $B$ et $C$ assurent qu'elle est de Schwartz en ces deux variables. D'apr\`es le lemme 3.16(iii), on a l'\'egalit\'e
$$\phi_{\star,M_{disc}}^{Q,T}(t,\nu,X;\lambda,\lambda-t\lambda'+\nu)=B(\lambda)\phi_{M_{disc}}^{Q,T}(t,\nu,X;\lambda,\lambda-t\lambda'+\nu) $$
pour tout $Q=LU_{Q}\in {\cal F}(M_{disc})$ et tout $X\in {\cal A}_{L,F}$. La d\'efinition de ${\bf E}^T_{t,\nu}$ se r\'ecrit
$${\bf E}^T_{t,\nu}=mes(i{\cal A}_{M_{disc},F}^*)^{-1}\sum_{Q=LU_{Q}\in {\cal F}(M_{disc})}\int_{{\cal A}_{L,F}\backslash {\cal A}_{A_{\tilde{G}},F}}S_Q(t;X-T[Q])$$
$$\int_{i{\cal A}_{M_{disc},F}^*}\phi_{\star,M_{disc}}^{Q,T}(t,\nu,X;\lambda,\lambda-t\lambda'+\nu)\,d\lambda\,dX.$$
En utilisant le lemme 1.8, on peut \'ecrire
$$\phi_{\star,M_{disc}}^{Q,T}(t,\nu,X;\lambda,\lambda-t\lambda'+\nu)=\sum_{Q'=L'U_{Q'}; M_{disc}\subset Q'\subset Q}\int_{{\cal A}_{L',F}}\int_{{\cal A}_{M_{disc},F}^L(X+H_{L}) }$$
$$\delta_{M_{disc}}^{Q'}(H')\Gamma_{Q'}^Q(H',H+T[Q'])\hat{\phi}_{\star}(t,\nu;\lambda,H,Q') e^{<\lambda-t\lambda'+\nu,H'>}\,dH'\,dH,$$
o\`u
$$\hat{\phi}_{\star}(t,\nu;\lambda,H,Q')=mes(i{\cal A}_{L',F}^*)^{-1}\int_{i{\cal A}_{L',F}^*}\phi_{\star}(t,\nu;\lambda,\Lambda,Q')e^{-<\Lambda,H>}\,d\Lambda.$$
Fixons $Q,Q'\in {\cal F}(M_{disc})$, avec $Q'\subset Q$. 

\ass{Lemme}{(i) Pour $X\in {\cal A}_{L,F}$, l'int\'egrale
$$\int_{i{\cal A}^*_{M_{disc},F}}\int_{{\cal A}_{L',F}}\int_{{\cal A}_{M_{disc},F}^L(X+H_{L})}$$
$$\vert\delta_{M_{disc}}^{Q'}(H') \Gamma_{Q'}^Q(H',H+T[Q'])\hat{\phi}_{\star}(t,\nu;\lambda,H,Q')\vert  \,dH'\,dH\,d\lambda$$
est convergente.

(ii) L'int\'egrale
$$\int_{{\cal A}_{L,F}\backslash {\cal A}_{A_{\tilde{G}},F}} \int_{{\cal A}_{L',F}}\int_{{\cal A}_{M_{disc},F}^L(X+H_{L})}\left\vert \delta_{M_{disc}}^{Q'}(H')S_Q(t;X-T[Q])\Gamma_{Q'}^Q(H',H+T[Q'])\right.$$
$$ \int_{i{\cal A}_{M_{disc},F}^*}\left.\hat{\phi}_{\star}(t,\nu;\lambda,H,Q') e^{<\lambda-t\lambda'+\nu,H'>}\,d\lambda\right\vert \,dH'\,dH\,dX$$
est convergente.}

Preuve de (i). Les deux conditions $\delta_{M_{disc}}^{Q'}(H')\Gamma_{Q'}^Q(H',H+T[Q'])=1$ et $H'\in {\cal A}_{M_{disc},F}^L(X+H_{L})$ entra\^{\i}nent une majoration $\vert H'\vert <<1+\vert H\vert$, les \'el\'ements $T$ et $X$ \'etant ici consid\'er\'es comme constants. L'int\'egrale en $H'$ est donc convergente et est essentiellement born\'ee par $(1+\vert H\vert )^D$ pour un certain entier $D$. Puisque $\phi_{\star}(t,\nu;\lambda,\Lambda,Q')$ est de Schwartz en les deux variables $\lambda$ et $\Lambda$, $\hat{\phi}_{\star}(t,\nu;\lambda,H,Q')$ est de Schwartz en les deux variables $\lambda$ et $H$. Alors l'int\'egrale
$$\int_{i{\cal A}^*_{M_{disc},F}}\int_{{\cal A}_{L',F}}(1+\vert H\vert )^D\vert \hat{\phi}_{\star}(t,\nu;t\lambda'-\nu,H,Q')\vert\,dH\,d\lambda$$
est convergente, ce qui prouve (i).

Preuve de (ii). L'int\'egrale int\'erieure en $\lambda$ est le produit de $e^{<\nu,H'>}$, qui est de valeur absolue $1$, et de la transform\'ee de Fourier en $\lambda$ de $\hat{\phi}_{\star}(t,\nu;\lambda,H,Q')$, \'evalu\'ee au point $(\theta t^{-1}-1)H'$. Puisque $\hat{\phi}_{\star}(t,\nu;\lambda,H,Q')$ est de Schwartz en $\lambda$ et $H$, on peut fixer des fonctions de Schwartz $h$ sur ${\cal A}_{L',F}$ et $h'$ sur ${\cal A}_{M_{disc},F}$, \`a valeurs positives, de sorte que la valeur absolue de cette int\'egrale int\'erieure soit major\'ee par $h'((\theta t^{-1}-1)H')h(H)$. Les conditions $H'\in {\cal A}_{M_{disc},F}^L(X+H_{L})$ et $\delta_{M_{disc}}^{Q'}(H')=1$ entra\^{\i}nent que $H'=X+H_{L}+(H')^L_{L'}$. La condition $\Gamma_{Q'}^Q(H',H+T[Q'])$ implique une majoration $\vert (H')^L_{L'}\vert <<1+\vert H^L\vert $, l'\'el\'ement $T$ \'etant consid\'er\'e comme constant. On en d\'eduit 
$$(\theta t^{-1}-1)H'=(\theta t^{-1}-1)X+Y(H'),$$
avec $\vert Y(H')\vert <<1+\vert H\vert $. D'autre part, on a une majoration
$$1+\vert U+V\vert >>(1+\vert U\vert )(1+\vert V\vert )^{-1}$$
pour tous $U,V\in {\cal A}_{M_{disc}}$. Pour tout r\'eel $r>0$, on a une majoration
$$h'((\theta t^{-1}-1)H')<<(1+\vert (\theta t^{-1}-1)H'\vert )^{-r}.$$
En utilisant les majorations pr\'ec\'edentes, on obtient
$$h'((\theta t^{-1}-1)H')<<(1+\vert (\theta t^{-1}-1)X\vert )^{-r}(1+\vert H\vert )^r.$$
A ce point, on a montr\'e que l'expression du (ii) de l'\'enonc\'e \'etait major\'ee par
$$\int_{{\cal A}_{L,F}/ {\cal A}_{A_{\tilde{G}},F}}\int_{{\cal A}_{L',F}}\int_{{\cal A}_{M_{disc},F}^L(X+H_{L})}\vert \delta_{M_{disc}}^{Q'}(H')S_Q(t;X-T[Q])\Gamma_{Q'}^Q(H',H+T[Q'])\vert $$
$$(1+\vert (\theta t^{-1}-1)X\vert )^{-r}(1+\vert H\vert )^rh(H)\,dH'\,dH\,dX.$$
Par des raisonnements d\'ej\`a faits, l'int\'egrale en $H'$ est convergente et essentiellement major\'ee par $(1+\vert H\vert )^D$ pour un certain entier $D$. L'int\'egrale en $H$ est convergente quel que soit $r$ puisque $h$ est de Schwartz. L'expression ci-dessus est donc essentiellement major\'ee par
$$\int_{{\cal A}_{L,F}/ {\cal A}_{A_{\tilde{G}},F}}\vert S_Q(t;X-T[Q])\vert (1+\vert (\theta t^{-1}-1)X\vert )^{-r} \,dX.$$
Il reste \`a montrer que l'on peut choisir $r$ tel que cette int\'egrale soit convergente. En revenant \`a la d\'efinition de $S_{Q}(t;X)$, on voit qu'il suffit de fixer un sous-groupe parabolique $R$ tel que $Q\subset R$ et $s_{Q}^R(t)\not=0$ et de prouver la  m\^eme assertion pour l'int\'egrale
$$\int_{{\cal A}_{L,F}/ {\cal A}_{A_{\tilde{G}},F}}\tilde{\sigma}_{Q}^R(X-T[Q])(1+\vert (\theta t^{-1}-1)X\vert )^{-r}\,dX.$$
Or il r\'esulte de 3.17(1) que, pour $\tilde{\sigma}_{Q}^R(X-T[Q])=1$, on a une majoration
$$1+\vert X^{\tilde{G}}\vert<<1+ \vert (\theta t^{-1}-1)X\vert .$$
L'int\'egrale ci-dessus est donc essentiellement major\'ee par
$$\int_{{\cal A}_{L,F}/ {\cal A}_{A_{\tilde{G}},F}}(1+\vert X^{\tilde{G}}\vert )^{-r}\,dX.$$
Pour $r$ assez grand, ceci est convergent. $\square$

Ce lemme  nous autorise \`a r\'ecrire ${\bf E}_{t,\nu}^T$ sous la forme
$$(1) \qquad {\bf E}_{t,\nu}^T=mes(i{\cal A}_{M_{disc},F}^*)^{-1}\sum_{Q=LU_{Q},Q'=L'U_{Q'}\in {\cal F}(M_{disc}); Q'\subset Q}\int_{{\cal A}_{L,F}/{\cal A}_{A_{\tilde{G}},F}}$$
$$\int_{{\cal A}_{L',F}}\int_{{\cal A}_{M_{disc},F}^L(X+H_{L})}\delta_{M_{disc}}^{Q'}(H')S_Q(t;X-T[Q])\Gamma_{Q'}^Q(H',H+T[Q'])$$
$$\int_{i{\cal A}_{M_{disc},F}^*}e^{<\lambda-t\lambda'+\nu,H'>}\hat{\phi}_{\star}(t,\nu;\lambda,H,Q')\,d\lambda\,dH'\,dH\,dX,$$
o\`u on peut permuter librement la somme en $Q$, $Q'$ et les trois premi\`eres int\'egrales. 

\bigskip

\subsection{Calcul de ${\bf E}_{t,\nu}^T$}

On a d\'efini en 3.16 une $(G,M)$-famille $(\phi_{reg}(t,\nu;\lambda,\Lambda,S'))_{S'\in {\cal P}(M_{disc})}$. Il s'en d\'eduit une $(\tilde{G},\tilde{M}_{t})$-famille $(\phi_{reg}(t,\nu;\lambda,\Lambda,\tilde{P}))_{\tilde{P}\in {\cal P}(\tilde{M}_{t})}$, cf. 2.3 . Les \'el\'ements $\lambda$ et $\Lambda$ appartiennent ici \`a $i{\cal A}_{M_{disc},F}^*$ et $i{\cal A}_{\tilde{M}_{t},F}^*$, $\lambda$ jouant un r\^ole de param\`etre. Pour $X\in {\cal A}_{\tilde{G},F}$, on construit la fonction
$$\phi^{\tilde{G},T}_{reg,\tilde{M}_{t}}(t,\nu,X;\lambda,\Lambda)=\sum_{\tilde{P}\in {\cal P}(\tilde{M}_{t})}\phi_{reg}(t,\nu;\lambda,\Lambda,\tilde{P})\epsilon_{\tilde{P}}^{\tilde{G},T[\tilde{P}]}(X;\Lambda).$$
C'est une fonction lisse en $\lambda$ et $\Lambda$. Si $F$ est archim\'edien, toutes ses d\'eriv\'ees sont \`a croissance mod\'er\'ee.  

\ass{Proposition}{On a l'\'egalit\'e
$${\bf E}_{t,\nu}^T=mes(i{\cal A}_{\tilde{M}_{t},F}^*)^{-1}\vert det((1-t\theta^{-1})_{\vert {\cal A}_{M_{disc}}^{\tilde{M}_{t},*}})\vert ^{-1}\sum_{\lambda\in \{\nu\}_{t}}\sum_{X\in {\cal A}_{\tilde{G},F}/{\cal A}_{A_{\tilde{G},F}}}$$
$$\int_{i{\cal A}_{\tilde{M}_{t},F}^*}B(\lambda+\xi)\epsilon(t,\nu;\lambda+\xi)\phi^{\tilde{G},T}_{reg,\tilde{M}_{t}}(t,\nu,X;\lambda+\xi,\Lambda(\lambda))\,d\xi.$$}

{\bf Remarque.} On v\'erifie que, pour $\lambda$ intervenant ci-desus, on a $e^{<\Lambda(\lambda),Z>}=1$ pour tout $Z\in {\cal A}_{A_{\tilde{G}},F}$. La fonction que l'on somme en $X$ est donc bien invariante par ${\cal A}_{A_{\tilde{G}},F}$.
 
\bigskip

Preuve. On part de la formule (1) du paragraphe pr\'ec\'edent.  Permutons les int\'egrales en $X$ et $H$. Effectuons ensuite le changement de variables $X\mapsto X-H_{L}$. L'int\'egrale en $H'$ devient une int\'egrale sur ${\cal A}_{M_{disc},F}^L(X)$. La composition de l'int\'egrale en $X$ et de cette int\'egrale en $H'$ devient une unique int\'egrale en $X\in   {\cal A}_{M_{disc},F}/{\cal A}_{A_{\tilde{G}},F}$.  Cela conduit \`a l'\'egalit\'e
$${\bf E}_{t,\nu}^T=mes(i{\cal A}_{M_{disc},F}^*)^{-1}\sum_{Q=LU_{Q},Q'=L'U_{Q'}\in {\cal F}(M_{disc}); Q'\subset Q}\int_{{\cal A}_{L',F}}\int_{{\cal A}_{M_{disc},F}/{\cal A}_{A_{\tilde{G}},F}}\delta_{M_{disc}}^{Q'}(X)$$
$$S_Q(t;X-H_{L}-T[Q])\Gamma_{Q'}^Q(X,H+T[Q'])\int_{i{\cal A}_{M_{disc},F}^*}e^{<\lambda-t\lambda'+\nu,X>}\hat{\phi}_{\star}(t,\nu;\lambda,H,Q')\,d\lambda\,dX\,dH.$$
On peut toujours permuter librement les sommes en $Q$ et $Q'$ et les deux premi\`eres int\'egrales, d'o\`u
 $${\bf E}_{t,\nu}^T=mes(i{\cal A}_{M_{disc},F}^*)^{-1}\int_{{\cal A}_{M_{disc},F}/{\cal A}_{A_{\tilde{G}},F}}\sum_{Q'=L'U_{Q'}\in {\cal F}(M_{disc})}\int_{{\cal A}_{L',F}}\delta_{M_{disc}}^{Q'}(X)\sum_{Q=LU_{Q}; Q'\subset Q}$$
 $$S_Q(t;X-H_{L}-T[Q])\Gamma_{Q'}^Q(X,H+T[Q'])\int_{i{\cal A}_{M_{disc},F}^*}e^{<\lambda-t\lambda'+\nu,X>}\hat{\phi}_{\star}(t,\nu;\lambda,H,Q')\,d\lambda\,dH\,dX.$$
Fixons $Q'$.  En se reportant \`a la d\'efinition de $S_{Q}(t,X)$ et en utilisant  3.20(3), on a l'\'egalit\'e
$$\sum_{Q=LU_{Q}; Q'\subset Q}S_Q(t;X-H_{L}-T[Q])\Gamma_{Q'}^Q(X,H+T[Q'])=$$
$$\sum_{\tilde{P}\in {\cal F}(\tilde{M}_{t}); Q'\subset P}(-1)^{a_{\tilde{P}}-a_{\tilde{G}}} \hat{\tau}_{\tilde{P}}(X-H-T[\tilde{P}])\sum_{Q=LU_{Q}; Q'\subset Q\subset P}\tau_{Q}^P(X-H-T[Q])\Gamma_{Q'}^Q(X,H+T[Q']).$$
D'apr\`es 1.3(4), la derni\`ere somme en $Q$ vaut $\tau_{Q'}^P(X)$. D'o\`u
$${\bf E}_{t,\nu}^T=\int_{{\cal A}_{M_{disc},F}/{\cal A}_{A_{\tilde{G}},F}}{\bf E}_{t,\nu}^T(X)\,dX,$$
o\`u
$${\bf E}_{t,\nu}^T(X)=mes(i{\cal A}_{M_{disc},F}^*)^{-1}\sum_{Q'=L'U_{Q'}\in {\cal F}(M_{disc})}\int_{{\cal A}_{L',F}}\delta_{M_{disc}}^{Q'}(X)\sum_{\tilde{P}=\tilde{M}U_{P}\in {\cal F}(\tilde{M}_{t}); Q'\subset P}$$
$$(-1)^{a_{\tilde{P}}-a_{\tilde{G}}} \hat{\tau}_{\tilde{P}}(X-H-T[\tilde{P}])\tau_{Q'}^P(X)\int_{i{\cal A}_{M_{disc},F}^*}e^{<\lambda-t\lambda'+\nu,X>}\hat{\phi}_{\star}(t,\nu;\lambda,H,Q')\,d\lambda\,dH.$$
Fixons $X$. Puisque $\hat{\phi}_{\star}(t,\nu;\lambda,H,Q')$ est de Schwartz en $\lambda$ et $H$,  l'expression  ${\bf E}^T_{t,\nu}(X)$ est absolument convergente. On peut donc \'ecrire
$${\bf E}^T_{t,\nu}(X)=mes(i{\cal A}_{M_{disc},F}^*)^{-1}\int_{i{\cal A}_{M_{disc},F}^*}e^{<\lambda-t\lambda'+\nu,X>}{\bf E}_{t,\nu}^T(X,\lambda)\,d\lambda,$$
o\`u
$${\bf E}_{t,\nu}^T(X,\lambda)=\sum_{Q'=L'U_{Q'}\in {\cal F}(M_{disc})}\int_{{\cal A}_{L',F}}\delta_{M_{disc}}^{Q'}(X)\sum_{\tilde{P}=\tilde{M}U_{P}\in {\cal F}(\tilde{M}_{t}); Q'\subset P}(-1)^{a_{\tilde{P}}-a_{\tilde{G}}} \hat{\tau}_{\tilde{P}}(X-H-T[\tilde{P}])$$
$$\tau_{Q'}^P(X)\hat{\phi}_{\star}(t,\nu;\lambda,H,Q')\,dH.$$
Fixons $\lambda$. On peut permuter les sommes en $\tilde{P}$ et en $Q'$, puis d\'ecomposer l'int\'egrale en $H\in {\cal A}_{L',F}$ en une int\'egrale sur $H\in {\cal A}_{\tilde{M},F}$ et une int\'egrale en $H'\in {\cal A}_{L',F}^{\tilde{M}}(H)$ et on peut permuter la premi\`ere avec la somme en $Q'$  On obtient
$${\bf E}_{t,\nu}^T(X,\lambda) =\sum_{\tilde{P}=\tilde{M}U_{P}\in {\cal F}(\tilde{M}_{t})}(-1)^{a_{\tilde{P}}-a_{\tilde{G}}}\int_{{\cal A}_{\tilde{M},F}}\hat{\tau}_{\tilde{P}}(X-H-T[\tilde{P}])$$
$$\sum_{Q'=L'U_{Q'}; M_{disc}\subset Q'\subset P}\delta_{M_{disc}}^{Q'}(X)\tau_{Q'}^P(X)\int_{{\cal A}_{L',F}^{\tilde{M}}(H)}\hat{\phi}_{\star}(t,\nu;\lambda,H',Q')\,dH'\,dH.$$
En se rappelant la d\'efinition de $\hat{\phi}_{\star}(t,\nu;\lambda,H',Q')$ et par inversion de Fourier partielle, on a
$$\int_{{\cal A}_{L',F}^{\tilde{M}}(H)}\hat{\phi}_{\star}(t,\nu;\lambda,H',Q')\,dH'\,=mes(i{\cal A}_{\tilde{M},F})^{-1}\int_{i{\cal A}_{\tilde{M},F}^*}\phi_{\star}(t,\nu;\lambda,\Lambda,Q')e^{-<\Lambda,H>}\,d\Lambda.$$
Puisque $\Lambda$ ne parcourt plus que $i{\cal A}_{\tilde{M},F}^*$, on  peut aussi bien remplacer 
$\phi_{\star}(t,\nu;\lambda,\Lambda,Q')$ par $\phi_{\star}(t,\nu;\lambda,\Lambda,\tilde{P})$ et l'expression ci-dessus est \'egale \`a $\hat{\phi}_{\star}(t,\nu;\lambda,H,\tilde{P})$. Notons que $Q'$ a ici disparu donc, dans l'expression pr\'ec\'edente de ${\bf E}_{\nu,t}^T$, la somme en $Q'$ n'est plus que
$$\sum_{Q'; M_{disc}\subset Q'\subset P}\delta_{M_{disc}}^{Q'}(X)\tau_{Q'}^P(X).$$
Celle-ci est \'egale \`a
$$\sum_{\tilde{P}'; \tilde{M}_{t}\subset \tilde{P}'\subset \tilde{P}}\delta_{\tilde{M}_{t}}^{\tilde{P}'}(X)\tau_{\tilde{P}'}^{\tilde{P}}(X).$$
En effet, les deux expressions sont \'egales \`a $1$ d'apr\`es 1.3(1) et sa variante tordue.
  On obtient
$${\bf E}_{t,\nu}^T(X,\lambda)= \sum_{\tilde{P}=\tilde{M}U_{P}\in {\cal F}(\tilde{M}_{t})}(-1)^{a_{\tilde{P}}-a_{\tilde{G}}}\int_{{\cal A}_{\tilde{M},F}}\hat{\tau}_{\tilde{P}}(X-H-T[\tilde{P}])$$
$$\sum_{\tilde{P}'=\tilde{M}'U_{P'};\tilde{M}_{t}\subset\tilde{P}'\subset \tilde{P}}\delta_{\tilde{M}_{t}}^{\tilde{P}'}(X)\tau_{\tilde{P}'}^{\tilde{P}}(X) \hat{\phi}_{\star}(t,\nu;\lambda,H,\tilde{P})\,dH.$$
On remonte maintenant le calcul ci-dessus: on permute les sommes en $\tilde{P}$ et $\tilde{P}'$;  on remplace $\hat{\phi}_{\star}(t,\nu;\lambda,H,\tilde{P})$ par l'expression \'egale
$$\int_{{\cal A}_{\tilde{M}',F}^{\tilde{M}}(H)}\hat{\phi}_{\star}(t,\nu;\lambda,H',\tilde{P}')\,dH';$$
\`a partir des int\'egrales en $H\in {\cal A}_{\tilde{M},F}$ et $H'\in {\cal A}_{\tilde{M}',F}^{\tilde{M}}(H)$, on reconstitue une int\'egrale en $H\in {\cal A}_{\tilde{M}',F}$. On obtient
$${\bf E}^T_{t,\nu}(X,\lambda)=\sum_{\tilde{P}'=\tilde{M}'U_{P'}\in {\cal F}(\tilde{M}_{t})}\int_{{\cal A}_{\tilde{M}',F}}\delta_{\tilde{M}_{t}}^{\tilde{P}'}(X)$$
$$\sum_{\tilde{P}; \tilde{P}'\subset \tilde{P}}(-1)^{a_{\tilde{P}}-a_{\tilde{G}}}\hat{\tau}_{\tilde{P}}(X-H-T[\tilde{P}'])\tau_{\tilde{P}'}^{\tilde{P}}(X)\hat{\phi}_{\star}(t,\nu;\lambda,H,\tilde{P}')\,dH.$$
Par d\'efinition,
$$\sum_{\tilde{P}; \tilde{P}'\subset \tilde{P}}(-1)^{a_{\tilde{P}}-a_{\tilde{G}}}\hat{\tau}_{\tilde{P}}(X-H-T[\tilde{P}'])\tau_{\tilde{P}'}^{\tilde{P}}(X)=\Gamma_{\tilde{P}'}^{\tilde{G}}(X,H+T[P']).$$
Cela fait dispara\^{\i}tre les $\tilde{P}$ de la formule ci-dessus et nous autorise \`a abandonner les $'$ des $\tilde{P}'$. D'o\`u
$${\bf E}^T_{t,\nu}(X,\lambda)=\sum_{\tilde{P}=\tilde{M}U_{P}\in {\cal F}(\tilde{M}_{t})}\int_{{\cal A}_{\tilde{M},F}}\delta_{\tilde{M}_{t}}^{\tilde{P}}(X)\Gamma_{\tilde{P}}^{\tilde{G}}(X,H+T[\tilde{P}])\hat{\phi}_{\star}(t,\nu;\lambda,H,\tilde{P})\,dH.$$

Revenons \`a ${\bf E}^T_{t,\nu}$ qui est l'int\'egrale en $\lambda$ puis $X$ de l'expression ci-dessus, multipli\'ee par $mes(i{\cal A}_{M_{disc},F}^*)^{-1}e^{<\lambda-t\lambda'+\nu,X>}$. On peut d\'ecomposer l'int\'egrale en $X$ en une int\'egrale en $X\in {\cal A}_{\tilde{M}_{t},F}/{\cal A}_{A_{\tilde{G}},F}$ et une int\'egrale en $Y\in {\cal A}_{M_{disc},F}^{\tilde{M}_{t}}(X)$. Remarquons que pour $Y$ dans cet ensemble, on a ${\bf E}^T_{t,\nu}(Y,\lambda)={\bf E}^T_{t,\nu}(X,\lambda)$. On a alors
$${\bf E}^T_{t,\nu}=mes(i{\cal A}_{M_{disc},F}^*)^{-1}\int_{{\cal A}_{\tilde{M}_{t},F}/{\cal A}_{A_{\tilde{G}},F}}\int_{{\cal A}_{M_{disc},F}^{\tilde{M}_{t}}(X)}\int_{i{\cal A}_{M_{disc},F}^*}e^{<\lambda-t\lambda'+\nu,Y>}{\bf E}^T_{t,\nu}(X,\lambda)\,d\lambda\,dY\,dX.$$
Pour tout $X$, la fonction $\lambda\mapsto {\bf E}^T_{t,\nu}(X,\lambda)$ est de Schwartz. La double int\'egrale int\'erieure est calcul\'ee par le lemme 3.20. On obtient
$${\bf E}^T_{t,\nu}=mes(i{\cal A}_{\tilde{M}_{t},F}^*)^{-1}\vert det((1-t\theta^{-1})_{\vert {\cal A}_{M_{disc}}^{\tilde{M}_{t}*}}\vert ^{-1}\int_{{\cal A}_{\tilde{M}_{t},F}/{\cal A}_{A_{\tilde{G}},F}}\sum_{\lambda\in \{\nu\}_{t}}e^{<\Lambda(\lambda),X>}$$
$$\int_{i{\cal A}_{\tilde{M}_{t},F}^*}{\bf E}^T_{t,\nu}(X,\lambda+\xi)\,d\xi\,dX.$$
En d\'eveloppant le dernier terme, on obtient
$${\bf E}^T_{t,\nu}=mes(i{\cal A}_{\tilde{M}_{t},F}^*)^{-1}\vert det((1-t\theta^{-1})_{\vert {\cal A}_{M_{disc}}^{\tilde{M}_{t}*}}\vert ^{-1}\int_{{\cal A}_{\tilde{M}_{t},F}/{\cal A}_{A_{\tilde{G}},F}}\sum_{\lambda\in \{\nu\}_{t}}e^{<\Lambda(\lambda),X>}\int_{i{\cal A}_{\tilde{M}_{t},F}^*}$$
$$\sum_{\tilde{P}=\tilde{M}U_{P}\in {\cal F}(\tilde{M}_{t})}\int_{{\cal A}_{\tilde{M},F}}\delta_{\tilde{M}_{t}}^{\tilde{P}}(X)\Gamma_{\tilde{P}}^{\tilde{G}}(X,H+T[\tilde{P}])\hat{\phi}_{\star}(t,\nu;\lambda+\xi,H,\tilde{P})\,dH\,d\xi\,dX.$$
Cette expression est absolument convergente. En effet, pour $\tilde{P}$ fix\'e, l'int\'egrale en $X$ est \`a support compact et est essentiellement born\'ee par $(1+\vert H\vert )^D$ pour un entier $D$ convenable. Les int\'egrales restantes en $H$ et $\xi$ sont convergentes puisque $\hat{\phi}_{\star}(t,\nu;\lambda+\xi,H,\tilde{P})$ est de Schwartz en $\xi$ et $H$. On peut donc commencer par int\'egrer en $\xi$, puis en $H$, puis en $X$. On d\'ecompose ensuite l'int\'egrale en $X$ en une somme sur $X\in {\cal A}_{\tilde{G},F}/{\cal A}_{A_{\tilde{G}},F}$ d'int\'egrales en $Y\in {\cal A}_{\tilde{M}_{t},F}^{\tilde{G}}(X+H_{\tilde{G}})$. Apr\`es encore quelques permutations, on obtient
$${\bf E}^T_{t,\nu}=mes(i{\cal A}_{\tilde{M}_{t},F}^*)^{-1}\vert det((1-t\theta^{-1})_{\vert {\cal A}_{M_{disc}}^{\tilde{M}_{t}*}}\vert ^{-1}\sum_{\lambda\in \{\nu\}_{t}}\sum_{X\in {\cal A}_{\tilde{G},F}/{\cal A}_{A_{\tilde{G}},F}}\int_{i{\cal A}_{\tilde{M}_{t},F}^*}\sum_{\tilde{P}=\tilde{M}U_{P}\in {\cal F}(\tilde{M}_{t})}$$
$$\int_{{\cal A}_{\tilde{M},F}}\int_{{\cal A}_{\tilde{M}_{t},F}^{\tilde{G}}(X+H_{\tilde{G}})}e^{<\Lambda(\lambda),Y>}\delta_{\tilde{M}_{t}}^{\tilde{P}}(Y)\Gamma_{\tilde{P}}^{\tilde{G}}(Y,H+T[\tilde{P}])\hat{\phi}_{\star}(t,\nu;\lambda+\xi,H,\tilde{P})\,dY\,dH\,d\xi.$$
L'expression int\'erieure (somme en $\tilde{P}$ et int\'egrales en $Y$ et $H$) est calcul\'ee par la variante tordue du  lemme 1.8. C'est $\phi_{\star,\tilde{M}_{t}}^{\tilde{G},T}(t,\nu,X;\lambda+\xi,\Lambda(\lambda))$. D'o\`u
$${\bf E}^T_{t,\nu}=mes(i{\cal A}_{\tilde{M}_{t},F}^*)^{-1}\vert det((1-t\theta^{-1})_{\vert {\cal A}_{M_{disc}}^{\tilde{M}_{t}*}}\vert ^{-1}\sum_{\lambda\in \{\nu\}_{t}}\sum_{X\in {\cal A}_{\tilde{G},F}/{\cal A}_{A_{\tilde{G}},F}}$$
$$\int_{i{\cal A}_{\tilde{M}_{t},F}^*}\phi_{\star,\tilde{M}_{t}}^{\tilde{G},T}(t,\nu,X;\lambda+\xi,\Lambda(\lambda))\,d\xi.$$
En revenant aux d\'efinitions des fonctions $\phi_{reg}(t,\nu;\lambda,\Lambda,S')$ et $\phi_{\star}(t,\nu;\lambda,\Lambda,S')$, on voit que la seconde est le produit de la premi\`ere et de $\epsilon(t,\nu;\lambda)B(\lambda)C( \Lambda-\lambda+t\lambda'-\nu)$. Mais, pour $\lambda\in \{\nu\}_{t}$ et $\xi\in i{\cal A}_{\tilde{M}_{t},F}^*$, on a par d\'efinition $\Lambda(\lambda)-\lambda-\xi+t\theta^{-1}(\lambda+\xi)-\nu\in i{\cal A}_{M_{disc},F}^{\vee}$. Donc $C( \Lambda(\lambda)-\lambda-\xi+t\theta^{-1}(\lambda+\xi)-\nu)=1$ et 
$$\phi_{\star,\tilde{M}_{t}}^{\tilde{G},T}(t,\nu,X;\lambda+\xi,\Lambda(\lambda))=B(\lambda+\xi)\epsilon(t,\nu;\lambda+\xi)\phi_{reg,\tilde{M}_{t}}^{\tilde{G},T}(t,\nu,X;\lambda+\xi,\Lambda(\lambda)).$$
La formule ci-dessus devient celle de l'\'enonc\'e. $\square$

\bigskip
\subsection{Le "terme constant" de ${\bf E}_{t,\nu}^T$}
Posons
$$j_{spec,t,\nu}=\vert det((1-t\theta^{-1})_{\vert {\cal A}_{M_{disc}}^{\tilde{M}_{t},*}})\vert ^{-1} \sum_{\lambda\in\{\nu\}_{t}; \Lambda(\lambda)=0}\int_{i{\cal A}_{\tilde{M}_{t},F}^*}B(\lambda+\xi)\phi_{reg,\tilde{M}_{t}}^{\tilde{G}}(t,\nu;\lambda+\xi,0)\,d\xi.$$
Cette expression est absolument convergente, la fonction $\xi\mapsto \phi_{reg,\tilde{M}_{t}}^G(t,\nu;\lambda+\xi,0)$ \'etant lisse et \`a croissance mod\'er\'ee. 

\ass{Lemme}{Il existe une unique fonction $T\mapsto f(T)$ qui appartient \`a $PolExp$ et qui co\"{\i}ncide avec $mes({\cal A}_{\tilde{G}}(F)_{c})^{-1}{\bf E}_{t,\nu}^T$ dans le c\^one o\`u celle-ci est d\'efinie. Si $F$ est archim\'edien, on a $c_{0}(f)=j_{spec,t,\nu}$. Si $F$ est non-archim\'edien, pour tout r\'eseau ${\cal R}\subset {\cal A}_{M_{0},F}\otimes_{{\mathbb Z}}{\mathbb Q}$, on a l'\'egalit\'e
$$lim_{k\to \infty}c_{\frac{1}{k}{\cal R},0}(f)=j_{spec,t,\nu}.$$}

{\bf Remarque.} Pour simplifier, on appellera "terme constant" de $f$ le terme $c_{0}(f)$ si $F$ est archim\'edien, $lim_{k\to \infty}c_{\frac{1}{k}{\cal R},0}(f)$ si $F$ est non-archim\'edien.
\bigskip

Preuve. On suppose $F$ non-archim\'edien, le cas archim\'edien \'etant similaire. Consid\'erons la formule de la proposition pr\'ec\'edente. Pour tous $\lambda$, $X$ et pour tout $\xi\in i{\cal A}_{\tilde{M}_{t},F}^*$, la fonction $T\mapsto f_{\lambda,X,\xi}(T)=\phi^T_{reg,\tilde{M}_{t}}(t,\nu,X;\lambda+\xi,\Lambda(\lambda))$ est d\'efinie  pour  $T\in {\cal A}_{M_{0},F}\otimes_{{\mathbb Z}}{\mathbb Q}$  et elle appartient \`a $PolExp$, cf. lemme 1.7. Plus pr\'ecis\'ement, elle appartient \`a un espace $PolExp_{\boldsymbol{\Xi},N}$, o\`u $\boldsymbol{\Xi}$ et $N$ ne d\'ependent pas de $\xi$. Autrement dit, pour tout r\'eseau ${\cal R}\subset {\cal A}_{M_{0},F}\otimes_{{\mathbb Z}}{\mathbb Q}$, on peut \'ecrire
$$ f_{\lambda,X,\xi}(T)=\sum_{\mu\in {\cal X}_{{\cal R}}}e^{<\mu,T>}p_{{\cal R},\mu}(T)$$
pour $T\in {\cal R}$, avec un ensemble ${\cal X}_{{\cal R}}$ ind\'ependant de $\xi$ et des polyn\^omes $p_{{\cal R},\mu}$ de degr\'e born\'e ind\'ependamment de $\xi$. Les coefficients de ces polyn\^omes se calculent par interpolation et v\'erifient donc les m\^emes propri\'et\'es que la fonction $ f_{\lambda,X,\xi}$ elle-m\^eme. Ils sont donc $C^{\infty}$ en $\xi$. Il en r\'esulte que le d\'eveloppement en $T$ commute \`a l'int\'egrale en $\xi$. Cela implique que, si on note $f(T)$ le membre de droite de l'\'egalit\'e de la proposition 3.22 multipli\'e par $mes({\cal A}_{\tilde{G}}(F)_{c})^{-1}$, la fonction $T\mapsto f(T)$ appartient \`a $PolExp$ et que son coefficient $c_{{\cal R},0}(f)$ se calcule en rempla\c{c}ant $ f_{\lambda,X,\xi}$ par son coefficient $c_{{\cal R},0}(f_{\lambda,X,\xi})$ dans la formule int\'egrale. Remarquons que la norme (au sens de 1.7) de la $(\tilde{G},\tilde{M}_{t})$-famille $(\phi_{reg}(t,\nu;\lambda+\xi,\Lambda,\tilde{P}))_{\tilde{P}\in {\cal P}(\tilde{M}_{t})}$ est born\'ee ind\'ependamment de $\xi$. Il r\'esulte donc du lemme 1.7 que $lim_{k\to \infty}c_{\frac{1}{k}{\cal R},0}(f)$ se calcule en rempla\c{c}ant $f_{\lambda,X,\xi}$ par $lim_{k\to \infty}c_{\frac{1}{k}{\cal R},0}(f_{\lambda,X,\xi})$ dans la formule int\'egrale  (multipli\'ee par $mes({\cal A}_{\tilde{G}}(F)_{c})^{-1}$). Cette derni\`ere limite est $0$ si $\Lambda(\lambda)\not\in ( i{\cal A}_{\tilde{M}_{t},F}^{\vee}+i{\cal A}_{\tilde{G}}^*)/ i{\cal A}_{\tilde{M}_{t},F}^{\vee}$. Si $\Lambda(\lambda)\in(\Lambda_{1}(\lambda)+i{\cal A}_{\tilde{M}_{t},F}^{\vee})/ i{\cal A}_{\tilde{M}_{t},F}^{\vee}$, avec $\Lambda_{1}(\lambda)\in i{\cal A}_{\tilde{G}}^*$, c'est 
$$mes({\cal A}_{\tilde{G}}(F)_{c})^{-1}mes(i{\cal A}_{\tilde{M}_{t},F}^*)mes(i{\cal A}_{\tilde{G},F}^*)^{-1}e^{<\Lambda_{1}(\lambda),X>}\phi_{reg,\tilde{M}_{t}}^{\tilde{G}}(t,\nu;\lambda+\xi,\Lambda_{1}(\lambda)).$$
La somme en $X$ devient simplement
$$\sum_{X\in {\cal A}_{\tilde{G},F}/{\cal A}_{A_{\tilde{G}},F}}e^{<\Lambda_{1}(\lambda),X>}.$$
Cette somme vaut $ [{\cal A}_{\tilde{G},F}:{\cal A}_{A_{\tilde{G}},F}]=mes(i{\cal A}_{\tilde{G},F}^*)mes(A_{\tilde{G}}(F)_{c})$ si $\Lambda_{1}(\lambda)\in i{\cal A}_{\tilde{G},F}^{\vee}$, $0$ sinon. La condition $\Lambda_{1}(\lambda)\in i{\cal A}_{\tilde{G},F}^{\vee}$ \'equivaut \`a $\Lambda(\lambda)=0$. Ces calculs conduisent \`a l'\'egalit\'e de l'\'enonc\'e. $\square$

\bigskip
\subsection{Le terme constant de $j^T$}
 Consid\'erons l'ensemble des $\lambda\in i{\cal A}_{M_{disc}}^*$ tels que $t(\sigma_{\lambda}\circ\theta)\simeq \omega \sigma_{\lambda}$. Il est invariant par translations par $i{\cal A}_{M_{disc},F}^{\vee}+i{\cal A}_{\tilde{M}_{t}}^*$. On note $[\sigma]_{t}$ son quotient par l'action de ce groupe. Ce quotient est fini. On v\'erifie que l'application $\lambda\mapsto (t\lambda'-\lambda,\lambda)$ est une bijection de $[\sigma]_{t}$ sur l'ensemble des couples $(\nu,\lambda)$, o\`u
 
 - $\nu\in [t\sigma',\omega\sigma]$;
 
 - $\lambda\in \{\nu\}_{t}$;
 
 - $\Lambda(\lambda)=0$.
  
  Posons
 $$j_{spec}= \sum_{t\in W^G(M_{disc}\vert M'_{disc})}\vert det((1-t\theta^{-1})_{\vert {\cal A}_{M_{disc}}^{\tilde{M}_{t},*}})\vert ^{-1}  \sum_{\lambda\in [\sigma]_{t}}$$
$$ \int_{i{\cal A}_{\tilde{M}_{t},F}^*}B(\lambda+\xi)\epsilon(t,t\lambda'-\lambda;\lambda+\xi)\phi_{reg,\tilde{M}_{t}}^{\tilde{G}}(t,t\lambda'-\lambda;\lambda+\xi,0)\,d\xi.$$
Ce terme est bien d\'efini d'apr\`es les propri\'et\'es ci-dessus.

\ass{Corollaire}{Il existe une unique fonction $T\mapsto f(T)$ qui appartient \`a $PolExp$ et qui v\'erifie la majoration
$$\vert j^T-f(T)\vert <<\vert T\vert ^{-r}$$
pour tout r\'eel $r$ et tout $T$ dans le c\^one o\`u $j^T$ est d\'efinie. Si $F$ est archim\'edien, on a $c_{0}(f)=j_{spec}$. Si $F$ est non-archim\'edien, pour tout r\'eseau ${\cal R}\subset {\cal A}_{M_{0},F}\otimes_{{\mathbb Z}}{\mathbb Q}$, on a l'\'egalit\'e
$$lim_{k\to \infty}c_{\frac{1}{k}{\cal R},0}(f)=j_{spec}.$$}

Preuve. L'existence de la fonction $f$ r\'esulte du lemme 3.19, de la d\'efinition de ${\bf E}^T$ et du lemme pr\'ec\'edent. L'unicit\'e est claire: un \'el\'ement de $PolExp$ est nul s'il est \`a d\'ecroissance rapide dans un c\^one. On calcule le terme constant en utilisant le lemme pr\'ec\'edent. Ce calcul conduit \`a une formule similaire \`a $j_{spec}$ ci-dessus, \`a ceci pr\`es que la somme en $\lambda\in [\sigma]_{t}$ y est remplac\'ee par une double somme sur $\nu\in [\sigma,\omega t\theta\sigma]$ et $\lambda\in \{\nu\}_{t}$ tel que $\Lambda(\lambda)=0$.  Comme on l'a dit ci-dessus, cette double somme co\"{\i}ncide avec une somme en $\lambda\in [\sigma]_{t}$. $\square$
  
 \bigskip
 
 \subsection{D\'efinition d'une expression spectrale}
 Soit $\tilde{M}$ un espace de Levi de $\tilde{G}$. On suppose que $M_{0}\subset M$.   Soit $\tau=(M_{disc},\sigma,r)$ un triplet form\'e d'un Levi semi-standard $M_{disc}\subset M$, d'une repr\'esentation $\sigma$ de $M_{disc}(F)$ irr\'eductible et de la s\'erie discr\`ete et d'un \'el\'ement $\tilde{r}\in R^{\tilde{M}}(\sigma)$. On fixe un rel\`evement $\boldsymbol{\tilde{r}}\in {\cal R}^{\tilde{M}}(\sigma)$ et on pose $\boldsymbol{\tau}=(M_{disc},\sigma,\boldsymbol{\tilde{r}})$. Fixons aussi un \'el\'ement $\tilde{P}\in {\cal P}(\tilde{M})$, puis un \'el\'ement $S\in {\cal P}(M_{disc})$ tel que $S\subset P$. Posons $\pi_{\tau}=Ind_{S\cap M}^M(\sigma)$ et $\Pi_{\tau}=Ind_{S}^G(\sigma)\simeq Ind_{P}^G(\pi_{\tau})$. A l'aide de $ \boldsymbol{\tau}$, on a d\'efini en  2.9 une repr\'esentation  de $\tilde{M}(F)$ dans l'espace $V_{\sigma,S\cap P}$ de $\pi_{\tau}$ et une repr\'esentation   de $\tilde{G}(F)$ dans l'espace $V_{\sigma,P}$ de $\Pi_{\tau}$. Notons-les respectivement $\tilde{\pi}_{\boldsymbol{\tau}}$ et $\tilde{\Pi}_{\boldsymbol{\tau}}$. On a $\tilde{\Pi}_{\boldsymbol{\tau}}=Ind_{\tilde{P}}^{\tilde{G}}(\tilde{\pi}_{\boldsymbol{\tau}})$. Soit $\tilde{\lambda}\in i\tilde{{\cal A}}_{\tilde{M},F}^*$ en position g\'en\'erale. On peut remplacer $\tau$ par $\tau_{\lambda}=(M,\sigma_{\lambda},\tilde{r})$ et $\boldsymbol{\tau}$ par $\boldsymbol{\tau}_{\tilde{\lambda}}=(M,\sigma_{\lambda},\boldsymbol{\tilde{r}})$ cf. 2.9. On d\'efinit comme en 2.7 la $(\tilde{G},\tilde{M})$-famille $({\cal M}(\pi_{\tau_{\lambda}};\Lambda,\tilde{Q}))_{\tilde{Q}\in {\cal P}(\tilde{M})}$ dont les fonctions prennent leurs valeurs dans l'espace des endomorphismes de $V_{\sigma,S\cap P}$. Rappelons que l'on a fix\'e deux fonctions $f_{1},f_{2}\in C_{c}^{\infty}(\tilde{G}(F))$.  On d\'efinit une $(\tilde{G},\tilde{M})$-famille $({\cal J}(\pi_{\tau_{\lambda}},f_{1},f_{2};\Lambda,\tilde{Q}))_{\tilde{Q}\in {\cal P}(\tilde{M})}$ par
 $${\cal J}(\pi_{\tau_{\lambda}},f_{1},f_{2};\Lambda,\tilde{Q})=\overline{trace({\cal M}(\pi_{\tau_{\lambda}};\Lambda,\tilde{\bar{Q}}))\tilde{\Pi}_{\boldsymbol{\tau}_{\tilde{\lambda}}}(f_{1})}trace({\cal M}(\pi_{\tau_{\lambda}};\Lambda,\tilde{Q})\tilde{\Pi}_{\boldsymbol{\tau}_{\tilde{\lambda}}}(f_{2})).$$
 On en d\'eduit une fonction ${\cal J}_{\tilde{M}}^{\tilde{G}}(\pi_{\tau_{\lambda}},f_{1},f_{2};\Lambda)$. On pose 
 $$J_{\tilde{M}}^{\tilde{G}}(\pi_{\tau_{\lambda}},f_{1},f_{2})={\cal J}_{\tilde{M}}^{\tilde{G}}(\pi_{\tau_{\lambda}},f_{1},f_{2};0).$$
   {\bf Remarque.}  A cause de la double apparition de $\tilde{\pi}_{\boldsymbol{\tau}_{\tilde{\lambda}}}$ dans les d\'efinitions ci-dessus, les fonctions ainsi d\'efinies ne d\'ependent ni  du rel\`evement   $\boldsymbol{\tau}$ de $\tau$, ni du rel\`evement $\tilde{\lambda}$ de $\lambda$.
   \bigskip
   
   On rappelle que l'on a d\'efini la notion de triplet essentiel, cf. 2.9. Cette notion d\'epend de l'espace ambiant. Ici, c'est l'espace $\tilde{M}$. On a
   
   (1) si le triplet $\tau$ n'est pas essentiel pour $\tilde{M}$, les fonctions ci-dessus sont nulles.
   
  Le caract\`ere de $\tilde{\pi}_{\boldsymbol{\tau}_{\tilde{\lambda}}}$ \'etant nul pour tout $\tilde{\lambda}$, la preuve est la m\^eme que celle de  2.7(3).
   
   On suppose d\'esormais $\tau$ essentiel, c'est-\`a-dire $\tau\in E(\tilde{M},\omega)$. Il est clair  que la fonction $\lambda\mapsto J_{\tilde{M}}^{\tilde{G}}(\pi_{\tau_{\lambda}},f_{1},f_{2})$ est la restriction \`a $i{\cal A}_{\tilde{M},F}^*$ d'une fonction m\'eromorphe.

 \ass{Lemme}{La fonction $\lambda\mapsto J_{\tilde{M}}^{\tilde{G}}(\pi_{\tau_{\lambda}},f_{1},f_{2})$ est r\'eguli\`ere  sur tout  $i {\cal A}_{\tilde{M},F}^*$. Si $F$ est archim\'edien, c'est une fonction de Schwartz sur cet espace.}

 Preuve. Convertissons tous les op\'erateurs d'entrelacement qui interviennent dans les d\'efinitions en produits d'op\'erateurs nomalis\'es et de facteurs de normalisation. On obtient une expression
$$ {\cal J}(\pi_{\tau_{\lambda}},f_{1},f_{2};\Lambda,\tilde{Q})={\bf r}(\pi_{\tau_{\lambda}};\Lambda,\tilde{Q}){\cal J}_{reg}(\pi_{\tau_{\lambda}},f_{1},f_{2};\Lambda,\tilde{Q}),$$
o\`u on a regroup\'e dans  ${\bf r}(\pi_{\tau_{\lambda}};\Lambda,\tilde{Q})$ les facteurs de normalisation.  On calcule ${\bf r}(\pi_{\tau_{\lambda}};\Lambda,\tilde{Q})$. C'est le produit de
$$r_{P\vert Q}(\sigma_{\lambda})r_{Q\vert P}(\sigma_{\lambda+\Lambda})r_{P\vert Q}(\sigma_{\lambda+\Lambda/2})^{-1}r_{Q\vert P}(\sigma_{\lambda+\Lambda/2})^{-1}$$
et du conjugu\'e de
$$r_{P\vert \bar{Q}}(\sigma_{\lambda})r_{\bar{Q}\vert P}(\sigma_{\lambda+\Lambda})r_{P\vert \bar{Q}}(\sigma_{\lambda+\Lambda/2})^{-1}r_{\bar{Q}\vert P}(\sigma_{\lambda+\Lambda/2})^{-1}.$$
On obtient
$${\bf r}(\pi_{\tau_{\lambda}};\Lambda,\tilde{Q})=r_{\bar{Q}\vert Q}(\sigma_{\lambda})r_{Q\vert\bar{ Q}}(\sigma_{\lambda+\Lambda})r_{Q\vert \bar{Q}}(\sigma_{\lambda+\Lambda/2})^{-1}r_{\bar{Q}\vert Q}(\sigma_{\lambda+\Lambda/2})^{-1}.$$
Le produit des deux derniers termes est ind\'ependant de $\tilde{Q}$. On peut donc \'ecrire
$${\bf r}(\pi_{\tau_{\lambda}};\Lambda,\tilde{Q})=C(\lambda,\Lambda) {\bf r}_{reg}(\pi_{\tau_{\lambda}};\Lambda,\tilde{Q}),$$
o\`u
$$C(\lambda,\Lambda)=r_{\bar{P}\vert P}(\sigma_{\lambda})r_{P\vert\bar{ P}}(\sigma_{\lambda+\Lambda})r_{P\vert \bar{P}}(\sigma_{\lambda+\Lambda/2})^{-1}r_{\bar{P}\vert P}(\sigma_{\lambda+\Lambda/2})^{-1}$$
et
$${\bf r}_{reg}(\pi_{\tau_{\lambda}};\Lambda,\tilde{Q})=r_{\bar{Q}\vert Q}(\sigma_{\lambda})r_{\bar{P}\vert P}(\sigma_{\lambda})^{-1}r_{Q\vert \bar{Q}}(\sigma_{\lambda+\Lambda})r_{P\vert \bar{P}}(\sigma_{\lambda+\Lambda})^{-1}.$$
En un point $\lambda$ g\'en\'eral, $J_{\tilde{M}}^{\tilde{G}}(\rho_{\tau_{\lambda}},f_{1},f_{2})$ est donc le produit de $C(\lambda,0)$ et d'un terme issu de la $(\tilde{G},\tilde{M})$-famille
$$({\bf r}_{reg}(\rho_{\tau_{\lambda}};\Lambda,\tilde{Q}){\cal J}_{reg}(\rho_{\tau_{\lambda}},f_{1},f_{2};\Lambda,\tilde{Q}))_{\tilde{Q}\in {\cal P}(\tilde{M})}.$$
Le terme $C(\lambda,0)$ vaut $1$. La $(\tilde{G},\tilde{M})$-famille ci-dessus est form\'ee de fonctions r\'eguli\`eres en $\lambda$ et $\Lambda$ (d'apr\`es 1.10(5)). Elle donne donc naissance \`a une fonction r\'eguli\`ere en $\lambda$. Dans le cas o\`u $F$ est archim\'edien, on d\'ecompose  plus finement cette derni\`ere $(\tilde{G},\tilde{M})$-famille en combinaison lin\'eaire de produit de familles similaires et de coefficients matriciels des op\'erateurs $\tilde{\Pi}_{\boldsymbol{\tau}_{\tilde{\lambda}}}(f_{1})$ et $\tilde{\Pi}_{\boldsymbol{\tau}_{\tilde{\lambda}}}(f_{2})$.  Ces coefficients sont de Schwartz tandis les autres termes donnent naissance \`a des fonctions \`a croissance mod\'er\'ee en $\lambda$ (lemme 1.4). D'o\`u le lemme. $\square$

{\bf Remarque.} Dans le cas o\`u $\tilde{M}=\tilde{G}$, on a simplement l'\'egalit\'e
$$J_{\tilde{G}}^{\tilde{G}}(\pi_{\tau_{\lambda}},f_{1},f_{2})=\overline{trace(\tilde{\Pi}_{\boldsymbol{\tau}_{\tilde{\lambda}}}(f_{1}))}trace(\tilde{\Pi}_{\boldsymbol{\tau}_{\tilde{\lambda}}}(f_{2})).$$
\bigskip
 
 D'apr\`es le lemme, le terme $J_{\tilde{M}}^{\tilde{G}}(\pi_{\tau_{\lambda}},f_{1},f_{2})$ est   d\'efini pour tout $\lambda$. On v\'erifie qu'il ne d\'epend pas de l'espace parabolique $\tilde{P}$ choisi et qu'il ne d\'epend que de la classe de conjugaison par $M(F)$ de l'\'el\'ement $\tau_{\lambda}$. On pose
 $$J_{\tilde{M},spec}^{\tilde{G}}(\omega,f_{1},f_{2})=\sum_{\tau\in (E_{disc}(\tilde{M},\omega)/conj)/i{\cal A}_{\tilde{M},F}^*}\vert {\bf Stab}(W^M\times i{\cal A}_{\tilde{M},F}^*,\tau)\vert ^{-1} \iota(\tau)$$
 $$\int_{i{\cal A}_{\tilde{M},F}^*}J_{\tilde{M}}^{\tilde{G}}(\pi_{\tau_{\lambda}},f_{1},f_{2})\,d\lambda.$$
 Pour toute classe $\tau\in (E_{disc}(\tilde{M},\omega)/conj)/i{\cal A}_{\tilde{M},F}^*$, on a choisi un point-base que l'on a \'egalement not\'e $\tau$. Les termes ${\bf Stab}(W^M\times i{\cal A}_{\tilde{M},F}^*,\tau)$ et $\iota(\tau)$ ont \'et\'e d\'efinis en 2.9 et 2.10. On v\'erifie que l'expression  $J_{\tilde{M},spec}^{\tilde{G}}(\omega,f_{1},f_{2})$ ne d\'epend que de la classe de conjugaison par $\tilde{G}(F)$ de l'espace de Levi $\tilde{M}$.
  
\bigskip

\subsection{La formule spectrale} 
 On pose $\tilde{W}^G=Norm_{G(F)}(\tilde{M}_{0})/M_{0}(F)$ (on ne confondra pas ce groupe avec le quotient par $M_{0}(F)$ du normalisateur de $M_{0}$ dans $\tilde{G}(F)$). Posons
$$J^{\tilde{G}}_{spec}(\omega,f_{1},f_{2})= \sum_{\tilde{M}\in {\cal L}(\tilde{M}_{0})}\vert \tilde{W}^M\vert \vert \tilde{W}^G\vert ^{-1}(-1)^{a_{\tilde{M}}-a_{\tilde{G}}}J_{\tilde{M},spec}^{\tilde{G}}(\omega,f_{1},f_{2}).$$

\ass{Proposition}{Il existe une unique fonction $T\mapsto \varphi(T)$ qui appartient \`a $PolExp$ et qui v\'erifie pour tout r\'eel $r$ la majoration
$$\vert J^T(\omega,f_{1},f_{2})-\varphi(T)\vert <<\vert T\vert ^{-r}$$
pour tout $T$ dans le c\^one o\`u $J^T(\omega,f_{1},f_{2})$ est d\'efinie. Si $F$ est archim\'edien, on a l'\'egalit\'e
$$c_{0}(\varphi)=J^{\tilde{G}}_{spec}(\omega,f_{1},f_{2}).$$
Si $F$ est non-archim\'edien, pour tout r\'eseau ${\cal R}\subset {\cal A}_{M_{0},F}\otimes_{{\mathbb Z}}{\mathbb Q}$, on a l'\'egalit\'e
$$lim_{k\to \infty}c_{\frac{1}{k}{\cal R},0}(\varphi)=J^{\tilde{G}}_{spec}(\omega,f_{1},f_{2}).$$}

Preuve. Il r\'esulte des formules 3.2(1) et 3.3(1) que $J^T(\omega,f_{1},f_{2})$ est somme finie de termes $j^T$ tels qu'en 3.4. L'existence de la fonction $\varphi$ r\'esulte donc du corollaire 3.24. Comme toujours, l'unicit\'e de $\varphi$ est \'evidente. Le terme constant de $\varphi$ se calcule en sommant les termes constants des diff\'erents $j^T$ intervenant, lesquels sont calcul\'es par le m\^eme corollaire 3.24. Commen\c{c}ons par calculer le terme constant de l'\'el\'ement de $PolExp$ asymptote \`a l'expression $J^T_{M_{disc},\sigma}(\omega,f_{1},f_{2})$ de 3.3(1). Il est \'egal \`a
$$(1) \qquad  \sum_{u,v\in {\cal B}, u',v'\in {\cal B}'} \sum_{t\in W^G(M_{disc}\vert M'_{disc})}\vert det((1-t\theta^{-1})_{\vert {\cal A}_{M_{disc}}^{\tilde{M}_{t},*}})\vert ^{-1}  \sum_{\lambda\in [\sigma]_{t}}$$
$$ \int_{i{\cal A}_{\tilde{M}_{t},F}^*}B_{u,v,u',v'}(\lambda+\xi)\epsilon(t,t\lambda'-\lambda;\lambda+\xi)\phi_{reg,\tilde{M}_{t}}^{\tilde{G}}(u,v,u',v';t,t\lambda'-\lambda;\lambda+\xi,0)\,d\xi.$$
On a pr\'ecis\'e la notation en r\'etablissant le quadruplet $(u,v,u',v')$ dans la fonction not\'ee simplement $\phi_{reg,\tilde{M}_{t}}^{\tilde{G}}(t,t\lambda'-\lambda;\lambda+\xi,0)$ dans le corollaire 3.24. Les sommes en $u$, $v$, $u'$ et $v'$ sont finies, on peut les faire entrer sous l'int\'egrale.  Fixons $t$, $\lambda$ et $\xi$, calculons
$$\sum_{u,v\in {\cal B}, u',v'\in {\cal B}'}B_{u,v,u',v'}(\lambda+\xi)\phi_{reg,\tilde{M}_{t}}^{\tilde{G}}(u,v,u',v';t,t\lambda'-\lambda;\lambda+\xi,0).$$
C'est la valeur en $\Lambda=0$ de la fonction $h_{\tilde{M}_{t}}^{\tilde{G}}(\Lambda)$ associ\'ee \`a la $(\tilde{G},\tilde{M}_{t})$-famille $(h(\Lambda,\tilde{P}))_{\tilde{P}\in {\cal P}(\tilde{M}_{t})}$ d\'efinie par
$$h(\Lambda,\tilde{P})=\sum_{u,v\in {\cal B}, u',v'\in {\cal B}'}B_{u,v,u',v'}(\lambda+\xi)\phi_{reg}(u,v,u',v';t,t\lambda'-\lambda;\lambda+\xi,\Lambda,\tilde{P}).$$
Fixons $S''\in {\cal P}(M_{disc})$ tel que $S''\subset P$. Posons pour simplifier $\nu=t\lambda'-\lambda$, $\underline{\sigma}=\sigma_{\lambda}$, $\underline{\pi}=\pi_{\lambda}$. En revenant aux d\'efinitions de 3.3 et 3.16, on a
$$h(\Lambda,\tilde{P})=\sum_{u,v\in {\cal B}, u',v'\in {\cal B}'}(\underline{\pi}_{\xi}(\varphi_{1})U_{\theta,\underline{\sigma}_{\xi}}u',v)(u,\underline{\pi}_{\xi}(\varphi_{2})U_{\theta,\underline{\sigma}_{\xi}}v') r_{S'',reg}(\underline{\sigma}_{\xi})^{-1}r_{S'',reg}(\underline{\sigma}_{\xi+\Lambda})$$
$$(A(t, \nu;\lambda+\xi)v',R_{\bar{S}''\vert S}(\underline{\sigma}_{\xi})^{-1}\circ R_{\bar{S}''\vert S}(\underline{\sigma}_{\xi+\Lambda})u)$$
$$(R_{S''\vert S}(\underline{\sigma}_{\xi})^{-1}\circ R_{S''\vert S}(\underline{\sigma}_{\xi+\Lambda})v,A(t,\nu;\lambda+\xi)u').$$
Les sommes en $u$ et $v$ se simplifient en
$$(2) \qquad h(\Lambda,\tilde{P})=\sum_{u',v'\in {\cal B}'} r_{S'',reg}(\underline{\sigma}_{\xi})^{-1}r_{S'',reg}(\underline{\sigma}_{\xi+\Lambda})$$
$$(A(t, \nu;\lambda+\xi)v',R_{\bar{S}''\vert S}(\underline{\sigma}_{\xi})^{-1}\circ R_{\bar{S}''\vert S}(\underline{\sigma}_{\xi+\Lambda})\underline{\pi}_{\xi}(\varphi_{2})U_{\theta,\underline{\sigma}_{\xi}}v')$$
$$(R_{S''\vert S}(\underline{\sigma}_{\xi})^{-1}\circ R_{S''\vert S}(\underline{\sigma}_{\xi+\Lambda})\underline{\pi}_{\xi}(\varphi_{1})U_{\theta,\underline{\sigma}_{\xi}}u',A(t,\nu;\lambda+\xi)u').$$
On a
$$\sum_{v'\in {\cal B}'}(A(t, \nu;\lambda+\xi)v',R_{\bar{S}''\vert S}(\underline{\sigma}_{\xi})^{-1}\circ R_{\bar{S}''\vert S}(\underline{\sigma}_{\xi+\Lambda})\underline{\pi}_{\xi}(\varphi_{2})U_{\theta,\underline{\sigma}_{\xi}}v')$$
$$=\sum_{v'\in {\cal B}'}( v',A(t, \nu;\lambda+\xi)^{-1}R_{\bar{S}''\vert S}(\underline{\sigma}_{\xi})^{-1}\circ R_{\bar{S}''\vert S}(\underline{\sigma}_{\xi+\Lambda})\underline{\pi}_{\xi}(\varphi_{2})U_{\theta,\underline{\sigma}_{\xi}}v')$$
$$=trace(A(t, \nu;\lambda+\xi)^{-1}R_{\bar{S}''\vert S}(\underline{\sigma}_{\xi})^{-1}\circ R_{\bar{S}''\vert S}(\underline{\sigma}_{\xi+\Lambda})\underline{\pi}_{\xi}(\varphi_{2})U_{\theta,\underline{\sigma}_{\xi}})$$
$$(3) \qquad =trace( R_{\bar{S}''\vert S}(\underline{\sigma}_{\xi})^{-1}\circ R_{\bar{S}''\vert S}(\underline{\sigma}_{\xi+\Lambda})\underline{\pi}_{\xi}(\varphi_{2})U_{\theta,\underline{\sigma}_{\xi}}A(t, \nu;\lambda+\xi)^{-1}).$$
La condition $\lambda\in [\sigma]_{t}$ \'equivaut \`a l'\'equivalence $t\underline{\sigma}_{\xi}'\simeq\omega \underline{\sigma}_{\xi}$, ou encore $\underline{\sigma}_{\xi}\circ ad_{\gamma}\simeq \omega\underline{\sigma}_{\xi}$, o\`u $\gamma=\gamma_{0}t^{-1}$ (en identifiant $t$ \`a un rel\`evement dans $G(F)$). On v\'erifie que le couple $(A_{\nu},\gamma)$ appartient \`a ${\cal N}^{\tilde{M_{t}}}(\underline{\sigma}_{\xi})\subset {\cal N}^{\tilde{G}}(\underline{\sigma}_{\xi})$. En  consid\'erant $(A_{\nu},\gamma)$ comme un \'el\'ement de ce dernier ensemble, on a d\'efini l'op\'erateur $\tilde{\nabla}_{S}(A_{\nu},\gamma)$ en 2.8. En comparant les d\'efinitions, on obtient l'\'egalit\'e
$$U_{\theta,\underline{\sigma}_{\xi}}A(t, \nu;\lambda+\xi)^{-1}=\tilde{\nabla}_{S}(A_{\nu},\gamma)\underline{\pi}_{\xi}(t).$$
Relevons $\xi$ en l'\'el\'ement $\tilde{\xi}\in i\tilde{{\cal A}}_{\tilde{M}_{t},F}^*$ tel que $<\tilde{\xi},\tilde{H}_{\tilde{M}_{t}}(\gamma)>=0$. A l'aide de cet \'el\'ement, on identifie ${\cal R}^{\tilde{M}_{t}}(\underline{\sigma})$ \`a ${\cal R}^{\tilde{M}_{t}}(\underline{\sigma}_{\xi})$. 
Notons $\boldsymbol{\tilde{r}}$ l'image de $(A_{\nu},\gamma)$ dans ${\cal R}^{\tilde{M_{t}}}(\underline{\sigma}_{\xi})={\cal R}^{\tilde{M_{t}}}(\underline{\sigma})$, posons $\boldsymbol{\tau}=(M_{disc},\underline{\sigma},\boldsymbol{\tilde{r}})$ et $\tau=(M_{disc},\underline{\sigma},\tilde{r})$, o\`u $\tilde{r}$ est l'image de $\boldsymbol{\tilde{r}}$ dans $R^{\tilde{M}_{t}}(\underline{\sigma})$.   Introduisons la repr\'esentation de $\tilde{G}(F)$ associ\'ee en 2.8 \`a $\boldsymbol{\tau}_{\tilde{\xi}}$, vu comme un triplet pour $\tilde{G}$, que l'on note ici simplement $\underline{\tilde{\Pi}}_{\xi}$. Alors 
$$\tilde{\nabla}_{S}(A_{\nu},\gamma)\underline{\pi}_{\xi}(t)=\omega(t)^{-1}\underline{\tilde{\Pi}}_{\xi}(\gamma_{0}).$$
En revenant \`a la d\'efinition de $\varphi_{2}$ en 3.2, on voit que
$$\underline{\pi}_{\xi}(\varphi_{2})U_{\theta,\underline{\sigma}_{\xi}}A(t, \nu;\lambda+\xi)^{-1}=\omega(t)^{-1}\underline{\tilde{\Pi}}_{\xi}(f_{2}).$$
L'expression (3) devient simplement
$$\omega(t)^{-1}trace( R_{\bar{S}''\vert S}(\underline{\sigma}_{\xi})^{-1}\circ R_{\bar{S}''\vert S}(\underline{\sigma}_{\xi+\Lambda})\underline{\tilde{\Pi}}_{\xi}(f_{2})).$$
On traite de m\^eme la somme en $u'$. Les op\'erateurs portant cette fois sur la premi\`ere variable des produits scalaires, on obtient une expression conjugu\'ee. En particulier, il sort un facteur $\overline{\omega(t)}^{-1}$, dont le produit avec $\omega(t)^{-1}$ ci-dessus vaut $1$. L'expression (2) devient
$$(4) \qquad h(\Lambda,\tilde{P})= r_{S'',reg}(\underline{\sigma}_{\xi})^{-1}r_{S'',reg}(\underline{\sigma}_{\xi+\Lambda})\overline{trace( R_{S''\vert S}(\underline{\sigma}_{\xi})^{-1}\circ R_{S''\vert S}(\underline{\sigma}_{\xi+\Lambda})\underline{\tilde{\Pi}}_{\xi}(f_{1}))}$$
$$trace( R_{\bar{S}''\vert S}(\underline{\sigma}_{\xi})^{-1}\circ R_{\bar{S}''\vert S}(\underline{\sigma}_{\xi+\Lambda})\underline{\tilde{\Pi}}_{\xi}(f_{2})).$$
Tous les termes sont relatifs \`a un parabolique $S$ de r\'ef\'erence, que l'on a fix\'e au d\'ebut du calcul. Fixons un \'el\'ement $\tilde{P}_{1}\in {\cal P}(\tilde{M}_{t})$ et un \'el\'ement $S_{1}\in {\cal P}(M_{disc})$ tel que $S_{1}\subset P_{1}$. Rempla\c{c}ons $S$ par $S_{1}$ dans les d\'efinitions  des termes ci-dessus. On obtient une nouvelle expression $h_{1}(\Lambda,\tilde{P})$. On a

(5) les valeurs en $\Lambda=0$ des fonctions $h_{\tilde{M}}^{\tilde{G}}(\Lambda)$ et $h_{1,\tilde{M}}^{\tilde{G}}(\Lambda)$ sont \'egales.

On a
$$trace( R_{\bar{S}''\vert S}(\underline{\sigma}_{\xi})^{-1}\circ R_{\bar{S}''\vert S}(\underline{\sigma}_{\xi+\Lambda})\underline{\tilde{\Pi}}_{\xi}(f_{2}))=$$
$$trace( R_{S_{1}\vert S}(\underline{\sigma}_{\xi})^{-1}R_{\bar{S}''\vert S_{1}}(\underline{\sigma}_{\xi})^{-1}\circ R_{\bar{S}''\vert S_{1}}(\underline{\sigma}_{\xi+\Lambda})R_{S_{1}\vert S}(\underline{\sigma}_{\xi+\Lambda})\underline{\tilde{\Pi}}_{\xi}(f_{2}))$$
$$=trace(  R_{\bar{S}''\vert S_{1}}(\underline{\sigma}_{\xi})^{-1}\circ R_{\bar{S}''\vert S_{1}}(\underline{\sigma}_{\xi+\Lambda})R_{S_{1}\vert S}(\underline{\sigma}_{\xi+\Lambda})\underline{\tilde{\Pi}}_{\xi}(f_{2})R_{S_{1}\vert S}(\underline{\sigma}_{\xi})^{-1}).$$
Le produit final 
$$R_{S_{1}\vert S}(\underline{\sigma}_{\xi+\Lambda})\underline{\tilde{\Pi}}_{\xi}(f_{2})R_{S_{1}\vert S}(\underline{\sigma}_{\xi})^{-1}$$
ne d\'epend pas de $S''$. Un principe g\'en\'eral valable pour toute $(\tilde{G},\tilde{M})$-famille dit que l'on ne modifie pas la valeur $h_{\tilde{M}}^{\tilde{G}}(0)$ si on remplace un facteur ind\'ependant de $S''$ par sa valeur en $\Lambda=0$. Mais, en $\Lambda=0$, le produit ci-dessus est \'egal \`a $ \underline{\tilde{\Pi}}_{1,\xi}(f_{2})$, o\`u $ \underline{\tilde{\Pi}}_{1,\xi}$ est l'analogue de $ \underline{\tilde{\Pi}}_{\xi}$ relatif au parabolique $S_{1}$. Un m\^eme raisonnement s'applique aux autres termes et (5) s'ensuit.

En oubliant ce changement de $S$ en $S_{1}$, on suppose pour simplifier que $S\subset P_{1}$. Pour calculer $h(\Lambda,\tilde{P})$, on a choisi $S''$ avec $S''\subset P$. On peut lui imposer de plus $S''\cap M_{t}=S\cap M_{t}$. Introduisons l'\'el\'ement $\underline{S}''\in {\cal P}(M_{disc})$ tel que $\underline{S}''\subset \bar{P}$ et $\underline{S}''\cap M_{t}=S\cap M_{t}$. On a
$$R_{\bar{S}''\vert S}(\underline{\sigma}_{\xi})^{-1}\circ R_{\bar{S}''\vert S}(\underline{\sigma}_{\xi+\Lambda})=R_{\underline{S}''\vert S}(\underline{\sigma}_{\xi})^{-1}\circ R_{\bar{S}''\vert \underline{S}''}(\underline{\sigma}_{\xi})^{-1}\circ R_{\bar{S}''\vert \underline{S}''}(\underline{\sigma}_{\xi+\Lambda})\circ R_{\underline{S}''\vert S}(\underline{\sigma}_{\xi+\Lambda}).$$
Mais $\bar{S}''$ et $\underline{S}''$ sont tous deux contenus dans $\bar{P}$. L'op\'erateur $R_{\bar{S}''\vert \underline{S}''}(\underline{\sigma}_{\xi+\Lambda})$ est induit de l'op\'erateur $R_{\bar{S}''\cap M_{t}\vert \underline{S}''\cap M_{t}}^{M_{t}}(\underline{\sigma}_{\xi+\Lambda})$ et celui-ci ne d\'epend pas de $\Lambda$ puisque $\Lambda\in i{\cal A}_{\tilde{M}_{t},F}^*$. D'o\`u
$$R_{\bar{S}''\vert \underline{S}''}(\underline{\sigma}_{\xi})^{-1}\circ R_{\bar{S}''\vert \underline{S}''}(\underline{\sigma}_{\xi+\Lambda})=1,$$
puis
$$R_{\bar{S}''\vert S}(\underline{\sigma}_{\xi})^{-1}\circ R_{\bar{S}''\vert S}(\underline{\sigma}_{\xi+\Lambda})=R_{\underline{S}''\vert S}(\underline{\sigma}_{\xi})^{-1}\circ R_{\underline{S}''\vert S}(\underline{\sigma}_{\xi+\Lambda}).$$
R\'etablissons maintenant les op\'erateurs d'entrelacement non normalis\'es. La formule (4) se transforme en
$$(6)\qquad h(\Lambda,\tilde{P})={\bf r}_{S''}(\underline{\sigma}_{\xi})^{-1}{\bf r}_{S''}(\underline{\sigma}_{\xi+\Lambda})\overline{trace(J_{S''\vert S}(\underline{\sigma}_{\xi})^{-1}\circ J_{S''\vert S}(\underline{\sigma}_{\xi+\Lambda})\underline{\tilde{\Pi}}_{\xi}(f_{1}))}$$
$$trace(J_{\underline{S}''\vert S}(\underline{\sigma}_{\xi})^{-1}\circ J_{\underline{S}''\vert S}(\underline{\sigma}_{\xi+\Lambda})\underline{\tilde{\Pi}}_{\xi}(f_{2})),$$
o\`u
$${\bf r}_{S''}(\underline{\sigma}_{\xi})=r_{S'',reg}(\underline{\sigma}_{\xi})r_{\underline{S}''\vert S}(\underline{\sigma}_{\xi})^{-1}r_{S\vert S''}(\underline{\sigma}_{\xi})^{-1}=r_{S'',reg}(\underline{\sigma}_{\xi})r_{\underline{S}''\vert S''}(\underline{\sigma}_{\xi})^{-1},$$
puisque les distances entre $\underline{S}''$ et $S$ et entre $S$ et $S''$ s'ajoutent.  Si l'on suppose $\xi$ en position g\'en\'erale, tous les termes ci-dessus sont bien d\'efinis. Introduisons la repr\'esentation $\tilde{\pi}_{\boldsymbol{\tau}_{\tilde{\xi}}}$ de $\tilde{M}_{t}(F)$ associ\'ee \`a $\boldsymbol{\tau}_{\tilde{\xi}}$. On note $\pi_{\tau_{\xi}}$ la repr\'esentation sous-jacente de $M_{t}(F)$. On a simplement $\underline{\tilde{\Pi}}_{\xi}=\tilde{\Pi}_{\boldsymbol{\tau}_{\tilde{\xi}}}=Ind_{\tilde{P}_{1}}^{\tilde{G}}(\tilde{\pi}_{\boldsymbol{\tau}_{\tilde{\xi}}})$. En identifiant ces repr\'esentations, les propri\'et\'es d'induction des op\'erateurs d'entrelacement entra\^{\i}nent que $J_{\underline{S}''\vert S}(\underline{\sigma}_{\xi})^{-1}\circ J_{\underline{S}''\vert S}(\underline{\sigma}_{\xi+\Lambda})$ co\"{\i}ncide avec
$J_{\bar{P}\vert P_{1}}(\pi_{\tau_{\xi}})^{-1}\circ J_{\bar{P}\vert P_{1}}(\pi_{\tau_{\xi+\Lambda}})$. En appliquant les d\'efinitions, la formule (6) se r\'ecrit
$$(7) \qquad h(\Lambda,\tilde{P})=C(\underline{\sigma}_{\xi},\Lambda,\tilde{P})\overline{trace({\cal M}(\pi_{\tau_{\xi}};\Lambda,\tilde{P})\underline{\tilde{\Pi}}_{\xi}(f_{1}))}trace({\cal M}(\pi_{\tau_{\xi}};\Lambda,\tilde{\bar{P}})\underline{\tilde{\Pi}}_{\xi}(f_{2})),$$
o\`u
$$C(\underline{\sigma}_{\xi},\Lambda,\tilde{P})={\bf r}_{S''}(\underline{\sigma}_{\xi})^{-1}{\bf r}_{S''}(\underline{\sigma}_{\xi+\Lambda})\mu_{\bar{P}\vert P_{1}}(\pi_{\tau_{\xi}})\mu_{\bar{P}\vert P_{1}}( \pi_{\tau_{\xi+\Lambda/2}})^{-1}\overline{\mu_{P\vert P_{1}}(\pi_{\tau_{\xi}})\mu_{P\vert P_{1}}(\pi_{\tau_{\xi+\Lambda/2}})^{-1}}.$$
 
On a l'\'egalit\'e
$$\mu_{\bar{P}\vert P_{1}}(\pi_{\tau_{\xi}})\overline{\mu_{P\vert P_{1}}(\pi_{\tau_{\xi}})}=r_{\bar{P}\vert P}(\underline{\sigma}_{\xi})^{-1}r_{P\vert \bar{P}}(\underline{\sigma}_{\xi})^{-1}.$$
Ce terme est ind\'ependant de $\tilde{P}$: c'est le produit des $r_{\alpha}(\underline{\sigma}_{\xi})^{-1}$ pour toutes les racines simples $\alpha$ de $A_{M_{disc}}$ intervenant dans $G$ et pas dans $M_{t}$. Pour la m\^eme raison, le produit
 $$\mu_{\bar{P}\vert P_{1}}( \pi_{\tau_{\xi+\Lambda/2}})^{-1}\overline{\mu_{P\vert P_{1}}(\pi_{\tau_{\xi+\Lambda/2}})^{-1}}$$
 est ind\'ependant de $\tilde{P}$. Comme dans la preuve de (5), on ne change pas le terme $h_{\tilde{M}_{t}}^{\tilde{G}}(0)$ en rempla\c{c}ant ces termes par leurs valeurs en $\Lambda=0$. Mais alors le produit de ces deux termes vaut $1$. On peut d\'efinir ${\bf r}_{S''}(\underline{\sigma}_{\xi+\mu})$ pour tout $\mu\in i{\cal A}_{M_{disc},F}^*$. Par d\'efinition  ${\bf r}_{S''}(\underline{\sigma}_{\xi})$ est la valeur de cette fonction en $\mu=0$. Pour $\mu$ en position g\'en\'erale, on a
 $${\bf r}_{S''}(\underline{\sigma}_{\xi+\mu})=r_{\bar{S}''\vert S''}(\underline{\sigma}_{\xi+\mu})r_{\underline{S}''\vert S''}(\underline{\sigma}_{\xi+\mu})^{-1}r_{\bar{S}\vert S}(\underline{\sigma}_{\xi+\mu})^{-1}.$$
 On introduit l'\'el\'ement $\underline{S}\in {\cal P}(M_{disc})$ tel que $\underline{S}\subset \bar{P}_{1}$ et $\underline{S}\cap M_{t}=S\cap M_{t}$. Les distances entre $\bar{S}''$ et $ \underline{S}''$ et entre $\underline{S}''$ et $S''$ s'ajoutent. De m\^eme, les distances entre $\bar{S}$ et $ \underline{S}$ et entre $\underline{S}$ et $S$ s'ajoutent. D'o\`u
 $${\bf r}_{S''}(\underline{\sigma}_{\xi+\mu})=r_{\bar{S}''\vert \underline{S}''}(\underline{\sigma}_{\xi+\mu})r_{\bar{S}\vert \underline{S}}(\underline{\sigma}_{\xi+\mu})^{-1}r_{\underline{S}\vert S}(\underline{\sigma}_{\xi+\mu})^{-1}.$$
 Par les propri\'et\'es habituelles d'induction,
 $$r_{\bar{S}''\vert \underline{S}''}(\underline{\sigma}_{\xi+\mu})=r^{M_{t}}_{(\bar{S}''\cap M_{t})\vert (\underline{S}''\cap M_{t})}(\underline{\sigma}_{\xi+\mu})=r^{M_{t}}_{(\bar{S}\cap M_{t})\vert (\underline{S}\cap M_{t})}(\underline{\sigma}_{\xi+\mu})=r_{\bar{S}\vert \underline{S}}(\underline{\sigma}_{\xi+\mu}).$$
 D'o\`u ${\bf r}_{S''}(\underline{\sigma}_{\xi+\mu})=r_{\underline{S}\vert S}(\underline{\sigma}_{\xi+\mu})^{-1}$. Ce terme est r\'egulier en $\mu=0$ (pour $\xi$ en position g\'en\'erale), d'o\`u  ${\bf r}_{S''}(\underline{\sigma}_{\xi})=r_{\underline{S}\vert S}(\underline{\sigma}_{\xi})^{-1}$. On obtient
 $${\bf r}_{S''}(\underline{\sigma}_{\xi})^{-1}{\bf r}_{S''}(\underline{\sigma}_{\xi+\Lambda})=r_{\underline{S}\vert S}(\underline{\sigma}_{\xi})r_{\underline{S}\vert S}(\underline{\sigma}_{\xi+\Lambda})^{-1}.$$
   Ce terme est ind\'ependant de $S''$. Comme ci-dessus, on ne change pas le terme $h_{\tilde{M}_{t}}^{\tilde{G}}(0)$ en rempla\c{c}ant l'expression pr\'ec\'edente par sa valeur en $\Lambda=0$. Mais celle-ci est $1$. Cela prouve qu'on peut supprimer le terme $C(\underline{\sigma}_{\xi},\Lambda,\tilde{P})$ de la formule (7) sans changer la valeur de $h_{\tilde{M}_{t}}^{\tilde{G}}(0)$. Quand on supprime ce terme $C(\underline{\sigma}_{\xi},\Lambda,\tilde{P})$, le membre de droite de (6) devient ${\cal J}(\pi_{\tau_{\xi}},f_{1},f_{2},\tilde{\bar{P}})$. Si c'\'etait plus simplement ${\cal J}(\pi_{\tau_{\xi}},f_{1},f_{2},\tilde{P})$, on en d\'eduirait $h_{\tilde{M}_{t}}^{\tilde{G}}(0)=J_{\tilde{M}_{t}}^{\tilde{G}}(\pi_{\tau_{\xi}},f_{1},f_{2})$. A cause du changement de $P$ en $\bar{P}$ et parce que $\epsilon_{\tilde{\bar{P}}}^{\tilde{G}}(\Lambda)=(-1)^{a_{\tilde{M}_{t}}-a_{\tilde{G}}}\epsilon_{\tilde{P}}^{\tilde{G}}(\Lambda)$, on obtient
 $$ h_{\tilde{M}_{t}}^{\tilde{G}}(0)=(-1)^{a_{\tilde{M}_{t}}-a_{\tilde{G}}}J_{\tilde{M}_{t}}^{\tilde{G}}(\pi_{\tau_{\xi}},f_{1},f_{2}).$$

 On doit aussi calculer le terme $\epsilon(t,\nu;\lambda+\xi)$ qui intervient dans (1). D'apr\`es sa d\'efinition en 3.16, c'est la valeur en $\mu=\lambda+\xi$ de $r_{\bar{S}\vert S}(\sigma_{\mu})r_{\bar{S}\vert S}(\sigma_{t\mu'-\nu})^{-1}$. Ou encore la valeur en $\mu=0$ de $r_{\bar{S}\vert S}(\underline{\sigma}_{\xi+\mu})r_{\bar{S}\vert S}(\underline{\sigma}_{\xi+t\mu'})^{-1}$. Le calcul est fait par Arthur en [A1] p.87. Le rapport pr\'ec\'edent est \'egal \`a
 $$(9)\qquad \prod_{\alpha>_{S}0}r_{\alpha}(\underline{\sigma}_{\xi+\mu})r_{\alpha}(\underline{\sigma}_{\xi+t\mu'})^{-1},$$
 o\`u $\alpha$ parcourt les racines simples de $A_{M_{disc}}$ dans $G$ qui sont positives pour $S$. Si $\alpha$ n'intervient pas dans $M_{t}$, $r_{\alpha}(\underline{\sigma}_{\xi+\mu})$ est r\'eguli\`ere en $\mu=0$ puisqu'on a suppos\'e $\xi$ en position g\'en\'erale. Le rapport $r_{\alpha}(\underline{\sigma}_{\xi+\mu})r_{\alpha}(\underline{\sigma}_{\xi+t\mu'})^{-1}$ vaut $1$ en $\mu=0$. Supposons que $\alpha$ intervient dans $M_{t}$. On d\'efinit un Levi $M_{\alpha}$ contenant $M_{disc}$ comme en 1.11. Si $m^{M_{\alpha}}(\underline{\sigma}_{\xi})\not=0$, $r_{\alpha}(\underline{\sigma}_{\xi+\mu})$ est encore  r\'eguli\`ere en $\mu=0$ et le  rapport $r_{\alpha}(\underline{\sigma}_{\xi+\mu})r_{\alpha}(\underline{\sigma}_{\xi+t\mu'})^{-1}$ vaut encore $1$ en $\mu=0$.  Supposons $m^{M_{\alpha}}(\underline{\sigma}_{\xi})=0$. On sait qu'alors $r_{\alpha}(\underline{\sigma}_{\xi+\mu})$ a un p\^ole d'ordre $1$ en $\mu=0$. Plus pr\'ecis\'ement,  $r_{\alpha}(\underline{\sigma}_{\xi+\mu})$ est  \'equivalente au voisinage de $\mu=0$ \`a  $c_{\alpha}<\mu,\check{\alpha}>$, o\`u $c_{\alpha}$ est un nombre r\'eel non nul. Le rapport $r_{\alpha}(\underline{\sigma}_{\xi+\mu})r_{\alpha}(\underline{\sigma}_{\xi+t\mu'})^{-1}$ est donc \'equivalent \`a $<\mu,\check{\alpha}><t\mu',\check{\alpha}>^{-1}$, ou encore \`a $<\mu,\check{\alpha}><\mu,\theta t^{-1}\check{\alpha}>^{-1}$. En notant $\Sigma$ l'ensemble des racines $\alpha$ intervenant dans $M_{t}$ pour lesquelles $m^{M_{\alpha}}(\underline{\sigma}_{\xi})=0$, l'expression (1) est donc \'equivalente \`a
 $$\prod_{\alpha\in \Sigma; \alpha>_{S}0}<\mu,\check{\alpha}><\mu,\theta t^{-1}\check{\alpha}>^{-1}.$$
 Ce produit est \'egal au signe $\epsilon_{\underline{\sigma}}(w)$ d\'efini en 2.9, o\`u $w$ est l'image de $\gamma=\gamma_{0}t^{-1}$ dans l'ensemble $W^{\tilde{M_{t}}}(\underline{\sigma}_{\xi})=W^{\tilde{M_{t}}}(\underline{\sigma})$. D'o\`u l'\'egalit\'e
 $$\epsilon(t,\nu;\lambda+\xi)=\epsilon_{\underline{\sigma}}(w).$$

 A ce point, on a transform\'e la formule (1) en 
 $$ \sum_{t\in W^G(M_{disc}\vert M'_{disc})}(-1)^{a_{\tilde{M}_{t}}-a_{\tilde{G}}} \vert det((1-t\theta^{-1})_{\vert {\cal A}_{M_{disc}}^{\tilde{M}_{t}}})\vert ^{-1}$$
 $$\sum_{\lambda\in [\sigma]_{t}}\epsilon_{\sigma_{\lambda}}(w)\int_{i{\cal A}_{\tilde{M}_{t},F}^*}J_{\tilde{M}_{t}}^{\tilde{G}}(\pi_{\tau_{\xi}},f_{1},f_{2})\,d\xi.$$
On rappelle que $w$ est l'image de $\gamma=\gamma_{0}t^{-1}$ dans $W^{\tilde{M_{t}}}(\sigma_{\lambda})$ et que $\tau=(M_{disc},\sigma_{\lambda},r)$, o\`u $r$ est l'image de $ w$ dans $R^{\tilde{M}_{t}}(\sigma_{\lambda})$. On d\'ecompose la formule ci-dessus selon les espaces de Levi $\tilde{M}_{t}$. Elle devient
$$(10)\qquad \sum_{\tilde{M};M_{disc}\subset M}X(\tilde{M},M_{disc},\sigma),$$
o\`u $X(\tilde{M},M_{disc},\sigma)$ est l'expression obtenue \`a partir de l'expression ci-dessus en limitant la somme aux $t$ tels que $\tilde{M}_{t}=\tilde{M}$. Fixons $\tilde{M}$. On peut remplacer la double somme en $t\in W^G(M_{disc},M'_{disc})$ tels que $\tilde{M}_{t}=\tilde{M}$ et en $\lambda\in [\sigma]_{t}$ par une double somme en $\lambda\in i{\cal A}_{M_{disc}}^*/(i{\cal A}_{M_{disc},F}^{\vee}+i{\cal A}_{\tilde{M}}^*)$ et en $t\in W^G(M_{disc},M'_{disc})$ tels que $\tilde{M}_{t}=\tilde{M}$ et $\sigma_{\lambda}\circ ad_{\gamma_{0}t^{-1}}\simeq \omega\sigma_{\lambda}$. Fixons $\lambda$. L'application $t\mapsto w$ est alors une bijection entre l'ensemble de sommation en $t$ et l'ensemble $W^{\tilde{M}}_{reg}(\sigma_{\lambda})$. On peut encore d\'ecomposer cette somme en une somme sur $r\in R^{\tilde{M}}(\sigma_{\lambda})$, que l'on \'ecrit $r=W_{0}^{M}(\sigma_{\lambda})w$, et une somme sur $w'\in W_{0}^{M}(\sigma_{\lambda})w\cap W_{reg}^{\tilde{M}}(\sigma_{\lambda})$. On voit que les seuls termes de l'expression qui d\'ependent de $w'$ sont $\epsilon_{\sigma_{\lambda}}(w')\vert det((1-w')_{\vert {\cal A}_{M_{disc}}^{\tilde{M}}}\vert ^{-1}$. Leur somme vaut $\vert W_{0}^M(\sigma_{\lambda})\vert \iota(\tau)$ o\`u $\tau$ est comme ci-dessus. D'o\`u
$$X(\tilde{M},M_{disc},\sigma)= (-1)^{a_{\tilde{M}}-a_{\tilde{G}}}\sum_{\lambda\in i{\cal A}_{M_{disc}}^*/(i{\cal A}_{M_{disc},F}^{\vee}+i{\cal A}_{\tilde{M}}^*)}\vert W_{0}^M(\sigma_{\lambda})\vert $$
$$\sum_{r\in R^{\tilde{M}}(\sigma_{\lambda})} \iota(\tau)\int_{i{\cal A}_{\tilde{M},F}^*}J_{\tilde{M}}^{\tilde{G}}(\pi_{\tau_{\xi}},f_{1},f_{2})\,d\xi.$$
Notons $\Pi_{disc}(M_{disc}(F))/i{\cal A}_{\tilde{M},F}^*$ l'ensemble des orbites dans $\Pi_{disc}(M_{disc}(F))$ pour l'action de $i{\cal A}_{\tilde{M},F}^*$. Pour $\underline{\sigma}\in \Pi_{disc}(M_{disc}(F))$, notons $Stab(i{\cal A}_{\tilde{M},F}^*,\underline{\sigma})$ son stabilisateur dans $i{\cal A}_{\tilde{M},F}^*$. Il ne d\'epend que de l'image de $\underline{\sigma}$ dans l'ensemble d'orbites pr\'ec\'edent. Notons ${\cal O}_{\sigma}$ l'orbite de $\sigma$ pour l'action de $i{\cal A}_{M_{disc},F}^*$. Elle se d\'ecompose en orbites pour l'action de $i{\cal A}_{\tilde{M},F}^*$, notons ${\cal O}_{\sigma}/i{\cal A}_{\tilde{M},F}^*$ l'ensemble de ces orbites. Notons que $Stab(i{\cal A}_{\tilde{M},F}^*,\underline{\sigma})$ est ind\'ependant de $\underline{\sigma}$ dans cet ensemble. L'application qui, \`a $\lambda$, associe l'orbite de $\sigma_{\lambda}$ pour  l'action de  $i{\cal A}_{\tilde{M},F}^*$,  a pour image  ${\cal O}_{\sigma}/i{\cal A}_{\tilde{M},F}^*$ et a un noyau d'ordre $\vert  Stab(i{\cal A}_{M_{disc},F}^*,\sigma)\vert \vert  Stab(i{\cal A}_{\tilde{M},F}^*,\sigma)\vert ^{-1} $. Alors
$$X(\tilde{M},M_{disc},\sigma)=  (-1)^{a_{\tilde{M}}-a_{\tilde{G}}}\vert Stab(i{\cal A}_{M_{disc},F}^*,\sigma)\vert$$
$$ \sum_{\underline{\sigma}\in {\cal O}_{\sigma}/i{\cal A}_{\tilde{M},F}^*}\vert  Stab(i{\cal A}_{\tilde{M},F}^*,\underline{\sigma})\vert ^{-1}\vert W_{0}^M(\underline{\sigma})\vert \sum_{r\in R^{\tilde{M}}(\underline{\sigma})} \iota(\tau)\int_{i{\cal A}_{\tilde{M},F}^*}J_{\tilde{M}}^{\tilde{G}}(\pi_{\tau_{\xi}},f_{1},f_{2})\,d\xi,$$
o\`u maintenant $\tau=(M_{disc},\underline{\sigma},\tilde{r})$.

Consid\'erons la formule 3.2(1). Le terme constant de son membre de droite est la somme sur $M_{disc}\in {\cal L}(M_{0})$ et sur $\sigma\in \Pi_{disc}(M_{disc}(F))/i{\cal A}_{M_{disc},F}^*$ de l'expression (10), multipli\'ee par $\vert W^{M_{disc}}\vert \vert W^G\vert ^{-1}\vert  Stab(i{\cal A}_{M_{disc},F}^*,\sigma)\vert^{-1}$. On peut l'\'ecrire
$$(11) \qquad \sum_{\tilde{M}; M_{0}\subset M}\vert W^M\vert \vert W^G\vert ^{-1}X(\tilde{M}),$$
o\`u
$$X(\tilde{M})=\sum_{M_{disc};M_{0}\subset M_{disc}\subset M}\vert W^{M_{disc}}\vert \vert W^M\vert ^{-1}$$
$$\sum_{\sigma\in \Pi_{disc}(M_{disc}(F))/i{\cal A}_{M_{disc},F}^*}\vert  Stab(i{\cal A}_{M_{disc},F}^*,\sigma)\vert^{-1} X(\tilde{M},M_{disc},\sigma).$$
Fixons $\tilde{M}$. La somme en $\sigma$ et celle en $\underline{\sigma}\in {\cal O}_{\sigma}/i{\cal A}_{\tilde{M},F}^*$ figurant dans $X(\tilde{M},M_{disc},\sigma)$ se simplifient en une somme sur l'ensemble $\Pi_{disc}(M_{disc}(F))/i{\cal A}_{\tilde{M},F}^*$.  Consid\'erons les triplets $\tau$ qui apparaissent. Ils ne sont pas forc\'ement essentiels, mais on peut se limiter aux essentiels (pour $\tilde{M}$) puisque la contribution des autres est nulle d'apr\`es 3.25(1).  Ils sont alors discrets, puisque sinon, le terme $\iota(\tau)$ est nul. On a donc
$$(12)\qquad X(\tilde{M})=  (-1)^{a_{\tilde{M}}-a_{\tilde{G}}}\sum_{\tau\in (E_{disc}(\tilde{M},\omega)/conj)/i{\cal A}_{\tilde{M},F}^*}C(\tau)\iota(\tau)\int_{i{\cal A}_{\tilde{M},F}^*}J_{\tilde{M}}^{\tilde{G}}(\pi_{\tau_{\xi}},f_{1},f_{2})\,d\xi,$$
o\`u $C(\tau)$ est la somme de 
$$(13)\qquad \vert W^{M_{disc}}\vert \vert W^M\vert ^{-1}\vert Stab(i{\cal A}_{\tilde{M},F}^*,\underline{\sigma}) \vert ^{-1}\vert W_{0}^M(\underline{\sigma})\vert $$
sur les $M_{disc}$ tels que $M_{0}\subset M_{disc}\subset M$, $\underline{\sigma}\in \Pi_{disc}(M_{disc}(F))/i{\cal A}_{\tilde{M},F}^*$, $\tilde{r}\in R^{\tilde{G}}(\underline{\sigma})$ tels que $(M_{disc},\underline{\sigma},\tilde{r})$ ait pour image $\tau$ dans $(E_{disc}(\tilde{M},\omega)/conj)/i{\cal A}_{\tilde{M},F}^*$. Si on fixe un repr\'esentant $(M_{disc},\underline{\sigma},\tilde{r})$ de $\tau$, l'ensemble de sommation est celui des diff\'erents $(w(M_{disc}),w\underline{\sigma}, w(\tilde{r}))$ pour $w\in W^M$. Ou encore l'ensemble de ces triplets pour $w\in W^M/Stab(W^M,\tau/i{\cal A}_{\tilde{M},F}^*)$, o\`u $Stab(W^M,\tau/i{\cal A}_{\tilde{M},F}^*)$ est le stabilisateur de $\tau\in E_{disc}(\tilde{M},\omega)/i{\cal A}_{\tilde{M},F}^*$ dans $W^M$. Le terme (13) est constant sur l'ensemble de sommation. D'o\`u
$$C(\tau)=\vert W^M/Stab(W^M,\tau/i{\cal A}_{\tilde{M},F}^*)\vert \vert W^{M_{disc}}\vert \vert W^M\vert ^{-1}\vert  Stab(i{\cal A}_{\tilde{M},F}^*,\underline{\sigma})\vert ^{-1}\vert W_{0}^M(\underline{\sigma})\vert$$
$$=\vert Stab(W^M,\tau/i{\cal A}_{\tilde{M},F}^*)\vert ^{-1} \vert W^{M_{disc}}\vert\vert Stab(i{\cal A}_{\tilde{M},F}^*,\underline{\sigma})\vert ^{-1}\vert W_{0}^M(\underline{\sigma})\vert.$$
On a d\'efini en 2.9 le groupe $Stab(W^M\times i{\cal A}_{\tilde{M},F}^*,\tau)$ des $(w,\lambda)\in W^M\times i{\cal A}_{\tilde{M},F}^*$ qui conservent $\tau$.  Prenons garde que dans la d\'efinition de ce groupe, $\tau$ est consid\'er\'e comme un \'el\'ement de $E_{disc}(\tilde{M},\omega)$ tandis que, dans celle ci-dessus de $Stab(W^M,\tau/i{\cal A}_{\tilde{M},F}^*)$, $\tau$ est consid\'er\'e comme un \'el\'ement du quotient $E_{disc}(\tilde{M},\omega)/i{\cal A}_{\tilde{M},F}^*$.   On voit qu'il y a une suite exacte
$$1\to Stab(i{\cal A}_{\tilde{M},F}^*,\underline{\sigma})\to  Stab(W^M\times i{\cal A}_{\tilde{M},F}^*,\tau)\to Stab(W^M,\tau/i{\cal A}_{\tilde{M},F}^*)\to 1$$
Donc 
$$\vert Stab(W^M,\tau/i{\cal A}_{\tilde{M},F}^*)\vert  \vert  Stab(i{\cal A}_{\tilde{M},F}^*,\underline{\sigma})\vert =\vert Stab(W^M\times i{\cal A}_{\tilde{M},F}^*,\tau)\vert .$$ 
En appliquant les d\'efinitions de 2.9, on a aussi
$$\vert Stab(W^M\times i{\cal A}_{\tilde{M},F}^*,\tau)\vert =\vert W^{M_{disc}}\vert \vert W_{0}^M(\underline{\sigma})\vert \vert  {\bf Stab}(W^M\times i{\cal A}_{\tilde{M},F}^*,\tau)\vert .$$
On conclut
$$C(\tau)=\vert {\bf  Stab}(W^M\times i{\cal A}_{\tilde{M},F}^*,\tau)\vert .$$
Alors (12) devient
$$X(\tilde{M})=(-1)^{a_{\tilde{M}}-a_{\tilde{G}}}J_{\tilde{M},spec}^{\tilde{G}}(\omega,f_{1},f_{2}).$$
Ce terme ne d\'ependant que de la classe de conjugaison de $\tilde{M}$, on v\'erifie que la somme (11) est \'egale \`a 
$$(14) \qquad \sum_{\tilde{M}\in{\cal L}( \tilde{M}_{0})}\vert \tilde{W}^M\vert \vert \tilde{W}^G\vert ^{-1}X(\tilde{M})$$
En effet, les deux termes sont \'egaux \`a
$$\sum_{\tilde{P}=\tilde{M}U_{P};P_{0}\subset P}\vert {\cal P}(\tilde{M})\vert ^{-1}X(\tilde{M}).$$ 
 On se rappelle que (11) est le terme constant de l'\'el\'ement de $PolExp$ asymptote \`a $J^T(\omega,f_{1},f_{2})$. D'apr\`es la formule ci-dessus calculant $X(\tilde{M})$, (14) est \'egal \`a $J^{\tilde{G}}_{spec}(\omega,f_{1},f_{2})$. Cela prouve la proposition. $\square$
 
 \bigskip
 
 \section{Le calcul g\'eom\'etrique}
 
 \bigskip
 
 \subsection{La formule de Weyl}
Pour quelques instants, nous allons consid\'erer des sous-groupes alg\'ebriques de $G$ qui ne sont pas forc\'ement d\'efinis sur $F$. Consid\'erons un sous-groupe de Borel $B$ de $G$ et un sous-tore maximal $T$ de $B$. On suppose $T$ d\'efini sur $F$ mais $B$ d\'efini seulement sur la cl\^oture alg\'ebrique $\bar{F}$ de $F$. On d\'efinit leur normalisateur commun $\tilde{T}=\{\gamma\in \tilde{G}; ad_{\gamma}(T)=T, ad_{\gamma}(B)=B\}$. On dit que $\tilde{T}$ est un tore tordu maximal de $\tilde{G}$ si et seulement si $\tilde{T}\cap \tilde{G}(F)\not=\emptyset$. Remarquons que cela entra\^{\i}ne que $\tilde{T}$ est d\'efini sur $F$: si on choisit $\gamma\in \tilde{T}\cap \tilde{G}(F)$, on a $\tilde{T} =T\gamma=\gamma T$. Remarquons aussi qu'\`a tout tore tordu maximal $\tilde{T}$ est associ\'e un sous-tore maximal $T$ de $G$. Pour $\gamma\in \tilde{T}(F)$, la restriction de $ad_{\gamma}$ \`a $\tilde{T}$ ne d\'epend pas de $\gamma$. On la note   simplement $\theta$. On note $A_{T}$ le plus grand sous-tore d\'eploy\'e de $T$ et $A_{\tilde{T}}$ le plus grand sous-tore contenu dans le sous-groupe $A_{T}^{\theta}$ des points fixes par $\theta$ (autrement dit, $A_{\tilde{T}}$ est la composante neutre $A_{T}^{\theta,0}$ de $A_{T}^{\theta}$). On dit que $\tilde{T}$ est elliptique dans $\tilde{G}$ si $A_{\tilde{T}}=A_{\tilde{G}}$. En g\'en\'eral, notons $M$ le commutant de $A_{\tilde{T}}$ dans $G$ et posons $\tilde{M}=M\tilde{T}$. Alors $\tilde{M}$ est un espace de Levi et $\tilde{T}$ est un tordu maximal et elliptique de $\tilde{M}$. 

Pour tout tore tordu maximal $\tilde{T}$, on munit le groupe $T^{\theta}(F)$ d'une mesure de Haar. On suppose que, si deux tores tordus maximaux sont conjugu\'es par un \'el\'ement de $G(F)$, cette conjugaison est compatible aux mesures.

Soit $\tilde{T}$ un tore tordu maximal. Le groupe $T(F)$ agit sur $\tilde{T}(F)$ par conjugaison. L'ensemble des classes de conjugaison est le quotient $\tilde{T}(F)/(1-\theta)(T(F))=(1-\theta)(T(F))\backslash \tilde{T}(F)$, o\`u $(1-\theta)(T(F))=\{t\theta(t^{-1}); t\in T(F)\}$.  Pour $\gamma\in \tilde{T}(F)$, l'application naturelle
$$\begin{array}{ccc}T^{\theta,0}(F)&\to& \tilde{T}(F)/(1-\theta)(T(F))\\ t&\mapsto& t\gamma\\ \end{array}$$
est un isomorphisme local.   On munit $\tilde{T}(F)/(1-\theta)(T(F))$ de la mesure invariante par translation \`a droite ou \`a gauche par $T(F)$ et telle que l'application ci-dessus pr\'eserve localement les mesures. Il est imm\'ediat que cette d\'efinition ne d\'epend pas du choix de $\gamma$. On pose 
$$W^{G}(\tilde{T})=Norm_{G(F)}(\tilde{T})/T(F).$$

L'ensemble des classes de conjugaison par $G(F)$ de tores tordus maximaux est fini. Fixons un ensemble de repr\'esentants $T(\tilde{G})$. On note $T_{ell}(\tilde{G})$ le sous-ensemble des \'el\'ements elliptiques de $T(\tilde{G})$.

Pour $\gamma\in \tilde{G}$, on note $Z_{G}(\gamma)$ son centralisateur dans $G$ (l'ensemble des $x\in G$ tels que $x\gamma x^{-1}=\gamma$) et $G_{\gamma}$ la composante neutre de $Z_{G}(\gamma)$. On dit que $\gamma$ est fortement r\'egulier si et seulement si $G_{\gamma}$ est un tore et $Z_{G}(\gamma)$ est commutatif. Un tel \'el\'ement est semi-simple. On note $G_{reg}$  l'ensemble des \'el\'ements fortement r\'eguliers. Un \'el\'ement $\gamma\in G_{reg}(F)$ appartient \`a un unique tore tordu maximal: c'est $\tilde{T}=T\gamma$, o\`u $T$ est le commutant de $G_{\gamma}$ dans $G$.  On d\'efinit une fonction $D^{\tilde{G}}$ sur $\tilde{G}(F)$ de la fa\c{c}on suivante. Si $\gamma\in \tilde{G}(F)$ est semi-simple, $D^{\tilde{G}}(\gamma)$ est la valeur absolue du d\'eterminant de $1-ad_{\gamma}$ agissant sur $\mathfrak{g}/\mathfrak{g}_{\gamma}$, o\`u on note par des lettres gothiques les alg\`ebres de Lie. Si $\gamma$ est quelconque, $D^{\tilde{G}}(\gamma)=D^{\tilde{G}}(\gamma_{ss})$, o\`u $\gamma_{ss}$ est la partie semi-simple de $\gamma$.

 La formule d'int\'egration de Weyl prend l'une ou l'autre des deux formes suivantes, pour une fonction $f\in C_{c}^{\infty}(\tilde{G}(F))$:
 
 $$\int_{\tilde{G}(F)}f(\gamma)\,d\gamma\,=\sum_{\tilde{T}\in T(\tilde{G})}\vert W^G(\tilde{T})\vert ^{-1} \int_{\tilde{T}(F)/(1-\theta)(T(F))}\int_{T^{\theta}(F)\backslash G(F)}f(x^{-1}\gamma x)\,dx\,D^{\tilde{G}}(\gamma)\,d\gamma$$
 
 $=\sum_{\tilde{M}\in {\cal L}(\tilde{M}_{0})}\vert \tilde{W}^M\vert \vert \tilde{W}^G\vert ^{-1}\sum_{\tilde{T}\in T_{ell}(\tilde{M})}\vert W^M(\tilde{T})\vert ^{-1}$
 $$ \int_{\tilde{T}(F)/(1-\theta)(T(F))}\int_{T^{\theta}(F)\backslash G(F)}f(x^{-1}\gamma x)\,dx\,D^{\tilde{G}}(\gamma)\,d\gamma.$$
 
 \bigskip
 
 \subsection{Quelques majorations}
 
 Soient $f\in C_{c}^{\infty}(\tilde{G}(F))$ et $\gamma\in \tilde{G}_{reg}(F)$ tel que $\omega$ soit trivial sur $Z_{G}(\gamma,F)$. On pose
 $$I_{\tilde{G}}(\gamma,\omega,f)=D^{\tilde{G}}(\gamma)^{1/2}\int_{Z_{G}(\gamma,F)\backslash G(F)}\omega(x) f(x^{-1}\gamma x)\,dx.$$
 Par d\'efinition, c'est l'int\'egrale orbitale de $f$ en $\gamma$.   Si $\omega$ n'est pas trivial sur $Z_{G}(\gamma,F)$,  on pose $I_{\tilde{G}}(\gamma,\omega,f)=0$.
 
 Dans le cas o\`u le caract\`ere $\omega$ est trivial, on  note simplement $I_{\tilde{G}}(\gamma,f)=I_{\tilde{G}}(\gamma,\omega,f)$.
   
 Soit $\tilde{T}$ un tore tordu maximal de $\tilde{G}$. Supposons $\omega$ trivial sur $T^{\theta}(F)$. On a
 
 (1) la fonction $\gamma\mapsto I_{\tilde{G}}(\gamma,\omega,f)$ est born\'ee sur $\tilde{T}(F)$. 
 
 Preuve. Il existe une fonction $f'\in C_{c}^{\infty}(\tilde{G}(F))$, \`a valeurs positives ou nulles, de sorte que $\vert f(\gamma)\vert \leq f'(\gamma)$ pour tout $\gamma\in \tilde{G}(F)$. Alors
 $$\vert I_{\tilde{G}}(\gamma,\omega,f)\vert\leq I_{\tilde{G}}(\gamma,f').  $$
 Il suffit de majorer le membre de droite.   En oubliant cette construction, on peut supposer $\omega=1$ et $f$ \`a valeurs positives ou nulles. La fonction $\gamma\mapsto I_{\tilde{G}}(\gamma,f)$ sur $\tilde{T}(F)$ se quotiente en une fonction sur $\tilde{T}(F)/(1-\theta)(T(F))$ qui est \`a support compact. Il suffit de v\'erifier qu'elle est localement born\'ee. Il suffit de fixer $\gamma\in \tilde{T}(F)$ et de trouver un voisinage $U$ de $0$ dans $\mathfrak{t}^{\theta}(F)$ tel que la fonction soit born\'ee sur $exp(U)\gamma$. Mais la th\'eorie de la descente vaut dans le cas tordu. Elle entra\^{\i}ne que l'on peut trouver un tel voisinage $U$ et une fonction $\varphi\in  C_{c}^{\infty}(\mathfrak{g}_{\gamma}(F))$ de sorte que
 $$I_{\tilde{G}}(exp(X)\gamma,f)=I_{G_{\gamma_{0}}}(X,\varphi)$$
 pour tout $X\in U$ tel que $exp(X)\gamma\in \tilde{G}_{reg}(F)$, avec une d\'efinition \'evidente de l'int\'egrale de droite. Le membre de droite est born\'e d'apr\`es Harish-Chandra. $\square$
 
 On a
 
 (2) il existe $\eta>0$ tel que la fonction $\gamma\mapsto D^{\tilde{G}}(\gamma)^{-\eta}$ soit localement int\'egrable sur $\tilde{T}(F)/(1-\theta)(T(F))$.
 
 Preuve. Il suffit de prouver l'existence de $\eta>0$ tel que, pour tout $\gamma\in \tilde{T}(F)$, la fonction $X\mapsto D^{\tilde{G}}(exp(X)\gamma)^{-\eta}$ est int\'egrable au voisinage de $0$ dans $\mathfrak{t}^{\theta}(F)$. Pour $X$ proche de $0$, on a $D^{\tilde{G}}(exp(X)\gamma)=D^{\tilde{G}}(\gamma)D^{G_{\gamma}}(exp(X))$. D'apr\`es Harish-Chandra, il existe $\eta_{\gamma}>0$ tel que $X\mapsto D^{G_{\gamma}}(exp(X))^{-\eta_{\gamma}}$ soit int\'egrable au voisinage de $0$. Le r\'eel $\eta_{\gamma}$ ne d\'epend en fait que du groupe $G_{\gamma}$. Or il n'y a qu'un nombre fini de tels commutants possibles. En prenant pour $\eta$ le plus petit des r\'eels $\eta_{\gamma}$, on obtient (2). $\square$

 On fixe une norme $\vert \vert .\vert \vert $ sur $G(F)$ selon la m\'ethode habituelle. On suppose qu'elle est biinvariante par $K$ et que $\vert \vert g\vert \vert \geq1$ pour tout $g\in G(F)$. On rappelle que, pour tout $g\in G(F)$, on a d\'efini $h_{0}(g)\in {\cal A}_{0}^{\geq}$. Les fonctions $1+log(\vert \vert g\vert \vert )$ et $1+\vert h_{0}(g)\vert $ sont \'equivalentes. Fixons $\gamma_{0}\in \tilde{G}(F)$. On d\'efinit une norme sur $\tilde{G}(F)$ par $\vert \vert \gamma\vert \vert =\vert \vert g\vert \vert $ si $\gamma=g\gamma_{0}$ avec $g\in G(F)$. Elle d\'epend du point-base $\gamma_{0}$ mais sa classe d'\'equivalence (en un sens plus ou moins clair) n'en d\'epend pas.  
  
 \ass{Lemme}{Soit $\tilde{T}$ un tore tordu maximal de $\tilde{G}$.
Il existe deux entiers $N,k>0$ et un r\'eel $c>0$ tels que l'on ait les majorations
 
 (i) $inf_{t\in T^{\theta,0}(F)}\vert \vert tx\vert \vert \leq c (inf_{t\in T^{\theta,0}(F)}\vert \vert t\gamma\vert \vert)^N\vert \vert x^{-1}\gamma x\vert \vert ^ND^{\tilde{G}}(\gamma)^{-k}$
 
 (ii)  $inf_{t\in T(F)}\vert \vert tx\vert \vert \leq c \vert \vert x^{-1}\gamma x\vert \vert ^ND^{\tilde{G}}(\gamma)^{-k}$
 
 \noindent pour tout $x\in G(F)$ et tout $\gamma\in \tilde{T}(F)\cap \tilde{G}_{reg}(F)$. }
 
 Preuve.  La preuve reprend celle du lemme 4.2 de [A1]. Consid\'erons une extension galoisienne finie $F'$ de $F$ que l'on pr\'ecisera plus tard.   On fixe comme pr\'ec\'edemment des normes sur $G(F')$ et $\tilde{G}(F')$, notons-les pour un instant $\vert \vert .\vert \vert_{F'}$. Il existe un entier $N_{1}\geq1$  tel que l'on ait les majorations
 $$\vert \vert x\vert \vert <<\vert \vert x\vert \vert _{F'}^{N_{1}},\,\,\vert \vert x\vert \vert _{F'}<<\vert \vert x\vert \vert ^{N_{1}}$$
 pour tout $x\in G(F)$ et on a des majorations similaires pour $\gamma\in \tilde{G}(F)$. On voit que l'\'enonc\'e est \'equivalent \`a celui obtenu en rempla\c{c}ant $\vert \vert .\vert \vert $ par $\vert \vert .\vert \vert _{F'}$ dans les in\'egalit\'es (i) et (ii). On peut donc travailler avec les seules normes $\vert \vert .\vert \vert _{F'}$ et  simplifier la notation pour le reste de la preuve en abandonnant les indices $F'$.
 
 Soit $S$ un sous-tore de $G$ d\'efini sur $F$, pas forc\'ement maximal. On montre d'abord
 
 (3) il existe un entier $N_{2}$ tel que l'on ait la majoration
$$inf_{t\in S(F)}\vert \vert tx\vert \vert <<( inf_{t\in S(F')}\vert \vert tx\vert \vert)^{N_{2}} $$
pour tout $x\in G(F)$.

Notons $d$ le degr\'e de l'extension $F'/F$ et $S(F')^d$ le sous-groupe des puissances $d$-i\`emes dans $S(F')$. Le quotient $S(F')/S(F')^d$ est compact, d'o\`u une majoration
$$(4) \qquad inf_{t\in S(F')^d}\vert \vert tx\vert \vert<<inf_{t\in S(F')}\vert \vert tx\vert \vert$$
pour tout $x\in G(F')$. Fixons $x\in G(F)$, posons $X=inf_{t\in S(F')^d}\vert \vert tx\vert \vert$ et choisissons $u\in S(F')$ tel que $\vert \vert u^dx\vert \vert =X$ (il est \`a peu pr\`es clair que la borne inf\'erieure est atteinte). Pour $\sigma\in Gal(F'/F)$, l'application $y\mapsto \vert \vert \sigma(y)\vert \vert $ est encore une norme sur $G(F')$, il y a donc un entier $N_{3}\geq1$ tel que $\vert \vert \sigma(y)\vert \vert <<\vert \vert y\vert \vert ^{N_{3}}$ pour tout $y\in G(F')$ et tout $\sigma\in Gal(F'/F)$. Donc $\vert \vert \sigma(u)^dx\vert \vert <<X^{N_{3}}$ pour tout $\sigma$. Parce que $\sigma(u)^du^{-d}=\sigma(u)^dx(u^dx)^{-1}$, on en d\'eduit l'existence de $N_{4}\geq1$ tel que $\vert \vert \sigma(u)^du^{-d}\vert \vert <<X^{N_{4}}$. On v\'erifie qu'il existe $N_{5}\geq1$ tel que $\vert \vert v\vert \vert <<\vert \vert v^d\vert \vert ^{N_{5}}$ pour tout $v\in S(F')$. D'o\`u $\vert \vert \sigma(u)u^{-1}\vert \vert<<X^{N_{4}N_{5}}$. Posons $v=\prod_{\sigma\in Gal(F'/F)}\sigma(u)$. On a $vx=(\prod_{\sigma\in Gal(F'/F)}\sigma(u)u^{-1}) u^dx$, d'o\`u l'existence de $N_{6}\geq1$ tel que $\vert \vert vx\vert \vert <<X^{N_{6}}$. Mais $v\in S(F)$. Donc $inf_{t\in S(F)}\vert \vert tx\vert \vert <<X^{N_{6}}$. Jointe \`a (4), cette relation d\'emontre (3).

On choisit $F'$ de sorte que $T$ soit d\'eploy\'e sur $F'$. En appliquant (3) \`a $S=T^{\theta,0}$ ou $S=T$, on voit que les assertions de l'\'enonc\'e sont \'equivalentes aux m\^emes assertions o\`u l'on remplace le corps de base $F$ par $F'$. En oubliant cette construction, on peut supposer que $T$ est d\'eploy\'e sur $F$. On fixe un sous-groupe de Borel $B=TU$ de $G$ contenant $T$. La d\'ecomposition d'Iwasawa montre que l'on peut se limiter \`a prouver l'\'enonc\'e pour $x\in B(F)$. Ecrivons $x=yu$ avec $y\in T(F)$ et $u\in U(F)$. Pour $S=T^{\theta,0}$ ou $S=T$, il existe $N_{7}\geq1$ tel que
$$(5)\qquad inf_{t\in S(F)}\vert \vert tx\vert \vert<<(inf_{t\in S(F)}\vert \vert ty\vert \vert)^{N_{7}}\vert \vert u\vert \vert ^{N_{7}}.$$
On a $x^{-1}\gamma x=u'\gamma' $, o\`u $\gamma'=y^{-1} \gamma y$ et $u'=u^{-1}ad_{\gamma'}(u)$. Il r\'esulte de la filtration habituelle du groupe unipotent $U$ que les coefficients de $u$ sont des fractions rationnelles en les coefficients de $\gamma'$ et $u'$, avec pour d\'enominateurs des puissances de $det((ad_{\gamma'}-1)_{\vert  \mathfrak{g}/\mathfrak{t}^{\theta}})$. Il y a donc des entiers $N_{8},k_{1}\geq1$ tels que $\vert \vert u\vert \vert <<\vert \vert x^{-1}\gamma x\vert \vert^{N_{8}}D^{\tilde{G}}(\gamma')^{-k_{1}}$.  Puisque $D^{\tilde{G}}(\gamma')=D^{\tilde{G}}(\gamma)$, la relation (5) nous ram\`ene \`a prouver l'existence d'un entier $N_{9}\geq1$ tel que
$$(6) \qquad inf_{t\in S(F)}\vert \vert ty\vert \vert<<(inf_{t\in S(F)}\vert \vert t\gamma\vert \vert)^{N_{9}}\vert \vert y^{-1}\gamma y\vert \vert ^{N_{9}}.$$
Le cas o\`u $S=T$ (qui concerne le (ii) de l'\'enonc\'e) est trivial puisque le terme de gauche vaut $1$. On suppose maintenant $S=T^{\theta,0}$. Parce que les deux groupes $(1-\theta)(T(F))$ et $T^{\theta,0}(F)$ sont d'intersection finie et que leur produit est un sous-groupe d'indice fini de $T(F)$, on v\'erifie qu'il existe $N_{10}\geq1$ tel que
$$ inf_{t\in T^{\theta,0}(F)}\vert \vert ty\vert \vert<<( inf_{t\in T^{\theta,0}(F)}\vert \vert ty^{-1}\theta(y)\vert \vert)^{N_{10}}.$$
Soit $t'\in T^{\theta,0}(F)$ tel que $\vert \vert (t')^{-1}\gamma\vert \vert $ soit minimal, posons $\gamma'=(t')^{-1}\gamma$. Alors
$$ inf_{t\in T^{\theta,0}(F)}\vert \vert ty\vert \vert<<\vert \vert t'y^{-1}\theta(y)\vert \vert ^{N_{10}}<<\vert \vert t'y^{-1}\theta(y)\gamma'\vert \vert ^{N_{11}}\vert \vert \gamma'\vert \vert ^{N_{11}}$$
pour un certain $N_{11}\geq1$. Mais $t'y^{-1}\theta(y)\gamma'=y^{-1}\gamma y$ et $\vert \vert \gamma'\vert \vert =inf_{t\in T^{\theta,0}(F)}\vert \vert t\gamma\vert \vert$. La relation pr\'ec\'edente n'est autre que (6). Cela ach\`eve la preuve. $\square$

  Le lemme entra\^{\i}ne \'evidemment
 
 (7) il existe un entier $k\geq0$ et, pour tous sous-ensembles compacts $\boldsymbol{\Omega}$ de $\tilde{G}(F)$ et $\Omega$ de $\tilde{T}(F)$,  il existe $c>0$ tel que l'on ait la majoration
 $$inf_{t\in T^{\theta,0}(F)}\vert \vert tx\vert \vert \leq c D^{\tilde{G}}(\gamma)^{-k}$$
 pour tout couple $(x,\gamma)$ tel que $x\in G(F)$, $\gamma\in  \Omega\cap \tilde{G}_{reg}(F)$ et $x^{-1}\gamma x\in \boldsymbol{\Omega}$.

 \bigskip
 \subsection{Application de la formule de Weyl}
 On a d\'efini l'int\'egrale $J^T(\omega,f_{1},f_{2})$ en 3.1. C'est une int\'egrale \`a support compact.  On applique la formule de Weyl sous sa deuxi\`eme forme \`a l'int\'egrale int\'erieure en $\gamma$. On obtient
 $$(1) \qquad J^T(\omega,f_{1},f_{2})=\sum_{\tilde{M}\in {\cal L}(\tilde{M}_{0})}\vert \tilde{W}^M\vert \vert \tilde{W}^G\vert ^{-1}\sum_{\tilde{S}\in T_{ell}(\tilde{M})}J_{\tilde{M},\tilde{S}}^T(\omega,f_{1},f_{2}),$$
 o\`u
 $$J_{\tilde{M},\tilde{S}}^T(\omega,f_{1},f_{2})=\vert W^M(\tilde{S})\vert ^{-1} \int_{\tilde{S}(F)/(1-\theta)(S(F))}\int_{S^{\theta}(F)\backslash G(F)}$$
 $$\int_{A_{\tilde{G}}(F)\backslash G(F)}\bar{f}_{1}(x^{-1}\gamma x)f_{2}(g^{-1}x^{-1}\gamma xg)\omega(g)\tilde{\kappa}^T(g)\,dg\,dx\,D^{\tilde{G}}(\gamma)\,d\gamma.$$
 On a not\'e nos tores $\tilde{S}$ plut\^ot que $\tilde{T}$ pour \'eviter les confusions avec le param\`etre $T$. 
 
 On fixe jusqu'en 4.8 un espace de Levi $\tilde{M}$ contenant $\tilde{M}_{0}$ et un tore tordu  elliptique $\tilde{S}$ de $\tilde{M}$. Par changement de variables $g\mapsto x^{-1}y$, on a
 $$J_{\tilde{M},\tilde{S}}^T(\omega,f_{1},f_{2})=\vert W^M(\tilde{S})\vert ^{-1} \int_{\tilde{S}(F)/(1-\theta)(S(F))}\int_{S^{\theta}(F)\backslash G(F)}$$
 $$\int_{A_{\tilde{G}}(F)\backslash G(F)}\bar{f}_{1}(x^{-1}\gamma x)f_{2}(y^{-1}\gamma y)\omega(x^{-1}y)\tilde{\kappa}^T(x^{-1}y)\,dy\,dx\,D^{\tilde{G}}(\gamma)\,d\gamma.$$
 Pour tout $\gamma$, la fonction $y\mapsto  f_{2}(y^{-1}\gamma y)$ est invariante par $S^{\theta}(F)$, a fortiori par $A_{\tilde{S}}(F)=A_{\tilde{M}}(F)$. D'o\`u
 $$J_{\tilde{M},\tilde{S}}^T(\omega,f_{1},f_{2})=\vert W^M(\tilde{S})\vert ^{-1} \int_{\tilde{S}(F)/(1-\theta)(S(F))}\int_{S^{\theta}(F)\backslash G(F)}$$
 $$\int_{A_{\tilde{M}}(F)\backslash G(F)}\bar{f}_{1}(x^{-1}\gamma x)f_{2}(y^{-1}\gamma y)\omega(x^{-1}y)u_{\tilde{M}}^T(x,y) \,dy\,dx\,D^{\tilde{G}}(\gamma)\,d\gamma,$$
 o\`u
 $$u_{\tilde{M}}^T(x,y)=\int_{A_{\tilde{G}}(F)\backslash A_{\tilde{M}}(F)}\omega(a)\tilde{\kappa}^T(x^{-1}ay)\,da.$$

 Soient $x,y\in G(F)$. Pour $\tilde{P}\in {\cal F}(\tilde{M})$, on a d\'efini $T[\tilde{P}]$ en 1.3. On pose $Y(x,y,T;\tilde{P})=T[\tilde{P}]+H_{\tilde{P}}(x)-H_{\tilde{\bar{P}}}(y)$. La famille ${\cal Y}(x,y,T)=(Y(x,y,T;\tilde{P}))_{\tilde{P}\in {\cal P}(\tilde{M})}$ est $(\tilde{G},\tilde{M})$-orthogonale. On introduit la  fonction $H\mapsto \Gamma_{\tilde{M}}(H,{\cal Y}(x,y,T))$ sur ${\cal A}_{\tilde{M}}$, cf. 2.2. Rappelons que l'on note $A_{\tilde{M}}(F)_{c}$ le plus grand sous-groupe compact de $A_{\tilde{M}}(F)$, c'est-\`a-dire le sous-groupe des $a\in A_{\tilde{M}}(F)$ tels que $H_{\tilde{M}}(a)=0$. Si la restriction de $\omega$ \`a ce groupe est non triviale, on pose
 $$v_{\tilde{M}}^T(x,y)=0.$$
 Supposons maintenant que la restriction de $\omega$ \`a $A_{\tilde{M}}(F)_{c}$ est triviale. Alors $\omega$ se factorise en un caract\`ere de ${\cal A}_{A_{\tilde{M}},F}$ que l'on note encore $\omega$. On pose alors
 $$v_{\tilde{M}}^T(x,y)=C\int_{{\cal A}_{A_{\tilde{G}},F}\backslash {\cal A}_{A_{\tilde{M}},F}}\Gamma_{\tilde{M}}(H,{\cal Y}(x,y,T))\omega(H)\,dH,$$
 o\`u
 $$C=mes(A_{\tilde{M}}(F)_{c})mes(A_{\tilde{G}}(F)_{c})^{-1}.$$
 
 Posons
 $$J_{v,\tilde{M},\tilde{S}}^T(\omega,f_{1},f_{2})=\vert W^M(\tilde{S})\vert ^{-1} \int_{\tilde{S}(F)/(1-\theta)(S(F))}\int_{S^{\theta}(F)\backslash G(F)}$$
 $$\int_{A_{\tilde{M}}(F)\backslash G(F)}\bar{f}_{1}(x^{-1}\gamma x)f_{2}(y^{-1}\gamma y)\omega(x^{-1}y)v_{\tilde{M}}^T(x,y) \,dy\,dx\,D^{\tilde{G}}(\gamma)\,d\gamma.$$
 
 \ass{Proposition}{L'expression ci-dessus est absolument convergente. Pour tout r\'eel $r$, on a la majoration
 $$\vert J_{\tilde{M},\tilde{S}}^T(\omega,f_{1},f_{2})-J_{v,\tilde{M},\tilde{S}}^T(\omega,f_{1},f_{2})\vert << \vert T\vert ^{-r}$$
 pour tout $T$.}
 
 Cette proposition sera prouv\'ee en 4.6.
 
  Dans les formules d\'efinissant $J_{\tilde{M},\tilde{S}}^T(\omega,f_{1},f_{2})$ et $J_{v,\tilde{M},\tilde{S}}^T(\omega,f_{1},f_{2})$, les variables d'int\'egration $\gamma$, $x$ et $y$ appartiennent \`a des ensembles quotients mais il convient de les relever dans $\tilde{S}(F)$, resp. $G(F)$. On fixe un sous-ensemble compact $\boldsymbol{\Omega}$ de $\tilde{G}(F)$ contenant les supports de $f_{1}$ et $f_{2}$. L'ensemble des \'el\'ements de $\tilde{S}(F)/(1-\theta)(S(F))$ dont la classe de conjugaison coupe $\boldsymbol{\Omega}$ est compact. On peut donc fixer un sous-ensemble compact $\Omega$ de $\tilde{S}(F)$ et supposer 
  
  (2) $\gamma\in \Omega$.
  
   Appliquons 4.2(7). Puisque $S^{\theta,0}(F)/A_{\tilde{M}}(F)$ est compact, on peut remplacer dans cette relation le groupe $S^{\theta,0}(F)$ par $A_{\tilde{M}}(F)$. Cela nous permet de supposer 
 $$(3) \qquad \vert \vert x\vert \vert ,\vert \vert y\vert \vert \leq cD^{\tilde{G}}(\gamma)^{-k}.$$

 \bigskip
 
 \subsection{Une majoration de $u_{\tilde{M}}^T(x,y)$ et $v_{\tilde{M}}^T(x,y)$}
 Fixons $\epsilon>0$. Notons $\chi^T_{\epsilon}$ la fonction caract\'eristique des $\gamma\in\tilde{S}(F)/(1-\theta)(S(F))$ tels que $D^{\tilde{G}}(\gamma)\leq e^{-\epsilon\vert T\vert }$.
 
\ass{Lemme}{Pour tout r\'eel $r$, on a la majoration
$$\int_{\tilde{S}(F)/(1-\theta)(S(F))}\int_{S^{\theta}(F)\backslash G(F)}\int_{A_{\tilde{M}}(F)\backslash G(F)}\vert f_{1}(x^{-1}\gamma x)f_{2}(y^{-1}\gamma y)\vert $$
$$(\vert u_{\tilde{M}}^T(x,y)\vert +\vert v_{\tilde{M}}^T(x,y)\vert  )\,dy\,dx\,\chi_{\epsilon}^T(\gamma)D^{\tilde{G}}(\gamma)\,d\gamma<<\vert T\vert ^{-r}$$
 pour tout $T$.}
 
 Preuve. Il suffit de traiter le cas o\`u $\omega=1$ et les fonctions $f_{1}$ et $f_{2}$ sont \`a valeurs positives ou nulles. Montrons qu'il existe un entier $D>0$ tel que, pour $x$ et $y$ v\'erifiant 4.3(3), on ait la majoration
 $$(1) \qquad \vert u_{\tilde{M}}^T(x,y)\vert <<(\vert T\vert +\vert log(D^{\tilde{G}}(\gamma))\vert)^D .$$
  On a une majoration
  $$log(\vert \vert g\vert \vert )<<\vert T\vert +\vert H_{G}(g)\vert $$
  pour tout $g\in G(F)$ tel que $\tilde{\kappa}^T(g)=1$. Le terme $u_{\tilde{M}}^T(x,y)$ est d\'efini par une int\'egration sur $A_{\tilde{M}}(F)/A_{\tilde{G}}(F)$.  Il existe un ensemble fini ${\bf b}\subset {\cal A}_{\tilde{G}}$ tel qu'en posant $A_{\tilde{M}}(F)^{{\bf b}}=\{a\in A_{\tilde{M}}(F); H_{\tilde{G}}(a)\in {\bf b}\}$, on ait
  $$\vert u_{\tilde{M}}^T(x,y)\vert<<\int_{A_{\tilde{M}}(F)^{{\bf b}}}\tilde{\kappa}^T(x^{-1}ay)\,da.$$
  Pour $a$ dans l'ensemble d'int\'egration, on a
  $$log(\vert \vert a\vert \vert )<<\vert T\vert +log(\vert \vert x\vert \vert )+log(\vert \vert y\vert \vert ),$$
  d'o\`u
   $$log(\vert \vert a\vert \vert )<<\vert T\vert +\vert log(D^{\tilde{G}}(\gamma))\vert $$
  d'apr\`es 4.3(3). L'int\'egrale est alors essentiellement born\'ee par la mesure des $H\in {\cal A}_{A_{\tilde{M}},F}$ tels que $\vert H\vert $ soit major\'e par le membre de droite ci-dessus. On en d\'eduit (1).
  
  Une preuve similaire montre que la fonction $v_{\tilde{M}}^T(x,y)$ v\'erifie la m\^eme propri\'et\'e.
 
D'apr\`es les explications de la fin du paragraphe pr\'ec\'edent, l'int\'egrale de l'\'enonc\'e est donc major\'ee par 
$$\int_{\tilde{S}(F)/(1-\theta)(S(F))}\int_{S^{\theta}(F)\backslash G(F)}\int_{A_{\tilde{M}}(F)\backslash G(F)} f_{1}(x^{-1}\gamma x)f_{2}(y^{-1}\gamma y)$$
 $$(\vert T\vert +\vert log(D^{\tilde{G}}(\gamma))\vert)^D \,dy\,dx\,\chi_{\epsilon}^T(\gamma)D^{\tilde{G}}(\gamma)\,d\gamma$$
 $$<<\int_{\tilde{S}(F)/(1-\theta)(S(F))}I_{\tilde{G}}(\gamma,f_{1})I_{\tilde{G}}(\gamma,f_{2})(\vert T\vert +\vert log(D^{\tilde{G}}(\gamma))\vert)^D \chi_{\epsilon}^T(\gamma)\,d\gamma.$$
 On choisit $\eta$ comme en 4.2(2). Dans le domaine o\`u les fonctions \`a int\'egrer ne sont pas nulles, on a
 $$(\vert T\vert +\vert log(D^{\tilde{G}}(\gamma))\vert)^D <<\vert T\vert ^DD^{\tilde{G}}(\gamma)^{-\eta/2},$$
 $$\chi_{\epsilon}^T(\gamma)\leq e^{-\epsilon\eta\vert T\vert /2}D^{\tilde{G}}(\gamma)^{-\eta/2}.$$
 En utilisant aussi 4.2(1), on voit que l'int\'egrale ci-dessus est donc essentiellement major\'ee par
 $$\vert T\vert ^De^{-\epsilon\eta\vert T\vert /2}\int_{\Omega}D^{\tilde{G}}(\gamma)^{-\eta}\,d\gamma.$$
 La derni\`ere int\'egrale est convergente et le lemme s'ensuit. $\square$

\bigskip

\subsection{Majoration de $u_{\tilde{M}}(x,y)-v_{\tilde{M}}(x,y)$}
\ass{Lemme}{Il existe $\epsilon>0$ tel que, pour tout r\'eel $r$, on ait la majoration
$$\vert u^T_{\tilde{M}}(x,y)-v^T_{\tilde{M}}(x,y)\vert <<\vert T\vert ^{-r}$$
pour tout $T$ et tous $x,y\in G(F)$ tels que $\vert \vert x\vert \vert, \vert \vert y\vert \vert  \leq e^{\epsilon\vert T\vert }$.}

Preuve. Fixons $\epsilon'>0$ que l'on pr\'ecisera plus tard. On note $\epsilon'{\cal T}$ la famille de points $(\epsilon'T[\tilde{P}])_{\tilde{P}\in {\cal P}(\tilde{M})}$. Pour $\tilde{Q}=\tilde{L}U_{Q}\in {\cal F}(\tilde{M})$, posons 
$$u^T_{\tilde{M}}(x,y,\tilde{Q})=\int_{A_{\tilde{M}}(F)/A_{\tilde{G}}(F)}\omega(a)\tilde{\kappa}^T(x^{-1}ay)\Gamma_{\tilde{M}}^{\tilde{Q}}(H_{\tilde{M}}(a),\epsilon'{\cal T})\tau_{\tilde{Q}}^{\tilde{G}}(H_{\tilde{L}}(a)-\epsilon'T[\tilde{Q}])\,da.$$
Si $\omega$ est non trivial sur $A_{\tilde{M}}(F)_{c}$, posons $v^T_{\tilde{M}}(x,y,\tilde{Q})=0$. Sinon, posons
$$v^T_{\tilde{M}}(x,y,\tilde{Q})=C\int_{{\cal A}_{A_{\tilde{M}},F}/{\cal A}_{A_{\tilde{G}},F}}\omega(H)\Gamma_{\tilde{M}}^{\tilde{G}}(H,{\cal Y}(x,y,T))\Gamma_{\tilde{M}}^{\tilde{Q}}(H,\epsilon'{\cal T})\tau_{\tilde{Q}}^{\tilde{G}}(H-\epsilon'T[\tilde{Q}])\,dH.$$
L'\'egalit\'e 1.3(6) nous permet d'\'ecrire
$$u^T_{\tilde{M}}(x,y)=\sum_{\tilde{Q}\in {\cal F}(\tilde{M})}u^T_{\tilde{M}}(x,y;\tilde{Q}),$$
$$v^T_{\tilde{M}}(x,y)=\sum_{\tilde{Q}\in {\cal F}(\tilde{M})}v^T_{\tilde{M}}(x,y;\tilde{Q}).$$
 On peut fixer $\tilde{Q}$ et majorer $\vert u^T_{\tilde{M}}(x,y,\tilde{Q})-v^T_{\tilde{M}}(x,y,\tilde{Q})\vert $. Ecrivons $x=ulk$, $y=\bar{u}'l'k'$, avec $l,l'\in L(F)$, $u\in U_{Q}(F)$, $\bar{u}'\in U_{\bar{Q}}(F)$, $k,k'\in K$. Soit $a\in A_{\tilde{M}}(F)$ tel que 
 
 (1) $\Gamma_{\tilde{M}}^{\tilde{Q}}(H_{\tilde{M}}(a),\epsilon'{\cal T})\tau_{\tilde{Q}}^{\tilde{G}}(H_{\tilde{M}}(a)-\epsilon'T[\tilde{Q}])=1$.
 
 On a $x^{-1}ay=k^{-1}hl^{-1}al'h'k'$, avec $h'=(a\bar{u}'l')^{-1}ua\bar{u}'l'$, $h=l^{-1}a\bar{u}'a^{-1}l$.  Les \'el\'ements $u$, $\bar{u}'$, $h'$ et $h$ sont unipotents, \'ecrivons  $u=exp(Y)$, $\bar{u}'=exp(Y')$, $h=exp(X)$, $h'=exp(X')$. Munissons $\mathfrak{g}(F)$ d'une norme $\vert .\vert $. Parce qu'il s'agit d'\'el\'ements unipotents,   on a des majorations $\vert Y\vert <<\vert \vert u\vert \vert ^D$, $ \vert Y'\vert <<\vert \vert \bar{u}'\vert \vert ^D$ pour un entier $D$ convenable.
 Fixons $\epsilon$ et supposons $\vert \vert x\vert \vert ,\vert \vert y\vert \vert \leq e^{\epsilon\vert T\vert }$. On a essentiellement la m\^eme majoration pour $l$, $u$, $l'$, $\bar{u}'$. On a donc aussi
 $\vert Y\vert ,\vert Y'\vert <<e^{D\epsilon\vert T\vert }$.  L'\'el\'ement $X$ se d\'eduit de $Y'$ par conjugaison par $l^{-1}a$. Dans la suite interviendront des r\'eels $c_{1}, c_{2}$ etc... que, pour simplifier, on n'introduira pas. A chaque fois, cela signifie qu'il existe un r\'eel $c_{1}, c_{2}...>0$, ind\'ependant des donn\'ees, tel que la relation soit v\'erifi\'ee.  Ainsi, la relation (1) entra\^{\i}ne
 $$(2) \qquad \vert \alpha(a)\vert >e^{c_{1}\epsilon'\vert T\vert }$$
 pour toute racine $\alpha$ de $A_{\tilde{M}}$ agissant dans $\mathfrak{u}_{Q} $. Ici $a$ agit dans $\mathfrak{u}_{\bar{Q}}$ et les racines y sont inf\'erieures \`a $e^{-c_{1}\epsilon'\vert T\vert }$. La conjugaison par $l^{-1}$ est de norme born\'ee par $e^{D'\epsilon\vert T\vert }$ pour un entier convenable. On en d\'eduit  des majorations
 $$(3)(a)\qquad \vert X\vert <<e^{(c_{2}\epsilon-c_{3}\epsilon')\vert T\vert }.$$
 De m\^eme,
  $$(3)(b)\qquad \vert Y\vert <<e^{(c_{2}\epsilon-c_{3}\epsilon')\vert T\vert }.$$
  
  Supposons d'abord $F$ non archim\'edien. On suppose $c_{2}\epsilon-c_{3}\epsilon'<0$. On peut toujours supposer $\vert T\vert $ assez grand (la majoration de l'\'enonc\'e \'etant triviale pour $\vert T\vert $ born\'e). Alors les relations ci-dessus entra\^{\i}nent que $h,h'\in K$. On a donc $\tilde{\kappa}^T(x^{-1}ay)=\tilde{\kappa}^T(l^{-1}al')$. Fixons un \'el\'ement $\tilde{P}\in {\cal P}(\tilde{M})$ tel que $\tilde{P}\subset \tilde{Q}$ et un \'el\'ement $w\in W^{\tilde{G}}$, relev\'e en un \'el\'ement de $K$, tel que $\tilde{P}_{0}\subset w(\tilde{P})$. Posons $\underline{\tilde{M}}=w(\tilde{M})$, $\underline{\tilde{L}}=w(\tilde{L})$, $\underline{\tilde{Q}}=w(\tilde{Q})$, $\underline{l}=wlw^{-1}$, $\underline{a}=waw^{-1}$, $\underline{l}'=wl'w^{-1}$. On a $\tilde{\kappa}^T(l^{-1}al')=\tilde{\kappa}^T(\underline{l}^{-1}\underline{a}\underline{l}')$. Fixons $u,u'\in K\cap \underline{L}(F)$ et $m\in M_{0}(F)^{\geq,\underline{Q}}$ tels que $\underline{l}^{-1}\underline{a}\underline{l}'=umu'$.
  On a
   
 (4)  si $\epsilon<c_{4}\epsilon'$ et $\vert T\vert $ est assez grand, alors $m\in M_{0}(F)^{\geq}$.
 
 Il s'agit de prouver que $\alpha(m)\geq1$ pour toute racine $\alpha$ de $A_{0}$ intervenant dans $\mathfrak{u}_{\underline{Q}}$. Il suffit de prouver que, pour tout $N\in \mathfrak{u}_{\underline{Q}}(F)$, on a $\vert ad_{m}(N)\vert \geq\vert N\vert $. En rempla\c{c}ant $N$ par $ad_{u'(\underline{l}')^{-1}}(N)$, il suffit que, pour tout $N\in \mathfrak{u}_{\underline{Q}}(F)$, on ait
 $$\vert ad_{u^{-1}\underline{l}^{-1}\underline{a}}(N)\vert \geq \vert ad_{u'(\underline{l}')^{-1}}(N)\vert .$$
  On a
 $$\vert ad_{u'(\underline{l}')^{-1}}(N)\vert<<\vert \vert l'\vert \vert^D \vert N\vert$$
 pour un entier $D$ convenable, d'o\`u
 $$\vert ad_{u'(\underline{l}')^{-1}}(N)\vert<<e^{c_{5}\epsilon\vert T\vert }\vert N\vert .$$
  Posons $N'=ad_{u^{-1}\underline{l}^{-1}\underline{a}}(N)$. On a $ad_{\underline{a}}(N)=ad_{\underline{l}u}(N')$ d'o\`u, par le m\^eme raisonnement
 $$\vert ad_{\underline{a}}(N)\vert <<e^{c_{6}\epsilon\vert T\vert }\vert N'\vert .$$
 Mais pour $\alpha$ intervenant dans $\mathfrak{u}_{\underline{Q}}$, $\alpha(\underline{a})$ est minor\'e par la relation (2). On en d\'eduit 
 $$\vert ad_{\underline{a}}(N)\vert >>e^{c_{1}\epsilon'\vert T\vert }\vert N\vert .$$
 Toutes ces relations entra\^{\i}nent
  $$\vert ad_{u'(\underline{l}')^{-1}}(N)\vert<<e^{((c_{5}+c_{6})\epsilon-c_{1}\epsilon')\vert T\vert }\vert N'\vert .$$
  En supposant $(c_{5}+c_{6})\epsilon-c_{1}\epsilon'<0$, cela entra\^{\i}ne la majoration cherch\'ee. Cela prouve (4).
  
  On a
  
  (5)  si $\epsilon<\epsilon'$, alors 
  $$\vert H_{0}(m)^{\underline{\tilde{L}}}\vert \leq c_{7}\epsilon'\vert T\vert.$$ 
  
  On fixe un ensemble fini ${\bf b}\subset {\cal A}_{\tilde{L},F}$ tel que ${\cal A}_{\tilde{L},F}\subset {\cal A}_{A_{\tilde{L}},F}+{\bf b}$. On peut choisir un \'el\'ement $b\in A_{\tilde{L}}(F)$ tel que $H_{\tilde{L}}(\underline{a})-H_{\tilde{L}}(b)\in {\bf b}$. Posons $\underline{a}'=\underline{a}b^{-1}$. La relation (1) entra\^{\i}ne $\Gamma_{\underline{\tilde{M}}}^{\underline{\tilde{Q}}}(H_{\underline{\tilde{M}}}(\underline{a}'),\epsilon'{\cal T})=1$. Puisque de plus $H_{\underline{\tilde{L}}}(\underline{a}')\in {\bf b}$, on a
  $$\vert H_{\underline{\tilde{M}}}(\underline{a}')\vert \leq c_{8}\epsilon'\vert T\vert .$$
  On en d\'eduit  une majoration
  $$\vert \vert a'\vert \vert <<e^{c_{9}\epsilon'\vert T\vert }.$$
  Posons $m'=mb^{-1}$. On a $H_{0}(m')^{\underline{\tilde{L}}}=H_{0}(m)^{\underline{\tilde{L}}}$. D'o\`u  une majoration
  $$\vert H_{0}(m)^{\underline{\tilde{L}}}\vert \leq c_{10}log(\vert \vert m'\vert \vert ).$$
  Mais $m'=\underline{l}^{-1}\underline{a}'\underline{l}'$. Par un calcul maintenant familier, on en d\'eduit 
  $$log(\vert \vert m'\vert \vert )\leq1+ (c_{9}\epsilon'+c_{11}\epsilon)\vert T\vert .$$
  La relation (5) r\'esulte des deux majorations pr\'ec\'edentes. 
  
  On suppose $\epsilon$, $\epsilon'$ et $T$  tels que (4) et (5) soient v\'erifi\'es. Montrons que

(6) si $\epsilon'$ est assez petit, on a l'\'egalit\'e $\phi_{\tilde{P}_{0}}^{\tilde{G}}(H_{0}(m)-T)=\phi_{\underline{\tilde{Q}}}^{\tilde{G}}(H_{\underline{\tilde{L}}}(m)-T)$.

\noindent Notons simplement $H=H_{\tilde{P}_{0}}(m)$. L'\'egalit\'e $\phi_{\tilde{P}_{0}}^{\tilde{G}}(H_{0}(m)-T)=1$ \'equivaut \`a 

(7) $<\varpi_{\tilde{\alpha}},H-T>\leq0$ pour tout $\tilde{\alpha}\in \Delta_{\tilde{P}_{0}}$. 

\noindent L'\'egalit\'e $\phi_{\underline{\tilde{Q}}}^{\tilde{G}}(H_{\underline{\tilde{L}}}(m)-T)=1$ \'equivaut \`a 

(8) $<\varpi_{\tilde{\alpha}},H-T>\leq0$ pour tout $\tilde{\alpha}\in \Delta_{\tilde{P}_{0}}-\Delta_{\tilde{P}_{0}}^{\underline{\tilde{Q}}}$. 

\noindent Evidemment (7) entra\^{\i}ne (8). Inversement, supposons (8) v\'erifi\'ee. Soit $\tilde{\alpha}\in \Delta_{\tilde{P}_{0}}^{\underline{\tilde{Q}}}$. On a
$$<\varpi_{\tilde{\alpha}},H-T>=<\varpi_{\tilde{\alpha}},H^{\underline{\tilde{L}}}-T^{\underline{\tilde{L}}}>+<\varpi_{\tilde{\alpha}},H_{\underline{\tilde{L}}}-T_{\underline{\tilde{L}}}>.$$
On va prouver que ces deux termes sont n\'egatifs pourvu que $\epsilon'$ soit assez petit. Cela prouvera que (8) entra\^{\i}ne (7) et ach\`evera la preuve de (6). Le terme $H^{\underline{\tilde{L}}}$ est une projection de $H_{0}(m)^{\underline{\tilde{L}}}$. D'apr\`es (5),  on a donc une majoration
$$\vert <\varpi_{\tilde{\alpha}},H^{\underline{\tilde{L}}}>\vert \leq c_{12}\epsilon'\vert T\vert .$$
On a aussi facilement
$$<\varpi_{\tilde{\alpha}},T>\geq c_{13}\vert T\vert .$$
Il r\'esulte de ces deux in\'egalit\'es que $<\varpi_{\tilde{\alpha}},H^{\underline{\tilde{L}}}-T^{\underline{\tilde{L}}}>$ est n\'egatif si $\epsilon'$ est assez petit. D'apr\`es (7), le terme $H_{\underline{\tilde{L}}}-T_{\underline{\tilde{L}}}$ est combinaison lin\'eaire \`a coefficients n\'egatifs ou nuls de termes $\check{\tilde{\beta}}_{\underline{\tilde{L}}}$ pour $\tilde{\beta}\in \Delta_{\tilde{P}_{0}}-\Delta_{\tilde{P}_{0}}^{\underline{\tilde{Q}}}$. Il reste \`a voir que, pour un tel $\tilde{\beta}$, on a $<\varpi_{\tilde{\alpha}},\check{\tilde{\beta}}_{\underline{\tilde{L}}}>\geq0$. On a 
$$<\varpi_{\tilde{\alpha}},\check{\tilde{\beta}}_{\underline{\tilde{L}}}>=<\varpi_{\tilde{\alpha}},\check{\tilde{\beta}}>-<\varpi_{\tilde{\alpha}},\check{\tilde{\beta}}^{\underline{\tilde{L}}}>.$$
Le premier terme est nul puisque $\tilde{\alpha}\in \Delta_{\tilde{P}_{0}}^{\underline{\tilde{Q}}}$ et $\tilde{\beta}\not\in \Delta_{\tilde{P}_{0}}^{\underline{\tilde{Q}}}$. Parce que $\Delta_{\tilde{P}_{0}}$ est une base aigu\" e, $\check{\tilde{\beta}}^{\underline{\tilde{L}}}$ appartient \`a la chambre n\'egative ferm\'ee associ\'ee \`a $\Delta_{\tilde{P}_{0}}^{\underline{\tilde{Q}}}$, a fortiori  est combinaison lin\'eaire \`a coefficients n\'egatifs ou nuls de $\check{\tilde{\alpha}}'$ pour $\tilde{\alpha}'\in \Delta_{\tilde{P}_{0}}^{\underline{\tilde{Q}}}$. Donc $<\varpi_{\tilde{\alpha}},\check{\tilde{\beta}}^{\underline{\tilde{L}}}>\leq0$. Cela ach\`eve la preuve de (6).

 On a
  
  (9) si $\epsilon'$ est assez petit, on a l'\'egalit\'e $\tilde{\kappa}^T(x^{-1}ay)=\phi_{\tilde{Q}}^{\tilde{G}}(H_{\tilde{L}}(a)-H_{\tilde{Q}}(x)+H_{\tilde{\bar{Q}}}(y)-T[\tilde{Q}])$.

On a vu que $\tilde{\kappa}^T(x^{-1}ay)=\tilde{\kappa}^T(\underline{l}^{-1}\underline{a}\underline{l}')=\tilde{\kappa}^T(m)$. Puisque $m\in M_{0}(F)^{\geq}$ d'apr\`es (4), on a $\tilde{\kappa}^T(m)=\phi_{\tilde{P}_{0}}^{\tilde{G}}(H_{0}(m)-T)$.  Gr\^ace \`a (5), on a donc $\tilde{\kappa}^T(x^{-1}ay)=
\phi_{\underline{\tilde{Q}}}^{\tilde{G}}(H_{\underline{\tilde{L}}}(m)-T)$ pourvu que $\epsilon'$ soit assez petit. Par conjugaison par $w$, cela entra\^{\i}ne $\tilde{\kappa}^T(x^{-1}ay)=\phi_{\tilde{Q}}^{\tilde{G}}(w^{-1}H_{\underline{\tilde{L}}}(m)-T[\tilde{Q}])$.  Par d\'efinition de $m$, on a l'\'egalit\'e 
$$H_{\underline{\tilde{L}}}(m)=-H_{\underline{\tilde{L}}}(\underline{l})+H_{\underline{\tilde{L}}}(\underline{a})+H_{\underline{\tilde{L}}}(\underline{l}').$$
On a aussi $w^{-1}H_{\underline{\tilde{L}}}(\underline{l})=H_{\tilde{L}}(l)$ etc... D'o\`u
$$w^{-1}H_{\underline{\tilde{L}}}(m)=-H_{\tilde{L}}(l)+H_{\tilde{L}}(a)+H_{\tilde{L}}(l').$$
Il r\'esulte des d\'efinitions que $H_{\tilde{L}}(l)=H_{\tilde{Q}}(x)$ et $H_{\tilde{L}}(l')=H_{\tilde{\bar{Q}}}(y)$. On obtient alors l'\'egalit\'e de (9).

Supposons $\epsilon'$ assez petit pour que (9) soit v\'erifi\'ee. Alors l'int\'egrale d\'efinissant $u^T_{\tilde{M}}(x,y;\tilde{Q})$ se factorise en 
$$(10)\qquad u^T_{\tilde{M}}(x,y;\tilde{Q})=C(\omega) \int_{{\cal A}_{A_{\tilde{M}},F}/{\cal A}_{A_{\tilde{G}},F}}\phi_{\tilde{Q}}^{\tilde{G}}(H-H_{\tilde{Q}}(x)+H_{\tilde{\bar{Q}}}(y)-T[\tilde{Q}])$$
$$\Gamma_{\tilde{M}}^{\tilde{Q}}(H,\epsilon'{\cal T})\tau_{\tilde{Q}}^{\tilde{G}}(H-\epsilon'T[\tilde{Q}])\,dH,$$
o\`u
$$C(\omega)=\int_{A_{\tilde{M}}(F)_{c}A_{\tilde{G}}(F)/A_{\tilde{G}}(F)}\omega(a)\,da.$$
Si $\omega$ n'est pas trivial sur $A_{\tilde{M}}(F)_{c}$, $C(\omega)$ est nul. Alors $u^T_{\tilde{M}}(x,y;\tilde{Q})=v_{\tilde{M}}^T(x,y;\tilde{Q})$ par d\'efinition de ce dernier terme. Supposons maintenant $\omega$ trivial sur $A_{\tilde{M}}(F)_{c}$.  Alors $C(\omega)=C$.  Fixons $H\in {\cal A}_{A_{\tilde{M}},F}$ tel que
$$(11)\qquad \Gamma_{\tilde{M}}^{\tilde{Q}}(H,\epsilon'{\cal T})\tau_{\tilde{Q}}^{\tilde{G}}(H-\epsilon'T[\tilde{Q}])=1.$$
On va prouver

(12) si $\epsilon'$ et $\epsilon<\epsilon'$ sont assez petits, on a l'\'egalit\'e
$$\Gamma_{\tilde{M}}^{\tilde{G}}(H,{\cal Y}(x,y,T))=\phi_{\tilde{Q}}^{\tilde{G}}(H-H_{\tilde{Q}}(x)+H_{\tilde{\bar{Q}}}(y)-T[\tilde{Q}]).$$

Fixons un sous-espace parabolique $\tilde{P}\in {\cal P}(\tilde{M})$ tel que $\tilde{P}\subset \tilde{Q}$ et que $H$ appartienne \`a la cl\^oture de la chambre positive associ\'ee \`a l'espace parabolique $\tilde{P}\cap \tilde{L}$ de $\tilde{L}$. On a alors l'\'egalit\'e $\Gamma_{\tilde{M}}^{\tilde{Q}}(H,\epsilon'{\cal T})=\phi_{\tilde{P}}^{\tilde{Q}}(H-\epsilon'T[\tilde{P}])$. L'\'egalit\'e (11) entra\^{\i}ne alors que $H$ appartient \`a la cl\^oture de la chambre positive associ\'ee \`a l'espace parabolique $\tilde{P}$. Par des arguments d\'ej\`a vus, on a des majorations
$$\vert H_{\tilde{P}'}(x)\vert, \vert H_{\tilde{P}'}(y)\vert \leq c_{14}\epsilon'\vert T\vert $$
pour tout $\tilde{P}'\in {\cal P}(\tilde{M})$. Si $\epsilon'$ est assez petit, il en r\'esulte que ${\cal Y}(x,y,T)$ est une famille positive. Donc  $H'\mapsto\Gamma_{\tilde{M}}^{\tilde{G}}(H',{\cal Y}(x,y,T))$ est la fonction caract\'eristique de l'ensemble des $H'$ tels que $(H')^{\tilde{G}}$ appartient \`a l'enveloppe convexe des $Y(x,y,T,\tilde{P}')^{\tilde{G}}$. Pour $H$ dans la cl\^oture de la chambre positive associ\'ee \`a $\tilde{P}$, on a simplement 
$$\Gamma_{\tilde{M}}^{\tilde{G}}(H,{\cal Y}(x,y,T))=\phi_{\tilde{P}}^{\tilde{G}}(H-H_{\tilde{P}}(x)+H_{\tilde{\bar{P}}}(y)-T[\tilde{P}]).$$
Posons $X=H-H_{\tilde{P}}(x)+H_{\tilde{\bar{P}}}(y)$. De la m\^eme fa\c{c}on que l'on a prouv\'e (5) et (6), on montre que si $\epsilon<\epsilon'$, on a $\vert X^{\tilde{L}}\vert \leq c_{15}\epsilon'\vert T\vert $, puis que, si $\epsilon'$ est assez petit, on a l'\'egalit\'e
$$\phi_{\tilde{P}}^{\tilde{G}}(H-H_{\tilde{P}}(x)+H_{\tilde{\bar{P}}}(y)-T[\tilde{P}])=\phi_{\tilde{Q}}^{\tilde{G}}(H-H_{\tilde{Q}}(x)+H_{\tilde{\bar{Q}}}(y)-T[\tilde{Q}]).$$
Cela prouve (12).

L'\'egalit\'e (10) jointe \`a (12) prouve l'\'egalit\'e $u_{\tilde{M}}^T(x,y;\tilde{Q})=v_{\tilde{M}}^T(x,y;\tilde{Q})$. On a d\^u supposer $\epsilon'$ et $\epsilon/\epsilon'$ assez petits. Si $\epsilon$ lui-m\^eme est assez petit, on peut choisir le param\`etre auxiliaire $\epsilon'$ satisfaisant ces conditions. Cela prouve le lemme. Cela prouve m\^eme plus, \`a savoir l'\'egalit\'e $u_{\tilde{M}}^T(x,y;\tilde{Q})=v_{\tilde{M}}^T(x,y;\tilde{Q})$. 

Rappelons que l'on avait suppos\'e $F$ non archim\'edien. Remarquons que cette hypoth\`ese n'a \'et\'e utilis\'ee que pour affirmer que $exp(X)$ et $exp(Y)$ appartenaient \`a $K$. Tout le reste est valable sans hypoth\`ese sur $F$. Ainsi, supposons maintenant $F$ archim\'edien et  d\'efinissons ${\bf u}^T_{\tilde{M}}(x,y,\tilde{Q})$ comme l'int\'egrale analogue \`a $u_{\tilde{M}}^T(x,y;\tilde{Q})$ o\`u l'on remplace $\tilde{\kappa}^T(x^{-1}ay)$ par $\tilde{\kappa}^T(l^{-1}al')$. On a pour cette int\'egrale les m\^emes r\'esultats d\'emontr\'es pr\'ec\'edemment pour $u_{\tilde{M}}^T(x,y;\tilde{Q})$.
  Revenons \`a l'\'egalit\'e $x^{-1}ay=k^{-1}exp(X)l^{-1}al'exp(Y)k'$, qui entra\^{\i}ne $\tilde{\kappa}^T(x^{-1}ay)=\tilde{\kappa}^T(exp(X)l^{-1}al'exp(Y))$. Fixons un \'el\'ement $m'\in M_{0}(F)^{\geq}$ tel que $exp(X)l^{-1}al'exp(Y)\in Km'K$ et d\'efinissons $m$ comme dans le cas non archim\'edien.  
On a encore (4): $m\in M_{0}(F)^{\geq}$ (sous hypoth\`eses sur $\epsilon$ et $\epsilon'$). Il r\'esulte alors de (3)(a) et (3)(b) et du lemme 5.2 de [A1] que l'on a une majoration
$$\vert H_{0}(m)-H_{0}(m')\vert \leq c_{15}exp^{-c_{16}\vert T\vert }.$$
Posons $\eta_{+}=1+exp^{-c_{16}\vert T\vert /2}$, $\eta_{-}=1-exp^{-c_{16}\vert T\vert /2}$.  De l'in\'egalit\'e ci-dessus r\'esulte les majorations
$$\tilde{\kappa}^{\eta_{-}T}(m)\leq\left\lbrace\begin{array}{c}\tilde{\kappa}^T(m)\\\tilde{\kappa}^T(m')\\ \end{array}\right\rbrace \leq\tilde{\kappa}^{\eta_{+}T}(m),$$
puis
$$\vert \tilde{\kappa}^T(m)-\tilde{\kappa}^T(m')\vert \leq\tilde{\kappa}^{\eta_{+}T}(m)-\tilde{\kappa}^{\eta_{-}T}(m).$$
Puisque $\tilde{\kappa}^T(m')=\tilde{\kappa}^T(x^{-1}ay)$, on en d\'eduit
$$\vert {\bf u}_{\tilde{M}}^T(x,y;\tilde{Q})-u_{\tilde{M}}^T(x,y;\tilde{Q})\vert \leq\int_{A_{\tilde{M}}(F)/A_{\tilde{G},F}}(\tilde{\kappa}^{\eta_{+}T}(l^{-1}al')-\tilde{\kappa}^{\eta_{-}T}(l^{-1}al'))$$
$$\Gamma_{\tilde{M}}^{\tilde{Q}}(H_{\tilde{M}}(a),\epsilon'{\cal T})\tau_{\tilde{Q}}^{\tilde{G}}(H_{\tilde{L}}(a)-\epsilon'T[\tilde{Q}]\,da.$$
Le membre de droite est la diff\'erence entre une expression analogue \`a ${\bf u}_{\tilde{M}}^T(x,y;\tilde{Q})$ o\`u l'on remplace $\omega$ par $1$, $T$ par $\eta_{+}T$ et $\epsilon'$ par $\eta_{+}^{-1}\epsilon'$, et l'expression similaire o\`u $\eta_{-}$ remplace $\eta_{+}$. Ces expressions sont calcul\'ees par (10). On obtient
$$\vert {\bf u}_{\tilde{M}}^T(x,y;\tilde{Q})-u_{\tilde{M}}^T(x,y;\tilde{Q})\vert \leq $$
$$ \int_{{\cal A}_{A_{\tilde{M}},F}/{\cal A}_{A_{\tilde{G}},F}}(\phi_{\tilde{Q}}^{\tilde{G}}(H-H_{\tilde{Q}}(x)+H_{\tilde{\bar{Q}}}(y)-\eta_{+}T[\tilde{Q}])-\phi_{\tilde{Q}}^{\tilde{G}}(H-H_{\tilde{Q}}(x)+H_{\tilde{\bar{Q}}}(y)-\eta_{-}T[\tilde{Q}]))$$
$$\Gamma_{\tilde{M}}^{\tilde{Q}}(H,\epsilon'{\cal T})\tau_{\tilde{Q}}^{\tilde{G}}(H-\epsilon'T[\tilde{Q}])\,dH.$$
Le membre de droite est la diff\'erence entre les volumes de deux poly\`edres dont l'un est inclus dans l'autre. Leurs sommets ont une norme essentiellement born\'ee par $\vert T\vert $ et tout point de la fronti\`ere du petit poly\`edre est \`a distance au plus $c_{17}(\eta_{+}-\eta_{-})\vert T\vert $ de la fronti\`ere du plus grand. La diff\'erence de leur volume est donc essentiellement born\'ee par $(\eta_{+}-\eta_{-})\vert T\vert ^D$ pour un entier $D$ convenable. Donc par $e^{-c_{16}\vert T\vert /2}\vert T\vert ^D$. On conclut que, pour tout r\'eel $r$,
$$\vert {\bf u}_{\tilde{M}}^T(x,y;\tilde{Q})-u_{\tilde{M}}^T(x,y;\tilde{Q})\vert << \vert T\vert ^{-r}.$$
Comme on l'a dit, la preuve du cas non-archim\'edien montre que ${\bf u}_{\tilde{M}}^T(x,y;\tilde{Q})=v_{\tilde{M}}^T(x,y;\tilde{Q})$. Cela prouve le lemme. $\square$

\bigskip
\subsection{Preuve de la proposition 4.3}
Pour prouver que $J_{v,\tilde{M},\tilde{S}}^T(\omega,f_{1},f_{2})$ est absolument int\'egrable, on consid\`ere $T$ comme fix\'e. L'expression $J_{v,\tilde{M},\tilde{S}}^T(\omega,f_{1},f_{2})$ est une int\'egrale en $\gamma$. Elle est \`a support compact. Au voisinage des \'el\'ements singuliers, elle est absolument int\'egrable d'apr\`es le lemme 4.4. En utilisant la relation 4.3(3), on v\'erifie que la fonction $v_{\tilde{M}}^T(x,y)$ reste born\'ee (pour $T$ fix\'e) si $\gamma$ reste hors d'un voisinage des \'el\'ements singuliers.  Alors l'int\'egrale hors d'un tel voisinage est essentiellement born\'ee par
$$(1) \qquad \int_{\tilde{S}(F)/(1-\theta)(S(F))}
 \int_{S^{\theta}(F)\backslash G(F)}\int_{A_{\tilde{M}}(F)\backslash G(F)}\vert f_{1}(x^{-1}\gamma x)f_{2}(y^{-1}\gamma y)\vert   \,dy\,dx\,D^{\tilde{G}}(\gamma)\,d\gamma,$$
 ou encore par
  $$\int_{\tilde{S}(F)/(1-\theta)(S(F))}I_{\tilde{G}}(\gamma,\vert f_{1}\vert )I_{\tilde{G}}(\gamma,\vert f_{2}\vert )d\gamma,$$
 qui est convergente d'apr\`es 4.2(1).

D\'emontrons maintenant la majoration de l'\'enonc\'e de la proposition 4.3. On fixe $\epsilon$ comme dans le lemme 4.5. On fixe ensuite $\epsilon'>0$ tel que la relation $D^{\tilde{G}}(\gamma)\geq e^{-\epsilon'\vert T\vert }$ et la relation 4.3(3):
$$\vert \vert x\vert \vert ,\vert \vert y\vert \vert \leq cD^{\tilde{G}}(\gamma)^{-k}$$
entra\^{\i}nent $\vert \vert x\vert \vert ,\vert \vert y\vert \vert \leq e^{\epsilon\vert T\vert }$, du moins si $T$ est grand (on peut supposer $T$ grand pour d\'emontrer notre majoration). On d\'ecompose $J_{\tilde{M},\tilde{S}}^T(\omega,f_{1},f_{2})-J_{v,\tilde{M},\tilde{S}}^T(\omega,f_{1},f_{2})$ en somme de deux int\'egrales. Dans la premi\`ere, on glisse le terme $\chi_{\epsilon'}^T(\gamma)$. Dans la seconde, on glisse le terme $1-\chi_{\epsilon'}^T(\gamma)$. La premi\`ere int\'egrale est major\'ee par le lemme 4.4: pour tout r\'eel $r$, elle est essentiellement born\'ee par $\vert T\vert ^{-r}$. Dans la seconde, gr\^ace \`a 4.3(3) et nos choix de $\epsilon'$ et $\epsilon$, on a la majoration du lemme 4.5. Pour tout $r$, cette int\'egrale est donc essentiellement major\'ee par le produit de $\vert T\vert ^{-r}$ et de l'int\'egrale (1). Celle-ci \'etant convergente, cela prouve la majoration cherch\'ee. $\square$

\bigskip

\subsection{Terme constant de $v_{\tilde{M}}^T(x,y)$}
Soient $x,y\in G(F)$. Pour $\tilde{P}\in {\cal P}(\tilde{M})$ et $\Lambda\in i{\cal A}_{\tilde{M}}^*$, posons
$$v(x,y;\Lambda,\tilde{P})=e^{<\Lambda, H_{\tilde{\bar{P}}}(y)-H_{\tilde{P}}(x)>}.$$
La famille $(v(x,y;\Lambda,\tilde{P}))_{\tilde{P}\in {\cal P}(\tilde{M})}$ est une $(\tilde{G},\tilde{M})$-famille, dont on d\'eduit une fonction $v^{\tilde{G}}_{\tilde{M}}(x,y;\Lambda)$. On pose $v^{\tilde{G}}_{\tilde{M}}(x,y)=v^{\tilde{G}}_{\tilde{M}}(x,y;0)$.
\ass{Lemme}{(i) La fonction $T\mapsto f(T)=v_{\tilde{M}}^T(x,y)$ appartient \`a $PolExp$. Plus pr\'ecis\'ement, elle appartient \`a un espace $PolExp_{\Xi,N}$ si $F$ est archim\'edien, resp. $PolExp_{\boldsymbol{\Xi},N}$ si $F$ est non-archim\'edien, o\`u $\Xi$ et $N$, resp. $\boldsymbol{\Xi}$ et $N$, sont ind\'ependants de $x$ et $y$. 

(ii) Supposons $F$ archim\'edien. On a 
$$c_{0}(f)=\left\lbrace\begin{array}{cc}(-1)^{a_{\tilde{M}}-a_{\tilde{G}}}v^{\tilde{G}}_{\tilde{M}}(x,y),&\text{ si }\omega\text{ est trivial sur }A_{\tilde{M}}(F),\\0, &\text{ sinon.}\\ \end{array}\right.$$

(iii) Supposons $F$ non-archim\'edien. Soit  ${\cal R}\subset {\cal A}_{M_{0},F}\otimes_{{\mathbb Z}}{\mathbb Q}$ un r\'eseau.   Alors il existe un entier $k_{1}>0$ et un r\'eel $c>0$ ne d\'ependant tous deux que de ${\cal R}$ de sorte que, pour tout entier $k\geq1$, on ait la majoration
 $$\vert c_{\frac{1}{k}{\cal R},0}(f)\vert\leq ck^{-1}(log(\vert \vert x\vert \vert )+log(\vert \vert y\vert \vert )) ^{a_{\tilde{M}}-a_{\tilde{G}}}$$
 si $\omega$ est non trivial sur $A_{\tilde{M}}(F)$, respectivement
  $$\vert c_{\frac{1}{k}{\cal R},0}(f)- (-1)^{a_{\tilde{M}}-a_{\tilde{G}}}v^{\tilde{G}}_{\tilde{M}}(x,y)\vert\leq ck^{-1}(log(\vert \vert x\vert \vert )+log(\vert \vert y\vert \vert )) ^{a_{\tilde{M}}-a_{\tilde{G}}}$$
  si $\omega$ est trivial sur $A_{\tilde{M}}(F)$.}
 
 Preuve. On traite le cas o\`u $F$ est non-archim\'edien, le cas archim\'edien \'etant plus simple. Si $\omega$ n'est pas trivial sur $A_{\tilde{M}}(F)_{c}$, $v_{\tilde{M}}^T(x,y)=0$ et les assertions sont \'evidentes. Supposons $\omega$ trivial sur $A_{\tilde{M}}(F)_{c}$. On a d\'eduit de $\omega$ un caract\`ere de ${\cal A}_{A_{\tilde{M}},F}$, encore not\'e $\omega$. Il y a un unique $\Lambda_{\omega}\in i{\cal A}_{\tilde{M}}^*/{\cal A}_{A_{\tilde{M}},F}^{\vee}$ tel que $\omega(H)=e^{<\Lambda_{\omega},H>}$ pour tout $H\in {\cal A}_{A_{\tilde{M}},F}$. Notons que, puisque $\omega$ est trivial sur $A_{\tilde{G}}(F)$, la projection $(\Lambda_{\omega})_{\tilde{G}}$ appartient \`a $i{\cal A}_{A_{\tilde{G}},F}^*$. On a
 $$v_{\tilde{M}}^T(x,y)=C\sum_{{\cal A}_{A_{\tilde{M}},F}/{\cal A}_{A_{\tilde{G}},F}}\Gamma_{\tilde{M}}^{\tilde{G}}(H,{\cal Y}(x,y,T))e^{<\Lambda_{\omega},H>}.$$
 On peut d\'ecomposer la somme en une somme sur $X\in {\cal A}_{\tilde{G},F}/{\cal A}_{A_{\tilde{G}},F}$ suivie d'une somme sur $H\in {\cal A}_{A_{\tilde{M}},F}\cap {\cal A}_{\tilde{M},F}^{\tilde{G}}(X)$. Par inversion de Fourier, cette derni\`ere somme peut \^etre remplac\'ee par la somme sur $H\in   {\cal A}_{\tilde{M},F}^{\tilde{G}}(X)$, suivie de
 $$[{\cal A}_{\tilde{M},F}:{\cal A}_{A_{\tilde{M}},F}]^{-1}\sum_{\nu\in i{\cal A}_{A_{\tilde{M}},F}^{\vee}/i{\cal A}_{\tilde{M},F}^{\vee}}e^{<\nu,H>}.$$
 Remarquons que $C[{\cal A}_{\tilde{M},F}:{\cal A}_{A_{\tilde{M}},F}]^{-1}=mes(A_{\tilde{G}}(F)_{c})^{-1}mes(i{\cal A}_{\tilde{M},F}^*)^{-1}$.
  Tout ceci est absolument convergent  puisque la fonction  $H\mapsto \Gamma_{\tilde{M}}^{\tilde{G}}(H,{\cal Y}(x,y,T))$ est \`a support compact. On obtient
 $$v_{\tilde{M}}^T(x,y)= mes(A_{\tilde{G}}(F)_{c})^{-1}mes(i{\cal A}_{\tilde{M},F}^*)^{-1} \sum_{\nu\in  i{\cal A}_{A_{\tilde{M}},F}^{\vee}/i{\cal A}_{\tilde{M},F}^{\vee} }\sum_{X\in {\cal A}_{\tilde{G},F}/{\cal A}_{A_{\tilde{G}},F}}I^T(x,y,X;\Lambda_{\omega}+\nu),$$
 o\`u, pour $\Lambda\in {\cal A}_{\tilde{M},{\mathbb C}}^*$, on a pos\'e
 $$I^T(x,y,X;\Lambda)=\sum_{H\in  {\cal A}_{\tilde{M},F}^{\tilde{G}}(X)}\Gamma_{\tilde{M}}^{\tilde{G}}(H,{\cal Y}(x,y,T))e^{<\Lambda,H>}.$$
  Cette derni\`ere somme est finie donc d\'efinit une fonction holomorphe en $\Lambda$. Fixons $\tilde{P}_{1}\in {\cal P}(\tilde{M})$ et supposons $<\Lambda,\check{\tilde{\alpha}}>>0$ pour tout $\tilde{\alpha}\in \Delta_{\tilde{P}_{1}}$. La variante tordue de 1.3(7) nous dit que
$$\Gamma_{\tilde{M}}^{\tilde{G}}(H,{\cal Y}(x,y,T))=\sum_{\tilde{P}\in {\cal P}(\tilde{M})}(-1)^{s(\tilde{P},\tilde{P}_{1})}\phi_{\tilde{P},\tilde{P}_{1}}^{\tilde{G}}(H-Y(x,y,T;\tilde{P}) ),$$
o\`u $Y(x,y,T;\tilde{P})=T[\tilde{P}]-H_{\tilde{P}}(x)+H_{\tilde{\bar{P}}}(y)$.   L'hypoth\`ese sur $\Lambda$ nous permet de permuter les deux sommes:
 $$I^T(x,y,X;\Lambda)=\sum_{\tilde{P}\in {\cal P}(\tilde{M})}(-1)^{s(\tilde{P},\tilde{P}_{1})}\sum_{H\in   {\cal A}_{\tilde{M},F}^{\tilde{G}}(X)}\phi_{\tilde{P},\tilde{P}_{1}}^{\tilde{G}}(H- Y(x,y,T;\tilde{P}))$$
 $$=\sum_{\tilde{P}\in {\cal P}(\tilde{M})}(-1)^{s(\tilde{P},\tilde{P}_{1})}\epsilon_{\tilde{P},\tilde{P}_{1}}^{\tilde{G},Y(x,y,T;\tilde{P})}(X;\Lambda).$$
 La variante tordue de 1.5(2) conduit \`a l'\'egalit\'e
 $$I^T(x,y,X;\Lambda)=\sum_{\tilde{P}\in {\cal P}(\tilde{M})}\epsilon_{\tilde{P}}^{\tilde{G},Y(x,y,T;\tilde{P})}(X;\Lambda).$$
 Introduisons la $(\tilde{G},\tilde{M})$-famille $(\varphi(\Lambda;\tilde{P}))_{\tilde{P}\in {\cal P}(\tilde{M})}$ form\'ee des fonctions constantes valant $1$. Alors
 $$I^T(x,y,X;\Lambda)=\varphi_{\tilde{M}}^{\tilde{G},{\cal Y}(x,y,T)}(X;\Lambda).$$
 On a prouv\'e cette \'egalit\'e sous des hypoth\`eses portant sur $\Lambda$ mais elle se prolonge \`a tout $\Lambda$, les deux membres \'etant m\'eromorphes (et en fait holomorphes).
 Notons simplement ${\cal Y}_{x,y}={\cal Y}(x,y,0)$. Avec la d\'efinition de 1.7 (dans sa version tordue), on a ${\cal Y}(x,y,T)={\cal Y}_{x,y}(T)$. On obtient
  $$v_{\tilde{M}}^T(x,y)=mes(A_{\tilde{G}}(F)_{c})^{-1}mes(i{\cal A}_{\tilde{M},F}^*)^{-1} \sum_{\nu\in  i{\cal A}_{A_{\tilde{M}},F}^{\vee}/i{\cal A}_{\tilde{M},F}^{\vee} }\sum_{X\in {\cal A}_{\tilde{G},F}/{\cal A}_{A_{\tilde{G}},F}}\varphi_{\tilde{M}}^{\tilde{G},{\cal Y}_{x,y}(T)}(X;\Lambda_{\omega}+\nu).$$
  La variante tordue du lemme 1.7 entra\^{\i}ne l'assertion (i) de l'\'enonc\'e. Pour un r\'eseau ${\cal R}$ et un entier $k$, le coefficient $c_{\frac{1}{k}{\cal R},0}(f)$ est le produit de $mes(A_{\tilde{G}}(F)_{c})^{-1}mes(i{\cal A}_{\tilde{M},F}^*)^{-1}$ et de la somme des coefficients analogues pour les fonctions $T\mapsto \varphi_{\tilde{M}}^{\tilde{G},{\cal Y}_{x,y}(T)}(X;\Lambda_{\omega}+\nu)$, lesquels sont calcul\'es par le m\^eme lemme. Notons ${\cal N}$ l'ensemble des $\nu \in i{\cal A}_{A_{\tilde{M}},F}^{\vee}/i{\cal A}_{\tilde{M},F}^{\vee}$ tel que $\Lambda_{\omega}+\nu\in i{\cal A}_{\tilde{M},F}^{\vee}+i{\cal A}_{\tilde{G}}^*$. Pour $\nu\in {\cal N}$, fixons $\Lambda_{\omega,\nu}\in i{\cal A}_{\tilde{G}}^*$ tel que $\Lambda_{\omega}+\nu\in\Lambda_{\omega,\nu} +i{\cal A}_{\tilde{M},F}^{\vee}$. Posons
  $$(1) \qquad c(f)= mes(A_{\tilde{G}}(F)_{c})^{-1}mes(i{\cal A}_{\tilde{G},F}^*)^{-1} \sum_{\nu\in {\cal N}}\sum_{X\in {\cal A}_{\tilde{G},F}/{\cal A}_{A_{\tilde{G}},F}}e^{<\Lambda_{\omega,\nu},X>}\varphi_{\tilde{M}}^{\tilde{G}}({\cal Y}_{x,y};\Lambda_{\omega,\nu}).$$
  Notons $N(x,y)$ la norme de la $(\tilde{G},\tilde{M})$-famille $(\varphi({\cal Y}_{x,y};\Lambda,\tilde{P}))_{\tilde{P}\in {\cal P}(\tilde{M})}$. 
  La variante du lemme 1.7 implique qu'il existe un entier $k_{1}>0$ et un r\'eel $c>0$, ne d\'ependant tous deux que de ${\cal R}$, de sorte que, si $k\geq k_{1}$, on ait la majoration 
  $$(2) \qquad \vert c_{\frac{1}{k}{\cal R},0}(f)-c(f)\vert \leq cN(x,y)k^{-1}.$$
  En comparant les d\'efinitions, on voit que $\varphi({\cal Y}_{x,y};\Lambda,\tilde{P})=v(x,y;-\Lambda,\tilde{P})$ pour tout $\tilde{P}$. Parce que $\epsilon_{\tilde{P}}^{\tilde{G}}(-\Lambda)=(-1)^{a_{\tilde{M}}-a_{\tilde{G}}}\epsilon_{\tilde{P}}^{\tilde{G}}(\Lambda)$, on en d\'eduit 
  $$\varphi_{\tilde{M}}^{\tilde{G}}({\cal Y}_{x,y};\Lambda_{\omega,\nu})=(-1)^{a_{\tilde{M}}-a_{\tilde{G}}}v^{\tilde{G}}_{\tilde{M}}(x,y;-\Lambda_{\omega,\nu}).$$
   Dans la formule (1), la somme en $X$ vaut $[{\cal A}_{\tilde{G},F}:{\cal A}_{A_{\tilde{G}},F}] $ si $\Lambda_{\omega,\nu}\in i{\cal A}_{\tilde{G},F}^{\vee}$, $0$ sinon. Notons que si $\Lambda_{\omega,\nu}\in i{\cal A}_{\tilde{G},F}^{\vee}$, on a $\Lambda_{\omega}+\nu\in i{\cal A}_{\tilde{M},F}^{\vee}$, d'o\`u $\Lambda_{\omega}\in i{\cal A}_{A_{\tilde{M}},F}^{\vee}$. Cela \'equivaut \`a ce que $\omega$ soit trivial sur $A_{\tilde{M}}(F)$. Donc, si $\omega$ n'est pas trivial sur $A_{\tilde{M}}(F)$, $c(f)=0$. Supposons $\omega$ trivial sur $A_{\tilde{M}}(F)$. On peut supposer $\Lambda_{\omega}=0$. La condition $\nu\in {\cal N}$ signifie que $\nu\in( i{\cal A}_{A_{\tilde{M}},F}^{\vee}\cap i{\cal A}_{\tilde{G}}^*)/i{\cal A}_{\tilde{G},F}^{\vee}$. La condition $\Lambda_{\omega,\nu}\in i{\cal A}_{\tilde{G},F}^{\vee}$ s\'electionne le terme $\nu=0$. En remarquant que $mes(A_{\tilde{G}}(F)_{c})^{-1}mes(i{\cal A}_{\tilde{G},F}^*)^{-1}[{\cal A}_{\tilde{G},F}:{\cal A}_{A_{\tilde{G}},F}] =1$, on obtient
  $$c(f)= (-1)^{a_{\tilde{M}}-a_{\tilde{G}}}v_{\tilde{M}}(x,y).$$
  La majoration
  $$N(x,y)<<sup_{\tilde{P}\in {\cal P}(\tilde{M})}(1+\vert H_{\tilde{P}}(x)\vert +\vert H_{\tilde{\bar{P}}}(y)\vert )^{a_{\tilde{M}}-a_{\tilde{G}}}$$
  r\'esulte des d\'efinitions. Il est clair que l'on a
$$(3) \qquad 1+\vert H_{\tilde{P}}(x)\vert +\vert H_{\tilde{\bar{P}}}(y)\vert<<1+log(\vert \vert x\vert \vert )+log(\vert \vert y\vert \vert ) .$$
Alors les assertions de l'\'enonc\'e r\'esultent de (2). $\square$

\bigskip
\subsection{Terme constant de $J_{v,\tilde{M},\tilde{S}}^T(\omega,f_{1},f_{2})$}
 Supposons d'abord $\omega$  trivial sur $S^{\theta}(F)$. Pour $\gamma\in \tilde{S}(F)$, posons
$$J_{\tilde{M}}^{\tilde{G}}(\gamma,\omega,f_{1},f_{2})= D^{\tilde{G}}(\gamma)\int_{S^{\theta}(F)\backslash G(F)}\int_{S^{\theta}(F)\backslash G(F)}\bar{f}_{1}(x^{-1}\gamma x)f_{2}(y^{-1}\gamma y)\omega(x^{-1}y)v^{\tilde{G}}_{\tilde{M}}(x,y) \,dy\,dx.$$
Posons 
$$J^{\tilde{G}}_{\tilde{M},\tilde{S}}(\omega,f_{1},f_{2})= \vert W^M(\tilde{S})\vert ^{-1} mes(A_{\tilde{M}}(F)\backslash S^{\theta}(F))\int_{\tilde{S}(F)/(1-\theta)(S(F))}J_{\tilde{M}}^{\tilde{G}}(\gamma,\omega,f_{1},f_{2})\,d\gamma.$$
Si $\omega$ n'est pas trivial sur $S^{\theta}(F)$, on pose simplement
$$J^{\tilde{G}}_{\tilde{M},\tilde{S}}(\omega,f_{1},f_{2})=0.$$

 \ass{Lemme}{(i) Les expressions ci-dessus sont absolument convergentes.
 
 (ii) La fonction $T\mapsto f(T)= J^T_{v,\tilde{M},\tilde{S}}(\omega,f_{1},f_{2})$ appartient \`a $PolExp$. Si $F$ est archim\'edien, on a l'\'egalit\'e $c_{0}(f)=(-1)^{a_{\tilde{M}}-a_{\tilde{G}}}J^{\tilde{G}}_{\tilde{M},\tilde{S}}(\omega,f_{1},f_{2})$. Si $F$ est non-archim\'edien, pour tout r\'eseau ${\cal R}\subset {\cal A}_{M_{0},F}\otimes_{{\mathbb Z}}{\mathbb Q}$, on a l'\'egalit\'e
 $$lim_{k\to \infty}c_{\frac{1}{k}{\cal R},0}(f)=(-1)^{a_{\tilde{M}}-a_{\tilde{G}}}J^{\tilde{G}}_{\tilde{M},\tilde{S}}(\omega,f_{1},f_{2}).$$}
 
 {\bf Remarque.} Comme en 3.23, on appellera "terme constant" de $f$ le terme $c_{0}(f)$ si $F$ est archim\'edien, $lim_{k\to \infty}c_{\frac{1}{k}{\cal R},0}(f)$ si $F$ est non-archim\'edien.
 
 Preuve. Pour l'assertion (i), on peut supposer $\omega$ trivial sur $S^{\theta}(F)$. Il r\'esulte  de la variante tordue de 1.7(1)  et de 4.7(3) que l'on a une majoration
 $$\vert v_{\tilde{M}}(x,y)\vert <<(1+log(\vert \vert x\vert \vert )+log(\vert \vert y\vert \vert ))^D$$
 pour tous $x,y\in G(F)$, o\`u $D=a_{\tilde{M}}-a_{\tilde{G}}$. En tenant compte de 4.3(3), il suffit de prouver la convergence de
$$(1)\qquad \int_{S^{\theta}(F)\backslash G(F)}\int_{A_{\tilde{M}}(F)\backslash G(F)}\vert  f_{1}(x^{-1}\gamma x)f_{2}(y^{-1}\gamma y) \vert (1+\vert log(D^{\tilde{G}}(\gamma))\vert )^D\,dy\,dx\,D^{\tilde{G}}(\gamma)\,d\gamma.$$
Gr\^ace \`a 4.2(1), on est ramen\'e \`a celle de l'int\'egrale
$$\int_{\Omega}(1+\vert log(D^{\tilde{G}}(\gamma))\vert )^D\,d\gamma.$$
Puisqu'on a sur $\Omega$ une majoration $(1+\vert log(D^{\tilde{G}}(\gamma))\vert )^D<<D^{\tilde{G}}(\gamma)^{-\eta}$, cela r\'esulte de 4.2(2).

Prouvons (ii), en supposant comme souvent $F$ non-archim\'edien. Le lemme pr\'ec\'edent nous dit que chaque fonction $T\mapsto v_{\tilde{M}}^T(x,y)$ appartient \`a $PolExp$.  Plus pr\'ecis\'ement, pour tout r\'eseau ${\cal R}\subset {\cal A}_{M_{0},F}\otimes_{{\mathbb Z}}{\mathbb Q}$,  on peut \'ecrire
$$v_{\tilde{M}}^T(x,y)=\sum_{\mu\in \chi_{{\cal R}}}e^{<\mu,T>}p_{{\cal R},\mu,x,y}(T),$$
o\`u $\chi_{{\cal R}}$ est ind\'ependant de $x$ et $y$ et les polyn\^omes $p_{{\cal R},\mu,x,y}$ ont un degr\'e born\'e ind\'ependamment de $x$ et $y$. Comme en 3.23, les coefficients de ces polyn\^omes se calculent par interpolation, donc ont les m\^emes propri\'et\'es de croissance que les fonctions $v_{\tilde{M}}^T(x,y)$ elles-m\^emes. Il en r\'esulte que, dans la formule int\'egrale d\'efinissant $J^T_{v,\tilde{M},\tilde{S}}(\omega,f_{1},f_{2})$, les int\'egrales commutent au d\'eveloppement en $T$. Cela entra\^{\i}ne que la fonction $T\mapsto f(T)= J^T_{v,\tilde{M},\tilde{S}}(\omega,f_{1},f_{2})$ appartient \`a $PolExp$. Cela entra\^{\i}ne aussi que  le terme $c_{\frac{1}{k}{\cal R},0}(f)$ se calcule en rempla\c{c}ant dans l'int\'egrale les fonctions $v_{\tilde{M}}^T(x,y)$ par leurs termes similaires. Ceux-ci sont calcul\'es par le lemme 4.7, \`a une erreur pr\`es. L'int\'egrale des termes d'erreurs est essentiellement major\'ee par le produit de $k^{-1}$ et de l'int\'egrale (1) ci-dessus. Puisque celle-ci est convergente, ce terme d'erreur tend vers $0$ quand $k$ tend vers l'infini. Le terme principal est nul si $\omega$ n'est pas trivial sur $A_{\tilde{M}}(F)$ et on obtient dans ce cas la formule cherch\'ee. Supposons $\omega$ trivial sur $A_{\tilde{M}}(F)$. On obtient dans ce cas pour terme constant
$$(2) \qquad (-1)^{a_{\tilde{M}}-a_{\tilde{G}}} \vert W^M(\tilde{S})\vert ^{-1} \int_{\tilde{S}(F)/(1-\theta)(S(F))}$$
 $$\int_{S^{\theta}(F)\backslash G(F)}\int_{A_{\tilde{M}}(F)\backslash G(F)}\bar{f}_{1}(x^{-1}\gamma x)f_{2}(y^{-1}\gamma y)\omega(x^{-1}y)v^{\tilde{G}}_{\tilde{M}}(x,y) \,dy\,dx\,D^{\tilde{G}}(\gamma)\,d\gamma.$$
 On remarque que la fonction $y\mapsto v_{\tilde{M}}(x,y)$ est invariante \`a gauche par $M(F)$, donc par $S(F)$. En factorisant l'int\'egrale en $y$ en une int\'egrale sur $S^{\theta}(F)\backslash G(F)$ et une int\'egrale en $A_{\tilde{M}}(F)\backslash S^{\theta}(F)$, on voit appara\^{\i}tre l'int\'egrale
 $$\int_{A_{\tilde{M}}(F)\backslash S^{\theta}(F)}\omega(y)\,dy.$$
 Elle est nulle (et donc aussi le terme constant), si $\omega$ n'est  pas trivial sur $S^{\theta}(F)$. Si $\omega$ est trivial sur $S^{\theta}(F)$, elle vaut $mes(A_{\tilde{M}}(F)\backslash S^{\theta}(F))$.  La formule (2) devient $(-1)^{a_{\tilde{M}}-a_{\tilde{G}}}J^{\tilde{G}}_{\tilde{M},\tilde{S}}(\omega,f_{1},f_{2})$. $\square$

\bigskip

\subsection{La formule g\'eom\'etrique}
Pour tout $\tilde{M}\in {\cal L}(\tilde{M}_{0})$, on note $T_{ell}(\tilde{M},\omega)$ l'ensemble des $\tilde{S}\in T_{ell}(\tilde{M})$ tels que $\omega$ soit trivial sur $S^{\theta}(F)$. On pose
$$J_{g\acute{e}om,\tilde{M}}^{\tilde{G}}(\omega,f_{1},f_{2})=\sum_{\tilde{S}\in T_{ell}(\tilde{M},\omega)}\vert W^M(\tilde{S})\vert ^{-1}mes(A_{\tilde{M}}(F)\backslash S^{\theta}(F))$$
$$\int_{\tilde{S}(F)/(1-\theta)(S(F))}J_{\tilde{M}}^{\tilde{G}}(\gamma,\omega,f_{1},f_{2})\,d\gamma.$$
On pose
$$J^{\tilde{G}}_{g\acute{e}om}(\omega,f_{1},f_{2})=\sum_{\tilde{M}\in {\cal L}(\tilde{M}_{0})}\vert \tilde{W}^M\vert \vert \tilde{W}^G\vert ^{-1}(-1)^{a_{\tilde{M}}-a_{\tilde{G}}}J_{g\acute{e}om,\tilde{M}}^{\tilde{G}}(\omega,f_{1},f_{2}).$$

\ass{Proposition}{Il existe une unique fonction $T\mapsto \varphi(T)$ qui appartient \`a $PolExp$ et qui v\'erifie pour tout r\'eel $r$ la majoration
$$\vert J^T(\omega,f_{1},f_{2})-\varphi(T)\vert <<\vert T\vert ^{-r}$$
pour tout $T$ dans le c\^one o\`u $J^T(\omega,f_{1},f_{2})$ est d\'efinie. Si $F$ est archim\'edien, on a l'\'egalit\'e
$$c_{0}(\varphi)=J_{g\acute{e}om}^{\tilde{G}}(\omega,f_{1},f_{2}).$$
Si $F$ est non-archim\'edien, pour tout r\'eseau ${\cal R}\subset {\cal A}_{M_{0},F}\otimes_{{\mathbb Z}}{\mathbb Q}$, on a l'\'egalit\'e
$$lim_{k\to \infty}c_{\frac{1}{k}{\cal R},0}(\varphi)=J_{g\acute{e}om}^{\tilde{G}}(\omega,f_{1},f_{2}).$$}

Preuve. Cela r\'esulte de 4.3(1) et du lemme pr\'ec\'edent. $\square$

\bigskip

\section{La formule des traces locale tordue, version non invariante}

\bigskip

\subsection{Le th\'eor\`eme}
On a d\'efini des expressions $J_{spec}^{\tilde{G}}(\omega,f_{1},f_{2})$ et $J_{g\acute{e}om}^{\tilde{G}}(\omega,f_{1},f_{2})$ en 3.26 et  4.9.

\ass{Th\'eor\`eme}{On a l'\'egalit\'e $J_{spec}^{\tilde{G}}(\omega,f_{1},f_{2})=J_{g\acute{e}om}^{\tilde{G}}(\omega,f_{1},f_{2})$.}

Il suffit de comparer les propositions 3.26 et 4.9. $\square$

 {\bf Remarque.} Nos mesures v\'erifiant les conditions de coh\'erence de 1.2, on voit que les expressions du th\'eor\`eme sont proportionnelles au carr\'e de la mesure sur $G(F)$, en un sens \'evident, et inversement proportionnelles \`a la mesure sur $A_{\tilde{G}}(F)$. Elles ne d\'ependent d'aucune autre mesure.

 \subsection{Extension de la formule des traces locale tordue aux fonctions $C^{\infty}$ \`a support compact}
 Dans ce paragraphe, $F$ est archim\'edien. Notons $\mathfrak{U}$ l'alg\`ebre enveloppante de l'alg\`ebre de Lie r\'eelle de $G(F)$. Elle agit \`a droite et \`a gauche sur les fonctions $C^{\infty}$ sur $\tilde{G}(F)$. Pour $X,Y\in \mathfrak{U}$ et $f$ une telle fonction, on note $XfY$ l'image de $f$ par les actions \`a gauche de $X$  et \`a droite de $Y$. Soit $\Omega$ un sous-ensemble compact de $\tilde{G}(F)$. Notons $C_{c}^{\infty}(\Omega)$ l'espace des \'el\'ements de $C_{c}^{\infty}(\tilde{G}(F))$ \`a support dans $\Omega$. On munit $C_{c}^{\infty}(\Omega)$ de la topologie d\'efinie par les semi-normes
 $$f\mapsto sup_{\gamma\in \Omega}\vert XfY(\gamma)\vert ,$$
 pour $X,Y\in \mathfrak{U}$. On a l'\'egalit\'e
 $$C_{c}^{\infty}(\tilde{G}(F))=\bigcup_{\Omega}C_{c}^{\infty}(\Omega)$$
 o\`u $\Omega$ d\'ecrit les sous-ensembles compacts de $\tilde{G}(F)$. On  munit $C_{c}^{\infty}(\tilde{G}(F))$ de la topologie limite inductive des topologies que l'on vient de d\'efinir sur les sous-espaces $C_{c}^{\infty}(\Omega)$. Autrement dit, une suite $(f_{n})_{n\in {\mathbb N}}$ d'\'el\'ements de $C_{c}^{\infty}(\tilde{G}(F))$ converge vers $f\in C_{c}^{\infty}(\tilde{G}(F))$ si et seulement si les deux conditions suivantes sont satisfaites:
 
 - il existe un sous-ensemble compact $\Omega\subset \tilde{G}(F)$ tel que $f$ et chaque $f_{n}$ soient \`a support dans $\Omega$;
 
 - pour tous $X,Y\in \mathfrak{U}$, on  a $lim_{n\to \infty}sup_{\gamma\in \tilde{G}(F)}\vert XfY(\gamma)-Xf_{n}Y(\gamma)\vert =0$.
 
 Comme dans le cas non tordu, le th\'eor\`eme de Peter-Weyl entra\^{\i}ne que $C_{c}^{\infty}(\tilde{G}(F),K)$ est dense dans $C_{c}^{\infty}(\tilde{G}(F))$.

 \ass{Proposition}{Soient $f_{1},f_{2}\in C_{c}^{\infty}(\tilde{G}(F))$.
 
 (i) Toutes les expressions intervenant dans la d\'efinition de $J^{\tilde{G}}_{spec}(\omega,f_{1},f_{2})$ et $J^{\tilde{G}}_{g\acute{e}om}(\omega,f_{1},f_{2})$ sont convergentes et ces expressions elles-m\^emes le sont.
 
 (ii) On a l'\'egalit\'e $J^{\tilde{G}}_{spec}(\omega,f_{1},f_{2})=J^{\tilde{G}}_{g\acute{e}om}(\omega,f_{1},f_{2})$.}
 
 Preuve. Le c\^ot\'e g\'eom\'etrique est \`a peu pr\`es trivial. Pour $i=1,2$, on peut majorer $\vert f_{i}\vert$ par une fonction $f'_{i}\in C_{c}^{\infty}(\tilde{G}(F),K)$ \`a valeurs positives ou nulles. Chaque terme de $J^{\tilde{G}}_{g\acute{e}om}(\omega,f_{1},f_{2})$ est major\'e par le m\^eme terme relatif aux fonctions $f'_{1}$ et $f'_{2}$ et au caract\`ere $\omega$ trivial. La convergence de ces termes entra\^{\i}ne celle des termes initiaux. On voit de m\^eme que $J^{\tilde{G}}_{g\acute{e}om}(\omega,f_{1},f_{2})$ est continu en $f_{1}$ et $f_{2}$.

Passons aux termes spectraux. Dans le cas o\`u $F={\mathbb C}$, on voit que notre probl\`eme est \'equivalent \`a celui concernant les objets d\'eduits de $G$ et $\tilde{G}$ par restriction des scalaires de ${\mathbb C}$ \`a ${\mathbb R}$. On suppose donc $F={\mathbb R}$. On peut fixer un sous-ensemble compact $\Omega\subset \tilde{G}({\mathbb R})$ et supposer nos fonctions $f_{1}$ et $f_{2}$ \`a support dans $\Omega$.
On a une formule de la forme
$$(1) \qquad J_{spec}^{\tilde{G}}(\omega,f_{1},f_{2})=\sum_{\tilde{M}\in {\cal L}(\tilde{M}_{0})}\sum_{\tau\in (E_{disc}(\tilde{M},\omega)/conj)/i{\cal A}_{\tilde{M}}^*}c(\tau)\int_{i{\cal A}_{\tilde{M}}^*} J_{\tilde{M}}^{\tilde{G}}(\pi_{\tau_{\lambda}},f_{1},f_{2})\,d\lambda,$$
o\`u les $c(\tau)$ sont des coefficients uniform\'ement born\'es.  On fixe $\gamma_{0}\in \tilde{M}_{0}({\mathbb R})$ v\'erifiant la condition 2.1(4)   et on d\'efinit des fonctions $\varphi_{i}$ sur $G({\mathbb R})$ pour $i=1,2$ par $\varphi_{i}(x)=f_{i}(x\gamma_{0})$. Fixons $\tilde{M}$ et $\tau$ intervenant dans la formule ci-dessus. On a une \'egalit\'e
$$J_{\tilde{M}}^{\tilde{G}}(\pi_{\tau_{\lambda}},f_{1},f_{2})=trace({\cal X}^{\tilde{G}}_{\tilde{M}}(\pi_{\tau_{\lambda}})(\overline{\tilde{\Pi}_{\boldsymbol{\tau}_{\tilde{\lambda}}}(f_{1})}\otimes \tilde{\Pi}_{\boldsymbol{\tau}_{\tilde{\lambda}}}(f_{2}))),$$
o\`u ${\cal X}^{\tilde{G}}_{\tilde{M}}(\pi_{\tau_{\lambda}})$ est l'op\'erateur d\'eduit de la $(\tilde{G},\tilde{M})$ famille $({\cal X}(\pi_{\tau_{\lambda}};\Lambda,\tilde{Q}))_{\tilde{Q}\in {\cal P}(\tilde{M})}$ d\'efinie par
$${\cal X}(\pi_{\tau_{\lambda}};\Lambda,\tilde{Q})=\overline{{\cal M}(\pi_{\tau_{\lambda}};\Lambda,\tilde{\bar{Q}})}\otimes {\cal M}(\pi_{\tau_{\lambda}};\Lambda,\tilde{Q}).$$
On a aussi
$$\overline{\tilde{\Pi}_{\boldsymbol{\tau}_{\tilde{\lambda}}}(f_{1})}\otimes \tilde{\Pi}_{\boldsymbol{\tau}_{\tilde{\lambda}}}(f_{2})=\left(\overline{\Pi_{\tau_{\lambda}}(\varphi_{1})}\otimes \Pi_{\tau_{\lambda}}(\varphi_{2})\right) U_{\tau_{\lambda}},$$
o\`u
$$U_{\tau_{\lambda}}=\overline{\tilde{\Pi}_{\boldsymbol{\tau}_{\tilde{\lambda}}}(\gamma_{0})}\otimes \tilde{\Pi}_{\boldsymbol{\tau}_{\tilde{\lambda}}}(\gamma_{0}).$$
On a donc
$$J_{\tilde{M}}^{\tilde{G}}(\pi_{\tau_{\lambda}},f_{1},f_{2})=trace(j_{\tilde{M}}^{\tilde{G}}(\pi_{\tau_{\lambda}},\varphi_{1},\varphi_{2})),$$
o\`u
$$j_{\tilde{M}}^{\tilde{G}}(\pi_{\tau_{\lambda}},\varphi_{1},\varphi_{2})={\cal X}^{\tilde{G}}_{\tilde{M}}(\pi_{\tau_{\lambda}})\left(\overline{\Pi_{\tau_{\lambda}}(\varphi_{1})}\otimes \Pi_{\tau_{\lambda}}(\varphi_{2})\right) U_{\tau_{\lambda}}.$$

   On peut bien s\^ur fixer $\tilde{M}$ et ne consid\'erer que la sous-somme de l'expression (1) index\'ee par $\tilde{M}$. Un \'el\'ement $\tau\in E_{disc}(\tilde{M},\omega)$ est de la forme $(M_{disc},\sigma,\tilde{r})$, o\`u $M_{disc}\subset M$ et $\sigma$ est de la s\'erie discr\`ete de $M_{disc}({\mathbb R})$. On peut encore fixer $M_{disc}$, que l'on suppose semi-standard, et ne consid\'erer que la sous-somme des $\tau$ issus de ce Levi.  On peut identifier la composante neutre de $A_{M_{disc}}({\mathbb R})$ \`a ${\cal A}_{M_{disc}}$. Puisque $\tau $ est discret, l'ensemble $W^{\tilde{M}}_{reg}(\sigma)$ n'est pas vide. La restriction du caract\`ere central de $\sigma$ \`a $ {\cal A}_{M_{disc}}$  est naturellement param\'etr\'ee par un \'el\'ement de $i{\cal A}_{M_{disc}}^*$. Cette restriction \'etant fix\'ee par tout \'el\'ement de  $W^{\tilde{M}}_{reg}(\sigma)$ et l'ensemble des points fixes de l'action d'un tel \'el\'ement dans $i{\cal A}_{M_{disc}}^*$ \'etant \'egal \`a $i{\cal A}_{\tilde{M}}^*$, le param\`etre en question appartient \`a $i{\cal A}_{\tilde{M}}^*$. Puisque seule intervient la classe de $\tau$ modulo ce groupe, on peut supposer que la restriction du caract\`ere central de $\sigma$ \`a $ {\cal A}_{M_{disc}}$  est triviale. Fixons comme en 1.10  une sous-alg\`ebre de Cartan $\mathfrak{h}_{{\mathbb C}}$ de l'alg\`ebre de Lie complexifi\'ee de $M_{disc}({\mathbb R})$. Pour  toute repr\'esentation irr\'eductible $\sigma$ de $M_{disc}({\mathbb R})$, on choisit un \'el\'ement   $\mu_{\sigma}\in \mathfrak{h}_{{\mathbb C}}$ dont la classe modulo un certain groupe de Weyl param\`etre le caract\`ere infinit\'esimal de $\sigma$. Notons $\hat{K}$ l'ensemble des repr\'esentations irr\'eductibles de $K$. De m\^eme, un \'el\'ement $\kappa\in \hat{K}$ a un param\`etre $\mu_{\kappa}\in \mathfrak{h}_{{\mathbb C}}$. Rappelons que les op\'erateurs qui interviennent dans l'int\'egrale de (3) agissent dans la repr\'esentation $Ind_{S}^G(\sigma_{\lambda})\otimes Ind_{S}^G(\sigma_{\lambda})$, o\`u $S$ est un \'el\'ement fix\'e de ${\cal P}(M_{disc})$. On d\'ecompose cette repr\'esentation selon les espaces isotypiques pour l'action de $K\times K$. Ces espaces sont donc index\'es par $(\kappa_{1},\kappa_{2})\in \hat{K}\times \hat{K}$. On note $p_{\kappa_{1},\kappa_{2}}$ le projecteur sur l'espace isotypique index\'e par $(\kappa_{1},\kappa_{2})$. Posons
   $$j_{\tilde{M}}^{\tilde{G}}(\pi_{\tau_{\lambda}},\varphi_{1},\varphi_{2};\kappa_{1},\kappa_{2})
   =p_{\kappa_{1},\kappa_{2}}j_{\tilde{M}}^{\tilde{G}}(\pi_{\tau_{\lambda}},\varphi_{1},\varphi_{2})p_{\kappa_{1},\kappa_{2}}.$$
   On a
   $$(2) \qquad J_{\tilde{M}}^{\tilde{G}}(\pi_{\tau_{\lambda}},f_{1},f_{2}) =\sum_{\kappa_{1},\kappa_{2}\in \hat{K}}trace(j_{\tilde{M}}^{\tilde{G}}(\pi_{\tau_{\lambda}},\varphi_{1},\varphi_{2};\kappa_{1},\kappa_{2})).$$
      Remarquons que l'op\'erateur ${\cal X}^{\tilde{G}}_{\tilde{M}}(\pi_{\tau_{\lambda}})$ conserve les espaces isotypiques par construction, tandis que $ U_{\tau_{\lambda}}$ envoie l'espace de type $(\kappa_{1},\kappa_{2})$ sur celui de type $(\kappa'_{1},\kappa'_{2})$, o\`u, pour $i=1,2$, $\kappa'_{i}=ad_{\gamma_{0}}(\kappa_{i}\otimes \omega)$. On a donc
  $$ j_{\tilde{M}}^{\tilde{G}}(\pi_{\tau_{\lambda}},\varphi_{1},\varphi_{2};\kappa_{1},\kappa_{2})= 
   {\cal X}^{\tilde{G}}_{\tilde{M}}(\pi_{\tau_{\lambda}};\kappa_{1},\kappa_{2})\Pi_{\tau_{\lambda}}(\varphi_{1},\varphi_{2};\kappa_{1},\kappa_{2})U_{\tau_{\lambda}}(\kappa_{1},\kappa_{2}),$$
   o\`u
  $${\cal X}^{\tilde{G}}_{\tilde{M}}(\pi_{\tau_{\lambda}};\kappa_{1},\kappa_{2})
  =p_{\kappa_{1},\kappa_{2}}{\cal X}^{\tilde{G}}_{\tilde{M}}(\pi_{\tau_{\lambda}})p_{\kappa_{1},\kappa_{2}};$$
  $$\Pi_{\tau_{\lambda}}(\varphi_{1},\varphi_{2};\kappa_{1},\kappa_{2})=p_{\kappa_{1},\kappa_{2}}
   \left(\overline{\Pi_{\tau_{\lambda}}(\varphi_{1})}\otimes \Pi_{\tau_{\lambda}}(\varphi_{2})\right)p_{\kappa'_{1},\kappa'_{2}};$$  
 $$U_{\tau_{\lambda}}(\kappa_{1},\kappa_{2})=p_{\kappa'_{1},\kappa'_{2}} U_{\tau_{\lambda}}p_{\kappa_{1},\kappa_{2}}.$$
 Puisque $\tilde{\Pi}_{\boldsymbol{\tau}_{\tilde{\lambda}}}$ est unitaire, ce dernier op\'erateur est de norme uniform\'ement born\'ee. Il est clair que les param\`etres associ\'es \`a $\kappa_{i}$ et $\kappa'_{i}$ sont de m\^eme norme, pour $i=1,2$. D'apr\`es [A8] p.174, pour tout r\'eel $r$, il existe une fonction  $c_{r}:C_{c}^{\infty}(G({\mathbb R}))\times C_{c}^{\infty}(G({\mathbb R}))\to {\mathbb R}_{\geq0}$ qui est born\'ee par le sup d'un ensemble fini de semi-normes et qui est telle que  l'on ait la majoration
 $$\vert \vert \Pi_{\tau_{\lambda}}(\varphi_{1},\varphi_{2};\kappa_{1},\kappa_{2})\vert \vert \leq c_{r}(\varphi_{1},\varphi_{2})(1+\vert \vert \mu_{\sigma}\vert \vert )^{-r}(1+\vert \vert \mu_{\kappa_{1}}\vert \vert )^{-r}(1+\vert \vert \mu_{\kappa_{2}}\vert \vert )^{-r}(1+\vert \vert \lambda\vert \vert )^{-r}.$$
On montrera ci-dessous qu'il existe $C\geq0$ et un entier $D\geq0$ de sorte que l'on ait la majoration
$$(3) \qquad \vert \vert {\cal X}^{\tilde{G}}_{\tilde{M}}(\pi_{\tau_{\lambda}};\kappa_{1},\kappa_{2})\vert \vert \leq C(1+\vert \vert \mu_{\sigma}\vert \vert )^{D}(1+\vert \vert \mu_{\kappa_{1}}\vert \vert )^{D}(1+\vert \vert \mu_{\kappa_{2}}\vert \vert )^{D}(1+\vert \vert \lambda\vert \vert )^{D}.$$
Donc $trace(j_{\tilde{M}}^{\tilde{G}}(\pi_{\tau_{\lambda}},\varphi_{1},\varphi_{2};\kappa_{1},\kappa_{2}))$ est born\'ee par le produit des deux expressions ci-dessus et des dimensions de $\kappa_{1}$ et $\kappa_{2}$. On sait que la dimension de $\kappa$ est essentiellement born\'ee par $(1+\vert \vert \mu_{\kappa}\vert \vert )^{D'}$ pour un entier $D'$ convenable. On en d\'eduit pour tout r\'eel $r$ une majoration
$$\vert trace(j_{\tilde{M}}^{\tilde{G}}(\pi_{\tau_{\lambda}},\varphi_{1},\varphi_{2};\kappa_{1},\kappa_{2}))\vert \leq c'_{r}(\varphi_{1},\varphi_{2})(1+\vert \vert \mu_{\sigma}\vert \vert )^{-r}(1+\vert \vert \mu_{\kappa_{1}}\vert \vert )^{-r}(1+\vert \vert \mu_{\kappa_{2}}\vert \vert )^{-r}(1+\vert \vert \lambda\vert \vert )^{-r},$$
o\`u $c'_{r}(\varphi_{1},\varphi_{2})$ est born\'ee par le sup d'un ensemble fini de semi-normes. Gr\^ace \`a (2), on obtient
 $$\vert J_{\tilde{M}}^{\tilde{G}}(\pi_{\tau_{\lambda}},f_{1},f_{2}) \vert \leq c'_{r}(\varphi_{1},\varphi_{2}) (1+\vert \vert \mu_{\sigma}\vert \vert )^{-r}(1+\vert \vert \lambda\vert \vert )^{-r}\left(\sum_{\kappa\in \hat{K}}(1+\vert \vert \mu_{\kappa}\vert \vert )^{-r}\right)^2.$$
 Si $r$ est assez grand, la derni\`ere s\'erie est convergente. L'int\'egrale en $\lambda$ de l'expression ci-dessus l'est aussi et on obtient simplement
 $$\int_{i{\cal A}_{\tilde{M}}^*}\vert J_{\tilde{M}}^{\tilde{G}}(\pi_{\tau_{\lambda}},f_{1},f_{2}) \vert \,d\lambda\leq C'c'_{r}(\varphi_{1},\varphi_{2})(1+\vert \vert \mu_{\sigma}\vert \vert )^{-r},$$
 pour une certaine constante absolue $C'$. Remarquons que chaque repr\'esentation $\sigma$ de la s\'erie discr\`ete  de $M_{disc}({\mathbb R})$ ne donne naissance qu'\`a un nombre uniform\'ement born\'e de triplets $\tau$. La sous-somme de $J_{spec}^{\tilde{G}}(\omega,f_{1},f_{2})$ o\`u on a fix\'e $\tilde{M}$ et $M_{disc}$, mais o\`u on a remplac\'e la fonction que l'on int\`egre par sa valeur absolue, est donc essentiellement major\'ee par
 $$c'_{r}(\varphi_{1},\varphi_{2})\sum_{\sigma}(1+\vert \vert \mu_{\sigma}\vert \vert )^{-r},$$
 o\`u on somme sur les repr\'esentations irr\'eductibles $\sigma$ de $M_{disc}({\mathbb R})$ de la s\'erie discr\`ete et dont la restriction du caract\`ere central \`a $
 {\cal A}_{M_{disc}}$ est triviale. Cette derni\`ere condition signifie que $\mu_{\sigma,M_{disc}}=0$. Comme on l'a dit en 1.10, la projection $\mu_{\sigma}^{M_{disc}}$ parcourt un r\'eseau de $i\mathfrak{h}^{M,*}$. De plus,  on sait qu'il n'y a qu'un nombre uniform\'ement born\'e de s\'eries discr\`etes d'un param\`etre donn\'e. Donc, si $r$ est assez grand, la s\'erie ci-dessus est convergente. Cela prouve les assertions du (i) de l'\'enonc\'e concernant $J_{spec}^{\tilde{G}}(\omega,f_{1},f_{2})$.  La propri\'et\'e du terme $c'_{r}(\varphi_{1},\varphi_{2})$ prouve en m\^eme temps la continuit\'e de cette expression en $f_{1}$ et $f_{2}$.
 
 Puique les deux membres de l'\'egalit\'e du (ii) de l'\'enonc\'e sont continus en $f_{1}$ et $f_{2}$, cette \'egalit\'e r\'esulte par continuit\'e du th\'eor\`eme 5.1, les fonctions lisses \`a support compacts et $K$-finies  \'etant denses dans $C_{c}^{\infty}(\tilde{G}({\mathbb R}))$.
 
 Il reste \`a prouver la majoration (3). Tout d'abord, en reprenant la preuve du  lemme 3.25, on voit que l'on peut remplacer la $(\tilde{G},\tilde{M})$-famile $({\cal X}(\pi_{\tau_{\lambda}};\Lambda,\tilde{Q}))_{\tilde{Q}\in {\cal P}(\tilde{M})}$ par 
 $$({\bf r}_{reg}(\pi_{\tau_{\lambda}};\Lambda,\tilde{Q}){\cal X}_{reg}(\pi_{\tau_{\lambda}};\Lambda,\tilde{Q}))_{\tilde{Q}\in {\cal P}(\tilde{M})},$$ o\`u ${\cal X}_{reg}(\pi_{\tau_{\lambda}};\Lambda,\tilde{Q})$ ne contient que des op\'erateurs d'entrelacement normalis\'es et
 $${\bf r}_{reg}(\pi_{\tau_{\lambda}}:\Lambda,\tilde{Q})=r_{\bar{Q}\vert Q}(\sigma_{\lambda})r_{\bar{P}\vert P}(\sigma_{\lambda})^{-1}r_{Q\vert \bar{Q}}(\sigma_{\lambda+\Lambda})r_{P\vert \bar{P}}(\sigma_{\lambda+\Lambda})^{-1}$$
 ($\tilde{P}$ est un \'el\'ement fix\'e de ${\cal P}(\tilde{M})$, cf. 3.25).  Notons que cette fonction est produit de fonctions $r_{\alpha}(\sigma_{\lambda})r_{-\alpha}(\sigma_{\lambda})^{-1}$ et $r_{\alpha}(\sigma_{\lambda+\Lambda})r_{-\alpha}(\sigma_{\lambda+\Lambda})^{-1}$ comme en 1.10. Evidemment, ${\cal X}^{\tilde{G}}_{\tilde{M}}(\pi_{\tau_{\lambda}};\kappa_{1},\kappa_{2})$ se d\'eduit de la $(\tilde{G},\tilde{M})$-famille $({\bf r}_{reg}(\pi_{\tau_{\lambda}}:\Lambda,\tilde{Q})p_{\kappa_{1},\kappa_{2}}{\cal X}_{reg}(\pi_{\tau_{\lambda}};\Lambda,\tilde{Q})p_{\kappa_{1},\kappa_{2}})_{\tilde{Q}\in {\cal P}(\tilde{M})}$. La formule 1.7(1) montre qu'il nous suffit de majorer les valeurs en $\Lambda=0$ des d\'eriv\'ees en $\Lambda$ d'ordre au plus $a_{\tilde{M}}-a_{\tilde{G}}$ de tous les termes intervenant dans la d\'efinition de la $(\tilde{G},\tilde{M})$-famille ci-dessus. Les op\'erateurs d'entrelacement normalis\'es \'etant unitaires, cela nous ram\`ene \`a majorer les op\'erateurs
 $$(4)\qquad D(p_{\kappa_{1},\kappa_{2}}R_{Q\vert P}(\sigma_{\lambda})p_{\kappa_{1},\kappa_{2}})$$
 et les fonctions
 $$(5)\qquad D(r_{\alpha}(\sigma_{\lambda})r_{-\alpha}(\sigma_{\lambda})^{-1})$$
 o\`u $D$ est une d\'erivation d'ordre au plus $a_{\tilde{M}}-a_{\tilde{G}}$ s'appliquant \`a la variable $\lambda$. Comme on l'a dit ci-dessus, nos hypoth\`eses entra\^{\i}nent $\mu_{\sigma,M_{disc}}=0$, donc $\mu_{\sigma}$ est orthogonal \`a $\lambda$. La majoration cherch\'ee de (4) r\'esulte alors du lemme 2.1 de [A5]. On a d\'ecrit les fonctions $r_{\alpha}(\sigma_{\lambda})r_{-\alpha}(\sigma_{\lambda})^{-1}$ en 1.10. Compte tenu de l'\'egalit\'e $\mu_{\sigma,M_{disc}}=0$, elles sont produit de termes
  $$(6)\qquad \frac{<\mu_{\sigma},\check{\beta}>-<\lambda,\check{\beta}>} {<\mu_{\sigma},\check{\beta}>+<\lambda,\check{\beta}>} $$
  et \'eventuellement d'un terme
  $$\frac{\Gamma(\frac{<\lambda,\check{\beta}>}{2})\Gamma(\frac{-<\lambda,\check{\beta}>+1}{2})}{\Gamma(\frac{-<\lambda,\check{\beta}>}{2})\Gamma(\frac{<\lambda,\check{\beta}>+1}{2})}.$$
  Ce dernier terme a des d\'eriv\'ees \`a croissance mod\'eree et ne d\'epend pas de $\mu_{\sigma}$. Il v\'erifie donc la majoration requise. Consid\'erons une fonction (6). Il r\'esulte de ce que l'on a rappel\'e en 1.10 que le produit $<\mu_{\sigma},\check{\beta}>$ appartient \`a un sous-groupe discret fixe de ${\mathbb R}$, tandis que $<\lambda,\check{\beta}>$ est imaginaire. Ou bien $<\mu_{\sigma},\check{\beta}>=0$ et la fonction est constante \'egale \`a $-1$. Ou bien $<\mu_{\sigma},\check{\beta}>\not=0$ et alors la valeur absolue du d\'enominateur est uniform\'ement minor\'ee par une constante strictement positive. La majoration (5) requise s'ensuit. Cela ach\`eve la preuve. $\square$

 \bigskip
 
 \subsection{Formules de descente pour les $(\tilde{G},\tilde{M})$-familles}
  On rappelle dans ce paragraphe plusieurs formules g\'en\'erales d'Arthur concernant les $(\tilde{G},\tilde{M})$-familles. Pour $i=1,2$, soient $\tilde{M}_{i}\in {\cal L}(\tilde{M}_{0})$ et $(x_{i}(\Lambda;\tilde{P}))_{\tilde{P}\in {\cal P}(\tilde{M}_{i})}$ une $(\tilde{G},\tilde{M}_{i})$-famille. Soit $\tilde{L}\in {\cal L}(\tilde{M}_{0})$ tel que $\tilde{M}_{i}\subset \tilde{L}$ pour $i=1,2$. De la $(\tilde{G},\tilde{M}_{i})$-famille $(x_{i}(\Lambda;\tilde{P}))_{\tilde{P}\in {\cal P}(\tilde{M}_{i})}$ se d\'eduit une $(\tilde{G},\tilde{L})$-famille $(x_{i}(\Lambda;\tilde{Q}))_{\tilde{Q}\in {\cal P}(\tilde{L})}$. Pour $\Lambda\in i{\cal A}_{\tilde{L}}^*$, on a $x_{i}(\Lambda;\tilde{Q})=x_{i}(\Lambda;\tilde{P})$ pour $\tilde{P}\in {\cal P}(\tilde{M}_{i}) $ tel que $P\subset Q$. Notons  $(y(\Lambda,\tilde{Q}))_{\tilde{Q}\in {\cal P}(\tilde{L})}$ la famille produit, c'est-\`a-dire $y(\Lambda,\tilde{Q})=x_{1}(\Lambda,\tilde{Q})x_{2}(\Lambda,\tilde{Q})$. De cette $(\tilde{G},\tilde{L})$-famille se d\'eduit une fonction $y_{\tilde{L}}^{\tilde{G}}(\Lambda)$ sur $i{\cal A}_{\tilde{L}}^*$. Rappelons que, pour $i=1,2$ et pour tout $\tilde{Q}_{i}=\tilde{L}_{i}U_{Q_{i}}\in {\cal F}(\tilde{M}_{i})$, on d\'eduit de  $(x_{i}(\Lambda;\tilde{P}))_{\tilde{P}\in {\cal P}(\tilde{M}_{i})}$ une $(\tilde{L}_{i},\tilde{M}_{i})$-famille  $(x_{i}^{\tilde{Q}_{i}}(\Lambda;\tilde{P}))_{\tilde{P}\in {\cal P}^{\tilde{Q}_{i}}(\tilde{M}_{i})}$, o\`u ${\cal P}^{\tilde{Q}_{i}}(\tilde{M}_{i})$ est l'ensemble des $\tilde{P}\in {\cal P}(\tilde{M}_{i})$ tels que $\tilde{P}\subset \tilde{Q}_{i}$. D'o\`u une fonction $x_{i,\tilde{M}_{i}}^{\tilde{Q}_{i}}(\Lambda)$ sur $i{\cal A}_{\tilde{M}_{i}}^*$. Dans l'espace ${\cal A}_{\tilde{M}_{1}}^{\tilde{G}}\oplus {\cal A}_{\tilde{M}_{2}}^{\tilde{G}}$, introduisons l'image $\Delta({\cal A}_{\tilde{L}}^{\tilde{G}})$ de ${\cal A}_{\tilde{L}}^{\tilde{G}}$ par le plongement diagonal. Notons $\Delta({\cal A}_{\tilde{L}}^{\tilde{G}})^{\bot}$ son orthogonal. Pour tout couple $(\tilde{L}_{1},\tilde{L}_{2})\in {\cal L}(\tilde{M}_{1})\times {\cal L}(\tilde{M}_{2})$, consid\'erons les conditions suivantes
  
  (1) $\Delta({\cal A}_{\tilde{L}}^{\tilde{G}})\oplus({\cal A}_{\tilde{L}_{1}}^{\tilde{G}}\oplus {\cal A}_{\tilde{L}_{2}}^{\tilde{G}})={\cal A}_{\tilde{M}_{1}}^{\tilde{G}}\oplus {\cal A}_{\tilde{M}_{2}}^{\tilde{G}}$;
 
 (2)  $\Delta({\cal A}_{\tilde{L}}^{\tilde{G}})^{\bot}\oplus({\cal A}_{\tilde{M}_{1}}^{\tilde{L}_{1}}\oplus {\cal A}_{\tilde{M}_{2}}^{\tilde{L}_{2}})={\cal A}_{\tilde{M}_{1}}^{\tilde{G}}\oplus {\cal A}_{\tilde{M}_{2}}^{\tilde{G}}$; 
 
 (3) l'application
 $$\begin{array}{ccc}{\cal A}_{\tilde{M}_{1}}^{\tilde{L}_{1}}\oplus {\cal A}_{\tilde{M}_{2}}^{\tilde{L}_{2}}&\to&{\cal A}_{\tilde{L}}^{\tilde{G}}\\ (H_{1},H_{2})&\mapsto&H_{1,\tilde{L}}+H_{2,\tilde{L}}\\ \end{array}$$
 est un isomorphisme.
 
 On v\'erifie qu'elles sont \'equivalentes. Si elles le sont, on note $d^{\tilde{G}}_{\tilde{M}_{1},\tilde{M}_{2}}(\tilde{L};\tilde{L}_{1},\tilde{L}_{2})$ le jacobien de l'isomorphisme (3), chaque espace \'etant muni des  mesures  fix\'ees en 1.2. Si les conditions ne sont pas v\'erifi\'ees, on pose $d^{\tilde{G}}_{\tilde{M}_{1},\tilde{M}_{2}}(\tilde{L};\tilde{L}_{1},\tilde{L}_{2})=0$. Fixons un point $\xi\in {\cal A}_{\tilde{M}_{1}}^{\tilde{G}}\oplus {\cal A}_{\tilde{M}_{2}}^{\tilde{G}}$ en position g\'en\'erale. Si les conditions ci-dessus sont v\'erifi\'ees, notons $(\xi_{1},\xi_{2})\in {\cal A}_{\tilde{L}_{1}}^{\tilde{G}}\oplus {\cal A}_{\tilde{L}_{2}}^{\tilde{G}}$ la projection (non orthogonale) de $\xi$ relativement \`a la d\'ecomposition (1). Puisque $\xi$ est en position g\'en\'erale, il existe pour $i=1,2$ un unique $\tilde{Q}_{i}\in {\cal P}(\tilde{L}_{i})$ tel que $\xi_{i}$ soit dans la chambre positive associ\'ee \`a cet espace parabolique.  Cela d\'efinit une application qui, \`a $(\tilde{L}_{1},\tilde{L}_{2})$ tel que $d^{\tilde{G}}_{\tilde{M}_{1},\tilde{M}_{2}}(\tilde{L};\tilde{L}_{1},\tilde{L}_{2})\not=0$, associe un couple $(\tilde{Q}_{1},\tilde{Q}_{2})\in {\cal P}(\tilde{L}_{1})\times {\cal P}(\tilde{L}_{2})$.
 
 On a l'\'egalit\'e
 
 $$(4) \qquad y_{\tilde{L}}^{\tilde{G}}(\Lambda)=\sum_{\tilde{L}_{1}\in {\cal L}(\tilde{M}_{1}),\tilde{L}_{2}\in {\cal L}(\tilde{M}_{2})}d^{\tilde{G}}_{\tilde{M}_{1},\tilde{M}_{2}}(\tilde{L};\tilde{L}_{1},\tilde{L}_{2})x_{1,\tilde{M}_{1}}^{\tilde{Q}_{1}}(\Lambda)x_{2,\tilde{M}_{2}}^{\tilde{Q}_{2}}(\Lambda)$$
 pour tout $\Lambda\in i{\cal A}_{\tilde{L}}^*$.
 
  Dans le cas non tordu, cela r\'esulte de [A6] proposition 7.1, appliqu\'ee au groupe $G\times G$, \`a son Levi $M_{1}\times M_{2}$ et \`a l'espace $\mathfrak{b}=\Delta({\cal A}_{L}^G)\oplus ({\cal A}_{G}\oplus {\cal A}_{G})$. La preuve s'\'etend au cas tordu. 
  
   On utilisera un cas particulier de la relation (4) o\`u les d\'efinitions se simplifient. Soient $\tilde{M},\tilde{M}_{1},\tilde{M}_{2}\in {\cal L}(\tilde{M}_{0})$ tels que $\tilde{M}\subset \tilde{M}_{i}$ pour $i=1=2$ et ${\cal A}_{\tilde{M}}^{\tilde{M}_{1}}\cap {\cal A}_{\tilde{M}}^{\tilde{M}_{2}}=0$. Il existe un unique $\tilde{L}\in {\cal L}(\tilde{M})$ tel que
  $${\cal A}_{\tilde{L}}={\cal A}_{\tilde{M}_{1}}\cap {\cal A}_{\tilde{M}_{2}},$$
  \`a savoir le commutant dans $\tilde{G}$ du  tore $(A_{\tilde{M}_{1}}\cap A_{\tilde{M}_{2}})^0$. On a les \'egalit\'es \'equivalentes
  $$ {\cal A}_{\tilde{M}}^{\tilde{L}}={\cal A}_{\tilde{M}_{1}}^{\tilde{L}}\oplus {\cal A}_{\tilde{M}_{2}}^{\tilde{L}},$$
  $$(5) \qquad    {\cal A}_{\tilde{M}}^{\tilde{L}}={\cal A}^{\tilde{M}_{1}}_{\tilde{M}}\oplus {\cal A}^{\tilde{M}_{2}}_{\tilde{M}}.$$
  On voit que les conditions \'equivalentes (1), (2) et (3) sont \'equivalentes aux deux \'egalit\'es \'equivalentes
  $$(6) \qquad {\cal A}_{\tilde{M}}^{\tilde{G}}={\cal A}_{\tilde{L}_{1}}^{\tilde{G}}\oplus {\cal A}_{\tilde{L}_{2}}^{\tilde{G}},$$
  $$(7)\qquad {\cal A}_{\tilde{M}}^{\tilde{G}}={\cal A}^{\tilde{L}_{1}}_{\tilde{M}}\oplus {\cal A}^{\tilde{L}_{2}}_{\tilde{M}}.$$
  Supposons ces conditions v\'erifi\'ees. On note $d_{\tilde{M}}^{\tilde{G}}(\tilde{L}_{1},\tilde{L}_{2})$ le jacobien de l'application somme
  $${\cal A}^{\tilde{L}_{1}}_{\tilde{M}}\oplus {\cal A}^{\tilde{L}_{2}}_{\tilde{M}}\to {\cal A}_{\tilde{M}}^{\tilde{G}}$$
  qui est un isomorphisme d'apr\`es (7). De m\^eme, gr\^ace (5), on d\'efinit le jacobien $d_{\tilde{M}}^{\tilde{L}}(\tilde{M}_{1},\tilde{M}_{2})$. On v\'erifie qu'alors
  $$(8) \qquad d^{\tilde{G}}_{\tilde{M}_{1},\tilde{M}_{2}}(\tilde{L};\tilde{L}_{1},\tilde{L}_{2})=d_{\tilde{M}}^{\tilde{G}}(\tilde{L}_{1},\tilde{L}_{2})d_{\tilde{M}}^{\tilde{L}}(\tilde{M}_{1},\tilde{M}_{2})^{-1}.$$
  Fixons $H\in {\cal A}_{\tilde{M}}^{\tilde{G}}$ en position g\'en\'erale. Gr\^ace \`a (6), on peut l'\'ecrire $H=H_{1}-H_{2}$, o\`u $H_{i}\in {\cal A}_{\tilde{L}_{i}}^{\tilde{G}}$ pour $i=1,2$. Comme pr\'ec\'edemment, $H_{i}$ d\'etermine un espace parabolique $\tilde{Q}_{i}\in {\cal P}(\tilde{L}_{i})$. On v\'erifie que $(\tilde{Q}_{1},\tilde{Q}_{2})$ co\"{\i}ncide avec le couple d\'etermin\'e pr\'ec\'edemment pour un choix convenable de $\xi$.
  
 Un cas encore plus particulier est celui o\`u $\tilde{M}_{1}=\tilde{M}_{2}=\tilde{M}$. Alors $\tilde{L}=\tilde{M}$ et la relation (4) prend la forme
 $$(9) \qquad y_{\tilde{M}}^{\tilde{G}}(\Lambda)=\sum_{\tilde{L}_{1},\tilde{L}_{2}\in {\cal L}(\tilde{M})}d_{\tilde{M}}^{\tilde{G}}(\tilde{L}_{1},\tilde{L}_{2})x_{1,\tilde{M}}^{\tilde{Q}_{1}}(\Lambda)x_{2,\tilde{M}}^{\tilde{Q}_{2}}(\Lambda).$$
 Cf. [A6] corollaire 7.4.
  
  Soient maintenant $\tilde{M},\tilde{L}\in {\cal L}(\tilde{M}_{0})$ tels que $\tilde{M}\subset \tilde{L}$ et soit  $(x(\Lambda;\tilde{P}))_{\tilde{P}\in {\cal P}(\tilde{M})}$ une $(\tilde{G},\tilde{M})$-famille. On en d\'eduit une $(\tilde{G},\tilde{L})$-famille puis une fonction $x^{\tilde{G}}_{\tilde{L}}(\Lambda)$ sur $i{\cal A}_{\tilde{L}}^*$. On a l'\'egalit\'e
$$(10) \qquad x^{\tilde{G}}_{\tilde{L}}(\Lambda)=\sum_{\tilde{L}'\in {\cal L}(\tilde{M})}d_{\tilde{M}}^{\tilde{G}}(\tilde{L},\tilde{L}')x_{\tilde{M}}^{\tilde{Q}'}(\Lambda),$$
 o\`u $\tilde{Q}'$ est le second terme de l'image de couple $(\tilde{L},\tilde{L}')$ par l'application d\'ecrite ci-dessus.   
  
  \subsection{Application des formules de descente}
 
 Soit $\tilde{M}\in {\cal L}(\tilde{M}_{0})$. Pour $x\in G(F)$, on d\'efinit   la $(\tilde{G},\tilde{M})$-famille $(v(x;\Lambda,\tilde{P}))_{\tilde{P}\in {\cal P}(\tilde{M})}$, o\`u $v(x;\Lambda,\tilde{P})=e^{-<\Lambda,H_{\tilde{P}}(x)>}$.  On en d\'eduit une fonction $v^{\tilde{G}}_{\tilde{M}}(x;\Lambda)$. On pose $v^{\tilde{G}}_{\tilde{M}}(x)=v^{\tilde{G}}_{\tilde{M}}(x;0)$. Pour $f\in C_{c}^{\infty}(\tilde{G}(F))$ et $\gamma\in \tilde{M}(F)\cap \tilde{G}_{reg}(F)$ tel que $\omega$ soit trivial sur $Z_{G}(\gamma,F)$, on d\'efinit l'int\'egrale orbitale pond\'er\'ee
 $$(1) \qquad J_{\tilde{M}}^{\tilde{G}}(\gamma,\omega,f)=D^{\tilde{G}}(\gamma)^{1/2}\int_{Z_{G}(\gamma,F)\backslash G(F)}\omega(x)f(x^{-1}\gamma x)v^{\tilde{G}}_{\tilde{M}}(x)\,dx.$$
 Si $\omega$ est non trivial sur $Z_{G}(\gamma,F)$, on pose $J_{\tilde{M}}^{\tilde{G}}(\gamma,\omega,f)=0$.
 
 D'autre part, pour une $\omega$-repr\'esentation temp\'er\'ee et de longueur finie $\tilde{\pi}$ de $\tilde{M}(F)$, on a d\'efini en 2.7 le caract\`ere pond\'er\'e $J_{\tilde{M}}^{\tilde{G}}(\tilde{\pi},f)$ dans le cas o\`u $f$ est $K$-finie. Les calculs du paragraphe 5.2 montrent que cette distribution s'\'etend contin\^ument  \`a toutes les fonctions $f\in C_{c}^{\infty}(\tilde{G}(F))$.

 Soient $f\in C_{c}^{\infty}(\tilde{G}(F))$ et $\tilde{P}=\tilde{M}U_{P}\in {\cal P}(\tilde{M}_{0})$.  On d\'efinit une fonction $f_{\tilde{P}}$ sur $\tilde{M}(F)$ par la formule habituelle
 $$f_{\tilde{P}}(x)=\delta_{\tilde{P}}(x)^{1/2}\int_{U_{P}(F)}\int_{K}\omega(k)f(k^{-1}xuk)\,dk\,du.$$
 Cette fonction appartient \`a $C_{c}^{\infty}(\tilde{M}(F))$.
 
 \ass{Lemme}{(i) Soient $f_{1},f_{2}\in C_{c}^{\infty}(\tilde{G}(F))$ et $\gamma\in \tilde{M}(F)\cap \tilde{G}_{reg}(F)$. On a l'\'egalit\'e
 $$ J_{\tilde{M}}^{\tilde{G}}(\gamma,\omega,f_{1},f_{2})=\sum_{\tilde{L}_{1},\tilde{L}_{2}\in {\cal L}(\tilde{M})}d_{\tilde{M}}^{\tilde{G}}(\tilde{L}_{1},\tilde{L}_{2})\overline{J_{\tilde{M}}^{\tilde{L}_{1}}(\gamma,\omega,f_{1,\tilde{\bar{Q}}_{1}})}J_{\tilde{M}}^{\tilde{L}_{2}}(\gamma,\omega,f_{2,\tilde{Q}_{2}}).$$
 
 (ii) Soient $f_{1},f_{2}\in C_{c}^{\infty}(\tilde{G}(F))$,  $\boldsymbol{\tau}\in {\cal E}_{disc}(\tilde{M},\omega)$ et $\tilde{\lambda}\in i\tilde{{\cal A}}_{\tilde{M},F}^*$. On a l'\'egalit\'e
 $$J_{\tilde{M}}^{\tilde{G}}(\pi_{\tau_{\lambda}},f_{1},f_{2})=\sum_{\tilde{L}_{1},\tilde{L}_{2}\in {\cal L}(\tilde{M})}d_{\tilde{M}}^{\tilde{G}}(\tilde{L}_{1},\tilde{L}_{2})\overline{J_{\tilde{M}}^{\tilde{L}_{1}}(\tilde{\pi}_{\boldsymbol{\tau}_{\tilde{\lambda}}},f_{1,\tilde{\bar{Q}}_{1}})}J_{\tilde{M}}^{\tilde{L}_{2}}(\tilde{\pi}_{\boldsymbol{\tau}_{\tilde{\lambda}}},f_{2,\tilde{Q}_{2}}).$$
 
 (iii) Soient $f\in C_{c}^{\infty}(\tilde{G}(F))$, $\gamma\in \tilde{M}(F)\cap \tilde{G}_{reg}(F)$  et  $\tilde{L}\in {\cal L}(\tilde{M})$. On a l'\'egalit\'e
 $$J_{\tilde{L}}^{\tilde{G}}(\gamma,\omega,f)=\sum_{\tilde{L}'\in {\cal L}(\tilde{M})}d_{\tilde{M}}^{\tilde{G}}(\tilde{L},\tilde{L}')J_{\tilde{M}}^{\tilde{L}'}(\gamma,\omega,f_{\tilde{Q}'}).$$
 
 (iv) Soient $f\in C_{c}^{\infty}(\tilde{G}(F))$, $\tilde{\pi}$ une $\omega$-repr\'esentation temp\'er\'ee de longueur finie de $\tilde{M}(F)$ et $\tilde{L}\in {\cal L}(\tilde{M})$. Posons $\tilde{\Pi}=Ind_{\tilde{P}}^{\tilde{L}}(\tilde{\pi})$, o\`u $\tilde{P}$ est un \'el\'ement de ${\cal P}^{\tilde{L}}(\tilde{M})$. Alors on a l'\'egalit\'e
 $$J_{\tilde{L}}^{\tilde{G}}(\tilde{\Pi},f)=\sum_{\tilde{L}'\in {\cal L}(\tilde{M})}d_{\tilde{M}}^{\tilde{G}}(\tilde{L},\tilde{L}')J_{\tilde{M}}^{\tilde{L}'}(\tilde{\pi},f_{\tilde{Q}'}).$$}
 
 La preuve est bas\'ee sur les formules 5.3(9) et (10) appliqu\'ees aux $(\tilde{G},\tilde{M})$-familles intervenant dans les d\'efinitions des membres de gauche. On la laisse au lecteur.
  
 On aura besoin plus tard des propri\'et\'es suivantes. Soient $\tilde{M}\in {\cal L}(\tilde{M}_{0})$ et $\tilde{T}$ un sous-tore tordu maximal de $\tilde{M}$ tel que $\omega$ soit trivial sur $T^{\theta}(F)$. Alors
 
 (2) pour tout $f\in C_{c}^{\infty}(\tilde{G}(F))$, la fonction $\gamma\mapsto J_{\tilde{M}}^{\tilde{G}}(\gamma,\omega,f)$ est lisse sur $\tilde{T}(F)\cap \tilde{G}_{reg}(F)$.
 
 En effet, au voisinage d'un point r\'egulier, $J_{\tilde{M}}^{\tilde{G}}(\gamma,\omega,f)$ est l'int\'egrale sur un compact d'une fonction lisse en $\gamma$ et en la variable d'int\'egration.

 (3) il existe un entier $N\geq0$ et, pour tout $f\in C_{c}^{\infty}(\tilde{G}(F))$, il existe $c>0$ tel que $\vert J_{\tilde{M}}^{\tilde{G}}(\gamma,\omega,f)\vert \leq c (1+\vert  log(D^{\tilde{G}}(\gamma))\vert )^N$.
 
 Preuve. Comme on l'a d\'ej\`a utilis\'e plusieurs fois, il r\'esulte de la variante tordue de 1.7(1) et de 4.2(7) qu'il existe $N\geq0$ ind\'ependant de $f$ et $c'>0$ tel que, pour $x$ contribuant \`a la formule (1), on ait $v^{\tilde{G}}_{\tilde{M}}(x)\leq c' (1+\vert  log(D^{\tilde{G}}(\gamma))\vert )^N$. Alors
 $$\vert J_{\tilde{M}}^{\tilde{G}}(\gamma,\omega,f)\vert \leq c' (1+\vert  log(D^{\tilde{G}}(\gamma))\vert )^ND^{\tilde{G}}(\gamma)^{1/2}\int_{Z_{G}(\gamma,F)\backslash G(F)}\vert f(x^{-1}\gamma x)\vert \,dx.$$
 Il reste \`a appliquer 4.2(1) pour obtenir (3).  $\square$

 \bigskip
 
 \subsection{Le th\'eor\`eme $0$}
 
 \ass{Th\'eor\`eme}{Soit $f\in C_{c}^{\infty}(\tilde{G}(F))$. Supposons $I_{\tilde{G}}(\tilde{\pi},f)=0$ pour toute $\omega$-repr\'esentation irr\'eductible temp\'er\'ee $\tilde{\pi}$ de $\tilde{G}(F)$. Alors $I_{\tilde{G}}(\gamma,\omega,f)=0$ pour tout $\gamma\in \tilde{G}_{reg}(F)$.}
 
 {\bf Remarques.} (1) On n'a d\'efini dans cet article que les int\'egrales orbitales relatives \`a des \'el\'ements $\gamma\in \tilde{G}_{reg}(F)$. Mais on sait bien qu'on peut les d\'efinir pour tout $\gamma\in \tilde{G}(F)$.   Par la th\'eorie de la descente, on montre que leur comportement local est le m\^eme que dans le cas non tordu. Cest-\`a-dire que les int\'egrales orbitales aux points singuliers s'expriment \`a l'aide d'int\'egrales orbitales aux points r\'eguliers, soit gr\^ace aux germes de Shalika dans le cas non-archim\'edien, soit par l'action d'op\'erateurs diff\'erentiels dans le cas archim\'edien. On peut donc renforcer la conclusion du th\'eor\`eme: $I_{\tilde{G}}(\gamma,\omega,f)=0$ pour tout $\gamma\in \tilde{G}(F)$ tel que $\omega$ soit trivial sur $Z_{G}(\gamma,F)$.
 
 (2) Dans le cas non-archim\'edien, on montre comme dans le cas non tordu qu'une fonction $f\in C_{c}^{\infty}(\tilde{G}(F))$ dont toutes les int\'egrales orbitales sont nulles est annul\'ee par toute forme lin\'eaire $l$ sur $C_{c}^{\infty}(\tilde{G}(F))$ qui v\'erifie la relation $l(f^x)=\omega(x)l(f)$ pour tout $f\in C_{c}^{\infty}(\tilde{G}(F))$ et tout $x\in G(F)$, o\`u $f^x$ est la fonction $f^x(\gamma)=f(x\gamma x^{-1})$.
 
 (3) Dans le cas non-archim\'edien, le th\'eor\`eme a \'et\'e prouv\'e dans [HL]. 
 
 \bigskip
 
 Preuve. On raisonne par r\'ecurence sur la dimension de $\tilde{G}$. Soit $\tilde{M} $ un espace
 de Levi semi-standard propre. Fixons $\tilde{P}\in {\cal P}(\tilde{M})$. Pour une $\omega$-repr\'esentation admisslble $\tilde{\pi}$ de $\tilde{M}(F)$, on a l'\'egalit\'e $I_{\tilde{M}}(\tilde{\pi},f_{\tilde{P}})=I_{\tilde{G}}(\tilde{\Pi},f)$, o\`u $\Pi=Ind_{\tilde{P}}^{\tilde{G}}(\tilde{\pi})$. La fonction $f_{\tilde{P}}$ v\'erifie donc la m\^eme hypoth\`ese que $f$. Donc, par l'hypoth\`ese de r\'ecurrence, $I_{\tilde{M}}(\gamma,\omega,f_{\tilde{P}})=0$ pour tout $\gamma\in \tilde{M}_{reg}(F)$.
Mais, pour $\gamma\in \tilde{M}(F)\cap \tilde{G}_{reg}(F)$, on a l'\'egalit\'e $I_{\tilde{G}}(\gamma,\omega,f)=I_{\tilde{M}}(\gamma,\omega,f_{\tilde{P}})$. Donc $I_{\tilde{G}}(\gamma,\omega,f)=0$. Cela prouve le r\'esultat pour tout $\gamma\in \tilde{G}_{reg}(F)$ qui n'est pas elliptique.  Pour traiter le cas des \'el\'ements elliptiques, consid\'erons une fonction $f_{2}\in C_{c}^{\infty}(\tilde{G}(F))$  \`a support dans l'ensemble $\tilde{G}_{ell}(F)$ des \'el\'ements elliptiques r\'eguliers de $\tilde{G}(F)$.  Soient $\tilde{M}$ un espace de Levi de $\tilde{G}$, $\tau\in E_{disc}(\tilde{M},\omega)$ et $\lambda\in i{\cal A}_{\tilde{M},F}^*$.  On rel\`eve $\tau$ en $\boldsymbol{\tau}\in {\cal E}_{disc}(\tilde{M},\omega)$ et $\lambda$ en $\tilde{\lambda}\in i\tilde{{\cal A}}_{\tilde{M},F}^*$. Le terme $J_{\tilde{M}}^{\tilde{G}}(\pi_{\tau_{\lambda}},f,f_{2})$ est calcul\'e par le lemme 5.4(ii). Si $\tilde{Q}_{2}$ est un espace parabolique propre de $\tilde{G}$, l'hypoth\`ese sur le support de $f_{2}$ entra\^{\i}ne que $f_{2,\tilde{Q}_{2}}=0$. La formule se simplifie donc en
$$J_{\tilde{M}}^{\tilde{G}}(\pi_{\tau_{\lambda}},f,f_{2})=\overline{J_{\tilde{M}}^{\tilde{M}}(\tilde{\pi}_{\boldsymbol{\tau}_{\tilde{\lambda}}},f_{\tilde{\bar{Q}}})}J_{\tilde{M}}^{\tilde{G}}(\tilde{\pi}_{\boldsymbol{\tau}_{\tilde{\lambda}}},f_{2}),$$
o\`u $\tilde{Q}$ est un certain \'el\'ement de ${\cal P}(\tilde{M})$. On a 
$$J_{\tilde{M}}^{\tilde{M}}(\tilde{\pi}_{\boldsymbol{\tau}_{\tilde{\lambda}}},f_{\tilde{\bar{Q}}})=I_{\tilde{M}}^{\tilde{M}}(\tilde{\pi}_{\boldsymbol{\tau}_{\tilde{\lambda}}},f_{\tilde{\bar{Q}}})=I_{\tilde{G}}(\tilde{\Pi},f),$$
o\`u $\tilde{\Pi}=Ind_{\tilde{\bar{Q}}}^{\tilde{G}}(\tilde{\pi}_{\boldsymbol{\tau}_{\tilde{\lambda}}})$. La repr\'esentation $\tilde{\Pi}$ est temp\'er\'ee, donc $I_{\tilde{G}}(\tilde{\Pi},f)=0$. Cela prouve que $J^{\tilde{G}}_{spec}(\omega,f,f_{2})=0$. Appliquons la proposition 5.2: elle entra\^{\i}ne $J^{\tilde{G}}_{g\acute{e}om}(\omega,f,f_{2})=0$. 
Pour un espace de Levi propre $\tilde{M}$ et un \'el\'ement $\tilde{S}\in T_{ell}(\tilde{M},\omega) $, on a $J_{\tilde{M},\tilde{S}}(\omega,f,f_{2})=0$: l'hypoth\`ese sur le support de $f_{2}$ entra\^{\i}ne que $f_{2}(y^{-1}\gamma y)=0$ pour tout $y\in G(F)$ et tout $\gamma\in \tilde{S}(F)$. On a donc
$$J_{g\acute{e}om}^{\tilde{G}}(\omega,f,f_{2})=\sum_{\tilde{S}\in T_{ell}(\tilde{G},\omega)}\vert W^G(\tilde{S})\vert ^{-1}mes(A_{\tilde{G}}(F)\backslash S^{\theta}(F))$$
$$\int_{\tilde{S}(F)/(1-\theta)(S(F))}\overline{I_{\tilde{G}}(\gamma,\omega,f)}I_{\tilde{G}}(\gamma,\omega,f_{2})\,d\gamma.$$
Soit maintenant $\gamma\in \tilde{G}_{ell}(F)$ tel que $\omega$ soit trivial sur $Z_{G}(\gamma,F)$. On veut montrer que $I_{\tilde{G}}(\gamma,\omega,f)=0$.   On peut conjuguer $\gamma$ et supposer que $\gamma\in \tilde{S}(F)$ pour un $\tilde{S}$ intervenant dans la formule ci-dessus. On fait maintenant parcourir \`a $f_{2}$ une suite de fonctions  $(f_{2,n})_{n\in {\mathbb N}}$ \`a valeurs positives ou nulles, non nulles en $\gamma$, et \`a supports dans des voisinages de plus en plus petits de $\gamma$. On voit que l'expression ci-dessus est de la forme
$$\int_{S^{\theta,0}(F)}\overline{I_{\tilde{G}}(x\gamma,\omega,f)}A_{n}(x)\,dx$$
o\`u $A_{n}$ parcourt une suite de fonctions  lisses sur $S^{\theta,0}(F))$  \`a valeurs positives ou nulles, non nulles en $x=1$, et \`a supports dans des voisinages de plus en plus petits de $1$. Puisque $x\mapsto I_{\tilde{G}}(x\gamma,\omega,f)$ est lisse, la nullit\'e de l'expression ci-dessus pour tout $n$ entra\^{\i}ne $I_{\tilde{G}}(\gamma,\omega,f)=0$. $\square$

  \bigskip
 
 \section{La formule invariante}
 
 \bigskip
 
 \subsection{Le th\'eor\`eme de Paley-Wiener}
 Notons ${\cal K}$ le groupe de Grothendieck de l'ensemble des classes d'isomorphisme de $\omega$-repr\'esentations temp\'er\'ees   de longueur finie de $\tilde{G}(F)$. 
 Notons ${\cal F}$ l'espace des fonctions lin\'eaires de ${\cal K}$ dans ${\mathbb C}$ qui sont nulles sur toute $\omega$-repr\'esentation temp\'er\'ee irr\'eductible et non $G$-irr\'eductible. Remarquons que ${\cal F}$ s'identifie \`a l'espace des fonctions \`a valeurs complexes sur   l'ensemble des classes d'isomorphisme de $\omega$-repr\'esentations  temp\'er\'ees $G$-irr\'eductibles de $\tilde{G}(F)$. Notons ${\cal PW}(\tilde{G},\omega)$ le sous-ensemble des \'el\'ements $\varphi\in {\cal F}$ v\'erifiant les conditions  suivantes:
 
 (1)  soient $\tilde{Q}=\tilde{L}U_{Q}$ un espace parabolique tel que $\omega$ soit trivial sur $Z_{L}(F)^{\theta}$ et $\tilde{\pi}$ une $\omega$-repr\'esentation temp\'er\'ee et $L$-irr\'eductible de $\tilde{L}(F)$;  alors la fonction $\tilde{\lambda}\mapsto \varphi(Ind_{\tilde{Q}}^{\tilde{G}}(\tilde{\pi}_{\tilde{\lambda}}))$ sur $i\tilde{{\cal A}}_{\tilde{L},F}^*$ est de Paley-Wiener (cf. 2.6);
 
  (2) si $F$ est non-archim\'edien, il existe un sous-groupe ouvert compact $H$ de $\tilde{G}(F)$ tel que $\varphi(\tilde{\pi})=0$ pour toute $\omega$-repr\'esentation temp\'er\'ee $\tilde{\pi}$ telle que le sous-espace des  invariants par $H$ dans $\pi$ soit nul;
 
 (3) si $F$ est archim\'edien, il existe un ensemble fini $\kappa_{1},...,\kappa_{n}$ de repr\'esentations irr\'eductibles de $K$ tel que $\varphi(\tilde{\pi})=0$ pour toute $\omega$-repr\'esentation temp\'er\'ee $\tilde{\pi}$ telle que, pour tout $i$, l'espace isotypique de type $\kappa_{i}$ dans $\pi$ soit nul.

 Notons $\underline{pw}_{\tilde{G}}:C_{c}^{\infty}(\tilde{G}(F),K)\to {\cal F}$ l'application lin\'eaire qui, \`a $f\in C_{c}^{\infty}(\tilde{G}(F),K)$, associe la fonction $\tilde{\pi}\mapsto I_{\tilde{G}}(\tilde{\pi},f)$. Notons $I(\tilde{G}(F),K,\omega)$ le quotient de $C_{c}^{\infty}(\tilde{G}(F),K)$ par le sous-espace des $f\in C_{c}^{\infty}(\tilde{G}(F),K)$ telles que $I_{\tilde{G}}(\gamma,\omega,f)=0$ pour tout $\gamma\in \tilde{G}_{reg}(F)$.
 
 \ass{Th\'eor\`eme}{L'application lin\'eaire $\underline{pw}_{\tilde{G}}$  a pour image ${\cal PW}(\tilde{G},\omega)$ et se quotiente en un isomorphisme de $I(\tilde{G}(F),K,\omega)$ sur cet espace.}
 
 Quand $F$ est non-archim\'edien, le th\'eor\`eme a \'et\'e d\'emontr\'e dans [R] pour $\omega=1$ et en g\'en\'eral dans [HL]. Si $F={\mathbb R}$, la premi\`ere assertion du th\'eor\`eme est d\'emontr\'ee dans [DM] pour $\omega=1$. Nous montrerons en 6.4 que le cas $\omega\not=1$ s'en d\'eduit. Si $F={\mathbb C}$, l'assertion est \'equivalente \`a celle pour l'espace tordu sur ${\mathbb R}$ d\'eduit  de $\tilde{G}$ par restriction des scalaires. La deuxi\`eme assertion est facile: le th\'eor\`eme 5.5 montre que tout \'el\'ement du noyau  a des int\'egrales orbitales nulles; la r\'eciproque provient de la locale int\'egrabilit\'e des caract\`eres des $\omega$-repr\'esentations $G$-irr\'eductibles, cf. 2.5(2). 
 
 {\bf Compl\'ement au th\'eor\`eme.} Supposons $F$ non-archim\'edien. Soit $H$ un sous-groupe ouvert compact de $G(F)$. Notons ${\cal PW}(\tilde{G},\omega)^H$ le sous-espace des fonctions $\varphi\in {\cal PW}(\tilde{G},\omega)$ tels que $\varphi(\tilde{\pi})=0$ pour toute $\omega$-repr\'esentation temp\'er\'ee $\tilde{\pi}$ telle que le sous-espace des  invariants par $H$ dans $\pi$ soit nul. Notons $C_{c}^{\infty}(H\backslash \tilde{G}(F)/H)$ le sous-espace des \'el\'ements de $C_{c}^{\infty}(\tilde{G}(F))$ qui sont biinvariants par $H$. Alors, $H$ \'etant donn\'e, il existe $H'$ tel que ${\cal PW}(\tilde{G},\omega)^H$ soit contenu dans $\underline{pw}_{\tilde{G}}(C_{c}^{\infty}(H'\backslash\tilde{G}(F)/H'))$. Cf. [HL].
 
 Supposons $F$ archim\'edien. Soit $\underline{\kappa}$ un ensemble fini de repr\'esentations irr\'eductibles de $K$. Notons ${\cal PW}(\tilde{G},\omega)^{\underline{\kappa}}$ le sous-espace des fonctions $\varphi\in PW(\tilde{G},\omega)$ tels que $\varphi(\tilde{\pi})=0$ pour toute $\omega$-repr\'esentation temp\'er\'ee $\tilde{\pi}$ telle que, pour tout $\kappa\in \underline{\kappa}$, l'espace isotypique de type $\kappa$ dans $\pi$ soit nul. Notons $C_{c}^{\infty}(\tilde{G}(F),\underline{\kappa})$ le sous-espace des $f\in C_{c}^{\infty}(\tilde{G}(F))$ tels que la repr\'esentation de $K\times K$ dans l'espace engendr\'e par les translat\'es de $f$ \`a droite et \`a gauche par des \'el\'ements de $K$ n'ait pour composantes irr\'eductibles que des \'el\'ements de $\underline{\kappa}\times \underline{\kappa}$. Alors, $\underline{\kappa}$ \'etant donn\'e, il existe $\underline{\kappa}'$ tel que ${\cal PW}(\tilde{G}(F))^{\underline{\kappa}}$ soit contenu dans $\underline{pw}_{\tilde{G}}(C_{c}^{\infty}(\tilde{G}(F),\underline{\kappa}'))$. Ceci n'est pas \'enonc\'e dans [DM], mais r\'esulte clairement de la preuve.
 
 \bigskip
 
 \subsection{Deuxi\`eme forme du th\'eor\`eme de Paley-Wiener}
 On rappelle que l'on a d\'efini l'ensemble ${\cal E}(\tilde{G},\omega)$ en 2.9. Il est muni d'une action de  $i\tilde{{\cal A}}_{\tilde{G},F}^*$ qui, \`a $\boldsymbol{\tau}\in {\cal E}(\tilde{G},\omega)$ et $\tilde{\lambda}\in i\tilde{{\cal A}}_{\tilde{G},F}^*$, associe $\boldsymbol{\tau}_{\tilde{\lambda}}$.  Il est aussi muni d'une action de ${\mathbb U}$ qui, \`a $\boldsymbol{\tau}=(M,\sigma,\boldsymbol{\tilde{r}})$ et $z\in {\mathbb U}$, associe $z\boldsymbol{\tau}=(M,\sigma,z\boldsymbol{\tilde{r}})$. Rappelons que l'action de $z\in i{\mathbb R}/2\pi i{\mathbb Z}\in i\tilde{{\cal A}}_{\tilde{G},F}^*$ co\"{\i}ncide avec celle de $e^z\in {\mathbb U}$. Notons ${\cal E}_{ell}(\tilde{G},\omega)$ l'ensemble des    triplets $(M,\sigma,\boldsymbol{\tilde{r}})\in {\cal E}(\tilde{G},\omega)$ tels que $(M,\sigma,\tilde{r})\in E_{ell}(\tilde{G},\omega)$, o\`u $\tilde{r}$ est l'image de $\boldsymbol{\tilde{r}}$ dans $R^{\tilde{G}}(\sigma)$. Notons ${\cal E}_{ell}(\tilde{G},\omega)/conj$ l'ensemble des classes de conjugaison par $G(F)$ dans ${\cal E}_{ell}(\tilde{G},\omega)$.     On a une surjection $({\cal E}_{ell}(\tilde{G},\omega)/conj)\to (E_{ell}(\tilde{G},\omega)/conj)$ dont les fibres sont toutes isomorphes \`a ${\mathbb U}$. Notons $PW_{ell}(\tilde{G},\omega)$ l'espace des fonctions $\varphi:({\cal E}_{ell}(\tilde{G},\omega)/conj)\to {\mathbb C}$  qui v\'erifient les conditions suivantes:

 (1) soient $\boldsymbol{\tau}\in {\cal E}_{ell}(\tilde{G},\omega)$ et $z\in {\mathbb U}$; alors $\varphi(z\boldsymbol{\tau})=z\varphi(\boldsymbol{\tau})$ (via la projection ${\cal E}_{ell}(\tilde{G},\omega)\to ({\cal E}_{ell}(\tilde{G},\omega)/conj)$, on a identifi\'e $\varphi$ \`a une fonction sur ${\cal E}_{ell}(\tilde{G},\omega)$);
 
 (2) soit $\boldsymbol{\tau}=(M,\sigma,\boldsymbol{\tilde{r}})\in {\cal E}_{ell}(\tilde{G},\omega)$;   alors la fonction $\tilde{\lambda}\mapsto \varphi(\boldsymbol{\tau}_{\tilde{\lambda}})$ est de Paley-Wiener sur $i\tilde{{\cal A}}_{\tilde{G},F}^*$;
 
 (3) le support de $\varphi$ est contenu dans un nombre fini d'orbites pour l'action de $i\tilde{{\cal A}}_{\tilde{G},F}^*$ dans ${\cal E}_{ell}(\tilde{G},\omega)/conj$.
 
 Rappelons (cf. 2.12) que l'on note ${\cal L}(\tilde{M}_{0},\omega)$ l'ensemble des $\tilde{L}\in {\cal L}(\tilde{M}_{0})$ tels que $\omega$ soit trivial sur $Z_{L}(F)^{\theta}$. Soient $\tilde{L}$, $\tilde{L}'$ deux  \'el\'ements de ${\cal L}(\tilde{M}_{0},\omega)$ et soit $x\in G(F)$ tel que $x\tilde{L}x^{-1}=\tilde{L}'$. De la conjugaison par $x$ se d\'eduit un isomorphisme $T_{ell}(\tilde{L},\omega)\simeq T_{ell}(\tilde{L}',\omega)$ (cf. 2.12), puis un isomorphisme $PW_{ell}(\tilde{L},\omega)\simeq PW_{ell}(\tilde{L}',\omega)$.  Notons que, si $\tilde{L}=\tilde{L}'$ et $x\in L(F)$, cet isomorphisme est l'identit\'e. En particulier, si on note $W^G(\tilde{L})$ le quotient par $L(F)$ du normalisateur de $\tilde{L}$ dans $G(F)$, ce groupe $W^G(\tilde{L})$ agit sur $PW_{ell}(\tilde{L})$.  On pose
 $$PW(\tilde{G},\omega)= (\oplus_{\tilde{L}\in {\cal L}(\tilde{M}_{0},\omega)}PW_{ell}(\tilde{L},\omega))^{\tilde{W}^G}$$
$$= \oplus_{\tilde{L}\in {\cal L}(\tilde{M}_{0},\omega)/\tilde{W}^G}PW_{ell}(\tilde{L},\omega)^{W^G(\tilde{L})},$$
 les exposants signifiant selon l'usage que l'on prend les invariants.
 
 {\bf Remarque.} Soient $\tilde{L}\in {\cal L}(\tilde{M}_{0},\omega)$, $\varphi\in PW_{ell}(\tilde{L},\omega)^{W^G(\tilde{L})}$ et $\boldsymbol{\tau}\in {\cal E}_{ell}(\tilde{L},\omega)$. Le triplet $\boldsymbol{\tau}$ est par d\'efinition essentiel dans $\tilde{L}$ mais on a d\'ej\`a dit qu'il ne l'\'etait pas forc\'ement dans $\tilde{G}$. S'il ne l'est pas, il existe $w\in W^G(\tilde{L})$ et $z\in {\mathbb U}$ tels que $w(\boldsymbol{\tau})=z\tau$, et $z\not=1$ (cf. preuve de la proposition 2.12). Cela entra\^{\i}ne $\varphi(\boldsymbol{\tau})=0$.
 \bigskip
 
 Soit $\psi\in {\cal PW}(\tilde{G},\omega)$. Pour un espace de Levi $\tilde{L}\in {\cal L}(\tilde{M}_{0},\omega)$ et un \'el\'ement $\boldsymbol{\tau}\in {\cal E}_{ell}(\tilde{L},\omega)$, rappelons que l'on a d\'efini en 2.9 une $\omega$-repr\'esentation $\tilde{\pi}_{\boldsymbol{\tau}}$ de $\tilde{L}(F)$. Fixons un espace parabolique $\tilde{Q}\in {\cal P}(\tilde{L})$, posons $\tilde{\Pi}_{\boldsymbol{\tau}}=Ind_{\tilde{Q}}^{\tilde{G}}(\tilde{\pi}_{\boldsymbol{\tau}})$. Cette repr\'esentation est temp\'er\'ee et sa classe d'isomorphisme ne d\'epend pas de $\tilde{Q}$. Posons $\varphi(\boldsymbol{\tau})=\psi(\tilde{\Pi}_{\boldsymbol{\tau}})$. On a ainsi une fonction $\varphi$ sur $\oplus_{\tilde{L}\in {\cal L}(\tilde{M}_{0},\omega)}T_{ell}(\tilde{L},\omega)$. Cette fonction appartient \`a $PW(\tilde{G},\omega)$. En effet,  les conditions (1) et (2) r\'esultent de 6.1(1) tandis que la condition  (3) r\'esulte de 6.1(2) et (3). La condition d'invariance par $\tilde{W}^G$ r\'esulte de l'\'egalit\'e $\tilde{\Pi}_{w(\boldsymbol{\tau})}=\tilde{\Pi}_{\boldsymbol{\tau}}$ pour tout $w\in\tilde{W}^G$. Notons 
 $$res:{\cal PW}(\tilde{G},\omega)\to PW(\tilde{G},\omega)$$
 l'application $\psi\mapsto \varphi$ et
 $$pw_{\tilde{G}}:C_{c}^{\infty}(\tilde{G}(F),K)\to PW(\tilde{G},\omega)$$
 la compos\'ee $pw_{\tilde{G}}=res\circ \underline{pw}_{\tilde{G}}$.

 \ass{Th\'eor\`eme}{L'application $pw_{\tilde{G}}$  se quotiente en un isomorphisme de $I(\tilde{G}(F),K,\omega)$ sur $PW(\tilde{G},\omega)$.}
 
 Preuve.  L'\'enonc\'e r\'esulte du th\'eor\`eme 6.1 pouvu que $res$ soit bijectif.   Cette bijectivit\'e r\'esulte ais\'ement de la proposition 2.12. $\square$
 
 Le th\'eor\`eme admet un compl\'ement similaire \`a celui du th\'eor\`eme 6.1.
 
 \bigskip
 
 \subsection{Extension du th\'eor\`eme de Delorme et Mezo au cas $\omega\not=1$}
 Dans ce paragraphe,  on suppose $F={\mathbb R}$. On veut prouver 
 
 (1) l'application $\underline{pw}_{\tilde{G}}$ a pour image ${\cal PW}(\tilde{G},\omega)$.  
 
  Delorme et Mezo ont prouv\'e l'assertion pour $\omega={\bf 1}$, o\`u ${\bf 1}$ d\'esigne le caract\`ere trivial de $G({\mathbb R})$. Supposons d'abord qu'il existe un caract\`ere unitaire $\mu$ de $G({\mathbb R})$ tel que $\omega=\mu\circ(1-\theta)$, cf. 2.4. Fixons un \'el\'ement $\gamma_{0}\in \tilde{G}({\mathbb R})$. Comme dans la preuve de 2.5(2), on associe \`a toute $\omega$-repr\'esentation $\tilde{\pi}$ de $\tilde{G}({\mathbb R})$ une ${\bf 1}$-repr\'esentation $\tilde{\pi}_{1}$ d\'efinie par $\tilde{\pi}_{1}(g\gamma_{0})=\mu(g)\tilde{\pi}(\gamma_{0})$ pour tout $g\in G({\mathbb R})$. L'application $\tilde{\pi}\mapsto \tilde{\pi}_{1}$ est bijective. Pour plus de pr\'ecision, on ajoute des indices  $\omega$ \`a certains objets d\'efinis en 6.1, par exemple ${\cal F}_{\omega}$ et $\underline{pw}_{\tilde{G},\omega }$. On d\'efinit une application 
 $$\begin{array}{ccc}{\cal F}_{\omega}&\to&{\cal F}_{{\bf 1}}\\ \varphi&\mapsto &\varphi_{1}\\ \end{array}$$
 par $ \varphi_{1}(\tilde{\pi}_{1})=\varphi(\tilde{\pi})$.  On d\'efinit une application
 $$\begin{array}{ccc}C_{c}^{\infty}(\tilde{G}({\mathbb R}),K)&\to&C_{c}^{\infty}(\tilde{G}({\mathbb R}),K)\\ f&\mapsto& f_{1} \\ \end{array}$$
 par $f_{1}(g\gamma_{0})=\mu(g)^{-1}f(g\gamma_{0})$ pour tout $g\in G({\mathbb R})$. Il est clair que $I_{\tilde{G}}(\tilde{\pi}_{1},f_{1})=I_{\tilde{G}}(\tilde{\pi},f)$ pour tous $f$, $\tilde{\pi}$. Donc le diagramme suivant est commutatif:
 $$\begin{array}{ccc}C_{c}^{\infty}(\tilde{G}({\mathbb R}),K)&\stackrel{f\mapsto f_{1}}{\to}&C_{c}^{\infty}(\tilde{G}({\mathbb R}),K)\\ \downarrow \underline{pw}_{\tilde{G},\omega}&&\downarrow \underline{pw}_{\tilde{G},{\bf 1}}\\ {\cal F}_{\omega}&\stackrel {\varphi\mapsto \varphi_{1}}{\to}&{\cal F}_{{\bf 1}}.\\ \end{array}$$
 L'application horizontale du haut est bijective. On v\'erifie imm\'ediatement que celle du bas se restreint en un isomorphisme de ${\cal PW}(\tilde{G},\omega)$ sur ${\cal PW}(\tilde{G},{\bf 1})$. Alors l'assertion (1) se d\'eduit de la m\^eme assertion pour le caract\`ere ${\bf 1}$.
 
 Dans le cas g\'en\'eral, on introduit des objets $G'$, $\tilde{G}'$, $C$, $p$, $\tilde{p}$ v\'erifiant la proposition 2.4. On pose $\omega'=\omega\circ p$. Ces termes v\'erifient l'hypoth\`ese pr\'ec\'edente: il existe un caract\`ere $\mu'$ de $G'({\mathbb R})$ tel que $\omega'=\mu'\circ(1-\theta')$. On note $K'$ l'unique sous-groupe compact maximal de $G'({\mathbb R})$ contenu dans $p^{-1}(K)$. Pour une $\omega$-repr\'esentation $\tilde{\pi}$ de $\tilde{G}({\mathbb R})$, on note $\tilde{\pi}'$ la $\omega'$-repr\'esentation $\tilde{\pi}\circ\tilde{p}$ de $\tilde{G}'({\mathbb R})$. Dualement, on en d\'eduit une application
 $$(2) \qquad \begin{array}{ccc}{\cal F}_{\tilde{G}',\omega'}&\to&{\cal F}_{\tilde{G},\omega}\\ \varphi'&\mapsto&\varphi\\ \end{array}$$
 par $\varphi(\tilde{\pi})=\varphi'(\tilde{\pi}')$ (on a  ajout\'e des indices $\tilde{G}$ et $\tilde{G}'$ aux notations pr\'ec\'edentes). On d\'efinit aussi une application
 $$\begin{array}{ccc}C_{c}^{\infty}(\tilde{G}'({\mathbb R}),K')&\to&C_{c}^{\infty}(\tilde{G}({\mathbb R}),K)\\ f'&\mapsto&f\\ \end{array}$$
 par 
 $$f(\gamma)=\int_{C({\mathbb R})}f'(c\gamma')\,dc,$$
 pour tout $\gamma\in \tilde{G}({\mathbb R})$, o\`u $\gamma'$ est un rel\`evement quelconque de $\gamma$ dans $\tilde{G}'({\mathbb R})$. La mesure sur $C(F)$ doit \^etre compatible aux mesures choisies sur $G({\mathbb R})$ et $G'({\mathbb R})$. On v\'erifie que le diagramme suivant est commutatif
  $$\begin{array}{ccc}C_{c}^{\infty}(\tilde{G}'({\mathbb R}),K')&\stackrel{f'\mapsto f}{\to}&C_{c}^{\infty}(\tilde{G}({\mathbb R}),K)\\ \downarrow \underline{pw}_{\tilde{G}',\omega'}&&\downarrow \underline{pw}_{\tilde{G},\omega}\\ {\cal F}_{\tilde{G}',\omega'}&\stackrel {\varphi'\mapsto \varphi}{\to}&{\cal F}_{\tilde{G},\omega}.\\ \end{array}$$
  L'application horizontale du haut est surjective. Puisque l'application $\underline{pw}_{\tilde{G}',\omega'}$ v\'erifie (1) d'apr\`es le cas d\'ej\`a trait\'e, il nous suffit pour conclure de montrer que l'application (2) se restreint en une surjection de ${\cal PW}(\tilde{G}',\omega')$ sur ${\cal PW}(\tilde{G},\omega)$. Que l'image de ${\cal PW}(\tilde{G}',\omega')$ par l'application (2) soit contenue dans ${\cal PW}(\tilde{G},\omega)$ r\'esulte ais\'ement du fait suivant. Soit $\tilde{Q}=\tilde{L}U_{Q}$ un espace parabolique de $\tilde{G}$, notons $\tilde{Q}'=\tilde{L}'U_{Q'}$ son image r\'eciproque dans $\tilde{G}'$. L'espace $i\tilde{{\cal A}}_{\tilde{L}}^*$ s'injecte naturellement dans $i\tilde{{\cal A}}_{\tilde{L}'}^*$ et la restriction \`a ce sous-espace d'une fonction de Paley-Wiener sur  $i\tilde{{\cal A}}_{\tilde{L}'}^*$ est encore de Paley-Wiener. Pour d\'emontrer la surjectivit\'e, on va construire une section $s$ de l'application (2) et montrer que $s$ envoie ${\cal PW}(\tilde{G},\omega)$ dans ${\cal PW}(\tilde{G}',\omega')$. On a une suite exacte naturelle
  $$0\to i\tilde{{\cal A}}_{\tilde{G}}^*\to i\tilde{{\cal A}}_{\tilde{G}'}^*\to i{\cal A}_{C}^*\to 0.$$
 En choisissant une section de la derni\`ere application, on identifie $i{\cal A}_{C}^*$ \`a  un suppl\'ementaire  de $i\tilde{{\cal A}}_{\tilde{G}}^*$ dans $i\tilde{{\cal A}}_{\tilde{G}'}^*$. On fixe une fonction $\varphi_{C}$ de Paley-Wiener sur $i{\cal A}_{C}^*$ telle que $\varphi_{C}(0)=1$. Soit $\varphi\in {\cal F}_{\tilde{G},\omega}$. Pour une $\omega'$-repr\'esentation $\tilde{\pi}'$ de $\tilde{G}'({\mathbb R})$, temp\'er\'ee et $G'$-irr\'eductible, ou bien il n'existe aucun $\tilde{\xi}\in i{\cal A}_{C}^*$ tel que $\tilde{\pi}'_{\tilde{\xi}}$ se factorise par $\tilde{p}$. On pose alors $(s(\varphi))(\tilde{\pi}')=0$. Ou bien il existe un unique $\tilde{\xi}\in i{\cal A}_{C}^*$ tel que $\tilde{\pi}'_{\tilde{\xi}}$ se factorise par $\tilde{p}$. Pour ce $\tilde{\xi}$, on note $\tilde{\pi}$ la $\omega$-repr\'esentation de $\tilde{G}({\mathbb R})$ telle que $\tilde{\pi}'_{\tilde{\xi}}=\tilde{\pi}\circ\tilde{p}$ et on pose alors $(s(\varphi))(\tilde{\pi}')=\varphi_{C}(\tilde{\xi})\varphi(\tilde{\pi})$. La fonction $s(\varphi)$ ainsi d\'efinie appartient \`a ${\cal F}_{\tilde{G}',\omega'}$. Cela d\'efinit une application $s:{\cal F}_{\tilde{G},\omega}\to {\cal F}_{\tilde{G}',\omega'}$ qui est clairement une section de l'application (2). On doit montrer que, si $\varphi\in {\cal PW}(\tilde{G},\omega)$, alors $s(\varphi)\in {\cal PW}(\tilde{G}',\omega')$. La condition 6.1(3) est imm\'ediate. Soient $\tilde{Q}'=\tilde{L}'U_{Q'}$ un espace parabolique et $\tilde{\sigma}'$ une $\omega'$-repr\'esentation temp\'er\'ee et $L'$-irr\'eductible de $\tilde{L}'(F)$. Posons $\tilde{Q}=\tilde{p}(\tilde{Q}')$, $\tilde{L}=\tilde{p}(\tilde{L}')$.  Supposons d'abord que, pour tout $\tilde{\lambda}'\in i\tilde{{\cal A}}_{\tilde{L}'}^*$ et pour toute composante $G'$-irr\'eductible $\tilde{\pi}'$ de $Ind_{\tilde{Q}'}^{\tilde{G}'}(\tilde{\sigma}'_{\tilde{\lambda}'})$, il n'existe pas de  $\tilde{\xi}\in i{\cal A}_{C}^*$ tel que $\tilde{\pi}'_{\tilde{\xi}}$ se factorise par $\tilde{p}$. Alors par d\'efinition $(s(\varphi))(Ind_{\tilde{Q}'}^{\tilde{G}'}(\tilde{\sigma}'_{\tilde{\lambda}'}))=0$ pour tout $\tilde{\lambda}'$ et cette fonction de $\tilde{\lambda}'$ est bien de Paley-Wiener. Supposons au contraire qu'il existe $\tilde{\lambda}'$, $\tilde{\pi}'$ et $\tilde{\xi}$ comme ci-dessus tel que $\tilde{\pi}'_{\tilde{\xi}}$ se factorise par $\tilde{p}$. Fixons de tels objets. Quitte \`a remplacer $\tilde{\sigma}'$ par $\tilde{\sigma}'_{\tilde{\lambda}'+\tilde{\xi}}$, on peut supposer que $\tilde{\pi}'$ est une composante de $Ind_{\tilde{Q}'}^{\tilde{G}'}(\tilde{\sigma}')$ et qu'elle se factorise par $\tilde{p}$. Cette condition \'equivaut \`a ce que le caract\`ere central de la repr\'esentation sous-jacente $\pi'$ se restreigne en le caract\`ere trivial de $C({\mathbb R})$. Mais cette restriction est la m\^eme que celle du caract\`ere central de $\sigma'$. Donc $\tilde{\sigma}'$ se factorise en une $\omega$-repr\'esentation $\tilde{\sigma}$ de $\tilde{L}({\mathbb R})$. Cela oblige $\omega$ \`a \^etre trivial sur $Z_{L}({\mathbb R})^{\theta}$. On a encore l'\'egalit\'e 
  $$i\tilde{{\cal A}}_{\tilde{L}'}^*=i\tilde{{\cal A}}_{\tilde{L}}^*\oplus i{\cal A}_{C}^*.$$
  Les constructions entra\^{\i}nent que, pour tout $\tilde{\lambda}'\in i\tilde{{\cal A}}_{\tilde{L}'}^*$, on a l'\'egalit\'e
  $$(s(\varphi))(Ind_{\tilde{Q}'}^{\tilde{G}'}(\tilde{\sigma}'_{\tilde{\lambda}'}))=\varphi_{C}(\tilde{\xi})\varphi(Ind_{\tilde{Q}}^{\tilde{G}}(\tilde{\sigma}_{\tilde{\lambda}})),$$
  o\`u $\tilde{\lambda}$ et $\tilde{\xi}$ sont les composantes de $\tilde{\lambda}'$ selon la d\'ecomposition ci-dessus. La fonction de $\tilde{\lambda}'$ ci-dessus est de Paley-Wiener. Cela ach\`eve la preuve. $\square$

  \bigskip  
 \subsection{L'application $\phi_{\tilde{M}}$}

 On note ${\cal H}_{ac}(\tilde{G}(F))$ l'espace des fonctions $f:\tilde{G}(F)\to {\mathbb C}$ qui v\'erifient les conditions suivantes:
 
 (1) si $F$ est non-archim\'edien, il existe un sous-groupe ouvert compact $H$ de $G(F)$ tel que $f$ soit biinvariante par $H$;
 
 (2) si $F$ est archim\'edien, il existe un ensemble fini $\underline{\kappa}$ de repr\'esentations irr\'eductibles  de $K$ tel que la repr\'esentation de $K\times K$ dans l'espace engendr\'e par les translat\'es de $f$ \`a droite et \`a gauche par des \'el\'ements de $K$ n'ait pour composantes irr\'eductibles que des \'el\'ements de $\underline{\kappa}\times \underline{\kappa}$;
 
 (3) pour toute fonction $b\in C_{c}^{\infty}( \tilde{{\cal A}}_{\tilde{G},F} )$, la fonction produit $f(b\circ \tilde{H}_{\tilde{G}})$  appartient \`a $C_{c}^{\infty}(\tilde{G}(F))$.
 
 Remarquons que toute forme lin\'eaire sur $C_{c}^{\infty}(\tilde{G}(F),K)$ dont le support a une projection compacte dans $ \tilde{{\cal A}}_{\tilde{G},F}$ s'\'etend \`a ${\cal H}_{ac}(\tilde{G}(F))$: pour une telle forme lin\'eaire $l$ et pour $f\in {\cal H}_{ac}(\tilde{G}(F))$, on pose $l(f)=l(f(b\circ \tilde{H}_{\tilde{G}}))$, o\`u $b$ est un \'el\'ement de $C_{c}^{\infty}(\tilde{{\cal A}}_{\tilde{G},F})$ qui vaut $1$ sur un voisinage de cette projection. En particulier, les int\'egrales orbitales ou les int\'egrales orbitales pond\'er\'ees sont d\'efinies sur ${\cal H}_{ac}(\tilde{G}(F))$. On note $I_{ac}(\tilde{G}(F),\omega)$ le quotient de ${\cal H}_{ac}(\tilde{G}(F))$ par le sous-espace des $f\in {\cal H}_{ac}(\tilde{G}(F))$ telles que $I_{\tilde{G}}(\gamma,\omega,f)=0$ pour tout $\gamma\in \tilde{G}_{reg}(F)$.
 
  Les caract\`eres de $\omega$-repr\'esentations admissibles ne s'\'etendent pas \`a l'espace ${\cal H}_{ac}(\tilde{G}(F))$ mais nous allons voir que leurs coefficients de Fourier s'y \'etendent. Pr\'ecis\'ement, soit $\tilde{\pi}$ une $\omega$-repr\'esentation admissible de $\tilde{G}(F)$, soit $f\in C_{c}^{\infty}(\tilde{G}(F),K)$ et soit $X\in \tilde{{\cal A}}_{\tilde{G},F}$. Pour $\tilde{\lambda}\in i\tilde{{\cal A}}_{\tilde{G},F}^*$, on d\'efinit $\tilde{\pi}_{\tilde{\lambda}}$ comme en 2.6. La fonction $\tilde{\lambda}\mapsto I_{\tilde{G}}(\tilde{\pi}_{\tilde{\lambda}},f)e^{-<\tilde{\lambda},X>}$ se descend en une fonction sur $i{\cal A}_{\tilde{G},F}^*$. Cette fonction est $C^{\infty}$, \`a d\'ecroissance rapide si $F$ est archim\'edien. Posons
 $$I_{\tilde{G}}(\tilde{\pi},X,f)=mes(i{\cal A}_{\tilde{G},F}^*)^{-1}\int_{i{\cal A}_{\tilde{G},F}^*}I_{\tilde{G}}(\tilde{\pi}_{\tilde{\lambda}},f)e^{-<\tilde{\lambda},X>}\,d\lambda.$$
 On a
 
 (4) la fonction $X\mapsto I_{\tilde{G}}(\tilde{\pi},X,f)$ est \`a support compact et est $C^{\infty}$ dans le cas o\`u $F$ est archim\'edien; son support est contenu dans la projection dans $\tilde{A}_{\tilde{G},F}$ de celui de $f$; plus pr\'ecis\'ement, pour une fonction lisse $b$ sur $\tilde{{\cal A}}_{\tilde{G},F}$, on a l'\'egalit\'e $ I_{\tilde{G}}(\tilde{\pi},X,f(b\circ\tilde{H}_{\tilde{G}}))=b(X) I_{\tilde{G}}(\tilde{\pi},X,f)$.
 
 Preuve. Par $K$-finitude, $I_{\tilde{G}}(\tilde{\pi}_{\tilde{\lambda}},f)$ est de la forme
 $$\sum_{i=1,...,n}<\check{v}_{i},\tilde{\pi}_{\tilde{\lambda}}(f)v_{i}>,$$
  pour des \'el\'ements $v_{i}\in V_{\pi}$ et $\check{v}_{i}\in V_{\pi^{\vee}}$.  Autrement dit
 $$I_{\tilde{G}}(\tilde{\pi}_{\tilde{\lambda}},f)=\int_{\tilde{G}(F)} B(\tilde{\lambda},\gamma)f(\gamma)\,d\gamma,$$
 o\`u
 $$B(\tilde{\lambda},\gamma)=\sum_{i=1,...,n}<\check{v}_{i},\tilde{\pi}_{\tilde{\lambda}}(\gamma)v_{i}>.$$
 L'espace topologique $\tilde{G}(F)$ est un fibr\'e au-dessus de $\tilde{{\cal A}}_{\tilde{G},F}$, de fibres isomorphes \`a $G(F)^1$ (le noyau de $H_{\tilde{G}}$). Pour tout $Y\in \tilde{{\cal A}}_{\tilde{G},F}$, notons $\tilde{G}(F;Y)$ la fibre au-dessus de $Y$.  Les r\'esultats habituels de d\'ecomposition d'int\'egrales sur des fibr\'es nous disent que l'on peut d\'efinir une int\'egrale 
 $$I'_{\tilde{G}}(\tilde{\pi}_{\tilde{\lambda}},Y,f)=\int_{\tilde{G}(F;Y)} B(\tilde{\lambda},\gamma)f(\gamma)\,d\gamma$$
 qui est  \`a support compact en $Y$ et $C^{\infty}$ si $F$ est archim\'edien, de sorte que 
 $$ I_{\tilde{G}}(\tilde{\pi}_{\tilde{\lambda}},f)=\int_{\tilde{{\cal A}}_{\tilde{G},F}}I'_{\tilde{G}}(\tilde{\pi}_{\tilde{\lambda}},Y,f)\,dY.$$
 Mais 
 $$I'_{\tilde{G}}(\tilde{\pi}_{\tilde{\lambda}},Y,f)=e^{<\tilde{\lambda},Y>}I'_{\tilde{G}}(\tilde{\pi},Y,f).$$ 
 Par inversion de Fourier, on en d\'eduit l'\'egalit\'e $I_{\tilde{G}}(\tilde{\pi},X,f)= I'_{\tilde{G}}(\tilde{\pi},X,f)$ et l'assertion (4) s'ensuit. $\square$
 
 Il r\'esulte de (4) que l'on peut d\'efinir $I_{\tilde{G}}(\tilde{\pi},X,f)$ pour $f\in {\cal H}_{ac}(\tilde{G}(F))$. Le th\'eor\`eme 5.5 s'\'etend: 
 
 (5) un \'el\'ement $f\in {\cal H}_{ac}(\tilde{G}(F))$ a une image nulle dans $I_{ac}(\tilde{G}(F),\omega)$ si et seulement si $I_{\tilde{G}}(\tilde{\pi},X,f)=0$ pour toute $\omega$-repr\'esentation irr\'eductible et temp\'er\'ee $\tilde{\pi}$ de $\tilde{G}(F)$ et tout $X\in \tilde{{\cal A}}_{\tilde{G},F}$.
 
 Plus g\'en\'eralement, soit $\tilde{M}\in {\cal L}(\tilde{M}_{0})$, soit $\tilde{\pi}$ une $\omega$-repr\'esentation temp\'er\'ee de longueur finie de $\tilde{M}(F)$ et soit $X\in \tilde{{\cal A}}_{\tilde{M},F}$. Pour $f\in C_{c}^{\infty}(\tilde{G}(F),K)$ et $\tilde{\lambda}\in i\tilde{{\cal A}}_{\tilde{M},F}^*$, on d\'efinit $J_{\tilde{M}}^{\tilde{G}}(\tilde{\pi}_{\tilde{\lambda}},f)$. La fonction $\tilde{\lambda}\mapsto J_{\tilde{M}}^{\tilde{G}}(\tilde{\pi}_{\tilde{\lambda}},f)e^{-<\tilde{\lambda},X>}$ se descend en une fonction sur $i{\cal A}_{\tilde{M},F}^*$. Cette fonction est $C^{\infty}$, \`a d\'ecroissance rapide si $F$ est archim\'edien. Posons
 $$J_{\tilde{M}}^{\tilde{G}}(\tilde{\pi},X,f)=mes(i{\cal A}_{\tilde{M},F}^*)^{-1}\int_{i{\cal A}_{\tilde{M},F}^*} J_{\tilde{M}}^{\tilde{G}}(\tilde{\pi}_{\tilde{\lambda}},f)e^{-<\tilde{\lambda},X>}\,d\lambda.$$
 Comme fonction de $X$, ce terme n'est pas \`a support compact dans $\tilde{{\cal A}}_{\tilde{M},F}$. Toutefois
 
 (6) la fonction $X\mapsto J_{\tilde{M}}^{\tilde{G}}(\tilde{\pi},X,f)$ est  de Schwartz; la projection dans $\tilde{{\cal A}}_{\tilde{G},F}$ de son support est contenue dans celle du support de $f$; plus pr\'ecis\'ement, pour une fonction lisse $b$ sur $\tilde{{\cal A}}_{\tilde{G},F}$, on a l'\'egalit\'e $ J^{\tilde{G}}_{\tilde{M}}(\tilde{\pi},X,f(b\circ\tilde{H}_{\tilde{G}}))=b(X_{\tilde{G}}) J^{\tilde{G}}_{\tilde{M}}(\tilde{\pi},X,f)$, o\`u $X_{\tilde{G}}$ est l'image naturelle de $X$ dans $\tilde{{\cal A}}_{\tilde{G},F}$.

 Preuve. La fonction est de Schwartz car c'est la transform\'ee de Fourier d'une fonction de Schwartz. La preuve de la deuxi\`eme assertion est la m\^eme que celle de (4), la fonction $B(\tilde{\lambda},\gamma)$ ayant maintenant la forme
 $$B(\tilde{\lambda},\gamma)=\sum_{i=1,...,n}<\check{v}_{i},{\cal M}_{\tilde{M}}^{\tilde{G}}(\pi_{\lambda})\tilde{\Pi}_{\tilde{\lambda}}(\gamma)v_{i}>,$$
 o\`u $\tilde{\Pi}_{\tilde{\lambda}}=Ind_{\tilde{P}}^{\tilde{G}}(\tilde{\pi}_{\tilde{\lambda}})$, cf. 2.7 pour les notations. $\square$
 
 Il r\'esulte de (6) que l'on peut d\'efinir $J_{\tilde{M}}^{\tilde{G}}(\tilde{\pi},X,f)$ pour $f\in {\cal H}_{ac}(\tilde{G}(F))$.
 
 \ass{Proposition}{Soit $\tilde{M}\in {\cal L}(\tilde{M}_{0},\omega)$. Pour tout $f\in {\cal H}_{ac}(\tilde{G}(F))$, il existe $\phi_{\tilde{M}}(f)\in {\cal H}_{ac}(\tilde{M}(F))$ telle que, pour toute $\omega$-repr\'esentation temp\'er\'ee et de longueur finie $\tilde{\pi}$ de $\tilde{M}(F)$ et pour tout $X\in \tilde{{\cal A}}_{\tilde{M},F}$, on ait l'\'egalit\'e
 $$I_{\tilde{M}}(\tilde{\pi},X,\phi_{\tilde{M}}(f))=J_{\tilde{M}}^{\tilde{G}}(\tilde{\pi},X,f).$$
 L'image de $\phi_{\tilde{M}}(f)$ dans $I_{ac}(\tilde{M}(F),\omega)$ est uniquement d\'etermin\'ee.}
 
 Preuve. Soient $f\in C_{c}^{\infty}(\tilde{G}(F),K)$ et $b\in C_{c}^{\infty}(\tilde{{\cal A}}_{\tilde{M},F})$. Pour une $\omega$-repr\'esentation temp\'er\'ee et de longueur finie $\tilde{\pi}$ de $\tilde{M}(F)$, posons
 $$\varphi_{f,b}(\tilde{\pi})=\int_{\tilde{{\cal A}}_{\tilde{M},F}}J_{\tilde{M}}^{\tilde{G}}(\tilde{\pi},X,f)\overline{b(X)}\,dX.$$
 Ceci est convergent d'apr\`es (6). Par transformation de Fourier, on a l'\'egalit\'e
 $$\varphi_{f,b}(\tilde{\pi})=mes(i{\cal A}_{\tilde{M},F}^*)^{-1}\int_{i{\cal A}_{\tilde{M},F}^*}J_{\tilde{M}}^{\tilde{G}}(\tilde{\pi}_{\tilde{\lambda}},f)\overline{\hat{b}(\tilde{\lambda})}\,d\lambda,$$
 o\`u
 $$\hat{b}(\tilde{\lambda})=\int_{\tilde{{\cal A}}_{\tilde{M},F}}b(X)e^{<\tilde{\lambda},X>}\,dX.$$
 On va montrer
 
 (7) la fonction $\varphi_{f,b}$ appartient \`a l'espace ${\cal PW}(\tilde{M},\omega)$.
 
 D'apr\`es 2.7(1), la fonction s'identifie bien \`a une fonction sur le groupe de Grothendieck ${\cal K}^{\tilde{M}}$ de 6.1. D'apr\`es 2.7(2), elle annule les repr\'esentations irr\'eductibles et non $M$-irr\'eductibles. Les conditions (2) et (3) de 6.1 sont \'evidentes puisque $f$ est $K$-finie. Il faut v\'erifier 6.1(1). D'apr\`es la bijectivit\'e de l'application $res$ de 6.2, il suffit de v\'erifier l'assertion suivante:
 
 (8) soit $\tilde{Q}=\tilde{L}U_{Q}\in {\cal F}^{\tilde{M}}(\tilde{M}_{0})$   et soit $\boldsymbol{\tau}\in {\cal E}_{ell}(\tilde{L},\omega)$; alors la fonction $\tilde{\mu}\mapsto \varphi_{f,b}(Ind_{\tilde{Q}}^{\tilde{M}}(\tilde{\pi}_{\boldsymbol{\tau}_{\tilde{\mu}}}))$ sur $i\tilde{{\cal A}}_{\tilde{L},F}^*$ est de Paley-Wiener.
 
 Posons simplement $\tilde{\pi}=\tilde{\pi}_{\boldsymbol{\tau}}$. Remarquons que, pour $\tilde{\lambda}\in i\tilde{{\cal A}}_{\tilde{M},F}^*$, on a l'\'egalit\'e $(Ind_{\tilde{Q}}^{\tilde{M}}(\tilde{\pi}_{\boldsymbol{\tau}_{\tilde{\mu}}}))_{\tilde{\lambda}}=Ind_{\tilde{Q}}^{\tilde{M}}(\tilde{\pi}_{\tilde{\mu}+\tilde{\lambda}})$.
 On utilise la formule de descente du lemme 5.4(iv) (o\`u les r\^oles de $\tilde{M}$ et $\tilde{L}$ sont \'echang\'es). C'est-\`a-dire
 $$J_{\tilde{M}}^{\tilde{G}}((Ind_{\tilde{Q}}^{\tilde{M}}(\tilde{\pi}_{\tilde{\mu}+\tilde{\lambda}}),f)=\sum_{\tilde{M}'\in {\cal L}(\tilde{L})}d_{\tilde{L}}^{\tilde{G}}(\tilde{M},\tilde{M}')J_{\tilde{L}}^{\tilde{M}'}(\tilde{\pi}_{\tilde{\mu}+\tilde{\lambda}},f_{\tilde{P}'}).$$
 On peut fixer $\tilde{M}'\in {\cal L}(\tilde{L})$ tel que $d_{\tilde{L}}^{\tilde{G}}(\tilde{M},\tilde{M}')\not=0$ et prouver que la fonction $\psi$ sur $i\tilde{{\cal A}}_{\tilde{L},F}^*$ d\'efinie par
 $$\psi(\tilde{\mu})=mes(i{\cal A}_{\tilde{M},F}^*)^{-1}\int_{i{\cal A}_{\tilde{M},F}^*}J_{\tilde{L}}^{\tilde{M}'}(\tilde{\pi}_{\tilde{\mu}+\tilde{\lambda}},f_{\tilde{P}'})\overline{\hat{b}(\tilde{\lambda})}\,d\lambda$$
 est de Paley-Wiener. Gr\^ace \`a (6), on peut exprimer $J_{\tilde{L}}^{\tilde{M}'}(\tilde{\pi}_{\tilde{\mu}+\tilde{\lambda}},f_{\tilde{P}'})$ par inversion de Fourier:
 $$J_{\tilde{L}}^{\tilde{M}'}(\tilde{\pi}_{\tilde{\mu}+\tilde{\lambda}},f_{\tilde{P}'})=\int_{\tilde{{\cal A}}_{\tilde{L},F}}J_{\tilde{L}}^{\tilde{M}'}(\tilde{\pi},Y,f_{\tilde{P}'})e^{<\tilde{\mu}+\tilde{\lambda},Y>}\,dY.$$
 D'o\`u
 $$\psi(\tilde{\mu})=mes(i{\cal A}_{\tilde{M},F}^*)^{-1}\int_{i{\cal A}_{\tilde{M},F}^*}\int_{\tilde{{\cal A}}_{\tilde{L},F}}J_{\tilde{L}}^{\tilde{M}'}(\tilde{\pi},Y,f_{\tilde{P}'})e^{<\tilde{\mu}+\tilde{\lambda},Y>}\overline{\hat{b}(\tilde{\lambda})}\,dY\,d\lambda.$$
Cette expression est absolument convergente. En int\'egrant d'abord en $\lambda$, on obtient par inversion de Fourier
$$\psi(\tilde{\mu})= \int_{\tilde{{\cal A}}_{\tilde{L},F}} \Psi(Y)e^{<\tilde{\mu},Y>}\,dY,$$
o\`u, en notant  $Y_{\tilde{M}}$  la projection naturelle de $Y$ dans $\tilde{{\cal A}}_{\tilde{M},F}$, 
$$\Psi(Y)=J_{\tilde{L}}^{\tilde{M}'}(\tilde{\pi},Y,f^{\tilde{P}'}) \bar{b}(Y_{\tilde{M}}).$$
 Ainsi, $\psi$ appara\^{\i}t comme la transform\'ee de Fourier de la fonction $\Psi$ sur $\tilde{{\cal A}}_{\tilde{L},F}$. Il s'agit de prouver que $\Psi$ est lisse et \`a support compact. Elle est lisse car c'est le produit de deux fonctions lisses. L'hypoth\`ese sur $\tilde{M}'$ implique que
 $${\cal A}_{\tilde{L}}^{\tilde{M}}\cap {\cal A}_{\tilde{L}}^{\tilde{M}'}=0,$$
 donc que la somme directe des projections
 $${\cal A}_{\tilde{L},F}\to {\cal A}_{\tilde{M},F}\oplus {\cal A}_{\tilde{M}',F}$$
 est injective. La r\'eplique pour les espaces affines est que le produit des projections
 $$\tilde{{\cal A}}_{\tilde{L},F}\to \tilde{{\cal A}}_{\tilde{M},F}\times \tilde{{\cal A}}_{\tilde{M}',F}$$
 est injectif. Il est clair que c'est une immersion ferm\'ee. La projection dans $\tilde{{\cal A}}_{\tilde{M}',F}$ du support de la fonction $Y\mapsto J_{\tilde{L}}^{\tilde{M}'}(\tilde{\pi},Y,f_{\tilde{P}'})$ est compacte d'apr\`es (6). La projection dans $\tilde{{\cal A}}_{\tilde{M},F}$ du support de la fonction $Y\mapsto \bar{b}(Y_{\tilde{M}})$ est compacte puisque $b$ est \`a support compact. Donc $\Psi$ est \`a support compact. Cela d\'emontre (8) et ach\`eve la preuve de (7).
 
 D'apr\`es le th\'eor\`eme 6.1, on peut choisir une fonction $\phi_{\tilde{M}}(f,b)\in C_{c}^{\infty}(\tilde{M}(F),K^M)$ (o\`u $K^M=K\cap M(F)$) de sorte que $I_{\tilde{M}}(\tilde{\pi},\phi_{\tilde{M}}(f,b))=\varphi_{f,b}(\tilde{\pi})$ pour toute $\omega$-repr\'esentation temp\'er\'ee et de longueur finie $\tilde{\pi}$ de $\tilde{M}(F)$. On peut pr\'eciser le comportement de $\phi_{\tilde{M}}(f,b)$ par translations \`a droite ou \`a gauche par $K$. Traitons le cas non-archim\'edien (le cas archim\'edien n'en diff\`ere que par les notations). Fixons un sous-groupe ouvert  compact $H$ de $G(F)$ tel que $f$ soit biinvariante par $H$. Il existe un tel sous-groupe $H'$ de $M(F)$, ne d\'ependant que de $H$, tel que $J_{\tilde{M}}^{\tilde{G}}(\tilde{\pi},f)=0$ si $\pi$ n'a pas d'invariant non nul par $H'$. A fortiori, dans ce cas,  $\varphi_{f,b}(\tilde{\pi})=0$ pour tout $b$. D'apr\`es le compl\'ement au th\'eor\`eme 6.1, on peut fixer un sous-groupe ouvert compact $H^M$ de $M(F)$, ne d\'ependant que de $H'$, donc ne d\'ependant que de $H$, et supposer que $\phi_{\tilde{M}}(f,b)$  est biinvariante par $H^M$.
 
 Remarquons que la fonction $\tilde{\lambda}\mapsto \varphi(\tilde{\pi}_{\tilde{\lambda}})$ sur $i\tilde{{\cal A}}_{\tilde{M},F}^*$ est par construction la transform\'ee de Fourier de la fonction $X\mapsto J_{\tilde{M}}(\tilde{\pi},X,f)\bar{b}(X)$ sur $\tilde{{\cal A}}_{\tilde{M},F}$. Par transformation de Fourier, l'\'egalit\'e $I_{\tilde{M}}(\tilde{\pi},\phi_{\tilde{M}}(f,b))=\varphi_{f,b}(\tilde{\pi})$ est donc \'equivalente \`a
$$(9) \qquad I_{\tilde{M}}(\tilde{\pi},X,\phi_{\tilde{M}}(f,b))=J_{\tilde{M}}^{\tilde{G}}(\tilde{\pi},X,f)\bar{b}(X)$$
pour tout $\tilde{\pi}$ et tout $X\in \tilde{{\cal A}}_{\tilde{M},F}$.
 
 Fixons une suite $(U_{n})_{n\in {\mathbb N}}$ de sous-ensembles ouverts relativement compacts de $\tilde{{\cal A}}_{\tilde{M},F}$ de sorte que $\bar{U}_{n}\subset U_{n+1}$ (o\`u $\bar{U}_{n}$ est la cl\^oture de $U_{n}$) et que $\tilde{{\cal A}}_{\tilde{M},F}=\cup_{n\in {\mathbb N}}U_{n}$. En posant $U'_{n}=U_{n}-\bar{U}_{n-2}$ pour $n\geq2$ et $U'_{1}=U_{1}$, on a aussi $\tilde{{\cal A}}_{\tilde{M},F}=\cup_{n\geq1}U'_{n}$. On peut choisir une partition de l'unit\'e $(b_{n})_{n\geq1}$ relative \`a ce dernier recouvrement, form\'ee de fonctions lisses (et forc\'ement \`a supports compacts). On choisit  une suite $(c_{n})_{n\geq1}$ de fonctions lisses \`a supports compacts sur $\tilde{{\cal A}}_{\tilde{M},F}$ de sorte que $c_{n}$ vaille $1$ sur $U'_{n}$ et, pour $n\geq3$, le support de $c_{n}$ soit contenu dans $U_{n+1}-\bar{U}_{n-3}$. On choisit enfin une suite $(d_{n})_{n\geq1}$ de fonctions lisses \`a supports compacts sur $\tilde{{\cal A}}_{\tilde{G},F}$ de sorte que $d_{n}$ vaille $1$ sur un voisinage de la projection de $U'_{n}$ dans $\tilde{{\cal A}}_{\tilde{G},F}$. Soit $f\in {\cal H}_{ac}(\tilde{G}(F))$. Pour tout $n\geq1$, la fonction $f(d_{n}\circ\tilde{H}_{\tilde{G}})$ est \`a support compact, on dispose donc de $\phi_{\tilde{M}}(f(d_{n}\circ\tilde{H}_{\tilde{G}}),b_{n})$. Posons
 $$\phi_{\tilde{M}}(f)=\sum_{n\geq1}(c_{n}\circ\tilde{H}_{\tilde{M}})\phi_{\tilde{M}}(f(d_{n}\circ\tilde{H}_{\tilde{G}}),b_{n}).$$
 Cette s\'erie est convergente: elle est localement finie d'apr\`es les propri\'et\'es de la suite $(c_{n})_{n\geq1}$. Pour la m\^eme raison, si $b\in C_{c}^{\infty}(\tilde{{\cal A}}_{\tilde{M},F})$, la fonction $\phi_{\tilde{M}}(f)(b\circ\tilde{H}_{\tilde{M}})$ est \`a support compact. Enfin, $\phi_{\tilde{M}}(f)$ v\'erifie les propri\'et\'es (1) ou (2). Par exemple, dans le cas non-archim\'edien, fixons un sous-groupe ouvert compact $H$ de $G(F)$ tel que $f$ soit biinvariante par $H$. Chaque fonction
$f(d_{n}\circ\tilde{H}_{\tilde{G}})$ v\'erifie la m\^eme propri\'et\'e. Comme on l'a dit ci-dessus, on peut supposer  $\phi_{\tilde{M}}(f(d_{n}\circ\tilde{H}_{\tilde{G}}),b_{n})$ biinvariante par $H^M$, o\`u $H^M$ est ind\'ependant de $n$. Alors $\phi_{\tilde{M}}(f)$ est aussi biinvariante par $H^M$. Cela prouve que $\phi_{\tilde{M}}(f)$ appartient \`a ${\cal H}_{ac}(\tilde{M})$. Soit $\tilde{\pi}$ une $\omega$-repr\'esentation temp\'er\'ee et de longueur finie de $\tilde{M}(F)$ et soit $X\in \tilde{{\cal A}}_{\tilde{M},F}$. En utilisant (4), on a
$$I_{\tilde{M}}(\tilde{\pi},X,\phi_{\tilde{M}}(f))=\sum_{n\geq1}c_{n}(X)I_{\tilde{M}}(\tilde{\pi},X,\phi_{\tilde{M}}(f(d_{n}\circ\tilde{H}_{\tilde{G}}),b_{n})).$$
D'o\`u, d'apr\`es (9),
$$I_{\tilde{M}}(\tilde{\pi},X,\phi_{\tilde{M}}(f))=\sum_{n\geq1}c_{n}(X)b_{n}(X)J_{\tilde{M}}^{\tilde{G}}(\tilde{\pi},X,f(d_{n}\circ\tilde{H}_{\tilde{G}})).$$
D'apr\`es les choix de nos fonctions, on a les \'egalit\'es
$$c_{n}(X)b_{n}(X)=b_{n}(X),$$
$$b_{n}(X)J_{\tilde{M}}^{\tilde{G}}(\tilde{\pi},X,f(d_{n}\circ\tilde{H}_{\tilde{G}}))=b_{n}(X)J_{\tilde{M}}^{\tilde{G}}(\tilde{\pi},X,f).$$
Donc
 $$I_{\tilde{M}}(\tilde{\pi},X,\phi_{\tilde{M}}(f))=\sum_{n\geq1}b_{n}(X)J_{\tilde{M}}^{\tilde{G}}(\tilde{\pi},X,f)=J_{\tilde{M}}^{\tilde{G}}(\tilde{\pi},X,f),$$
 puisque $(b_{n})_{n\geq1}$ est une partition de l'unit\'e. C'est la propri\'et\'e requise, ce qui d\'emontre la premi\`ere partie de la proposition. Que l'image de $\phi_{\tilde{M}}(f)$ dans $I_{ac}(\tilde{M}(F),\omega)$ soit uniquement d\'etermin\'ee r\'esulte de (5). $\square$
 
 Remarquons que, pour $f\in {\cal H}_{ac}(\tilde{G}(F))$ et pour une fonction lisse $b$ sur $\tilde{{\cal A}}_{\tilde{G},F}$, 
 
 (10) les images dans $I_{ac}(\tilde{M}(F),\omega)$ de $\phi_{\tilde{M}}(f(b\circ\tilde{H}_{\tilde{G}}))$ et $(b\circ\tilde{H}_{\tilde{G}})\phi_{\tilde{M}}(f)$ co\"{\i}ncident.
 
 Cela r\'esulte de (4) et (6).

 Pour $f\in {\cal H}_{ac}(\tilde{G}(F))$, on note encore $\phi_{\tilde{M}}(f)$ l'image dans $I_{ac}(\tilde{M}(F),\omega)$ de la fonction ainsi not\'ee dans l'\'enonc\'e. Alors $\phi_{\tilde{M}}$ devient une application lin\'eaire bien d\'efinie de ${\cal H}_{ac}(\tilde{G}(F))$ dans $I_{ac}(\tilde{M}(F),\omega)$.
 
 Soit $\tilde{P}=\tilde{M}U_{P}\in {\cal F}(\tilde{M}_{0})$. On dispose de l'application $f\mapsto f_{\tilde{P}}$ de $C_{c}^{\infty}(\tilde{G}(F),K)$ dans $C_{c}^{\infty}(\tilde{M}(F),K^M)$. On sait que l'image de $f_{\tilde{P}}$ dans $I(\tilde{M}(F),K,\omega)$ ne d\'epend que de $\tilde{M}$ et pas de $\tilde{P}$. On note cette image $f_{\tilde{M}}$. Ceci s'\'etend \`a l'espace ${\cal H}_{ac}(\tilde{G}(F))$: pour $f\in {\cal H}_{ac}(\tilde{G}(F))$, on d\'efinit $f_{\tilde{M}}\in I_{ac}(\tilde{M}(F),\omega)$. Soient maintenant $\tilde{M}\in {\cal L}(\tilde{M}_{0},\omega)$ et $\tilde{L}\in {\cal L}(\tilde{M})$. Pour $f\in {\cal H}_{ac}(\tilde{G}(F))$, on a l'\'egalit\'e dans $I_{ac}(\tilde{M}(F),\omega)$:
 
 $$(11)\qquad (\phi_{\tilde{L}}^{\tilde{G}}(f))_{\tilde{M}}=\sum_{\tilde{L}'\in {\cal L}(\tilde{M})}d_{\tilde{M}}^{\tilde{G}}(\tilde{L},\tilde{L}')\phi_{\tilde{M}}^{\tilde{L}'}(f_{\tilde{Q}'}).$$
 
 La preuve est formelle \`a partir de la formule de descente du lemme 5.4(iv).

 \bigskip

 \subsection{Int\'egrales orbitales pond\'er\'ees \'equivariantes}
 Appelons forme lin\'eaire $\omega$-\'equivariante sur $C_{c}^{\infty}(\tilde{G}(F),K)$, resp. ${\cal H}_{ac}(\tilde{G}(F))$, une forme lin\'eaire qui se factorise en une  forme lin\'eaire  sur $I(\tilde{G}(F),K,\omega)$, resp. $I_{ac}(\tilde{G}(F),\omega)$. Dans le cas o\`u $F$ est non-archim\'edien, une forme lin\'eaire $l$ disons sur $C_{c}^{\infty}(\tilde{G}(F))$ est  $\omega$-\'equivariante si et seulement si elle v\'erifie la relation $l(^gf)=\omega(g)^{-1}l(f)$ pour toute $f\in C_{c}^{\infty}(\tilde{G}(F))$ et tout $g\in G(F)$ (cf. 5.5 remarque (2)).  On identifiera souvent une forme lin\'eaire $\omega$-\'equivariante \`a une forme lin\'eaire sur $I(\tilde{G}(F),K,\omega)$, resp. $I_{ac}(\tilde{G}(F),\omega)$.
 
 \ass{Proposition}{Soient $\tilde{M}\in {\cal L}(\tilde{M}_{0})$ et $\gamma\in \tilde{M}(F)\cap \tilde{G}_{reg}(F)$. Il existe une unique forme lin\'eaire $\omega$-\'equivariante $f\mapsto I_{\tilde{M}}^{\tilde{G}}(\gamma,\omega,f)$ sur ${\cal H}_{ac}(\tilde{G}(F))$ qui v\'erifie l'\'egalit\'e
 $$I_{\tilde{M}}^{\tilde{G}}(\gamma,\omega,f)=J_{\tilde{M}}^{\tilde{G}}(\gamma,\omega,f)-\sum_{\tilde{L}\in {\cal L}(\tilde{M}), \tilde{L}\not=\tilde{G}}I_{\tilde{M}}^{\tilde{L}}(\gamma,\omega,\phi_{\tilde{L}}(f))$$
 pour tout $f\in {\cal H}_{ac}(\tilde{G}(F))$.}
 
   Pour donner un sens \`a cet \'enonc\'e, on doit raisonner par r\'ecurrence sur  le rang semi-simple  $rg_{ss}(G)$ de $G$. Si ce rang est nul, ou plus g\'en\'eralement si  $\tilde{M}=\tilde{G}$, l'\'enonc\'e est  tautologique: on a simplement  
   $$I_{\tilde{G}}^{\tilde{G}}(\gamma,\omega,f)=J_{\tilde{G}}^{\tilde{G}}(\gamma,\omega,f).$$
   Si $rg_{ss}(G)>0$, on suppose par r\'ecurrence que les formes lin\'eaires $I_{\tilde{M}}^{\tilde{L}}(\gamma,\omega,.)$ sont d\'efinies pour $\tilde{L}\not=\tilde{G}$ (auquel cas $ rg_{ss}(L)<rg_{ss}(G)$) et qu'elles sont  $\omega$-\'equivariantes. D'apr\`es la derni\`ere assertion de la proposition pr\'ec\'edente, le terme $I_{\tilde{M}}^{\tilde{L}}(\gamma,\omega,\phi_{\tilde{L}}(f))$ est bien d\'efini. Il en est donc de m\^eme du membre de droite de l'\'egalit\'e de l'\'enonc\'e. Cette \'egalit\'e d\'efinit la forme lin\'eaire $f\mapsto I_{\tilde{M}}^{\tilde{G}}(\gamma,\omega,f)$. L'assertion de la proposition est que celle-ci est $\omega$-\'equivariante. Dans le cas o\`u $F$ est non-archim\'edien, on montre que l'\'egalit\'e $I_{\tilde{M}}^{\tilde{G}}(\gamma,\omega,^gf)=\omega(g)^{-1}I_{\tilde{M}}^{\tilde{G}}(\gamma,\omega,f)$ est v\'erifi\'ee pour tout $g\in G(F)$: 
   la d\'emonstration, essentiellement formelle, est la m\^eme que dans le cas non tordu. On se contente de renvoyer \`a [A2] proposition 4.1 pour ce cas.  Cela suffit pour conclure. Si $F$ est archim\'edien, cette relation n'a plus de sens si on se limite aux fonctions $K$-finies. On peut l'adapter \`a de telles fonctions, mais elle ne suffit de toute fa\c{c}on pas \`a conclure.  Pour l'instant, nous laissons la preuve inachev\'ee (on la compl\`etera en 7.1). 
   
   {\bf On conserve pour  ce paragraphe et jusqu'\`a la fin de 7.1 l'hypoth\`ese de r\'ecurrence ci-dessus, \`a savoir que la proposition est v\'erifi\'ee si l'on remplace $\tilde{G}$ par $\tilde{G}'$ avec $rg_{ss}(G')<rg_{ss}(G)$}. 
   
   Comme on l'a dit, cela suffit \`a d\'efinir la forme lin\'eaire $f\mapsto I_{\tilde{M}}^{\tilde{G}}(\gamma,\omega,f)$ sur ${\cal H}_{ac}(\tilde{G}(F))$.
   
   Il est clair que $I_{\tilde{M}}^{\tilde{G}}(\gamma,\omega,f)=0$ si $\omega$ n'est pas trivial sur $Z_{G}(\gamma,F)$.
   
   On aura besoin des propri\'et\'es suivantes. On ne d\'emontrera pas les deux premi\`eres, leurs preuves \'etant essentiellement formelles. Soient $\tilde{M},\tilde{M}'\in {\cal L}(\tilde{M}_{0})$ et $g\in G(F)$. Supposons $\tilde{M}'=g\tilde{M}g^{-1}$. Alors
  
  (1) pour tout $\gamma\in \tilde{M}(F)\cap \tilde{G}_{reg}(F)$ et tout $f\in {\cal H}_{ac}(\tilde{G}(F))$, on a l'\'egalit\'e
  $$I_{\tilde{M}'}^{\tilde{G}}(g\gamma g^{-1},f)=\omega(g)I_{\tilde{M}}^{\tilde{G}}(\gamma,\omega,f).$$

  Soient $\tilde{M}\in {\cal L}(\tilde{M}_{0})$ et $\tilde{L}\in {\cal L}(\tilde{M})$. Soit $\gamma\in \tilde{M}(F)\cap \tilde{G}_{reg}(F)$. Pour $f\in {\cal H}_{ac}(\tilde{G}(F))$, on a l'\'egalit\'e
  $$(2) \qquad I_{\tilde{L}}^{\tilde{G}}(\gamma,\omega,f)=\sum_{\tilde{L}'\in {\cal L}(\tilde{M})}d_{\tilde{M}}^{\tilde{G}}(\tilde{L},\tilde{L}')I_{\tilde{M}}^{\tilde{L}'}(\gamma,\omega,f_{\tilde{L}'}).$$
  
  Soient $\tilde{M}\in {\cal L}(\tilde{M}_{0})$, $\gamma\in \tilde{M}(F)\cap \tilde{G}_{reg}(F)$ et $b$ une fonction lisse sur $\tilde{{\cal A}}_{\tilde{G},F}$. On a l'\'egalit\'e
  $$(3) \qquad I_{\tilde{M}}^{\tilde{G}}(\gamma,\omega,f(b\circ\tilde{H}_{\tilde{G}}))=b(\tilde{H}_{\tilde{G}}(\gamma))I_{\tilde{M}}^{\tilde{G}}(\gamma,\omega,f).$$
  
  Preuve. La propri\'et\'e analogue pour les int\'egrales pond\'er\'ees $J_{\tilde{M}}^{\tilde{G}}(\gamma,\omega,.)$ r\'esulte des d\'efinitions. D'autre part, pour $\tilde{L}\in {\cal L}(\tilde{M})$, on peut supposer $\phi_{\tilde{L}}( f(b\circ\tilde{H}_{\tilde{G}}))=(b\circ\tilde{H}_{\tilde{G}})\phi_{\tilde{L}}(f)$, cf. 6.4(10). Pour $\tilde{L}\not=\tilde{G}$, on peut supposer par r\'ecurrence que (3) est vrai quand on remplace $(\tilde{G},\tilde{M})$ par $(\tilde{L},\tilde{M})$. Alors, quand on remplace $f$ par $f(b\circ\tilde{H}_{\tilde{G}})$, le membre de droite de l'\'egalit\'e de l'\'enonc\'e est multipli\'e par $b(\tilde{H}_{\tilde{G}}(\gamma))$. Donc le membre de gauche aussi. $\square$
  
  Soient $\tilde{M}\in {\cal L}(\tilde{M}_{0})$ et $\tilde{T}$ un tore tordu maximal de $\tilde{M}$. Alors
  
  (4) il existe un entier $N\geq0$ et, pour $f\in {\cal H}_{ac}(\tilde{G}(F))$ et $\gamma_{0}\in \tilde{T}(F)$, il existe $c>0$ et un voisinage $\Omega$ de $\gamma_{0}$ dans $\tilde{T}(F)$ tel que l'on ait la majoration
  $$\vert I_{\tilde{M}}^{\tilde{G}}(\gamma,\omega,f)\vert\leq c(1+\vert log(D^{\tilde{G}}(\gamma))\vert )^N$$
  pour tout $\gamma\in \Omega\cap \tilde{G}_{reg}(F)$;
  
  (5) pour tout $f\in {\cal H}_{ac}(\tilde{G}(F))$,   la fonction $\gamma\mapsto I_{\tilde{M}}^{\tilde{G}}(\gamma,\omega,f)$ est   lisse sur $\tilde{T}(F)\cap \tilde{G}_{reg}(F)$.
  
  Preuve. En raisonnant par r\'ecurrence, il suffit de d\'emontrer les m\^emes propri\'et\'es pour $J_{\tilde{M}}^{\tilde{G}}(\gamma,\omega,f)$. Au voisinage d'un point $\gamma_{0}$, cette fonction co\"{\i}ncide avec $J_{\tilde{M}}^{\tilde{G}}(\gamma,\omega,f(b\circ\tilde{H}_{\tilde{G}}))$, o\`u $b\in C_{c}^{\infty}(\tilde{{\cal A}}_{\tilde{G},F})$ vaut $1$ sur un voisinage de $\tilde{H}_{\tilde{G}}(\gamma_{0})$. Cela nous ram\`ene au cas o\`u $f\in C_{c}^{\infty}(\tilde{G}(F),K)$, lequel cas est trait\'e par   5.4(2) et (3). $\square$

 \bigskip
 
 \subsection{Le th\'eor\`eme }
 
 Soient $f_{1},f_{2}\in C_{c}^{\infty}(\tilde{G}(F),K)$. Pour $\tilde{M}\in {\cal L}(\tilde{M}_{0})$, $\tilde{S}\in T_{ell}(\tilde{M},\omega)$   et pour $\gamma\in \tilde{S}(F)\cap \tilde{G}_{reg}(F)$, on pose
 $$I_{\tilde{M}}^{\tilde{G}}(\gamma,\omega,f_{1},f_{2})=\sum_{\tilde{L}_{1},\tilde{L}_{2}\in {\cal L}(\tilde{M})}d_{\tilde{M}}^{\tilde{G}}(\tilde{L}_{1},\tilde{L}_{2})\overline{I_{\tilde{M}}^{\tilde{L}_{1}}( \gamma,\omega,f_{1,\tilde{L}_{1}})}I_{\tilde{M}}^{\tilde{L}_{2}}(\gamma,\omega,f_{2,\tilde{L}_{2}}).$$
 
 {\bf Remarque.} Pour $i=1,2$, l'\'el\'ement $f_{i,\tilde{L}_{i}}$ appartient \`a $I(\tilde{L}_{i}(F),K^L,\omega)$. Si $\tilde{L}_{i}\subsetneq \tilde{G}$, le terme $I_{\tilde{M}}^{\tilde{L}_{i}}(\gamma,\omega,f_{i,\tilde{L}_{i}})$ est bien d\'efini d'apr\`es la proposition 6.5. Si $\tilde{L}_{i}=\tilde{G}$, on ne sait pas encore que $f\mapsto I_{\tilde{M}}^{\tilde{G}}(\gamma,\omega,f)$ se factorise par $I(\tilde{G}(F),K,\omega)$. Par convention, on suppose dans ce cas que $f_{\tilde{G}}=f$ et $I_{\tilde{M}}^{\tilde{G}}(\gamma,\omega,f_{\tilde{G}})=I_{\tilde{M}}^{\tilde{G}}(\gamma,\omega,f)$.
 
 On pose
 $$I_{\tilde{M},\tilde{S}}^{\tilde{G}}(\omega,f_{1},f_{2})=\vert W^M(\tilde{S})\vert ^{-1} mes(A_{\tilde{M}}(F)\backslash S^{\theta}(F))\int_{\tilde{S}(F)/(1-\theta)(S(F))}I_{\tilde{M}}^{\tilde{G}}(\gamma,\omega,f_{1},f_{2})\,d\gamma.$$
 Montrons que
 
 (1) cette int\'egrale est absolument convergente. 
 
 Preuve. D'apr\`es 6.5(4) et (5) et 4.2(2), la fonction \`a int\'egrer est localement int\'egrable. Il suffit de prouver qu'elle est \`a support compact. Puisque $\tilde{S}$ est elliptique dans $\tilde{M}$, il suffit de prouver que la projection de ce support sur $\tilde{{\cal A}}_{\tilde{M}}$ l'est. On peut encore fixer $\tilde{L}_{1},\tilde{L}_{2}$ tels que $d_{\tilde{M}}^{\tilde{G}}(\tilde{L}_{1},\tilde{L}_{2})\not=0$ et remplacer la fonction par 
 $$\overline{I_{\tilde{M}}^{\tilde{L}_{1}}( \gamma,\omega,f_{1,\tilde{L}_{1}})}I_{\tilde{M}}^{\tilde{L}_{2}}(\gamma,\omega,f_{2,\tilde{L}_{2}}).$$
 Puisque $f_{1,\tilde{L}_{1}}$ est \`a support compact, la relation 6.5(3) entra\^{\i}ne que  le support de la premi\`ere fonction ci-dessus a une projection compacte dans $\tilde{{\cal A}}_{\tilde{L}_{1}}$. De m\^eme, le support de la seconde fonction a une projection compacte dans $\tilde{{\cal A}}_{\tilde{L}_{2}}$. La non-nullit\'e de $d_{\tilde{M}}^{\tilde{G}}(\tilde{L}_{1},\tilde{L}_{2})$ entra\^{\i}ne que le produit des projections
  $$\tilde{{\cal A}}_{\tilde{M}}\to\tilde{{\cal A}}_{\tilde{L}_{1}}\times\tilde{{\cal A}}_{\tilde{L}_{2}}$$
  est injective. Donc le support du produit des deux fonctions a bien les propri\'et\'es requises. $\square$
  
  Posons
  $$I_{\tilde{M},g\acute{e}om}^{\tilde{G}}(\omega,f_{1},f_{2})=\sum_{\tilde{S}\in T_{ell}(\tilde{M},\omega)}I_{\tilde{M},\tilde{S}}^{\tilde{G}}(\omega,f_{1},f_{2}),$$
  puis
  $$I_{g\acute{e}om}^{\tilde{G}}(\omega,f_{1},f_{2})=\sum_{\tilde{M}\in {\cal L}(\tilde{M}_{0})}\vert \tilde{W}^M\vert \vert \tilde{W}^G\vert ^{-1}(-1)^{a_{\tilde{M}}-a_{\tilde{G}}}I_{\tilde{M},g\acute{e}om}^{\tilde{G}}(\omega,f_{1},f_{2}).$$
  
  D'autre part, posons
  $$I_{disc}^{\tilde{G}}(\omega,f_{1},f_{2})=J_{\tilde{G},spec}^{\tilde{G}}(\omega,f_{1},f_{2}),$$
  cf. 3.25.  
  
  \ass{Th\'eor\`eme}{Pour tous $f_{1},f_{2}\in C_{c}^{\infty}(\tilde{G}(F),K)$, on a l'\'egalit\'e
  $$I_{g\acute{e}om}^{\tilde{G}}(\omega,f_{1},f_{2})=I_{disc}^{\tilde{G}}(\omega,f_{1},f_{2}).$$}
  
  Preuve. Cette preuve est formelle. Compte tenu de l'importance du th\'eor\`eme, nous la traitons en d\'etail.  On a besoin de quelques constructions, formelles comme on vient de le dire. Pour deux espaces vectoriels complexes $V$ et $V'$, notons $V\boxtimes V'$ leur produit tensoriel "sesquilin\'eaire", pr\'ecis\'ement le produit tensoriel $\bar{V}\otimes V'$, o\`u $\bar{V}$ est le conjugu\'e complexe de $V$. Posons simplement
  $$\underline{C}(\tilde{G}(F))=C_{c}^{\infty}(\tilde{G}(F),K)\boxtimes C_{c}^{\infty}(\tilde{G}(F),K),$$
  $$\underline{H}_{ac}(\tilde{G}(F))={\cal H}_{ac}(\tilde{G}(F))\boxtimes {\cal H}_{ac}(\tilde{G}(F)),$$
  $$\underline{I}_{ac}(\tilde{G}(F),\omega)=I_{ac}(\tilde{G}(F),\omega)\boxtimes I_{ac}(\tilde{G}(F),\omega).$$

     Pour $\tilde{L}\in {\cal L}(\tilde{M}_{0})$,  $\tilde{L}\not=\tilde{G}$, on d\'efinit une application
  $$\underline{\phi}_{\tilde{L}}: \underline{H}_{ac}(\tilde{G}(F))\to \underline{I}_{ac}(\tilde{L}(F),\omega)$$
  par la formule
  $$(2) \qquad \underline{\phi}_{\tilde{L}}(f_{1}\boxtimes f_{2})=\sum_{\tilde{L}_{1},\tilde{L}_{1}\in {\cal L}(\tilde{L})}d_{\tilde{L}}^{\tilde{G}}(\tilde{L}_{1},\tilde{L}_{2})\phi_{\tilde{L}}^{\tilde{L}_{1}}(f_{1,\tilde{\bar{Q}}_{1}})\boxtimes \phi_{\tilde{L}}^{\tilde{L}_{2}}(f_{2,\tilde{Q}_{2}}).$$
  Le membre de droite d\'epend d'un choix de param\`etre auxiliaire d\'efinissant l'application $(\tilde{L}_{1},\tilde{L}_{2})\mapsto (\tilde{Q}_{1},\tilde{Q}_{2})$. Pour que la d\'efinition soit correcte, on doit montrer qu'elle est ind\'ependante de ce choix.  Fixons des $\omega$-repr\'esentations $L$-irr\'eductibles et temp\'er\'ees $ \tilde{\pi}_{1}$ et $\tilde{\pi}_{2}$ de $\tilde{L}(F)$ et des \'el\'ements $X_{1},X_{2}\in \tilde{{\cal A}}_{\tilde{L},F}$. Consid\'erons la forme lin\'eaire sur $\underline{I}_{ac}(\tilde{L}(F),\omega)$ d\'efinie par
  $$(\varphi_{1},\varphi_{2})\mapsto \overline{I_{\tilde{L}}(\tilde{\pi}_{1},X_{1},\varphi_{1})}I_{\tilde{L}}(\tilde{\pi}_{2},X_{2},\varphi_{2}).$$
  D'apr\`es 6.4(5), il suffit de prouver que cette forme lin\'eaire prend sur le membre de droite de (2) une valeur qui ne d\'epend pas du choix du param\`etre auxiliaire. D'apr\`es les d\'efinitions des applications $\phi_{\tilde{L}}^{\tilde{L}_{i}}$ pour $i=1,2$, cette valeur est
 $$\sum_{\tilde{L}_{1},\tilde{L}_{1}\in {\cal L}(\tilde{L})}d_{\tilde{L}}^{\tilde{G}}(\tilde{L}_{1},\tilde{L}_{2})\overline{J_{\tilde{L}}^{\tilde{L}_{1}}(\tilde{\pi}_{1},X_{1},f_{1,\tilde{\bar{Q}}_{1}}) }J_{\tilde{L}}^{\tilde{L}_{2}}(\tilde{\pi}_{2},X_{2},f_{2,\tilde{Q}_{2}}).$$
 On peut remplacer $f_{1}$ et $f_{2}$ par leurs produits avec $b\circ\tilde{H}_{\tilde{G}}$, o\`u $b\in C_{c}^{\infty}(\tilde{{\cal A}}_{\tilde{G},F})$ vaut $1$ sur les images de $X_{1}$ et $X_{2}$ dans $\tilde{{\cal A}}_{\tilde{G},F}$. On peut donc supposer $f_{1}$ et $f_{2}$ \`a support compacts. L'expression ci-dessus est alors d\'eduite par transformation de Fourier de la fonction
 $$(\tilde{\lambda}_{1},\tilde{\lambda}_{2})\mapsto\sum_{\tilde{L}_{1},\tilde{L}_{1}\in {\cal L}(\tilde{L})}d_{\tilde{L}}^{\tilde{G}}(\tilde{L}_{1},\tilde{L}_{2})\overline{J_{\tilde{L}}^{\tilde{L}_{1}}(\tilde{\pi}_{1,\tilde{\lambda}_{1}},f_{1,\tilde{\bar{Q}}_{1}}) }J_{\tilde{L}}^{\tilde{L}_{2}}(\tilde{\pi}_{2,\tilde{\lambda}_{2}},f_{2,\tilde{Q}_{2}})$$ 
 sur $i\tilde{{\cal A}}_{\tilde{L},F}^*\times i\tilde{{\cal A}}_{\tilde{L},F}^*$. Il suffit de prouver que cette  fonction ne d\'epend pas du choix du param\`etre auxiliaire. Quitte \`a tordre $\tilde{\pi}_{1}$ et $\tilde{\pi}_{2}$, on est ramen\'e \`a montrer que l'expression  
$$(3) \qquad  \sum_{\tilde{L}_{1},\tilde{L}_{1}\in {\cal L}(\tilde{L})}d_{\tilde{L}}^{\tilde{G}}(\tilde{L}_{1},\tilde{L}_{2})\overline{J_{\tilde{L}}^{\tilde{L}_{1}}(\tilde{\pi}_{1},f_{1,\tilde{\bar{Q}}_{1}}) }J_{\tilde{L}}^{\tilde{L}_{2}}(\tilde{\pi}_{2},f_{2,\tilde{Q}_{2}})$$ 
est ind\'ependante de ce choix. Fixons $\tilde{P}\in {\cal P}(\tilde{L})$.  Pour $i=1,2$, introduisons la $(\tilde{G},\tilde{L})$-famille \`a valeurs op\'erateurs $({\cal M}(\pi_{i};\Lambda,\tilde{Q}))_{\tilde{Q}\in {\cal P}(\tilde{L})}$. Posons
$${\cal M}(\pi_{1}\boxtimes \pi_{2};\Lambda,\tilde{Q})={\cal M}(\pi_{1};\Lambda,\tilde{\bar{Q}})\boxtimes {\cal M}(\pi_{2};\Lambda,\tilde{Q}).$$
De la $(\tilde{G},\tilde{L})$-famille $({\cal M}(\pi_{1}\boxtimes \pi_{2};\Lambda,\tilde{Q}))_{\tilde{Q}\in {\cal P}(\tilde{L})}$ se d\'eduit un op\'erateur ${\cal M}^{\tilde{G}}_{\tilde{L}}(\pi_{1}\boxtimes \pi_{2})$ comme en 2.7. On pose
$$J_{\tilde{L}}^{\tilde{G}}(\tilde{\pi}_{1}\boxtimes \tilde{\pi}_{2},f_{1}\boxtimes f_{2})=trace({\cal M}^{\tilde{G}}_{\tilde{L}}(\pi_{1}\boxtimes \pi_{2})(Ind_{\tilde{P}}^{\tilde{G}}(f_{1})\boxtimes Ind_{\tilde{P}}^{\tilde{G}}(f_{2}))).$$
   Le lemme 5.4(ii) se g\'en\'eralise \`a cette situation: $J_{\tilde{L}}^{\tilde{G}}(\tilde{\pi}_{1}\boxtimes \tilde{\pi}_{2},f_{1}\boxtimes f_{2})$ est \'egal \`a (3). Puisque $J_{\tilde{L}}^{\tilde{G}}(\tilde{\pi}_{1}\boxtimes \tilde{\pi}_{2},f_{1}\boxtimes f_{2})$ ne d\'epend d'aucun choix, il en est de m\^eme de (3), ce que l'on voulait d\'emontrer.
  
  Remarquons que pour $f\in {\cal H}_{ac}(\tilde{G}(F))$, l'image par $\tilde{H}_{\tilde{G}}$ du support $Supp(f)$ de $f$ est ferm\'ee: localement, c'est l'image de $Supp(f(b\circ\tilde{H}_{\tilde{G}}))$ pour une fonction $b$ \`a support compact convenable et ce dernier support est compact.
  On note $ \underline{H}_{ac}(\tilde{G}(F))^1$, resp. $\underline{H}_{ac}(\tilde{G}(F))_{0}$, le sous-espace de $\underline{H}_{ac}(\tilde{G}(F))$ engendr\'e par les fonctions $f_{1}\otimes f_{2}$ telles que
  $$\tilde{H}_{\tilde{G}}(Supp(f_{1}))\cap \tilde{H}_{\tilde{G}}(Supp(f_{2}))$$
  soit compact, resp. vide. On a l'\'egalit\'e
  $$(4) \qquad \underline{H}_{ac}(\tilde{G}(F))^1=\underline{C}(\tilde{G}(F))+\underline{H}_{ac}(\tilde{G}(F))_{0}.$$
   En effet, soient $f_{1},f_{2}\in {\cal H}_{ac}(\tilde{G}(F))$ tels que $\tilde{H}_{\tilde{G}}(Supp(f_{1}))\cap \tilde{H}_{\tilde{G}}(Supp(f_{2}))$ soit  compact. Choisissons une fonction $b'\in C_{c}^{\infty}(\tilde{{\cal A}}_{\tilde{G}})$ valant $1$ sur un voisinage de ce compact. Posons $b''=1-b'$ et, pour $i=1,2$, $f'_{i}=f_{i}(b'\circ\tilde{H}_{\tilde{G}})$, $f''_{i}=f_{i}(b''\circ\tilde{H}_{\tilde{G}})$. On a l'\'egalit\'e
  $$f_{1}\boxtimes f_{2}=(f'_{1}\boxtimes f'_{2})+(f'_{1}\boxtimes f''_{2})+(f''_{1}\boxtimes f'_{2})+(f''_{1}\boxtimes f''_{2}).$$
  Le premier terme appartient \`a $\underline{C}(\tilde{G}(F))$, les trois autres \`a $\underline{H}_{ac}(\tilde{G}(F))_{0}$. D'o\`u (4).

   Notons $\underline{I}_{ac}(\tilde{G}(F),\omega)^1$ l'image de $\underline{H}_{ac}(\tilde{G}(F))^1$ dans $\underline{I}_{ac}(\tilde{G}(F),\omega)$. Montrons que
  
  (5) l'application $\underline{\phi}_{\tilde{L}}$ envoie $ \underline{H}_{ac}(\tilde{G}(F))^1$ dans  $ \underline{I}_{ac}(\tilde{L}(F),\omega)^1$.
  
 Preuve. D'apr\`es (4), il suffit de montrer que $\underline{\phi}_{\tilde{L}}(f_{1}\boxtimes f_{2})\in \underline{I}_{ac}(\tilde{L}(F),\omega)^1$ dans les deux cas suivants
  
  (6) $f_{i}\in C_{c}^{\infty}(\tilde{G}(F),K)$ pour $i=1,2$;
  
  (7) $f_{i}\in {\cal H}_{ac}(\tilde{G}(F))$ pour $i=1,2$ et  $\tilde{H}_{\tilde{G}}(Supp(f_{1}))\cap \tilde{H}_{\tilde{G}}(Supp(f_{2}))=\emptyset$.
   
    Soit $(\tilde{L}_{1},\tilde{L}_{2})$ intervenant dans la formule (2). On pose $\varphi_{1}=\phi_{\tilde{L}}^{\tilde{L}_{1}}(f_{1,\tilde{\bar{Q}}_{1}})$ et $\varphi_{2}=\phi_{\tilde{L}}^{\tilde{L}_{2}}(f_{2,\tilde{Q}_{2}})$. Dans le  cas (6), les fonctions $f_{1,\tilde{\bar{Q}}_{1}} $ et $f_{2,\tilde{Q}_{2}}$ sont \`a supports compacts. D'apr\`es 6.4(10), on peut supposer que  $\tilde{H}_{\tilde{L}_{i}}(Supp(\varphi_{i}))$ est compact pour $i=1,2$.   Il en r\'esulte que $ \tilde{H}_{\tilde{L}}(Supp(\varphi_{1}))\cap \tilde{H}_{\tilde{L}}(Supp(\varphi_{2}))$ a une projection compacte dans $\tilde{{\cal A}}_{\tilde{L}_{i}}$ pour $i=1,2$. Comme dans la preuve de (1), la non-nullit\'e de $d_{\tilde{L}}^{\tilde{G}}(\tilde{L}_{1},\tilde{L}_{2})$ entra\^{\i}ne que l'ensemble $ \tilde{H}_{\tilde{L}}(Supp(\varphi_{1}))\cap \tilde{H}_{\tilde{L}}(Supp(\varphi_{2}))$ lui-m\^eme est compact. Donc $\varphi_{1}\boxtimes \varphi_{2}\in  \underline{H}_{ac}(\tilde{L}(F))^1$. Dans le cas (7), puisque $\tilde{H}_{\tilde{G}}(Supp(f_{i}))$ est ferm\'e pour $i=1,2$, on peut fixer des fonctions $b_{i}$ sur $\tilde{{\cal A}}_{\tilde{G},F}$, lisses, \`a supports disjoints, telles que $b_{i}$ vaille $1$   sur $\tilde{H}_{\tilde{G}}(Supp(f_{i}))$. D'apr\`es 6.4(10), on peut remplacer $\varphi_{i}$ par $\varphi_{i}(b_{i}\circ\tilde{H}_{\tilde{G}})$. Mais alors $ \tilde{H}_{\tilde{G}}(Supp(\varphi_{1}))\cap \tilde{H}_{\tilde{G}}(Supp(\varphi_{2}))$ est vide, a  fortiori $ \tilde{H}_{\tilde{L}}(Supp(\varphi_{1}))\cap \tilde{H}_{\tilde{L}}(Supp(\varphi_{2}))$ l'est et $\varphi_{1}\boxtimes \varphi_{2}$ appartient \`a $ \underline{H}_{ac}(\tilde{L}(F))^1$. Cela prouve (5).
    
    On a d\'efini $\underline{\phi}_{\tilde{L}}$ pour $\tilde{L}\not=\tilde{G}$. On note simplement $\underline{\phi}_{\tilde{G}}$ l'identit\'e de  $\underline{H}_{ac}(\tilde{G}(F))$. 
    
   Les distributions $(f_{1},f_{2})\mapsto J_{g\acute{e}om}^{\tilde{G}}(\omega,f_{1},f_{2})$, $(f_{1},f_{2})\mapsto I_{g\acute{e}om}^{\tilde{G}}(\omega,f_{1},f_{2})$ et celles qui les constituent peuvent \^etre consid\'er\'ees comme des formes lin\'eaires sur  $ \underline{C}(\tilde{G}(F))$.  En fait elles se prolongent \`a l'espace $ \underline{H}_{ac}^1$: pour une telle distribution $D$ et pour $\underline{f}\in \underline{H}_{ac}^1$, on \'ecrit $\underline{f}=\underline{f}_{c}+\underline{f}_{0}$ avec $\underline{f}_{c}\in \underline{C}(\tilde{G}(F))$ et $\underline{f}_{0}\in \underline{H}_{ac}(\tilde{G}(F))_{0}$ et on pose $D(\underline{f})=D(\underline{f}_{c})$. Pour que cette d\'efinition soit loisible, il faut \'evidemment montrer que $D(\underline{f})=0$ si $\underline{f}\in \underline{C}(\tilde{G}(F))\cap \underline{H}_{ac}(\tilde{G}(F))_{0}$. Les distributions en question s'expriment \`a l'aide des distributions basiques $\underline{f}\mapsto J_{\tilde{M}}^{\tilde{G}}(\gamma,\omega,\underline{f})$ ou $\underline{f}\mapsto I_{\tilde{M}}^{\tilde{G}}(\gamma,\omega,\underline{f})$, il suffit donc de traiter celles-ci. Chacune d'elles s'\'etend naturellement \`a $\underline{H}_{ac}(\tilde{G}(F))$ tout entier, il suffit donc de montrer que ces distributions \'etendues annulent $\underline{H}_{ac}(\tilde{G}(F))_{0}$. Soient donc $f_{1}$ et $f_{2} $ v\'erifiant (7). On choisit des fonctions $b_{i}$ comme dans la preuve de (5). Il r\'esulte des d\'efinitions et de 6.5(3) que, pour chacune de nos deux distributions $D$ ci-dessus, on a les \'egalit\'es
   $$D(f_{1}\boxtimes f_{2})=D(f_{1}(b_{1}\circ\tilde{H}_{\tilde{G}})\boxtimes f_{2}(b_{2}\circ\tilde{H}_{\tilde{G}}))=b_{1}(\tilde{H}_{\tilde{G}}(\gamma))b_{2}(\tilde{H}_{\tilde{G}}(\gamma))D(f_{1}\boxtimes f_{2})=0,$$
   ce qu'on voulait d\'emontrer. 
   
   Si on admet la proposition 6.5, le prolongement \`a $\underline{H}_{ac}(\tilde{G}(F))^1$ de la  distribution $(f_{1},f_{2})\mapsto I_{g\acute{e}om}^{\tilde{G}}(\omega,f_{1},f_{2})$ se quotiente en une forme lin\'eaire sur $\underline{I}_{ac}(\tilde{G}(F),\omega)^1$ parce les distributions $\underline{f}\mapsto I_{\tilde{M}}^{\tilde{G}}(\gamma,\omega,\underline{f})$ qui le constituent se quotientent ainsi.
   
   La distribution $(f_{1},f_{2})\mapsto I_{disc}^{\tilde{G}}(\omega,f_{1},f_{2})$ peut elle-aussi \^etre vue comme une forme lin\'eaire sur $\underline{C}(\tilde{G}(F))$. Montrons qu'elle se prolonge \`a $\underline{H}_{ac}(\tilde{G}(F))^1$ par le m\^eme proc\'ed\'e que ci-dessus et que ce prolongement se quotiente en une forme lin\'eaire sur  l'espace $\underline{I}_{ac}(\tilde{G}(F),\omega)^1$. La distribution en question est combinaison lin\'eaire de distributions
   $$(f_{1},f_{2})\mapsto \int_{i{\cal A}_{\tilde{G},F}^*}\overline{I_{\tilde{G}}(\tilde{\pi}_{\tilde{\lambda}},f_{1})}I_{\tilde{G}}(\tilde{\pi}_{\tilde{\lambda}},f_{2})\,d\lambda,$$
o\`u $\tilde{\pi}$ est une $\omega$-repr\'esentation temp\'er\'ee de $\tilde{G}(F)$. On transforme celles-ci par inversion de Fourier en
$$(f_{1},f_{2})\mapsto mes(i{\cal A}_{\tilde{G},F}^*)\int_{\tilde{{\cal A}}_{\tilde{L},F}}\overline{I_{\tilde{G}}(\tilde{\pi},X,f_{1})}I_{\tilde{G}}(\tilde{\pi},X,f_{2})\,dX.$$
 On montre gr\^ace \`a 6.4(4) qu'une distribution
$$(8) \qquad (f_{1},f_{2})\mapsto  \overline{I_{\tilde{G}}(\tilde{\pi},X,f_{1})}I_{\tilde{G}}(\tilde{\pi},X,f_{2})$$
annule $\underline{C}(\tilde{G}(F))\cap \underline{H}_{ac}(\tilde{G}(F))_{0}$. Comme pour les distributions "g\'eom\'etriques", cela permet de prolonger la distribution $(f_{1},f_{2})\mapsto I_{disc}^{\tilde{G}}(\omega,f_{1},f_{2})$ \`a $\underline{H}_{ac}(\tilde{G}(F))^1$. Ce prolongement se quotiente en une forme lin\'eaire sur $\underline{I}_{ac}(\tilde{G}(F),\omega)^1$ parce les distributions (8) qui le constituent se quotientent ainsi.

Remarquons que, si l'\'enonc\'e du th\'eor\`eme est vrai, il s'\'etend en l'\'egalit\'e
$$I_{g\acute{e}om}^{\tilde{G}}(\omega,\underline{f})=I_{disc}^{\tilde{G}}(\omega,\underline{f})$$
pour tout $\underline{f}\in \underline{H}_{ac}(\tilde{G}(F))^1$ (ou $\underline{f}\in \underline{I}_{ac}(\tilde{G}(F),\omega)^1$ si on admet la proposition 6.5) puisque chaque terme est par d\'efinition le m\^eme terme \'evalu\'e sur $\underline{f}_{c}$ o\`u, comme plus haut $\underline{f}_{c}$ est un \'el\'ement de $\underline{C}(\tilde{G}(F))$ tel que $\underline{f}\in \underline{f}_{c}+\underline{H}_{ac}(\tilde{G}(F))_{0}$.  

Venons-en \`a la preuve du th\'eor\`eme. Elle se fait par r\'ecurrence sur $rg_{ss}(G)$. On suppose v\'erifi\'es le th\'eor\`eme et la proposition 6.5 pour les espaces $\tilde{G}'$ tels que $rg_{ss}(G')<rg_{ss}(G)$. Soient $f_{i}\in C_{c}^{\infty}(\tilde{G}(F),K)$ pour $i=1,2$, posons $\underline{f}=f_{1}\boxtimes f_{2}$. Soit $\tilde{L}\in {\cal L}(\tilde{M}_{0})$. Montrons que
$$(9)\qquad J_{\tilde{L},spec}^{\tilde{G}}(\omega,f_{1},f_{2})=I^{\tilde{L}}_{disc}(\omega,\underline{\phi}_{\tilde{L}}(\underline{f})).$$
C'est tautologique si $\tilde{L}=\tilde{G}$. Supposons $\tilde{L}\not=\tilde{G}$. 
Les deux c\^ot\'es sont des combinaisons lin\'eaires index\'ees par  $\tau\in(E_{disc}(\tilde{L},\omega)/conj)/i{\cal A}_{\tilde{L},F}^*$ de produits des m\^emes coefficients et de certaines distributions. La distribution qui intervient dans le membre de gauche est
$$(10)\qquad \int_{i{\cal A}_{\tilde{L},F}^*}J_{\tilde{L}}^{\tilde{G}}(\pi_{\tau_{\lambda}},f_{1},f_{2})\,d\lambda.$$
Relevons $\tau$ en un \'el\'ement $\boldsymbol{\tau}\in {\cal E}(\tilde{L},\omega)$. Pour $X\in \tilde{{\cal A}}_{\tilde{L},F}$, notons $\underline{\varphi}\mapsto I_{\tilde{L}}(\pi_{\tau},X,\underline{\varphi})$ le prolongement \`a $\underline{H}_{ac}(\tilde{G})^1$ de la distribution 
$$(\varphi_{1},\varphi_{2})\mapsto \overline{I_{\tilde{L}}(\tilde{\pi}_{\boldsymbol{\tau}},X,\varphi_{1})}I_{\tilde{L}}(\tilde{\pi}_{\boldsymbol{\tau}},X,\varphi_{2}).$$
Alors la distribution intervenant dans le membre de droite de (9) est
$$(11)\qquad mes(i{\cal A}_{\tilde{L},F}^*)\int_{\tilde{{\cal A}}_{\tilde{L},F}}I_{\tilde{L}}(\pi_{\tau},X,\underline{\phi}_{\tilde{L}}(\underline{f}))\,dX.$$
Il faut montrer que les expressions (10) et (11) sont \'egales. Par transformation de Fourier, le lemme 5.4(ii) entra\^{\i}ne que (10) est \'egal \`a
$$mes(i{\cal A}_{\tilde{L},F}^*)\sum_{\tilde{L}_{1},\tilde{L}_{2}\in {\cal L}(\tilde{L})}d_{\tilde{L}}^{\tilde{G}}(\tilde{L}_{1},\tilde{L}_{2})\int_{\tilde{{\cal A}}_{\tilde{L},F}}\overline{J_{\tilde{L}}^{\tilde{L}_{1}}(\tilde{\pi}_{\boldsymbol{\tau}},X,f_{1,\tilde{\bar{Q}}_{1}})}J_{\tilde{L}}^{\tilde{L}_{2}}(\tilde{\pi}_{\boldsymbol{\tau}},X,f_{2,\tilde{Q}_{2}})\,dX.$$
D'apr\`es les d\'efinitions des applications $\phi_{\tilde{L}}^{\tilde{L}_{i}}$ pour $i=1,2$, c'est aussi
 $$mes(i{\cal A}_{\tilde{L},F}^*)\sum_{\tilde{L}_{1},\tilde{L}_{2}\in {\cal L}(\tilde{L})}d_{\tilde{L}}^{\tilde{G}}(\tilde{L}_{1},\tilde{L}_{2})\int_{\tilde{{\cal A}}_{\tilde{L},F}}\overline{I_{\tilde{L}}(\tilde{\pi}_{\boldsymbol{\tau}},X,\phi_{\tilde{L}}^{\tilde{L}_{1}}(f_{1,\tilde{\bar{Q}}_{1}}))}I_{\tilde{L}}(\tilde{\pi}_{\boldsymbol{\tau}},X,\phi_{\tilde{L}}^{\tilde{L}_{2}}(f_{2,\tilde{Q}_{2}}))\,dX.$$
 L'\'egalit\'e de cette expression avec (11) r\'esulte alors de la d\'efinition de $\underline{\phi}_{\tilde{L}}$. Cela prouve (9).
   
 Pour $\tilde{L}\not=\tilde{G}$, on peut par l'hypoth\`ese de r\'ecurrence utiliser le th\'eor\`eme prolong\'e comme indiqu\'e ci-dessus: on a  
   $$I^{\tilde{L}}_{disc}(\omega,\underline{\phi}_{\tilde{L}}(\underline{f}))=I^{\tilde{L}}_{g\acute{e}om}(\omega,\underline{\phi}_{\tilde{L}}(\underline{f})).$$
  Alors (9) implique
  $$J_{spec}^{\tilde{G}}(\omega,f_{1},f_{2})=I^{\tilde{G}}_{disc}(\omega,f_{1},f_{2})-I_{g\acute{e}om}^{\tilde{G}}(\omega,f_{1},f_{2})+X,$$
  o\`u
  $$X=\sum_{\tilde{L}\in {\cal L}(\tilde{M}_{0})}\vert \tilde{W}^{L}\vert \vert \tilde{W}^G\vert ^{-1}(-1)^{a_{\tilde{L}}-a_{\tilde{G}}}I^{\tilde{L}}_{g\acute{e}om}(\underline{\phi}_{\tilde{L}}(\underline{f}))$$
  (rappelons que, par convention,  $\underline{\phi}_{\tilde{G}}$ l'identit\'e de  $\underline{H}_{ac}(\tilde{G}(F))$). 
   En utilisant le th\'eor\`eme 5.1, la relation cherch\'ee
  $$I_{g\acute{e}om}^{\tilde{G}}(\omega,f_{1},f_{2})=I_{disc}^{\tilde{G}}(\omega,f_{1},f_{2})$$
  \'equivaut \`a l'\'egalit\'e
  $$(12) \qquad J_{g\acute{e}om}^{\tilde{G}}(\omega,f_{1},f_{2})=X.$$
 Par d\'efinition
 $$X=\sum_{\tilde{L}\in {\cal L}(\tilde{M}_{0})}\vert \tilde{W}^{L}\vert \vert \tilde{W}^G\vert ^{-1}(-1)^{a_{\tilde{L}}-a_{\tilde{G}}}\sum_{\tilde{M}\in {\cal L}^{\tilde{L}}(\tilde{M}_{0})}\vert \tilde{W}^M\vert \vert \tilde{W}^L\vert ^{-1}(-1)^{a_{\tilde{M}}-a_{\tilde{L}}}I_{\tilde{M},g\acute{e}om}^{\tilde{L}}(\omega,\underline{\phi}_{\tilde{L}}(\underline{f}))$$
 $$=\sum_{\tilde{M}\in {\cal L}(\tilde{M}_{0})}\vert \tilde{W}^M\vert \vert \tilde{W}^G\vert ^{-1}(-1)^{a_{\tilde{M}}-a_{\tilde{G}}}X_{\tilde{M}},$$
 o\`u
 $$X_{\tilde{M}}=\sum_{\tilde{L}\in {\cal L}(\tilde{M})}I_{\tilde{M},g\acute{e}om}^{\tilde{L}}(\omega,\underline{\phi}_{\tilde{L}}(\underline{f})).$$
 En se rappelant la d\'efinition de $J_{g\acute{e}om}^{\tilde{G}}(\omega,f_{1},f_{2})$, il suffit pour d\'emontrer (12) de fixer $\tilde{M}\in {\cal L}(\tilde{M}_{0})$ et de prouver l'\'egalit\'e
  $$(13)\qquad  J_{\tilde{M},g\acute{e}om}^{\tilde{G}}(\omega,f_{1},f_{2})=X_{\tilde{M}}.$$
  Les deux membres de cette \'egalit\'e sont des sommes sur $\tilde{S}\in T_{ell}(\tilde{M},\omega)$  de coefficients (qui sont les m\^emes pour les deux membres) et d'int\'egrales sur $\tilde{S}(F)/(1-\theta)(S(F))$ de certaines fonctions. Il suffit de fixer $\tilde{S}$ et de prouver que les fonctions sont les m\^emes. Fixons donc $\tilde{S}$ et un point $\gamma\in \tilde{S}(F)\cap \tilde{G}_{reg}(F)$. La valeur en $\gamma$ de la fonction relative au membre de gauche de (13) est
  $$\sum_{\tilde{M}_{1},\tilde{M}_{2}\in{\cal L}(\tilde{M})}d_{\tilde{M}}^{\tilde{G}}(\tilde{M}_{1},\tilde{M}_{2})\overline{J_{\tilde{M}}^{\tilde{M}_{1}}(\gamma,\omega,f_{1,\tilde{\bar{P}}_{1}})}J_{\tilde{M}}^{\tilde{M}_{2}}(\gamma,\omega,f_{2,\tilde{\bar{P}}_{2}}),$$
  cela d'apr\`es le lemme 5.4(i). On peut exprimer les int\'egrales orbitales pond\'er\'ees \`a l'aide d'int\'egrales \'equivariantes gr\^ace \`a la proposition 6.5. On obtient
  $$\sum_{\tilde{M}_{1},\tilde{M}_{2}\in{\cal L}(\tilde{M})}d_{\tilde{M}}^{\tilde{G}}(\tilde{M}_{1},\tilde{M}_{2})\sum_{\tilde{L}_{1}\in {\cal L}^{\tilde{M}_{1}}(\tilde{M})}\sum_{\tilde{L}_{2}\in {\cal L}^{\tilde{M}_{2}}(\tilde{M})}\overline{I_{\tilde{M}}^{\tilde{L}_{1}}(\gamma,\omega,\phi_{\tilde{L}_{1}}^{\tilde{M}_{1}}(f_{1,\tilde{\bar{P}}_{1}})}I_{\tilde{M}}^{\tilde{L}_{2}}(\gamma,\omega,\phi_{\tilde{L}_{2}}^{\tilde{M}_{2}}(f_{2,\tilde{\bar{P}}_{2}})),$$
  en posant la convention que $\phi_{\tilde{G}}^{\tilde{G}}$ est l'identit\'e de $C_{c}^{\infty}(\tilde{G}(F),K)$.
Si $\tilde{L}_{1}$ et $\tilde{L}_{2}$ sont tous deux diff\'erents de $\tilde{G}$, la distribution 
 $$(\varphi_{1},\varphi_{2})\to \overline{I_{\tilde{M}}^{\tilde{L}_{1}}(\gamma,\omega, \varphi_{1})}I_{\tilde{M}}^{\tilde{L}_{2}}(\gamma, \omega,\varphi_{2})$$
 peut \^etre consid\'er\'ee comme une forme lin\'eaire sur $I_{ac}(\tilde{L}_{1}(F),\omega)\boxtimes I_{ac}(\tilde{L}_{2}(F),\omega)$. Notons-la $\underline{\varphi}\mapsto
 I_{\tilde{M}}^{\tilde{L}_{1},\tilde{L}_{2}}(\gamma,\omega,\underline{\varphi})$. Dans le cas o\`u l'un des $\tilde{L}_{i}$ est \'egal \`a $\tilde{G}$ (ou les deux), on utilise la m\^eme notation en rempla\c{c}ant la composante $I_{ac}(\tilde{L}_{i}(F),\omega)$ de l'espace de d\'epart par ${\cal H}_{ac}(\tilde{G}(F))$.
 Alors l'expression pr\'ec\'edente s'\'ecrit
 $$(14) \qquad \sum_{\tilde{L}_{1},\tilde{L}_{2}\in {\cal L}(\tilde{M})}I_{\tilde{M}}^{\tilde{L}_{1},\tilde{L}_{2}}(\gamma,\omega,\underline{\varphi}(\tilde{L}_{1},\tilde{L}_{2})),$$
 o\`u
 $$\underline{\varphi}(\tilde{L}_{1},\tilde{L}_{2})=\sum_{\tilde{M}_{1}\in {\cal L}(\tilde{L}_{1}),\tilde{M}_{2}\in {\cal L}(\tilde{L}_{2})}d_{\tilde{M}}^{\tilde{G}}(\tilde{M}_{1},\tilde{M}_{2})\phi_{\tilde{L}_{1}}^{\tilde{M}_{1}}(f_{1,\tilde{\bar{P}}_{1}})\boxtimes \phi_{\tilde{L}_{2}}^{\tilde{M}_{2}}(f_{2,\tilde{P}_{2}}).$$
 
 La valeur de la fonction intervenant dans le membre de droite de (13) est
 $$\sum_{\tilde{L}\in {\cal L}(\tilde{M})}I_{\tilde{M}}^{\tilde{L}}(\gamma,\omega,\underline{\phi}_{\tilde{L}}(\underline{f})).$$
 D'apr\`es les d\'efinitions, cela s'\'ecrit encore
 $$\sum_{\tilde{L}\in {\cal L}(\tilde{M})}\sum_{\tilde{L}_{1},\tilde{L}_{2}\in {\cal L}^{\tilde{L}}(\tilde{M})}d_{\tilde{M}}^{\tilde{L}}(\tilde{L}_{1},\tilde{L}_{2})I_{\tilde{M}}^{\tilde{L}_{1},\tilde{L}_{2}}(\gamma,\omega,(\underline{\phi}_{\tilde{L}}(\underline{f}))_{\tilde{L}_{1},\tilde{L}_{2}}),$$
 o\`u on a not\'e $\underline{\varphi}\mapsto \underline{\varphi}_{\tilde{L}_{1},\tilde{L}_{2}}$ l'application lin\'eaire de $\underline{H}_{ac}(\tilde{L}(F))$ dans $I_{ac}(\tilde{L}_{1}(F),\omega)\boxtimes I_{ac}(\tilde{L}_{2}(F),\omega)$ qui envoie $\varphi_{1}\boxtimes \varphi_{2}$ sur $ \varphi_{1,\tilde{L}_{1}}\boxtimes \varphi_{2,\tilde{L}_{2}}$. Ici encore, on doit modifier la d\'efinition si l'un des $\tilde{L}_{i}$ est \'egal \`a $\tilde{G}$: on remplace $I_{ac}(\tilde{L}_{i}(F),\omega)$ par ${\cal H}_{ac}(\tilde{G}(F))$ et $\varphi_{i,\tilde{L}_{i}}$ par $\varphi_{i}$. L'expression ci-dessus s'\'ecrit encore
 $$(15)\qquad \sum_{\tilde{L}_{1},\tilde{L}_{2}\in {\cal L}(\tilde{M})}I_{\tilde{M}}^{\tilde{L}_{1},\tilde{L}_{2}}(\gamma,\omega,\underline{\varphi}'(\tilde{L}_{1},\tilde{L}_{2})),$$
 o\`u
 $$\underline{\varphi}'(\tilde{L}_{1},\tilde{L}_{2})=\sum_{\tilde{L}\in {\cal L}(\tilde{M}); \tilde{L}_{1},\tilde{L}_{2}\subset \tilde{L}}d_{\tilde{M}}^{\tilde{L}}(\tilde{L}_{1},\tilde{L}_{2})(\underline{\phi}_{\tilde{L}}(\underline{f}))_{\tilde{L}_{1},\tilde{L}_{2}}.$$
 On doit prouver que les expressions (14) et (15) sont \'egales. Il suffit de fixer $\tilde{L}_{1}$ et $\tilde{L}_{2}$ et de prouver l'\'egalit\'e
 $$(16) \qquad  \underline{\varphi}(\tilde{L}_{1},\tilde{L}_{2})= \underline{\varphi}'(\tilde{L}_{1},\tilde{L}_{2}).$$
 Fixons donc $\tilde{L}_{1}$ et $\tilde{L}_{2}$.  On voit d'abord que les deux membres sont nuls si la condition
 $$(17) \qquad {\cal A}_{\tilde{M}}^{\tilde{L}_{1}}\cap {\cal A}_{\tilde{M}}^{\tilde{L}_{2}}=\{0\}$$
 n'est pas v\'erifi\'ee. En effet, dans ce cas, il n'y a pas de couples $(\tilde{M}_{1},\tilde{M}_{2})$  intervenant dans la d\'efinition de $\underline{\varphi}(\tilde{L}_{1},\tilde{L}_{2})$ pour lesquels on ait $d_{\tilde{M}}^{\tilde{G}}(\tilde{M}_{1},\tilde{M}_{2})\not=0$ et il n'y a pas de $\tilde{L}$ intervenant dans la d\'efinition de $\underline{\varphi}'(\tilde{L}_{1},\tilde{L}_{2})$ pour lequel $d_{\tilde{M}}^{\tilde{L}}(\tilde{L}_{1},\tilde{L}_{2})\not=0$. On suppose donc que  (17) est v\'erifi\'ee. Il y a alors un unique $\tilde{L}$ intervenant dans la d\'efinition de $\underline{\varphi}'(\tilde{L}_{1},\tilde{L}_{2})$ pour lequel $d_{\tilde{M}}^{\tilde{L}}(\tilde{L}_{1},\tilde{L}_{2})\not=0$: celui pour lequel 
 $${\cal A}_{\tilde{L}}={\cal A}_{\tilde{L}_{1}}\cap {\cal A}_{\tilde{L}_{2}},$$
 cf. 5.3. Dans la suite $\tilde{L}$ d\'esigne cet ensemble de Levi. 
 
 Supposons $\tilde{L}_{1}=\tilde{L}_{2}=\tilde{G}$. L'\'egalit\'e (17) entra\^{\i}ne $\tilde{M}=\tilde{G}$ et on v\'erifie que les deux membres de (16) sont simplement \'egaux \`a $f_{1}\boxtimes f_{2}$. Supposons par exemple $\tilde{L}_{1}=\tilde{G}$ et $\tilde{L}_{2}\not=\tilde{G}$. L'\'egalit\'e (17) entra\^{\i}ne $\tilde{L}_{2}=\tilde{M}$. On a aussi $\tilde{L}=\tilde{G}$. D'o\`u $\underline{\varphi}'(\tilde{G},\tilde{M})=f_{1}\boxtimes f_{2,\tilde{M}}$. Seul le couple $(\tilde{M}_{1},\tilde{M}_{2})=(\tilde{G},\tilde{M})$ contribue \`a la d\'efinition de $\underline{\varphi}(\tilde{G},\tilde{M})$, d'o\`u $\underline{\varphi}(\tilde{G},\tilde{M})=f_{1}\boxtimes \phi_{\tilde{M}}^{\tilde{M}}(f_{2,\tilde{P}_{2}})=f_{1}\boxtimes f_{2,\tilde{M}}$, d'o\`u (16) dans ce cas.

On suppose maintenant $\tilde{L}_{1}$ et $\tilde{L}_{2}$ tous deux diff\'erents de $\tilde{G}$.  Les deux termes  de (16) appartiennent \`a  $I_{ac}(\tilde{L}_{1}(F),\omega)\boxtimes I_{ac}(\tilde{L}_{2}(F),\omega)$. Pour prouver (16), on peut fixer pour $i=1,2$  une $\omega$-repr\'esentation temp\'er\'ee  $\tilde{\pi}_{i}$ de $\tilde{L}_{i}(F)$ et un \'el\'ement $X_{i}\in \tilde{{\cal A}}_{\tilde{L}_{i}}$ et prouver que la forme lin\'eaire
 $$(18)\qquad \varphi_{1}\boxtimes \varphi_{2}\mapsto\overline{ I_{\tilde{L}_{1}}(\tilde{\pi}_{1},X_{1},\varphi_{1})}I_{\tilde{L}_{2}}(\tilde{\pi}_{2},X_{2},\varphi_{2})$$
 prend la m\^eme valeur sur les deux membres. Sa valeur sur $\underline{\varphi}(\tilde{L}_{1},\tilde{L}_{2})$ est
 $$\sum_{\tilde{M}_{1}\in {\cal L}(\tilde{L}_{1}),\tilde{M}_{2}\in {\cal L}(\tilde{L}_{2})}d_{\tilde{M}}^{\tilde{G}}(\tilde{M}_{1},\tilde{M}_{2})\overline{ I_{\tilde{L}_{1}}(\tilde{\pi}_{1},X_{1},\phi_{\tilde{L}_{1}}^{\tilde{M}_{i}}(f_{1,\tilde{\bar{P}}_{1}}))}I_{\tilde{L}_{2}} (\tilde{\pi}_{2},X_{2},\phi_{\tilde{L}_{2}}^{\tilde{M}_{i}}(f_{2,\tilde{\bar{P}}_{2}})).$$
 En utilisant les d\'efinitions des applications $\phi_{\tilde{L}_{i}}^{\tilde{M}_{i}}$, on obtient
$$\sum_{\tilde{M}_{1}\in {\cal L}(\tilde{L}_{1}),\tilde{M}_{2}\in {\cal L}(\tilde{L}_{2})}d_{\tilde{M}}^{\tilde{G}}(\tilde{M}_{1},\tilde{M}_{2})\overline{ J_{\tilde{L}_{1}}^{\tilde{M}_{i}}(\tilde{\pi}_{1},X_{1}, f_{1,\tilde{\bar{P}}_{1}})}J_{\tilde{L}_{2}}^{\tilde{M}_{i}}(\tilde{\pi}_{2},X_{2}, f_{2,\tilde{\bar{P}}_{2}}).$$
 C'est la valeur en $(-X_{1},X_{2})$ (le signe $-$ provenant de la conjugaison complexe figurant dans la formule ci-dessus) de la transform\'ee de Fourier de la fonction sur $i\tilde{{\cal A}}_{\tilde{L}_{1},F}^*\times i\tilde{{\cal A}}_{\tilde{L}_{2},F}^*$ qui, \`a $(\tilde{\lambda}_{1},\tilde{\lambda}_{2})$, associe
$$ \sum_{\tilde{M}_{1}\in {\cal L}(\tilde{L}_{1}),\tilde{M}_{2}\in {\cal L}(\tilde{L}_{2})}d_{\tilde{M}}^{\tilde{G}}(\tilde{M}_{1},\tilde{M}_{2})\overline{ J_{\tilde{L}_{1}}^{\tilde{M}_{i}}(\tilde{\pi}_{1,\tilde{\lambda}_{1}}, f_{1,\tilde{\bar{P}}_{1}})}J_{\tilde{L}_{2}}^{\tilde{M}_{i}}(\tilde{\pi}_{2, \tilde{\lambda}_{2}},f_{2,\tilde{\bar{P}}_{2}}).$$
Fixons des espaces paraboliques $\tilde{Q}_{i}\in {\cal P}(\tilde{L}_{i})$ pour $i=1,2$  contenus dans des espaces paraboliques d'espaces de Levi $\tilde{L}$ et introduisons les repr\'esentations induites $\tilde{\Pi}_{i,\tilde{\lambda}_{i}}=Ind_{\tilde{Q}_{i}}^{\tilde{G}}(\tilde{\pi}_{i,\tilde{\lambda}_{i}})$ dans l'espace $V_{i}=V_{\pi_{i},\tilde{Q}_{i}}$. On dispose des $(\tilde{G},\tilde{L}_{i})$-familles \`a valeurs op\'erateurs $({\cal M}(\pi_{i,\lambda_{i}};\Lambda,\tilde{Q}))_{\tilde{Q}\in {\cal P}(\tilde{L}_{i})}$. On v\'erifie que l'expression pr\'ec\'edente n'est autre que
$$(19)\qquad  trace({\cal N}(\lambda_{1},\lambda_{2})(\tilde{\Pi}_{1,\tilde{\lambda}_{1}}(f_{1})\boxtimes \tilde{\Pi}_{2,\tilde{\lambda}_{2}}(f_{2})),$$
o\`u ${\cal N}(\lambda_{1},\lambda_{2})$ est l'op\'erateur de $V_{1}\boxtimes V_{2}$ d\'efini par
$${\cal N}(\lambda_{1},\lambda_{2})= \sum_{\tilde{M}_{1}\in {\cal L}(\tilde{L}_{1}),\tilde{M}_{2}\in {\cal L}(\tilde{L}_{2})}d_{\tilde{M}}^{\tilde{G}}(\tilde{M}_{1},\tilde{M}_{2}){\cal M}_{\tilde{L}_{1}}^{\tilde{\bar{P}}_{1}}(\pi_{1,\lambda_{1}})\boxtimes {\cal M}_{\tilde{L}_{2}}^{\tilde{P}_{2}}(\pi_{2,\lambda_{2}}).$$
Pour calculer la valeur sur $\underline{\varphi}'(\tilde{L}_{1},\tilde{L}_{2})$ de la forme lin\'eaire (18), on doit d'abord calculer $I_{\tilde{L}_{i}}(\tilde{\pi}_{i},X_{i},\varphi_{i,\tilde{L}_{i}})$ pour $\varphi_{i}\in {\cal H}_{ac}(\tilde{L}(F))$. On introduit une fonction $b\in C_{c}^{\infty}(\tilde{{\cal A}}_{\tilde{L},F})$ telle que $b(\tilde{H}_{\tilde{L}}(X_{i}))=1$. On a alors
$$I_{\tilde{L}_{i}}(\tilde{\pi}_{i},X_{i},\varphi_{i,\tilde{L}_{i}})=I_{\tilde{L}_{i}}(\tilde{\pi}_{i},X_{i},\varphi_{i,\tilde{L}_{i}}(b\circ\tilde{H}_{\tilde{L}}))=I_{\tilde{L}_{i}}(\tilde{\pi}_{i},X_{i},(\varphi_{i} (b\circ\tilde{H}_{\tilde{L}}))_{\tilde{L}_{i}}).$$
Maintenant, $\varphi_{i} (b\circ\tilde{H}_{\tilde{L}})$ est \`a support compact, d'o\`u
$$I_{\tilde{L}_{i}}(\tilde{\pi}_{i},X_{i},(\varphi_{i} (b\circ\tilde{H}_{\tilde{L}}))_{\tilde{L}_{i}})=mes( i{\cal A}_{\underline{\tilde{M}}_{i},F}^*)^{-1}\int_{i{\cal A}_{\tilde{L}_{i},F}^*}
I_{\tilde{L}_{i}}(\tilde{\pi}_{i,\tilde{\lambda}_{i}},(\varphi_{i} (b\circ\tilde{H}_{\tilde{L}}))_{\tilde{L}_{i}})e^{-<\tilde{\lambda}_{i},X_{i}>}\,d\lambda_{i}.$$
En introduisant la repr\'esentation $\tilde{\pi}_{i,\tilde{\lambda}_{i}}^{\tilde{L}}=Ind_{\tilde{Q}_{i}\cap \tilde{L}}^{\tilde{L}}(\tilde{\pi}_{i,\tilde{\lambda}_{i}})$, l'expression ci-dessus devient
$$mes( i{\cal A}_{\underline{\tilde{M}}_{i},F}^*)^{-1}\int_{i{\cal A}_{\tilde{L}_{i},F}^*}
I_{\tilde{L}}(\tilde{\pi}^{\tilde{L}}_{i,\tilde{\lambda}_{i}},\varphi_{i} (b\circ\tilde{H}_{\tilde{L}}))e^{-<\tilde{\lambda}_{i},X_{i}>}\,d\lambda_{i}.$$
On peut aussi l'\'ecrire
$$ mes( i{\cal A}_{\underline{\tilde{M}}_{i},F}^*)^{-1}\int_{i{\cal A}_{\tilde{L}_{i},F}^*}\int_{\tilde{{\cal A}}_{\tilde{L},F}}I_{\tilde{L}}(\tilde{\pi}^{\tilde{L}}_{i,\tilde{\lambda}_{i}},Y_{i},\varphi_{i} (b\circ\tilde{H}_{\tilde{L}}))\,dY_{i}\,e^{-<\tilde{\lambda}_{i},X_{i}>}\,d\lambda_{i},$$
puis, par inversion de Fourier,
$$ mes( i{\cal A}_{\underline{\tilde{M}}_{i},F}^*)^{-1}mes(i{\cal A}_{\tilde{L},F}^*)\int_{i{\cal A}_{\tilde{L}_{i},F}^*/i{\cal A}_{\tilde{L},F}^*}I_{\tilde{L}}(\tilde{\pi}^{\tilde{L}}_{i,\tilde{\lambda}_{i}},X_{i,\tilde{L}},\varphi_{i} (b\circ\tilde{H}_{\tilde{L}}))e^{-<\tilde{\lambda}_{i},X_{i}>}\,d\lambda_{i},$$
o\`u $X_{i,\tilde{L}}$ est la projection de $X_{i}$ dans $\tilde{{\cal A}}_{\tilde{L},F}$. Maintenant, on peut faire dispara\^{\i}tre la fonction $(b\circ\tilde{H}_{\tilde{L}})$ et on obtient simplement
 $$ mes( i{\cal A}_{\underline{\tilde{M}}_{i},F}^*)^{-1}mes(i{\cal A}_{\tilde{L},F}^*)\int_{i{\cal A}_{\tilde{L}_{i},F}^*/i{\cal A}_{\tilde{L},F}^*}I_{\tilde{L}}(\tilde{\pi}^{\tilde{L}}_{i,\tilde{\lambda}_{i}},X_{i,\tilde{L}},\varphi_{i}  )e^{-<\tilde{\lambda}_{i},X_{i}>}\,d\lambda_{i}.$$
 En appliquant ce calcul et les d\'efinitions, on obtient que la valeur sur $\underline{\varphi}'(\tilde{L}_{1},\tilde{L}_{2})$ de la forme lin\'eaire (18) est
 $$(20) \qquad mes( i{\cal A}_{\underline{\tilde{M}}_{1},F}^*)^{-1} mes( i{\cal A}_{\underline{\tilde{M}}_{2},F}^*)^{-1}mes(i{\cal A}_{\tilde{L},F}^*)^2\int_{i{\cal A}_{\tilde{L}_{1},F}^*/i{\cal A}_{\tilde{L},F}^*}\int_{i{\cal A}_{\tilde{L}_{2},F}^*/i{\cal A}_{\tilde{L},F}^*}$$
 $$I_{\tilde{L}}(\tilde{\pi}^{\tilde{L}}_{1,\tilde{\lambda}_{1}}\boxtimes \tilde{\pi}^{\tilde{L}}_{2,\tilde{\lambda}_{2}},X_{1,\tilde{L}},X_{2,\tilde{L}},\underline{\phi}_{\tilde{L}}(\underline{f}))e^{<\tilde{\lambda}_{1},X_{1}>-<\tilde{\lambda}_{2},X_{2}>}\,d\lambda_{1}\,d\lambda_{2},$$
 o\`u on a not\'e $\underline{\varphi}\mapsto I_{\tilde{L}}(\tilde{\pi}^{\tilde{L}}_{1,\tilde{\lambda}_{1}}\boxtimes \tilde{\pi}^{\tilde{L}}_{2,\tilde{\lambda}_{2}},X_{1,\tilde{L}},X_{2,\tilde{L}},\underline{\varphi})$ la forme lin\'eaire
 $$\varphi_{1}\boxtimes \varphi_{2}\mapsto \overline{I_{\tilde{L}}(\tilde{\pi}^{\tilde{L}}_{1,\tilde{\lambda}_{1}},X_{1,\tilde{L}},\varphi_{1}  )}I_{\tilde{L}}(\tilde{\pi}^{\tilde{L}}_{2,\tilde{\lambda}_{2}},X_{2,\tilde{L}},\varphi_{2}  )$$
 sur $\underline{H}_{ac}(\tilde{L}(F))$. Mais on a calcul\'e la compos\'ee de $\underline{\phi}_{\tilde{L}}$ avec cette forme lin\'eaire dans la preuve suivant la relation (2). Le r\'esultat est le suivant. On d\'eduit des  $(\tilde{G},\tilde{L}_{i})$-familles  $({\cal M}(\pi_{i,\lambda'_{i}};\Lambda,\tilde{Q}))_{\tilde{Q}\in {\cal P}(\tilde{L}_{i})}$ (o\`u $\lambda'_{i}\in i{\cal A}_{\tilde{L}_{i},F}^*$) des $(\tilde{G},\tilde{L})$-familles et on introduit leur produit sesquilin\'eaire
 $${\cal M}(\pi_{1,\lambda'_{1}}\boxtimes \pi_{2,\lambda'_{2}};\Lambda,\tilde{Q})={\cal M}(\pi_{1,\lambda'_{1}};\Lambda,\tilde{Q})\boxtimes {\cal M}(\pi_{2,\lambda'_{2}};\Lambda,\tilde{Q})$$
 pour $\tilde{Q}\in {\cal P}(\tilde{L})$. On pose
 $${\cal N}'(\lambda'_{1},\lambda'_{2})={\cal M}^{\tilde{G}}_{\tilde{L}}(\pi_{1,\lambda'_{1}}\boxtimes\pi_{2,\lambda'_{2}};0).$$
 Alors $I_{\tilde{L}}(\tilde{\pi}^{\tilde{L}}_{1,\tilde{\lambda}_{1}}\boxtimes \tilde{\pi}^{\tilde{L}}_{2,\tilde{\lambda}_{2}},X_{1,\tilde{L}},X_{2,\tilde{L}},\underline{\phi}_{\tilde{L}}(\underline{f}))$ est la valeur en $(-X_{1},X_{2})$ de la transform\'ee de Fourier de la fonction
 $$(\tilde{\mu}_{1},\tilde{\mu}_{2})\mapsto trace({\cal N}'(\lambda_{1}+\mu_{1},\lambda_{2}+\mu_{2})(\tilde{\Pi}_{1,\tilde{\lambda}_{1}+\tilde{\mu}_{1}}(f_{1})\boxtimes \tilde{\Pi}_{2,\tilde{\lambda}_{2}+\tilde{\mu}_{2}}(f_{2})))$$
 sur $i\tilde{{\cal A}}_{\tilde{L},F}^*\times i\tilde{{\cal A}}_{\tilde{L},F}^*$.  En ins\'erant cette valeur dans l'expression (20), on voit que la valeur de la forme lin\'eaire (18) sur $\underline{\varphi}'(\tilde{L}_{1},\tilde{L}_{2})$ est la valeur en $(-X_{1},X_{2})$ de la transform\'ee de Fourier de la fonction sur  sur $i\tilde{{\cal A}}_{\tilde{L}_{1},F}^*\times i\tilde{{\cal A}}_{\tilde{L}_{2},F}^*$ qui, \`a $(\tilde{\lambda}_{1},\tilde{\lambda}_{2})$, associe
 $$(21)\qquad   trace({\cal N}'(\lambda_{1},\lambda_{2})(\tilde{\Pi}_{1,\tilde{\lambda}_{1}}(f_{1})\boxtimes \tilde{\Pi}_{2,\tilde{\lambda}_{2}}(f_{2}))$$
 sur $i\tilde{{\cal A}}_{\tilde{L}_{1},F}^*\times i\tilde{{\cal A}}_{\tilde{L}_{2},F}^*$.
 
 On a montr\'e que les valeurs de la forme lin\'eaire (18) sur les deux membres de (16) \'etaient les valeurs en $(-X_{1},X_{2})$ des transform\'ees de Fourier des fonctions en $(\tilde{\lambda}_{1},\tilde{\lambda}_{2})$ d\'efinies par (19) et (21). Pour d\'emontrer leur  \'egalit\'e, il suffit de prouver l'\'egalit\'e de ces derni\`eres expressions. Il suffit encore de prouver l'\'egalit\'e
 $${\cal N}(\lambda_{1},\lambda_{2})={\cal N}'(\lambda_{1},\lambda_{2}).$$
 Mais c'est ce qu'affirme la relation 5.3(1) simplifi\'ee dans notre situation comme expliqu\'e dans ce paragraphe (et \'etendue aux familles \`a valeurs op\'erateurs, ce qui ne pose pas de difficult\'e). Cela ach\`eve la preuve. $\square$

\bigskip

\subsection{Variante avec caract\`ere central}

Pour tout sous-groupe $H$ de $G$ contenant $A_{\tilde{G}}$, on pose $\underline{H}=A_{\tilde{G}}\backslash H$. Pour tout sous-vari\'et\'e $\tilde{H}$ de $\tilde{G}$ invariante par translations par $A_{\tilde{G}}$, on pose $\underline{\tilde{H}}=A_{\tilde{G}}\backslash \tilde{H}$.  On a simplement $\underline{H}(F)=A_{\tilde{G}}(F)\backslash H(F)$, $\underline{\tilde{H}}(F)=A_{\tilde{G}}(F)\backslash \tilde{H}(F)$ puisque $A_{\tilde{G}}$ est d\'eploy\'e donc cohomologiquement trivial.
 
 Pour une fonction $f$ sur $\tilde{G}(F)$ et pour $z\in A_{\tilde{G}}(F)$, d\'efinissons la fonction $f^{[z]}$ par $f^{[z]}(\gamma)=f(z\gamma)=f(\gamma z)$. Soit $\mu$ un caract\`ere unitaire de $A_{\tilde{G}}(F)$. On note $C_{\mu}^{\infty}(\tilde{G}(F),K)$ l'espace des fonctions $f:\tilde{G}(F)\to {\mathbb C}$ qui sont lisses, $K$-finies \`a droite et \`a gauche, qui v\'erifient la relation $f^{[z]}=\mu(z)^{-1}f$ et dont le support est compact modulo $A_{\tilde{G}}(F)$. Soit $\pi$ une repr\'esentation irr\'eductible et unitaire de $G(F)$.  Pour $z\in A_{\tilde{G}}(F)$, $\pi(z)$ est une homoth\'etie dont on note $\mu_{\pi}(z)$ le rapport. L'application $z\mapsto \mu_{\pi}(z)$ est un caract\`ere unitaire.  Appelons-le le $A$-caract\`ere central de $\pi$. Soit maintenant $\tilde{\pi}$ une $\omega$-repr\'esentation unitaire de $\tilde{G}(F)$, de longueur finie. Supposons que  toutes les composantes irr\'eductibles de la repr\'esentation sous-jacente $\pi$ aient le m\^eme $A$-caract\`ere central. On dit alors que c'est le $A$-caract\`ere central de $\tilde{\pi}$. Supposons qu'il en soit ainsi et que ce caract\`ere soit $\mu$. On peut alors d\'efinir l'op\'erateur $\tilde{\pi}(f)$ pour $f\in  C_{\mu}^{\infty}(\tilde{G}(F),K)$  par
 $$\tilde{\pi}(f)=\int_{ A_{\tilde{G}}(F)\backslash \tilde{G}(F)}f(g)\tilde{\pi}(g)\,dg.$$
 On d\'efinit  aussi sa trace $I_{\tilde{\underline{G}}}(\tilde{\pi},f)=trace(\tilde{\pi}(f))$. Pour $\boldsymbol{\tau}\in {\cal E}(\tilde{G},\omega)$, la repr\'esentation $\tilde{\pi}_{\boldsymbol{\tau}}$ n'est pas $G$-irr\'eductible, mais admet un $A$-caract\`ere central, qui ne d\'epend bien s\^ur que de l'image $\tau$ de $\boldsymbol{\tau}$ dans $E(\tilde{G},\omega)$. Notons-le $\mu_{\tau}$. On note $E_{disc,\mu}(\tilde{G},\omega)$   l'ensemble des $\tau\in E_{disc}(\tilde{G},\omega)$ tels que $\mu_{\tau}=\mu$. On note $E_{disc,\mu}(\tilde{G},\omega)/conj$ l'ensemble des classes de conjugaison par $G(F)$ dans $E_{disc,\mu}(\tilde{G},\omega)$. Pour $\tau=(M,\sigma,\tilde{r})\in E_{disc,\mu}(\tilde{G},\omega)$ (ou $E_{disc,\mu}(\tilde{G},\omega)/conj$), notons $Stab(W^G,\tau)$ le stabilisateur de $\tau$ dans $W^G$, puis, comme en 2.9,
 $${\bf Stab}(W^G,\tau)=(Stab(W^G,\tau)/W^M)/W_{0}^{G}(\sigma).$$
   Pour $f_{1},f_{2}\in C_{\mu}^{\infty}(\tilde{G},K)$, posons
  $$I_{disc}^{\tilde{\underline{G}}}(\omega,f_{1},f_{2})=\sum_{\tau\in E_{disc,\mu}(\tilde{G},\omega)/conj}\vert {\bf Stab}(W^G,\tau)\vert ^{-1}\iota(\tau)\overline{I_{\tilde{\underline{G}}}(\tilde{\pi}_{\boldsymbol{\tau}},f_{1})}I_{\tilde{\underline{G}}}(\tilde{\pi}_{\boldsymbol{\tau}},f_{1}).$$
  
  Soient $\tilde{M}\in {\cal L}(\tilde{M}_{0})$ et $\tilde{S}\in T_{ell}(\tilde{M})$.   On dispose de mesures sur $S^{\theta,0}(F)$ et sur $A_{\tilde{G}}(F)$, donc aussi sur $A_{\tilde{G}}(F)\backslash S^{\theta,0}(F)=\underline{S}^{\theta,0}(F)$. Comme en 4.1, on munit  $\tilde{S}(F)/A_{\tilde{G}}(F)(1-\theta)(S(F))=\underline{\tilde{S}}(F)/(1-\theta)(\underline{S}(F))$ de la mesure telle que, pour tout $\gamma\in \tilde{S}(F)$, l'application
  $$\begin{array}{ccc} \underline{ S}^{\theta,0}(F)&\to&\underline{\tilde{S}}(F)/(1-\theta)(\underline{S}(F))\\ t&\mapsto &t\gamma\\ \end{array}$$
  pr\'eserve localement les mesures au voisinage de l'origine.  Supposons $\omega$ trivial sur $S^{\theta}(F)$.  Remarquons que $C_{\mu}^{\infty}(\tilde{G},K)\subset {\cal H}_{ac}(\tilde{G})$. Les int\'egrales pond\'er\'ees \'equivariantes de 6.5 sont donc d\'efinies pour $f_{1}$ et $f_{2}$ et $\gamma\in \tilde{S}(F)\cap \tilde{G}_{reg}(F)$, ainsi que les termes $I_{\tilde{M}}^{\tilde{G}}(\gamma,\omega,f_{1},f_{2})$ de 6.6.  On v\'erifie que
  $$I_{\tilde{M}}^{\tilde{G}}(z\gamma,\omega,f_{1},f_{2})=I_{\tilde{M}}^{\tilde{G}}(\gamma,\omega,f_{1},f_{2})$$
  pour tout $z\in A_{\tilde{G}}(F)$.   Remarquons que l'inclusion $A_{\tilde{G}}(F)\backslash S^{\theta}(F)\subset \underline{S}^{\theta}(F)$ n'est pas forc\'ement une \'egalit\'e. On pose
  $$I_{\tilde{\underline{M}}}^{\tilde{\underline{G}}}(\gamma,\omega,f_{i})=[\underline{S}^{\theta}(F):(A_{\tilde{G}}(F)\backslash S^{\theta}(F))]^{-1}I_{\tilde{M}}^{\tilde{G}}(\gamma,\omega,f_{i})$$
  pour $i=1,2$ et 
  $$I_{\tilde{\underline{M}}}^{\tilde{\underline{G}}}(\gamma,\omega,f_{1},f_{2})=[\underline{S}^{\theta}(F):(A_{\tilde{G}}(F)\backslash S^{\theta}(F))]^{-2}I_{\tilde{M}}^{\tilde{G}}(\gamma,\omega,f_{1},f_{2}).$$
  
  {\bf Remarque.} On a d\'efini les int\'egrales orbitales comme des int\'egrales sur $S^{\theta}(F)\backslash G(F)$. On doit ici les convertir en int\'egrales sur $\underline{S}^{\theta}(F)\backslash\underline{G}(F)$, ce qui justifie les facteurs introduits ci-dessus.
  \bigskip
  
    On pose, au moins formellement,
  $$(1) \qquad I_{\tilde{\underline{M}},\tilde{\underline{S}}}^{\tilde{\underline{G}}}(\omega,f_{1},f_{2})=\vert W^M(\tilde{S})\vert ^{-1} mes(A_{\underline{\tilde{M}}}(F)\backslash \underline{ S}^{\theta}(F))$$
  $$\int_{\tilde{S}(F)/A_{\tilde{G}}(F)(1-\theta)(S(F))}I_{\tilde{\underline{M}}}^{\tilde{\underline{G}}}(\gamma,\omega,f_{1},f_{2})\,d\gamma.$$
   
    Comme en 6.6, on d\'efinit ensuite
  $$I_{\tilde{\underline{M}},g\acute{e}om}^{\tilde{\underline{G}}}(\omega,f_{1},f_{2})=\sum_{\tilde{S}\in T_{ell}(\tilde{M},\omega)}I_{\tilde{\underline{M}},\tilde{\underline{S}}}^{\tilde{\underline{G}}}(\omega,f_{1},f_{2})$$
  et
  $$I_{g\acute{e}om}^{\tilde{\underline{G}}}(\omega,f_{1},f_{2})=\sum_{\tilde{M}\in {\cal L}(\tilde{M}_{0})}\vert \tilde{W}^M\vert \vert \tilde{W}^G\vert ^{-1}(-1)^{a_{\tilde{M}}-a_{\tilde{G}}}I_{\tilde{\underline{M}},g\acute{e}om}^{\tilde{\underline{G}}}(\omega,f_{1},f_{2}).$$
  
  \ass{Th\'eor\`eme}{Pour $f_{1},f_{2}\in C_{\mu}^{\infty}(\tilde{G}(F),K)$, les expressions (1) sont absolument convergentes. On a l'\'egalit\'e
  $$I_{g\acute{e}om}^{\tilde{\underline{G}}}(\omega,f_{1},f_{2})=I_{spec}^{\tilde{\underline{G}}}(\omega,f_{1},f_{2}).$$}
  
  Preuve.   On ne d\'etaillera que l'aspect combinatoire de la preuve. On  note
  $$p_{\mu}:C_{c}^{\infty}(\tilde{G}(F),K)\to C_{\mu}^{\infty}(\tilde{G}(F),K)$$
  l'application qui, \`a $\varphi\in C_{c}^{\infty}(\tilde{G}(F),K)$, associe la fonction $f$ d\'efinie par
  $$f(\gamma)=\int_{A_{\tilde{G}}(F)}f(z\gamma)\mu(z)\,dz.$$
  Elle est surjective. Ainsi on peut introduire pour $i=1,2$ une fonction $\varphi_{i}\in C_{c}^{\infty}(\tilde{G}(F),K)$ telle que $p_{\mu}(\varphi_{i})=f_{i}$. Le th\'eor\`eme 6.6 nous fournit l'\'egalit\'e
  $$I_{g\acute{e}om}^{\tilde{G}}(\omega,\varphi_{1},\varphi_{2}^{[z]})=I_{spec}^{\tilde{G}}(\omega,\varphi_{1},\varphi_{2}^{[z]})$$
  pour tout $z\in A_{\tilde{G}}(F)$. On v\'erifie que le membre de gauche (donc aussi celui de droite) est \`a support compact en $z$ et born\'e. On a donc l'\'egalit\'e
  $$(2) \qquad \int_{A_{\tilde{G}}(F)}I_{g\acute{e}om}^{\tilde{G}}(\omega,\varphi_{1},\varphi_{2}^{[z]})\mu(z)\,dz\,=\int_{A_{\tilde{G}}(F)}I_{spec}^{\tilde{G}}(\omega,\varphi_{1},\varphi_{2}^{[z]})\mu(z)\,dz.$$
  On va montrer que   le membre de gauche est \'egal \`a $I_{g\acute{e}om}^{\tilde{\underline{G}}}(\omega,f_{1},f_{2})$, tandis que celui de droite est \'egal \`a $I_{spec}^{\tilde{\underline{G}}}(\omega,f_{1},f_{2})$.
  
  Commen\c{c}ons par les termes g\'eom\'etriques. On peut fixer $\tilde{M}\in {\cal L}(\tilde{M}_{0})$ et $\tilde{S}\in T_{ell}(\tilde{M},\omega)$ et prouver l'\'egalit\'e
  $$(3) \qquad mes(A_{\tilde{M}}(F)\backslash S^{\theta}(F))\int_{A_{\tilde{G}}(F)}\int_{\tilde{S}(F)/(1-\theta)(S(F))}I_{\tilde{M}}^{\tilde{G}}(\gamma,\omega,\varphi_{1},\varphi_{2}^{[z]})\,d\gamma\,dz=$$
  $$mes(A_{\underline{\tilde{M}}}(F)\backslash \underline{S}^{\theta}(F))\int_{\underline{\tilde{S}}(F)/(1-\theta)(\underline{S}(F))}I_{\tilde{\underline{M}}}^{\tilde{\underline{G}}}(\gamma,\omega,f_{1},f_{2})\,d\gamma.$$
  Il existe une constante $C>0$ telle que l'on ait la formule d'int\'egration
  $$(4) \qquad \int_{\tilde{S}(F)/(1-\theta)(S(F))}\phi(\gamma)\,d\gamma\,=C\int_{\underline{\tilde{S}}(F)/(1-\theta)(\underline{S}(F))}\int_{A_{\tilde{G}}(F)}\phi(z'\gamma)\,dz'\,d\gamma$$
  pour toute fonction int\'egrable $\phi$ sur $\tilde{S}(F)/(1-\theta)(S(F))$. Le membre de gauche de (3) devient
  $$ C mes(A_{\tilde{M}}(F)\backslash S^{\theta}(F))\int_{\underline{\tilde{S}}(F)/(1-\theta)(\underline{S}(F))}\int_{A_{\tilde{G}}(F)}\int_{A_{\tilde{G}}(F)}I_{\tilde{M}}^{\tilde{G}}(z'\gamma,\omega,\varphi_{1},\varphi_{2}^{[z]})\,dz'\,dz\,d\gamma.$$
  On v\'erifie que
$$I_{\tilde{M}}^{\tilde{G}}(z'\gamma,\omega,\varphi_{1},\varphi_{2}^{[z]})=I_{\tilde{M}}^{\tilde{G}}(\gamma,\omega,\varphi_{1}^{[z']},\varphi_{2}^{[z'z]}),$$
puis que   
$$\int_{A_{\tilde{G}}(F)}\int_{A_{\tilde{G}}(F)}I_{\tilde{M}}^{\tilde{G}}(z'\gamma,\omega,\varphi_{1},\varphi_{2}^{[z]})\,dz'\,dz= [\underline{S}^{\theta}(F):(A_{\tilde{G}}(F)\backslash S^{\theta}(F))]^2 I_{\tilde{\underline{M}}}^{\tilde{\underline{G}}}(\gamma,\omega,f_{1},f_{2}).$$
Alors le membre de gauche de (3) devient celui de droite, \`a ceci pr\`es que la constante $mes(A_{\underline{\tilde{M}}}(F)\backslash \underline{S}^{\theta}(F))$ y est remplac\'ee par $Cmes(A_{\tilde{M}}(F)\backslash S^{\theta}(F))[\underline{S}^{\theta}(F):(A_{\tilde{G}}(F)\backslash S^{\theta}(F))]^2$. Il reste \`a montrer que ces deux termes sont \'egaux, autrement dit, on doit calculer $C$.  Fixons un point base $\gamma_{0}\in \tilde{S}(F)$, ce qui permet d'identifier $S(F)$ \`a $\tilde{S}(F)$ par $x\mapsto x\gamma_{0}$ et de m\^eme $\underline{S}(F)$ \`a $\underline{\tilde{S}}(F)$. Consid\'erons le diagramme
$$\begin{array}{ccccc}&&S^{\theta,0}(F)&&\\ &\swarrow&&\searrow&\\\underline{S}^{\theta,0}(F)&&&&S(F)/(1-\theta)(S(F))\\ &\searrow&&\swarrow&\\&&\underline{S}(F)/(1-\theta)(\underline{S}(F))&&\\ \end{array}$$
Pour chaque fl\`eche
$$\begin{array}{ccc}&&X\\&\swarrow&\\Y&&\\ \end{array}$$
de ce diagramme, on d\'efinit une application 
$$\begin{array}{ccc}C_{c}^{\infty}(X)&\to&C_{c}^{\infty}(Y)\\ \psi_{X}&\mapsto&\psi_{Y}\\ \end{array}$$
par 
$$\psi_{Y}(y)=\int_{A_{\tilde{G}}(F)}\psi_{X}(zy_{X})\,dz,$$
 o\`u $y_{X}$ est une image r\'eciproque de $y$ dans $X$. Pour chaque fl\`eche
$$\begin{array}{ccc}X&&\\&\searrow&\\&&Y\\ \end{array}$$
de ce diagramme, on d\'efinit une application 
$$\begin{array}{ccc}C_{c}^{\infty}(X)&\to&C_{c}^{\infty}(Y)\\ \psi_{X}&\mapsto&\psi_{Y}\\ \end{array}$$
par $\psi_{Y}(y)=\sum_{x}\psi_{X}(x)$, o\`u $x$ parcourt l'image r\'eciproque de $y$ dans $X$. Soit $\psi\in C_{c}^{\infty}(S^{\theta,0}(F))$. En appliquant ces constructions au c\^ot\'e gauche,  du diagramme, on d\'eduit de $\psi$ une fonction $\psi_{g}$  sur $\underline{S}(F)/(1-\theta)(\underline{S}(F))$. En appliquant ces constructions au c\^ot\'e droit du diagramme, on obtient d'abord une fonction $\phi$ sur $S(F)/(1-\theta)(S(F))$, puis une fonction $\psi_{d}$ sur $\underline{S}(F)/(1-\theta)(\underline{S}(F))$. Appliquons (4) \`a $\phi$. Par d\'efinition de la mesure sur $S(F)/(1-\theta)(S(F)) $,  on a
$$\int_{S(F)/(1-\theta)(S(F))}\phi(\gamma)\,d\gamma\,=\int_{S^{\theta,0}(F)}\psi(\gamma)\,d\gamma.$$
Par d\'efinition de la mesure sur $\underline{S}(F)/(1-\theta)(\underline{S}(F))$, c'est aussi
$$\int_{\underline{ S}(F)/(1-\theta)(\underline{S}(F))}\psi_{g}(\gamma)\,d\gamma.$$
L'int\'egrale int\'erieure du membre de droite de (4) n'est autre que $\psi_{d}(\gamma)$ et ce membre de droite est \'egal \`a
$$\int_{\underline{ S}(F)/(1-\theta)(\underline{S}(F))}\psi_{d}(\gamma)\,d\gamma.$$
Les fonctions $\psi_{d}$ et $\psi_{g}$ sont proportionnelles et les calculs ci-dessus montrent que la constante $C$ est celle pour laquelle $\psi_{g}=C\psi_{d}$. L'action de $A_{\tilde{G}}(F)$ sur $S^{\theta,0}(F)$ est libre, tandis que celle sur $S(F)/(1-\theta)(S(F))$ se quotiente en une action libre de $A_{\tilde{G}}(F)/(A_{\tilde{G}}(F)\cap (1-\theta)(S(F)))$. Il r\'esulte alors des d\'efinitions que
$C=\vert A_{\tilde{G}}(F)\cap (1-\theta)(S(F))\vert^{-1} $. Pour $x\in \underline{S}^{\theta}(F)$, soit $y\in S(F)$ relevant $x$, posons $z=(1-\theta)y$. L'application $x\mapsto z$ envoie 
$\underline{S}^{\theta}(F)$ dans $ A_{\tilde{G}}(F)\cap (1-\theta)(S(F))$. Parce que l'application $S^{\theta,0}(F)\to \underline{S}^{\theta,0}(F)$ est surjective, l'application pr\'ec\'edente  se quotiente en une application de $\underline{S}^{\theta}(F)/\underline{S}^{\theta,0}(F)$ dans  $ A_{\tilde{G}}(F)\cap (1-\theta)(S(F))$. On obtient une suite d'applications
$$1\to S^{\theta}(F)/S^{\theta,0}(F)\to \underline{S}^{\theta}(F)/\underline{S}^{\theta,0}(F)\to A_{\tilde{G}}(F)\cap (1-\theta)(S(F))\to 1.$$
On v\'erifie ais\'ement que cette suite est exacte. Donc
$$\vert A_{\tilde{G}}(F)\cap (1-\theta)(S(F))\vert=[\underline{S}^{\theta}(F):(A_{\tilde{G}}(F)\backslash S^{\theta}(F))]=[\underline{S}^{\theta}(F):\underline{S}^{\theta,0}(F)][ S^{\theta}(F):S^{\theta,0}(F)]^{-1} .$$
 D'autre part, on a les \'egalit\'es
 $$mes(A_{\underline{\tilde{M}}}(F)\backslash\underline{S}^{\theta}(F))=[\underline{S}^{\theta}(F):\underline{S}^{\theta,0}(F)]mes(A_{\underline{\tilde{M}}}(F)\backslash \underline{S}^{\theta,0}(F))$$
 $$=
[\underline{S}^{\theta}(F):\underline{S}^{\theta,0}(F)]mes(A_{\tilde{M}}(F)\backslash S^{\theta,0}(F))$$
$$
=[\underline{S}^{\theta}(F):\underline{S}^{\theta,0}(F)][ S^{\theta}(F):S^{\theta,0}(F)]^{-1} mes(A_{\tilde{M}}(F)\backslash S^{\theta}(F)).$$
De ces calculs r\'esulte l'\'egalit\'e
$$C=mes(A_{\underline{\tilde{M}}}(F)\backslash\underline{S}^{\theta}(F))mes(A_{\tilde{M}}(F)\backslash S^{\theta}(F))^{-1}[\underline{S}^{\theta}(F):(A_{\tilde{G}}(F)\backslash S^{\theta}(F))]^{-2}$$
que l'on voulait prouver. Cela ach\`eve la preuve de (3) et l'identification du membre de gauche de (2) avec $I_{g\acute{e}om}^{\tilde{\underline{G}}}(\omega,f_{1},f_{2})$.

  Montrons maintenant que le membre de droite de (2) est \'egal \`a $I_{spec}^{\tilde{\underline{G}}}(\omega,f_{1},f_{2})$. Ce membre de droite est \'egal \`a
  $$\sum_{\tau\in (E_{disc}(\tilde{G},\omega)/conj)/i{\cal A}_{\tilde{G},F}^*}\vert {\bf Stab}(W^G\times i{\cal A}_{\tilde{G},F}^*,\tau)\vert ^{-1}\iota(\tau)$$
  $$\int_{A_{\tilde{G}}(F)}\int_{i{\cal A}_{\tilde{G},F}^*}\overline{I_{\tilde{G}}(\tilde{\pi}_{\boldsymbol{\tau}_{\tilde{\lambda}}},\varphi_{1})}I_{\tilde{G}}(\tilde{\pi}_{\boldsymbol{\tau}_{\tilde{\lambda}}},\varphi_{2}^{[z]})\,d\lambda\,\mu(z)\,dz.$$
 On identifie $(E_{disc}(\tilde{G},\omega)/conj)/i{\cal A}_{\tilde{G},F}^*$ \`a un ensemble de repr\'esentants $\underline{E}$ dans $E_{disc}(\tilde{G},\omega)$. On a
 $$I_{\tilde{G}}(\tilde{\pi}_{\boldsymbol{\tau}_{\tilde{\lambda}}},\varphi_{2}^{[z]})=I_{\tilde{G}}(\tilde{\pi}_{\boldsymbol{\tau}_{\tilde{\lambda}}},\varphi_{2})\mu_{\tau}(z)^{-1}e^{-<\lambda,H_{\tilde{G}}(z)>}.$$
 L'int\'egrale sur le plus grand sous-groupe compact de $A_{\tilde{G}}(F)$ s\'electionne les $\tau\in \underline{E}$ tels que $\mu_{\tau}$ co\"{\i}ncide avec $\mu$ sur ce sous-groupe. Notons $\underline{E}'$ ce sous-ensemble. L'expression ci-dessus devient
 $$mes(A_{\tilde{G}}(F)_{c})\sum_{\tau\in \underline{E}'}\vert {\bf Stab}(W^G\times i{\cal A}_{\tilde{G},F}^*,\tau)\vert ^{-1}\iota(\tau)$$
 $$\int_{{\cal A}_{A_{\tilde{G}},F}}\int_{i{\cal A}_{\tilde{G},F}^*}\overline{I_{\tilde{G}}(\tilde{\pi}_{\boldsymbol{\tau}_{\tilde{\lambda}}},\varphi_{1})}I_{\tilde{G}}(\tilde{\pi}_{\boldsymbol{\tau}_{\tilde{\lambda}}},\varphi_{2})\,d\lambda\,(\mu_{\tau}^{-1}\mu)(H)e^{-<\lambda,H>}\,dH,$$
 o\`u on  a quotient\'e  $\mu_{\tau}^{-1}\mu$ en un caract\`ere de ${\cal A}_{A_{\tilde{G}},F}$. Par inversion de Fourier, on obtient le produit de
 $$(5)\qquad mes(A_{\tilde{G}}(F)_{c})mes(i{\cal A}_{\tilde{G},F}^*)[{\cal A}_{G,F}:{\cal A}_{A_{G},F}]^{-1}$$ et de
 $$\sum_{\tau\in \underline{E}'}\vert {\bf Stab}(W^G\times i{\cal A}_{\tilde{G},F}^*,\tau)\vert ^{-1}\iota(\tau)\sum_{\lambda\in \Lambda_{\mu}(\tau)}\overline{I_{\tilde{G}}(\tilde{\pi}_{\boldsymbol{\tau}_{\tilde{\lambda}}},\varphi_{1})}I_{\tilde{G}}(\tilde{\pi}_{\boldsymbol{\tau}_{\tilde{\lambda}}},\varphi_{2}),$$
 o\`u $\Lambda_{\mu}(\tau)$ est l'ensemble des $\lambda\in i{\cal A}_{\tilde{G},F}^*$ tels que $\mu_{\tau_{\lambda}}=\mu$, ou encore tels que $\tau_{\lambda}\in E_{disc,\mu}(\tilde{G},\omega)$.  L'expression (5) vaut $1$ d'apr\`es les normalisations de 1.2. L'ensemble $\Lambda_{\mu}(\tau)$ est un espace principal homog\`ene sous le groupe $i{\cal A}_{A_{\tilde{G}},F}^{\vee}/i{\cal A}_{\tilde{G},F}^{\vee}$. On  v\'erifie que, pour $\lambda$ dans cet ensemble, on a l'\'egalit\'e
 $$I_{\tilde{G}}(\tilde{\pi}_{\boldsymbol{\tau}_{\tilde{\lambda}}},\varphi_{i})=I_{\tilde{\underline{G}}}(\tilde{\pi}_{\boldsymbol{\tau}_{\tilde{\lambda}}},f_{i})$$
 pour $i=1,2$.  Notons $\underline{X}$ l'ensemble des couples $(\tau,\lambda)$, o\`u $\tau\in \underline{E}'$ et $\lambda\in \Lambda_{\mu}(\tau)$. Consid\'erons l'application de $ \underline{X}$ dans $ E_{disc,\mu}(\tilde{G},\omega)/conj$ qui, \`a $(\tau,\lambda)\in \underline{X}$, associe la classe de conjugaison de $\tau_{\lambda}$.
 On v\'erifie qu'elle est surjective. Notons $\underline{X}_{\tau}$ sa fibre au-dessus d'un \'el\'ement $\tau\in  E_{disc,\mu}(\tilde{G},\omega)/conj$. A ce point, on a transform\'e  le membre de droite de (2) en
 $$(6) \qquad \sum_{\tau\in E_{disc,\mu}(\tilde{G},\omega)/conj}C(\tau)\overline{I_{\tilde{\underline{G}}}(\tilde{\pi}_{\boldsymbol{\tau}},f_{1})}I_{\tilde{\underline{G}}}(\tilde{\pi}_{\boldsymbol{\tau}},f_{1}),$$
 o\`u
 $$C(\tau)=\sum_{(\underline{\tau},\lambda)\in \underline{X}_{\tau}}\vert {\bf Stab}(W^G\times i{\cal A}_{\tilde{G},F}^*,\underline{\tau})\vert ^{-1}\iota(\underline{\tau}).$$
Fixons $\tau\in  E_{disc,\mu}(\tilde{G},\omega)/conj$  que l'on rel\`eve  en un \'el\'ement de $E_{disc}(\tilde{G},\omega)$. L'ensemble $\underline{X}_{\tau}$ est celui des $(\underline{\tau},\lambda)$ tels que $\underline{\tau}\in \underline{E}'$, $\lambda\in \Lambda_{\mu}(\underline{\tau})$ et $\underline{\tau}_{\lambda}$ est conjugu\'e \`a $\tau$.  Cette relation entra\^{\i}ne que les images de $\tau$ et $\underline{\tau}$ dans   $(E_{disc}(\tilde{G},\omega)/conj)/i{\cal A}_{\tilde{G},F}^*$ sont \'egales. Donc $\underline{\tau}$ est bien d\'etermin\'e: c'est l'\'el\'ement de $\underline{E}'$ qui repr\'esente l'image de $\tau$ dans $(E_{disc}(\tilde{G},\omega)/conj)/i{\cal A}_{\tilde{G},F}^*$. Cela simplifie $C(\tau)$ en
$$C(\tau)=\vert \underline{X}_{\tau}\vert \vert {\bf Stab}(W^G\times i{\cal A}_{\tilde{G},F}^*,\tau)\vert ^{-1}\iota(\tau).$$
On peut fixer $\underline{\lambda}\in i{\cal A}_{\tilde{G},F}^*$ tel que $\underline{\tau}_{\underline{\lambda}}$ soit conjugu\'e \`a $\tau$. Notons $\Lambda'(\tau)$ l'ensemble des $\lambda\in i{\cal A}_{\tilde{G},F}^*$ tels que $\tau_{\lambda}$ soit conjugu\'e \`a $\tau$.  Alors $\underline{X}_{\tau}$ est l'ensemble des $(\underline{\tau},\underline{\lambda}+\lambda)$ pour $\lambda\in \Lambda'(\tau)$. D'o\`u $\vert \underline{X}_{\tau}\vert=\vert \Lambda'(\tau)\vert $. Remarquons  que $\tau$ et $\tau_{\lambda}$ sont conjugu\'es par $G(F)$  si et seulement s'ils le sont par $W^G$.   On a d\'efini le groupe $Stab(W^G\times i{\cal A}_{\tilde{G},F}^*,\tau)$ en 2.9, form\'e des $(w,\lambda)\in Stab(W^G\times i{\cal A}_{\tilde{G},F}^*,\tau)$ tels que $w\tau=\tau_{\lambda}$. L'application qui, \`a $(w,\lambda)$, associe $\lambda$ fournit une suite exacte
 $$1\to Stab(W^G,\tau)\to Stab(W^G\times i{\cal A}_{\tilde{G},F}^*,\tau)\to \Lambda'(\tau)\to 1.$$
 On en d\'eduit ais\'ement
 $$\vert \Lambda'(\tau)\vert =\vert {\bf Stab}(W^G\times i{\cal A}_{\tilde{G},F}^*,\tau)\vert \vert {\bf Stab}(W^G,\tau)\vert ^{-1}.$$
 D'o\`u 
 $$C(\tau)=\vert {\bf Stab}(W^G,\tau)\vert ^{-1}\iota(\tau).$$
 Alors (6) devient $I_{disc}^{\tilde{\underline{G}}}(\omega,f_{1},f_{2})$. $\square$

 \bigskip
 
 \section{Cons\'equences}
 
 \bigskip
 
 \subsection{Fonctions cuspidales}
 Rappelons que l'on note $\tilde{G}(F)_{ell}$ l'ensemble des \'el\'ements semi-simples r\'eguliers et elliptiques de $\tilde{G}(F)$. Une fonction $f\in C_{c}^{\infty}(\tilde{G}(F),K)$ est dite cuspidale si et seulement si $I_{\tilde{G}}(\gamma,\omega,f)=0$ pour tout $\gamma\in \tilde{G}_{reg}(F)$ tel que $\gamma\not\in \tilde{G}(F)_{ell}$.  Cela \'equivaut \`a ce que, pour tout espace parabolique $\tilde{P}=\tilde{M}U_{P}\in {\cal F}(\tilde{M}_{0})$, avec $\tilde{P}\not=\tilde{G}$, l'image de $f_{\tilde{P}}$ dans $I(\tilde{M}(F),K^M,\omega)$ soit nulle (en convenant que cet espace lui-m\^eme est nul si  $\omega$ n'est pas trivial sur $Z_{M}(F)^{\theta}$). Gr\^ace au th\'eor\`eme 5.5, cela \'equivaut aussi \`a $I_{\tilde{G}}(Ind_{\tilde{P}}^{\tilde{G}}(\tilde{\pi}),f)=0$ pour tout $\tilde{P}$ comme ci-dessus et toute $\omega$-repr\'esentation $\tilde{\pi}$ de $\tilde{M}(F)$,  temp\'er\'ee et de longueur finie. On note $I_{cusp}(\tilde{G}(F),K,\omega)$ l'image dans $I(\tilde{G}(F),K,\omega)$ de l'espace des fonctions cuspidales. D'apr\`es le th\'eor\`eme 6.2, de l'application $pw_{\tilde{G}}$ se d\'eduit un isomorphisme $I_{cusp}(\tilde{G}(F),K,\omega)\simeq PW_{ell}(\tilde{G},\omega)$.
 
 Dans le cas o\`u $f_{2}$ est cuspidale, le th\'eor\`eme 6.6 se simplifie. Pour $\tilde{M}\in {\cal L}(\tilde{M}_{0})$, $\tilde{S}\in T_{ell}(\tilde{M},\omega)$  et  $\gamma\in \tilde{S}(F)\cap \tilde{G}_{reg}(F)$, on a simplement
 $$I_{\tilde{M}}^{\tilde{G}}(\gamma,\omega,f_{1},f_{2})=\overline{I_{\tilde{G}}(\gamma,\omega,f_{1})}I_{\tilde{M}}^{\tilde{G}}(\gamma,\omega,f_{2}).$$
 En effet, dans la somme d\'efinissant le membre de gauche (cf. 6.6), les termes pour $\tilde{L}_{2}\not=\tilde{G}$ sont nuls parce que $f_{2,\tilde{L}_{2}}=0$. Si $\tilde{L}_{2}=\tilde{G}$, on a $\tilde{L}_{1}=\tilde{M}$ et $I_{\tilde{M}}^{\tilde{M}}(\gamma,\omega,f_{1,\tilde{M}})=I_{\tilde{G}}(\gamma,\omega,f_{1})$.  D'autre part, on a
 $$I_{disc}^{\tilde{G}}(\omega,f_{1},f_{2})=\sum_{\tau\in (E_{ell}(\tilde{G},\omega)/conj)/i{\cal A}_{\tilde{G},F}^*}\vert {\bf Stab}(W^G\times i{\cal A}_{\tilde{G},F}^*,\tau)\vert ^{-1}\iota(\tau)$$
 $$\int_{i{\cal A}_{\tilde{G},F}^*}\overline{I_{\tilde{G}}(\tilde{\pi}_{\boldsymbol{\tau}_{\tilde{\lambda}}},f_{1})}I_{\tilde{G}}(\tilde{\pi}_{\boldsymbol{\tau}_{\tilde{\lambda}}},f_{2})\,d\lambda.$$
 En effet, si $\tau\in E_{disc}(\tilde{G}(F))-E_{ell}(\tilde{G}(F))$, une repr\'esentation $\tilde{\pi}_{\boldsymbol{\tau}_{\tilde{\lambda}}}$ est induite \`a partir d'un espace parabolique propre (cf. lemmes 2.10 et 2.11), donc son caract\`ere annule $f_{2}$. Remarquons que, pour un triplet $\tau=(M_{disc},\sigma,\tilde{r})\in E_{ell}(\tilde{G}(F))$, on a simplement
 $$\iota(\tau)=\vert det((1-\tilde{r})_{\vert {\cal A}_{M}^{\tilde{G}}})\vert ^{-1},$$
 puisque $W_{0}^G(\sigma)=\{1\}$. 
 
 Rappelons que les caract\`eres de $\omega$-repr\'esentations de longueur finie de $\tilde{G}(F)$ sont des distributions  localement int\'egrables. Pour une telle $\omega$-repr\'esentation $\tilde{\pi}$, notons $\gamma\mapsto \Theta(\tilde{\pi},\gamma)$ la fonction sur $\tilde{G}(F)$ telle que
 $$I_{\tilde{G}}(\tilde{\pi},f)=\int_{\tilde{G}(F)}\Theta(\tilde{\pi},\gamma)f(\gamma)\,d\gamma.$$
 On suppose, ainsi qu'il est loisible, que cette fonction  est lisse sur $\tilde{G}_{reg}(F)$.
 
 \ass{Th\'eor\`eme}{Soient $f\in C_{c}^{\infty}(\tilde{G}(F),K)$, $\tilde{M}\in {\cal L}(\tilde{M}_{0})$ et $\gamma\in \tilde{M}(F)\cap \tilde{G}_{reg}(F)$. On suppose $f$ cuspidale. Alors
 
 (i) si $\gamma\not\in \tilde{M}(F)_{ell}$, $I_{\tilde{M}}^{\tilde{G}}(\gamma,\omega,f)=0$;
 
 (ii) si $\gamma\in \tilde{M}(F)_{ell}$,
 $$mes(A_{\tilde{M}}(F)\backslash Z_{G}(\gamma,F))I_{\tilde{M}}^{\tilde{G}}(\gamma,\omega,f)=(-1)^{a_{\tilde{M}}-a_{\tilde{G}}}D^{\tilde{G}}(\gamma)^{1/2}$$
 $$\sum_{\tau\in (E_{ell}(\tilde{G},\omega)/conj)/i{\cal A}_{\tilde{G},F}^*}\vert {\bf Stab}(W^G\times i{\cal A}_{\tilde{G},F}^*,\tau)\vert ^{-1}\iota(\tau)
\int_{i{\cal A}_{\tilde{G},F}^*}\overline{ \Theta(\tilde{\pi}_{\boldsymbol{\tau}_{\tilde{\lambda}}},\gamma)}I_{\tilde{G}}(\tilde{\pi}_{\boldsymbol{\tau}_{\tilde{\lambda}}},f)\,d\lambda.$$}

 Preuve. Si $\gamma\not\in \tilde{M}(F)_{ell}$, on peut conjuguer $\gamma$ de sorte que $\gamma\in \tilde{L}(F)$, o\`u $\tilde{L}\in {\cal L}(\tilde{M}_{0})$ et $\tilde{L}\subsetneq \tilde{M}$. La formule de descente 6.5(2) (o\`u on \'echange les r\^oles de $\tilde{M}$ et $\tilde{L}$) et la cuspidalit\'e de $f$ entra\^{\i}nent la conclusion de (i).

 D\'efinissons une fonction $\varphi$  sur $\tilde{G}_{reg}(F)$ de la fa\c{c}on suivante. Soit $\gamma\in \tilde{G}_{reg}(F)$. Si $\omega$ est non trivial sur $Z_{G}(\gamma,F)$, on pose $\varphi(\gamma)=0$. Sinon, choisissons $g\in G(F)$ et $\tilde{M}\in {\cal L}(\tilde{M}_{0})$ de sorte que $g\gamma g^{-1}\in \tilde{M}(F)_{ell}$. Posons 
 $$\varphi(\gamma)=(-1)^{a_{\tilde{M}}-a_{\tilde{G}}}mes(A_{\tilde{M}}(F)\backslash Z_{G}(g\gamma g^{-1},F))D^{\tilde{G}}(\gamma)^{1/2}\omega(g)^{-1}I_{\tilde{M}}^{\tilde{G}}(g\gamma g^{-1},\omega,f).$$
  Cette d\'efinition est loisible d'apr\`es 6.5(1). Soit $f'\in C_{c}^{\infty}(\tilde{G}(F),K)$. Comme on l'a dit ci-dessus, le terme $I^{\tilde{G}}_{g\acute{e}om}(\omega,f',f)$ se simplifie puisque $f$ est cuspidale. Gr\^ace \`a la formule de Weyl (cf. 4.1), on a simplement
 $$I_{g\acute{e}om}^{\tilde{G}}(\omega,f',f)=\int_{\tilde{G}(F)}\overline{f'(\gamma)}\varphi(\gamma)\,d\gamma.$$
 D'autre part
 $$I_{spec}^{\tilde{G}}(\omega,f',f)=\sum_{\tau\in (E_{ell}(\tilde{G},\omega)/conj)/i{\cal A}_{\tilde{G},F}^*}\vert {\bf Stab}(W^G\times i{\cal A}_{\tilde{G},F}^*,\tau)\vert ^{-1}\iota(\tau)$$
 $$\int_{i{\cal A}_{\tilde{G},F}^*}I_{\tilde{G}}(\tilde{\pi}_{\boldsymbol{\tau}_{\tilde{\lambda}}},f)\int_{\tilde{G}(F)}\overline{\Theta(\tilde{\pi}_{\boldsymbol{\tau}_{\tilde{\lambda}}},\gamma)f'(\gamma)}\,d\gamma \,d\lambda.$$
 Cette formule est absolument convergente, d'o\`u
 $$I_{spec}^{\tilde{G}}(\omega,f',f)=\int_{\tilde{G}(F)}\overline{f'(\gamma)}\varphi'(\gamma)\,d\gamma,$$
 o\`u
 $$\varphi'(\gamma)=\sum_{\tau\in (E_{ell}(\tilde{G},\omega)/conj)/i{\cal A}_{\tilde{G},F}^*}\vert {\bf Stab}(W^G\times i{\cal A}_{\tilde{G},F}^*,\tau)\vert ^{-1}\iota(\tau)\int_{i{\cal A}_{\tilde{G},F}^*}\overline{\Theta(\tilde{\pi}_{\boldsymbol{\tau}_{\tilde{\lambda}}},\gamma)}I_{\tilde{G}}(\tilde{\pi}_{\boldsymbol{\tau}_{\tilde{\lambda}}},f)\,d\lambda.$$
 Le th\'eor\`eme 6.6 entra\^{\i}ne l'\'egalit\'e
 $$\int_{\tilde{G}(F)}\overline{f'(\gamma)}\varphi(\gamma)\,d\gamma=\int_{\tilde{G}(F)}\overline{f'(\gamma)}\varphi'(\gamma)\,d\gamma.$$
 Si $F$ est non-archim\'edien, $f'$ est n'importe quel \'el\'ement de $C_{c}^{\infty}(\tilde{G}(F))$ et cette \'egalit\'e entra\^{\i}ne $\varphi(\gamma)=\varphi'(\gamma)$ pour tout $\gamma$, ce qui est l'assertion (2) du th\'eor\`eme. Si $F$ est archim\'edien, il y a un petit probl\`eme car $f'$ est suppos\'ee $K$-finie. Quelques calculs similaires \`a ceux du paragraphe 5.2  montrent que l'\'egalit\'e ci-dessus se prolonge contin\^ument \`a $f'\in C_{c}^{\infty}(\tilde{G}(F))$ tout entier. La conclusion est alors la m\^eme que dans le cas non-archim\'edien. $\square$

 {\bf Preuve de la proposition 6.5} Soient $\tilde{M}\in {\cal L}(\tilde{M}_{0})$, $\gamma\in \tilde{M}(F)\cap \tilde{G}_{reg}(F)$ et $f\in {\cal H}_{ac}(\tilde{G}(F))$ dont l'image dans $I_{ac}(\tilde{G}(F),\omega)$ soit nulle. On veut montrer que $I_{\tilde{M}}^{\tilde{G}}(\gamma,\omega,f)=0$. La relation 6.5(3) nous ram\`ene au cas o\`u $f\in C_{c}^{\infty}(\tilde{G}(F),K)$. Son image dans $I(\tilde{G}(F),K,\omega)$ est nulle, a fortiori $f$ est cuspidale. Si $\gamma$ n'est pas elliptique dans $\tilde{M}(F)$, l'assertion cherch\'ee r\'esulte du (i) du th\'eor\`eme ci-dessus. Si $\gamma$ est elliptique dans $\tilde{M}(F)$, elle r\'esulte du (ii) puisque $I_{\tilde{G}}(\tilde{\pi}_{\boldsymbol{\tau}_{\tilde{\lambda}}},f)=0$ pour tous $\tau$ et $\tilde{\lambda}$. $\square$
 
 \bigskip
 
 \subsection{Fonctions cuspidales, variante avec caract\`ere central}
 On introduit le groupe $\underline{G}$ de 6.7 et les notations aff\'erentes.
 Soit $\mu$ un caract\`ere unitaire de $A_{\tilde{G}}(F)$. On d\'efinit comme au paragraphe pr\'ec\'edent la notion de cuspidalit\'e pour une fonction $f\in C_{\mu}^{\infty}(\tilde{G}(F),K)$. On note $C_{\mu,cusp}^{\infty}(\tilde{G}(F),K)$ le sous-espace des fonctions cuspidales. On note $I_{\mu,cusp}(\tilde{G}(F),K,\omega)$ le quotient de $C_{\mu,cusp}^{\infty}(\tilde{G}(F),K)$ par le sous-espace des fonctions dont toutes les int\'egrales orbitales r\'eguli\`eres sont nulles.  Une variante avec caract\`ere central du th\'eor\`eme 6.2 conduit \`a l'assertion suivante. On d\'efinit de fa\c{c}on \'evidente l'ensemble $E_{ell,\mu}(\tilde{G},\omega)$. Pour tout $\tau\in E_{ell,\mu}(\tilde{G},\omega)/conj$, fixons un \'el\'ement $\boldsymbol{\tau}\in {\cal E}(\tilde{G},\omega)$ qui rel\`eve $\tau$.
 Consid\'erons l'application
 $$\begin{array}{ccc}C_{\mu,cusp}^{\infty}(\tilde{G}(F),K)&\to&\sum_{\tau\in E_{ell,\mu}(\tilde{G},\omega)/conj}{\mathbb C}\\f&\mapsto& \oplus_{\tau\in E_{ell,\mu}(\tilde{G},\omega)/conj}I_{\tilde{\underline{G}}}(\tilde{\pi}_{\boldsymbol{\tau}},f).\\ \end{array}$$
 Alors
 
 (1) cette application se quotiente en un isomorphisme 
  $$I_{\mu,cusp}(\tilde{G},K,\omega)\simeq \oplus_{\tau\in E_{ell,\mu}(\tilde{G},\omega)/conj}{\mathbb C}.$$
 
 Autrement dit, pour tout $\tau\in E_{ell,\mu}(\tilde{G},\omega)/conj$, il existe un unique pseudo-coefficient $f_{\tau}\in I_{\mu,cusp}(\tilde{G}(F),K,\omega)$ tel que, pour $\tau'\in E_{ell,\mu}(\tilde{G},\omega)/conj$, on ait
 $$I_{\tilde{\underline{G}}}(\tilde{\pi}_{\boldsymbol{\tau}'},f_{\tau})=\left\lbrace\begin{array}{cc}1,& \text{ si }\tau'=\tau,\\0,&\text{ si }\tau'\not=\tau.\\ \end{array}\right.$$
 Et la famille $(f_{\tau})_{\tau\in E_{ell,\mu}(\tilde{G},\omega)/conj}$ est une base de $I_{\mu,cusp}(\tilde{G}(F),K)$.

 \ass{Th\'eor\`eme}{Soient $\tau\in E_{ell,\mu}(\tilde{G},\omega)/conj$, $\tilde{M}\in {\cal L}(\tilde{M}_{0})$ et $\gamma\in \tilde{M}(F)_{ell}\cap \tilde{G}_{reg}(F)$. On a l'\'egalit\'e
 $$mes(A_{\underline{\tilde{M}}}(F)\backslash Z_{\underline{G}}(\gamma,F))I_{\tilde{\underline{M}}}^{\tilde{\underline{G}}}(\gamma,\omega,f_{\tau})=(-1)^{a_{\tilde{M}}-a_{\tilde{G}}}D^{\tilde{G}}(\gamma)^{1/2} \vert {\bf Stab}(W^G,\tau)\vert ^{-1}\iota(\tau)\overline{ \Theta(\tilde{\pi}_{\boldsymbol{\tau}},\gamma)} .$$}
 
 C'est la version avec caract\`ere central du (ii) du th\'eor\`eme pr\'ec\'edent, appliqu\'e au pseudo-coefficient $f_{\tau}$. Remarquons que, pour $\tau=(M_{disc},\sigma,\tilde{r})\in E_{ell,\mu}(\tilde{G},\omega)$, ${\bf Stab}(W^G,\tau)$ n'est autre que le commutant de $\tilde{r}$ dans $R^G(\sigma)$.

 \bigskip
 
 \subsection{Produit scalaire elliptique}
   Notons $\underline{\tilde{G}}(F)_{ell}/conj$ l'ensemble des classes de conjugaison par $G(F)$ dans $\underline{\tilde{G}}(F)_{ell} =\tilde{G}(F)_{ell}/A_{\tilde{G}}(F)$. On munit cet ensemble d'une structure de vari\'et\'e analytique sur $F$ de sorte que, pour tout $\gamma\in \underline{\tilde{G}}(F)_{ell}$, l'application qui, \`a $x\in \underline{G}_{\gamma}(F)$, associe la classe de conjugaison de $x\gamma$, se quotiente en un isomorphisme local de $\underline{G}_{\gamma}(F)$ sur $\underline{\tilde{G}}(F)_{ell}/conj$ au voisinage de $x=1$. De m\^eme, on munit $\underline{\tilde{G}}(F)_{ell}/conj$ de la mesure pour laquelle ces applications pr\'eservent la mesure au voisinage de $x=1$.  On d\'efinit  une fonction $m$ sur $\underline{\tilde{G}}(F)_{ell}/conj$ de la fa\c{c}on suivante. Pour $\gamma\in \underline{\tilde{G}}(F)$, soit $\underline{\tilde{S}}$ l'unique tore tordu maximal de $\underline{\tilde{G}}$ contenant $\gamma$. Alors
  $$m(\gamma)=mes(\underline{S}^{\theta}(F)).$$
  Pour une fonction $\varphi\in C_{c}^{\infty}(\underline{G}(F)_{ell}/conj)$, l'int\'egrale
  $$\int_{\underline{G}(F)_{ell}/conj}m(\gamma)^{-1}\varphi(\gamma)\,d\gamma$$
  ne d\'epend d'aucune mesure. Elle est \'egale \`a
  $$\sum_{S\in T_{ell}(\tilde{G})}\vert W^G(\tilde{S})\vert ^{-1}mes(\underline{S}^{\theta}(F))^{-1}\int_{\underline{\tilde{S}}(F)/(1-\theta)(\underline{S}(F))}\varphi(\gamma)\,d\gamma.$$

 Soient $\mu$ un caract\`ere unitaire de $A_{\tilde{G}}(F)$. Soient $\tilde{\pi}_{1}$ et $\tilde{\pi}_{2}$ deux $\omega$-repr\'esentations de $\tilde{G}(F)$, de longueur finie et de $A$-caract\`ere central $\mu$. La fonction $\gamma\mapsto \overline{\Theta(\tilde{\pi}_{1},\gamma)}\Theta(\tilde{\pi}_{2},\gamma)$ sur $\tilde{G}(F)_{ell}$ se quotiente en une fonction sur $\underline{\tilde{G}}(F)_{ell}/conj$. On pose
 $$(\tilde{\pi}_{1},\tilde{\pi}_{2})_{ell}=\int_{\underline{\tilde{G}}(F)_{ell}/conj} m(\gamma)^{-1}D^{\tilde{G}}(\gamma)\overline{\Theta(\tilde{\pi}_{1},\gamma)}\Theta(\tilde{\pi}_{2},\gamma)\,d\gamma.$$
 Cette expression ne d\'epend d'aucune mesure. 
 
 Rappelons que, pour tout $\tau\in E_{ell,\mu}(\tilde{G},\omega)$, on a choisi un rel\`evement $\boldsymbol{\tau}$ de $\tau$ dans ${\cal E}_{ell}(\tilde{G},\omega)$. 
 
 \ass{Th\'eor\`eme}{(i) Soient $\tau_{1},\tau_{2}\in {\cal E}_{ell,\mu}(\tilde{G},\omega)$.   Supposons qu'ils  ne sont pas conjugu\'es par $G(F)$. Alors on a l'\'egalit\'e $(\tilde{\pi}_{\boldsymbol{\tau}_{1}},\tilde{\pi}_{\boldsymbol{\tau}_{2}})_{ell}=0$. 
 
 (ii) Soit $\tau\in {\cal E}_{ell,\mu}(\tilde{G},\omega)$. On a l'\'egalit\'e
 $$(\tilde{\pi}_{\boldsymbol{\tau}},\tilde{\pi}_{\boldsymbol{\tau}})_{ell}=\vert {\bf Stab}(W^G,\tau)\vert \iota(\tau)^{-1}.$$}

 Preuve. En appliquant le th\'eor\`eme 7.2, on voit que
 $$(\tilde{\pi}_{\boldsymbol{\tau}_{1}},\tilde{\pi}_{\boldsymbol{\tau}_{2}})_{ell}=\vert {\bf Stab}(W^G,\tau_{1})\vert \vert {\bf Stab}(W^G,\tau_{2})\vert \iota(\tau_{1})^{-1}\iota(\tau_{2})^{-1} X,$$
 o\`u
 $$X=\int_{\underline{\tilde{G}}(F)_{ell}/conj} m(\gamma)\overline{I_{\tilde{\underline{G}}}(\gamma,\omega,f_{\tau_{1}})}I_{\tilde{\underline{G}}}(\gamma,\omega,f_{\tau_{2}})\,d\gamma.$$
 Il r\'esulte des d\'efinitions que $X=I_{\tilde{\underline{G}},g\acute{e}om}^{\tilde{\underline{G}}}(\omega,f_{\tau_{1}},f_{\tau_{2}})$. Par le m\^eme argument qu'en 7.1 et parce que les fonctions $f_{\tau_{1}}$ et $f_{\tau_{2}}$ sont toutes deux cuspidales, les termes $I_{\tilde{\underline{M}},g\acute{e}om}^{\tilde{\underline{G}}}(\omega,f_{\tau_{1}},f_{\tau_{2}})$ sont nuls pour $\tilde{M}\not=\tilde{G}$. Donc $X=I^{\tilde{\underline{G}}}_{g\acute{e}om}(\omega,f_{\tau_{1}},f_{\tau_{2}})$. Appliquons le th\'eor\`eme 6.7: $X=I^{\tilde{\underline{G}}}_{disc}(\omega,f_{\tau_{1}},f_{\tau_{2}})$. Mais cette expression se calcule gr\^ace \`a la d\'efinition des pseudo-coefficients. On obtient $I^{\tilde{\underline{G}}}_{disc}(\omega,f_{\tau_{1}},f_{\tau_{2}})=0$
 si $\tau_{1}$ et $\tau_{2}$ ne sont pas conjugu\'es. Si $\tau_{1}=\tau_{2}=\tau$, on obtient
 $$I^{\tilde{\underline{G}}}_{disc}(\omega,f_{\tau},f_{\tau})=\vert {\bf Stab}(W^G,\tau)\vert^{-1} \iota(\tau). $$
  Le th\'eor\`eme en r\'esulte. $\square$
 
 \bigskip
 
 \subsection{Produit elliptique pour les $\omega$-repr\'esentations irr\'eductibles}
Soit $\mu$ un caract\`ere unitaire de $A_{\tilde{G}}(F)$. Pour $i=1,2$, soient $M_{i}\in {\cal P}(M_{0})$ et $\sigma_{i}$ une repr\'esentation irr\'eductible et de la s\'erie discr\`ete de $M_{i}(F)$. On suppose que la restriction \`a $A_{\tilde{G}}(F)$ du caract\`ere central de $\sigma_{i}$ est \'egale \`a $\mu$. On suppose ${\cal N}^{\tilde{G}}(\sigma_{i})$ non vide. Soit $\tilde{\rho}_{i}$ une repr\'esentation ${\cal R}^G(\sigma_{i})$-irr\'eductible de ${\cal R}^{\tilde{G}}(\sigma_{i})$.  D\'efinissons un terme $(\tilde{\rho}_{1},\tilde{\rho}_{2})_{ell}$ de la fa\c{c}on suivante. Supposons d'abord que les couples $(M_{1},\sigma_{1})$ et $(M_{2},\sigma_{2})$ soient \'egaux. On les note simplement $(M,\sigma)$. Supposons de plus $W^G_{0}(\sigma)=\{1\}$. La fonction $\boldsymbol{\tilde{r}}\mapsto \overline{trace(\tilde{\rho}_{1}(\boldsymbol{\tilde{r}}))}trace(\tilde{\rho}_{2}(\boldsymbol{\tilde{r}}))$ sur ${\cal R}^{\tilde{G}}(\sigma)$ se descend en une fonction sur $R^{\tilde{G}}(\sigma)$. On pose
$$(\tilde{\rho}_{1},\tilde{\rho}_{2})_{ell}=\vert R^G(\sigma)\vert ^{-1}\sum_{\tilde{r}\in R^{\tilde{G}}(\sigma)\cap W^{\tilde{G}}_{reg}(\sigma)}\vert det((1-\tilde{r})_{\vert {\cal A}_{M}^{\tilde{G}}})\vert  trace(\tilde{\rho}_{1}(\tilde{r}))\overline{trace(\tilde{\rho}_{2}(\tilde{r}))}.$$
Supposons maintenant que  $W^G_{0}(\sigma_{i})=\{1\}$ pour $i=1,2$ et que
les couples $(M_{1},\sigma_{1})$ et $(M_{2},\sigma_{2})$ sont conjugu\'es par un \'el\'ement de $G(F)$. En effectuant une telle conjugaison, on remplace $(M_{2},\sigma_{2},\tilde{\rho}_{2})$ par $(M_{1},\sigma_{1},\tilde{\rho}'_{2})$. On pose 
$$(\tilde{\rho}_{1},\tilde{\rho}_{2})_{ell}=(\tilde{\rho}_{1},\tilde{\rho}'_{2})_{ell}.$$
Cela ne d\'epend pas de la conjugaison choisie. Dans les cas restants, c'est-\`a-dire si  $W^G_{0}(\sigma_{1})$ ou $W^G_{0}(\sigma_{2})$ est non trivial, ou si $(M_{1},\sigma_{1})$ et $(M_{2},\sigma_{2})$ ne sont pas conjugu\'es, on pose $(\tilde{\rho}_{1},\tilde{\rho}_{2})_{ell}=0$.  

Pour $i=1,2$, on a associ\'e en 2.8 \`a $(M_{i},\sigma_{i},\tilde{\rho}_{i})$ une repr\'esentation $G$-irr\'eductible $\tilde{\pi}_{\tilde{\rho}_{i}}$ de $\tilde{G}(F)$. 

\ass{Corollaire}{On a l'\'egalit\'e $(\tilde{\pi}_{\tilde{\rho}_{1}},\tilde{\pi}_{\tilde{\rho}_{2}})_{ell}=(\tilde{\rho}_{1},\tilde{\rho}_{2})_{ell}$.}

Preuve. Pour $i=1,2$ et pour $\tilde{r}\in R^{\tilde{G}}(\sigma_{i})$, posons $\tau_{i}(\tilde{r})=(M_{i},\sigma_{i},\tilde{r})$. Notons $R_{ess}^{\tilde{G}}(\sigma_{i})$ l'ensemble des $\tilde{r}\in R^{\tilde{G}}(\sigma_{i})$ tels que  $\tau_{i}(\tilde{r})$ soit essentiel.   Pour $\tilde{r}\in R_{ess}^{\tilde{G}}(\sigma_{i})$, on fixe un rel\`evement $\boldsymbol{\tau}_{i}(\tilde{r})=(M_{i},\sigma_{i},\boldsymbol{\tilde{r}})$ de $\tau_{i}$ dans ${\cal E}(\tilde{G},\omega)$. On a vu dans la preuve de la proposition 2.9 une formule d'inversion exprimant $\Theta(\tilde{\pi}_{\tilde{\rho}_{i}})$ en fonction des $\Theta(\tilde{\pi}_{\boldsymbol{\tau}_{i}(\tilde{r})})$ pour $\tilde{r}\in R^{\tilde{G}}_{ess}(\sigma_{i})$. Dans cette formule, on sommait sur les classes de conjugaison dans cet ensemble. On peut la r\'ecrire comme une somme sur tous les \'el\'ements de cet ensemble sous la forme
$$\Theta(\tilde{\pi}_{\tilde{\rho}_{i}})=\vert R^G(\sigma_{i})\vert ^{-1}\sum_{\tilde{r}\in R^{\tilde{G}}_{ess}(\sigma_{i})} \overline{trace(\tilde{\rho}_{i}(\boldsymbol{\tilde{r}}))}\Theta(\tilde{\pi}_{\boldsymbol{\tau}_{i}(\tilde{r})}).$$
Si $W_{0}^G(\sigma_{i})\not=\{1\}$, il r\'esulte des lemmes 2.10 et 2.11 que toutes les repr\'esentations $\tilde{\pi}_{\boldsymbol{\tau}_{i}(\tilde{r})}$ sont des induites \`a partir d'espaces paraboliques propres. Leur caract\`ere est donc nul sur $\tilde{G}(F)_{ell}$. Il en r\'esulte que, si $W^G_{0}(\sigma_{1})$ ou $W^G_{0}(\sigma_{2})$ est non trivial, $(\tilde{\pi}_{\tilde{\rho}_{1}},\tilde{\pi}_{\tilde{\rho}_{2}})_{ell}=0$. D'o\`u l'\'egalit\'e de l'\'enonc\'e dans ce cas. On suppose maintenant $W^G_{0}(\sigma_{1})=\{1\}$ et $W^G_{0}(\sigma_{2})=\{1\}$. On a alors $R^{\tilde{G}}(\sigma_{i})=W^{\tilde{G}}(\sigma_{i})$ pour $i=1,2$. Toujours d'apr\`es les lemmes 2.10 et 2.11, les repr\'esentations $\tilde{\pi}_{\boldsymbol{\tau}_{i}(\tilde{r})}$ pour $\tilde{r}\not\in W^{\tilde{G}}_{reg}(\sigma_{i})$ ont un caract\`ere nul sur $\tilde{G}(F)_{ell}$. On obtient l'\'egalit\'e
$$(\tilde{\pi}_{\tilde{\rho}_{1}},\tilde{\pi}_{\tilde{\rho}_{2}})_{ell}=\vert R^G(\sigma_{1})\vert ^{-1}\vert R^G(\sigma_{2})\vert ^{-1}\sum_{\tilde{r}_{1}\in R^{\tilde{G}}_{ess}(\sigma_{1})\cap W^{\tilde{G}}_{reg}(\sigma_{1})}\sum_{\tilde{r}_{2}\in R^{\tilde{G}}_{ess}(\sigma_{2})\cap W^{\tilde{G}}_{reg}(\sigma_{2})} $$
$$trace(\tilde{\rho}_{1}(\boldsymbol{\tilde{r}_{1}}))\overline{trace(\tilde{\rho}_{2}(\boldsymbol{\tilde{r}_{2}}))}(\tilde{\pi}_{\boldsymbol{\tau}_{1}(\tilde{r}_{1})},\tilde{\pi}_{\boldsymbol{\tau}_{2}(\tilde{r}_{2})})_{ell}.$$
Si $(M_{1},\sigma_{1})$ et $(M_{2},\sigma_{2})$ ne sont pas conjugu\'es, le (i) du th\'eor\`eme 7.3 dit que le membre de droite ci-dessus est nul, d'o\`u encore l'\'egalit\'e de l'\'enonc\'e. Si les deux couples sont conjugu\'es, on peut les supposer \'egaux, on les note simplement $(M,\sigma)$. Le (i) du th\'eor\`eme 7.3 dit que seuls les couples $(\tilde{r}_{1},\tilde{r}_{2})$ form\'es d'\'el\'ements conjugu\'es contribuent \`a la formule ci-dessus. Pour un couple d'\'el\'ements conjugu\'es, les valeurs des diff\'erents termes intervenant sont les m\^emes pour $\tilde{r}_{1}$ et $\tilde{r}_{2}$. On a d\'ej\`a dit que, pour un triplet elliptique $\tau=(M,\sigma,\tilde{r})$, on a l'\'egalit\'e ${\bf Stab}(W^G,\tau)=Stab(R^G(\sigma),\tilde{r})$, ce dernier groupe \'etant le commutant de $\tilde{r}$ dans $R^G(\sigma)$.  On a aussi $\iota(\tau)^{-1}=\vert det((1-\tilde{r})_{\vert {\cal A}_{M}^{\tilde{G}}}\vert $. Gr\^ace au (ii) du th\'eor\`eme 7.3, la contribution d'un couple   $(\tilde{r}_{1},\tilde{r}_{2})$ d'\'el\'ements conjugu\'es \`a la formule ci-dessus est donc
$$\vert Stab(R^G(\sigma),\tilde{r}_{1})\vert \vert det((1-\tilde{r}_{1})_{\vert {\cal A}_{M}^{\tilde{G}}}\vert trace(\tilde{\rho}_{1}(\boldsymbol{\tilde{r}_{1}}))\overline{trace(\tilde{\rho}_{2}(\boldsymbol{\tilde{r}_{1}}))}.$$
Pour $\tilde{r}_{1}$ fix\'e, le nombre de $\tilde{r}_{2}$ conjugu\'es \`a $\tilde{r}_{1}$ est $\vert R^G(\sigma)\vert \vert Stab(R^G(\sigma),\tilde{r}_{1})\vert^{-1}$. Notre formule devient
$$(\tilde{\pi}_{\tilde{\rho}_{1}},\tilde{\pi}_{\tilde{\rho}_{2}})_{ell}=\vert R^G(\sigma)\vert ^{-1}\sum_{\tilde{r}\in R^{\tilde{G}}_{ess}(\sigma)\cap W^{\tilde{G}}_{reg}(\sigma)}\vert det((1-\tilde{r})_{\vert {\cal A}_{M}^{\tilde{G}}}\vert trace(\tilde{\rho}_{1}(\boldsymbol{\tilde{r}}))\overline{trace(\tilde{\rho}_{2}(\boldsymbol{\tilde{r}}))}.$$
C'est presque la d\'efinition de $(\tilde{\rho}_{1},\tilde{\rho}_{2})_{ell}$, \`a ceci pr\`es qu'ici, la somme est limit\'ee \`a $\tilde{r}\in R^{\tilde{G}}_{ess}(\sigma)\cap W^{\tilde{G}}_{reg}(\sigma)$. Mais, si $\tilde{r}\not\in R^{\tilde{G}}_{ess}(\sigma)$, les termes  $ trace(\tilde{\rho}_{1}(\boldsymbol{\tilde{r}}))$ et $ trace(\tilde{\rho}_{2}(\boldsymbol{\tilde{r}}))$ sont nuls. On peut donc remplacer la somme sur $R^{\tilde{G}}_{ess}(\sigma)\cap W^{\tilde{G}}_{reg}(\sigma)$ par celle sur tout $W^{\tilde{G}}_{reg}(\sigma)$. Cela ach\`eve la preuve. $\square$
 
 \bigskip
 
 {\bf Bibliographie}

[A1]  J. Arthur: {\it A local trace formula}, Publ. Math. de l'IHES 73 (1991), p.5-96 

[A2]  -----------: {\it The trace formula in invariant form}, Annals of Math. 114 (1981), p.1-74

[A3]  -----------: {\it Intertwining operators and residues I. Weighted characters}, J. of Functional Analysis 84 (1989), p.19-84

[A4] -----------: {\it Canonical normalization of weighted characters and a transfer conjecture}, C. R. Math. Rep. Acad. Sci. Canada 20 (1998), p.35-52

[A5] -----------: {\it On the Fourier transforms of weighted orbital integrals}, J. reine angew. Math. 452 (1994), p.163-217

[A6] -----------: {\it The invariant trace formula I. Local theory}, J. AMS 1 (1988), p.323-383

[A7] -----------: {\it On elliptic tempered characters}, Acta Math. 171 (1993), p.73-138

[A8] -----------:{\it The trace Paley-Wiener theorem for Schwartz functions}, Contemporary Math. 177 (1994), p. 171-180

[BT] F. Bruhat, J. Tits: {\it Groupes r\'eductifs sur un corps local I. Donn\'ees radicielles valu\'ees}, Publ. Math. IHES 41 (1972), p.5-251

[B] A. Bouaziz: {\it Sur les caract\`eres des groupes de Lie r\'eductifs non connexes}, J. of Functional Analysis 70 (1987), p.1-79

[C] L. Clozel: {\it Characters of non-connected, reductive $p$-adic groups} p.149-

[DM] P. Delorme, P. Mezo: {\it A twisted invariant Paley-Wiener theorem for real reductive groups}, Duke Math. J. 144 (2008), p.341-380

[HC] Harish-Chandra: {\it Spherical functions on a semi-simple Lie group I}, Amer. J. of Math. 80 (1958), p.241-316

[HL] G. Henniart, B. Lemaire: {\it La transform\'ee de Fourier tordue pour les groupes r\'eductifs $p$-adiques}, en pr\'eparation

[LW] J.-P. Labesse, J.-L. Waldspurger: {\it La formule des traces tordue d'apr\`es le Friday morning seminar}, pr\'epublication 2012, Arxiv RT 12042888

[L] B. Lemaire: {\it Caract\`eres tordus des repr\'esentations admissibles}, pr\'epublication 2010.

[R] J. Rogawski: {\it The trace Paley-Wiener theorem in the twisted case}, Trans. Amer. Math. Soc. 309 (1988), p.215-229

[W] J.-L. Waldspurger: {\it La formule de Plancherel pour les groupes $p$-adiques d'apr\`es Harish-Chandra}, Journal de l'Inst. de Math. Jussieu 2 (2003), p.235-333

\bigskip

\noindent Institut de math\'ematiques de Jussieu- CNRS

\noindent 2 place Jussieu, 75005 Paris

\noindent e-mail: waldspur@math.jussieu.fr

 \end{document}